\documentclass[a4paper]{article}

\usepackage{amssymb}
\usepackage{amsmath,amsthm}      
\usepackage{graphicx}            
\usepackage{microtype}           
\usepackage{booktabs}            
\usepackage{multicol}            
\usepackage[hidelinks]{hyperref} 
\usepackage{physics}
\usepackage{mathtools}
\usepackage{relsize}
\usepackage{setspace}

\theoremstyle{definition}
\newtheorem{example}{Example}
\newtheorem{definition}{Definition}

\newcommand{\obs}{\mathcal{O}}

\newcommand{\lbl}[1]{\label{eq: #1}}
\newcommand{\rf}[1]{(\ref{eq: #1})}

\newcommand{\dt}[1]{\frac{d #1}{dt}}
\newcommand{\dst}[1]{\frac{d^2 #1}{dt^2}}

\newcommand{\dx}[1]{\frac{d #1}{dx}}
\newcommand{\dsx}[1]{\frac{d^2 #1}{dx^2}}

\newcommand{\prt}[1]{\partial_{#1}}
\newcommand{\inv}[1]{\frac{1}{#1}}

\newcommand{\ttx}[1]{\textit{#1}}

\newcommand{\intx}{\int_{x_0}^{x_1}}
\newcommand{\intxi}{\int_{x_0}^{\xi}}
\newcommand{\intix}{\int_{\xi}^{x_1}}
\newcommand{\intt}{\int^{t_1}_{t_0}}
\newcommand{\intinf}{\int^{\infty}_{-\infty}}

\newcommand{\isp}{\frac{1}{\sqrt{2 \pi}}}
\newcommand{\eps}{\epsilon}
\newcommand{\ra}{\rightarrow}
\newcommand{\Ra}{\Rightarrow}

\newcommand{\Lra}{\Leftrightarrow}
\newcommand{\tp}[1]{10^{#1}}
\newcommand{\vg}{\vb{\gamma}}

\newcommand{\vx}{\vb{x}}
\newcommand{\vxuv}{\vb{x(u,v)}}
\newcommand{\vT}{\vb{T}}
\newcommand{\vTu}{\vb{T}_u}
\newcommand{\vTv}{\vb{T}_v}
\newcommand{\prtdvxu}{\frac{\partial\vb{x}}{\partial{} u}}
\newcommand{\prtdvxv}{\frac{\partial\vb{x}}{\partial{} v}}
\newcommand{\prtdLdy}{\frac{\partial \mathcal{L} }{\partial y}}
\newcommand{\prtdLdyd}{\frac{\partial \mathcal{L} }{\partial y'}}
\newcommand{\prtdLy}{\frac{\partial L }{\partial y}}
\newcommand{\prtdLyd}{\frac{\partial L }{\partial y'}}
\newcommand{\intRn}{\int_{\mathbf{R}^n}}
\newcommand{\mathsym}[1]{{}}
\newcommand{\unicode}[1]{{}}
\newcommand{\ba}{\begin{align}}
\newcommand{\ea}{\end{align}}

\begin{document}

\title{Topics in Applied Mathematics and Nonlinear Waves}
\author{ Per Kristen Jakobsen\\
\\
Department of Mathematics and Statistics,\\
 the Arctic University of Norway, Troms\o, Norway}
\maketitle
\tableofcontents

\section{Introduction}
Applied mathematics, in the widest sense, is possibly as old as humanity itself. In this widest sense, applied mathematics denote an activity where one classify objects and events in external reality  using symbols of some sort, and subsequently, through manipulations of the symbols, try to predict, influence or control events in the same external reality. It can be argued that humanity's ability to represent elements of reality using symbols  is the defining feature of the cognitive explosion that our species underwent between 80 and 60 thousand years ago, and which set us on the path to become the dominating species on our planet. Most of this early use of symbols by our species has been lost to history. The oldest evidence we have for our species practicing applied mathematics,   is the Ishango bone, which is a tally stick from central Africa and which may be as old as 35000 years. The oldest written accounts of applied mathematics, is the Rhind Mathematical Papyrus, which  dates to around 4000 years before the present.

  The well known split between pure and applied mathematics arguably occurred around 2600 years ago through the work originating from the Pythagorean school. However, even if applied mathematics predates pure mathematics, and pure mathematics originated from applied mathematics, today,  pure mathematics encompass a vast domain of human thought that is deep, subtle, important,  and whose continuing evolution is driven by it's own internal motivations and aspirations. Together, pure and applied mathematics has weaved  the  fabric underlying human civilization.

No written account of applied mathematics in this widest sense can ever be produced. Applied mathematics is too vast and varied for this to be possible, and the current text, or any written text for that matter, cannot  claim to be defining  what applied mathematics is. At best, a written text  in applied mathematics, like the current one, can only give an account of a very small corner of the vast tapestry that is applied mathematics.

 Given this, the best I can hope for is that the text contained in the following pages,  just perhaps, has a focus, and a level of generality, which makes it a worthwhile study.
The selection of topics in the text has formed the core of a one semester course in applied mathematics at the Arctic University of Norway that has been running continuously since the 1970s. The class has, during it's existence, drawn participants from both applied mathematics and physics, and also to some extent from pure mathematics, analysis in particular. The material in these lecture notes can be covered by one semester's worth of five lecture hours a week. The work requirements for the students consists of seven  obligatory projects whose content are taken from exercises and computational projects  included in these lecture notes.

\section{Cartesian tensors}

Physical quantities, as measured by some observer, in general consist of indexed collections of components. Thus, if $t$ is a physical quantity, its components are 
\begin{align*}
t = ( t_{i_1 i_2 ...i_p}),
\end{align*}	
where each index runs over 1, 2 ..,$n$. The number $n$ can be any value, but for the case of interest for us, $n=3$. This is related to the fact that physical space has three dimensions. \\
The components, $t_{i_1 .. i_p}$, can be numbers, vectors, matrices, operators, etc,  depending on the context. For us, they will mainly be functions
\begin{align*}
t_{i_1 i_2 ...i_p} =  t_{i_1 i_2 ...i_p}(\vb{x}).
\end{align*}
The important point is that if another observer describes the same physical quantity, her components might be entirely different. The obvious example of this  is the velocity of a particle. Two different observers moving with respect to one another will observe different velocity components even though they are describing the same physical quantity. The same is true for two different observers whose frames of reference are rotated with respect to one another. 

Let us denote the frame of reference for a given observer $\mathcal{O}$ by $x_i$, and  let $t$ be a physical quantity whose components with respect to two observers $\mathcal{O}$ and $\mathcal{O'}$ are 
\begin{align*}
(t_{i_1 .. i_p}) , (t'_{i_1 .. i_p}).
\end{align*}
It is obviously of interest to know what  the relation between these two sets of components is. This must clearly depend on the relation between the frames of reference $x_i$ and $x'_i$ for $\obs$ and $\obs'$. In these notes, a frame or reference will be a unique labeling of points in physical space using three numbers $x_{i}, \; i=1,2,3$.  Many such frames are in use by physicists and applied mathematicians. Cartesian coordinates, spherical coordinates and cylindrical coordinates are three such frames. Since the frames $x_i$ and $x_i'$ label the \ttx{same} set of points we must have 
\begin{align*}
x_i' = x_i' (x_1,x_2,x_3) && i = 1,2,3.
\end{align*}
In these notes we will only discuss Cartesian frames related by a rotation. Thus we will always have 
\begin{align}
x_i' = a_{ij} \; x_j  \; \; \; \; \; \text{[Einstein summation convention in use!]},\label{eq: 1.12}
\end{align}
where $a_{ij}$ is a matrix representing rotation about some axis in space. Recall that such matrices were called \ttx{unitary} in linear algebra. Unitary matrices have the interesting property that their inverse can be found by taking their transpose. Let $b_{ij}$ be the components of the inverse of a unitary matrix whose components are $a_{ij}$. Then, using the Kronecker delta we have 
\begin{align}
b_{ik} \; a_{kj} = \delta_{ij}\nonumber,\\
a_{ik} \; b_{kj} = \delta_{ij}.\lbl{6.12}
\end{align}
The fact the  of $a_{ij}$ is unitary  is expressed by the relation
\begin{align}
b_{ij} = a_{ji}, \lbl{7.12}
\end{align}
and using \rf{7.12} in \rf{6.12} we find the identity 
\begin{align}
a_{ki} \; a_{kj} = \delta_{ij}\nonumber,\\
a_{ik} \; a_{jk} = \delta_{ij}.\lbl{6.1}
\end{align}
Let now $x$ and $x'$ be two Cartesian frames related by a rotation matrix whose unitary matrix has components $a_{ij}$. Thus we have 
\begin{align*}
x_i' = a_{ij} \; x_j.
\end{align*}
Let $t$ be a physical quantity whose components with respect to the frames $x$ and $x'$ are 
\begin{align*}
t_{i_1 .. i_p} (x) , t'_{i_1 .. i_p} (x').
\end{align*}	
Then $t$ is a \textit{tensor field} if 
\begin{align}
t'_{i_1 .. i_p} (x') = a_{i_1 j_1} .. a_{i_p j_p } t_{j_1 .. j_p} (x). \label{eq: 11.12} 
\end{align}
The integer $p$ is the \textit{rank} of the tensor field. Note that if a tensor field $t_i$ has only one index  it is said to be a {\it vector field} and if a tensor field has no indices, and is thus a single function, it is called a {\it scalar field}.

 In order to be precise,  $t=(t_{i_1 ,..i_n}(x))$ should be called a {\it Cartesian} tensor field since tensor fields can also be defined in terms of more general change of coordinates than rotations. In general, the frames can even mix space and time. This is what occurs in relativity theory, where the frames are related through a \textit{Lorentz transformation}. The resulting tensors have indices running over $ \mu = 0, 1,2,3$, where the index value zero refers to the time coordinate. Such space time tensors are called \textit{Lorentzian tensors}. Tensors can also be defined for spaces of dimension higher than 3. In fact, in the current hottest contender for the "Theory of Everything", a space with 10 spatial dimensions is involved! \\
Finally, you should be aware of the fact that in modern physics there are important physical quantities that are \textit{not} tensors, meaning that they do not transform like in equation (\ref{eq: 11.12}) when changing coordinate frames. \textit{Connections} and \textit{Spinors} are two such non-tensorial quantities.
However, here we will stick with Cartesian tensor fields in physical space.

  Our aim is  not to give a comprehensive introduction to all important aspects of Cartesian tensors, but rather to introduce enough of the machinery to be able to use these tensors as a tool for deriving all sorts of useful vector calculus identities.

  Note that if the components of a tensor field are constant in one Cartesian frame of reference they will be constant in all frames, and we have 
\begin{align*}
t'_{i_1 .. i_p} = a_{i_1 j_1} .. a_{i_p j_p } t_{j_1 .. j_p}.
\end{align*}
This quantity we will usually call a \textit{tensor}, not a tensor field.

Even if we have a way of defining a set of components, $t_{i_1 .. i_p}$,  with respect to all Cartesian frames, these components might still not define a tensor. In order to verify that we have defined a tensor we must change coordinates according to equation (\ref{eq: 1.12}) and verify that relations (\ref{eq: 11.12}) hold. 

\subsection{Tensors}
\paragraph{The Kronecker delta}
In any given Cartesian frame we can define a collection of components $\delta_{ij}$ by 
\begin{align}\label{eq: 13.12}
\delta_{ij} = \begin{cases}
1 \; \; \; \; \;  & i = j  \\
0  & i \ne j.
\end{cases}
\end{align}
Observe that this collection of components is the \textit{same} in all Cartesian frames, there is nothing in definition (\ref{eq: 13.12}) that refers directly to the frame. Thus 
\begin{align*}
\delta'_{ij} = \delta_{ij}, 
\end{align*}
But if $\delta'_{ij}$ is to define a tensor $\delta$, we must have 
\begin{align*}
\delta'_{ij} = a_{ik} \; a_{jl} \; \delta_{kl}.
\end{align*}
We thus get the condition
\begin{align}
a_{ik} \; a_{jl} \; \delta_{kl} &= \delta_{ij},\nonumber \\
&\Updownarrow \nonumber \\  a_{ik} \; a_{jk} &= \delta_{ij}. \lbl{17.12}
\end{align}
But \rf{17.12} holds because $a_{ij}$ are the matrix elements  of a unitary matrix. Thus $\delta_{ij}$ are the components of a tensor. This tensor is called the \ttx{Kronecker delta}. 

\paragraph{A three component quantity that is not a vector}
In any Cartesian frame define a collection of components $t_i$ by 
\begin{align}
t_1 = t_2 = t_3 = 1. \lbl{18.13}
\end{align}
Let us check if \rf{18.13} defines a tensor. Observe that, as in example 1, the collection of components is the same in all Cartesian frames 
\begin{align*}
t'_i = t_i.
\end{align*}
If $t_i$ is to define a tensor we must have 
\begin{align}
t'_i = a_{ij} \; t_j.\lbl{18.15}
\end{align}
We thus have the condition
\begin{align*}
a_{ij} \; t_j = t_i.
\end{align*}
If we recall that all unitary matrices are rotations around some axis in space, it is evident that \rf{18.15} is only true for the subset of matrices that defines rotations around the axis determined by the given vector $t_i$. Therefore, $t_i$ does not define a tensor.

\paragraph{The Levi-Civita tensor}
In any Cartesian frame there is defined a collection of components $\epsilon_{ijk}$ by 
\begin{align}
\epsilon_{ijk} =\begin{cases}
\; 1 \; \; \; \; \; \text{if the permutation (1,2,3) $\rightarrow$ (i,j,k) is even}\\
-1 \; \; \; \; \text{if the permutation (1,2,3) $\rightarrow$ (i,j,k) is odd}\\
\; 0 \; \; \; \; \; \text{if two indices are equal} \lbl{25.12}
\end{cases}.
\end{align}
Recall that a permutation is even if it is composed of an even number of pair switches, and odd if it is composed of an odd number of pair switches.

 \noindent  Since the definition \rf{25.12} does not refer to the frame, we have that $\epsilon'_{ijk} = \epsilon_{ijk}$. If $\epsilon_{ijk}$ is to define a tensor, we must according to \rf{11.12} have 
\begin{align*}
\epsilon_{ijk}' = a_{ip} \; a_{jq} \; a_{kr} \; \epsilon_{pqr}.
\end{align*}
Since $\epsilon'_{ijk} = \epsilon_{ijk}$, we have defined a tensor only if 
\begin{align}
\epsilon_{ijk} = a_{ip} \; a_{jq} \; a_{kr} \; \epsilon_{pqr}. \lbl{27.122}
\end{align}
Using the definition of the determinant from linear algebra, we can show that \rf{27.122}  is true, so \rf{25.12} in fact defines a tensor. This is the celebrated \ttx{Levi-Civita} tensor. \\
For example, choosing $i=1, \; j=2, \; k=3$ in \rf{27.122} we have 
\begin{align*}
a_{1p} \; a_{2q} \; a_{3r} \; \epsilon_{pqr} = \mathlarger{\sum}_{\sigma \in S_3} a_{1 \sigma (1)} \; a_{2 \sigma(2)} a_{3 \sigma (3)} \; \epsilon_{\sigma(1) \sigma(2) \sigma(3)},
\end{align*}
 where $S_3$ is the collection of all permutation of (1,2,3). From the definition of $\epsilon_{ijk}$ we have
\begin{align*}
\epsilon_{\sigma(1) \sigma(2) \sigma(3)} = (-1)^\sigma \; \epsilon_{123},
\end{align*}
where $(-1)^\sigma$ is $\pm 1$ according to if $\sigma$ is even or odd. But then we get 
\begin{align*}
a_{1p} \; a_{2q} \; a_{3r} \; \epsilon_{pqr} &= \Big\{ \mathlarger{\sum}_{\sigma \in S_3} (-1)^\sigma \; a_{1 \sigma(1)} \; a_{2 \sigma(2)} \; a_{3 \sigma(3)} \Big\} \; \epsilon_{123} \\ 
&= \det(a) \; \epsilon_{123} = \epsilon_{123},
\end{align*}
since the determinant of a unitary matrix is equal to 1. Thus \rf{27.122} holds for this case. The rest is proved in a similar way.

\subsection{Tensor properties}
Let $t_{ij}$ be the components of a tensor of rank 2, with respect to some observer $\mathcal{O}$. Let us assume that with respect to this observer, the components of the tensor obeys the identity
\begin{align}
t_{ij}=t_{ji}.\nonumber
\end{align}
\noindent Thus, in terms of linear algebra, the components form a symmetric matrix. Let the components of the tensor $t$, with respect to some other observer,$\mathcal{O'}$, be $t'_{ij}$. Is it true that we also have $t'_{ij}=t'_{ji}$? Thus, does the components of $t$ with respect to the observer $\mathcal{O'}$ also form a symmetric matrix?

Using the fact the $t$ is a tensor, we have
\begin{align}
t'_{ji}&=a_{jk}a_{il}t_{kl}\nonumber\\
&=a_{jk}a_{il}t_{lk}\nonumber\\
&=a_{il}a_{jk}t_{lk}\nonumber\\
&=t'_{ij}.\nonumber
\end{align}
\noindent Thus, if the components of a tensor form a symmetric matrix with respect to one observer, then the components with respect to any other observer also forms a symmetric matrix. It therefore make sense to say that the \textit{tensor}, $t$, is symmetric. Thus symmetry is a property of the components of a tensor that can be elevated to the tensorial level where it defines a property of the underlying tensor. Properties like these are called \textit{tensorial}. The great utility of tensors in theoretical physics is that they make it possible to express relations between observed quantities that are independent of the observers. Only such relations are actual \textit{physical} relations, and as a consequence only tensorial equations are physical. All fundamental physical equations are equations relating tensors, and they express tensorial relations between the relevant tensors. 

As an example of this, the ultimate theory describing gravitational interactions is the General Theory of Relativity discovered by Albert Einstein 100 years ago. The fundamental equation in this theory, the \textit{Einstein field equation} is a tensorial equation
\begin{align*}
R(g)_{\mu\nu}-\frac{1}{2}R(g)g_{\mu\nu}+\Lambda g_{\mu\nu}=\frac{8\pi G}{c^4}T_{\mu\nu},
\end{align*}
relating the geometric properties of space-time, described by the \textit{metric tensor}, $g=(g_{\mu\nu})$, and the density of mass and energy, described by the \textit{Stress-Energy tensor}, $T=(T{\mu\nu})$. Here $G$ is the gravitational constant, $c$ is the speed of light and $R(g)$ is the \textit{Ricci tensor} whose components are nonlinear functions of the metric tensor and all its partial derivatives up to second order. The Einstein field equations is thus a system of sixteen highly nonlinear coupled partial differential equations of second order. Describing the whole universe by one equation makes for an equation that is hard to solve, no surprise there!

\subsection{Tensor operations}

Let $\phi_i$ and $\psi_i$ be the components of two tensors $\phi$ and $\psi$. In any Cartesian frame define a number 
\begin{align}
(\phi, \psi) = \delta_{ij} \; \phi_i \; \psi_j .\lbl{22.12}
\end{align}
If \rf{22.12} is to define a tensor, then that tensor has only a single component. It is what we call a \ttx{scalar}. The transformation rule for such tensors is 
\begin{align*}
(\phi , \psi)' = (\phi, \psi),
\end{align*}
or in other words, the single component must be the same in any Cartesian frame. But according to \rf{22.12}, we have 
\begin{align*}
(\phi, \psi)' &= \delta_{ij}' \; \phi_i' \; \psi_j' = \delta_{ij} \; a_{ik} \; \phi_k \; a_{jl} \; \psi_l \\
&= \delta_{ij} \; a_{ik} \; a_{jl} \; \phi_k \; \psi_l = a_{ik} \; a_{il} \; \phi_k \; \psi_l\\ 
&= \delta_{kl} \; \phi_k \; \psi_l = (\phi, \psi),
\end{align*}
where we have used the identity \rf{6.1}. Recall that tensors like $\phi_i$ and $\psi_i$ which are of rank one are called \ttx{vectors}. We recognize that ($\phi, \psi$) is nothing but the scalar product of vectors in $\vb{R}^3$.  Operations defined on components of tensors that produce a new tensor are called \textit{Tensor operations}. We will now introduce several common and useful tensor operations.

\subsubsection{Contraction}
Let $\{ \alpha_i \}$ and $\{ \beta_{ij} \}$ be collections of components defining tensors $\alpha$ and $\beta$ where $\alpha$ is a vector and $\beta$ some tensor of rank 2. These kind of tensor are of great importance in fluid dynamics. Define a collection of components $\{ c_i \}$ by
\begin{align*}
c_i &= \alpha_j \; \beta_{ji}.
\end{align*}
Then we have
\begin{align*}
c_i' &= \alpha'_j \; \beta_{ji}' = a_{jk} \; \alpha_k \; a_{jl} \; a_{ir} \; \beta_{lr} \\ &= a_{jk} \; a_{jl} \; a_{ir} \; \alpha_k \; \beta_{lr} \\
&= \delta_{kl} \; a_{ir} \; \alpha_k \; \beta_{lr}\\ 
&= a_{ir} \; \alpha_k \; \beta_{kr}\\ 
&= a_{ir} \; c_r.
\end{align*}
Therefore, $\{ c_i \}$ defines a tensor $c$. The tensor $c$ is the \ttx{contraction} of $\alpha$ and $\beta$. In a similar way, more general contractions of tensors can be defined. These contractions always produce tensors. Examples are 
\begin{align*}
&\{ \alpha_{ij} \} , \{ \beta_{ij} \} \rightarrow \{ \alpha_{ij} \; \beta_{ij} \} &&\text{ rank 0},\\
&\{ \alpha_{ij} \} , \{ \beta_{ijkl} \} \rightarrow \{ \alpha_{ij} \; \beta_{ijkl} \} &&\text{ rank 2},\\
&\{ \epsilon_{ijk} \}, \{ \alpha_{i}\} , \{ \beta_{i} \} \rightarrow \{ \epsilon_{ijk} \; \alpha_{j} \; \beta_{k} \} &&\text{ rank 1}, \\
\end{align*}
\subsubsection{Dyadic notation}\label{cross}
Let $\{ \alpha_i \} , \; \{ \beta_i \}$ be components of tensors of rank 1, they are thus by definition vectors. Define a collection of components
 \begin{equation*}
 c_i = \epsilon_{ijk} \; \alpha_j \; \beta_k.
 \end{equation*}
These components define a tensor since we get them by a contraction involving three tensors, $\epsilon, \; \alpha$ and $\beta$. 

 This is an opportune point to introduce an older, commonly used notation for Cartesian tensors. 
In this notation, called \ttx{dyadic},  vectors are denoted by boldface letters $\boldsymbol{\alpha}, \; \boldsymbol{\beta}, \; \boldsymbol{\gamma}$ etc. The cross product of the vectors  $\boldsymbol{\alpha}, \; \boldsymbol{\beta}$  is written as $\boldsymbol{\alpha} \cp \boldsymbol{\beta}$. 

\noindent The reader should verify that the tensor introduced in the last example is  in fact the well known cross product of vectors. Thus we have
\begin{equation*}
  (\boldsymbol{\alpha} \cp \boldsymbol{\beta})_i = \epsilon_{ijk} \; \alpha_j \; \beta_k.
\end{equation*}
  The contraction 
$$ c_i = \alpha_j \; \beta_{ji}, $$ 
is written using the \ttx{dyadic} notation as
$$\vb{c} = \boldsymbol{\alpha} \vdot \boldsymbol{\beta}, $$ 
and the scalar contraction 
$$ c = \alpha_i \; \beta_{ij} \; \gamma_j, $$ 
is written as
$$c= \boldsymbol{\alpha} \vdot\boldsymbol{\beta}\vdot \boldsymbol{\gamma}. $$ 

\noindent Vector calculus formulas are usually displayed using the dyadic notation. For calculations and derivations I however find the component formulas more effective. For tensors of higher rank the dyadic notation become cumbersome and also ambiguous.

In order to use the formalism for Cartesian tensors in an effective way, we need some identities connecting $\epsilon_{ijk}$ and $\delta_{ij}$. There are many such, here are some of them 
\begin{align}
&\vb{1)} \; \; \epsilon_{ijk} \; \epsilon_{lmn} =  \text{det}\mqty(\delta_{il} && \delta_{im} && \delta_{in} \\ \delta_{jl} && \delta_{jm} && \delta_{jn} \\ \delta_{kl} && \delta_{km} && \delta_{kn}), \nonumber \\ 
&\vb{2)} \; \; \epsilon_{ijk} \; \epsilon_{lmk} = \delta_{il} \; \delta_{jm} - \delta_{im} \; \delta_{jl}, \nonumber \\ 
&\vb{3)} \; \; \epsilon_{ijk} \; \epsilon_{ljk} = 2 \; \delta_{il}, \nonumber \\
&\vb{4)} \epsilon_{ijk} \; \epsilon_{ijk} = 6. \nonumber  
\end{align}
Let $ \{ \alpha_i \} , \; \{ \beta_i \}$ and $ \{ \gamma_i \}$ be the components of three vectors. Then 
\begin{align} 
[\boldsymbol{\alpha} \cp (\boldsymbol{\beta} \cp \boldsymbol{\gamma})]_i &= \epsilon_{ijk} \; \alpha_j \; (\vb{\beta} \cp \vb{\gamma})_k \nonumber \\ 
&= \epsilon_{ijk} \; \alpha_j \; \epsilon_{kln} \; \beta_l \; \gamma_n = \epsilon_{ijk} \; \epsilon_{kln} \; \alpha_j \; \beta_l \; \gamma_n \nonumber \\ 
&= \epsilon_{ijk} \; \epsilon_{lnk} \; \alpha_j \; \beta_l \; \gamma_n \nonumber \\ 
&= \delta_{il} \; \delta_{jn} \; \alpha_j \; \beta_l \; \gamma_n - \delta_{in} \; \delta_{jl} \; \alpha_j \; \beta_l \; \gamma_n \nonumber \\
&= \beta_i \; \alpha_n \; \gamma_n - \gamma_i \; \beta_n \; \alpha_n \nonumber \\
&= [ \boldsymbol{\beta} \; ( \boldsymbol{\alpha} \cdot \boldsymbol{\gamma}) - \boldsymbol{\gamma} \; ( \boldsymbol{\alpha} \cdot \boldsymbol{\beta})]_i, \nonumber 
\end{align}
and we get the well known formula $$ \boldsymbol{\alpha} \cp (\boldsymbol{\beta} \cp \boldsymbol{\gamma}) = \boldsymbol{\beta} \; (\boldsymbol{\alpha} \cdot \boldsymbol{\gamma}) - \boldsymbol{\gamma} \; (\boldsymbol{\alpha} \cdot \boldsymbol{\beta}). $$

\subsubsection{Sum}\label{sum}
Let $\{ \alpha_{i_1 .. i_p}\}$, $\{ \beta_{i_1..i_p} \} $ be the components of two tensors $\alpha, \; \beta$, both of rank $p$. Define a collection of components $\{ \gamma_{i_1 .. i_p}\}$ by 
$$ \gamma_{i_1 .. i_p} = \alpha_{i_1..i_p} + \beta_{i_1 .. i_p}. $$
Show that $\{ \gamma_{i_1 .. i_p} \}$ are the components of a tensor $\gamma$ of rank $p$. The resulting tensor $\gamma$ is the \ttx{sum} of the two tensors $\alpha$ and $\beta$. 

\subsubsection{Product}\label{product}
Let $\{ \alpha_{i_1 .. i_p}\}$ , $\{ \beta_{i_1 .. i_q}\}$ be the components of two tensors $\alpha, \; \beta$ of rank $p$ and $q$. Define a collection of components $\{ \gamma_{i_1 .. i_{p+q}} \}$ by 
$$\gamma_{i_1 .. i_{p+q}} = \alpha_{i_1 .. i_p} \; \beta_{i_{p+1} .. i_{p+q}}. $$ 
Let us show that  that  $\{ \gamma_{i_1 .. i_{p+q}} \}$ are the components of a tensor $\gamma$ of rank $p+q$. 
\begin{align*}
\gamma'_{i_1 .. i_{p+q}} &= \alpha'_{i_1 .. i_p} \; \beta'_{i_{p+1} .. i_{p+q}}\\
&=a_{i_1,j_1}...a_{i_p,j_p}\alpha_{j_1 .. j_p} \; a_{i_{p+1},k_1}...a_{i_{p+q},k_q}\beta_{k_1 .. k_q}\\
&=a_{i_1,j_1}...a_{i_p,j_p} a_{i_{p+1},k_1}...a_{i_{p+q},k_q}\alpha_{j_1 .. j_p}\beta_{k_1 .. k_q}\\
&=a_{i_1,j_1}...a_{i_p,j_p} a_{i_{p+1},j_{p+1}}...a_{i_{p+q},j_{p+q}}\alpha_{j_1 .. j_p}\beta_{j_{p+1} .. j_{p+q}}\\
&=a_{i_1,j_1}...a_{i_p,j_p} a_{i_{p+1},j_{p+1}}...a_{i_{p+q},j_{p+q}}\gamma_{j_1 .. j_{p+q}}.
\end{align*}
Thus, $\gamma$ is a tensor that is called the  \ttx{tensor product} of the two tensors $\alpha$ and $\beta$ and is denoted by $\gamma=\alpha\beta$. In dyadic notation, a tensor $\beta$, of rank two, that is the tensor product of two vectors $\vb{u}$ and $\vb{v}$,  is written 
\begin{equation*}
\boldsymbol{\beta} = \vb{u} \; \vb{v}
\end{equation*}
The full contraction,$c$, of this tensor with two vectors $\boldsymbol{\alpha}$ and $\boldsymbol{\gamma}$ can be calculated using the dyadic notation in the following way
\begin{equation*}
c=\boldsymbol{\alpha} \cdot (\vb{u}\vb{v}) \cdot \boldsymbol{\gamma}  = (\boldsymbol{\alpha} \cdot \vb{u})(\vb{v} \cdot \boldsymbol{\gamma}),\\
\end{equation*}
where $\boldsymbol{\alpha} \cdot \vb{u}$ is the usual dot product of vectors.

\subsubsection{Gradient}\label{gradient}
Let $\phi$ be a function, and define a collection of component functions $\{ c_i \}$ by 
$$ c_i = \prt{x_i} \; \phi \nonumber. $$ 
Then $\{ c_i \}$ are the components of a vector because 
$$c_i' = \prt{x_i'} \; \phi' = \partial_{x_i'} x_r \; \prt{x_r} \; \phi, $$
and $$ x_r = a_{sr} \; x_s' \Rightarrow \partial_{x_i'} x_r = a_{ir},$$
thus
$$ c_i' = a_{ir} \; \prt{x_r} \; \phi = a_{ir} \; c_r .$$ 
The vector $c$ is the \ttx{gradient} of $\phi$. $c$ is a vector but in this context we call it a \ttx{vector field}. In the same way we call $\phi$ a scalar field. Note that the transformation rule for the vector field $c_i$ can more precisely be written as 
$$c'_i (x') = a_{ir} \; c_r (x), $$
where $x_i' = a_{ij} \; x_j$. For a scalar field we have the transformation rule $$ \phi'(x') = \phi(x),$$
the function value is thus the same at points that corresponds under the transformation of coordinates. 

In dyadic notation the gradient of a scalar field is written $\grad{\phi}$. Thus $$(\grad{\phi})_i = \prt{x_i} \; \phi. $$
In a similar way we can define the gradient of a tensor field $t$ with components $ \{ t_{i_1 .. i_m} \} $ by $$ (\grad{t})_{ii_1 .. i_m} = \prt{x_i} \; t_{i_1 .. i_m}, $$ where we have used the dyadic notation for the gradient of a tensor.
Taking gradients can clearly be repeated. For example, the tensor of rank two whose dyadic notation is $\grad\grad\phi$, is defined by
\begin{equation*}
(\grad\grad\phi)_{i,j}=\partial_{x_i}\partial_{x_j}\phi.
\end{equation*}

\subsubsection{Divergence}
Let $\{ t_{i_1 .. i_m} \}$ be the components of a tensor, $t$, of rank $m$. Define the components $$c_{i_1 .. i_{m-1}} = \prt{x_i} t_{ii_1 .. i_{m-1}}.$$
Let us show that the components $c_{i_1 .. i_{m-1}}$ defines a tensor.
\begin{align*}
c'_{i_1 .. i_{m-1}} &= \prt{x'_i}\; t'_{ii_1 .. i_{m-1}}\\
&=\partial_{x_i'} x_r \partial_{x_r}\;a_{i,j_1}a_{i_1,j_2}..a_{i_{m-1},j_m} t_{j_1,j_2, .. j_{m}}\\
&=a_{i,r} a_{i,j_1}a_{i_1,j_2}..a_{i_{m-1},j_m}\partial_{x_r} t_{j_1,j_2, .. j_{m}}\\
&=\delta_{r,j_1}a_{i_1,j_2}..a_{i_{m-1},j_m}\partial_{x_r} t_{j_1,j_2, .. j_{m}}\\
&=a_{i_1,j_2}..a_{i_{m-1},j_m}\partial_{x_j} t_{j,j_2, .. j_{m}}\\
&=a_{i_1,j_2}..a_{i_{m-1},j_m}c_{j_2, .. j_{m}}\\
&=a_{i_1,j_1}..a_{i_{m-1},j_{m-1}}c_{j_1, .. j_{m-1}}.
\end{align*}
This tensor is called the \ttx{divergence} of $t$ and is written $\div{t}$ in dyadic notation. Thus in this notation 
$$(\div{t})_{i_1 .. i_{m-1}} = \prt{x_i}  t_{ii_1 .. i_{m-1}}.$$
For tensors of rank $\geq 2$ we can define more than one divergence operation. For example, for rank two tensor fields the two divergence operations are
\begin{eqnarray*}
&  ( t \vdot \grad)_i = \prt{x_j} \; t_{ij},\\
 & (\grad\vdot t)_i = \prt{x_j} \; t_{ji}.
 \end{eqnarray*}
 Thus, as the formulas indicate, the dyadic notation for the two divergence operations are $t\vdot\grad$ and $\grad\vdot t$.
  
\subsubsection{Curl}\label{curl}

Let $\boldsymbol{\alpha}$ be a tensor field of rank one, thus by definition a vector field, then we have
$$ (\curl{\boldsymbol{\alpha}})_i = \epsilon_{ijk} \; \prt{x_j} \; \alpha_k. $$
In a similar way we can define the curl of higher rank tensor fields. For a tensor field, $t=(t_{i_1,i_2,..i_n})$ of rank $n$, we define a tensor field of rank $n$, which in dyadic notation is written as $ \curl{t} $,  by the components
$$ (\curl{t})_{i_1,i_2,..i_n} = \epsilon_{i_1kl} \; \prt{x_k} \; t_{l,i_2,..i_n}. $$
 
\noindent For tensors of rank greater than one, we can define more than one curl operation. For example,  a second curl of the tensor field $t$ can be defined by 
 $$(t \cp \grad)_{i_1,i_2,..i_n}= \epsilon_{i_nkl} \; \prt{x_k} \; t_{i_1,i_2,..l}. $$ 
  As indicated, the dyadic notation for this second curl is  $ t \cp \grad $.

\paragraph{Example 1}

Let $\{ \alpha_i \} , \{ \beta_i \} $ be the components of vectors denoted by $\boldsymbol{\alpha}, \; \boldsymbol{\beta}$ in dyadic notation. Then we have 
\begin{align}
[\div{(\boldsymbol{\alpha} \cp \boldsymbol{\beta})}] &= \prt{x_i} (\vb{\alpha} \cp \vb{\beta})_i \nonumber \\
&= \prt{x_i} \; \epsilon_{ijk} \; \alpha_j \; \beta_k \nonumber \\ 
&= \epsilon_{ijk} \; (\prt{x_i} \; \alpha_j) \; \beta_k + \epsilon_{ijk} \; \alpha_j \; (\prt{_{x_i} \; \beta_k}) \nonumber \\
&= \beta_k \; \epsilon_{kij} \; \prt{x_i} \; \alpha_j - \alpha_j \; \epsilon_{jik} \; \prt{x_i} \; \beta_k \nonumber \\ 
&= \beta_k \; [ \grad \cp \boldsymbol{\alpha}]_k - \alpha_j \; [ \curl{\boldsymbol{\beta}}]_j \nonumber \\
&= \boldsymbol{\beta} \cdot (\curl{\boldsymbol{\alpha}}) - \boldsymbol{\alpha} \vdot (\curl{\boldsymbol{\beta}}). \nonumber 
\end{align}

\paragraph{Example 2}
Let $\{ \alpha_i \} $ be the components of a vector. Then 
\begin{align} 
(\curl{(\curl{\boldsymbol{\alpha}})})_i &= \epsilon_{ijk} \; \prt{x_j} \; (\curl{\boldsymbol{\alpha}})_k \nonumber \\
&= \epsilon_{ijk} \; \prt{x_j} \; \epsilon_{klm} \; \prt{x_l} \; \alpha_m \nonumber \\ 
&= \epsilon_{ijk} \; \epsilon_{lmk} \; \prt{x_j} \; \prt{x_l} \; \alpha_m \nonumber \\ 
&= \delta_{il} \; \delta_{jm} \; \prt{x_j} \; \prt{x_l} \; \alpha_m - \delta_{im} \; \delta_{jl} \; \prt{x_j} \; \prt{x_l} \; \alpha_m \nonumber \\ 
&= \prt{x_i} \; \prt{x_j} \; \alpha_j - \prt{x_j} \; \prt{x_j} \; \alpha_i \nonumber \\ 
&= [ \grad \; (\div{\boldsymbol{\alpha}}) - \laplacian{\boldsymbol{\alpha}}]_i. \nonumber 
\end{align}
Thus we get the well known formula $$ \curl{(\curl{\boldsymbol{\alpha}})} = \grad{(\div{\boldsymbol{\alpha}})} - \laplacian{\boldsymbol{\alpha}}. $$ Where $\laplacian = \prt{x_i} \prt{x_i} $ is the 3D Laplace operator.

\subsection{The Divergence theorem for Cartesian tensors} 

Let $D$ be a domain in 3D, and let the boundary of $D$ be $S$. Let $\{n_i\}$ be the unit normal vector field defined on $S$ and pointing out of $D$. Then for any tensor of rank $m$ we have
\begin{equation} 
\int_D dV \prt{x_i} \; t_{i_1 .. i_m} = \int_S dS \; n_i \; t_{i_1 .. i_m}. \label{divthm1}
\end{equation} 
Since (\ref{divthm1}) is an identity between tensors we get a special case of (\ref{divthm1}) by contracting the first index
\begin{equation}
\int_D dV \; \prt{x_i} \; t_{ii_2 .. i_m} = \int_S dS \; n_i \; t_{ii_2 .. i_m}\label{divthm2}.
\end{equation}
Both (\ref{divthm1}) and (\ref{divthm2}) are called the divergence theorem even if (\ref{divthm1}) is more general than (\ref{divthm2}).
 
\paragraph{Example 1}
Let us consider (\ref{divthm1}) with a scalar, $\phi$ we get 
$$ \int_D dV \; \prt{x_i} \phi = \int_S dS \; n_i \; \phi, $$ 
or in dyadic notation 
$$ \int_D dV \; \grad{\phi} = \int_S dS \; \phi \; \vb{n}. $$ 

\paragraph{Example 2}
Let us a consider (\ref{divthm2}) for a vector field $\vb{a}$ with components $ \{ a_i \} $. 
$$ \int_D dV \; \prt{i} \; a_i = \int_S dS \; n_i \; a_i, $$ 
which in dyadic notation is 
$$ \int_D dV \; \div{\vb{a}} = \int_S dS \; \vb{n} \vdot \vb{a}. $$ 
This is the usual divergence theorem from vector calculus.
\paragraph{Example 3}
We use the tensor with components $\{ \epsilon_{ijk} \; a_j \}$ in (\ref{divthm2}) and get 
\begin{align*}
\int_D dV \; \epsilon_{ijk} \; \prt{x_i} \; a_j &= \int_S dS \; n_i \; \epsilon_{ijk} \; a_j , \\
&\Updownarrow  \\  \int_D dV \; \epsilon_{kij} \; \prt{x_i} \; a_j &= \int_S dS \; \epsilon_{kij} \; n_i \; a_j, \\ 
\end{align*}
which in dyadic notation is 
\begin{align} 
\int_D dV \; \curl{\vb{a}} = \int_S dS \; \vb{n} \cp \vb{a}. \nonumber 
\end{align}
\subsection*{Example 4}
Use the tensor with components $ \{ \epsilon_{ijk} \; a_j \; t_{kl} \} $ in (\ref{divthm2}). We get 
\begin{align} 
\int_S dS \; n_i \; \epsilon_{ijk} \; a_j \; t_{kl} &= \int_D dV \; \prt{x_i} \epsilon_{ijk} a_j  t_{kl} \nonumber\\
&= \int_D dV \; \{ \; \epsilon_{ijk} (\prt{x_i} a_j) \; t_{kl} + \epsilon_{ijk} \; a_j \; \prt{x_i} t_{kl} \; \} \nonumber \\
&= \int_D dV \; \{ \; \epsilon_{kij} \; (\prt{x_i} a_j) \; t_{kl} - a_j \; \epsilon_{jik} \; \prt{x_i} \; t_{kl} \; \}, \nonumber
\end{align} 
which in dyadic notation is 
$$ \int_S dS\; \vb{n} \vdot (\vb{a} \cp \vb{t}) = \int_D dV \; \{ \; ( \curl{\vb{a}}) \vdot \vb{t} - \vb{a} \vdot ( \curl{\vb{t}}) \; \}, $$ 
where the cross product of the vector $\vb{a}$ and the rank 2 tensor $t$ is defined by 
\begin{align}
(\vb{a} \cp \vb{t})_{ij} = \epsilon_{ikl} \; a_k \; t_{lj}. \nonumber
\end{align}

\subsection{Stoke's theorem for Cartesian tensors}

Let $S$ be a surface bounded by a closed curve $C$. Let $\{ t_{i_1 .. i_m} \}$ be the components of a Cartesian tensor of rank $m$. Then Stoke's theorem is 
\begin{equation}
\int_C dx_k \; t_{i_1 .. i_m} = \int_S dS \; n_i \; \epsilon_{ijk} \; \prt{x_j} t_{i_1 .. i_m}.\label{stokesthm1}
\end{equation}
Since (\ref{stokesthm1}) is an identity for tensors, we can contract the first index and get 
\begin{equation}
\int_C dx_i \; t_{ii_2 .. i_m} = \int_S dS \; n_i \; \epsilon_{ijk} \; \prt{x_j} \; t_{ki_2 .. i_m}.\label{stokesthm2}
\end{equation} 
The orientation of $C$ is related to the choice of unit normal for $S$ in the same way as for the usual Stoke's theorem in vector calculus.
\paragraph{Example 1}
Use a scalar, $\phi$, in (\ref{stokesthm1}). We get 
$$ \int_C dx_k \; \phi = \int_S dS \; n_i \; \epsilon_{ijk} \; \prt{x_j} \phi, $$ 
which in dyadic notation is 
$$ \int_C d\vb{l} \; \phi = \int_S dS \; \vb{n} \cp \grad{\phi}. $$
\paragraph{Example 2}
Using the components of a vector $\{ a_i \}$ in (\ref{stokesthm2}) we get 
$$ \int_C dx_i \; a_i = \int_S dS \; n_i \; \epsilon_{ijk} \; \prt{x_j} a_k, $$
which in dyadic notation is
$$ \int_C d\vb{l} \vdot \vb{a} = \int_S dS \; \vb{n} \vdot (\curl{\vb{a}}). $$
This is the usual Stoke's theorem from vector calculus.
\paragraph{Example 3}
Using a tensor with components $\{ \epsilon_{ijk} \; a_j \} $ in (\ref{stokesthm2}). We get 
\begin{align*}
\int_C dx_i \; \epsilon_{ijk} \; a_j &= \int_S dS \; n_i \; \epsilon_{ijn} \; \prt{x_j} \epsilon_{nlk} a_l \\ &= \int_S dS \; \epsilon_{ijn} \; \epsilon_{lkn} \; n_i \; \prt{x_j} a_l \\ 
&= \int_S dS \; \{ \; \delta_{il} \; \delta_{jk} \; n_i \; \prt{x_j} a_l - \delta_{ik} \; \delta_{jl} \; n_i \; \prt{x_j} a_l \; \} \\
&= \int_S dS \; \{ \; (\prt{x_k} a_l) \; n_l - n_k \; \prt{x_j} a_j \; \},
\end{align*}
which in dyadic notation is 
$$ \int_C d\vb{l} \cp \vb{a} = \int_S dS \; \{ \; (\grad{\vb{a}}) \vdot \vb{n} - \vb{n} \; \div{\vb{a}} \}. $$ 

\noindent This ends our exposition of Cartesian tensors.

\subsection{Exercises}
\paragraph{Problem 1}

Let $t_{ij}$ be the components of a tensor of rank 2, with respect to some observer $\mathcal{O}$. Let us assume that with respect to this observer, the components of the tensor obeys the identity
\begin{align}
t_{ij}=-t_{ji}.\nonumber
\end{align}
\noindent Thus, in terms of linear algebra, the components form an anti-symmetric matrix. Show that anti-symmetry is a tensorial property.

\paragraph{Problem 2}

Let $t_{i_1 .. i_p}(x)$ be the components of a tensor field of rank p with respect to some observer $\mathcal{O}$. Let us assume that with respect to this observer, we have 
\begin{align}
t_{i_1 .. i_p}(x_0)=0,\nonumber
\end{align}
at some point $x_0$.  Show that the same equation holds for the components of the tensor field $t'_{i_1 .. i_p}(x')$ with respect to any observer $\mathcal{O'}$ at the point corresponding to $x_0$ under the change of coordinates $x'=a_{ij}x_j$. Thus a tensor field being zero at some point is a tensorial property of tensor fields.

\paragraph{Problem 3}
Prove that the sum of tensors and product of tensors, as defined in section \ref{sum} and \ref{product}, are tensors.
\paragraph{Problem 4}
Show that the gradient and curl of tensors, as defined in section \ref{gradient} and \ref{curl}, are tensors.
\paragraph{Problem 5}
Prove the following vector identifies using Cartesian tensors. 
\begin{align*}
&\vb{a)} \; \div(\curl\boldsymbol{\alpha}) = 0, \\ 
&\vb{b)} \;\curl(\grad\alpha) = 0,  \\ 
&\vb{c)} \; \grad(\boldsymbol{\alpha} \vdot \boldsymbol{\beta}) = \grad\boldsymbol{\alpha} \vdot \boldsymbol{\beta} + \grad\boldsymbol{\beta} \vdot \boldsymbol{\alpha} \\ 
& \qquad \qquad \; \; \; \; = (\boldsymbol{\alpha} \vdot \grad) \; \boldsymbol{\beta} + (\boldsymbol{\beta} \vdot \grad) \; \boldsymbol{\alpha} + \boldsymbol{\alpha} \cp (\curl\boldsymbol{\beta}) + \boldsymbol{\beta} \cp (\curl\boldsymbol{\alpha}),  \\
&\vb{d)} \; \curl(\boldsymbol{\alpha}\cp\boldsymbol{\beta}) = \boldsymbol{\alpha} \; \div{\boldsymbol{\beta}} - \boldsymbol{\beta} \; \div{\boldsymbol{\alpha}} + (\boldsymbol{\beta} \vdot \grad) \; \boldsymbol{\alpha} - (\boldsymbol{\alpha} \vdot \grad) \; \boldsymbol{\beta}.
\end{align*}

\setcounter{equation}{0}
\section{Fluid dynamics}

Liquids and gases are similar, in the sense that they have no fixed shape like solids do. A liquid or a gas will shape themselves to fit perfectly to any container we pour them into. This similarity of liquids and gases makes it possible to present their mathematical description in a unified way. This unified way is called \textit{fluid dynamics}.

  Liquids and gases are certainly different too. The first one is for example very hard to compress, whereas the second one is easy to compress.
  
\noindent The part of fluid dynamics that concerns itself with easily compressible substances is called \textit{gas dynamics}. Water is a liquid and its theoretical description is called \textit{hydrodynamics}. Hydrodynamics is a subfield of fluid dynamics.

A fluid appears to us to be a continuous substance, it fills space smoothly. However, we know that this continuity is only apparent. Underneath it all, we know that fluids consist of discrete entities in the form of atoms or molecules. These atoms or molecules move around, sometimes at great speed, and interact with each other. This immense activity, at the microscopic level, appears to us, at the macroscopic level, as a continuously moving fluid. Exactly how this happen is not understood in all detail, even today, but the overall picture is clear. 

We will \textit{not} try to do a detailed derivation of the equations of fluid dynamics from the motion of atoms and molecules. What we will do, is to use some of the descriptive tools from this derivation in order to gain insight into the various terms occurring in the equations of fluid dynamics. The exposition of fluid dynamics in this section is inspired by the books written on the subject by Landau and Lifschitz\cite{Landau} and G. B. Whitham\cite{Whitham}.

\noindent We will start by assuming that the particles underlying the fluid are simple mass points. We will also assume that they are all of the same type,  and have a common mass $m$. 
What we really are assuming here is that the particles in the fluid only collide elastically, so that any internal degrees of freedom, like rotations and vibrations, are never excited in collisions, and are invisible from our macroscopic  point of view. Not all fluids are like this, so this is a real restriction. It does however hold approximately for many liquids and gases. 

We know that the particles in the fluid move along curves that can be described using Newton's law. However, the number of particles in the fluid is so immense that the equations for the particles  can not be solved, even on today's largest machines. And even if we \textit{could} solve them, making measurements precise enough to supply the equations with initial conditions is beyond our means for the foreseeable future.
In situations like this, one resorts to a coarser description using a particle distribution function $f(\mathbf{x},\mathbf{u},t)$. Here 
\begin{align}
f(\mathbf{x},\mathbf{u},t)\; d\mathbf{x} \; d\mathbf{u}, \label{eq: 1}
\end{align}
is the number of particles in a domain in phase space of volume $d\mathbf{x}\;d\mathbf{u}$ centred on position $\mathbf{x}$ and velocity $\mathbf{u}$.

The \textit{mass density} at position $\mathbf{x}$ and time $t$ is 
\begin{align}
\rho(\mathbf{x}, t) = \int d\mathbf{u}\;m\;f(\mathbf{x}, \mathbf{u}, t), \label{eq: 2}
\end{align}
and the \textit{mean velocity},  $\mathbf{v}$, of the particle distribution at $\mathbf{x}$, $t,$ is defined by the expression
\begin{align}
\rho(\mathbf{x}, t)\;\mathbf{v}(x,t) = \int d\mathbf{u}\;m\;\mathbf{u}\;f(\mathbf{x}, \mathbf{u}, t). \label{eq: 3}
\end{align}

\noindent A key assumption in the theory of fluid dynamics is that there is a "mesoscale", much larger than the microscale and much smaller than the macroscale. On the mesoscale $\rho, \mathbf{v},$ and other quantities we will define as moments of the particles distribution, are assumed to be constant.

\begin{figure}[htbp]
\centering
\includegraphics{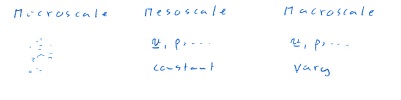}
\label{fig1}
\caption{}
\end{figure}

\noindent Traditionally, a mesoscale piece of the fluid is called a \textit{fluid element}. The assumption is that the fluid element is small on a macroscopic scale but large on a microscopic scale, in the sense that it contains a large number of particles. The dominant, short time, behavior of a fluid element, is a translation along the velocity field. On a larger timescale the fluid element also deforms. Imagine marking out a fluid element by a tiny drop of ink in the fluid. This drop of ink is a fluid element and initially it just translates along the velocity field of the fluid. On larger timescales it deforms and eventually smears out and vanishes.

If this split into micro-, meso- and macroscale is not possible, we can not hide the particle nature of the system, and fluid dynamics does not apply. If this is the case one must use \textit{kinetic theory}, which is much more challenging than fluid dynamics. In these notes we will not discuss kinetic theory.

Let $\mathbf{v}(\mathbf{x},t)$ be the fluid velocity field as defined in (\ref{eq: 3}) and let $\mathbf{x}(t)$ be the position of a fluid element. Then we have 
\begin{align}
\frac{d\mathbf{x}}{dt}=\mathbf{v}(\mathbf{x},t). \label{eq: 4}
\end{align}
This is the equation of motion for a fluid element. 

\noindent Let now $A(\mathbf{x},t)$ be some local quantity associated with the fluid, like $\rho(\mathbf{x},t), \mathbf{v}(\mathbf{x}, t)$... Then the function $A(t)$, defined by 
\begin{align}
A(t)=A(\mathbf{x}(t),t), \label{eq: 5}
\end{align}
where $\mathbf{x}(t)$ is the position of a fluid element, will describe how $A$ changes for a fluid element following the fluid flow. Using the chain rule we have
\begin{align}
\frac{dA}{dt}(t) =\partial_t \; A(\mathbf{x}(t),t)+\frac{d\mathbf{x}}{dt}(t)\vdot \grad{A(\mathbf{x},t)} = (\partial_t + \mathbf{v} \vdot \grad)\; A(\mathbf{x} , t)\; \mathlarger{|}_{\mathbf{x} = \mathbf{x}(t)}. \label{eq: 6}
\end{align}
The operator
\begin{align}
\frac{D}{Dt} = \partial_t + \mathbf{v} \vdot \grad, \label{eq: 7}
\end{align}
plays an important role in fluid dynamics, and is called the {\it material derivative}. We have for example 

\begin{align}
\frac{D\rho}{Dt} = \partial_t \; \rho + \mathbf{v} \vdot \grad \rho, \label{eq: 8}\\
\frac{D\mathbf{v}}{Dt} = \partial_t \; \mathbf{v} + \mathbf{v} \vdot \grad \mathbf{v}. \nonumber
\end{align}

\noindent For the particle dynamics underlying the fluid motion we have three fundamental conservation laws

\begin{align*}
&1. \; \text{Conservation} \; \text{of}\;  \text{mass.}\\
&2. \; \text{Conservation} \; \text{of}\; \text{momentum.}\\
&3. \; \text{Conservation} \; \text{of}\;\text{ energy.}\\
\end{align*}
We are going to find macroscopic analogues for the microscopic quantities, mass, momentum and energy and postulate that they are conserved. This will give us the equations of fluid dynamics. 

\subsection{Conservation of mass}

Let $V$ be some volume of fluid, with bounding surface $S$.

\begin{figure}[htbp]
\centering
\includegraphics{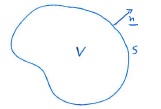}
\label{fig1}
\caption{: A fluid volume $V$, with boundings surface $S$ and unit normal $\vb{n}$.}
\end{figure}
\noindent The total mass in $V$ is a function of $t$, defined by
\begin{align}
M(t) = \int_V dV \; \rho (\mathbf{x}, t). \label{eq: 9}
\end{align}
Since mass is postulated to be conserved, $M(t)$ can vary only if mass enters or leaves the volume $V$ by crossing the boundary. Let us analyze this crossing carefully, we will use it several times.

Let us consider a particle that is moving at velocity $\mathbf{u}$ and is close to the boundary surface at time $t$. During the time from $t$ to $t + dt$ the particle is moving a distance normal to the surface given by
\begin{align*}
dl_{\mathbf{u}} = \mathbf{n} \vdot \mathbf{u} \: dt.
\end{align*}
Here $dl_{\mathbf{u}}$ is a \textit{signed distance}, $dl_{\mathbf{u}} >0$ when the particle has a velocity component parallel to $\mathbf{n}$ pointing in the direction of $\mathbf{n}$. If this velocity component points opposite to $\mathbf{n}$ we have $dl_{\mathbf{u}} < 0 $.  Let the signed volume $d V_{\mathbf{u}}$ be defined by
\begin{align*}
d V_{\mathbf{u}} = dA\;\mathbf{n} \vdot \mathbf{u} \; dt, 
\end{align*}
where $dA$ is a surface area element. If $d V_{\mathbf{u}} > 0$, all particles of velocity $\mathbf{u}$, that are inside the volume $d V_{\mathbf{u}}$ at time $t$, will cross the boundary of $V$ from the inside, and leave $V$ during the time interval between $t$ and $t + dt$.
\begin{figure}[htbp]
\centering
\includegraphics{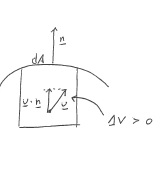}
\label{fig1}
\caption{}
\end{figure}

\noindent If $d V_{\mathbf{u}} < 0$,  all particles of velocity $\mathbf{u}$, that are inside the volume $d V_{\mathbf{u}}$ at time $t$, will enter $V$ during the time interval between $t$ and $t + dt$. Therefore, the number of particles with velocity in the range $d\mathbf{u}$ around $\mathbf{u}$  leaving or entering the volume $V$ through the surface element $dA$, is 
\begin{align*}
d N_{\mathbf{u}}&=d V_{\mathbf{u}}\;d\mathbf{u}\; f(\mathbf{x}, \mathbf{u}, t)\\
&=dA\;\mathbf{n} \vdot \mathbf{u} \; dt \;d\mathbf{u}\; f(\mathbf{x}, \mathbf{u}, t),
\end{align*}
and the number of  particles  with velocity in the range $d\mathbf{u}$ around $\mathbf{u}$  leaving or entering  the volume $V$ through the surface element $dA$ per unit time, $dP_{\mathbf{u}}$,  is given by 
\begin{align}
d P_{\mathbf{u}}=dA\;\mathbf{n} \vdot \mathbf{u}\;d\mathbf{u}\; f(\mathbf{x}, \mathbf{u}, t)\label{Pu}.
\end{align}
This expression will be used repeatedly in the following.

Using (\ref{Pu}), we can conclude that the net mass  of the particles with velocity in the range $d\mathbf{u}$ around $\mathbf{u}$, entering or leaving the volume $V$ through the surface element $dA$, per unit time, is
\begin{align}
m\;dP_{\mathbf{u}}=dA\;d\mathbf{u} \; m \; \mathbf{n} \vdot \; \mathbf{u} \; f(\mathbf{x}, \mathbf{u}, t), \label{eq: 12}
\end{align}
Integrating this expression over all possible velocities, we get the mass transport, per unit time, through the surface element $dA$. This is what we call the \textit{mass flux}.
But 
\begin{align}
dA\;\int d\mathbf{u} \; m \; \mathbf{n} \vdot \mathbf{u} \; f &=dA\; \mathbf{n} \vdot \int d\mathbf{u} \; m \; \mathbf{u} \; f 
&=dA\; \mathbf{n} \vdot (\rho \; \mathbf{v}) =dA\; \rho \; \mathbf{v}\vdot  \mathbf{n} , \label{eq: 13}
\end{align}
so that the mass flux through $dA$ is
\begin{align}
dA\;\rho \;  \mathbf{v} \vdot \mathbf{n}. \label{eq: 14}
\end{align}
The total mass passing through the bounding surface $S$ is therefore
\begin{align}
\int_S dA \; \rho \; \mathbf{v} \vdot \mathbf{n}. \label{eq: 15}
\end{align}
Conservation of mass then imposes the identity
\begin{align}
\frac{dM}{dt} &= - \int_S dA \; \rho \; \mathbf{v} \vdot \mathbf{n}, \; \nonumber \\ 
&\Updownarrow \nonumber \\\; \int_V dV \;  \partial_t \; \rho &= - \int_S dA \; \rho \; \mathbf{v} \vdot \mathbf{n},\; \nonumber\\
&\Updownarrow \nonumber \\
 \int_V dV \; \{ \partial_t \; \rho + \grad \vdot (\rho \; \mathbf{v})\} &= 0. \label{eq: 16}
\end{align}
This identity holds for all volumes $V$. We can therefore conclude that
\begin{align}
\partial_t \; \rho + \grad \vdot (\rho \; \mathbf{v}) = 0. \label{eq: 17}
\end{align}
This is the conservation of mass in differential form, and is the first fundamental equation of fluid dynamics. Using the material derivative, (\ref{eq: 17}) can be written in the form
\begin{align}
 \frac{D\rho}{Dt} &= - \rho \; \grad \vdot \mathbf{v}. \label{eq:18}
\end{align}

\subsection{Conservation of momentum}

Let, as before, $V$ be a volume of fluid with bounding surface $S$. The total momentum inside the volume at a time $t$ is 
\begin{align*}
\mathbf{P}(t) = \int_V dV \; \rho \; \mathbf{v}.
\end{align*}
Using (\ref{Pu}), we can conclude that the net momentum  of the particles with velocity in the range $d\mathbf{u}$ around $\mathbf{u}$, entering or leaving the volume $V$ through the surface element $dA$, per unit time, is
\begin{align*}
m\mathbf{u}\;dP_{\mathbf{u}}=dA\;d\mathbf{u} \; m\mathbf{u}\; \mathbf{n} \vdot \; \mathbf{u} \; f(\mathbf{x}, \mathbf{u}, t),
\end{align*}
Integrating this expression over all possible velocities, we get the momentum transport, per unit time, through the surface element $dA$
\begin{align}
dA \; \int d\mathbf{u} \; m \;  \mathbf{u} \; \mathbf{n} \vdot \mathbf{u} \; f(\mathbf{x}, \mathbf{u}, t). \label{eq:20}
\end{align}
This is what we call the \textit{momentum flux}.

Define a velocity $\mathbf{c}$ by
\begin{align}
\mathbf{u} = \mathbf{v} + \mathbf{c},\label{eq:21}
\end{align}
thus, $\mathbf{c}$ describes the deviation of the particle velocity from the local mean velocity, $\mathbf{v}(\mathbf{x}, t)$.

Using (\ref{eq:21}) in (\ref{eq:20}) we get
\begin{align*}
&\int d\mathbf{u} \; m \; \mathbf{u} \; \mathbf{n} \vdot \mathbf{u} \; f = \int d\mathbf{c} \; m \; (\mathbf{v} +\mathbf{c} )\; \mathbf{n} \vdot (\mathbf{v} + \mathbf{c}) \; f
\\ &= \int d\mathbf{c} \; m  \; \mathbf{v} \; \mathbf{n} \vdot \mathbf{v} \; f + \int d \mathbf{c} \; m \; \mathbf{v} \; \mathbf{n} \vdot \mathbf{c} \; f
\\ &+ \int d \mathbf{c} \; m \; \mathbf{c} \; \mathbf{n} \vdot \mathbf{v} \; f + \int d\mathbf{c} \; m \; \mathbf{c} \; \mathbf{n} \vdot \mathbf{c} \;f. 
\end{align*}
From the definition of $\mathbf{c}$ as the deviation from the mean velocity we must have
\begin{align*}
\int d\mathbf{c} \; \mathbf{c} \; f = 0. 
\end{align*}
Therefore
\begin{align}
\int d\mathbf{c} \; m \; \mathbf{v} \; \mathbf{n} \vdot \mathbf{c} \; f = m \; \mathbf{v} \; \mathbf{n} \vdot (\int d\mathbf{c} \; \mathbf{c} \; f) &= 0, \nonumber\\
\int d\mathbf{c} \; m \; \mathbf{c} \; \mathbf{n} \vdot \mathbf{v} \; f = m \; \mathbf{n} \vdot \mathbf{v} \; (\int d\mathbf{c} \; \mathbf{c} \; f) &= 0,\nonumber\\
\int d\mathbf{c} \; m \; \mathbf{v} \; \mathbf{n} \vdot \mathbf{v} \; f &= (\int d \mathbf{c} \; m \; f) \; \mathbf{v} \; \mathbf{n} \vdot \mathbf{v} \nonumber
\\ &= \rho \; \mathbf{v} \; \mathbf{n} \vdot \mathbf{v} = (\rho \; \mathbf{v} \mathbf{v}) \vdot \mathbf{n}, \label{eq:22}
\\ \int d \mathbf{c} \; m \; \mathbf{c} \; \mathbf{n} \vdot \mathbf{c}\; f &= ( \int d\mathbf{c}\;m\;\mathbf{c}  \mathbf{c} \; f) \vdot \mathbf{n} \nonumber \\ &= -\mathbf{{\cal P}} \vdot \mathbf{n}, \label{eq:23}
\end{align}
where $\mathbf{{\cal P}}$ is a tensor of rank 2, called the \textit{stress tensor} for the fluid. Note that we are using the dyadic notation for Cartesian tensors in (\ref{eq:22}) and (\ref{eq:23}). Explicitly, the stress tensor is 
\begin{align*}
\mathbf{{\cal P}} = - \int d\mathbf{c} \; m \; \mathbf{c} \vb{c} \; f.
\end{align*}
The choice of minus sign is conventional. We observe that ${\cal P}_{ij}$ describes the correlation between fluctuating velocities in the direction of the $i$ and $j$ axis. Also note that $\mathbf{{\cal P}}$ is \textit{symmetric}
\begin{align*}
{\cal P}_{ij} ={\cal P}_{ji}.
\end{align*}
Using the above expressions we have found that the momentum flux through the surface element $dA$ is given by 
\begin{align*}
dA \; \int d\mathbf{u} \; m \;  \mathbf{u} \; \mathbf{n} \vdot \mathbf{u} \; f(\mathbf{x}, \mathbf{u}, t)
= dA\;\left(\rho \; \mathbf{v} \mathbf{v} \vdot \mathbf{n} -\mathbf{ {\cal P}} \vdot \mathbf{n}\right) =dA \;( \rho \; \mathbf{v} \mathbf{v} -\mathbf{ {\cal P}}) \vdot \mathbf{n}
 \end{align*}

 Let us allow for the possibility that there is a volume force acting on the fluid. Gravity is such a force, as is the electromagnetic force. The last one would act if the fluid consisted of charged particles. For geophysical applications, volume forces will enter as gravitational forces and inertial forces, like the Coriolis force.  Denote the volume force density by $\mathbf{F}_V$. Recalling that force is change in momentum per unit time, the law of conservation of momentum implies that
\begin{align*}
\frac{d\mathbf{P}}{dt} = - \int_S dA \; ( \rho \; \mathbf{v} \mathbf{v} -\mathbf{ {\cal P}}) \vdot \mathbf{n} + \int_V dV \; \mathbf{F}_V.
\end{align*}
Using the divergence theorem for 2-tensors, and letting the volume approach zero, we get the second fundamental equation of fluid dynamics
\begin{align}
\partial_t \; (\rho \; \mathbf{v}) + \grad \vdot(\rho\; \mathbf{v} \mathbf{v} - \mathbf{ {\cal P}}) = \mathbf{F}_V.\label{eq:23}
\end{align}
This equation can be simplified using the equation of mass conservation
\begin{align*}
\partial_t \; (\rho \; \vb{v}) &= \partial_t \; \rho \; \vb{v} + \rho \; \partial_t \; \vb{v} \\ &= - \vb{v} \vdot \grad \rho \; \vb{v} - \rho \; \grad \vdot \vb{v} \vb{v} + \rho \; \partial_t \; \vb{v}. 
\end{align*}
Using the component notation for Cartesian tensors, we have
\begin{align*} 
(\div{\rho \; \vb{v} \vb{v}})_j &= \partial_{x_i} \; (\rho \; v_i \; v_j) \\
&= \partial_{x_i} \; \rho \; v_i \; v_j + \rho \; \partial_{x_i} \; v_i \; v_j + \rho \; v_i \; \partial_{x_i} \; v_j.
\end{align*}
Thus in dyadic notation we have the identity
\begin{align*}
\div{\rho \; \vb{v} \vb{v}} = \vb{v} \vdot \grad \rho \; \vb{v} + \rho \; \div{\vb{v}} \; \vb{v} + \rho \; \vb{v} \vdot \grad \vb{v}.
\end{align*}
Therefore, (\ref{eq:23}) turns into
\begin{align*}
\rho \; \partial_t \; \vb{v} - \vb{v} \vdot \grad  \rho \; \vb{v} - \rho \; \grad \vdot \vb{v} \vb{v} + \vb{v} \vdot \grad  \rho \; \vb{v} \nonumber \\ + \rho \; \grad \vdot \vb{v} \vb{v} + \rho \; \vb{v} \vdot \grad \vb{v} - \div{\mathbf{ {\cal P}}} = \vb{F}_V.
\end{align*}
Thus
\begin{align}
\rho \; \partial_t \; \vb{v} + \rho \; \vb{v} \vdot \grad \vb{v} = \div{\mathbf{ {\cal P}}} + \vb{F}_V. \label{eq:36}
\end{align}
This is the second fundamental equation of fluid dynamics in simplified form. Using the material derivative, equation (\ref{eq:36}) can be written compactly as 
\begin{align}
\rho\frac{D \vb{v}}{Dt} = \div{\mathbf{ {\cal P}}} + \vb{F}_V. \label{eq:36.1}
\end{align}
This is clearly Newton's law for a fluid element where the force is given by 
\begin{align}
\vb{F} = \div{\mathbf{ {\cal P}}} + \vb{F}_V. \label{eq:36.2}
\end{align}

\subsection{Conservation of energy}

For a particle of mass $m$ and velocity $\vb{u}$, the kinetic energy is given as 
\begin{align*}
E = \frac{1}{2} \; m \; \vb{u}^2. 
\end{align*}
We will assume that the particle energy is dominated by the kinetic part. The macroscopic energy is then
\begin{align*}
\int d\vb{u} \; \frac{1}{2} \; m \; \vb{u}^2 \; f &= \int d\vb{c} \; \frac{1}{2} \; m \; (\vb{v} + \vb{c})^2 \; f \\
&= \int d\vb{c} \; \frac{1}{2} \; m \; \vb{v}^2 \;f + \int d\vb{c} \; m \; \vb{v} \vdot \vb{c} \; f + \int d\vb{c} \; \frac{1}{2} \; m \; \vb{c}^2 \; f  \\
&= \frac{1}{2} \; \rho \; \vb{v}^2 + \rho \; e, 
\end{align*}
where we have defined $e$ by 
\begin{align*}
\rho \; e = \int d\vb{c} \; \frac{1}{2} \; m \; \vb{c}^2 \; f.
\end{align*}
Clearly, $e$ measures the kinetic energy in the fluctuating part of the particle motion. We call $e$ the \textit{internal energy} of the fluid.

Transport of kinetic  energy across the bounding surface element $dA$, is given by
\begin{align*}
 \int \frac{1}{2} \; m \; \vb{u}^2 \;dP_{\vb{u}} =dA \; \int d\vb{u} \; \frac{1}{2} \; m \; \vb{u}^2 \; \vb{n} \vdot \vb{u} \; f, 
\end{align*}
where we have argued as in (\ref{eq: 12}). This is what we call the \textit{energy flux}

Introducing mean and fluctuating velocity as before, we have
\begin{align}
\int d\vb{u} \; \frac{1}{2} \; m \; \vb{u}^2 \; \vb{n} \vdot \vb{u} \; f
&= \int d\vb{c} \; \frac{1}{2} \; m \; (\vb{v} + \vb{c} )^2 \; \vb{n}  \vdot (\vb{v} + \vb{c}) \; f \nonumber\\ &= \int d\vb{c} \; \frac{1}{2} \; m \; \vb{v}^2 \; \vb{n} \vdot \vb{v} \; f + \int d\vb{c} \; m \; \vb{v} \vdot \vb{c} \; \vb{n} \vdot \vb{v} \; f \nonumber \\ &+ \int d\vb{c} \; \frac{1}{2} \; m \; \vb{c}^2 \; \vb{n} \vdot \vb{v} \; f + \int d\vb{c} \; \frac{1}{2} \; m \; \vb{v}^2 \; \vb{n} \vdot \vb{c} \; f \nonumber \\ &+ \int d\vb{c} \; m \; \vb{v} \vdot \vb{c} \; \vb{n} \vdot \vb{c} \; f + \int d\vb{c} \; \frac{1}{2} \; m \; \vb{c}^2 \; \vb{n} \vdot \vb{c} \; f,  \label{eq:37}
\end{align}
and
\begin{align}
\int d \vb{c} \; \frac{1}{2} \; m \; \vb{v}^2 \; \vb{n} \vdot \vb{v} \; f = \frac{1}{2} \; \rho \; \vb{v}^2 \;\vb{n} \vdot \vb{v},  \nonumber\\
\int d\vb{c} \; m \; \vb{v} \vdot \vb{c} \; \vb{n} \vdot \vb{v} \; f = m \; \vb{n} \vdot \vb{v} \; \vb{v} \vdot (\int d\vb{c} \; \vb{c} \; f) = 0, \nonumber \\
\int d\vb{c} \; \frac{1}{2} \; m \; \vb{c}^2 \; \vb{n} \vdot \vb{v} \; f = (\int d\vb{c} \; \frac{1}{2} \; m \; \vb{c}^2 \; f) \; \vb{n} \vdot \vb{v}  \nonumber\\ = \rho \; e \; \vb{n} \vdot \vb{v}, \nonumber \\
\int d\vb{c} \; \frac{1}{2} \; m \; \vb{v}^2 \; \vb{n} \vdot \vb{c} \; f = \frac{1}{2} \; m \; \vb{v}^2 \; \vb{n} \vdot ( \int d\vb{c} \; \vb{c} \; f ) = 0, \nonumber\\
\int d \vb{c} \; m \; \vb{v} \vdot \vb{c} \; \vb{n} \vdot \vb{c} \; f = \vb{v} \vdot ( \int d\vb{c} \; m \; \vb{c} \vb{c} \; f) \vdot \vb{n}  \nonumber\\ = - \vb{v} \vdot  \mathbf{ {\cal P}} \vdot \vb{n}, \nonumber \\
\int d\vb{c} \; \frac{1}{2} \; m \; \vb{c}^2 \; \vb{n} \vdot \vb{c} \; f = (\int d\vb{c} \; \frac{1}{2} \; m \; \vb{c}^2 \vb{c} \; f) \vdot \vb{n}  \nonumber\\ = \vb{q} \vdot \vb{n}, \label{eq:38}
\end{align}
where we have introduced the vector $\vb{q}$ by 
\begin{align*}
\vb{q} = \int d \vb{c} \; \frac{1}{2} \; m \; \vb{c}^2 \vb{c} \; f. 
\end{align*}
This vector describes the transport of internal energy by the fluctuating velocity field $\vb{c}$. The vector $\vb{q}$ is called the \textit{heat flux} vector.

\noindent Using (\ref{eq:38})  in (\ref{eq:37}) we get the following formula for the {\it energy flux} through the surface element $dA$
\begin{align*}
 dA \; \{ \frac{1}{2} \; \rho \; \vb{v}^2 \; \vb{n} \vdot \vb{v} + \rho \; e \; \vb{n} \vdot \vb{v} - \vb{v} \vdot \mathbf{ {\cal P}} \vdot \vb{n} + \vb{q} \vdot \vb{n} \}.
\end{align*}
The work per unit time and unit volume performed by the volume force $\vb{F}_V$ is 
\begin{align*}
\vb{F_V} \vdot \vb{v}.
\end{align*}
Postulating conservation of energy we get
\begin{align*}
&\frac{d}{dt} \int_V dV \; \{ \frac{1}{2} \; \rho \; \vb{v}^2 + \rho \; e \}  \\ = - \int_S dA \; \{(\frac{1}{2} \rho \; \vb{v}^2& + \rho \; e) \; \vb{v} - \mathbf{ {\cal P}} \vdot \vb{v}  + \vb{q} \} \vdot \vb{n} + \int_V dV \; \vb{F} \vdot \vb{v}. 
\end{align*}
Using the divergence theorem, and letting the volume approach zero, we get the third fundamental equation of fluid dynamics
\begin{align}
\partial_t \; (\frac{1}{2} \; \rho \; \vb{v}^2 + \rho \; e) + \div{((\frac{1}{2} \; \rho \; \vb{v}^2 + \rho \; e) \; \vb{v} - \mathbf{ {\cal P}} \vdot \vb{v} + \vb{q})} = \vb{F}_V \vdot \vb{v}. \label{eq:39}
\end{align}
We can simplify (\ref{eq:39}) using the mass and momentum conservation equations. We have 
\begin{align}
	&\partial_t \; (\frac{1}{2} \; \rho \; \vb{v}^2 + \rho \; e) = \frac{1}{2} \; \partial_t \; \rho \; \vb{v}^2 + \rho \; \partial_t \; \vb{v} \vdot \vb{v} + \partial_t \; \rho \; e + \rho \; \partial_t \; e \nonumber\\
	&= \frac{1}{2} \; (-\vb{v} \vdot \grad \rho - \rho \; \grad \vdot \vb{v}) \; \vb{v}^2 + \rho \; (- \vb{v} \vdot \grad \vb{v} + \frac{\div{\mathbf{ {\cal P}}}}{\rho} + \frac{1}{\rho} \; \vb{F}_V ) \vdot \vb{v}
	\nonumber\\ &+ ( - \vb{v} \vdot \grad \rho - \rho \; \grad \vdot \vb{v}) \; e + \rho \; \partial_t \; e \nonumber \\
	&= - \frac{1}{2} \; \vb{v} \vdot \grad \rho \; \vb{v}^2 - \frac{1}{2} \; \rho \; \div{\vb{v}} \; \vb{v}^2 - \rho \; \vb{v} \vdot \grad \vb{v} \vdot \vb{v} \nonumber
	+ \div{\mathbf{ {\cal P}}} \vdot \vb{v} \nonumber \\ &+ \vb{F}_V \vdot \vb{v} - \vb{v} \vdot \grad \rho \; e - \rho \; \div{\vb{v}} \; e + \rho \; \partial_t \; e, \label{eq:40}
\end{align}
and
\begin{align}
	&\div{((\frac{1}{2} \; \rho \; \vb{v}^2 + \rho \; e) \; \vb{v} - \mathbf{ {\cal P}} \vdot \vb{v} + \vb{q})}\nonumber \\
	&=\frac{1}{2} \; \div{(\rho \; \vb{v}^2 \; \vb{v})} + \div{(\rho \; e \; \vb{v})} - \div{(\mathbf{ {\cal P}} \vdot \vb{v})} + \div{\vb{q}} \nonumber \\
	&=\frac{1}{2} \; \grad (\rho \; \vb{v}^2) \vdot \vb{v} + \frac{1}{2}\; \rho \; \vb{v}^2 \; \div{\vb{v}} + \grad (\rho \; e) \vdot \vb{v} \nonumber \\
	&+ \rho \; e \; \div{\vb{v}} - \div{(\mathbf{ {\cal P}} \vdot \vb{v})} + \div{\vb{q}} \nonumber \\
	&= \frac{1}{2} \; \vb{v} \vdot \grad \rho \; \vb{v}^2 + \rho \; \vb{v} \vdot \grad \vb{v} \vdot \vb{v} + \frac{1}{2} \; \rho \; \vb{v}^2 \; \grad \vdot \vb{v} \nonumber \\
	&+ \vb{v} \vdot \grad \rho \; e + \rho \; \vb{v} \vdot \grad e + \rho \; e \; \div{\vb{v}} - \div{(\mathbf{ {\cal P}} \vdot \vb{v})} + \div{\vb{q}}. \label{eq:41}
\end{align}
Using (\ref{eq:40}) and (\ref{eq:41}) in (\ref{eq:39}), gives us the third fundamental equation of fluid dynamics in the form
\begin{align*} 
\rho \; \partial_t \; e + \rho \; \vb{v} \vdot \grad e = \div{(\mathbf{ {\cal P}} \vdot \vb{v})} - (\div{\mathbf{ {\cal P}}}) \vdot \vb{v} - \div{\vb{q}}. 
\end{align*}
But
\begin{align}
\{ \div{({\cal P}\vdot \vb{v})} - (\div{ {\cal P}}) \vdot \vb{v} \}
&= \partial_{x_i} \; ( {\cal P}_{ij} \; v_j) -\partial_{x_i} \;  {\cal P}_{ij} \; v_j \nonumber \\
&= \partial_{x_i} \; {\cal P}_{ij} \; v_j + {\cal P}_{ij} \; \partial_{x_i} \; v_j - \partial_{x_i} \; {\cal P}_{ij} \; v_j \nonumber \\
&= {\cal P}_{ij} \; \partial_{x_i} \; v_j. \nonumber
\end{align}
Using the symbol
\begin{align*}
A:B = a_{ij} \; b_{ij},
\end{align*} 
for the full contraction of the tensors $A$ and $B$ we finally get the equation 
\begin{align}
\rho \; \partial_t \; e + \rho \; \vb{v} \vdot \grad e = \mathbf{ {\cal P}} : \grad \vb{v} - \div{\vb{q}}. \label{eq:42}
\end{align}
Using the material derivative,equation (\ref{eq:42}) can be written compactly as 
\begin{align}
\rho \; \frac{De}{Dt} = &- \div{\vb{q}} + \mathbf{ {\cal P}} : \grad \vb{v}\label{eq:43}
\end{align}
The  fundamental equations for fluid dynamics consists of (\ref{eq:18}),(\ref{eq:36.1}) and (\ref{eq:43})
\begin{align}
\frac{D\rho}{Dt} = &- \rho \; \div{\vb{v}}, \label{eq:44} \\
\rho \; \frac{D\vb{v}}{Dt} = &\div{\mathbf{ {\cal P}}} + \vb{F}_V, \label{eq:45}\\
\rho \; \frac{De}{Dt} = &- \div{\vb{q}} + \mathbf{ {\cal P}} : \grad \vb{v}. \label{eq:46}
\end{align}

\subsection{Closures}
The system (\ref{eq:44})-(\ref{eq:46}), consist of 5 equations for  14 unknowns, which are $\rho$, $e$, $\vb{v}$, $\vb{q}$ and the 6 independent components of $\mathbf{ {\cal P}}$.

In order to apply (\ref{eq:44}) - (\ref{eq:46}) to a particular problem, we must, for that particular problem, specify how $\mathbf{{\cal P}}$ and $\vb{q}$ depends on $\rho$, $\vb{v}$ and $e$. This will give us a closed system of equations that can be solved. Equations (\ref{eq:44}) - (\ref{eq:46}) are thus not the end of the story, further modeling is required to close the system. We will discuss two general closures leading to what are called \textit{ideal} and \textit{non-ideal} fluids.

\subsubsection{Ideal fluid}

Let us make the following three assumptions
\begin{align} 
&\mathbf{ {\cal P}}=-p \; I \; \; \; \; \; \; \; \; \; \; \; \; p = p(\vb{x}, t),  \label{eq:47}\\
&\vb{q} = 0, \label{eq:48}\\
 &\text{The fluid is in local thermodynamic equilibrium}. \label{eq:49}
\end{align}
Here $I$ is the $3\times 3$ identity matrix. Recall that by definition, the stress tensor ${\cal P}$ is 
\begin{align}
\mathbf{ {\cal P}} = - \int d\vb{c} \; m \; \vb{c} \vb{c} \; f,\label{eq:49.1}
\end{align}
so that (\ref{eq:47}) implies
\begin{align*}
p ={\cal P}_{ii} = \int d\vb{c} \; m \; c_{i}^2 \; f \geq 0 \; ,\; \; \; \; i = 1,2,3\;\;,
\end{align*}
where $p$, by definition, is the \textit{pressure} in the fluid. Since change in momentum per unit time is by definition \textit{force}, we observe from equation (\ref{eq:36.2}) that, using the assumption (\ref{eq:47}), the pressure induced force on a volume $V$ with surface $S$ is

\begin{figure}[htbp]
\centering
\includegraphics{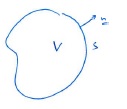}
\label{fig4}
\caption{}
\end{figure}

\begin{align*}
\int_V dV \;\vb{F}=\int_V dV \;\div{\mathbf{ {\cal P}}}=\int_S dA \; \mathbf{ {\cal P}} \vdot \vb{n} = - \int_S dA \; p \; \vb{n}.
\end{align*}
Note that this formula say that a {\it positive} pressure in the fluid acts as a {\it compressional} force, as it should. This fact is the motivation for the sign convention in the definition (\ref{eq:49.1}) of the stress tensor  for a fluid.

From statistical mechanics we know that a physical system is in thermodynamical equilibrium when the phase space distribution has settled down to the canonical distribution, or the \textit{Gibbs ensemble}, which it is also called. When the system is in thermodynamical equilibrium, it can in most cases be described in terms of five  parameters. These are the energy, entropy, volume, temperature and pressure. These parameters are usually denoted by $E$, $S$, $V$, $T$ and $p$. Only two of these parameters are independent, the rest can be determined by two chosen ones. For example if we let $p$ and $V$ be independent then $T=T(p,V)$, $E = E(p, V)$ and $S = S(p,V)$. These functional relationships constitute the equation of state for the system and contains everything there is to know about the system from a thermodynamical point of view. Equations of state are sometimes calculated from the Gibbs ensemble, but often just postulated on theoretical or empirical grounds.

Whichever way the equation of state is constructed it \textit{must} satisfy the following {\it fundamental thermodynamical relation}, which can be derived from the Gibbs ensemble
\begin{align}
T\; dS = dE + p \; dV. \label{eq:70}
\end{align}
If someone present you with an equation of state that does not obey (\ref{eq:70}), you just ask him or her to go back to the drawing board, their equation of state is wrong! 

In the assumption (\ref{eq:49}),  we do not assume that the fluid is in thermodynamical equilibrium, we rather ask that each small piece of it is. Thus the physical systems that are in thermodynamical equilibrium are the \textit{fluid elements}. I will not argue why this is a good assumption, but believe me; it is a very good assumption that almost always applies.

We have up to this point been a little vague about what a fluid element is, but we will now have to be a little more precise. We imagine a macroscopically small, but microscopically large volume, $V$, of fluid bounded by a surface, $S$. The important point here is that as time goes by, the fluid element and its bounding surface moves as determined by the velocity field $\vb{v}(\vb{x},t)$. The surface will deform, but no mass will flow across the surface exactly because the surface goes along with the flow. We say that $S$ is a \textit{material surface}.

\begin{figure}[htbp]
\centering
\includegraphics{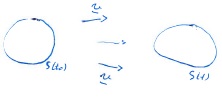}
\label{fig5}
\caption{}
\end{figure}

\noindent The fact that S is a material surface implies that the mass of the fluid element does not change. 

We will now derive the fundamental thermodynamical relation for fluid elements by applying (\ref{eq:70})
 to a fluid element.

Let $M$ be the mass of the fluid in the fluid element. Since $M$ does not change we have $dM=0$. We now introduce \textit{mass density}, $\rho$, \textit{entropy density}, $s$ and \textit{energy density}, $e$, through 
\begin{align}
\rho = \frac{M}{V},&& s=\frac{S}{M}, && e= \frac{E}{M}.  \label{eq:72}
\end{align}
Thus $\rho$ is mass per unit volume, $s$ is entropy per unit mass and $e$ is energy per unit mass. Note that $\rho$ and $e$ here coincide with quantities with the same names occurring in the fundamental equations of fluid dynamics (\ref{eq:44}) - (\ref{eq:46}). 

From (\ref{eq:70}), (\ref{eq:72}) and $dM=0$ we get
\begin{align}
T\;ds &= T\; d(\frac{S}{M}) = \frac{1}{M} \; T \; dS = \frac{1}{M} \; dE + \frac{1}{M} \; p \; dV  \nonumber\\ &= de + p \; d(\frac{1}{\rho}).\label{eq:73} 
\end{align}
Equation (\ref{eq:73}) is the fundamental thermodynamical relation for a fluid element. 

Observe that for an ideal fluid we have
\begin{align*}
\mathbf{ {\cal P}}: \grad \vb{v} = {\cal P}_{ij} \; \partial_{x_i} \; v_j = - p \; \delta_{ij} \; \partial_{x_i} \; v_j = - p \; \partial_{x_i} \; v_i = - p \; \div{\vb{v}}. 
\end{align*}
This, taken together with assumption (\ref{eq:48}), which states that there is not heat flow, implies that, for the ideal case, equation (\ref{eq:46}) simplifies into
\begin{align}
\rho \; \frac{De}{Dt} = - p \; \div{\vb{v}},  \label{eq:75}
\end{align}
But from equation (\ref{eq:44}) we have
\begin{align}
\frac{D\rho}{Dt} = - \rho \; \div{\vb{v}}.  \label{eq:76}
\end{align}
Combining (\ref{eq:75}) and (\ref{eq:76}) we get
\begin{align}
\frac{De}{Dt} + p \; \frac{D}{Dt} (\frac{1}{\rho}) = 0.  \label{eq:77}
\end{align}
Recall that from the definition of the material derivative, and the fact that the fluid element follows the fluid flow, we have that during a short time interval $dt$, between $t,t+dt$
\begin{align}
&\frac{De}{Dt} \; dt = \frac{d}{dt} \; e(\vb{x}(t),t) \; dt = de,  \label{eq: 78}\\
&\frac{D}{Dt}(\frac{1}{\rho})\; dt = \frac{d}{dt} \; (\frac{1}{\rho}) \; (\vb{x}(t), t) \; dt=d(\frac{1}{\rho}). \nonumber
\end{align}

\begin{figure}[htbp]
\centering
\includegraphics{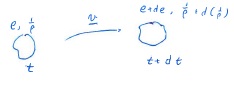}
\label{fig6}
\caption{}
\end{figure}

Therefore, for a fluid element in an ideal fluid, equations (\ref{eq:77}) and (\ref{eq:73}) imply that
\begin{align*}
T\; ds = de + p \; d(\frac{1}{\rho}) = (\frac{De}{Dt} + p \; \frac{D}{Dt} \; (\frac{1}{\rho}))\; dt = 0 \Rightarrow ds = 0. 
\end{align*}
Thus the entropy of the fluid element stays fixed. 

This kind of fluid flow is called \textit{isentropic}. Since each piece of fluid preserves whatever entropy it had at $t=0$, if the entropy was constant in space  $s(\vb{x}, 0) = s_0$ at $t=0$ then it will remain constant in space for all time. We will in general assume that the entropy is constant in space at $t=0$. Equation (\ref{eq:46}) will after this play no role in the dynamics and our fundamental system is reduced to 
\begin{align}
\frac{D\rho}{Dt} = - \rho \; \div{\vb{v}}  \label{eq:80} \\
\frac{D\vb{v}}{Dt} = - \frac{\grad p}{\rho} + \frac{1}{\rho} \; \vb{F}_V  \label{eq:81}
\end{align}
because
\begin{align*}
(\div{\mathbf{ {\cal P}}})_j = \partial_{x_i} \;{\cal  P}_{ij} = - \partial_{x_i} \; (p \; \delta_{ij}) = -\partial_{x_j} \; p = - (\grad p)_j.
\end{align*}
Using $p$ and $s$ as independent variables we have in general $\rho = \rho(p, s)$, but since $s=s_0$ is a constant we can simplify this into
\begin{align}
\rho = \rho(p), \label{eq:83}
\end{align}
and (\ref{eq:80}), (\ref{eq:81}) and (\ref{eq:83}) is a closed system of five equations for five unknowns, $\rho$, $\vb{v}$ and $p$. These are the equations for an \textit{ideal fluid}. 

\noindent If the fluid is a liquid we know that it is essentially incompressible. This means that the mass density of each fluid element is constant in time as it is transported along with the flow, thus $\frac{D\rho}{Dt}=0$. This is equivalent to the mathematical condition $\div{\vb{v}}=0$. If the density of the fluid is constant in space for $t=0$, incompressibility will imply that it will remain constant in space for all time, thus $\rho = \rho_0$ is a constant. For this case (\ref{eq:80})-(\ref{eq:83}) reduces to the system
\begin{align}
\frac{D\vb{v}}{Dt} = - \frac{\grad p}{\rho_0} + \frac{1}{\rho_0} \; \vb{F}_V, \label{eq:84}\\
\div{\vb{v}} = 0. \label{eq:85}
\end{align}
This is a closed system of four equations for the four unknowns $\rho$ and $\vb{v}$. These are the equations for an \textit{ideal liquid}. Equations (\ref{eq:84}), (\ref{eq:85}) were first published by Leonard Euler in 1757, and are in his honor called the \textit{Euler equations}. These equations are in most situations a very good model for ordinary water.

  Note that an important case when incompressibility does not imply constant density is if the fluid consists of two different components, each with a different density. For this case the density at $t=0$ does vary in space.

  You might at this point be confused about what looks like a contradiction. Equation (\ref{eq:83}) says that $\rho$ is a function of $p$ and for a liquid $\rho = \rho_0$ is constant but $p$ is not! How can this be? The explanation is simple, $\rho$ is not actually constant but depends on $p$ so weakly that even in the deepest part of the ocean, where the pressure is immense, $\rho$ is still very close to $\rho_0$. In effect, the equation of state has the form 
\begin{align*}
\rho = \rho_0 + \epsilon \; p && \epsilon\;p << \rho_0.
\end{align*}
Observe that, through the material derivative, the Euler equations are \textit{nonlinear}. The fact that the basic equations for an ideal liquid are nonlinear makes this theory \textit{hard}.

\subsubsection{Non-ideal fluid}

We make the following closure assumptions
\begin{align}
& \mathcal{P} = -p \; I + \tau, \label{eq:86}\\
&\vb{q} = - \kappa \; \grad T \; \; \; \; \; \; \; \; \; \kappa > 0, \label{eq:87}\\
&\text{The fluid is in local thermodynamic equilibrium}.\label{eq:88}
\end{align}
In (\ref{eq:87}) the parameter $\kappa$ is positive and is called the \textit{heat conductivity} of the fluid. The equation expresses the observational fact that heat flows from hot to cold objects. 

In (\ref{eq:86}) one could at the outset assume that 
\begin{align} 
\tau = \tau(\vb{v}, \grad{\vb{v}}, ..). \label{eq:89}
\end{align}
However physical laws must look the same to all inertial observers, and this means that $\tau$ can not depend on $\vb{v}$. If we assume that the gradients are not too large, we can linearize (\ref{eq:89}) and write
\begin{align*}
\tau_{ij} = A_{ijkl} \; \partial_{x_k} \; v_l,
\end{align*}
where $A_{ijkl}$ is a tensor of rank 4. For all gases and most liquids, the relation between $\tau$ and $\grad v$ must be rotationally invariant. There is no preferred direction for such fluids and they are said to be \textit{isotropic}. At the mathematical level this means that the tensor $A_{ijkl}$ is \textit{invariant}. It has the same components in all rotated Cartesian  frames. Thus for any change of coordinates
\begin{align*}
x_i' = a_{ir} \; x_r, 
\end{align*}
we have
\begin{align*}
A_{ijkl}' = A_{ijkl}.
\end{align*}
and the requirement for invariance is
\begin{align}
 A_{ijkl}=a_{it} \; a_{ju} \; a_{kn} \; a_{lm} \; A_{tunm}.  \label{eq:90}
\end{align}
We are not going to solve the tensor equation (\ref{eq:90}), this kind of analysis belongs in a more specialized course on tensors. Here we merely state the result. The solutions of (\ref{eq:90}) are of the form
\begin{align*}
A_{ijkl} = \eta \; (\delta_{ik} \; \delta_{jl} + \delta_{il} \; \delta_{jk} -  \frac{2}{3} \; \delta_{ij} \; \delta_{kl}) + \xi \; \delta_{ij} \; \delta_{kl},
\end{align*}
where $\xi=\xi(\vb{x})$ and $\eta=\eta(\vb{x})$ are two free parameters that are scalar fields in the tensorial sense. 

With this solution for $A_{ijkl}$, the relation between $\tau$ and $\grad \vb{v}$, using dyadic notation, is 
\begin{align*}
\tau = \eta \; (\grad \vb{v} + \grad  \vb{v}^t -  \frac{2}{3} \; I \; \div{\vb{v}}) + \xi \; I \; \div{\vb{v}}.
\end{align*}
The equation for $\rho$ is the same as for the ideal case
\begin{align}
\frac{D\rho}{Dt} = - \rho \; \div{\vb{v}}. \label{eq:91.1}
\end{align}
In order to write down the equation for $\vb{v}$ we need $\div{\tau}$. If we assume that the scalar fields $\xi(\vb{x})$ and $\eta(\vb{x})$ are constants, independent of $\vb{x}$, which is a good approximation in many cases, we have
\begin{align}
&(\div{\tau})_j = \partial_{x_i} \; \tau_{ij} = \eta \; \{ \partial_{x_i} \; \partial_{x_i} \; v_j + \partial_{x_i} \; \partial_{x_j} \; v_i - \frac{2}{3} \; \delta_{ij} \; \partial_{x_i} \; \partial_{x_k} \; v_k \} \nonumber \\ &+ \xi \; \delta_{ij} \; \partial_{x_i} \; \partial_{x_k} \; v_k \nonumber \\
&= \frac{1}{3} \; \eta \; \partial_{x_j} \; \partial_{x_k} \; v_k + \eta \; \partial_{x_i} \; \partial_{x_i} \; v_j + \xi \; \partial_{x_j} \; \partial_{x_k} \; v_k. \nonumber
\end{align}
Thus in dyadic notation we have
\begin{align*}
\div{\tau} = \eta \; \grad^2 \vb{v} + ( \frac{1}{3} \; \eta+\xi) \; \grad (\div{\vb{v}}),
\end{align*}
and the equation for $\vb{v}$ in the non-ideal case can be written as
\begin{align*}
\rho \; \frac{D\vb{v}}{Dt} = - \grad p + \eta \; \grad^2 \vb{v} + (\frac{1}{3} \; \eta+\xi) \; \grad (\div{\vb{v}}) + \vb{F}_V.
\end{align*}
Let us finally consider the equation for the energy (\ref{eq:46}). We need to calculate the full contraction $\tau:\grad \vb{v}$.

We have 
\begin{align}
&\tau : \grad v = \tau_{ij} \; \partial_{x_i} \; v_j = \eta \; (\partial_{x_i} \; v_j  \\ & + \partial_{x_j} \; v_i - \frac{2}{3} \; \delta_{ij} \; \partial_{x_k} \; v_k)\; \partial_{x_i} \; v_j + \xi \; \delta_{ij} \; \partial_{x_k} \; v_k \; \partial_{x_i} \; v_j \nonumber\\
&= \eta \; (\partial_{x_i} \; v_j \; \partial_{x_i} \; v_j + \partial_{x_j} \; v_i \; \partial_{x_i} \; v_j -  \frac{2}{3} \; \partial_{x_i} \; v_i \; \partial_{x_k} \; v_k) \nonumber \\ &+ \xi \; \partial_{x_i} \; v_i \; \partial_{x_k} \; v_k, \label{eq:91}
\end{align}
For any number $a$, define a 2-tensor $Q$  by
\begin{equation*}
Q_{ij}=\partial_{x_i} \; v_j+\partial_{x_j} \; v_i  - a \; \delta_{ij} \; \partial_{x_k} \; v_k.
\end{equation*}
Observe that for any $a$ we have $Q^2\equiv Q:Q=Q_{ij}Q_{ij}\geq 0$.

We have
\begin{align}
&Q^2=
 \partial_{x_j} \; v_i \; \partial_{x_j} \; v_i + 2 \; \partial_{x_i} \; v_j \; \partial_{x_j} \; v_i  \nonumber \\ &- 2 \; a \; \delta_{ij} \; \partial_{x_j} \; v_i \; \partial_{x_k} \; v_k \nonumber \\
&+ \partial_{x_i} \; v_j \; \partial_{x_i} \; v_j - 2 \; a \; \delta_{ij} \; \partial_{x_i} \; v_j \; \partial_{x_k} \; v_k + a^2 \; \delta_{ij} \; \delta_{ij} \; \partial_{x_k} \; v_k \; \partial_{x_l} \; v_l \nonumber  \\
&= 2 \; (\partial_{x_j} \; v_i \; \partial_{x_j} \; v_i + \partial_{x_i} \; v_j \; \partial_{x_j} \; v_i - ( 2 \; a - \frac{3}{2} \; a^2) \; \partial_{x_k} \; v_k \; \partial_{x_l} \; v_l), \label{eq:92}
\end{align}
Choose $a$ such that $2\; a - \frac{3}{2}\; a^2 =  \frac{2}{3}$ or in other words choose $a = \frac{1}{3} \; (2+ \sqrt{2})$. Then (\ref{eq:91}) and (\ref{eq:92}) imply that
\begin{align*} 
\tau : \grad \vb{v} =\frac{1}{2}\eta\;  Q^2 + \xi \;(\div{\bf v})^2.
\end{align*}
Using the same thermodynamic arguments as under the ideal case, we get the energy equation, in the non-ideal case, on the form
\begin{align}
\rho \; T \; \frac{Ds}{Dt} = \div{( \kappa \;\grad T)}+ \frac{1}{2}\eta\;  Q^2 +  \xi \;(\div{\bf v})^2. \label{eq:93}
\end{align}
Let $V$ be a fluid volume with bounding surface $S$. The entropy per unit volume is $\rho \; s$. Therefore the total entropy in $V$ is 
\begin{align*}
S(t) = \int_V \rho \; s \; dV,
\end{align*}
which implies that
\begin{align}
\frac{dS}{dt} = \int_V dV \; \partial_t \; (\rho \; s). \label{eq:94}
\end{align}
Using the definition of the material derivative and equation (\ref{eq:91.1}) we have
\begin{align}
\partial_t \; (\rho \; s) = - \div{(\rho \; s \; \vb{v})} + \rho \; \frac{Ds}{Dt}. \label{eq:95}
\end{align}
Therefore, upon using (\ref{eq:93}) and (\ref{eq:95}) in (\ref{eq:94}), we get
\begin{align}
&\frac{dS}{dt} = \int_V dV \{ - \div{(\rho \; s \; \vb{v})} + \div{(\frac{\kappa}{T} \; \grad T)} + \frac{\kappa}{T^2} \; (\grad T)^2 \nonumber\\
&+ \frac{1}{2} \; \frac{\eta}{T} \; Q^2 + \frac{\xi}{T} \; ( \div{\vb{v}}) ^2 \} \nonumber \\
&= \int_S dA \; \frac{\kappa}{T} \; \grad T \vdot \vb{n} - \int_S dA \; \rho s \; \vb{v} \vdot \vb{n} + \int_V dV\; \theta, \label{eq:96}
\end{align}
where $\theta \geq 0$.

  Let us assume that $S$ is a material surface that is perfectly insulated. This means that no mass or energy can pass through $S$. Said in another way; the volume $V$ is a \textit{closed} system. For a closed system we have 
\begin{align*}
&\kappa \; \grad T \vdot \vb{n} = \vb{q} \vdot \vb{n} = 0 && \text{(no energy flow)}, \\
&\vb{v} \vdot \vb{n} = 0 && \text{(no mass flow)}, 
\end{align*}
and therefore, (\ref{eq:96}) implies that 
\begin{align*}
\frac{dS}{dt} = \int_V dV \; \theta \geq 0. 
\end{align*}
The entropy in a closed body of fluid can never decrease. This expresses the fact that non-ideal fluids are \textit{irreversible systems} in the thermodynamic sense.

  Summing up, the fundamental system of equations for a non-ideal fluid is 
\begin{align}
\frac{D\rho}{Dt} &= - \rho \; \div{\vb{v}}, \label{eq:97} \\
\rho \; \frac{D\vb{v}}{Dt} &= - \grad  p + \eta \; \grad^2 \vb{v} + (\xi + \frac{1}{3} \; \eta) \; \grad  (\div{\vb{v}}) + \vb{F}_V, \label{eq:98}\\
\rho \; T \; \frac{Ds}{Dt}& = \div{( \kappa \;\grad T)}+ \frac{1}{2}\eta\;  Q^2 +  \xi \;(\div{\bf v})^2,\label{eq:99}\\
\rho &= \rho(p, s), \; \; \; \; \; \; \; \; \; \; \; \; \; \; \; \; \; T = T(p, s), \label{eq:100}
\end{align}
where (\ref{eq:100}), are the equations of state for a non-ideal fluid which must be supplied in order to close the system; six equations for six unknowns. A single component  non-ideal liquid satisfies the additional requirement of incompressibility, $\div{\vb{v}}=0$, which for a single component fluid implies that $\rho = \rho_0$. For such fluids  (\ref{eq:97}) and (\ref{eq:98}) reduce to 
\begin{align}
\rho_0 \; \frac{D\vb{v}}{Dt} = - \grad p + \eta \; \grad^2 \vb{v} + \vb{F}_V. \label{eq:101}\\
\div{\vb{v}} = 0. \label{eq:102}
\end{align}
They decouple from (\ref{eq:99}) and (\ref{eq:100}), which can be used to calculate $T$ and $s$ after we have found $\vb{v}$ and $p$ from (\ref{eq:101}) and (\ref{eq:102}). The system
(\ref{eq:101}) and (\ref{eq:102}) is called the \textit{Navier-Stokes equations} and where first published by Claude Navier in 1822. 

It is worth noting that if $\vb{v}(\vb{x},t),p(\vb{x},t)$ is a solution to the Navier-Stokes equations, then $\vb{v}(\vb{x},t),p(\vb{x},t)+\alpha(t)$ is also a solution for any function $\alpha(t)$. This arbitrariness in the definition of the pressure is not a feature of the general equations for non-ideal fluids,  but is a consequence of the assumption of incompressibility for non-ideal liquids.

The Navier-Stokes equations are believed to be a very precise model for water. The equations describe an incredibly rich array of phenomena in water and other liquids. From a mathematical point of view they are \textit{very hard} equations to solve. In fact it is not even known if they are well posed. Prize money in the order of a million dollars goes to the first person that can solve this problem.

\subsection{Simple fluid systems}

We will now apply the equations of fluid dynamics to some simple systems. 

\subsubsection{Static fluid}

Let us consider the case of a static fluid in a constant gravitational field. We assume thus that 
\begin{align*}
\vb{v} = 0,\\
\vb{F}_V = \rho\; \vb{g}, 
\end{align*}
where $\vb{g}$ is a constant vector.

The first two equations of fluid dynamics (\ref{eq:97}), (\ref{eq:98}) decouple from the third (\ref{eq:99}) and we get 
\begin{align}
\rho = \rho(\vb{x}),\nonumber\\
-\grad  p + \rho \; \vb{g} = 0. \label{eq:103}
\end{align}
Let us apply (\ref{eq:103}) to a volume $V$ with bounding surface $S$

\begin{figure}[htbp]
\centering
\includegraphics{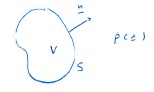}
\label{fig7}
\caption{}
\end{figure}

We have 
\begin{align}
\int_V\;dV\grad p &= \int_V dV \; \rho \; \vb{g}, \nonumber \\
&\Updownarrow\nonumber\\
\int_S dA \; p \; \vb{n} &= \int_V dV \; \rho \; \vb{g}, \label{eq:104}
\end{align}
where we have used the divergence theorem and the fact that 
\begin{align*}
\grad  p = \div{(p \; I)}.
\end{align*}
Equation (\ref{eq:104}) tells us that, in a stationary fluid, the pressure-induced force on a volume of the fluid is exactly equal to the weight of the fluid in the volume. Observe that the pressure-induced force on the volume $V$ does not actually depend on whatever is inside the volume $V$, it only depends on the pressure in the surrounding fluid. Thus the \textit{same} pressure-induced force acts on \textit{any} object with bounding surface $S$ that is immersed in the fluid. 

This is the well known \textit{Archimedean principle}. It explains why helium balloons rise in the air, why stones sink in water and how come a boat made of iron can float.

Let the surface of a calm sea be defined by $z=0$ 

\begin{figure}[htbp]
\centering
\includegraphics{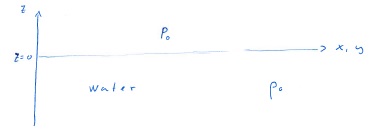}
\label{fig8}
\caption{}
\end{figure}

\noindent The pressure at the surface is $p_o$. Since pure water is an incompressible single component fluid, it has constant density $\rho_0$ independent of the depth. The force of gravity is constant and points vertically down
\begin{align*}
\vb{g} = - g \; e_z && g>0, 
\end{align*}
where $g$ is the strength of the gravitational force at sea-level. From equation (\ref{eq:103}) we get 
\begin{align}
\partial_x \; p &= 0, \nonumber\\
\partial_y \; p &= 0, \nonumber\\
\partial_z \; p &= -\rho_0 \; g. \label{eq:105}
\end{align}
Thus $ p = p(z)$ and (\ref{eq:105}) implies that
\begin{align*}
p(z) = - \rho_0 \; g \; z + c.
\end{align*}
But $p=p_0$ at $z=0$. Thus $c=p_0$ and we have
\begin{align}
p(z) = p_0 - \rho_0 \; g \; z.\label{eq:106}
\end{align}
The pressure is a linear function of the water depth and increases as we descend into deeper water. Anyone that has done any amount of diving knows this. 

Defining the pressure at sea-level to be one atmosphere (1 atm), and letting $z_n$ be the depth at which the pressure is $n$ atm, we have from (\ref{eq:106})
\begin{align*} 
n \; p_0 = p_0 - \rho \; g \; z_n \Rightarrow z_n = - (n-1) \; \frac{p_0}{g\; \rho_0}.
\end{align*}
For water we have 
\begin{align*}
\frac{p_0}{g\; \rho_0} \approx 10m.
\end{align*}
Thus the pressure increases by 1 atm for every 10 meters of depth. This shows how heavy water is compared to air. The 1 atm of pressure at sea level is the weight of an air column that is more than 120 km tall. This weight is the same as for a cylinder of water that is merely 10m tall! 

\subsubsection{The Bernoulli equation}

Let us consider an incompressible fluid of constant density $\rho_0$ in a constant gravity field, which we without loss of generality can assume to point vertically down. Choosing the z-axis in the positive vertical direction we have 
\begin{align*}
\vb{g} = - g \; e_z.
\end{align*}
We will assume that the velocity field of the fluid is \textit{stationary} 
\begin{align*}
\vb{v}(\vb{x}, t) = \vb{v}(x).
\end{align*}
The equation for the velocity field (\ref{eq:98}) then simplifies into
\begin{align}
\vb{v} \vdot \grad \vb{v} = - \frac{1}{\rho_0} \; \grad  p - g \; e_z. \label{eq:107}
\end{align}
We have here assumed that the fluid is ideal.
A \textit{streamline} is the path of a fluid element and plays an important role in the theory of stationary flows. A streamline is then by definition determined by 
\begin{align*}
\frac{d\vb{x}}{ds} = \vb{v}(\vb{x}), 
\end{align*}
where $s$ parametrizes the streamline. Note that we have
\begin{align}
\vb{v} \vdot \grad  \vb{v}& = \grad (\frac{1}{2} \; \vb{v}^2) - \vb{v} \cp (\curl{\vb{v}}), \label{eq:108} \\
g \; e_z& = \grad  (g \; z). \label{eq:109}
\end{align}
Using (\ref{eq:108}) and (\ref{eq:109}), equation (\ref{eq:107}) can be rewritten as 
\begin{align*}
\grad ( \frac{1}{2} \; \vb{v}^2 + \frac{p}{\rho_0} + g \; z) = \vb{v} \cp (\curl{\vb{V}}). 
\end{align*}
Let
\begin{align*}
A(\vb{x}) = \frac{1}{2} \; \vb{v}^2 + \frac{p}{\rho_0} + g \; z. 
\end{align*}
Along a streamline, $A(s) = A(\vb{x}(s))$ changes according to 
\begin{align*}
\frac{dA}{ds} = \frac{d \vb{x}}{ds} \vdot \grad  A = \vb{v} \vdot (\vb{v} \cp (\curl{\vb{v}})) = 0,
\end{align*}
and therefore $A$ is \textit{constant} along streamlines. 
\begin{align*}
\frac{1}{2} \; \vb{v}^2 + \frac{p}{\rho_0} + g \; z = const 
\end{align*} 
This is the \textit{Bernoulli's law}. It explains (naively) many simple fluid phenomena, for example why a metal plane that is heavier than air very rarely falls out of the sky!

\begin{figure}[htbp]
\centering
\includegraphics{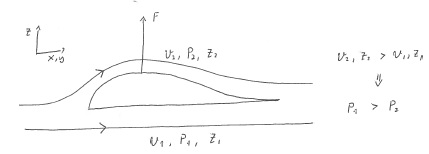}
\label{fig9}
\caption{}
\end{figure}

\subsubsection{Sound waves}

For an ideal fluid, equations (\ref{eq:80}), (\ref{eq:81}), (\ref{eq:83}) defines the correct model.
\begin{align}
\partial_t \; \rho + \vb{v} \vdot \grad  \rho + \rho \; \grad \vdot \vb{v} = 0, \label{eq:110} \\
\partial_t \; \vb{v} + \vb{v} \vdot \grad  \vb{v} + \frac{\grad  p}{\rho} = 0, \label{eq:111}\\
p = p(\rho). \label{eq:112}
\end{align}
Note that, using (\ref{eq:110}), we have
\begin{align*}
\frac{Dp}{Dt} = \frac{d}{dt} \; p(\rho(t, \vb{v}(t))) = \frac{dp}{d\rho} \; \frac{d}{dt} \; \rho(t, \vb{x}(t)) = \frac{dp}{d\rho} \;\frac{D\rho}{Dt}&= - \rho \; \frac{dp}{d\rho} \; \div{\vb{v}}, \\ &\Updownarrow \\  \partial_t \; p + \vb{v} \vdot \grad  p + \rho \; a^2(\rho) \; \div{\vb{v}} &= 0, \nonumber 
\end{align*} 
where we have defined
\begin{align*}
a^2(\rho) = \frac{dp}{d\rho}. 
\end{align*}
For all reasonable equations of state
$\frac{dp}{d\rho} > 0$. Thus,  pressure increases when density increases.

We can now dispense with the equation of state (\ref{eq:112}) and rather consider the system 
\begin{align} 
\partial_t \; \rho + \vb{v} \vdot \; \grad  \rho + \rho \; \div{\vb{v}} = 0, \label{eq:113}\\
\partial_t \vb{v} + \vb{v} \vdot \grad  \vb{v} + \frac{1}{\rho} \; \grad  p = 0, \label{eq:114}\\
\partial_t \; p + \vb{v} \vdot \grad \ p + \rho \; a^2(\rho) \; \div{\vb{v}} = 0.  \label{eq:115}
\end{align}
This system has a simple solution 
\begin{align*} 
\vb{v} = 0, && \rho = \rho_0, && p = p_0,
\end{align*}
corresponding to a static homogeneous fluid. We will now investigate small disturbances of this state
\begin{align*} 
\vb{v} = \vb{v}',&& \rho = \rho_0 + \rho', && p = p_0 + p',
\end{align*}
where 
\begin{align}
\norm{\vb{v}^2} << 1, && \frac{\abs{\rho - \rho_0}}{\rho_0} << 1, && \frac{\abs{ p - p_0}}{p_0} << 1. \label{eq:116}
\end{align}
Inserting (\ref{eq:116}) into the fluid equations (\ref{eq:113}) - (\ref{eq:115}) we get 
\begin{align*} 
\partial_t \; \rho' + \vb{v}' \vdot \grad \rho' + ( \rho_0 + \rho') \; \div\vb{v}' = 0, \\
\partial_t \; \vb{v}' + \vb{v}' \vdot \grad \vb{v}' + \frac{\grad  p'}{\rho_0 + \rho'} = 0, \\ \partial_t \; p' + \vb{v}' \vdot \grad p' + (\rho_0 + \rho') \; a^2(\rho_0 + \rho') \; \div{\vb{v}'} = 0. 
\end{align*}
We now \textit{linearize} this system by dropping terms that contain products of small quantities. This gives a linear system of equations for the small disturbances 
\begin{align}
\partial_t \; \rho' + \rho_0 \; \div{\vb{v}'} = 0,  \nonumber\\
\partial_t \; \vb{v}' + \frac{1}{\rho_0} \; \grad p' = 0, \label{eq:117}\\ 
\partial_t \; p' + \rho_0 \; a_0^2 \; \div{\vb{v}'} = 0, \label{eq:118}
\end{align}
where
\begin{align*} 
a_0^2 = a^2(\rho_0) = \eval{\frac{dp}{d\rho}}_{\rho = \rho_0}.
\end{align*}
These equations are the starting point for most work in \textit{acoustics}, which is the science of sound. Small disturbances like these, in air, is perceived as sound by humans and other animals.
Observe that using (\ref{eq:117}), (\ref{eq:118}), we get 
\begin{align*} 
\partial_{tt} \; p' = - \rho_0 \; a_0^2 \; \partial_t \; (\div{\vb{v}'}) &= - \rho_0 \; a_0^2 \; \div{(\partial_t \; \vb{v}')}  \\ &= - \rho_0 \; a_0^2 \; \div{(-\frac{1}{\rho_0} \; \grad p')}, 
\end{align*}
and thus we get 
\begin{align*}
\partial_{tt} \; p' - a_0^2 \; \grad^2 p' = 0. 
\end{align*}
This is the wave equation. We know that this equation has a speed limit of 
\begin{align*}
c = a_0. 
\end{align*}
Therefore, small pressure disturbances will propagate through the fluid at speeds less or equal to $c = a_0$. 

Thus the speed of such disturbances, which is called the \textit{sound speed} depends on the physical properties of the fluid through the equations of state. 

\subsubsection{Potential flow}

The motion of an ideal incompressible fluid is described by the Euler equations 
\begin{align}
\partial_t \; \vb{v} + \vb{v} \vdot \grad  \vb{v} = - \frac{1}{\rho_0} \; \grad  p + \frac{1}{\rho_0} \; \vb{F}_V, \nonumber\\
\div{\vb{v}} = 0. \label{eq:167}
\end{align}
The volume force is the force of gravity, which is assumed to be constant
\begin{align*} 
\vb{F}_V = \rho_0 \; \vb{g}. 
\end{align*}
This is for example true for oceans on the surface of the earth. 

We choose a coordinate system with z-axis along the vertical, where the vertical is determined by the force of gravity. 
\begin{figure}[htbp]
\centering
\includegraphics{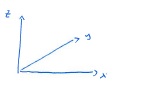}
\label{fig10}
\caption{}
\end{figure}

\noindent  So in these coordinates 
\begin{align*}
\vb{F}_V = - \rho_0 \; g \; e_z, 
\end{align*}
where $e_z$ is the unit vector in the z-direction. Let the \textit{vorticity} of the fluid velocity field be defined by 
\begin{align*} 
\vb{w} = \curl{\vb{v}}. 
\end{align*}
Using the vector identity 
\begin{align*} 
(\curl{\vb{v}}) \cp \vb{v} = -\frac{1}{2} \; \grad \vb{v}^2 + \vb{v} \vdot \grad  \vb{v},
\end{align*}
the Euler equation (\ref{eq:167}) can be rewritten as 
\begin{align}
\partial_t \; \vb{v} + \grad  (\frac{1}{2} \; \vb{v}^2) + (\curl{\vb{v}}) \cp \vb{v} = - \grad (\frac{p}{\rho_0}) -\grad( g \; z). \label{eq:172}
\end{align}
Taking the curl of (\ref{eq:172}), and using the fact that the curl of a gradient is zero we get 
\begin{align}
\partial_t \; \vb{w} + \curl{(\vb{w} \cp \vb{v})} = 0. \label{eq:173}
\end{align}
Using cartesian tensors we have 
\begin{align}
&[\curl{(\vb{w} \cp \vb{v})}]_i = \epsilon_{ijk} \; \partial_{x_j} \; (\vb{w} \cp \vb{v})_k \nonumber \\
&=\epsilon_{ijk} \; \partial_{x_j} \; (\epsilon_{kln} \; w_l \; v_n) \nonumber \\
&= \epsilon_{ijk} \; \epsilon_{kln}\; \partial_{x_j} \; (w_l \; v_n) \nonumber \\
&= \epsilon_{ijk} \; \epsilon_{lnk} \; ( \partial_{x_j} \; w_l ) \; v_n + \epsilon_{ijk} \; \epsilon_{lnk} \; w_l \; \partial_{x_j} \; v_n \nonumber \\
&= \delta_{il} \; \delta_{jn} \; \partial_{x_j} \; w_l \; v_n - \delta_{in} \; \delta_{jl} \; \partial_{x_j} \; w_l \; v_n \nonumber \\ &+ \delta_{il} \; \delta_{jn} \; w_l \; \partial_{x_j} \; v_n - \delta_{in} \; \delta_{jl} \; w_l \; \partial_{x_j} \; v_n \nonumber \\ &= \partial_{x_n} \; w_i \; v_n - \partial_{x_j} \; w_j \; \; v_i + w_i \; \partial_{x_j} \; v_j - w_j \; \partial_{x_j} \; v_i. \nonumber
\end{align}
Thus in dyadic notation we have 
\begin{align} 
\curl{(\vb{w} \cp \vb{v})} = \vb{v} \vdot \grad \vb{w} - \vb{w} \vdot \grad  \vb{v} - (\div{\vb{w}}) \; \vb{v} + (\div{\vb{v}}) \; \vb{w}. \label{eq:174}
\end{align}
But $\div{\vb{v}} = 0$ because the flow is incompressible and $ \div{\vb{w}} = 0$ because $\vb{w} = \curl{\vb{v}}$ is a curl. Therefore, (\ref{eq:174}) simplifies into 
\begin{align} 
\curl{(\vb{w} \cp \vb{v})} = \vb{v} \vdot \grad  \vb{w} - \vb{w} \vdot \grad \vb{v}. \label{eq:175}
\end{align}
Inserting (\ref{eq:175}) into (\ref{eq:173}), and using the definition of the material derivative, we get
\begin{align*} 
\frac{D \vb{w}}{Dt} = \vb{w} \vdot \grad \vb{v}. 
\end{align*}
This equation says that if $\vb{w} = 0 $ for a fluid element at $t=0$ then it will remain so for all $t>0$. It is therefore consistent with the Euler equation to seek solutions that are vorticity free, or \textit{irrotational}, 
\begin{align*}
\curl{\vb{v}} = 0.
\end{align*}
For such solutions, the velocity field can be described in terms of a \textit{velocity potential}, $\phi$,
\begin{align}
\vb{v} = \grad \phi. \label{eq:178}
\end{align}
For this reason, irrotational flows are also \textit{potential flows}.
Inserting (\ref{eq:178}) into the Euler equation (\ref{eq:172}) we get 
\begin{align}
\partial_t \; \grad  \phi + \grad (\frac{1}{2} \; (\grad  \phi)^2) &= - \grad  (\frac{p}{\rho_0} + g \; z), \nonumber \\ 
&\Updownarrow \nonumber\\
\grad (\partial_t \; \phi + \frac{1}{2} \; (\grad  \phi)^2 + \frac{p}{\rho_0} + g \; z) &= 0, \nonumber \\
&\Updownarrow \nonumber\\
\partial_t \; \phi + \frac{1}{2} \; (\grad \phi)^2 + \frac{p}{\rho_0} + g \; z &= \alpha(t), \label{eq:179}
\end{align}
where $\alpha(t)$ is an arbitrary function of time only. Equation (\ref{eq:179}) determines the pressure in the fluid in terms of $\phi$ and $\alpha$. 
\begin{align} 
p = \rho_0 ( - \partial_t \; \phi  - \frac{1}{2} \; (\grad \phi) ^2 - g \; z + \alpha (t)). \label{eq:180}
\end{align}
Equation (\ref{eq:178}), taken together with the condition of incompressibility in (\ref{eq:167}), implies that 
\begin{align*}
\div{\grad \phi} &= 0,\nonumber\\
 &\Updownarrow \nonumber \\ \grad^2 \phi &= 0.
\end{align*}
The space-time dependent function $\phi(\vb{x}, t)$ thus satisfy the Laplace equation. For later use, we separate out an arbitrary constant from $\alpha(t)$ and write (\ref{eq:180}) in the form 
\begin{align*}
\frac{p-p_0}{\rho_0} = - \phi_t - \frac{1}{2} \; (\grad \phi)^2 - g\; z + \alpha(t).
\end{align*}
Note that from (\ref{eq:178}) ,it is clear that $\phi$ is not uniquely determined by $\vb{v}$, we can add an arbitrary constant of time, $\rho(t)$, to $\phi$ without changing $\vb{v}$
\begin{align*}
\vb{v} = \grad \phi' && \phi' = \phi + \rho(t).
\end{align*}
With this choice of the potential we have 
\begin{align*} 
\frac{p-p_0}{\rho_0} = - \phi_t' - \frac{1}{2} \; (\grad \phi')^2 - g \; z + \alpha(t) = -\phi_t - \frac{1}{2} \; (\grad \phi)^2 - g \; z - \rho'+ \alpha.
\end{align*}
Choosing $\rho' = \alpha$ we see that the potential can always be chosen so as to eliminate $\alpha(t)$. We will usually do this and conclude that irrotational solutions to the Euler equations are determined by 
\begin{align}
&\grad^2 \phi(\vb{x}, t ) = 0,  \label{eq:185}\\
&\frac{p - p_0}{\rho_0} = - \phi_t (\vb{x}, t) - \frac{1}{2} \; (\grad \phi)^2(\vb{x}, t) - g \; z, \label{eq:186} \\
&\vb{v}(\vb{x}, t) = \grad \phi(\vb{x}, t). \label{eq:187}
\end{align}
The only equation we actually need to solve is the Laplace equation (\ref{eq:185}). Given $\phi$, (\ref{eq:186}) and \ref{eq:187} determine $p$ and $\vb{v}$ in terms of $\phi$. The equations (\ref{eq:185}), (\ref{eq:186}) and (\ref{eq:187}) are the fundamental equations for potential flows.

\subsection{Surface waves}

Ocean surface waves are the best known wave phenomena there is, in fact our basic intuition on the behavior of waves comes historically from the formal and informal study of ocean waves.
\subsubsection{Surface waves for irrotations flows}
\begin{figure}[htbp]
\centering
\includegraphics{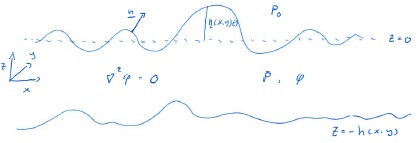}
\caption{}
\label{fig11.1}
\end{figure}

  In order to get our modeling off the ground, certain simplifying assumptions will have to be made. The geometry of the situation is as illustrated in figure \ref{fig11.1}

\noindent Our first and most important restriction is to consider only potential flows. The second restriction is that there are no breaking waves. This means that the fluid surface is the graph of a function
\begin{align*}
z = \eta(x,y,t). 
\end{align*}
For breaking waves this would not be the case, and a different treatment is needed. 

\begin{figure}[htbp]
\centering
\includegraphics{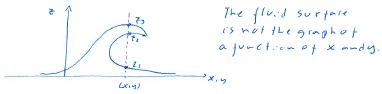}
\label{fig12.1}
\caption{}
\end{figure}

\noindent Inside the fluid volume we know that the velocity potential satisfies the Laplace equation. In order to pick out the unique solution describing ocean surface waves, we must impose boundary conditions at the bottom and the surface of the ocean. What is special here, is that the upper boundary surface is not fixed, but will move. This is thus an example of a \textit{moving boundary value problem}.

Let us start with the boundary conditions at the bottom of the ocean. The bottom is defined as the graph of a function $h(x,y)$.

\begin{equation*}
 z  = -h(x,y) \nonumber
\end{equation*} 
The basic assumption is that the bottom is impenetrable to water. This means that that the normal component of the fluid velocity field must be zero at the bottom.
\begin{align}
\vb{v} \vdot \vb{n} = 0 \; \; \text{at} \; \; z=-h(x,y). \label{eq:189}
\end{align}
Let
\begin{align*}
F(x,y,z) = z + h(x,y). 
\end{align*}
Then the bottom is a level surface for $F$ 
\begin{align*}
F(x,y,z) &= 0,\\
 &\Updownarrow  \\ 
 z &= - h(x,y). 
\end{align*}
A unit normal for the bottom is therefore given by
\begin{align*}
\vb{n} = \frac{\grad{F}}{\norm{\grad{F}}}.
\end{align*}
But
\begin{align*}
\grad{F} = (h_x,h_y,1), 
\end{align*}
and therefore the unit normal is
\begin{align*}
\vb{n} = \frac{(h_x, h_y, 1)}{(h_x^2 + h_y^2 + 1)^{\frac{1}{2}}},
\end{align*}
and the boundary condition at the bottom, (\ref{eq:189}), can be written as 
\begin{align*}
(\phi_x, \phi_y, \phi_z)\vdot\frac{(h_x,h_y,1)}{(h_x^2 + h_y^2 + 1)^{\frac{1}{2}}} &= 0, \\ 
&\Updownarrow \nonumber \\  \phi_x \; h_x + \phi_y \; h_y + \phi_z = 0 \; \; \text{at} \; \; z &= - h(x,y),
\end{align*}
where we have used the fact that for potential flows, the velocity field is determined by a potential $\phi$ through $\vb{v}=\grad{\phi}$.

At the surface we have two conditions, one \textit{kinematic} and one \textit{dynamic}. Let us first consider the kinematic condition. 

The surface is defined as an interface between air and water. As a consequence of this, water does not pass through the surface. This means that the fluid velocity at the surface must be equal to the velocity of the surface. This is the kinematic boundary condition at the surface.

Let
\begin{align}
\vb{x} = (x(t), y(t), z(t)), \label{eq:196}
\end{align}
be the position vector for a point on the surface. Since the point (\ref{eq:196}) is on the surface and the surface is the graph of the function $\eta(x,y)$, we must have
\begin{align}
z(t) = \eta(x(t) , y(t),t). \label{eq:197}
\end{align}
Differentiating (\ref{eq:197}) with respect to time we have 
\begin{align}
z' &= x'\; \eta_x + y'\; \eta_y + \eta_t.  \label{eq:198a}
\end{align}
We now use the kinematic boundary condition
\begin{align*}
\vb{x}'(t) &= \vb{v}(\vb{x}(t),t), \\
&\Updownarrow \nonumber \\  x'(t) &= v_x = \phi_x, \\ 
y'(t) &= v_y = \phi_y,  \\ 
z'(t) &= v_z = \phi_z, 
\end{align*}
where now $\phi_x\equiv\partial_x\phi$, etc. Thus, (\ref{eq:198a}) is transformed into the equation
\begin{align*}
  \eta_t + \eta_x \; \phi_x + \eta_y \; \phi_y &= \phi_z. 
\end{align*}
This is the final form of the kinematic boundary condition at the surface.

Let us next look at the dynamic boundary condition. 
We will assume that the surface has no mass and that there is no surface tension. Then, the net force acting on a small piece of the fluid surface is equal to $p-p_0$. The fact that the surface is mass-less implies then, through Newton's law, that 
\begin{align*}
p - p_0 = 0 \; \; at \; \; z= \eta (x,y,t).
\end{align*}
This is the dynamic boundary condition. Using the basic equation (\ref{eq:186}) for potential flow, our ocean surface wave problem is in summary:
\begin{align}
&\laplacian{\phi} = 0 && -h(x,y) < z<\eta(x,y,t), \label{eq:201}\\ 
&\phi_z + h_x \; \phi_x + h_y \; \phi_y = 0 && z = - h(x,y), \label{eq:202}\\  
&\eta_t + \eta_x \; \phi_x + \eta_y \; \phi_y = \phi_z &&  z = \eta (x, y, t),  \label{eq:203}\\ 
&\phi_t + \frac{1}{2} \; (\grad{\phi})^2 + y \; z = 0 &&  z = \eta (x, y, t). \label{eq:204}
\end{align}

\subsubsection{Low amplitude surface waves for irrotational flows}

A smooth undisturbed ocean is characterized by the solution 
\begin{align}
\eta (x,y,t) = \eta_0, \nonumber\\
\phi (x,y,t) = 0. \label{eq:205}
\end{align}
We will now linearize the system (\ref{eq:201}) - (\ref{eq:204}) around the simple solution (\ref{eq:205}). This will describe a situation where the ocean waves are of low amplitude. Introduce $\eta', \phi'$ by 
\begin{align*}
\eta = \eta_0 + \eta', \\
\phi = \phi'. 
\end{align*}
The Laplace equation is already linear so we have 
\begin{align*}
\laplacian{\phi'} = 0. 
\end{align*}
Note that, without loss of generality, we can assume that $\eta_0 = 0$ by choosing the origin of our coordinate system in an appropriate way. We assume that this has been done.

The boundary condition at the bottom is also linear and we get 
\begin{align*}
\phi'_z + h_x \; \phi'_x + h_y \; \phi'_y = 0 && z = - h(x,y).
\end{align*}
On the surface we get
\begin{equation}
\begin{rcases}
\eta'_t + \eta'_x \; \phi'_x + \eta'_y \; \phi'_y = \phi'_z \\
\phi'_t + \frac{1}{2} \; (\grad{\phi'})^2 + g\;z = 0 
\end{rcases} z = \eta'(x,y,t). \nonumber
\end{equation}
Linearizing by dropping products of small quantities we get
\begin{equation*}
\begin{rcases}
\eta'_t  = \phi'_z \\
\phi'_t + g \; \eta' = 0 
\end{rcases} z = 0. 
\end{equation*}
Note that the linearized boundary conditions are evaluated at $z=0$. This is because $\phi$ in general depends on $z$ in a nonlinear way, so that for example 
\begin{align*}
\phi_z'(x,y,\eta',t),
\end{align*}
will be nonlinear in $\eta$. To extract the linear part we Taylor expand
\begin{align*}
\phi'_z(x,y,\eta' , t) = \phi'_z(x,y,0,t) + \phi'_{zz}(x,y,0,t) \; \eta' + ... \;.
\end{align*}
The second term is a product of small quantities and can be dropped.

Summing up, the linearized surface wave problem is
\begin{align*}
\laplacian{\phi} = 0, && -h(x,y) < z <0, \\
\phi_z + h_x \; \phi_x + h_y \; \phi_y = 0, && z = -h(x,y)  \\
\eta_t = \phi_z, && z=0, \\ 
\phi_t + g \; \eta = 0, && z = 0, 
\end{align*} 
 where we are now dropping the primes from the variables.
This is a linear problem, but still hard to solve for a bottom of variable depth.

We will simplify the problem further by assuming that the bottom is perfectly flat
\begin{align*}
h(x,y) = h_0. 
\end{align*}
For this simplified problem we have
\begin{align}
\laplacian{\phi} = 0 && -h_0 < z <0, \label{eq:217}\\
\phi_z = 0 && z = -h_0, \label{eq:218}\\
\eta_t = \phi_z && z=0, \label{eq:219}\\ 
\phi_t + g \; \eta = 0 && z = 0. \label{eq:220}
\end{align}
We can solve this problem using Fourier modes of the form 
\begin{align}
\phi(x,y,z,t) = a(z) \; e^{i(\vb{k} \vdot \vb{x} - \omega \; t)}, \label{eq:221}
\end{align}
where $\vb{k} = (k_x, k_y),\; \vb{x} = (x,y)$. \\  \noindent Equation (\ref{eq:220}) determines $\eta$ in terms of $a(z)$
\begin{align*}
\eta(x,y,t) =  \Re {\frac{i \; \omega}{g} \; a(0) \; e^{i(\vb{k} \vdot \vb{x} - \omega \; t)} }. 
\end{align*}
From the Laplace equation (\ref{eq:217}), we get 
\begin{align*}
a''(z) - k^2\; a(z) = 0 && -h_0 < z <0, 
\end{align*}
where $k^2 = \vb{k} \vdot \vb{k}$. \\ 
The boundary conditions (\ref{eq:219}) and (\ref{eq:220}) implies that 
\begin{align*} g\phi_z+\phi_{tt} &=0,  \\ 
&\Updownarrow  \\ 
 g \; a'(z)-\omega^2 \; a(z)  &= 0 && z=0, 
\end{align*}
and the boundary condition (\ref{eq:218}) implies that 
\begin{align*}
a'(-h_0) = 0.
\end{align*} 
We thus get the following boundary value problem for the function $a(z)$ 
\begin{align}
a''(z) - k^2 \; a(z) = 0 && -h_0 < z <0, \label{eq:226}\\
g \; a'(0) - \omega^2 \; a(0) = 0, \label{eq:227}\\
a'(-h_0) = 0. \label{eq:228}
\end{align}
This problem we can easily solve. The general solution of (\ref{eq:226}) is 
\begin{align}
a(z) = A \; \cosh{(k \;z)} + B \; \sinh \: (k\;z). \label{eq:229}
\end{align}
The boundary conditions (\ref{eq:227}), (\ref{eq:228}) give 
\begin{align}
g\;a'(0) - \omega^2 \; a(0) = 0 \Leftrightarrow k\; g \; B - \omega^2 \; A = 0, \nonumber \\a'(-h_0) = 0 \Leftrightarrow -A \; \sinh \: (k \; h_0) + B \; \cosh \: (k\;h_0) = 0. \nonumber
\end{align}
We thus have the $2\times 2$ linear system 
\begin{align}
\mqty[-\omega^2 & k\;g \\ -\sinh \: (k\;h_0) & \cosh \: (k \; h_0)] \; \mqty[A \\ B] = 0.\label{eq:230}
\end{align}
A non-trivial solution exists only if the determinant of the matrix is zero 
\begin{align*}
-\omega^2 \; \cosh \: (k \; h_0) + k \; g \; \sinh \: (k \; h_0) = 0, 
\end{align*}
which can be written as 
\begin{align}
\omega^2 = g \; k \; \tanh \: (k \; h_0). \label{eq:232}
\end{align}
This is the \textit{dispersion relation} for small amplitude ocean surface waves. For a given $\vb{k}$ the surface elevation is 
\begin{align}
\eta(\vb{x} , t) = \Re{\frac{i \; \omega(k)}{g} \;a(0) \;  e^{i(\vb{k} \vdot \vb{x} - \omega(k) \; t)}}, \label{eq:233}
\end{align}
where 
\begin{align}
\omega(k) =\sqrt{g \; k \; \tanh \: (k\; h_0)}, \nonumber
\end{align}
is a solution of the dispersion relation (\ref{eq:232}).
The formula for the surface elevation (\ref{eq:233}) describes a two dimensional plane wave with phase speed,
\begin{align*}
v_f = \frac{\omega(k)}{k},
\end{align*}
 moving across the ocean surface.
\begin{figure}[!h]
\centering
\includegraphics{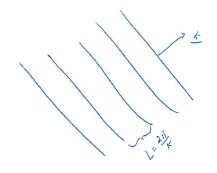}
\label{fig13}
\caption{}
\end{figure}
 This is illustrated in figure 13.
 
  Let us consider two limiting cases. Observe that the wavelength $L$, for a given plane wave with wavenumber $k$, is $L=\frac{2\pi}{k}$.

\noindent  The first case is defined by the inequality
 \begin{equation*}
 \frac{h_0 }{L}\gg 1.
 \end{equation*}
This means that the wavelength of the wave is much smaller than the depth of the ocean. This case is called \textit{waves in deep water}.
For this case we have
\begin{equation*}
k \: h_0 = \frac{2 \; \pi \; h_0}{L} \gg 1,
\end{equation*}
and therefore 
\begin{align*}
&\omega^2 = g \; k \; \tanh \: (k \; h_0) \approx g \; k, \\ 
&\Updownarrow  \\ 
&\omega \approx \sqrt{g \; k}. 
\end{align*}
Surface waves in deep water are therefore \textit{dispersive}; the phase velocity depends on $k$ 
\begin{align}
v_f = \frac{\omega}{k} \approx \sqrt{\frac{g}{k}}. \nonumber
\end{align}
Waves with long wavelengths move faster than waves of short wavelengths. Observe that the phase speed does not depend on the ocean depth for deep water waves. 

The second case is defined by the inequality
 \begin{equation*}
 \frac{h_0 }{L} \ll 1.
 \end{equation*}
This means that the wavelength of the surface waves are much larger than the ocean depth. This case is called \textit{waves in shallow water}. \\
For the dispersion relation we now get
\begin{align*}
\omega = \sqrt{ g \; k \; \tanh(k \; h_0)}, \\ 
\approx \sqrt{g \; k^2  \; h_0} = \sqrt{g \; h_0} \; k. 
\end{align*}
For this case the phase speed does \textit{not} depend on the wave number and waves in shallow water are \textit{non-dispersive}.  \\ 
Note that the phase speed 
\begin{align*}
v_f = \frac{\omega}{k} = \sqrt{g \; h_0}, 
\end{align*}
of shallow water waves depends on the depth of the ocean. The waves move faster in a deep ocean than in a shallow one. This explains why waves break over a reef.

\begin{figure}[htbp]
\centering
\includegraphics{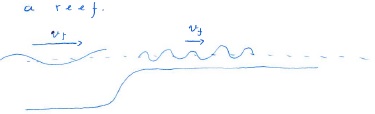}
\label{fig14}
\caption{}
\end{figure}
The waves "pile up" as they enter the shallow part of the reef. This piling up leads to larger wave amplitude and eventually wave breaking. Recall however, that wave breaking is not included in our modelling. \\
Shallow water waves differ from deep water waves in another important way that we are now going to explain. \\ 
Using (\ref{eq:229}) and (\ref{eq:230}) and assuming that the dispersion relation (\ref{eq:232}) holds, a solution of the boundary value problem (\ref{eq:226})-(\ref{eq:228}) is 
\begin{align*}
a(z) = \cosh \: (k \; z) + \tanh \: (k \; h_0) \; \sinh \: (k \; z). 
\end{align*}
Using the fact that
\begin{align*}
\cosh \: (x + y) = \cosh \: (x) \; \cosh \: (y) + \sinh \: (y) \; \sinh \: (y), 
\end{align*}
we can rewrite $a(z)$ into the form 
\begin{align*}
a(z) = \frac{\cosh \: (k \; (z + h_0))}{\cosh \: (k \; h_0)}, 
\end{align*}
and the corresponding velocity potential is, using (\ref{eq:221}), 
\begin{align}
\phi(\vb{x},z,t) = \frac{\cosh \: (k \; (z + h_0))}{\cosh \: (k \; h_0)} \; e^{i(\vb{k} \vdot \vb{x} - \omega \; t)}. \label{eq:241}
\end{align}
These are the modes for the linearized surface wave problem for the case when the ocean bottom is flat. \\ Let $\vb{u}$ be the horizontal and $w$ the vertical part of the fluid velocity field 
\begin{align*}
\vb{v} = ( \vb{u}, w).
\end{align*}
From (\ref{eq:241}) we have 
\begin{align*}
w = \partial_z \; \phi = \frac{k \; \sinh \: (k \; (z + h_0))}{\cosh \: (k \; h_0)} \; e^{i(\vb{k} \vdot \vb{x} - \omega \; t)}, \\
\vb{u} = ( \phi_x, \phi_y) = \vb{k} \; \frac{\cosh \: (k \; (z + h_0))}{\cosh \: (k \; h_0)} \; e^{i(\vb{k} \vdot \vb{x} - \omega \; t)}. 
\end{align*}
For shallow water waves we therefore have 
\begin{align*}
\frac{\abs{w}}{\norm{\vb{u}}} = \frac{\sinh \: (k \; (z + h_0))}{\cosh \: (k \; (z + h_0))} = \tanh \: (k \; (z + h_0)) << 1, 
\end{align*}
here we have used the fact that $\mid k(z+h_0)\mid\ <\;2kh_0\ll\;1$. Thus the vertical flow of the fluid is negligible compared to the horizontal one. The fluid flow is essentially \textit{horizontal}. \\
Also observe that 
\begin{align}
\frac{\norm{\partial_z  \vb{u}}}{\norm{\partial_x  \vb{u}}} = \frac{k}{k_x} \; \tanh(k \; (z + h_0)) << 1, \label{eq:246}
\end{align}
here we assume that the mode does not move along the y-axis because then $k_x = 0$ and the flow velocity of the mode in fact depends on y only. In a similar way we find 
\begin{align}
\frac{\norm{\partial_z  \vb{u}}}{\norm{\partial_y  \vb{u}}} << 1. \label{eq:247}
\end{align}
What (\ref{eq:246}) and (\ref{eq:247}) tell us, is that for shallow water waves, the flow velocity, in addition to being essentially horizontal, depends very weakly on the z-coordinate. This means that for shallow water waves the whole water column, from the bottom to the surface, will move. This is not true for deep water waves where there will be a significant fluid flow only in the surface layers. \\
The fact that the whole water column flows in shallow water waves, can pose an extreme danger to anyone living close to the shore. \\ 
Imagine a wave that satisfies the shallow water condition $\frac{h_0}{L} << 1$, in the deep ocean. Let us say that the depth is 5km. Such a wave would have to have a wavelength of around 5 -10km, say. For such a wave we have a water column of height 5km moving horizontally. The phase speed is approximately 
$$
v_f \approx \text{800 km/h !}
$$
It is thus clear that a wave of this type can carry an enormous amount of momentum and this momentum is transported at extreme speed. Eventually, such a wave will hit land and the momentum will be deposited in the shore area. Clearly this can be catastrophic. And it is has been, many times. This kind of wave is called a Tsunami when it gets close to the shore. For the 2004 Tsunami in the Indian ocean some stretches of coastline (close to the city of Banda Aceh) experienced a wave of a height close to 30 meters that hit the shore moving at approximately 50 km/h. The material damage was total and hundreds of thousands of lives were lost. \\
In order to initiate a shallow water wave in the deep ocean, something must happen that can set the whole water column in motion. Wind cannot do this but earthquakes can. 

\subsubsection{The shallow water equations}

Our aim is now to derive a simplified system of equations describing shallow water waves by using the two special properties of such waves. This derivation will not assume that the flow is vorticity free and the resulting equations are thus of wider generality than the equations from the previous section that was based on the general equations for potential flow. From our investigations in the previous section we have found that if we write 
\begin{align}
\vb{v} = ( \vb{u} , w), \nonumber
\end{align}
then 
\begin{align}
\vb{(i)} && &\frac{\abs{w}}{\norm{\vb{u}}} << 1,\nonumber \\ 
\vb{(ii)} && &\vb{u}(\vb{x},z,t) \approx \vb{u}(\vb{x},t)  \; \; \; \; \vb{x} = (x,y). \nonumber
\end{align}
We start the derivation from the Euler equations for an incompressible fluid
\begin{align*}
&\vb{v}_t + \vb{v} \vdot \grad{\vb{v}} = - \frac{\grad{p}}{\rho_0} + \vb{g}, \\
&\div{\vb{v}} = 0, 
\end{align*}
where, as usual, $\vb{g}$ is the constant force of gravity, and the coordinate system is chosen such that 
\begin{align*}
\vb{g} = - g \; e_z. 
\end{align*}
The vertical part of the Euler equation is 
\begin{align*}
 w_t + \vb{u} \vdot \grad_\perp{w} + w \; w_z = - \frac{p_z}{\rho_0} - g, 
\end{align*}
where $\nabla_{\perp}=(\partial_x,\partial_y)$ is the horizontal gradient operator. The left side of this equation is small by assumption so we must have
\begin{align*}
p_z \approx - \rho_0 \; g.
\end{align*}
Integrating this equation and applying the boundary condition $ p = p_0$ at $ z = \eta(x,y)$ we get 
\begin{align}
p = p_0 + \rho_0 \; g \; ( \eta - z). \label{eq:252}
\end{align}
The horizontal part of the Euler equation is 
\begin{align}
&\vb{u}_t + \vb{u} \vdot \grad_\perp{\vb{u}} + w \; \vb{u}_z = - \frac{\grad_\perp{p}}{\rho_0} = - g \; \grad_\perp{\eta} \nonumber \\
&\Updownarrow \nonumber \\
&\vb{u}_t+ \vb{u} \vdot \grad_\perp{u} + g \; \grad_\perp{\eta} = 0, \label{eq:253}
\end{align}
where we have used (\ref{eq:252}) and the fact that $w$ is small. I should really have used approximate equality in all expressions but it is customary to use equality in derivations like the one I am doing here. From the incompressibility condition we have 
\begin{align}
\div{\vb{v}} = \div_\perp{\vb{u}} + w_z = 0. \label{eq:254}
\end{align}
We now integrate (\ref{eq:254}) with respect to $z$ from the bottom $z = h(x,y)$ to the surface $z=\eta(x,y,t)$ 
\begin{align}
\int_{-h(x,y)}^{\eta(x,y,t)} dz \; \div_\perp{\vb{u}} = - (w \mathlarger{|}_{\eta} - w \mathlarger{|}_{-h}). \label{eq:255}
\end{align}
But we have the following boundary conditions 
\begin{align}
\eta_t + u \; \eta_x + v \; \eta_y &= w && z=\eta,\nonumber \\
 u \; h_x + v \; h_y + w &= 0 && z = -h, r\nonumber \\ 
&\Updownarrow  \nonumber\\ 
 w\mathlarger{|}_{z=\eta} &= \eta_t + \eta_x \; u\mathlarger{|}_{\eta} + \eta_y \; v\mathlarger{|}_{\eta}, \label{eq:257}\\
 w \mathlarger{|}_{z=-h} &= -h_x\; u\mathlarger{|}_{-h} - h_y \; v\mathlarger{|}_{-h}, \nonumber
\end{align}
where we have introduced the notation 
\begin{align*}
\vb{u} = ( u, v).
\end{align*}
Inserting (\ref{eq:257}) into (\ref{eq:255}) we get 
\begin{align}
\int_{-h(x,y)}^{\eta(x,y,t)} dz \; \div_\perp{\vb{u}} = -\eta_t - ( \grad_\perp{\eta} \vdot \vb{u} \mathlarger{|}_{\eta} + \grad_\perp{h} \vdot \vb{u} \mathlarger{|}_{-h}). \label{eq:259}
\end{align}
But we also have
\begin{align}
\partial_x \; ((\eta + h) \; u) &= \partial_x \; \int_{-h(x,y)}^{\eta(x,y,t)} u \; dz \nonumber\\
&= \eta_x \; u\mathlarger{|}_{\eta} - (-h_x) \; u\mathlarger{|}_{-h} + \int_{-h(x,y)}^{\eta(x,y,t)} dz \; \partial_x \; u, \label{eq:260}\\
\partial_y \; ((\eta + h) \; v) &= \partial_y \; \int_{-h(x,y)}^{\eta(x,y,t)} v \; dz \nonumber\\
&= \eta_y \; v\mathlarger{|}_{\eta} - (-h_y) \; v\mathlarger{|}_{-h} + \int_{-h(x,y)}^{\eta(x,y,t)} dz \; \partial_y \; v. \label{eq:261}
\end{align}
Adding (\ref{eq:260}) and (\ref{eq:261}) we get 
\begin{align}
\div_\perp{((\eta + h) \; \vb{u})} = \grad_\perp{\eta} \vdot \vb{u} \mathlarger{|}_{\eta} + \grad_\perp{h} \vdot \vb{u} \mathlarger{|}_{-h} + \int_{-h(x,y)}^{\eta(x,y,t)} dz \; \div_\perp{\vb{u}}, \label{eq:262}
\end{align}
and combining (\ref{eq:262}) and (\ref{eq:259}) we get 
\begin{align}
\div_\perp{((\eta + h) \; \vb{u})} = - \eta_t. \label{eq:263}
\end{align}
The system (\ref{eq:253}) and (\ref{eq:263}) are the \textit{shallow water equations}
\begin{align}
\vb{u}_t + \vb{u} \vdot \grad_\perp{\vb{u}} + g \; \grad_\perp{\eta} = 0, \label{eq:264}\\
\eta_t + \div_\perp{(\vb{u} \; ( \eta + h))} = 0 \nonumber. 
\end{align}
The shallow water equations is a closed system of 3 equations for the 3 unknowns $\eta$ and $\vb{u}$, and are very important in the science of surface waves, in particular for the modeling of Tsunami waves. They however have a flaw that I will discuss now. \\
Let us for simplicity assume that the bottom is horizontal, so that $h=h_0$ and that the flow is 1D, $ u = u(x,t)$, $v=0$. Then the system (\ref{eq:264}) can be written as 
\begin{align}
\partial_t \; \mqty(u \\ \eta) + \mqty( u && g \\ h_0 + \eta && u) \; \partial_x \; \mqty(u \\ \eta) = 0. \label{eq:265}
\end{align}
Define
\begin{align}
\vb{U} = \mqty(u \\ \eta), && A(\vb{U}) = \mqty(u && g \\ \eta + h_0 && u). \label{eq:266}
\end{align}
Then (\ref{eq:265}) is of the general form 
\begin{align}
\partial_t \; \vb{U} + A(U) \; \partial_x \; \vb{U} = 0 \label{eq:267}
\end{align}
Let us look for a solution to (\ref{eq:267}) of the form 
\begin{align}
\vb{U}(x,t) = \vb{S}( \phi(x,t)), \label{eq:268}
\end{align}
for some scalar function $\phi$. Inserting (\ref{eq:268}) into (\ref{eq:267}) we get 
\begin{align}
\frac{d\vb{S}}{d\phi} \; \partial_t \; \phi + A(\vb{S}) \; \frac{d\vb{S}}{d \phi} \; \partial_x  \; \phi = 0. \label{eq:269}
\end{align}
Let us assume that $ \frac{d \vb{S}}{d \phi}$ is an eigenvector of $A( \vb{S} (\phi))$ with corresponding eigenvalue $ \lambda(\phi)$. Then from (\ref{eq:269}) we get 
\begin{align*}
 ( \partial_t \; \phi + \lambda (\phi) \; \partial_x \; \phi ) \; \frac{d \vb{S}}{d \phi} = 0. 
\end{align*}
Thus (\ref{eq:268}) is a solution to (\ref{eq:267}) if 
\begin{align}
\vb{R}(\phi)\equiv\frac{d \vb{S}}{d \phi}, \label{eq:271}
\end{align}
is an eigenvector of $A(\vb{S} ( \phi))$ with eigenvalue $\lambda (\phi)$ and where the following equation must hold
\begin{align*}
\partial_t \; \phi + \lambda (\phi) \; \partial_x \; \phi = 0.
\end{align*}
These kinds of solutions can be found for many nonlinear systems of equations and are called \textit{simple solutions}. They are a kind of \textit{nonlinear mode}. For the shallow water equations where $A(\vb{S})$ are defined in (\ref{eq:266}) we find two eigenvalues and eigenvectors 
\begin{align*}
\lambda_1 = u - \sqrt{g \; h} && \vb{R}_1 = \mqty(1 \\ - \frac{1}{\sqrt{\frac{h}{g}}}), \\
\lambda_2 = u + \sqrt{g \; h} && \vb{R}_2 = \mqty(1 \\ \frac{1}{\sqrt{\frac{h}{g}}}), 
\end{align*}
where $ h = \eta + h_0$. \\ 
For the first choice we must solve the equation 
\begin{align*}
\partial_t \; \phi + \lambda_1 (\phi) \; \partial_x \; \phi,
\end{align*}
where 
\begin{align}
\lambda_1(\phi) = u (\phi) - \sqrt{g \; h(\phi)} = u(\phi) - \sqrt{g \; (\eta (\phi) + h_0)}, \label{eq:276}
\end{align}
and where 
\begin{align*}
\frac{d\vb{S}}{d\phi}=\frac{d}{d \phi} \; \mqty( u \\ \eta) =\vb{R}_1 = \mqty(1 \\ - \frac{1}{\sqrt{\frac{\eta + h_0}{g}}}). 
\end{align*}
Thus 
\begin{align*}
\frac{du}{d \phi} = 1, \\
\frac{d \eta}{d \phi} = - \sqrt{ \frac{\eta + h_0}{g}}.
\end{align*}
A solution of the first equation is 
\begin{align}
u = \phi, \label{eq:279}
\end{align}
and an implicit solution of the second equation is 
\begin{align}
2 \; \sqrt{\eta + h_0} = - \frac{\phi}{\sqrt{g}}.  \label{eq:280}
\end{align}
Inserting (\ref{eq:279}) and (\ref{eq:280}) into (\ref{eq:276}) we get 
\begin{align*}
\lambda_1 (\phi) = \phi - \sqrt{g} \; ( - \frac{\phi}{2\; \sqrt{g}}) = \frac{3}{2} \; \phi. 
\end{align*}
Thus the equation for $\phi = \phi (x,t)$ is 
\begin{align}
\partial_t \; \phi + \frac{3}{2} \; \phi \; \partial_x \; \phi = 0. \label{eq:282}
\end{align}
But we know that solutions to (\ref{eq:282}) break down and form vertical waves at a finite break-down time. This strongly indicates that solutions of the shallow water equations (\ref{eq:264}) will tend to form singularities in the form of breaking waves. When this happens the assumptions underlying the shallow water equations also break down and the equations are no longer valid.

The equations can be repaired by doing a more careful and less heuristic derivation for shallow water waves. This will add higher order spatial derivatives to (\ref{eq:264}), which will remove the breaking waves, and thus regularize the equations. The approach used to regularize the shallow water equations is part of a large domain of applied mathematics called {\it perturbation methods}. We will give an introduction to some part of this domain in section five of these lecture notes.

\subsection{Computational project}

In this project we will simulate surface waves in a narrow channel of finite
length. The waves will be generated by a time dependent deformation of the
bottom of the channel

\bigskip

\bigskip

\begin{description}
\item[a)] Show that surface waves in a narrow channel is modelled by the
following system of partial differential equations%
\begin{align*}
\varphi_{xx}+\varphi_{zz}  & =0,\text{ \ \ \ \ \thinspace}-h(x,t)<z<\eta
(x,t)\\
\varphi_{z}+h_{x}\varphi_{x}+h_{t}  & =0,\text{ \ \ \ \ \ \ }z=-h(x,t)\\
\varphi_{z}-\eta_{x}\varphi_{x}-\eta_{t}  & =0,\text{ \ \ \ \ \ \ }%
z=\eta(x,t)\\
\varphi_{t}+\frac{1}{2}\nabla\varphi\cdot\nabla\varphi+gz  & =0,\text{
\ \ \ \ \ \ }z=\eta(x,t)\\
\varphi_{x}(-L,z,t)  & =0,\\
\varphi_{x}(L,z,t)  & =0\text{ \ \ \ \ \ \ }%
\end{align*}

\item Narrow channel here means that we disregard all dependence on the y
coordinate. I want you to present a detailed derivation of the equations for
this case. Imagine you are writing the presentation for a person that is not familiar with the theory of surface waves.

\item[b)] Linearize the system by assuming that $\varphi$ and $\eta$ are small
and that the time dependent bottom is given by
\[
h(x,t)=h_{0}+\xi(x,t)
\]
where $\xi(x,t)$ is small compared to the constant $h_{0}$. \ Elliminate
the surface elevation $\eta(x,t)$ from the system by differentiation. You now
have a linear system for the unknown function $\varphi(x,z,t)$.

\item[c)] Show that eigenvalues and corresponding eigenfunctions defined by
\begin{align*}
\chi^{^{\prime\prime}}(x)  & =\lambda\chi(x)\\
\chi^{\prime}(-L)  & =\chi^{\prime}(L)=0
\end{align*}
are
\[
\left.
\begin{array}
[c]{c}%
\chi_{n}(x)=\cos(k_{n}(x+L))\\
\lambda_{n}=-k_{n}^{2}%
\end{array}
\right\}  ,\text{ \ \ \ \ \ }n=0,1,....
\]
where $k_{n}=\frac{n\pi}{2L}$.

\item[d)] Introduce expansions%
\begin{align*}
\varphi(x,z,t)  & =%
{\displaystyle\sum_{n=0}^{\infty}}
\varphi_{n}(z,t)\chi_{n}(x)\\
\xi(x,t)  & =%
{\displaystyle\sum_{n=0}^{\infty}}
\xi_{n}(t)\chi_{n}(x)
\end{align*}
 in the linearized equations from b) and derive equations for the coefficient
functions $\varphi_{n}(z,t)$.

\item[e)] Solve the equations from d) and show that
\begin{align*}
\varphi_{0}(z,t)  & =-\xi_{0}^{\prime}(t)z+B_{0}(t)\\
\varphi_{n}(z,t)  & =A_{n}(t)e^{k_{n}z}+B_{n}(t)e^{-k_{n}z},\text{
\ \ \thinspace}n=1,2,...
\end{align*}
where
\begin{align*}
B_{0}^{\prime}(t)  & =g(\xi_{0}(t)+C)\\
B_{n}(t)  & =A_{n}(t)e^{-2k_{n}h_{0}}+\frac{e^{-k_{n}h_{0}}}{k_{n}}\xi
_{n}^{\prime}(t)\\
A_{n}^{^{\prime\prime}}(t)+\omega_{n}^{2}A_{n}(t)  & =f_{n}(t)
\end{align*}
and where $f_{n}(t)$ is a certain function determined by $\xi_{n}(t)$
and where
\[
\omega_{n}^{2}=\kappa_{n}\tanh(k_{n}h_{0})
\]
is the dispersion relation.

\item[f)] The initial conditions for the system are
\begin{align*}
\eta(x,0)  & =0\\
\varphi(x,z,0)  & =0
\end{align*}
Show that the initial conditions are satisfied if
\begin{align*}
A_{n}(0)  & =A_{n}^{\prime}(0)=0\\
C  & =-\xi_{0}(0)\\
\xi_{n}^{\prime}(0)  & =\xi_{n}^{\prime\prime}(0)=0
\end{align*}
The last condition is a constraint on how we can let the bottom deform.

\item[g)] We will consider bottom deformations of the general form\
\[
\xi(x,t)=f(t)g(x)
\]
Show that the condition $\xi_{n}^{\prime\prime}(0)=0$ from f) is satisfied if
$f^{\prime\prime}(0)=0$. Solve the equations for $A_{n}(t)$ numerically and
plot the time evolution of the surface elevation and the fluid velocity field
for at least two different shapes, $g(x)$ and some choise for $f(t)$. Let one
of the shapes be a Gaussian and at least one of the other some nonsymmetric
shape, for example a wedge. Use the plots to argue that the generated waves
are shallow water waves.
\end{description}

\setcounter{equation}{0}

\section{Calculus Of Variation}

\subsection{Generalized extremal problems}
Extremal problems, like minimum and maximum problems, have played a major role in the development of calculus. In fact, calculus was more or less invented to solve such problems. 

\noindent In the language of calculus, the quantity we need to maximize or minimize is a \ttx{function} of a real variable, $x$, and the challenge is to find an $x_0$ such that 
\begin{align}
f(x) \leq f(x_0) \; \; (f(x) \geq f(x_0)) && \forall x \ne x_0, \lbl{1} 
\end{align}
If $f(x)$ is a well behaved function, calculus tells us that we only need to look at points $x^*$ such that
\begin{align}
f'(x^*) = 0. \lbl{2}
\end{align}
All maximum and minimum points will be found among the set of points that satisfy \rf{2}. By well behaved functions, we mean here that $f(x)$ is continuously differentiable and defined on the whole real axis. 

There are however many important extremal problems that do not fall into the category described above. In fact, some of these problems are much older than calculus itself. 

The so-called \ttx{isoperimetric problem} was clearly stated already 200 BC by the Greek mathematician Zenodorus. The problem consists of finding, among all curves of a fixed length $L$, the curve that encloses the largest area. 

We will start our study of the calculus of variation by introducing several other interesting extremal problems that are beyond the bounds of standard calculus, and then move on to describing mathematical tools that we can use to solve them.

\subsubsection{Curve of shortest length in a plane}\label{VariationalExample1}

Let $p$ and $q$ be two points on the plane. The challenge is to find a curve, $C$, connecting $p$ and $q$, and which is of shortest possible length. 
\begin{figure}[htbp]
\centering
\includegraphics{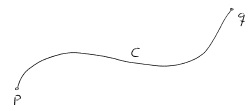}
\caption{Curve connecting two points $p$ and $q$  in the plane}
\label{fig1}
\end{figure}
\\\\
In order to state the problem in precise mathematical terms we, introduce a parametrization $\vb{\gamma}(t)$ for $C$. 
\\
\begin{align*}
\vb{\gamma}(t) &= (x(t) , y(t)), && 0 \leq\;t\leq 1,  \\ 
\vb{\gamma} ([0,1]) &= C,  \\
\vb{\gamma}(0) &= p, && \vb{\gamma} (1) = q.  
\end{align*}
Using this parametrization, the length of the curve, $L(C)$, can be written as 
\begin{align*}
L(C) = \int^1_0 dt \; \norm{\vb{\gamma}'(t)} = \int^1_0 dt \; \sqrt{x'(t)^2 + y'(t)^2}. 
\end{align*}
The challenge is then to find a curve, $C_0$, such that
\begin{align*}
L(C) \geq L(C_0),
\end{align*}
for all $C$ connecting $p$ and $q$. This looks just like the minimum problem \rf{1} from elementary calculus. 

The only difference is that $L$ is \ttx{not} a function of a real variable, but is rather a function defined on the set of smooth curved connecting the points $p$ and $q$. 
\begin{figure}[htbp]
\centering
\includegraphics{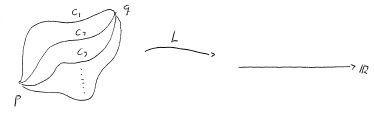}
\caption{A functional defined on plane curves}
\label{fig2}
\end{figure}

\noindent Such a function is called a \ttx{functional}. We will in this section of the  lecture notes encounter many other functions of this type, with ever larger domains of definition. All such functions will be called functionals.

\subsubsection{Curve of shortest length on a surface}\label{VariationalExample2}

Let $S$ be a surface in $\mathbf{R}^3$ and let $p$ and $q$ be points on $S$. 
The challenge is to find a curve, $C$, on the surface $S$, connecting $p$ and $q$, and that is of the shortest possible length. 
\begin{figure}[htbp]
\centering
\includegraphics{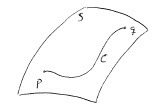}
\caption{A curve on the surface, $S$ connecting two points $p$ and $q$ on the surface}
\label{fig3}
\end{figure}

\noindent Let $\vb{\gamma}(t)$ be a parametrization for $C$ 
\begin{align*}
\vb{\gamma}(t) &= (x(t), y(t), z(t)) && 0 \leq t \leq 1,  \\
\vg ([0,1]) &= C \subset S,  \\
\vg(0) &= p, && \vg(1) = q. 
\end{align*}
The length of $C$ is
\begin{align*}
L(C) = \int^1_0 dt \; \sqrt{x'(t)^2 + y'(t)^2 + z'(t)^2}, 
\end{align*}
and the challenge is to find a curve, $C_0$, such that 
\begin{align*}
L(C) \geq L(C_0) \; \; \; \; \; \text{for all $C \subset S$ connecting $p$ and $q$.}
\end{align*}
This looks exactly like an extremal problem from regular calculus where we have a \ttx{constraint}. In regular calculus such problems are solved using Lagrange multipliers. We will see that a similar approach will work for functionals. 

For this particular problem we can remove the constraint using some vector calculus.

Let the surface be parametrized by $\vb{x}(u,v)$,
\begin{figure}[htbp]
\centering
\includegraphics{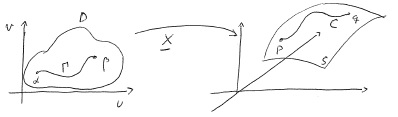}
\caption{A parametrization for the surface S}
\label{fig4}
\end{figure}

\noindent where 
\begin{align*}
\vb{x}(u,v) &= (x(u,v), \; y(u,v), \; z(u,v)), &&  (u,v) \in D,  \\
\vb{x} (D) &= S,  \\
\vb{x} (\Gamma) &= C \;, \; \; \; \vb{x(\alpha)} = p, \;  \; \; \; \vb{x} (\beta) = q. 
\end{align*}
Let the curve $\Gamma$ be parametrized by $\vb{\gamma}(t)$ 
\begin{align*}
\vg (t) &= (u(t), v(t) ), && 0 \leq t \leq 1,  \\ 
\vg([0,1]) &= \Gamma,  \\ 
\vg(0) &= \alpha, && \vg(1) = \beta. 
\end{align*}
Then then the curve $C$ is parametrized by 
\begin{align*}
\vb{\xi} (t) = \vb{x}(\vg(t)), && 0 \leq t \leq 1.
\end{align*}
In component form we have 

\begin{align*}
x(t) &= x(u(t), \; v(t)),  \\ 
y(t) &= y(u(t), \; v(t)),  \\
z(t) &= z(u(t), \; v(t)),
\end{align*}
so that 
\begin{align*}
\vb{\xi}(t) & = (x(t), \; y(t), \; z(t)).
\end{align*}
 Using the chain rule we have 
\begin{align*}
x'(t) &= \prt{u}x(u(t), v(t)) u'(t) + \prt{v} x(u(t), v(t)) v'(t),  \\
y'(t) &= \prt{u} y(u(t),  v(t))  u'(t) + \prt{v} y(u(t),  v(t)) v'(t),  \\ 
z'(t) &= \prt{u} z(u(t),  v(t))  u'(t) + \prt{v} z(u(t),  v(t))  v'(t). 
\end{align*}
In vector form this can be written as 
\begin{align*}
\vb{\xi}'(t) = \vb{T}_u \; u' + \vb{T}_v \; v',
\end{align*}
where the vectors $\vb{T}_u$ and $\vb{T}_v$ are defined by 
\begin{align*}
\vb{T}_u &= \vb{T}_u(t) =\vb{T}_u(u(t), v(t))  \\ 
&= (\prt{u}x(u(t),  v(t)) , \prt{u} y(u(t),  v(t)) ,  \prt{u} z(u(t),  v(t))),  \\
\vb{T}_v &= \vb{T}_v(t) = \vb{T}_v(u(t), v(t))  \\ 
&= (\prt{v}x(u(t), v(t)) , \prt{v} y(u(t), v(t)) ,  \prt{v} z(u(t), v(t))). 
\end{align*}
The two vectors $\vb{T}_u(u,v)$. $\vb{T}_v(u,v)$ are tangent to the surface $S$. Using $\vb{T}_u$ and $\vb{T}_v$ we have 
\begin{align}
x'(t)^2 + y'(t)^2 + z'(t)^2 &= \norm{\vb{\xi}'(t)}^2 \nonumber \\ 
&= \vb{\xi}'(t) \vdot \vb{\xi}'(t) = (\vb{T}_u \; u' + \vb{T}_v \; v' ) \vdot (\vb{T}_u \; u' + \vb{T}_v \; v') \nonumber \\
&= A \; u'^2 + 2 \; B \; u' \; v' + C \; v'^2, \lbl{20.1}
\end{align}
where
\begin{align}
A &= A(t) = \vb{T}_u \vdot \vb{T}_u, && B = B(t) = \vb{T}_u \vdot \vb{T}_v,\nonumber\\ 
C &= \vb{T}_v \vdot \vb{T}_v.  \lbl{21.1} 
\end{align}
Observe that the functions $A, \; B, \; C$ depends only on the structure of the surface $S$ 
\begin{align*}
A(u,  v) &= \vb{T}_u(u, v) \vdot \vb{T}_u (u, v), \\
B(u, v) &= \vb{T}_u(u, v) \vdot \vb{T}_v (u, v),  \\
C(u, v) &= \vb{T}_v(u, v) \vdot \vb{T}_v (u, v). 
\end{align*}
We are here abusing the notation in the usual calculus way.
\begin{align*}
A = A(t) = A ( u(t), \; v(t)).  
\end{align*}
Using \rf{20.1} and \rf{21.1}, our minimum problem is now to find a curve $\Gamma_0 \subset D$, connecting $\alpha$ and $\beta$, such that 
\begin{align*}
L(\Gamma) \geq L(\Gamma_0) \;\text{for all curves in $D$ connecting $\alpha$ and $\beta$}, 
\end{align*}
where the functional $L$ is 
\begin{align*}
L(\Gamma) = \int^1_0 dt \; [A(t) \; u'(t)^2 + 2 \; B(t) \; u'(t) \; v'(t) + C(t) \; v'(t)^2]^{\inv{2}}.  
\end{align*}
This is now an unconstrained minimum problem. \\
Observe that the three functions, $A, \; B$ and $C$ determine the length of \ttx{all} curves on the surface $S$. Thus the three functions determine the \ttx{geometry} of the surface $S$. 

The functions $A, \; B$ and $C$ are called the \ttx{metric} coefficients for $S$ and the curve of minimal length, $C_0$, connecting $p$ and $q$ is called a \ttx{geodesics} for the surface. 

\noindent These ideas and their generalization play a fundamental role in theoretical physics. 

In fact, Einstein's fundamental contribution to gravitational physics was to merge gravitational and inertial forces into a collection of 10 metric coefficients determining the geometry of the four dimensional space-time continuum. From the metric coefficients we can calculate the length of curves connecting points, here called \ttx{events}, in space-time. The length of a curve, in this space-time context, is the time it would take an observer,  to follow the curve from an event $p$ to another event $q$. The time in question is the one measured by the observer following the curve. This time is called \ttx{proper} time in Einstein's gravitational theory (General Relativity). 

The physical postulate is that an observer following the curve $C_0$ from $p$ to event $q$, will feel no inertial forces if the curve has maximal length. 
Since length is proper time, this means that if an observer wants to move from event $p$ to another $q$, without feeling inertial forces, she should choose a curve that takes as much time as possible as measured by her clock. 

These curves are also called geodiscs, but here, in the space-time context, they are curves of maximal, not minimal, length. 
The metric coefficients are determined by the distribution of mass and energy through the \ttx{Einstein field equations}.

\subsubsection{The isoperimetric problem} \label{VariationalExample3}
Let us consider the truly ancient isoperimetric problem.
Let $C$ be a curve enclosing a domain $D$.
\begin{figure}[htbp]
\centering
\includegraphics{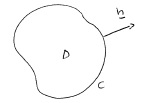}
\caption{Curve enclosing a domain $D$}
\label{fig5}
\end{figure}

\noindent  We want to express the area of $D$ in terms of $C$, and for this purpose introduce a vector field $\vb{f}$, on the plane given by 
\begin{align*}
\vb{f}(x,y) = \inv{2} (x,y).
\end{align*}
Evidently we have 
\begin{align*}
\div{\vb{f}} = \prt{x}(\inv{2} \;x) + \prt{y} (\inv{2} \; y) = \inv{2} + \inv{2} = 1. 
\end{align*}
But then, using the divergence form of Green's theorem, we have
\begin{align*}
A = \iint\limits_{D} dA = \iint\limits_D dA \; \div{\vb{f}} = \oint\limits_{C} dl \; \vb{f} \vdot \vb{n}. 
\end{align*}
Thus the area of $D$ is a functional of $C$ 
\begin{align}
A = A(C) = \oint\limits_{C} dl \; \vb{f} \vdot \vb{n}. \lbl{30.1} 
\end{align}
The length of the curve $C$ is determined by the functional 
\begin{align}
L = L(C) = \oint\limits_{C} dl .\lbl{31.1} 
\end{align}
The isoperimetric problem consists of, for a fixed length $L^*$, finding the curve $C_0$ such that 
\begin{align*}
A(C) \leq A(C_0) && \forall C \; , \; \; L(C) = L^*.
\end{align*}
We will now parametrize this problem, and therefore introduce a counter clockwise orientation for the curve $C$. 
\begin{figure}[htbp]
\centering
\includegraphics{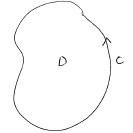}
\caption{Defining orientation for the curve $C$}
\label{fig6}
\end{figure}

\noindent Let $\vb{\gamma}$ be a parametrization of the orientated curve $C$. We thus have 
\begin{align}
\vg(t) &= (x(t) , \; y(t)) && 0 \leq t \leq 1, \lbl{34.1} \\ 
\vg([0,1]) &= C, \nonumber \\
\vg(0) &= \vg(1), \nonumber 
\end{align}
and $\vg'(t)$ points along the orientation of $C$. The choice of orientation ensures that 
\begin{align}
\vb{n} (t) = \frac{(y'(t), - x'(t))}{\norm{\vg'(t)}}, \lbl{35.1} 
\end{align}
is a unit normal defined on $C$ that points $\it{out}$ of the domain $D$. Using \rf{34.1} and \rf{35.1} in \rf{30.1} and \rf{31.1} we have
\begin{align*}
L(C) &= \oint\limits_{C} dl = \int\limits_0^1 dt \; \norm{\vg'(t)} = \int\limits^1_0 dt \; \sqrt{x'(t)^2 + y'(t)^2},  \\ 
A(C) &= \oint\limits_C dl \; \vb{f} \vdot \vb{n} = \int\limits_0^1 dt \; \norm{\vg'(t)} \; \inv{2} \; (x(t), \; y(t)) \cdot \frac{(y'(t), \; -x'(t))}{\norm{\vg'(t)}}  \\ 
&= \inv{2} \; \int^1_0 dt \; (x(t) \; y'(t) - y(t) \; x'(t)). 
\end{align*}
Thus our problem consists in finding functions $x_0(t),  y_0(t)$ such that 
\begin{align*}
&\mathbf{i)} \; \; x_0(0) = x_0(1),\;\; \; y_0(0) = y(1),\\ 
&\mathbf{ii)} \; \; \int^1_0 dt \; \sqrt{x'(t)^2 + y'(t)^2} = L^*,  \\
&\mathbf{iii)} \; \; \inv{2} \; \int^1_0 dt \; (x(t) \; y'(t) - y(t) \; x'(t)) \;\;\; \text{is maximal}. 
\end{align*}

\subsubsection{Surface of revolution of minimal area}\label{VariationalExample4}
Let $y(t)$ be a function defined on the interval $(x_1, x_2)$. Assume $y(x) > 0$ for all $x \in (x_1, x_2)$.
\begin{figure}[htbp]
\centering
\includegraphics{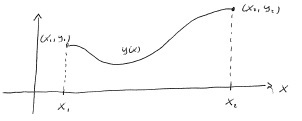}
\caption{Curve defining surface of revolution}
\label{fig7}
\end{figure}

\noindent Let $A(y)$ be the area of the surface of revolution that we get by rotating the curve $y(x)$ around the $x$-axis. The challenge is to find a curve $y(x)$ such that $A(y)$ is minimal for given fixed points $(x_1, y_1)$ and $(x_2, y_2)$. 

Recall from elementary Calculus that the formula for the area $A(y)$ is 
\begin{align*}
A(y) = 2 \; \pi \; \int^{x_2}_{x_1} dx \; y(x) \; \sqrt{1 + y'(x)^2}.  
\end{align*}
Thus the challenge is to find a curve $y(x)$, defined on $(x_1, x_2)$, such that 
\begin{align*}
&\mathbf{i)} \; \;  y(x_1) = y_1,\;\, y(x_2) = y_2,  \\ 
&\mathbf{ii)} \; \;  2 \; \pi \; \int_{x_1}^{x_2} dx \; y(x) \; \sqrt{1 + y'(x)^2} \; \; \; \text{is minimal}.
\end{align*}
\\\\
\subsubsection{General surface of minimal area}\label{VariationalExample5}
Let a curve $C$, in $\text{R}^3$ be given. The challenge is to find a surface $S \subset \textbf{R}^3$ such that 
\begin{align*}
&\mathbf{i)} \; \; \partial S = C,  \\
&\mathbf{ii)} \; A(S ) \; \text{is minimal}, 
\end{align*}
where $A(S)$ is the area of the surface $S$. Such a surface is called a \textit{minimal surface}. Many important problems in theoretical physics and applied mathematics can be reduced to the problem of finding a minimal surface. 

\noindent For example, if we dip a piece of string, described by a closed curve in $\textbf{R}^3$, into a bucket of soap water, the resulting soap film, clinging to the string, will form a minimal surface. 

\noindent Let us parametrize this problem  by  introducing a parametrization $\vb{x}$ of the surface $S$. Thus  
\begin{figure}[htbp]
\centering
\includegraphics{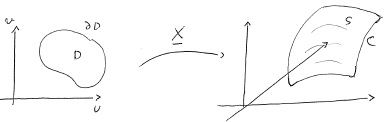}
\caption{Parametrization of surface $S$ bounded by the curve $C$}
\label{fig8}
\end{figure}

\noindent
\begin{align*}
&\mathbf{i)} \; \; \vb{x} : \mathbf{R}^2 \rightarrow \mathbf{R}^3, && \vx = \vxuv,  \\ 
&\mathbf{ii)} \; \; \vx(D) = S,  \\ 
&\mathbf{iii)} \; \; \vx(\prt{} D) = C.  
\end{align*}
From calculus we know that the area of $S$, $A(S)$, is given by the formula 
\begin{align*}
A(S) = \iint\limits_D du \; dv \; \norm{\vTu \cp \vTv},  
\end{align*}
where as before 
\begin{align*}
\vTu = \prtdvxu,\;\;\; \vTv = \prtdvxv.  
\end{align*}
The challenge is to choose functions $x(u,v), \; y(u,v), \; z(u,v)$ with
\begin{align*}
\vb{x}(u,v) = (x(u,v),  y(u,v),  z(u,v)), 
\end{align*}
such that $A(S)$ is minimal under the constraint 
\begin{align*}
\vx(\prt{} D) = C. 
\end{align*}
Let us restrict to the case when $C$ is given by the graph of a function $h=h(x,y)$ defined on the boundary, $\prt{} D$, of $D$. A surface of this type is displayed in figure \ref{fig9.1}.
\begin{figure}[htbp]
\centering
\includegraphics{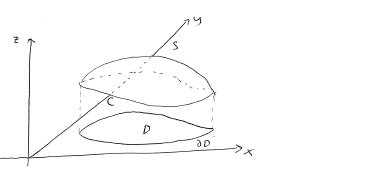}
\caption{Surface whose boundary is the graph of a function.\label{fig9.1}}
\end{figure}
\noindent For such surfaces we can use a parametrization of the form 
\begin{align*}
\vx(x,y) = (x,y,f(x,y)).  
\end{align*}
For this parametrization we have 
\begin{align*}
\vT_x = \frac{\prt{}\vx}{\prt{} x} = (1 , \; 0, f_x),  \\
\vT_y = \frac{\prt{}\vx}{\prt{} y} = (0, \; 1, f_y).  
\end{align*}
Thus 
\begin{align*}
\norm{\vT_x \cp \vT_y} = \sqrt{1 + f_x^2 + f_y^2}.  
\end{align*}
Therefore the challenge is to find a function $f(x,y)$, defined on a domain $D$ in $\mathbf{R}^2$, such that 
\begin{align*}
&\mathbf{i)} \; f(x,y) = h(x,y), \;\;\;  (x,y) \in \prt{} D,  \\ 
&\mathbf{ii)} \; \iint\limits_D dx \; dy \; \sqrt{1 + f_x^2 + f_y^2} \; \; \; \; \; \text{is minimal}. \nonumber 
\end{align*}
Note that in this case the functional is defined on a domain consisting of functions of two variables. This is clearly a very large domain.

\subsubsection{The Fermat Principle}\label{VariationalExample6}

Let $c$ be the speed of light in vacuum. Recall that the speed of light typically depends on the physical properties of the medium it is travelling through. \\
The ratio between the speed of light in vacuum and the speed of light in a medium, $v$, is a dimensionless number 
\begin{align}
n = \frac{c}{v}, \lbl{54.1} 
\end{align}
which is called the refractive index of the material. Under normal circumstances $v<c$ so that $n>1$. \\ 
The refractive index typically depends on the frequency of the light, but we will disregard this effect here. Unless the material is homogeneous, the refractive index will depend on position
\begin{align*}
n = n ( \vx).  
\end{align*}
Let us now consider a light ray passing through a medium that has a refractive index $n(\vx)$. 
\begin{figure}[htbp]
\centering
\includegraphics{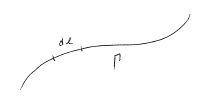}
\caption{Light ray in a refractive medium.\label{fig10.1}}
\end{figure}

\noindent Using \rf{54.1} and figure \ref{fig10.1}, we have 
\begin{align*}
n &= \frac{c}{\frac{dl}{dt}}, \; \;  \\ &\Updownarrow  \\ \; \; dt &= \frac{n}{c} \; dl,  
\end{align*}
where now $dt$ is the time it takes light to propagate the distance $dl$ along the curve $\Gamma$. The total time it takes light to propagate along a curve $\Gamma$ is then 
\begin{align*}
T(\Gamma) = \frac{1}{c} \; \int_{\Gamma} dl \; n.  
\end{align*}
Fermat's principle say that light follows the path through a medium of refractive index $n(\vx)$, that takes the shortest time. \\
Thus in order to find the path followed by light we must minimize $T(\Gamma)$ over all paths $\Gamma$. \\
Let us parametrize this problem. Let $\vg$ be a parametrization for $\Gamma$. 
\begin{align*}
\vg(t) &= (x(t), y(t),z(t)), \;\;\;0\leq t \leq 1,  \\ 
\vg(0) &= p, \;\;\; \vg(1) = q.  
\end{align*}
Thus, in order to find the path followed by light from $p$ to $q$ we must find functions $x(t), \; y(t), \; z(t) $ such that 
\begin{align*}
T(\Gamma) = \int^1_0 dt \; n(x(t), y(t), z(t)) \; \sqrt{x'(t)^2 + y'(t)^2 + z'(t)^2},  
\end{align*}
subject to the constraints 
\begin{align*}
(x(0), y(0), z(0)) &= p,  \\
(x(1), y(1), z(1)) &= q,  
\end{align*}
is as small as possible.

\subsubsection{The brachistochrone problem.}\label{VariationalExample7}

(brachistochrone - shortest time in Greek)

The challenge is to find the arc $(x, y(x))$ that a particle of mass $m$ must follow from $(x_1,y_1)$ to $(x_2,y_2)$ in order to use as little time as possible. The particle is influenced by a constant gravitational field pointing vertically down. 
\begin{figure}[htbp]
\centering
\includegraphics{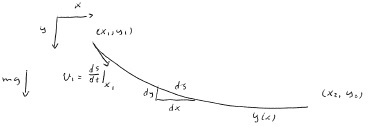}
\caption{The brachistochrone problem}
\label{fig11}
\end{figure}

\noindent We have 
\begin{align*}
v &= \frac{ds}{dt}, \; \;  \\ &\Updownarrow \\ \; \; dt &= \frac{ds}{v}. 
\end{align*}
Thus, the time it takes the particle to move from $(x_1, y_1)$ to $(x_2, y_2)$ along the arc $(x , y(x))$ is 
\begin{align*}
T(y) = \int^{x_2}_{x_1} dx \; \inv{v} \; \sqrt{1 + y'^2}.  
\end{align*}
We find a relation between $v$ and $y$ using the conservation of energy. We set the zero for potential energy at $y=y_1$. Then we have 
\begin{align*}
\frac{1}{2} \; m \; v_1^2 = \inv{2} \; m \; v^2 - m \; g \; (y - y_1),  
\end{align*}
which leads to 
\begin{align*}
v = \sqrt{2 \; g} \; \sqrt{y-y_0}, && y_0 = y_1 - \frac{v_1^2}{2g}, 
\end{align*}
and thus
\begin{align*}
T(y) = \inv{\sqrt{2g}} \; \int_{x_1}^{x_2} dx \; \sqrt{\frac{1+y'(x)^2}{y(x) - y_0}}.  
\end{align*}
The challenge is then to minimize $T(y)$ subject to the constraints 
\begin{align*}
y(x_1) = y_1, \;\;\; y(x_2) = y_2. 
\end{align*}

\subsubsection{The Action Principle.}\label{VariationalExample8}

Let us consider a system consisting of $N$ mass-points with positions $\{ \vx_i \}^N_{i=1}$, velocities $\{ \vx'_i \}^N_{i=1}$ and masses $m_i$. Let us assume that the mass-points are moving under the influence of a conservative force. Let $V(\vx_1,...,\vx_N)$ be the potential of this conservative force. Then by definition
\begin{align*}
\vb{f}_i = \frac{\prt{}V}{\prt{}\vx_i}, 
\end{align*}
 is the force acting on particle number $i$.  The kinetic energy of the system of particles is 
\begin{align*}
T(\vx'_1,...,\vx'_N)= \mathlarger{\sum}^N_{i=1} \inv{2} \; m_i \; \vx_i'^2,  
\end{align*}
and the \ttx{Lagrangian} of the system is by definition 
\begin{align*}
L(\vx_1,...,\vx_N,\vx'_1,...,\vx'_N) = T(\vx'_1,...,\vx'_N) - V(\vx_1,...,\vx_N).  
\end{align*}
A position vector 
\begin{align*}
\vb{P} = (\vx_1, ... , \vx_N) \in \mathbf{R}^{3N},  
\end{align*}
is a \ttx{configuration} for the system. Let a parametrized curve 
\begin{align*}
\vb{P}(t) = (\vx_1(t),...,\vx_N(t)), && t_1 \leq t \leq t_2, 
\end{align*}
in \ttx{configuration space} $\mathbf{R}^{3N}$ be given. The \ttx{action} of the parametrized curve is by definition
\begin{align*}
\mathcal{S}(\vb{P}) = \int^{t_2}_{t_1} dt \; \mathcal{L}.  
\end{align*}
The action principle says that the path, $\vb{P}(t)$, traced out in configuration space by a system of mass-points under the influence of conservative forces,  is the one that is \ttx{stationary} for the action. This means, by definition, that for all curves $\vb{q}(t)$ whose size is of order one, $\parallel\vb{q}\parallel=\obs(1)$, we have
\begin{equation*}
\mathcal{S}(\vb{P}+\epsilon\vb{q})=\mathcal{S}(\vb{P})+ \obs(\eps^2).
\end{equation*}
This approach to dynamics is very different from the usual one where we solve Newton's equations subject to given initial conditions. We will show that they are in fact equivalent. 
The action principle is also called \ttx{Hamilton's principle}. 

The action principle is the single most important idea in theoretical physics. \ttx{All} fundamental physical models are derived from the action principle, this is true for both classical physics and quantum physics. 

\subsubsection{The Maximum Entropy Principle}\label{VariationalExample9}

Let $x_1, ... , x_n$ be random variables with an associated probability distribution $\rho(x_1,...,x_n)$. Let $f_1(x_1,...,x_n),...,f_p(x_1,...,x_n)$ be functions defined on the space of random variables $\mathbf{R}^n$.
The functions $f_j$ are our \ttx{observables}. Their $\it{expectation}$ values are as usual defined by 
\begin{align*}
\left<f_j\right> \; = \int_{\mathbf{R}^n} dV \; f_j(x_1,...,x_n) \; \rho(x_1,...,x_n).  
\end{align*}
The expectation value of a given observable of course depends on which probability distribution $\rho$, we use. The challenge in statistics is to figure out which probability distribution one should use in any given situation. Let us say that we for some reason, (guesswork, hearsay, ...) believe that a probability distribution $\rho_0$,  accurately represents what we currently know about a given system. The probability distribution $\rho_0$ is called the \textit{prior distribution}, or just the \textit{prior}.

Let us next assume that we measure the mean values of the observables $f_1, ...,f_p$ and find the values $c_1,...,c_p$. If 
\begin{align*}
\left<f_j\right>_0 \;  = \int_{\mathbf{R}^n} dV \; f_j \; (x_1,...,x_n) \; \rho_0 (x_1,...,x_n) = c_j,  
\end{align*}
for $j=1,\cdots,p$, we are satisfied with our choice of prior. It predicts exactly the mean values that are observed. \\
But we might not be so lucky. Perhaps 
\begin{align*}
\left<f_j\right>_0 \; \ne c_j,  
\end{align*}
for at least one $j$. Our selected $\rho_0$ is then not the correct one, it predicts expectation values that are not observed. The challenge is to modify $\rho_0$ into a new distribution $\rho$ that is consistent with \ttx{all} the observed mean values. 

For this purpose we define a functional $S(\rho)$ by 
\begin{align*}
S(\rho) = - \int_{\mathbf{R}^n} dV \; \rho \; \ln (\frac{\rho}{\rho_0}).  
\end{align*}
$S$ is by definition the \ttx{relative entropy} of the probability distribution $\rho$ with respect  to $\rho_0$. We will see later that our use of the word entropy here is consistent with its usage in thermodynamics. 

The \ttx{maximum entropy} principle states that one should choose the probability distribution that maximizes the functional 
\begin{align*}
S(\rho) = - \int_{\mathbf{R}^n} dV \; \rho \; \ln (\frac{\rho}{\rho_0}),  
\end{align*}
subject to the constraints 
\begin{align*}
\left<f_j\right> \; = \int_{\mathbf{R}^n} dV \; f_j \; \rho = c_j. 
\end{align*}
\\

\subsection{The Euler-Lagrange Equations}
From examples given in the previous subsections, it should be clear that extremal problems for functionals play an important role in our description of nature. It is now time to find a way to solve such problems. 
\subsubsection {One dependent variable}
Several of the examples involved a functional of the form
\begin{align}
T(y) = \intt dt \; L (t, y , y'), \lbl{81.1}  
\end{align}
with constraints of the form
\begin{align}
y(t_0)&= y_0,\nonumber  \\
y(t_1) &= y_1. \lbl{82.1}
\end{align}
Note that the integrand, defining the functional \rf{81.1}, is in the context of the calculus of variation sometimes called an {\it integral density}, but more often a {\it Lagrangian}. The inspiration for the second name came originally from  from applications of the calculus of variation to mechanical systems. Joseph-Louis Lagrange reshaped the subject of particle mechanics in the 1700's and has had a large influence on how we today  think about the subject of mechanics. The name Lagrangian for the integral densities defining functionals has subsequently migrated to field theory, which forms the foundation for fundamental physics in general, and particle physics in particular. In these subjects  the functionals of interest are time integrals over a Lagrangian, which itself is the space integral over a function, that in this context is called a {\it Lagrangian density}. Theories in fundamental physics are {\it defined} in terms of Lagrangian  densities.

  In order to be consistent with the standard usage of the terms Lagrangian and Lagrangian density, we call the integrand in one of our functional for the Lagrangian, with symbol $L$, if the integral is over time, and Lagrangian density, with symbol $\mathcal{L}$,  if the functional is defined by an integral over space and  time, or over space alone.

Our challenge is to find a function $y(x)$ that satisfies the constraints \rf{82.1} and that is extremal for \rf{81.1}. 

Let us for a moment return to ordinary calculus. Let $f(x)$ be a function of a real variable, and let $x=x_0$ be some fixed value of $x$. The function $f(x)$ is by definition differentiable at $x=x_0$ if for $\epsilon << 1\; \text{and for all}\; h$ we have
\begin{align}
f(x_0 + \epsilon \; h) = f(x_0) + \epsilon \; A(x_0) \; h + \obs(\epsilon^2). \lbl{83.1}
\end{align}
$A(x_0)$ is by definition the derivative of $f$ at $x=x_0$. The point $x_0$ is $\it{stationary}$ for $f(x)$ if 
\begin{align}
f'(x_0) \equiv A(x_0) = 0.\lbl{84.12}
\end{align}
Thus at a stationary point we have
\begin{align}
f(x_0 + \epsilon \; h) = f(x_0) + \obs(\epsilon^2). \lbl{85.1} 
\end{align}
We know that for a differentiable function defined on the real line, all extremal points of a function $f(x)$ are found among the list of stationary points. Recall that extremal points are  local or global maximum points or minimum points.

Inspired by \rf{83.1}, \rf{84.12} and \rf{85.1} we say that a functional $T$ is differentiable at $y$ if for $\epsilon << 1$ and all $\eta(x)$ of order one, we have 
\begin{align}
T (y + \epsilon \; \eta) = T(y) + \epsilon \; A(y,\eta) + \obs(\epsilon^2), \lbl{86.1}
\end{align}
where the map 
\begin{align*}
\eta \rightarrow A(y, \eta),  
\end{align*}
is \ttx{linear} in $\eta$. This linear map is called the \ttx{variational derivative} of $T$ at $y(x)$, and in these lecture notes we denote it by $\delta \; T(y)$, thus
\begin{equation*}
\delta T(y)(\eta)=A(y,\eta).
\end{equation*}
The function $y=y(x)$ is stationary for $T(y)$ if 
\begin{align*}
\delta \; T(y) &= 0, \; \; \nonumber \\ &\Updownarrow \nonumber \\ \; \; \delta \; T(y)(\eta) &= 0 \; \; \; \; \; \; \; \forall \eta(x).
\end{align*}
Thus at a stationary point we have 
\begin{align*}
T(y + \epsilon \; \eta) = T(y) + \obs(\epsilon^2) && \forall \eta(x).
\end{align*}
Just like in ordinary calculus, extremal points for a functional are found among the list of stationary points.

Let us now calculate  the variational derivative of the functional \rf{81.1}. 
\begin{align}
T(y+\epsilon \; \eta) &= \intt dt \; L (t,y+ \epsilon \; \eta, y' + \epsilon \; \eta') \nonumber\\
&= \intt dt \; \{ L (t,y,y') + \frac{\prt{} L}{\prt{} y} \; \epsilon \eta + \frac{\prt{}L}{\prt{} y'} \; \epsilon \; \eta' + \obs(\epsilon^2) \} \nonumber \\
&= T(y) + \epsilon \; \intt dt \; \{ \frac{\prt{} L}{\prt{}y} \; \eta + \frac{\prt{} L}{\prt{} y'} \; \eta' \} + \obs(\epsilon^2). \lbl{91.1}
\end{align}
The functions $y(x)$ and $y(x) + \epsilon \; \eta(x)$ must satisfy the constraints \rf{82.1}, only such functions are relevant. Thus we have 
\begin{align}
y(t_0) + \epsilon \; \eta(t_0) = y_0 = y(t_0) \; \; \; \Rightarrow \; \; \eta(t_0) = 0,\nonumber \\ y(t_1) + \epsilon \; \eta (t_1) = y_1 = y(t_1) \; \; \; \Rightarrow \; \; \eta(t_1) = 0.  \lbl{92.1}
\end{align}
Using integration by parts and \rf{92.1} we have 
\begin{align}
\intt dt \; \frac{\prt{}L}{\prt{} y'} \; \eta' &= \frac{\prt{}L}{\prt{} y'} \; \eta \mathlarger{|}^{t_1}_{t_0} - \intt dt \; \dt{} (\frac{\prt{} L}{\prt{} y'}) \; \eta \lbl{93.1} \\ 
&=- \intt dt \; \dt{} (\frac{\prt{} L}{\prt{} y'}) \; \eta.  \nonumber 
\end{align}
Inserting \rf{93.1} into \rf{91.1} we get 
\begin{align*}
T (y+ \epsilon \; \eta) = T (y) + \epsilon \; \intt dt \; \{ \prtdLy - \dt{} (\prtdLyd) \} \; \eta + \obs(\eps^2) && \forall \eta(t). 
\end{align*} 
Thus, by the definition of the variational derivative\;\rf{86.1}, we have
\begin{equation*}
\delta\;T(y)(\eta)=\intt dt \; \{ \prtdLy - \dt{} (\prtdLyd) \} \; \eta.
\end{equation*}
We have the following result; the so-called the fundamental lemma of variational calculus. 
\\\\
\textbf{lemma} 
\\
Let $a(t)$ be a continuous function and assume that 
\begin{align}
\intt dt \; a(t) \; \eta(t) = 0, \lbl{95.1} 
\end{align}
for \ttx{all} continuous functions $\eta$. Then 
\begin{align}
a(t) = 0 && t_0 \leq t \leq t_1.  \nonumber 
\end{align}
If we apply the lemma \rf{95.1} we find that $y(x)$ is a stationary point for $T$ if $y(x)$ satisfies the equation 
\begin{align*}
\prtdLy - \dt{} (\prtdLyd) = 0.
\end{align*}
This equation is called the \ttx{Euler-Lagrange} equation. 

\begin{example}\label{Example1}
In section \ref{VariationalExample4} we discussed the problem of finding a curve which, upon rotation around the x-axis, generates a surface of minimal area. 

We observed that such a curve, $y(x)$,  minimized the functional 
\begin{align*}
A(y) = 2 \; \pi \; \int^{x_2}_{x_1} dx \; y(x) \; \sqrt{1 + y'(x)^2}, 
\end{align*}
subject to the constraints
\begin{align*}
y(x_1) &= y_1, \\
y(x_2) &= y_2. 
\end{align*}
Just as in ordinary calculus the minimum will be found among the stationary points for $A$. 

Let us therefore find the Euler-Lagrange equation for the functional $A$. For $A$ we evidently the Lagrangian density is given by
\begin{align}
\mathcal{L} = 2 \; \pi \; y \; \sqrt{1 + y'^2}.  \lbl{99.12}
\end{align}
Thus 
\begin{align*}
\prtdLdy &= 2 \; \pi \; \sqrt{1 + y'^2},\\
\prtdLdyd &= 2 \; \pi \; y \; \frac{y'}{\sqrt{1+y'^2}}, 
\end{align*}
and therefore the Euler-Lagrange equation is 
\begin{align}
2 \; \pi \; \sqrt{1 + y'^2} - \dt{} (2 \; \pi \; y \; \frac{y'}{\sqrt{1+y'^2}}) = 0. \lbl{101.12} 
\end{align}
\end{example}
\noindent Any stationary point $y(x)$ is a solution to this equation. 
Equation \rf{101.12} is a highly nonlinear second order differential equation. In fact, the Euler-Lagrange equation will always be second order for functionals of the form \rf{81.1}. \\
Note that the the Lagrange density $\mathcal{L}$, in formula  \rf{99.12}, defining the functional \rf{81.1},  is not totally general since it does not depend explicitly on $x$. \\
For such Lagrangians densities we have in general 
\begin{align*}
\dx{}[y' \; \prtdLdyd - \mathcal{L}] &= y'' \; \prtdLdyd + y' \; \dx{} (\prtdLdyd)\\
- \prtdLdy \; y' - \prtdLdyd \; y'' &= y' \; (\dx{} (\prtdLdyd) - \prtdLdyd) = 0,
\end{align*}
where we in the last step used the Euler-Lagrange equation. \\
This calculation shows that for an $\mathcal{L}$ that does not explicitly depends on $x$, any extremal of $A$ is a solution to the following \ttx{first} order differential equation 
\begin{align}
y' \; \prtdLdyd - \mathcal{L} = c, \lbl{103.12} 
\end{align}
where $c$ is a constant. This reduction of order for such special Lagrangian densities, is not just a lucky break. Behind this result there is a hugely important mathematical machine called \ttx{Noether's theorem}. It will give similar reductions of order in many other unrelated situations. 

Noether's theorem has for almost a century been at the center of the action  in theoretical physics, both quantum and classical. 
In quantum theory Noether's theorem, in the form of \ttx{the gauge principle}, is used to derive \ttx{all} the fundamental equations in elementary particle physics. These equations, taken together, form what is modestly called \ttx{the Standard Model}. This model has predicted the outcome of \ttx{all} experiments in fundamental physics since the 1970's. It is the most accurate theory of nature that has ever been created, and it all flows from Noether's theorem. \\
Using \rf{103.12} for the Lagrangian density \rf{99.12}, we get 
\begin{align}
y' \; 2 \; \pi \; y \; \frac{y'}{\sqrt{1+y'^2}} - 2 \; \pi \; y \; \sqrt{1+y'^2} &= c, \nonumber \\
&\Updownarrow \nonumber \\  -\frac{y}{\sqrt{1+y'^2}} = c_1 &\equiv \frac{c}{2 \; \pi}. \lbl{104.12} 
\end{align}
Since $y(x)$ from section \ref{VariationalExample4} is positive, equation \rf{104.12} can have solutions only if $c_1 < 0$. We can therefore, without loss og generality, write $c_1 = -\alpha$, where $\alpha > 0$. Using this we have
\begin{align}
\frac{y}{\sqrt{1+y'^2}} &= \alpha, \nonumber \\
&\Updownarrow \nonumber \\ y^2 &= \alpha^2 \; (1+y'^2), \nonumber \\ 
&\Updownarrow \nonumber \\  y'^2 &= (\frac{y}{\alpha})^2 - 1, \nonumber \\ 
&\Updownarrow \nonumber \\ y' &= \pm \; \sqrt{(\frac{y}{\alpha})^2 - 1}. \lbl{105.12} 
\end{align}
Equation \rf{105.12} is separable and can be solved. The general solution is 
\begin{align}
y(x) = \alpha \; \cosh (\frac{x}{\alpha} + \beta) && \alpha > 0, \; \; \beta \in \mathbf{R}, \lbl{106.12}
\end{align}
where the constants $\alpha$ and $\beta$ are determined from the conditions 
\begin{align*}
y(x_1) = y_1, && y(x_2) = y_2. 
\end{align*}
A graph of functions of the type \rf{106} is called a \ttx{Catenary}.

By construction, functions of the type \rf{106.12} are stationary points for the area functional $A(y)$. In ordinary calculus we use the second derivative test in order to decide if a stationary point is a local maximum or a local minimum. Global maximum or minimum will be found among the local maximums or minimums \ttx{or} at points where the function is singular (not differentiable). This is true assuming there are no boundary points. If there are boundary points, maximum and minimum can occur there also. 

\noindent The same rules apply here in the \ttx{calculus of variation}, which is the official name for what we are doing. We will however not develop the theory further in this direction and will not discuss the second variational derivative. 

The problem of deciding whether a given stationary point gives a global minimum or a global maximum must be investigated in each separate case. Solving this problem can be highly non-trivial. \\
For the stationary points \rf{106.12}, it is easy easy to verify that a Catenary passing through the point $(x_1,y_1)$ is given by 
\begin{align*}
y(x) = \alpha \; \cosh (\frac{x - x_1}{\alpha} + \cosh^{-1} (\frac{y_1}{\alpha})). 
\end{align*}
In order to make sure that $y(x)$ also pass through the point $(x_2, y_2)$ we must find an $\alpha$ such that 
\begin{align*}
\alpha \; \cosh (\frac{x_2 - x_2}{\alpha} + \cosh^{-1}(\frac{y_1}{\alpha})) = y_2.
\end{align*}
This is a transcendental equation for $\alpha$ and depending on coordinates, $(x_1,y_1),(x_2,y_2)$ . It might have no positive solutions or several positive solutions. \\
The following cases are known to occur 
\begin{align}
& \mathbf{i)} \; \; \; \text{There exists no Catenary connecting the points $(x_1, y_1)$ and $(x_2,y_2)$}, \nonumber \\
& \mathbf{ii)} \;  \; \text{There exists exactly one Catenary connecting $(x_1, y_1)$ and $(x_2,y_2)$}, \nonumber \\
& \mathbf{iii)} \;  \text{There exists exactly two Catenaries connecting $(x_1, y_1)$ and $(x_2,y_2)$}. \nonumber
\end{align}
In case $\mathbf{ii)}$ the unique Catenary is neither a global nor  local minimum. In case $\mathbf{iii)}$ one of the solutions is a local minimum, and this is also a global minimum, if $y_1$ is large enough in a certain sense. 
Thus a surface of revolution of minimum area exists only in case $\mathbf{iii)}$ and only if $y_1$ is large enough. 

Clearly, the Euler-Lagrange equation only gives stationary points and these points are in general not even \ttx{local} maximum or minimum. \\
However, recall that in some cases, like the action principle, we are only looking for stationary points. 

\begin{example}\label{Example2}
The Brachistochrone, from section \ref{VariationalExample7},  is defined to be the minimum of the functional 
\begin{align*}
T(y) = \inv{\sqrt{2 \; g}} \; \int^{x_2}_{x_1} dx \; (\frac{1+y'^2}{y - y_0})^{\inv{2}}.
\end{align*}
Thus, the Lagrangian  for this problem is
\begin{align}
\mathcal{L} = \inv{\sqrt{2 \; g}} \; (1 + y'^2)^{\inv{2}} \; (y - y_0)^{-\inv{2}}. \lbl{111.12} 
\end{align}
Observe that the Lagrangian density\rf{111.12} does not depend explicitly on the independent variable $x$. Therefore, using the result from \rf{103.12}, we know that any extremal of the functional $T$ is a solution to the following first order ODE
\begin{align}
y' \; \prtdLdyd - \mathcal{L} = c, \lbl{121.1} 
\end{align}
Differentiating $\mathcal{L}$ we get 
\begin{align}
\prtdLdyd &= \inv{\sqrt{2 \; g}} \; y' \; (1 + y'^2)^{- \inv{2}} \; (y - y_0)^{-\inv{2}},\lbl{112.1}
\end{align}
and inserting \rf{112.1} into \rf{121.1} we get
\begin{align}
y' \; \prtdLdyd - \mathcal{L} &= c, \nonumber \\ 
&\Updownarrow \nonumber \\  y'^2 \; (1 + y'^2)^{-\inv{2}} \; (y - y_0)^{-\inv{2}} - (1+y'^2)^{\inv{2}} \; (y - y_0)^{-\inv{2}} = c_1 &\equiv c \; \sqrt{2 \; g}, \nonumber \\
&\Updownarrow \nonumber \\  y'^2 = 1 + y'^2 + c_1 \; (1+y'^2)^{\inv{2}} \; (y-y_0)^{\inv{2}}, \nonumber \\ 
&\Updownarrow \nonumber \\ c_1 \; (1+y'^2)^{\inv{2}} \; (y - y_0)^{\inv{2}} &= -1. \lbl{113.12} 
\end{align}
For there to be a solution we must have $c_1 < 0$. Let $c_1 = -\beta, \; \beta > 0$. From \rf{113.12} we then get 
\begin{align}
(1+y'^2)^{\inv{2}} \; (y-y_0)^{\inv{2}} &= \inv{\beta}, \nonumber\\ 
&\Updownarrow \nonumber \\  1 + y'^2 &= \inv{\beta^2 \; (y -y_0)}, \nonumber \\ 
&\Updownarrow \nonumber \\  y'^2 &= \frac{1 - \beta^2 \; (y-y_0)}{\beta^2 \; (y-y_0)}, \nonumber \\
&\Updownarrow \nonumber \\ y' &= (\frac{\alpha^2 - (y-y_0)}{y-y_0})^{\inv{2}} && \alpha^2 = \inv{\beta^2} > 0, \lbl{115.12}
\end{align}
where we have used the positive root since we know that for the current example $y' > 0$.\\ The equation \rf{115.12} is  separable, and we get the following implicit solution 
\begin{align}
\int dy \; (\frac{y - y_0}{\alpha^2 - (y-y_0)})^{\inv{2}} = x + c, \lbl{116.12} 
\end{align}
where $c$ is an arbitrary integration constant. \\ 
The integral in \rf{116.12} can be solved using a trigonometric substitution 
\begin{align}
y - y_0 = \alpha^2 \; \sin^2 \frac{\theta}{2}. \lbl{117.12}
\end{align}
Inserting \rf{117.12} into \rf{116.12} give us 
\begin{align*}
\alpha^2 \; \int d\theta \; \sin^2 \frac{\theta}{2} = x + c, \\
&\Updownarrow \nonumber \\  \inv{2} \; \alpha^2 \; (\theta - \sin \theta ) = x+c.
\end{align*}
Thus our solution is 
\begin{align*}
x &= \inv{2} \; \alpha^2 \; (\theta - \sin \theta) - c,  \\
y &= y_0 + \inv{2} \; \alpha^2 \; (1-\cos \theta).
\end{align*}
\end{example}
This is a parametric representation of a type of curve called a \ttx{cycloid}. One can show that there is a unique cycloid passing through any pair of points $(x_1, y_1), \; (x_2, y_2)$ with $y_1 < y_2$ and that the unique cycloid is a global minimum for the propagation time functional $T(y)$. \\

\subsubsection{Several dependent variables}
Several of our examples involved functionals of the general form 
\begin{align*}
T(y_1,...,y_n) = \intt dt \; L (t, y_1,...,y_n.y_1',...,y_n'). 
\end{align*}
By analogy with multi-variable calculus, we consider small independent perturbations of all the functions $\{ y_i(x) \}_{i=1}^n$. Proceeding like for the case of one function, we get 
\begin{align*}
&T(y_1 + \eps\; \eta_1, ... , y_n + \eps \; \eta_n) \\
&= \intt dt L(t_1, y_1 + \eps \; \eta_1, ..., y_n + \eps \; \eta_n, y'_1 + \eps \; \eta_1', ..., y_n' + \eps \; \eta_n') \nonumber \\ 
&= \intt dt \; \{ L(t_1 y_1,...,y_n,y_1',...,y_n') + \frac{\prt{} L}{\prt{} y_1} \eps \; \eta_1  \nonumber \\ 
& +... + \frac{\prt{} L}{\prt{} y_n} \; \eps \; \eta_n + \frac{\prt{} L}{\prt{} y_1'} \; \eps \; \eta_1' + ... + \frac{\prt{} L}{\prt{} y_n'} \; \eps \; \eta_n' \} \nonumber \\ 
&+ \obs(\eps^2) \nonumber \\ 
&= T (y_1,...,y_n) \nonumber \\ 
&+ \eps \; \intt dt \; \{ (\frac{\prt{} L}{\prt{} y_1} - \dt{} (\frac{\prt{} L}{\prt{} y_1'})) \; \eta_1 + ... + (\frac{\prt{} L}{\prt{} y_n} - \dt{} (\frac{\prt{} L}{\prt{} y_n'})) \; \eta_n \} + \obs(\eps^2), \nonumber  
\end{align*}
where we in the last step has used integration by parts and the fact that 
\begin{align}
\eta_j(t_0) = \eta_j(t_1) = 0 && j = 1,...,n, \lbl{122.12} 
\end{align}
The relations \rf{122.12} follows from the fact we have constraints 
\begin{align*}
y_j(t_0) = y_j^0, && y_j(t_1) = y_j^1, 
\end{align*}
by an argument that is identical to the one in  equation \rf{92.1}. A stationary point for the functional $T(y_1, ... , y_n)$ is determined by the condition that 
\begin{align*}
\intt dt \; \{ (\frac{\prt{}L}{\partial y_1} - \dt{} (\frac{\partial L}{\partial y_1'})) \; \eta_1 + ... + (\frac{\partial L}{\partial y_n} - \dt{} (\frac{\partial L}{\partial y_n'})) \; \eta_n \}  \\
=0, \nonumber 
\end{align*}
for \ttx{all} functions $\{ n_i \}^n_{i=1}$. Since there by assumption are no dependencies among the functions $\eta_i(t)$, we conclude, using the fundamental lemma \rf{95.1}, that $(y_1,...,y_n)$ is a stationary point for $T(y_1,...,y_n)$ only if 
\begin{align}
\frac{\partial L}{\partial y_i} - \dt{} (\frac{\partial L}{\partial y_i'}) = 0 && i = 1...n. \lbl{125.12} 
\end{align}
These are the Euler-Lagrange equations for the functional $T(y_1,...,y_n)$. Observe that, in general, \rf{125.12} are $n$ coupled non-linear second order differential equations.

\begin{example}\label{Example3}
In \ref{VariationalExample1}. we found that the problem of finding the shortest curve in the plane, connecting two fixed points $p=(x_0,y_0)$ and $q=(x_1,y_1)$, amounted to minimizing the functional 
\begin{align*}
T(x,y) = \int^1_0 dt \; \sqrt{x'^2 + y'^2},  
\end{align*}
subject to the constraints 
\begin{align*}
(x(0), y(0)) = p, && (x(1), y(1)) = q.  
\end{align*}
The Lagrangian density is 
\begin{align*}
\mathcal{L}(x,y) = \sqrt{x'^2 + y'^2}.  
\end{align*}
Observe that
\begin{align*}
\frac{\partial \mathcal{L}}{\partial x} &= 0, && \frac{\partial \mathcal{L}}{\partial y} = 0,  \\ 
\frac{\partial \mathcal{L}}{\partial x'} &= \frac{x'}{\sqrt{x'^2 + y'^2}}, && \frac{\partial \mathcal{L}}{\partial y'} = \frac{y'}{\sqrt{x'^2 + y'^2}}. \nonumber 
\end{align*}
Thus the Euler-Lagrange equations are 
\begin{align}
\dt{} (\frac{x'}{\sqrt{x'^2 + y'^2}}) &= 0, \; \; \nonumber \\ &\Updownarrow \nonumber \\  \; \; \frac{x'}{\sqrt{x'^2 + y'^2}} &= c_1, \nonumber\\ 
\dt{} (\frac{y'}{\sqrt{x'^2 + y'^2}}) &= 0, \; \; \nonumber \\ &\Updownarrow \nonumber \\  \; \; \frac{y'}{\sqrt{x'^2 + y'^2}} &= c_2. \lbl{130.12}
\end{align}
From \rf{130.12} we find that 
\begin{align*}
\frac{dy}{dx} = \frac{\dt{y}}{\dt{x}} &= \frac{c_2}{c_1}, \; \;\nonumber \\  &\Updownarrow \nonumber \\ \; \; y &= \frac{c_2}{c_1} \; x + c_3,  
\end{align*}
and the curve passing through $p$ and $q$ is 
\begin{align*}
y(x) = \frac{y_1 - y_0}{x_1 - x_0} \; (x - x_0) + y_0.  
\end{align*}
\end{example}
\noindent This is a straight line and is evidently the curve of minimal length. We of course already knew that the straight line is the shortest curve connecting two points in the plane. \\
In  \ref{VariationalExample2} we asked the same question for two points on a surface in $\mathbf{R}^3$. Here the answer is not obvious, a straight line in $\mathbf{R}^3$ will not work unless $S$ is a plane. We are in some sense seeking a curve on a general curved surface that is the analogy of straight lines in $\mathbf{R}^3$. \\
When the family of such curves has been found we can use them to construct analogues of triangles, squares etc.,  on the curved surface and ask geometrical questions like: What is the sum of the internal angles of a triangle on a given surface $S$? We can in fact develop a whole analogue to Euclidean geometry for plane figures, for any surface, not only plane surfaces. 

This has been done for many kind of surfaces and is of obvious practical importance for the case of a sphere. (Fuel efficient long distance transport)

\begin{example}\label{Example4}
In section \ref{VariationalExample6} we introduced the Fermat principle determining the light-rays in a material of variable refractive index. \\
The functional was 
\begin{align*}
T(x,y,z) = \int^1_0 dt \; n(x(t), y(t), z(t))\sqrt{x'^2(t) + y'^2(t) + z'^2(t)}.  
\end{align*}
The Lagrangian is 
\begin{align*}
L = n(x,y,z) \; \sqrt{x'^2 + y'^2 + z'^2} = n\norm{\vb{x}'}  
\end{align*}
Observe that 

\begin{align*}
\frac{\partial L}{\partial x} &= \prt{x}n\norm{\vb{x}'}, && \frac{\partial L}{\partial x'} = \frac{n \; x'}{\norm{\vb{x}'}},  \\
\frac{\partial L}{\partial y} &= \prt{y}n\norm{\vb{x}'}, && \frac{\partial L}{\partial y'} = \frac{n \; y'}{\norm{\vb{x}'}}, \nonumber \\
\frac{\partial L}{\partial z} &= \prt{z}n\norm{\vb{x}'}, && \frac{\partial L}{\partial z'} = \frac{n \; z'}{\norm{\vb{x}'}}. \nonumber 
\end{align*}
The Euler-Lagrange equations are 
\begin{align}
\prt{x} n\norm{\vb{x}'} - \dt{} (\frac{n \; x'}{\norm{\vb{x}'}}) &= 0, \lbl{136.12} \\
\prt{y} n\norm{\vb{x}'} - \dt{} (\frac{n \; y'}{\norm{\vb{x}'}}) &= 0, \nonumber \\
\prt{z} n\norm{\vb{x}'} - \dt{} (\frac{n \; z'}{\norm{\vb{x}'}}) &= 0. \nonumber
\end{align}
Consider the case of a homogeneous medium where $n(x,y,z) = n_0$. For this case \rf{136.12} simplifies into 
\begin{align*}
\dt{} (\frac{x'}{\norm{\vb{x}'}}) &= 0, \; \; \nonumber \\  &\Updownarrow \nonumber \\ \; \; \dt{x} &= c_1 \; \norm{\vb{x}'},  \\
\dt{} (\frac{y'}{\norm{\vb{x}'}}) &= 0, \; \; \nonumber \\  &\Updownarrow \nonumber \\  \; \; \dt{y} &= c_2 \; \norm{\vb{x}'}, \nonumber \\
\dt{} (\frac{z'}{\norm{\vb{x}'}}) &= 0, \; \; \nonumber \\  &\Updownarrow \nonumber \\  \; \; \dt{z} &= c_3 \; \norm{\vb{x}'}. \nonumber
\end{align*}
Using $x$ as a new independent variables assuming $c_1 \ne 0$ we get in the usual way
\begin{align*}
\frac{dy}{dx} = \frac{\dt{y}}{\dt{x}} &= \frac{c_2}{c_1}, \; \; \nonumber \\ &\Updownarrow \nonumber \\  \; \; y &= \frac{c_2}{c_1} \; x + c_4,  \\
\frac{dz}{dx} = \frac{\dt{z}}{\dt{x}} &= \frac{c_3}{c_1}, \; \; \nonumber \\  &\Updownarrow \nonumber \\  \; \; y &= \frac{c_3}{c_1} \; x + c_5. \nonumber
\end{align*}
Thus the light-rays in a homogeneous medium are straight lines. We of course know this from elementary physics. \\
For the general case, we introduce path-length as a new parameter 
\begin{align*}
s(t) = \int_{t_0}^t dt' \; \norm{\vb{x}'(t')}.  
\end{align*}
Using this as our curve parameter we have, using the chain rule, that
\begin{equation*}
\frac{d}{dt}=\norm{\vb{x}'}\frac{d}{ds}.
\end{equation*}
Thus we have
\begin{align*}
\partial_x n \; \norm{\vb{x}'}- \dt{} (\frac{n \; x'}{\norm{\vb{x}'}}) &= 0,\nonumber \\
&\Updownarrow\nonumber\\ \partial_x n \; \norm{\vb{x}'}- \dt{} (\frac{n}{\norm{\vb{x}'}}\frac{dx}{dt}) &= 0,\nonumber \\
&\Updownarrow \nonumber \\ \partial_x n \; \norm{\vb{x}'}- \norm{\vb{x}'}\frac{d}{ds}(\frac{n}{\norm{\vb{x}'}}\norm{\vb{x}'}\frac{d}{ds}) &= 0, \nonumber \\ 
&\Updownarrow \nonumber \\ \prt{x} n \norm{\vb{x}'}- \norm{\vb{x}'}\frac{d}{ds} (n \; \frac{dx}{ds}) &= 0, \nonumber \\
&\Updownarrow \nonumber \\  \prt{x} n - \frac{d}{ds} (n \; \frac{dx}{ds}) &= 0.  
\end{align*}
We rewrite the other two Euler-Lagrange equations in the same way. Collecting the three scalar equations into one vector equation we have 
\begin{align}
\frac{d}{ds} (n \; \frac{d \vb{x}}{ds}) = \grad{n}. \lbl{141.12} 
\end{align}
This is the fundamental equation for \ttx{ray optics}. \\
We are not going to solve this equation but will make a general observation. Since we are using path-length parametrization, the tangent vector
\begin{align*}
\vT = \frac{d\vb{x}}{ds},   
\end{align*}
is a \ttx{unit} vector. This means that 
\begin{align*}
\vT \vdot \vT &= 1,  \\
&\Updownarrow \nonumber \\ \frac{d\vT}{ds} \vdot \vT &= 0, \nonumber \\
&\Updownarrow \nonumber \\  \vb{\alpha} \vdot \vT &= 0, \nonumber 
\end{align*}
where $\vb{\alpha} = \frac{d\vT}{ds}$ is the \ttx{curvature vector} for the curve $\vb{x}(s)$. It is normal to the curve and points in the direction in which the curve bends
\begin{figure}[htbp]
\centering
\includegraphics{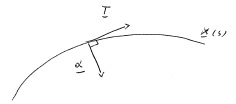}
\caption{Unit tangent $\vb{T}$ and curvature vector $\vb{\alpha}$,  for a light ray.\label{fig12.12}}
\end{figure}

From the fundamental equation for ray optics \rf{141.12} we have 
\begin{align*}
\frac{dn}{ds} \; \vT + n \; \frac{d \vT}{ds} &= \grad{n},  \\
\Downarrow\nonumber\\
 \vb{\alpha} &= \frac{\grad{n}}{n} - \frac{d}{ds} (\ln n) \; \vT, \nonumber \\
\Downarrow\nonumber\\
 0 < \vb{\alpha} \vdot \vb{\alpha} &= \inv{n} \; \vb{\alpha} \vdot \grad{n}, \nonumber \\
\Downarrow\nonumber\\
 \vb{\alpha} \vdot \grad{n} &> 0. \nonumber 
\end{align*}

\noindent Thus the curvature vector for  the light-ray points in the direction of increasing refractive index.
\begin{figure}[htbp]
\centering
\includegraphics{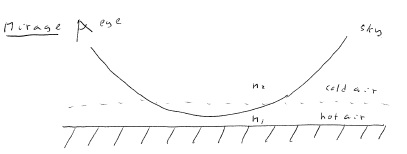}
\caption{Explaining the occurrence of a Mirage using ray optics.\label{fig13.12}}
\end{figure}
This means that light-rays in an inhomogeneous medium always curve towards regions of higher refractive index. This fact has many important physical consequences. 
Hot air close to the ground is less dense and has lower refractive index than the higher colder air. The light coming from the sky will therefore bend upwards, The sky will appear as "water" on the ground to an observer. This is illustrated in figure \ref{fig13.12}.

Another consequence of practical importance is that light rays will be confined to move within a region of higher refractive index: The optical fibers which forms the physical backbone of the internet works because of this fact. 
\end{example}

\begin{example}\label{Example5}\label{MassPoints}
In \ref{VariationalExample8}.  we introduced the action principle for a system of mass-points moving under the influence of conservative forces. \\
The action functional was 
\begin{align}
\mathcal{S}(\vx_1,...,\vx_N)= \intt dt \; \{T(\vx'_1(t),...,\vx'_N(t))-V(\vx_1(t),...,\vx_N(t))\}, \lbl{148.12} 
\end{align}
where $T$ is the kinetic energy of the system of mass-points and $V$ is the potential determining the conservative forces. \\
The Lagrangian is thus 
\begin{align}
L = \inv{2} \; \mathlarger{\sum}^n_{i=1} m_i \; \vb{x}_i'^2 - V(\vx_1,...,\vx_n). \lbl{149.12} 
\end{align}
Observe that 
\begin{align*}
\frac{\partial L}{\partial \vx_i} = - \frac{\partial V}{\partial \vx_i}, && \frac{\partial L}{\partial \vx_i'} = m_i \; \vx_i', 
\end{align*}
and thus the Euler-Lagrange equations are 
\begin{align*}
\frac{\partial L}{\partial \vx_i} - \dt{} (\frac{\partial L}{\partial \vx_i'}) &= 0 && i = 1,...,n\;,  \\ 
&\Updownarrow \nonumber \\  - \frac{\partial V}{\partial \vx_i} - m_i \; \vx_i'' &= 0, \nonumber \\ 
&\Updownarrow \nonumber \\  m_i \; \vx_i'' &= - \frac{\partial V}{\partial \vx_i} \equiv \vb{f}_i. \nonumber 
\end{align*}
This is Newton's law for $n$ mass-points, $m_i$, moving under the influence of conservative forces 
\begin{align*}
\vb{f}_i = - \frac{\partial V}{\partial \vx_i}.  
\end{align*}
Note that the Lagrangian \rf{149.12} does not depend explicitly on time. Let 
\begin{align*}
E = \mathlarger{\sum}^n_{i=1} \vx_i' \vdot \frac{\partial L}{\partial \vx_i'} - \mathcal{L}.  
\end{align*}
Then, using the Euler-Lagrange equation we have 
\begin{align*}
\frac{dE}{dt} &= \mathlarger{\sum}^n_{i=1} \{ \vx_i'' \vdot \frac{\partial L}{\partial \vx_i'} + \vx_i' \vdot \dt{} (\frac{\partial L}{\partial \vx_i'}) \}  \\
&- \mathlarger{\sum}^n_{i=1} \{ \frac{\partial L}{\partial \vx_i} \vdot \vx_i' - \frac{\partial L}{\partial \vx_i'} \vdot \vx_i'' \} \nonumber \\ 
&= \mathlarger{\sum}^n_{i= 1} \vx_i' \cdot\{ \dt{} (\frac{\partial L}{\partial \vx_i'}) - \frac{\partial L}{\partial \vx_i} \} = 0, \nonumber 
\end{align*}
and thus $E$ is a constant of the motion for any solution of the Euler-Lagrange equation, and therefore for any stationary point for the action functional \rf{148.12}. \\
Using \rf{150} we find that 
\begin{align}
E &= \mathlarger{\sum}^n_{i=1} \vx_i' \vdot (m_i \; \vx_i') - \mathlarger{\sum}^n_{i=1} \inv{2} \; m_i \; \vx_i^2 +V\nonumber\\ 
 &= \inv{2} \; \mathlarger{\sum}^n_{i=1} m_i \; \vx_i'^2 + V = T+V. \lbl{155.12}  
\end{align}
The conserved quantity is as expected the total energy of the system of mass-points. The function $E$ in \rf{155.12} is called the \ttx{Hamiltonian} for the system and is traditionally written using the letter $H$. Thus 
\begin{align*}
L &= T- V && \text{-Lagrangian},  \\
H &= T + V && \text{-Hamiltonian}. \nonumber 
\end{align*}
\end{example}
\noindent As we have just seen, the action principle leads to Euler-Lagrange equations that are equivalent to Newton's law. Thus from this point of view, nothing appears to have been gained from using the action principle. \\
However, the action principle, when it is taken together with the theory of canonical transformations, is a much more efficient tool for solving mechanical problems that Newton's law. 

Also, insights gained from the use of the action principle in mechanics played a major role in the development of quantum mechanics. Because of this, the action principle plays an important and ever increasing role in theoretical (and practical!) physics today. 
\\

\subsubsection{Constraints} 
Let us consider the general problem of finding the extremum of a functional 
\begin{align}
T (y) = \intt dt \; L(t, y, y'), \lbl{157.121} 
\end{align}
subject to the constraint 
\begin{align}
J(y) = \intt dt \; j(t, y, y') = c. \lbl{158.12} 
\end{align}
We solve this problem using the Lagrange multiplier approach. Let us recall how this approach works in the case of function on $\textbf{R}^n$.
 
  For this case, $f: \textbf{R}^n \rightarrow \mathbf{R} , \; \; g: \mathbf{R}^n \rightarrow \mathbf{R}$ are functions on $\textbf{R}^n$. The challenge is to find a stationary  point of $f$ subject to the constraint 
\begin{align}
g(x_1,...,x_n) = c. \label{159.12}
\end{align} 
In calculus one proves that stationary points for the function $f$ under the constraint (\ref{159.12}) determined by the function $g$,  can be found by introducing a new function 
\begin{align*}
h(x_1,...,x_n) = f(x_1,...,x_n) + \lambda \; (g(x_1,...,x_n) - c),  
\end{align*}
and then seek an unconstrained stationary point for $h$. The parameter $\lambda$ is called a \ttx{Lagrange multiplier}. \\
We thus solve the equation
\begin{align*}
\grad{h} &= 0,  \\
&\Updownarrow \nonumber \\ \grad{(f + \lambda \; (g-c))} &= 0, \nonumber \\
&\Updownarrow \nonumber \\ \grad{f} &= - \lambda \; \grad{g}. \nonumber 
\end{align*}
The value for $\lambda$ is chosen so that the constraint 
\begin{align*}
g(x_1,...,x_n) = c.  
\end{align*}
is satisfied. \\
The exact same approach works for constrained variational problems. We will not prove this fact. \\
Thus in order to find a stationary point for the functional \rf{157.121} subject to the constraint \rf{158.12}, we introduce the new functional 
\begin{align*}
K(y) = T(y) + \lambda \; (J(y) - c),  
\end{align*}
and find the unconstrained stationary points for $k$. Using the notation $\delta R$ for the variational derivatives of a functional $R$, we have
\begin{align}
\delta K &= 0,\nonumber \\
&\Updownarrow \nonumber \\  \delta \; T + \lambda \; \delta \; J &= 0. \lbl{164.121}
\end{align}
The value of the Lagrange multiplier is chosen so that 
\begin{align*}
J(y^*) = c,  
\end{align*}
where $y^*$ is a solution of \rf{164.121}. Observe that we have 
\begin{align*}
K(y) &= \intt dt \; L(t,y,y') + \lambda \; ( \intt dt \; j(t, y, y') - c)  \\
&= \intt dt \; \{ L(t, y, y') + \lambda \; (j(t,y,y') - \frac{c}{T}) \}, \nonumber 
\end{align*}
where $T=t_1-t_0$. Thus the Lagrangian for $K$ is 
\begin{align*}
L^* = L + \lambda \; (j - \frac{c}{T}).  
\end{align*}
Stationary points for the functional  $K$ are as usual solutions to the corresponding  Euler-Lagrange equations 
\begin{align*}
\frac{\partial L^{*}}{\partial y} - \dt{} (\frac{\partial L^{*}}{\partial y'}) = 0.  
\end{align*}
The Euler-Lagrange equation is here a scalar ODE of second order, and therefore the general solution contains two integration constants. These two integration constants, together with $\lambda$, are determined by the following three constraints 
\begin{align*}
y(t_0) &= y_0,  \\
y(t_1) &= y_1, \nonumber \\
J(y) &= c. \nonumber
\end{align*}
This procedure can obviously be extended to any number of constraints and any number of dependent variables: If we want to find stationary points for the functional 
\begin{align*}
T (y_1, ..., y_n) = \intt dt \; L(t,y_1,...,y_n,y_1',...,y_n'),  
\end{align*}
subject to the constraints 
\begin{align*}
J_p(y_1,...,y_n) = c_p && p = 1,...,m\;\;,  
\end{align*}
where 
\begin{align*}
J_p(y_1,...,y_n) = \intt dt \; j_p (t, y_1,...,y_n,y_1',...,y_n'),  
\end{align*}
we solve the Euler-Lagrange equation for the functional 
\begin{align*}
K(y_1,...,y_n) = T(y_1,...,y_n) + \mathlarger{\sum}_p \;  \lambda_p \; (J_p(y_1,...,y_n) - c_p).  
\end{align*}
The Lagrangian corresponding to $K$ is 
\begin{align*}
L^* = L + \mathlarger{\sum}_p  \; \lambda_p  \; (j_p - \frac{c_p}{T}),  
\end{align*}
where $T=t_1-t_0$. The $2n$ integration constants, together with the $m$ Lagrange multipliers $\{ \lambda_p \}^m_{p=1}$,  are determined by the conditions 
\begin{align*}
&\begin{rcases}
y_i(t_0) = y_i^0 \\
y_i(t_1) = y^1_i 
\end{rcases}i = 1,...,n\;\;,  \\ 
\nonumber \\
& J_p(y_1,...,y_n) = c_p \qquad \qquad \qquad \qquad p = 1,...,m\;\;. \nonumber  
\end{align*}

\begin{example}\label{Example6}
Let us consider the isoperimetric problem from section \ref{VariationalExample3}. The solution to this problem was there reduced to maximizing the functional 
\begin{align*}
2 \; T = \int^1_0 dt \; (x(t) \; y'(t) - y(t) \; x'(t)),
\end{align*} 
subject to the constraint 
\begin{align*}
J = \int^1_0 dt \; \sqrt{x'(t)^2 + y'(t)^2} = L^*. 
\end{align*}
We introduce the Lagrangian
\begin{align*}
L^* = x\; y' - y \; x' + \lambda \; (\sqrt{x'^2 + y'^2} - L^*), 
\end{align*}
and observe that 
\begin{align*}
\frac{\partial L^*}{\partial x} &= y',\;\;\frac{\partial L^*}{\partial y} = - x',  \\ 
\frac{\partial L^*}{\partial x'} &= -y + \frac{\lambda \; x'}{\sqrt{x'^2 + y'^2}}, \nonumber \\
\frac{\partial L^*}{\partial y'} &= x + \frac{\lambda \; y'}{\sqrt{x'^2 + y'^2}}. \nonumber 
\end{align*}
The Euler-Lagrange equations are then 
\begin{align}
\frac{\partial L^*}{\partial x} - \dt{} (\frac{\partial L^*}{\partial x'}) &= 0, \nonumber \\
&\Updownarrow \nonumber \\  y ' - \dt{} (-y + \frac{\lambda \; x'}{\sqrt{x'^2 + y'^2}}) &= 0, \lbl{180.12}
\end{align}
and 
\begin{align}
\frac{\partial L^*}{\partial y} - \dt{} (\frac{\partial L^*}{\partial y'}) &= 0, \nonumber \\
&\Updownarrow \nonumber \\  - x' - \dt{} (x + \frac{\lambda \; y'}{\sqrt{x'^2 + y'^2}}) &= 0. \lbl{181.12} 
\end{align}
Equations \rf{180.12} and \rf{181.12} can be integrated once to yield 
\begin{align}
2 \; y - \frac{\lambda \; x'}{\sqrt{x'^2 + y'^2}} &= c_1, \nonumber \\ 
2 \; x + \frac{\lambda \; y'}{\sqrt{x'^2 + y'^2}} &= c_2, \nonumber \\ 
&\Updownarrow \nonumber \\ y - \frac{c_1}{2} &= \inv{2} \; \frac{\lambda \; x'}{\sqrt{x'^2 + y'^2}}, \lbl{183.12} \\ 
x - \frac{c_2}{2} &= - \inv{2} \; \frac{\lambda \; y'}{\sqrt{x'^2 + y'^2}}. \nonumber 
\end{align}
Squaring and adding the two equations \rf{183.12} we get 
\begin{align}
(y - \frac{c_1}{2})^2 + (x- \frac{c_2}{2})^2 = \inv{4} \; \lambda^2, \lbl{184.12} 
\end{align}
which we recognize to be the equation for a circle. Thus all extremals are circles. The radius of the circle \rf{184.12} is $R = \inv{2} \; \lambda$. Thus the constraint is satisfied if 
\begin{align*}
2 \; \pi \; R &= L^*, \nonumber \\  &\Updownarrow \nonumber \\ \lambda &= \frac{L^*}{\pi}.  
\end{align*}

\end{example}\label{Example6}
\subsubsection{Several independent variables}
Consider a functional of the form 
\begin{align*}
T(u) = \int_D dx \; dy \; \mathcal{L}(x,y,u,u_x,u_y).  
\end{align*}
The challenge is to find the stationary points for $T$ subject to the constraint 
\begin{align*}
u\mid_{\partial D} = f.  
\end{align*}
We proceed like before by introducing a variation 
\begin{align*}
v = u + \eps \; \eta, && \eta = \eta (x,y).  
\end{align*}
Observe that since the boundary condition has to be fixed we have 
\begin{align}
v\mid_{\partial D} = u\mid_{\partial D} + \eps \; \eta\mid_{\partial D} = f \; \Rightarrow \; \eta\mid_{\partial D} = 0. \lbl{189.12} 
\end{align}
For the functional $T$ we now have 
\begin{align*}
T (u + \eps \; \eta) &= \int_D dx\; dy \; \mathcal{L}(x,y,u+\eps \; \eta, u_x + \eps \; \eta_x, u_y + \eps \; \eta_y)  \\ 
&= \int_D dx \; dy \; \{ \mathcal{L}(x,y,u,u_x,u_y) + \frac{\partial \mathcal{L}}{\partial u} \; \eps \; \eta + \frac{\partial \mathcal{L}}{\partial u_x} \; \eps \; \eta_x \nonumber \\
&+ \frac{\partial \mathcal{L}}{\partial u_y} \; \eps \; \eta_y  \} + \obs(\eps^2) \nonumber \\
&= T (u) + \eps \; \int_D dx \; dy \; \{ \frac{\partial \mathcal{L}}{\partial u} \; \eta + \frac{\partial \mathcal{L}}{\partial u_x} \; \eta_x + \frac{\partial \mathcal{L}}{\partial u_y} \; \eta_y \} \nonumber \\ 
&+ \obs(\eps^2) \nonumber \\
&= \vT(u) + \eps \; \int_D dx \; dy \; \{ \frac{\partial \mathcal{L}}{\partial u} - \prt{x} (\frac{\partial \mathcal{L}}{\partial u_x}) - \partial_y (\frac{\partial \mathcal{L}}{\partial u_y}) \} \; \eta \nonumber \\
&+ \obs(\eps^2), \nonumber 
\end{align*}
where we in the last step have used Green's theorem in divergence form and the boundary condition \rf{189.12} on $\eta$. 

Using the fundamental lemma we conclude that $u$ is a stationary point for $T$ if it satisfies the following Euler-Lagrange equation. 
\begin{align*}
\frac{\partial \mathcal{L}}{\partial u} - \prt{x} (\frac{\partial \mathcal{L}}{\partial u_x}) - \prt{y} (\frac{\partial \mathcal{L}}{\partial u_y}) = 0.  
\end{align*}

\begin{example}\label{Example7}
In section \ref{VariationalExample5} we introduced the notion of a minimal surface. We consider here the simplified situation where the surface is the graph of a function $u(x,y)$ over a domain $D$ in the plane.
The challenge is to minimize the functional 
\begin{align*}
T (u) = \int_D dx \; dy \; \sqrt{1+u_x^2 + u_y^2},  
\end{align*}
subject to the constraint 
\begin{align*}
u\mid_{\partial D} = h(x,y).  
\end{align*}
The Lagrangian density is 
\begin{align*}
\mathcal{L} = \sqrt{1+u_x^2 + u_y^2},  
\end{align*}
and we have 
\begin{align*}
\frac{\partial \mathcal{L}}{\partial u} &= 0,\nonumber\\
 \frac{\partial \mathcal{L}}{\partial u_x} &=\frac{ u_x}{\sqrt{1 + u_x^2 + u_y^2}},\nonumber \\
\frac{\partial \mathcal{L}}{\partial u_y} &=\frac{ u_y}{\sqrt{1 + u_x^2 + u_y^2}}. 
\end{align*}
The Euler-Lagrange equation is thus 
\begin{align*}
- \prt{x} \frac{u_x}{\sqrt{1 + u_x^2 + u_y^2}} - \prt{y}\frac{u_y}{\sqrt{1+u_x^2 + u_y^2}} = 0,  
\end{align*}
which can be rewritten as 
\begin{align}
u_{xx} + u_{yy} =2 u_xu_yu_{xy}-u_x^2u_{yy}-u_y^2u_{xx} \lbl{197.12} 
\end{align}
The boundary condition for $u$ is 
\begin{align*}
u\mid_{\partial D} = h.  
\end{align*}
The equation \rf{197.12} is a non-linear second order partial differential equation and is not by any means easy to solve in general. 

However,  if the boundary curve is horizontal
\begin{align*}
h(x,y) = h_0,  
\end{align*} 
the boundary value problem clearly has the unique solution 
\begin{align*}
u(x,y) = h_0, && (x,y) \in D.  
\end{align*}
This is flat and is obviously of minimum area among all surfaces with flat boundary curve $h=h_0$. If the boundary curve is not constant but varies little on the scale of $h_0$ 
\begin{align*}
h = h_0 + \eps \; k && \eps << 1,  
\end{align*}
we seek a solution that is a small modification of $u(x,y) = h_0$
\begin{align*}
u(x,y) = h_0 + \eps \; v && \eps << 1.  
\end{align*}
The function $v$ must then satisfy the equation
\begin{align*}
v_{xx} + v_{yy} = \eps^2 (2 v_xv_yv_{xy}-v_x^2v_{yy}-v_y^2v_{xx}).    
\end{align*}
This equation we can solve approximately using a perturbation expansion. Expansions of this type will be discussed in section five of these lecture notes.
\end{example}

\begin{example}\label{Example8}
In section \ref{VariationalExample9}  we introduced the maximum entropy principle. The aim of this section is to derive the maximum entropy distribution by solving the corresponding Euler-Lagrange equation. It turns out that special cases of the resulting probability distribution form the foundation for statistical mechanics and thermodynamics, information theory and probably also elementary particle physics through it’s  mathematical grounding in quantum field theory. The current section is an excerpt of a more detailed treatment of the maximum entropy principle given in Appendix B.

As you recall, the maximum entropy principle states that,  if what we known about  a system $S$, prior to a measurement, is described by a probability distribution $\rho_0$, and we measure the mean $c_j$ of $p$ observables $f_j$ for $S$, then the probability distribution that best incorporates this new  information about the system, is the one that maximizes the functional 
\begin{align}
S(\rho) = - \int_{\mathbf{R}^n} dV \; \rho \; \ln (\frac{\rho}{\rho_0}), \lbl{204.12} 
\end{align}
under the constraints 
\begin{align*}
\left<f_j\right> \; = \int_{\mathbf{R}^n} dV \; f (x_1, ... , x_n) \; \rho(x_1,...,x_n) = c_j.  
\end{align*}
Since $\rho$ must be a probability distribution we have one more constraint 
\begin{align*}
\left<1\right> \; = \int_{\mathbf{R}^n} dV \; \rho (x_1,...,x_n) = 1,  
\end{align*}
and we thus have $p+1$ constraints and therefore introduce an extended functional 
\begin{align*}
T(\rho) = S(\rho) + \lambda_0 \; \left<1\right> + \mathlarger{\sum}^p_{j=1} \lambda_j \; \left<f_j\right>  
\end{align*}
Note that we could have written 
\begin{align*}
T (\rho) = S(\rho) + \lambda_0 \; ( \left<1\right> - 1) + \mathlarger{\sum}^p_{j=1} \lambda_j \; ( \left<f_j\right> - c_j),  
\end{align*}
in order to make the values of the constraints explicit, like we have done on previous occasions. However, all constant terms vanish when we take variational derivative, so we might as well drop the constant terms. 

The Lagrangian  density corresponding to the extended functional $T (\rho)$ is then
\begin{align*}
\mathcal{L} = - \rho \; \ln (\frac{\rho}{\rho_0}) + \lambda_0 \; \rho + \mathlarger{\sum}^p_{j=1} \lambda_j \; f_j \; \rho  
\end{align*}
Observe that $\mathcal{L}$ does not depend on any derivatives of $\rho$. The Euler-Lagrange equation for $T$ is therefore simply 
\begin{align*}
\frac{\partial \mathcal{L}}{\partial \rho} &= 0,  \\ 
&\Updownarrow \nonumber \\
& - \ln (\frac{\rho}{\rho_0}) - 1 + \lambda_0 + \mathlarger{\sum}^p_{j=1} \lambda_j f_j , \nonumber \\ 
&\Updownarrow \nonumber \\ \rho &= \frac{\rho_0}{Z} \; \exp(\mathlarger{\sum}_j \lambda_j \; f_j), \nonumber 
\end{align*}
where we have defined $Z=\exp(1- \lambda_0 )$. In order for the constraint $\left<1\right> = 1$ to be satisfied, we must have 
\begin{align*}
\left<1\right> &= 1, \nonumber \\
&\Updownarrow \nonumber \\  \intRn dV \; \frac{\rho_0}{Z} \; \exp(\mathlarger{\sum}_j \lambda_j \; f_j) &= 1, \nonumber \\ 
&\Updownarrow \nonumber \\ Z = Z(\lambda_1,...,\lambda_p) &= \intRn dV \; \rho_0 \; \exp{\mathlarger{\sum}_j \; \lambda_j \; f_j},  
\end{align*}
and the stationary distribution is 
\begin{align}
\rho(x_1,...,x_n) = \frac{\rho_0(x_1,...,x_n)}{Z(\lambda_1,...,\lambda_p)} \; \exp(\mathlarger{\sum}^p_{j=1} \lambda_j \; f_j (x_1, ...,x_n)). \lbl{212.12} 
\end{align}
$\rho$ is called the \ttx{maximum entropy distribution} and $Z$ is the $partition function$. Note that we have not proved that it in fact is the distribution that gives a maximum for $S$, but this can be done. \\
The Lagrange multipliers $\lambda_1,...,\lambda_p$ are chosen so that all the constraints are satisfied 
\begin{align}
\left<f_j\right> \; = \intRn dV \; f_j (x_1,...,x_n) \; \rho(x_1,...,x_n) = c_j && j=1,...,p\;. \lbl{213.12} 
\end{align}
The system of equations \rf{213.12} consists of  $p$ equations for the $p$ quantities $\lambda_j$. As it turns out, we almost never need to know the distribution $\rho$ from \rf{212.12}, it is enough to know the partition function. Observe that 
\begin{align}
\left<f_j\right> \; &= \intRn dV \; f_j \; \rho \lbl{214.12} \\ 
&= \inv{Z} \; \intRn dV \; f_j \; \rho_0 \; \exp(\mathlarger{\sum}^p_{i=1} \lambda_{i} f_{i}) \nonumber \\ 
&= \inv{Z} \; \intRn dV \; \partial_{\lambda_j}\{ \rho_0 \; \exp(\mathlarger{\sum}^p_{i=1} \lambda_{i} f_{i})\} \nonumber \\ 
&= \inv{Z} \; \partial_{\lambda_j}  \intRn dV \;  \rho_0 \; \exp(\mathlarger{\sum}^p_{i=1} \lambda_{i} f_{i}) \nonumber \\
&= \inv{Z} \; \prt{\lambda_j} Z = \prt{\lambda_j} \ln Z, \nonumber 
\end{align}
and thus we can find the mean of all the quantities $f_j$ by taking partial derivatives of the partition function with respect to the Lagrangian multipliers. Moreover, we also have 
\begin{align*}
\prt{\lambda_j \lambda_k}\ln Z &= \prt{\lambda_j} (\inv{Z} \; \prt{\lambda_k} Z)  \\
&= - \inv{Z^2} \; \prt{\lambda_j} Z \; \prt{\lambda_k} Z + \inv{Z} \; \prt{\lambda_j \lambda_k} Z \nonumber \\
&= - \prt{\lambda_j} \ln Z \; \prt{\lambda_k} \ln Z + \inv{Z} \; \intRn dV \; f_j \; f_k \; \rho_0 \; \exp(\mathlarger{\sum}_i \lambda_i \; f_i) \nonumber \\
&= - \prt{\lambda_j}\ln Z \; \prt{\lambda_k} \ln Z + \left< f_j \; f_k\right>. \nonumber 
\end{align*}
Thus 
\begin{align*}
\left<f_j \; f_k\right> \; = \prt{\lambda_j} \ln Z \; \prt{\lambda_k} \; \ln Z + \prt{\lambda_j \; \lambda_k} \ln Z  
\end{align*}
\end{example}
In a similar way \ttx{all} correlation coefficients $\left<f_1^{n_1} ... f_p^{n_p}\right>$ can be expressed through derivatives of the partition function.\\ Let us consider the special case when our underlying space is the classical state space for a mechanical system with $n$ degrees of freedom. This could for example consist of $n$ mass points. We thus have a state space $\mathbf{R}^{6n}$ since we need 3 position coordinates $\vb{x}=(x_1,x_2,x_3)$,  and 3 velocity coordinates $\vb{v}=(v_1,v_2,v_3)$, or equivalently three momentum coordinates $\vb{p}=m\vb{v}$, for each particle in order to uniquely specify the state of the system. \\
Let $H = H(\vx_1,...,\vx_n, \vb{p}_1,...,\vb{p}_n)$ be the Hamiltonian for the system of mass   points. Recall that the value of the Hamiltonian on any given state  \\$(\vx_1,...,\vx_n,\vb{p}_1,...,\vb{p}_n)$, is the energy of that state. \\
When $n$ is large it is very hard, and also mostly useless, to try to track the exact state $(\vx_1(t),...,\vx_n(t), \vb{p}_1(t),...,\vb{p}_n(t))$ of a system of mass points. \\
For such a large system it is more useful to consider a probability distribution $\rho(\vx_1,...,\vx_n , \vb{p}_1,...,\vb{p}_n)$ on the state-space. We have seen how useful this point of view is in fluid dynamics. \\
Here we will assume that we have some prior distribution $\rho_0$ and the observation of the mean value of the total energy, $H$, of the system 
\begin{align*}
\left<H\right> \; = E.  
\end{align*}
Using the maximal entropy principle we are lead to select the following probability distribution 
\begin{align*}
\rho (\vx_1,...,\vx_n, \vb{p}_1,..., \vb{p}_n) = \frac{\rho_0}{Z} \; \exp(\lambda \; H).  
\end{align*}
In this situation one usually redefines $\lambda$ by writing 
\begin{align*}
\lambda = - \inv{k \; T},  
\end{align*}
where $k$ is the Boltzmann constant and $T$ is a new parameter. \\
What we then we get is the well known \ttx{Gibb's ensemble} from statistical physics 
\begin{align}
\rho = \frac{\rho_0}{Z(T)} \; \exp{- \inv{k T} H}. \lbl{219.12} 
\end{align}
The parameter $T$ is determined by 
\begin{align}
E&=\left<H\right>, \lbl{220.12} \\
&\Updownarrow \nonumber \\
E&= k \; T^2 \; \prt{T} \ln Z , \nonumber 
\end{align}
where we have used the chain rule 
\begin{align*}
\prt{\lambda} = k \; T^2 \; \prt{T},  
\end{align*}
in the general formula \rf{214.12}. \\
The Gibb's ensemble is the foundation of statistical physics. All results in statistical physics flows from formula \rf{219.12}. Statistical physics is also the foundation of thermodynamics so all conclusions from that subject also flow from formula \rf{219.12}. In this context, $T$ is the temperature of the system of mass-points and \rf{220.12} is nothing but the \ttx{equation of state}. 

An interesting insight here is that the temperature of a thermodynamic system is in fact a Lagrange multiplier!! This is a profound insight that to this day has not been fully understood or explored.

An extended discussion of the maximum entropy principle and how it relates to foundational problems in statistical physics  is included in appendix B.

From this example, it appears as if it might be useful to think of any application of the maximal entropy principle as an extension of the methods of statistical mechanics to systems that has absolutely nothing to do with the motion of mass points. 

This wide general applicability of the methods of statistical physics has also lead to deep questions and insights into the nature and significance of the assumption of equilibrium that appears to underline the application of the Gibb's ensemble in statistical physics. 

There is also the intriguing fact that the very same functional \rf{204.12} used in the maximum entropy principle,  is also the foundation of information theory which was discovered by Shannon in 1948. This connection between information theory and statistical mechanics (and thermodynamics) has lead to deep insights into the role of information in our fundamental physical theories. \\
The general nature and wide applicability of the maximum entropy principle has been described well by E.T. Jaynes in many papers and the monumental book "Probability theory: The Logic of Science". 

As if all this is not impressive enough for one single principle, it is also a very intriguing fact that when one looks deep into the heart of fundamental physics, in the form of quantum field theory, one again finds an appropriately generalized form of the Gibb's ensemble! The whole computational engine in the theory of quantum fields revolve around this generalized Gibb's ensemble. 

What on earth is going on...
\\
\subsection{Equations of Variational Type}

We have seen that stationary points for  functionals are solutions to the Euler-Lagrange equations corresponding to the functional. The exact structure of the Euler-Lagrange equations and their number depends on the functional. We have seen several examples of differential equations of the Euler-Lagrange type in the previous sections of these lecture notes. For example have we found that the differential equations \rf{101.12}, \rf{136.12} and \rf{197.12} are of Euler-Lagrange type.

\noindent In this section we will ask which (systems of) differential equations are Euler-Lagrange equations for some functional. This is an important question to ask, because many important structural properties of differential equations can be decided if we know that they are Euler-Lagrange equations for some functional. 

Equations that are Euler-Lagrange equations for a functional are said to be \ttx{variational}.
\subsubsection{Real valued functions}
\begin{example}\label{Example9}
Let $T(u)$ be the functional 
\begin{align*}
T(u) = \iint_D dx \; dy \; ( \inv{2} \; u_x^2 + \inv{2} \; u_y^2),  
\end{align*}
with a constraint 
\begin{align*}
u\mid{\partial D} = f.  
\end{align*}
The Lagrangian density is 
\begin{align*}
\mathcal{L} = \inv{2} \; u_x^2 + \inv{2} \; u_y^2,  
\end{align*}
and we have 
\begin{align*}
\frac{\partial \mathcal{L}}{\partial u} = 0, \qquad \qquad \qquad \frac{\partial \mathcal{L}}{\partial u_x} = u_x, \qquad \qquad \qquad \frac{\partial \mathcal{L}}{\partial u_y}= u_y,  
\end{align*}
so the Euler-Lagrange equation is 
\begin{align*}
\frac{\partial \mathcal{L}}{\partial u} - \prt{x} (\frac{\partial \mathcal{L}}{\partial u_x}) - \prt{y} (\frac{\partial \mathcal{L}}{\partial u_y}) &= 0,  \\ 
&\Updownarrow \nonumber \\  u_{xx} + u_{yy} &= 0. \nonumber 
\end{align*}
\end{example}
\noindent Thus the 2D Laplace equation is variational. The same is true in 3D or in any dimension for that matter. The same is also true for Poisson's equation in any dimension.

\begin{example}\label{Example10}
Consider the functional 
\begin{align*}
T (u) = \int_D dx \; dy \; \intt dt \; \{ \inv{2} \; u_t^2 - \inv{2} \; c^2 \; u_x^2 - \inv{2} \; c^2 \; u_y^2.  \}  
\end{align*}
We are looking for functions $u(x,y, t)$ that are stationary points with respect to variations that vanish on the boundaries to the domain of integration 
\begin{figure}[htbp]
\centering
\includegraphics{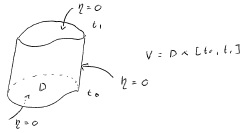}
\caption{The integration domain of the functional whose Euler-Lagrange equation is the wave equation.}
\label{fig14}
\end{figure}
\begin{align*}
T (u + \eps \; \eta) &= \int_D dx \; dy \; \intt dt \; \{ \inv{2} \; (u_t + \eps \; \eta_t)^2 - \inv{2} \; c^2 \; (u_x + \eps \; \eta_x)^2 \nonumber \\ 
&- \inv{2} \; c^2 \; (u_y + \eps \; \eta_y)^2 \} \nonumber \\
&= \int_D dx \; dy \; \intt dt \; \{ \inv{2} \; u_t^2 - \inv{2} \; c^2 \; u_x^2 - \inv{2} \; c^2 \; u_y^2 \nonumber \\
&+ \eps \; u_t \; \eta_t - c^2 \; \eps \; u_x \; \eta_x - c^2 \; \eps u_y \; \eta_y \} + \obs(\eps) \nonumber \\
&= T(u) + \eps \; \int_D dx \; dy \; \intt dt \; \{ u_t \; \eta_t - c^2 \; u_x \; \eta_x - c^2 \; u_y \; \eta_y \} \nonumber \\ &+ \obs(\eps^2) \nonumber\\
&= T(u) + \eps \; [ \int_D dx \; dy \; \intt dt \; ((u_t \; \eta)_t - u_{tt} \; \eta) \nonumber \\
&- c^2\intt dt \; \int_D dx \; dy \; (\div(\grad{u} \eta) - \laplacian{u \; \eta}) ] + \obs(\eps^2) \nonumber 
\end{align*}
\begin{align*}
&\Downarrow\\
T (u + \eps \; \eta)&= T(u) + \eps [ \int_D dx \; dy \; (u_t \; \eta)\mathlarger{|}^{t_1}_{t_0} - \iint_D dx \; dy \; \intt dt \; u_{tt} \; \eta \nonumber \\ 
&-c^2 \; \intt dt \; \int_{\partial D} dl \; \eta \; \grad{u} \vdot \vb{n} + c^2 \; \intt dt \; \int_D dx \; dy \; \laplacian{u \; \eta} ] \nonumber \\ &+ \obs(\eps^2) \nonumber\\ 
&= T(u) + \eps \; \int_D dx \; dy \; \intt dt \; \{ - u_{tt} + c^2 \; \laplacian{u} \} \; \eta + \obs(\eps^2),
\end{align*}
where we in the last step have used the boundary conditions for $\eta$. Since $\eta(x,y,t)$ is arbitrary within the domain $D\cross[t_0,t_1]$ we can use the fundamental lemma to conclude that $u$ is a stationary point for $T$ only if 
\begin{align*}
u_{tt} - c^2 \; \laplacian{u} = 0.  
\end{align*}
\end{example}
\noindent This proves that the 2D wave equation is variational. The same is true for wave equations on any number of spatial dimensions. 

\begin{example}\label{Example11}
Consider the functional 
\begin{align*}
T(u) &= \int_D dx \; dy \; \intt dt \; \{\inv{2} \; u_t^2 - \inv{2} \; c^2 \; u_x^2 - \inv{2} \; c^2 \; u_y^2 - \inv{2} \;m\; u^2 \} && m>0  
\end{align*}
Using the exact same approach as in example \ref{Example10}, it is easy to show that $u$ is extremal for variation $\eta$ vanishing on the boundary of $D\cross[t_0, t_1]$ only if 
\begin{align*}
u_{tt} -c^2 \laplacian{u} + m \; u = 0 
\end{align*}
This is the (real) Klein-Gordon equation. Thus the Klein-Gordon equation is variational.
\end{example}
\noindent This is all fine, you might say, but where did the functionals in examples \ref{Example9}-\ref{Example11} come from?  
Is there some practical, useful, general method for deciding if a given (system of) differential equations is variational and find an appropriate functional if they are variational? 

Sadly, no such general method exists. The functionals must be constructed by trial and error.

\subsubsection{Complex valued functions}\label{ComplexValuedFields}

Many important equations in theoretical physics and applied mathematics involve complex valued functions. The prime example is of course the quantum mechanical Schr$\ddot{o}$dinger equation. 

In order to decide whether such equations are variational or not, we need to extend the calculus of variation to functions whose values are complex. 
\begin{example}
Let us consider the complex second order differential equation
\begin{align}
A_{tt} = \abs{A}^2 \; A && A(t) \in \mathbf{C}. \lbl{232.12} 
\end{align}
Introduce the real and imaginary parts of $A$ 
\begin{align}
A(t) = u(t) + i \; v(t). \lbl{233.12} 
\end{align}
Inserting \rf{233.12} into \rf{232.12}, we get 
\begin{align*}
u_{tt} + i \; v_{tt} = (u^2 + v^2) \; (u + i \; v). 
\end{align*}
Separating real and imaginary parts we get two coupled real equations 
\begin{align}
u_{tt} &= (u^2 + v^2) \; u, \lbl{235.12} \\ 
v_{tt} &= (u^2 + v^2) \; v. \nonumber 
\end{align}
Consider the functional 
\begin{align}
T (u,v) = \intt dt \; (\inv{2} \; u_t^2 + \inv{2} \; v_t^2 + \inv{4} \; u^4 + \inv{2} \; u^2 \; v^2 + \inv{4} \; v^4). \lbl{236.12} 
\end{align}
The Lagrangian is 
\begin{align}
L = \inv{2} \; u_t^2 + \inv{2} \; v_t^2 + \inv{4} \; u^4 + \inv{2} \; u^2 \; v^2 + \inv{4} \; v^2, \lbl{237.12} 
\end{align}
and the corresponding Euler-Lagrange equations are 
\begin{align*}
\frac{\partial L}{\partial u} - \dt{} (\frac{\partial L}{\partial u_t}) &= 0,  \\
&\Updownarrow \nonumber \\ u^3 + v^2 \; u - u_{tt} &= 0, \nonumber \\
&\Updownarrow \nonumber \\  u_{tt} &= (u^2 + v^2) \; u, \nonumber \\
\frac{\partial L}{\partial v} - \dt{} (\frac{\partial L}{\partial v_t}) &= 0,  \\ 
&\Updownarrow \nonumber \\  v^3 + u^2 \; v - v_{tt} &= 0,  \nonumber \\
&\Updownarrow \nonumber \\ v_{tt} &= (u^2 + v^2) \; v. \nonumber 
\end{align*}
This shows that \rf{235.12}  is variational with the corresponding functional defined in \rf{236.12}. We could now {\it define} a complex equation like \rf{232.12} to be variational if the system of real equations we get when separating into real and imaginary parts is variational. It would however be better if we could define what it means for a complex equation to be variational {\it without} first separating it  into real and imaginary parts.

Observe that the Lagrangian \rf{237.12} can be written more compactly as 
\begin{align}
L = \inv{2} \; u_t^2 + \inv{2} \; v_t^2 + \inv{4} \; (u^2 + v^2)^2.  \lbl{240.12}
\end{align}
We have 
\begin{align}
& u + i \; v = A, && u - i \; v = A^*,\nonumber  \\ 
& u = \inv{2} \; (A + A^*), && v = \inv{2 i} \; (A - A^*). \lbl{241.12}
\end{align}
Inserting \rf{241.12} into the Lagrangian  \rf{240.12}, we get 
\begin{align*}
L &= \inv{8} \; (A_t + A_t^*)^2 - \inv{8} \; (A_t - A_t^*)^2  \\
&+ \inv{4} \; (A \; A^*)^2 \nonumber \\
&= \inv{8} \; A_t^2 + \inv{4} \; A_t \; A_t^* +\inv{8} \; A_t^{*2}  - \inv{8} \; A_t^2 + \inv{4} \; A_t \; A_t^* \nonumber \\
&- \inv{8} \; A_t^{*2} + \inv{4} \; A^2 \; A^{*2} \nonumber \\
&= \inv{2} \; A_t \; A_t^* + \inv{4} \; A^2 \; A^{*2}. \nonumber
\end{align*}
Since a factor of 2 makes no difference for the Euler-Lagrange equation, we might as well use the Lagrangian 
\begin{align}
L = A_t \; A_t^* + \inv{2} \; A^2 \; A^{*2}. \lbl{243.12} 
\end{align}
In this Lagrangian density, $A$ is the only dependent variable. $A^*$ is of course calculated by taking the complex conjugate of $A$. However let us disregard this fact and assume that $A$ and $A^*$ can be varied \ttx{independently}. \\
Then, any functional of the form 
\begin{align*}
T(A,A^*) = \intt dt \; L (A,A^*,A_t,A_t^*),   
\end{align*}
for some Lagrangian  $L$, will in the usual way lead to Euler-Lagrange equations 
\begin{align*}
\frac{\partial L}{\partial A^*} - \dt{} (\frac{\partial L}{\partial A_t^*}) &= 0,  \\ 
\frac{\partial L}{\partial A} - \dt{} (\frac{\partial L}{\partial A_t}) &= 0. \nonumber 
\end{align*}
For our particular Lagrangian \rf{243.12} we have 
\begin{align*}
\frac{\partial L}{\partial A} &= A\; A^{*2}, && \frac{\partial L}{\partial A_t} = A_t^*,  \\ 
\frac{\partial L}{\partial A^*} &= A^2 \; A^*, && \frac{\partial L}{\partial A_t^*} = A_t. \nonumber
\end{align*}
The Euler-lagrange equations are thus 
\begin{align*}
A \; A^{*2} - A_{tt}^* &= 0,  \\ 
&\Updownarrow \nonumber \\ A_{tt}^* &= \abs{A}^2 \; A^*, \nonumber \\ 
A^2 \; A^* - A_{tt} &= 0,  \\
&\Updownarrow \nonumber \\  A_{tt} &= \abs{A}^2 \; A, \nonumber 
\end{align*}
which are our original equation \rf{232.12} and its complex conjugate. \\
Thus following this formal procedure, where we assume that  $A$ and $A^*$ can be varied independently, we have again proved that \rf{232.12} is variational, and we have done this without separating the problem into real and imaginary parts. 

This procedure will always work and its conclusions are equivalent to what we get by separating into real and imaginary part. \\
However, an important caveat is that we must use this formal procedure only Lagrangians that are \ttx{real valued}. We must thus always make sure that 
\begin{align*}
L^* = L  
\end{align*}
\end{example}

\begin{example}
Consider the functional 
\begin{align*}
T (A,A^*) = \int_D dx \; dy \; \intt dt \; \mathcal{L}(A,A^*, A_t, A_t^*, A_x,A_x^*, A_y, A_y^*),  
\end{align*}
where $\mathcal{L} = \mathcal{L}^*$. Following the procedure from example \ref{Example10} with dependent variables $A$ and $A^*$ we get the Euler-Lagrange equations 
\begin{align}
\frac{\partial \mathcal{L}}{\partial A} - \prt{t} (\frac{\partial \mathcal{L}}{\partial A_t} ) - \prt{x} (\frac{\partial \mathcal{L}}{\partial A_x} ) - \prt{y} (\frac{\partial \mathcal{L}}{\partial A_y} ) &= 0, \nonumber \\ 
\frac{\partial \mathcal{L}}{\partial A^*} - \prt{t} (\frac{\partial \mathcal{L}}{\partial A^*_t} ) - \prt{x} (\frac{\partial \mathcal{L}}{\partial A^*_x} ) - \prt{y} (\frac{\partial \mathcal{L}}{\partial A^*_y} ) &= 0. \lbl{251.12}
\end{align}
As a matter of fact, we only need one of the equations \rf{251.12} since the first one is just the complex conjugate of the second one. 

Her we choose to use the second equation. Consider the special Lagrangian density
\begin{align}
\mathcal{L} = A_t \; A_t^* - c^2 \; A_x \; A_x^* - c^2 \; A_y \; A_y^* - m \; A \; A^*. \lbl{252.12} 
\end{align}
We have 
\begin{align*}
&\frac{\partial \mathcal{L}}{\partial A^*} = - m \; A, && \frac{\partial \mathcal{L}}{\partial A_t^*} = A_t, \nonumber \\ &\frac{\partial \mathcal{L}}{\partial A_x^*} = - c^2 \; A_x, && \frac{\partial \mathcal{L}}{\partial A_y^*} = - c^2 \; A_y, 
\end{align*}
and the Euler -Lagrange equation is 
\begin{align*}
- m \; A - A_{tt} + c^2 \; A_{xx} + c^2 \; A_{yy} &= 0,  \\
&\Updownarrow \nonumber \\  A_{tt} - c^2 \; \laplacian{A} + m \; A &= 0. \nonumber 
\end{align*}
This shows that the 2D complex Klein-Gordon equation is variational with Lagrangian density defined in  \rf{252.12}. The same is true for the 1D and 3D cases. 
\end{example}

\begin{example}\label{SchrodingerEquation}
Consider a Lagrangian density
\begin{align}
\mathcal{L} &=   i\;\frac{\hbar}{2} \; (\psi^* \;\psi_t - \psi \; \psi^*_t) - \frac{\hbar^2}{2m} \; (\psi_x \; \psi_x^* + \psi_y \; \psi_y^*) - V(x,y) \; \psi \; \psi^*. \lbl{255.12} 
\end{align}
We evidently have $\mathcal{L}= \mathcal{L}^*$. Observe that 
\begin{align*}
\frac{\partial \mathcal{L}}{\partial \psi^*} &=i\; \frac{\hbar}{2} \; \psi_t - V \; \psi, &&  &\frac{\partial \mathcal{L}}{\partial \psi_t^*} = - i\;\frac{\hbar}{2}\psi,\nonumber \\ 
\frac{\partial \mathcal{L}}{\partial \psi_x^*}&= - \frac{\hbar^2}{2m} \; \psi_x, && &\frac{\partial \mathcal{L}}{\partial \psi_y^*} = - \frac{\hbar^2}{2m} \; \psi_y,
\end{align*}
and the Euler-Lagrange equation is 
\begin{align*}
i\;\frac{\hbar}{2} \; \prt{t} \psi - V \; \psi - \prt{t} (- \frac{\hbar}{2} \; i \; \psi )  \\ 
- \prt{x} ( - \frac{\hbar^2}{2m} \; \psi_x) - \prt{y} (- \frac{\hbar^2}{2m} \; \psi_y) &= 0, \nonumber \\
&\Updownarrow \nonumber \\  i \; \hbar \; \prt{t} \psi &= - \frac{\hbar^2}{2m} \; \laplacian{\psi} + V \; \psi. \nonumber
\end{align*}
This shows that the Sch$\ddot{o}$dinge equation is variational with Lagrangian density defined in  \rf{255.12}. It is easy to show that the stationary Sch$\ddot{o}$dinger equation is also variational.
\\
\subsection{Noether's Theorem}

Noether's theorem creates a one-to-one correspondence between \ttx{conserved quantities} of variational equations and \ttx{symmetries} of the corresponding functionals. 

The theorem was proved by the mathematician Emmy Noether in 1915. It has been described as: \\

"One of the most important mathematical theorems ever proved in guiding the development of modern physics"\\
\subsubsection{One dependent variable}
\noindent Let us introduce the theorem in the simplest possible context. We consider a functional of the form
\begin{align}
T (y) = \intt dt \; L(t,y,y'). \lbl{258.121}
\end{align}
Let us consider some variation
\begin{align}
y \rightarrow y + \eps \; \eta, \lbl{259.121} 
\end{align}
where $\eta=\eta(t)$ now is some \ttx{specific} function. Inserting \rf{259.121} into \rf{258.121} we get 
\begin{align*}
T(y + \eps \; \eta) &= \intt dt \; L(t,y+\eps \; \eta, y' + \eps \; \eta')  \\
&= \intt dt \; \{ L(t,y,y') + \frac{\partial L}{\partial y} \; \eps \; \eta + \frac{\partial L}{\partial y'} \; \eps \; \eta' \} + \obs(\eps^2) \nonumber \\ 
&= T(y) + \eps \; \intt dt \; \{ \frac{\partial L}{\partial y} \; \eta + \frac{\partial L}{\partial y'}\eta' \} + \obs(\eps^2). \nonumber 
\end{align*}
We now introduce the key idea of invariance. 

\noindent The functional $T(y)$ is {\it invariant} under the variation \rf{259.121} if there exists a function $F(t)$ such that 
\begin{align}
\frac{\partial L}{\partial y} \; \eta + \frac{\partial L}{\partial y'} \; \eta' = \frac{dF}{dt}. \lbl{261.121} 
\end{align}
If $F=0$ we say that the Lagrangian $L$ is invariant.

\noindent Let us next consider a more general variation
\begin{align}
y(t) \rightarrow y(t) + \eps(t) \; \eta(t), \lbl{262.121} 
\end{align}
where $T$ is invariant with respect to the variation \rf{259.121} and $\eps(t)$ is a function that is numerically small, $\abs{\eps(t)} << 1$, vanishes at the boundary points $t_0$ and $t_1$, but which is otherwise arbitrary. 
Inserting the variation \rf{262.121} into the functional \rf{258.121} we get 
\begin{align}
T(y + \eps \; \eta) &= \intt dt \; L(t,y+\eps \; \eta, y' + \eps \; \eta' + \eps' \; \eta) \nonumber \\ 
&= \intt dt \; \{L+ \eps \; ( \frac{\partial L}{\partial y} \; \eta + \frac{\partial L}{\partial y'} \; \eta') + \eta \; \frac{\partial L}{\partial y'} \; \eps' \} + \obs(\eps^2) \nonumber \\
&=T(y)+ \intt dt \; \{ \eps \; \frac{dF}{dt} - \dt{} ( \eta \; \frac{\partial L}{\partial y}) \; \eps \} + \obs(\eps^2) \nonumber \\ 
&=T(y)+\intt dt \; \eps \; \{ \frac{dF}{dt} - \dt{} (\eta \; \frac{\partial L}{\partial y'}) \} + \obs(\eps^2) \nonumber \\ 
&=T(y)+ \intt dt \; \eps \; \dt{j} + \obs(\eps^2), \lbl{263.121}
\end{align}
where we have used \rf{261.121} in line three,  and where the \ttx{Noether current} $j(t)$ is defined to be  
\begin{align}
j = F - \eta \; \frac{\partial L}{\partial y'}\lbl{Ncurrent}.  
\end{align}
Equation \rf{263.121} is true for \ttx{any} $y(t)$. In particular it is true for a $y(t)$ that is a stationary point for the functional $T$. But if $y$ is stationary  we must have 
\begin{align}
T (y + \tilde{\eps} \; \tilde{\eta}) = T(y) + \obs(\tilde{\eps}^2), \lbl{264.12} 
\end{align}
for any $\tilde{\eta}$ of order one and small number $\tilde{\eps}$. 

If we let $\tilde{\eps}$ measure the size of $\eps(t)$ and define 
\begin{align*}
\tilde{\eta}(t) = \frac{\eps(t)}{\tilde{\eps}} \; \eta(t),  
\end{align*}
then the variation \rf{262.121} is exactly of the form 
\begin{align}
y(t) \rightarrow y(t) + \tilde{\eps} \; \tilde{\eta} (t). \lbl{266.12} 
\end{align}
Therefore when $y$ is a stationary point for the functional $T(y)$, we must from \rf{263.121} have 
\begin{align}
\intt dt \; \eps(t) \; \dt{j} = 0, \lbl{267.12} 
\end{align}
and this holds for all $\eps(t)$ that vanishes on the boundaries. The fundamental lemma then implies that 
\begin{align*}
\dt{j} = 0,  
\end{align*}
or in other words, the Noether current corresponding to an invariant for a functional is conserved for any stationary point of the functional. Such stationary $y(t)$'s satisfy, as we recall, the Euler- Lagrange equations corresponding to the functional. This is one instance of Noethers theorem 
\end{example}

\begin{example}\label{Nexample1}
Let us consider a functional 
\begin{align*}
T(y) = \intt dt \; L(y,y'),  
\end{align*}
thus $L$ does not depend explicitly on $t$. \\
We now consider an infinitesimal translation of the variable $t$ 
\begin{align}
t \rightarrow t + \eps. \lbl{270.12} 
\end{align}
The translation \rf{270.12} induces a corresponding variation of $y(t)$ that we find using Taylor's formula
\begin{align}
y(t) \rightarrow y(t+ \eps) = y(t) + \eps \; y'(t) + ...\;\;. \lbl{271.121}
\end{align}
Thus we have a variation of the form \rf{259.121} with $\eta(t) = y'(t)$. \\
Observe,  that using this particular variation we have 
\begin{align*}
\frac{\partial L}{\partial y} \; \eta + \frac{\partial L}{\partial y'} \; \eta' &= \frac{\partial L}{\partial y} \; y' + \frac{\partial L}{\partial y'} \; y''  \\ 
&= \frac{\partial L}{\partial t} + \frac{\partial L}{\partial y} \; y' + \frac{\partial L}{\partial y'} \; y'' = \frac{dL}{d t}. \nonumber 
\end{align*}
Thus the functional is invariant under the variation \rf{271.121} \ttx{because} $L$ does not depend explicitly on $t$. \\The conserved Noether current corresponding to the variation \rf{271.121} is then from \rf{Ncurrent} 
\begin{align*}
j = L - y' \; \frac{\partial L}{\partial y'}.  
\end{align*}
Let us verify directly that $j$ is indeed conserved. Using the Euler-Lagrange equations we have 
\begin{align*}
\dt{j} &= \frac{\partial L}{\partial y} \; y' + \frac{\partial L}{\partial y'} \; y'' - y'' \; \frac{\partial L}{\partial y'} - y' \; \dt{} (\frac{\partial L}{\partial y'})  \\
&= y' \; (\frac{\partial L}{\partial y} - \dt{} (\frac{\partial L}{\partial y'})) = 0. \nonumber 
\end{align*}
\end{example}

\begin{example}\label{Nexample2}
The motion of a mass-point $m$ under the influence of a conservative, time  invariant force
\begin{align*}
F = - \frac{\partial V}{\partial x},\;\;\; V = V(x),  
\end{align*}
is determined by Newton's law 
\begin{align*}
m \; x'' = - \frac{\partial V}{\partial x}.  
\end{align*}
We have seen that this equation is variational with Lagrangian 
\begin{align}
L = \inv{2} \; m \; x'^2 - V(x), \lbl{277.12} 
\end{align}
and we observe that $L$ is invariant under time translation because $V$ does not depend on time.

\noindent By applying the general result from example \ref{Nexample1}, we have the following conserved Noether current 
\begin{align*}
j &= L - x' \; \frac{\partial L}{\partial x'}  \\ 
&= \inv{2} \; m \; x'^2 - V(x) - m \; x'^2 \nonumber \\ 
&= - ( \inv{2} \; m \; x'^2 + V(x)) = - E(t), \nonumber 
\end{align*}
where $E(t)$ is the total \ttx{energy} of the mass-point. \\
Thus the energy is conserved because the Lagrangian \rf{277.12} of the action functional does not depend on time and is thus invariant under time translation. 

This link between energy conservation and invariance under time translation for the Lagrangians holds in general.

The reason why energy conservation plays such a prominent role in our description of nature is because we insist that our natural laws should look the same for all observers, even if they live at different times. Thus, energy conservation is not actually a part of nature, but is rather a consequence of how we choose to describe nature.

\subsubsection{Several dependent variables}
Let us next consider the case when we have several dependent variables 
\begin{align}
T(y_1,y_2,...,y_n) = \intt dt \; L(t , y_1,...,y_n,y_1',...,y_n'). \lbl{279.121} 
\end{align}
We introduce a variation
\begin{align}
y_i \rightarrow y_i + \eps \; \eta_i && i=1,...,n\;\;, \lbl{280.121} 
\end{align}
where $\eta_i = \eta_i(t)$ is a specific set of $n$ functions and $\eps << 1$. \\
Inserting the variation \rf{280.121} into the functional \rf{279.121} we get 
\begin{align*}
&T(y_1 + \eps \; \eta_1,...,y_n + \eps \; \eta_n)  \\
&= \intt dt \; L (t,y_1 + \eps \; \eta_1,...,y_n + \eps \; \eta_n, y_1' + \eps \; \eta_1',...,y_n' + \eps \; \eta_n') \nonumber \\ &= \intt dt \; \{ L(t, y_1, ... , y_n, y_1', ..., y_n') + \frac{\partial L}{\partial y_1} \; \eps \; \eta_1 + ... + \frac{\partial L}{\partial y_n} \; \eps \; \eta_n \nonumber \\ 
&+ \frac{\partial L}{\partial y_1'} \; \eps \; \eta_1' + ... + \frac{\partial L}{\partial y_n'} \; \eps \; \eta_n' \} + \obs(\eps^2) \nonumber \\ 
&= T(y_1,...,y_n) + \eps \; \intt dt\; \{ \frac{\partial L}{\partial y_1} \; \eta_1 + ... + \frac{\partial L}{\partial y_n} \; \eta_n + \frac{\partial L}{\partial y_1'} \; \eta_1' \nonumber \\
&+ ... + \frac{\partial L}{\partial y_n'} \; \eta_n' \} + \obs(\eps^2). \nonumber 
\end{align*}
We now define the functional $T$ to be invariant under the variation \rf{280.121} if there exists a function $F(t)$ such that 
\begin{align}
\frac{\partial L}{\partial y_1} \; \eta_1 + ...+ \frac{\partial L}{\partial y_n} \; \eta_n + \frac{\partial L}{\partial y_1'} \; \eta_1' + ... + \frac{\partial L}{\partial y_n'} \; \eta_n' = \dt{F}.  \lbl{282.121}
\end{align}
If $F=0$ we say that the Lagrangian $L$  is invariant. \\ 
Let us next introduce the more general variation 
\begin{align}
y_i(t) \rightarrow y_i(t) + \eps(t) \; \eta_i(t) \lbl{283.121} 
\end{align}
where the functions $\eps(t)$ satisfies the properties preceding \rf{263.121}. \\
Inserting the variation \rf{283.121} into the functional \rf{279.121} we get
\begin{align*}
&T (y_1 + \eps \; \eta_1,...,y_n + \eps \; \eta_n)  \\ 
&= \intt dt \; L(t,y_1 + \eps \; \eta_1,...,y_n + \eps \; \eta_n, y_1' \nonumber \\  &+ \eps \; \eta_1' + \eps' \; \eta_1, ..., y_n' + \eps \; \eta_n' + \eps' \; \eta_n) \nonumber \\ 
&= \intt dt \; \{ L(t, y_1, ... , y_n, y_1', ..., y_n') + \eps(t) \; \{ \frac{\partial L}{\partial y_1} \eta_1+\frac{\partial L}{\partial y'_1} \eta'_1  \nonumber \\ &+ ...+ \frac{\partial L}{\partial y_n} \; \eta_n+ \frac{\partial L}{\partial y'_n} \; \eta'_n \} +\eps'(t)\{ \eta_1 \; \frac{\partial L}{\partial y_1'}+ \; \nonumber \\ 
&... + \eta_n \; \frac{\partial L}{\partial y_n'} \;\} + \obs(\eps^2) \nonumber \\ 
&= T(y_1,...,y_n) + \intt dt \; \eps(t) \; \{ \frac{dF}{dt} - \dt{} (\eta_1 \; \frac{\partial L}{\partial y_1'}+ \nonumber \\ 
&... + \eta_n \; \frac{\partial L}{\partial y_n'}) \} + \obs(\eps^2) \nonumber \\ 
&= T (y_1,...,y_n) + \intt dt \; \eps(t) \; \dt{j} + \obs(\eps^2), \nonumber 
\end{align*}
where the Noether current is 
\begin{align}
j = F - \eta_1 \; \frac{\partial L}{\partial y_1'} - ... - \eta_n \; \frac{\partial L}{\partial y_n'}\lbl{285.12}  
\end{align}
By an argument identical to \rf{264.12}-\rf{267.12} we conclude that the Noether current is conserved 
\begin{align*}
\dt{j} = 0,  
\end{align*}
for any solutions to the Euler - Lagrange equations 
\begin{align*}
\frac{\partial L}{\partial y_i} - \dt{} (\frac{\partial L}{\partial y'_i}) = 0 && i = 1,...,n\;\;.  
\end{align*}
This is Noether's theorem for the functional \rf{279.121} 
\end{example}

\begin{example}
Let us assume that the Lagrangian $L$, does not depend explicitly on time. It is thus invariant under an infinitesimal time translation
\begin{align*}
t \rightarrow t + \eps  
\end{align*}
The infinitesimal  time translation induce, like in \rf{271.121}, variations of the form 
\begin{align}
y_i(t) \rightarrow y_i(t) + \eps \; y_i'(t), \lbl{289.121} 
\end{align}
and we observe that the functional \rf{279.121} is invariant under the variation \rf{289.121} 
\begin{align*}
& \frac{\partial L}{\partial y_1} \; \eta_1 + ... + \frac{\partial \mathcal{L}}{\partial y_n} \; \eta_n + \frac{\partial \mathcal{L}}{\partial y_1'} \; \eta_1' + ... + \frac{\partial \mathcal{L}}{\partial y_n'} \; \eta_n'  \\ 
&= \frac{\partial \mathcal{L}}{\partial y_1} \; y_1' + ... + \frac{\partial \mathcal{L}}{\partial y_n} \; y_n' + \frac{\partial \mathcal{L}}{\partial y_1'} \; y_1'' + ... + \frac{\partial \mathcal{L}}{\partial y_n'} \; y_n'' \nonumber \\ 
&= \dt{\mathcal{L}}. \nonumber
\end{align*}
Thus we have the following conserved Noether current.
\begin{align}
j = \mathcal{L} - y_1' \; \frac{\partial \mathcal{L}}{\partial y_1'} - ... - y_n' \; \frac{\partial \mathcal{L}}{\partial y_n'}. \lbl{291.12}
\end{align}
\end{example}

\begin{example}
In example \ref{MassPoints} we discussed the action principle for systems of mass  points moving under the influence of conservative forces.
\begin{align*}
T(\vx_1,...,\vx_n) = \intt dt \; L(t , \vx_1 , ... , \vx_n , \vx_1' , ... ,\vx_n'), 
\end{align*}
where the Lagrangian is 
\begin{align}
L = \inv{2} \; \mathlarger{\sum}^n_{i=1} m_i \; \vx_i'^2 - V(t, \vx_1,...,\vx_n). \lbl{293.12} 
\end{align}
Let us assume that the potential does not depend explicitly on time. Then the Lagrangian \rf{293.12} is invariant with respect to translation of time and according to \rf{291.12} we have the following conserved Noether current. 
\begin{align*}
j &= L - \mathlarger{\sum}^n_{i=1} \vx_i' \vdot \frac{\partial L}{\partial \vx_i'}  \\ 
&= \inv{2} \; \mathlarger{\sum}^n_{i=1} m_i \; \vx_i'^2 - V - \mathlarger{\sum}^n_{i=1} m_i \; \vx_i'^2 \nonumber \\ 
&= - (\inv{2} \; \mathlarger{\sum}^n_{i=1} m_i \; \vx_i'^2 + V) = - E, \nonumber 
\end{align*}
where $E$ is the total energy of the system of mass-points. \\
Thus, we see again that energy conservation exists because we insist on natural laws that appear the same for all observers, independently of when they live. 

Let us next assume that the Lagrangian in the functional \rf{279.121} is invariant with respect to the variation 
\begin{align*}
y_i(t) \rightarrow y_i(t) + \eps \; a_i && a_i \in \mathbf{R} \qquad \qquad \qquad i =1,...,n\;\;.  
\end{align*}
Since it is the Lagrangian that is invariant we have that $F=0$ in definition \rf{282.121}, and according to \rf{285.12} we have the conserved Noether current 
\begin{align}
j = \mathlarger{\sum}_i a_i \; \frac{\partial L}{\partial y_i'}. \lbl{296.12} 
\end{align}
\end{example}
\begin{example}
We return to the system of mass-points with Lagrangian 
\begin{align*}
L = \inv{2} \; \mathlarger{\sum}^n_{i=1} m_i \; \vx_i'^2 - V(t,\vx_1,...,\vx_n).  
\end{align*}
Let us assume that the potential is invariant under a translation
\begin{align*}
\vx_i \rightarrow \vx_i +\eps\; \vb{a}_i, 
\end{align*}
then the Lagrangian is invariant under the variation
\begin{align}
\vx_i(t) \rightarrow \vx_i(t) + \eps \; \vb{a}_i, \lbl{299.12} 
\end{align}
and we get, according to \rf{296.12}, the conserved Noether current 
\begin{align*}
j = -\mathlarger{\sum}_i \vb{a}_i \vdot \frac{\partial L}{\partial \vx_i'} = -\mathlarger{\sum}_i m_i \; \vb{a}_i \vdot \vx_i'.  
\end{align*}
The most common situation is when $\vb{a}_i = \vb{a} \; \; \; \forall i$. For this case the Noether current is 
\begin{align*}
j =- \vb{a} \vdot \mathlarger{\sum}_i m_i \; \vx_i'.  
\end{align*}
Thus the component of the total momentum in the direction of $\vb{a}$ is conserved. If the invariance \rf{299.12} holds for three vectors $\vb{a}, \vb{b}, \vb{c},$ that span $\mathbf{R}^3$, we can conclude that the total momentum of the system of mass points is conserved 
\begin{align*}
\vb{P} = \mathlarger{\sum}_i m_i \; \vx_i'.  
\end{align*}
This conservation law holds for example if the potential only depends on differences of the vectors $\vx_i$. \\
In this case, the laws of motion for mass-points looks the same for all observers, independently of where in space they are located. \\
This invariance with respect to location in space is something we \ttx{choose} to impose on \ttx{all} our fundamental natural laws. The consequence of this \ttx{choice} is that we will have conservation of momentum
\end{example}

\begin{example}
Let us assume that the potential for a system of mass-points is invariant under rotation of coordinates around some axis $\vb{k}$.

Recall that rotations of an angle $\theta$ around some axis $\vb{k}$ can be written in the following way
\begin{align*}
\vx \rightarrow \vx \; \cos \theta + (\vb{k} \cp \vx ) \; \sin \theta + \vb{k} \; (\vb{k} \vdot \vx) \; (1 - \cos \theta).  
\end{align*}
This is \textit{Rodrigue's formula}. For an infinitesimal rotation angle we get 
\begin{align}
\vx \rightarrow \vx + \eps (\vb{k} \cp \vx). \lbl{304.12} 
\end{align}
The rotation of coordinates \rf{304.12}, induces a corresponding variation
\begin{align*}
\vx_i(t) \rightarrow \vx_i(t) + \eps \; (\vb{k} \cp \vx_i(t)).  
\end{align*}
Since the Lagrangian is conserved by this variation we have $F=0$ in \rf{282.121} and the conserved Noether current is from \rf{285.12} 
\begin{align*}
j &= - \mathlarger{\sum}^n_{i=1} (\vb{k} \cp \vx_i) \vdot \frac{\partial L}{\partial \vx_i'}  \\ 
&= - \mathlarger{\sum}^n_{i=1} (\vb{k} \cp \vx_i) \vdot (m_i \; \vx_i') \nonumber \\
&= - \mathlarger{\sum}^n_{i=1} \vb{k} \vdot (\vx_i \cp m_i \; \vx_1') \nonumber \\ 
&= - \vb{k} \vdot \mathlarger{\sum}^n_{i=1} \vx_i \cp m_i \; \vx_i' = - \vb{k} \vdot \vb{L},
\end{align*}
where $\vb{L}$ is the total angular momentum of the system of mass-points. \\
Thus, if the Lagrangian is invariant with respect to rotations around some axis $\vb{k}$, then the component of the total angular momentum along $\vb{k}$ is conserved. \\
If the system is invariant with respect to rotations around three separate axes, that together span $\mathbf{R}^3$, the total angular momentum is conserved. This is for example true if the Lagrangian is invariant with respect to all possible rotations in $\mathbf{R}^3$.

Like for the location in time and space, this invariance with respect to direction, it is called \ttx{isotropy}, is something we \ttx{choose} to impose on all our natural laws. 

The consequence of this choice is that we will have conservation of angular momentum. \\
Our natural laws in general involve fields and partial differential equations. Maxwell's equations for the electromagnetic field is an example that comes to mind. 

In all these cases, invariance with respect to location in space and time and with respect to direction in space will lead to conserved Noether currents that we, by convention, call energy, momentum and angular momentum. 
\end{example}

\subsubsection{Several independent variables}

We consider a functional of the form 
\begin{align}
T(u) = \int_D dx \; dy \; \mathcal{L}(x,y,u,u_x,u_y) \lbl{307.12} 
\end{align}
Let us consider a variation
\begin{align}
u(x,y) \rightarrow u(x,y) + \eps \; \eta(x,y) \lbl{308.12} 
\end{align}
Where $\eta$ is some specific function.
Inserting \rf{308.12} into \rf{307.12} we get
\begin{align*}
T(u + \eps \; \eta) &= \int_D dx \; dy \; \mathcal{L}(x,y,u + \eps \; \eta, u_x + \eps \; \eta_x, u_y + \eps \; \eta_y)  \\
&= \int_D dx \; dy \; \{ \mathcal{L}(x,y,u,u_x,u_y) + \frac{\partial \mathcal{L}}{\partial u} \; \eps \; \eta + \frac{\partial \mathcal{L}}{\partial u_x} \; \eps \; \eta_x \nonumber \\
&+ \frac{\partial \mathcal{L} }{\partial u_y} \; \eps \; \eta_y  \} + \obs(\eps^2) \nonumber \\
&= T(u) + \eps \; \int_D dx \; dy \; \{ \frac{\partial \mathcal{L}}{\partial u} \; \eta + \frac{\partial \mathcal{L}}{\partial u_x} \; \eta_x + \frac{\partial \mathcal{L}}{\partial u_y} \; \eta_y \} + \obs(\eps). \nonumber 
\end{align*}
The functional \rf{307.12} is said to be invariant under the variation \rf{308.12}, if there exists functions $F_1(x,y), \; F_2(x,y)$ such that 
\begin{align}
\frac{\partial \mathcal{L}}{\partial u} \; \eta + \frac{\partial \mathcal{L}}{\partial u_x} \; \eta_x + \frac{\partial \mathcal{L}}{\partial u_y} \; \eta_y = \prt{x} F_1 + \prt{y} F_2 \lbl{310.121} 
\end{align}
If $F_1 = F_2 = 0$ we say that the Lagrangian density $\mathcal{L}$, is invariant under the variation \rf{308.12}.

Let us now assume that the functional $T$ is invariant under a variation of the form \rf{308.12} for some specific choice of $\eta(x,y)$. Using this $\eta(x,y)$, we consider the more general variation
\begin{align}
u(x,y) \rightarrow u(x,y) + \eps(x,y) \; \eta(x,y), \lbl{311.12} 
\end{align}
 where $\eps(x,y)$ is a numerically small function that vanished on the boundary of $D$ but is otherwise arbitrary. 
Inserting \rf{311.12} into the functional \rf{307.12} we get 
\begin{align}
T (u + \eps \; \eta) &= \int_D dx \; dy \; \mathcal{L}(x,y,u + \eps \; \eta_, u_x + \eps \; \eta_x + \eps_x \; \eta, u_y + \eps \; \eta_y + \eps_y \; \eta) \nonumber \\ 
&= \int_D dx \; dy \{ \mathcal{L}(x,y,u,u_x,u_y) + \frac{\partial \mathcal{L}}{\partial u} \; \eps \; \eta + \frac{\partial \mathcal{L}}{\partial u_x} \; (\eps \; \eta_x + \eps_x \; \eta) \nonumber \\
&+ \frac{\partial \mathcal{L}}{\partial u_y} \; (\eps \; \eta_y + \eps_y \; \eta)  \} + \obs(\eps^2) \nonumber \\ 
&= T(u) + \int_D dx \; dy \; \{ \eps(\frac{\partial \mathcal{L}}{\partial u} \eta + \frac{\partial \mathcal{L}}{\partial u_x} \; \eta_x + \frac{\partial \mathcal{L}}{\partial u_y} \; \eta_y) \nonumber \\
&+ \eta \; \frac{\partial \mathcal{L}}{\partial u_x} \; \eps_x + \eta \; \frac{\partial \mathcal{L}}{\partial u_y} \; \eps_y  \} + \obs(\eps^2) \nonumber \\ 
&= T(u) + \int_D dx \; dy \; \{ \eps \; (\prt{x} F_1 + \prt{x} F_2) - \prt{x} \; (\eta \; \frac{\partial \mathcal{L}}{\partial u_x}) \; \eps \nonumber \\
&- \prt{y} (\eta \; \frac{\partial \mathcal{L}}{\partial u_y}) \; \eps   \} + \obs(\eps^2), \lbl{312.12}
\end{align}
where we have used the divergence theorem and the boundary conditions on $\eps$ in addition to the definition \rf{310.121} of  invariance of $T$ with respect to the variation \rf{308.12}. \\
Thus from \rf{312.12} we have 
\begin{align}
T (u + \eps \; \eta) &= T (u)  + \int_D dx \; dy \; \eps \; \{ \prt{x}j_1 + \prt{y} j_2  \} + \obs(\eps^2), \lbl{313.12} 
\end{align}
where the Noether current $\vb{j} = (j_1, j_2)$ is 
\begin{align}
j_1 &= F_1 - \eta \; \frac{\partial \mathcal{L}}{\partial u_x}, \nonumber \\ 
j_2 &= F_2 - \eta \; \frac{\partial \mathcal{L}}{\partial u_y}. \lbl{314.121}
\end{align}
Equation \rf{313.12} holds for all $u$. If $u$ is extremal, all variations of $T(u)$ vanish, in particular they vanish for variations of the form \rf{311.12}. Using a 2D version of the fundamental lemma we conclude that the Noether current satisfies
\begin{align}
\prt{x} j_1  + \prt{y} j_2 = 0. \lbl{315.12}
\end{align}
This kind of identity is in general called a \ttx{conservation law}.

\begin{example}
In example \ref{Example9} we discussed the functional 
\begin{align*}
T(u) = \iint\limits_D dx \; dy \; \left\{\inv{2} \; u_x^2 + \inv{2} \; u_y^2\right\} . 
\end{align*}
The Lagrangian density $\mathcal{L}$, is here
\begin{align*}
\mathcal{L} = \inv{2} \; u_x^2 + \inv{2} \; u_y^2, 
\end{align*}
and we found that the corresponding Euler-Lagrange equation is the Laplace equation
\begin{align*}
u_{xx} + u_{yy} = 0.  
\end{align*}
Observe that $\mathcal{L}$ invariant under the variation 
\begin{align*}
u(x,y) \rightarrow u(x,y) + \eps . 
\end{align*}
Thus $\eta = 1$ and the components of the Noether current are 
\begin{align*}
j_1 &= - \frac{\partial \mathcal{L}}{\partial u_x} = - u_x,  \\ 
j_2 &= - \frac{\partial \mathcal{L}}{\partial u_y} = - u_y.                                                                                                                                                                                                                                                                                                                                                    \nonumber
\end{align*}
The conservation law \rf{315.12} is thus 
\begin{align*}
\prt{x} j_1 + \prt{y} j_1 &= 0, \nonumber \\
&\Updownarrow \nonumber \\  u_{xx} + u_{yy} &= 0.  
\end{align*}
This conserved current does not tell us anything new since it's conservation law  is just the Laplace equation itself. 

We also observe that $\mathcal{L}$ is independent of $x$ and $y$. \\
Let us consider an arbitrary infinitesimal translation in the plane
\begin{align}
x \rightarrow x + \eps \; a, \lbl{322.12} \\ 
y \rightarrow y + \eps \; b, \nonumber \\
\nonumber \\ 
a^2 + b^2 = 1. \nonumber 
\end{align}
Thus $\vb{n} = (a,b)$ is a unit vector determining the direction of the translation.
The infinitesimal translation \rf{322.12} induces the following variation in $u$
\begin{align*}
u(x,y) \rightarrow u(x,y) + \eps \; (a \; u_x + b \; u_y).  
\end{align*}
Observe that 
\begin{align*}
& \frac{\partial \mathcal{L}}{\partial u} \; (a \; u_x + b \; u_y) + \frac{\partial \mathcal{L}}{\partial u_x} \; (a \; u_{xx} + b \; u_{xy}) + \frac{\partial \mathcal{L}}{\partial u_y} \; (a \; u_{xy} + b \; u_{yy})  \\ 
&= a \; (\frac{\partial \mathcal{L}}{\partial u} \; u_x + \frac{\partial \mathcal{L}}{\partial u_x} \; u_{xx} + \frac{\partial \mathcal{L}}{\partial u_y} \; u_{xy} ) \nonumber \\ 
&= b \; (\frac{\partial \mathcal{L}}{\partial u} \; u_y + \frac{\partial \mathcal{L}}{\partial u_x} \; u_{xy} + \frac{\partial \mathcal{L}}{\partial u_y} \; u_{yy} ) \nonumber \\ 
&= \prt{x} (a \; \mathcal{L}) + \prt{y}(b \; \mathcal{L}), \nonumber 
\end{align*}
since $\mathcal{L}$ does not depend explicitly on $x$ and $y$. \\
Thus from \rf{310.121} we conclude that 
\begin{align*}
F_1 = a\; \mathcal{L}, \;\;\;F_2 =b \; \mathcal{L},  
\end{align*}
and the components of the Noether current are 
\begin{align*}
j_1 &= a \; \mathcal{L} - (a \; u_x + b \; u_y) \;\frac{\partial \mathcal{L}}{\partial u_x}  \\ 
&= a \; \mathcal{L} - (a \; u_x + b \; u_y) \; u_x, \nonumber \\ 
j_2 &= b \; \mathcal{L} - (a \; u_x + b \; u_y) \; \frac{\partial \mathcal{L}}{\partial u_y} \nonumber \\ 
&= b \; \mathcal{L} - (a \; u_x + b \; u_y) \; u_y, \nonumber 
\end{align*}
and thus the conservation law \rf{315.12} is 
\begin{align}
\prt{x} j_1 + \prt{y} j_2 &= 0 \nonumber \\ 
&\Updownarrow \nonumber \\  a \; \prt{x} (\inv{2} \; u_x^2 + \inv{2} \; u_y^2) - \prt{x}(a \; u_x^2 + b \; u_x \; u_y)& \nonumber \\
+ b \; \prt{y} (\inv{2} \; u_x^2 + \inv{2} \; u_y^2) - \prt{y} (a \; u_x \; u_y + b \; u_y^2)& = 0 \nonumber \\
&\Updownarrow \nonumber \\  \; a \; u_x \; u_{xx} + a \; u_y \; u_{xy} - 2 \; a \; u_x \; u_{xx} - b \; u_{xx} \; u_y& \nonumber \\
 - b \; u_x \; u_{xy} + b \; u_x \; u_{xy} + b \; u_y \; u_{yy}& \nonumber \\
 - a \; u_{xy} \; u_y - a \; u_x \; u_{yy} - 2 \; b \; u_y \; u_{yy}& = 0 \nonumber \\
&\Updownarrow \nonumber \\   - (a \; u_x + b \; u_y) \;( u_{xx} + u_{yy})& = 0. \lbl{327.121} 
\end{align}
The conservation law \rf{327.121} clearly holds for any solution to the Laplace equation and does not tell us anything new.
\end{example}

\begin{example}\label{Example22}
From example \ref{Example10} we have seen that the 1D wave equation is the Euler-Lagrange equation for the functional 
\begin{align*}
T (u) = \intx dx \; \intt dt \; \{ \inv{2} \; u_t^2 - \inv{2} \; c^2 \; u_x^2 \}.  
\end{align*}
We introduce an infinitesimal time translation 
\begin{align*}
t \rightarrow t + \eps,  
\end{align*}
and this induces a variation of the form 
\begin{align*}
u(t,x) \rightarrow u(t,x) + \eps \; u_t(t,x).  
\end{align*}
For this variation $\eta = u_t$ and we have 
\begin{align*}
& \frac{\partial \mathcal{L}}{\partial u} \; \eta + \frac{\partial \mathcal{L}}{\partial u_t} \; \eta_t + \frac{\partial \mathcal{L}}{\partial u_x} \; \eta_x  \\ 
&= u_t \; (u_{tt}) +  (-c^2 \; u_x)\;u_{xt} \nonumber \\
&= u_t \; u_{tt} - c^2 \; u_x \; u_{xt} \nonumber \\ 
&= \prt{t} (\inv{2} \; u_t^2 - \inv{2} \; c^2 \; u_x^2), \nonumber 
\end{align*}
and thus, according to \rf{310.121}, the functional $T$ is invariant with 
\begin{align*}
F_1 &= \inv{2} \; u_t^2 - \inv{2} \; c^2 \; u_x^2 , \\ 
F_2 &= 0. \nonumber 
\end{align*}
The components of the conserved Noether current are according to \rf{314.121} 
\begin{align*}
j_1 &= \inv{2} \; u_t^2 - \inv{2} \; c^2 \; u_x^2 - u_t \; u_t  \\
&= - (\inv{2} \; u_t^2 + \inv{2} \; c^2 \; u_x^2), \nonumber \\
j_2 &= -u_t\; (-c^2 \; u_x)= c^2 \; u_t \; u_x, \nonumber
\end{align*}
and the conservation law is 
\begin{align}
\prt{t} j_1 + \prt{x} j_2=0.\label{ConservationLaw}  
\end{align}
Inserting the expressions for the components of the Noether current, we find the following form for the conservation law.
\begin{align*}
u_t (u_{tt}-c^2u_{xx})=0.
\end{align*}
This is clearly satisfied for any solution to the wave equation.
In order to see why it is natural to call equation (\ref{ConservationLaw}), for the Noether current, a convervation law, define
\begin{align*}
E(t) = - \intx dx \; j_1 (t,x).  
\end{align*}
For $E(t)$ we have 
\begin{align}
\dt{E} = - \intx dx \; \prt{t} j_1 = \intx dx \; \prt{x} j_2 = j_2\mathlarger{|}^{x_1}_{x_0}. \lbl{336.12} 
\end{align}
Let us assume that $u$ satisfies one of the following boundary conditions 
\begin{align*}
u(t,x_0)& = u(t,x_1) = 0,\;\;\;\;\text{Dirichlet}\\
\\
u_x(t,x_0)& = u_x(t,x_1) = 0,\;\;\;\;\text{Neumann}\\
\\
x_0= -\infty&, \; \; \; x_1 = + \infty,\;\;\;\;\;\;\;\text{and $u(x,t)$ vanish at $\pm \infty$}
\end{align*}

\noindent If this the case we, conclude that 
\begin{align*}
\dt{E} &= 0, \; \; \; \nonumber \\ &\Updownarrow \nonumber \\  \; \; E &= \text{const},  
\end{align*}
and thus 
\begin{align*}
E = \intx dx \; (\inv{2} \; u_t^2 + \inv{2} \; c^2 \; u_x^2),  
\end{align*}
is a \ttx{conserved quantity} for any solution to the wave equation. \\
It is in general true that, whenever we have a functional of time dependent fields, and the Lagrangian of the functional does not depend explicitly on time, the space integral of the time component of the Noether current will be a conserved quantity for appropriate boundary conditions at the spatial boundary. 
Inspired by the situation for systems of mass points, we \ttx{define} the space integral of the time component of the Noether current to be the total energy of the system of time dependent fields. 

The time component of the Noether current is defined to be the \ttx{energy density} and the space component(s) are called the \ttx{energy flux density}. 

Thus for the current example the energy density is 
\begin{align}
e = \inv{2} \; u_t^2 + \inv{2} \; c^2 \; u_x^2, \lbl{339.12} 
\end{align}
and the energy flux density is 
\begin{align*}
f = c^2 \; u_t \; u_x.  
\end{align*}
Actually, here $e = - j_1$. We use a standard sign convention that ensures that the energy density is positive.

In some cases the formula for the energy density of a field system is known from physical modelling. In all such cases the energy density derived from the modeling is equal or proportional to the energy density defined using time translation invariance and Noethers theorem. 

For the current example we know that the 1D wave equation is a model for small vibrations of a string. With respect to this example we recognize the first term in \rf{339.12} to be proportional to the kinetic energy, $K$,  and the second term to be proportional to the potential energy, $V$,  of a small piece of the string. 

Thus for this example we clearly we have
\begin{align*}
e \propto K + V.  
\end{align*}

\noindent Let us next introduce an infinitesimal space translation
\begin{align*}
x \rightarrow x + \eps.  
\end{align*}
This introduces a variation of the form 
\begin{align*}
u(t,x) \rightarrow u(t,x) + \eps \; u_x (t,x). 
\end{align*}
For this variation $\eta = u_x$  and we have 
\begin{align*}
&\frac{\partial \mathcal{L}}{\partial u} \; \eta + \frac{\partial \mathcal{L}}{\partial u_t} \; \eta_t + \frac{\partial \mathcal{L}}{\partial u_x} \; \eta_x  \\ 
&= u_t \; u_{xt} + (-c^2 \; u_x) \; u_{xx} \nonumber \\ 
&= u_t \; u_{xt} - c^2 \; u_x \; u_{xx} \nonumber \\
&= \prt{x} (\inv{2} \; u_t^2 - \inv{2} \; c^2 \; u_x^2), \nonumber   
\end{align*}
and thus according to \rf{310.121}, the functional $T$ is invariant with 
\begin{align*}
F_1 &= 0,  \\
F_2 &= \inv{2} \; u_t^2 - \inv{2} \; c^2 \; u_x^2, \nonumber 
\end{align*}
The components of the conserved Noether current are according to \rf{314.121} 
\begin{align*}
j_1 &= - u_x \; u_t,  \\ 
j_2 &= \inv{2} \; u_t^2 - \inv{2} c^2 \; u_x^2 - u_x \; (- c^2 \; u_x) \nonumber \\
&= \inv{2} \; u_t^2 + \inv{2} \; c^2 \; u_x^2, \nonumber 
\end{align*}
and we get the conservation law 
\begin{align*}
\prt{t} (- u_t \; u_x) + \prt{x} (\inv{2} \; u_t^2 + \inv{2} \; c^2 \; u_x^2) = 0.  
\end{align*}
This holds for any solution to the wave equation. Let us verify this directly 
\begin{align*}
&\prt{t} (- u_t \; u_x) + \prt{x} (\inv{2} \; u_t^2 + \inv{2} \; c^2 \; u_x^2)  \\ 
&= - u_{tt} \; u_x  - u_t \; u_{xt} + u_t \; u_{xt} + c^2 \; u_x \; u_{xx} \nonumber \\ 
&= - u_x \; (u_{tt} - c^2 \; u_{xx}) = 0. 
\end{align*}
Defining 
\begin{align*}
P = \intx dx \; (- u_t \; u_x)  
\end{align*}
We get, with the same caveat about boundary conditions as after \rf{336.12},
\begin{align*}
\frac{dP}{dt} &= \intx dx \; \prt{t} (-u_t \; u_x)  \\ 
&= - \intx dx \; \prt{x} (\inv{2} \; u_t^2 + \inv{2} \; c^2 \; u_x^2) \nonumber \\ 
&= - (\inv{2} \; u_t^2 + \inv{2} \; c^2 \; u_x^2) \mathlarger{|}^{x_1}_{x_0} = 0, 
\end{align*}
and thus $P$ is a conserved quantity for the wave equation. Since this conservation law comes from the space translation through the use of Noether's theorem, the time component of the Noether current is defined to be the \ttx{momentum density} and the spatial component is defined to be the \ttx{momentum flux density}. 

\noindent Let us next consider functionals of the form
\begin{align}
T(u) = \int_D dx \; dy \; \intt dt \; \mathcal{L}(t,x,y,u,u_t,u_x,u_y). \lbl{351.12} 
\end{align}
Following the by now familiar procedure, we define \rf{351.12} to be invariant under a variation
\begin{align*}
u(t,x,y) \rightarrow u(t,x,y) + \eps \; \eta(t,x,y),  
\end{align*}
if there exists functions $F_1, F_2$ and $F_3$ such that 
\begin{align*}
\frac{\partial \mathcal{L}}{\partial u} \; \eta + \frac{\partial \mathcal{L}}{\partial u_t} \; \eta_t + \frac{\partial \mathcal{L}}{\partial u_x} \; \eta_x + \frac{\partial \mathcal{L}}{\partial u_y} \; \eta_y = \prt{t} F_1 + \prt{x} F_2 + \prt{y} F_3. \nonumber    
\end{align*}
The conserved Noether current components are 
\begin{align}
j_1 &= F_1 - \eta \; \frac{\partial \mathcal{L}}{\partial u_t}, \lbl{354.121} \\ 
j_2 &= F_2 - \eta \; \frac{\partial \mathcal{L}}{\partial u_x}, \nonumber \\ 
j_3 &= F_3 - \eta \; \frac{\partial \mathcal{L}}{\partial u_y}, \nonumber    
\end{align}
and the conservation law is 
\begin{align*}
\prt{t} j_1 + \prt{x} j_2 + \prt{y} j_3 = 0.  
\end{align*}
This is Noether's theorem for the functional \rf{351.12}. 
\end{example}

\begin{example}
The 2D wave equation has been seen, in example \ref{Example10}, to be variational with Lagrangian
\begin{align*}
\mathcal{L} = \inv{2} \; u_t^2 - \inv{2} \; c^2 \; u_x^2 - \inv{2} \; c^2 \; u_y^2.  
\end{align*}
This Lagrangian density  is translation invariant with respect to time and space and this will lead to conserved Noether currents whose spatial integral of the time component of the current will be the total energy and the total momentum. Let us start with time translation 
\begin{align*}
t \rightarrow t+\eps,  
\end{align*}
which leads to the variation
\begin{align*}
u(t,x,y) \rightarrow u(t,x,y) + \eps \; u_t(t,x,y).  
\end{align*}
In a calculation entirely similar to the one in example \ref{Example22} we find that the functional i invariant with $F_1 = \mathcal{L}, \; F_2 = F_3 = 0$. \\
The components of the Noether current are then from \rf{354.121} 
\begin{align*}
j_1 &= \inv{2} \; u_t^2 - \inv{2} \; c^2 \; u_x^2 - \inv{2} \; c^2 \; u_y^2  - u_t \; (u_t) \\
&= - (\inv{2} \; u_t^2 + \inv{2} \; c^2 \; u_x^2 + \inv{2} \; c^2 \; u_y^2), \nonumber \\ 
j_2 &= - u_t \; (-c^2 \; u_x) = c^2 \; u_t \; u_x, \nonumber \\
j_3 &= - u_t \; (-c^2 \; u_y) = c^2 \; u_t \; u_y, \nonumber  
\end{align*}
and the conservation law is 
\begin{align*}
\prt{t} j_1 + \prt{x} j_2 + \prt{y} j_3 = 0.  
\end{align*}
Defining the energy density to be 
\begin{align*}
e = \inv{2} \; u_t^2 + \inv{2} \; c^2 \; u_x^2 + \inv{2} \; c^2 \; u_y^2,  
\end{align*}
and the energy flux density to be 
\begin{align*}
\vb{f} = - c^2 \; u_t \; \grad{u},  
\end{align*}
where $\grad{} $ is the 2D gradient operator, we have 
\begin{align}
\prt{t} e + \div{\vb{f}} = 0. \lbl{363.12} 
\end{align}
For the total energy inside a domain $D \subset \mathbf{R}^2$ we have 
\begin{align*}
E = \int_D dx \; dy \; e, 
\end{align*}
and using \rf{363.12} we get 
\begin{align}
\dt{E} &= \int_D dx \; dy \; \prt{t}e = - \int_D dx \; dy \; \div{\vb{f}} \nonumber \\ 
&= - \int_{\partial D} dl \; \vb{f} \vdot \vb{n} \lbl{365.12} 
\end{align}
The sign convention chosen for $\vb{f}$ ensures that $\vb{f} \vdot \vb{n} > 0$ means that energy is \ttx{leaving} the domain $D$.
\begin{figure}[htbp]
\centering
\includegraphics{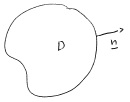}
\caption{The integration domain for the functional and convention for unit normal on the boundary}
\label{fig15}
\end{figure}

This convention is standard. 
Equation \rf{365.12} clearly expresses conservation of energy; the only way energy can change in $D$ is if energy is escaping across the boundary. If the flux is zero on $\partial D$, which would happen if for example $u\mathlarger{|}_D = 0$, or if $D = \mathbf{R}^2$ and solutions decay at infinity, we get conservation of total energy
\begin{align*}
\dt{} \int_D dx \; dy \; \{ \inv{2} \; u_t^2 + \inv{2} \; c^2 \; u_x^2 + \inv{2} \; c^2 \; u_y^2 \} = 0.  
\end{align*}
Let us next look at translation along a direction $\vb{a} = (a,b)$, $a^2 + b^2 = 1$ 
\begin{align*}
\vx = (x,y) \rightarrow \vx + \eps \; \vb{a}.  
\end{align*}
The corresponding variation is clearly 
\begin{align*}
u(t,x,y) \rightarrow u(t,x,y) + \eps(a \; u_x + b \; u_y),  
\end{align*}
and for this variation we have 
\begin{align*}
&\frac{\partial \mathcal{L}}{\partial u} \; \eta + \frac{\partial \mathcal{L}}{\partial u_t} \ \eta_t  + \frac{\partial \mathcal{L}}{\partial u_x} \; \eta_x+ \frac{\partial \mathcal{L}}{\partial u_y} \; \eta_y  \\ 
&=\frac{\partial \mathcal{L}}{\partial u_t} \; (a \; u_{xt} + b \; u_{yt} ) + \frac{\partial \mathcal{L}}{\partial u_x} \; (a \; u_{xx} + b \; u_{xy}) + \frac{\partial \mathcal{L}}{\partial u_y} \; (a \; u_{xy} + b \; u_{yy}) \nonumber \\
&= a \; (\frac{\partial \mathcal{L}}{\partial u_t} \; u_{xt} + \frac{\partial \mathcal{L}}{\partial u_x} \; u_{xx} + \frac{\partial \mathcal{L}}{\partial u_y} \; u_{xy}) \nonumber \\
&+ b \; (\frac{\partial L}{\partial u_t} \; u_{yt} + \frac{\partial \mathcal{L}}{\partial u_x} \; u_{yx} + \frac{\partial \mathcal{L}}{\partial u_y} \; u_{yy}) \nonumber \\ 
&= a \; \prt{x} \mathcal{L} + b \; \prt{y} \mathcal{L}. \nonumber 
\end{align*}
So the functional is invariant with $F_1 = 0, \; F_2 = a \; \mathcal{L}, \; F_3 = b\; \mathcal{L}$. 
The components of the Noether current are from \rf{354.121} 
\begin{align*}
j_1 &= - (a \; u_x + b \; u_y) \; u_t = - u_t \; (a \; u_x + b \; u_y),  \\
j_2 &= a \; (\inv{2} \; u_t^2 - \inv{2} \; c^2 \; u_x^2 - \inv{2} \; c^2 \; u_y^2) - (a \; u_x + b \; u_y) \; (-c^2 \; u_x) \nonumber \\ 
&= a \; \inv{2} \; u_t^2 - a \; \inv{2} \; c^2 \; u_x^2 - a \; \inv{2} \; c^2 \; u_y^2 + a \; c^2 \; u_x^2 + b \; c^2 \; u_x \; u_y \nonumber \\
&= \inv{2} \; a \; u_t^2  + \inv{2} \; a \; c^2 \; u_x^2 - \inv{2} \; a \; c^2 \; u_y^2 + b \; c^2 \; u_x \; u_y, \nonumber \\
j_3 &= b\; (\inv{2} \; u_t^2 - \inv{2} \; c^2 \; u_x^2 - \inv{2} \; c^2 \; u_y^2) - (a \; u_x + b\; u_y) \; (-c^2 \; u_y) \nonumber \\ 
&= \inv{2} \; b \; u_t^2 - \inv{2} \; b \; c^2 \; u_x^2 + \inv{2} \; b \; c^2 \; u_y^2 + a \; c^2 \; u_x \; u_y, \nonumber 
\end{align*}
and the conservation law is 
\begin{align*}
\prt{t} j_1 + \prt{x} j_2 + \prt{y} j_3 = 0,  
\end{align*}
and looks kind of messy. However, if we introduce a vector $\vb{P}$ and a Cartesian tensor $f$ of rank 2 by 
\begin{align*}
\vb{P} &= u_t \; \grad{u},  \\
f &= c^2 \; \grad{u} \; \grad{u} - \inv{2} \; c^2 \; \Tr (\grad{u} \; \grad{u}) \; I + \inv{2} \; u_t^2 \; I,  
\end{align*}
where $I$ is the identity matrix, we can write the conservation law as 
\begin{align}
\prt{t} \vb{a} \vdot \vb{P} + \div{(\vb{a} \vdot f)} = 0. \lbl{375.12} 
\end{align}
The argument leading up to \rf{375.12} is true for all vectors $\vb{a}$. Therefore we have the conservation law 
\begin{align}
\prt{t} \vb{P} + \div{f} = 0. \lbl{376.12}
\end{align}
By definition, $\vb{P}$ is the momentum density and $f$ is the momentum flux density. \\
The total momentum inside some domain $D$ is 
\begin{align*}
\mathcal{P} = \int_D dx \; dy \; \vb{P},  
\end{align*}
and from \rf{376.121} we get 
\begin{align}
\frac{d \mathcal{P}}{dt} = \int_D dx \; dy \; \prt{t} \vb{P} = - \int_{\partial D} dl \; \vb{f} \vdot \vb{n}. \lbl{378.121} 
\end{align}
All of this can be generalized to the 3D wave equation. Then \rf{378.121} will involve a surface integral of the momentum flux over a 2D surface $\partial D$ bounding a 3D domain $D \subset \mathbf{R}^3$.
\end{example}
Formula \rf{378.121}, and formulas like it for other field systems, have important practical applications.

Let us for example assume that an object, filling a domain $D \subset \mathbf{R}^3$, is embedded in a wave field satisfying the 2D-wave equation. Let us apply formula \rf{378.121} to $\mathbf{R}^2 - D$ and assume that the wave field is localized so that we get no contribution from the boundary at infinity. The only boundary to $\mathbf{R}^2 - D$ is then $\partial D$. The outward normal to $\mathbf{R}^2 - D$ points into $D$ 
\begin{figure}[htbp]
\centering
\includegraphics{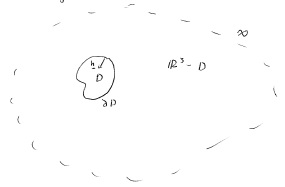}
\caption{}
\label{fig16}
\end{figure}

\noindent Since total momentum is conserved, we know that momentum lost from $\mathbf{R}^3 - D$ to $D$ through the boundary $\partial D$, must be gained by $D$. Similarly momentum gained by $\mathbf{R} - D$ must be lost by $D$. 

If $\vb{Q}$ is the total momentum inside $D$ we thus have, using \rf{378.121}, that
\begin{align*}
\dt{\vb{Q}} = - \dt{\vb{P}} = \int_{\partial D} dl \; \vb{f} \vdot \vb{n} = - \int_{\partial D} dl \; \vb{f} \vdot \vb{s},  
\end{align*}
where now $\vb{s} = -\vb{n}$ is the normal to $\partial D$ pointing \ttx{out} of $D$. \\
Recall that change of momentum per. unit time is in fact \ttx{force}. Thus the force, $\vb{F} = \dt{\vb{Q}}$, acting on the embedded object $D$ in the wave field is 
\begin{align*}
\vb{F} = - \int_{\partial D} dl \; \vb{f} \vdot \vb{s}.  
\end{align*}
This is an important result since it allows us, for example, to calculate the surface-wave induced mechanical forces on ship hulls and stationary structures like oil platforms or structures that generate electric  power from surface waves. 

We end our exposition of the calculus of variation and Noether's theorem here. This is not because this is the end of the story for Noether's theorem. 
For Noether's theorem, and also for the calculus of variations, we have barely scratched the surface. Beyond these notes there is a whole universe to explore.

\subsection{Exercises}

\begin{enumerate}
\item Find all extremals for the functional%
\[
I(y)=\int_{t_{0}}^{t_{1}}dte^{t}\sqrt{1+(y^{\prime})^{2}}
\]

\item Find the unique extremal for
\[
I(y)=\int_{0}^{1}dt(y^{\prime2}-y^{2}+2ty)
\]

satisfying the constrains $y(0)=0,y(1)=1$.

\item Show that there is no extremal for the functional%
\[
I(y)=\int_{-1}^{1}dtt^{4}(y^{\prime})^{2}
\]

that satisfy the constraints $y(-1)=-1,y(1)=1$.

\item Find a unique smooth extremal for
\[
I(y)=\int_{2}^{3}dty^{2}(1-y^{\prime})^{2}
\]

that satisfy the constraints $y(2)=1,y(3)=\sqrt{3}$.

\item Derive the Euler-Lagrange equation for a functional of the general form%
\[
I(y)=\int_{t_{0}}^{t_{1}}dt\mathcal{L}(t,y,y^{\prime},y^{\prime\prime})
\]

Find the unique extremal to the functional%
\[
I(y)=\int_{0}^{1}dt(y^{\prime\prime})^{2}
\]

that satisfy the constraints%
\begin{align*}
y(0)  & =y^{\prime}(0)=0\\
y(1)  & =y^{\prime}(1)=0
\end{align*}

\item Find the extremals for the functional%
\[
I(y)=\int_{0}^{1}dt(y^{\prime\prime})^{2}
\]

subject to the constraints%
\[
J(y)=\int_{0}^{1}dty=L
\]

and $y(0)=0,y(1)=2,y'(0)=1,y'(1)=-1$.

\item Find the extremals of the functional%
\[
I(y)=\int_{0}^{1}dt(y^{\prime})^{2}
\]

subject to the constraints%
\begin{align*}
J_{1}(y)  & =\int_{0}^{1}dty=2\\
J_{2}(y)  & =\int_{0}^{1}dtyt=\frac{1}{2}%
\end{align*}

and $y(0)=y(1)=0$.

\item On page \pageref{ComplexValuedFields} in these lecture notes,  we
discussed extremals of real valued functionals whose domain of definition
consisted of a complex valued function $A$. We argued that the Euler-Lagrange
equations could be found by varying the field $A$ and it's complex conjugate
$A^{\ast}$ as if they were independent variables. Let us consider such
functionals of the form%
\begin{equation}
I(A,A^{\ast})=\int_{t_{0}}^{t_{0}}dt\mathcal{L}(t,A,A^{\ast},A_{t},A_{t}%
^{\ast})\label{Functional1}%
\end{equation}
Show that the Euler-Lagrange equations for this type of functional are%
\begin{align*}
\frac{\partial\mathcal{L}}{\partial A}-\frac{d}{dt}\left(  \frac
{\partial\mathcal{L}}{\partial A_{t}}\right)   & =0\\
\frac{\partial\mathcal{L}}{\partial A^{\ast}}-\frac{d}{dt}\left(
\frac{\partial\mathcal{L}}{\partial A_{t}^{\ast}}\right)   & =0
\end{align*}
Since these two equations are complex conjugate of each other we only need one
of them. Find the Euler-Lagrange equations corresponding to the functionals

\begin{description}
\item[a)]
\begin{equation}
I(A,A^{\ast})=\int_{t_{0}}^{t_{0}}dt(A_{t}A_{t}^{\ast}+AA^{\ast}AA^{\ast
})\label{Functional2}%
\end{equation}

\item[b)]
\begin{equation}
I(A,A^{\ast})=\int_{t_{0}}^{t_{0}}dt(\frac{1}{2}i(AA_{t}^{\ast}-A^{\ast}%
A_{t})+\frac{1}{2}AA^{\ast}AA^{\ast})\label{Functional3}%
\end{equation}

\end{description}

\item A functional of the type (\ref{Functional1}) is said to be invariant
under variations of the form%
\begin{align}
A(t) & \mapsto A(t)+\varepsilon\eta(t)\label{Variation1}\\
A^{\ast}(t) & \mapsto A^{\ast}(t)+\varepsilon\eta(t)\nonumber
\end{align}

if there exists a function $F(t)$ such that%
\[
\frac{\partial\mathcal{L}}{\partial A}\eta+\frac{\partial\mathcal{L}}{\partial
A^{\ast}}\eta^{\ast}+\frac{\partial\mathcal{L}}{\partial A_{t}}\eta_{t}%
+\frac{\partial\mathcal{L}}{\partial A_{t}^{\ast}}\eta_{t}^{\ast}=\frac
{dF}{dt}
\]

If $F=0$ we say that the Lagrangian is invariant. Show that the conserved
Noether current corresponding to the variation (\ref{Variation1}) is%
\[
j=F-\eta\frac{\partial\mathcal{L}}{\partial A_{t}}-\eta^{\ast}\frac
{\partial\mathcal{L}}{\partial A_{t}^{\ast}}
\]

Show that the Lagrangian in both functionals (\ref{Functional2}) and
(\ref{Functional3}) are invariant under an infinitesimal rotation of the
complex phase of $A$.%
\begin{align}
A(t)  & \mapsto A(t)+i\varepsilon A(t)\label{Variation2}\\
A^{\ast}(t)  & \mapsto A^{\ast}(t)-i\varepsilon A^{\ast}(t)\nonumber
\end{align}

Find the conserved Noether current corresponding to the variation
(\ref{Variation2}) for the two functionals (\ref{Functional2}) and
(\ref{Functional3}) and show directly using the Euler-Lagrange equations that
the Noether currents are indeed conserved.

\item Consider real valued functionals of the form%
\begin{equation}
I(A,A^{\ast})=\int_{D}dxdy\int_{t_{0}}^{t_{0}}dt\mathcal{L}(t,A,A^{\ast}%
,A_{t},A_{t}^{\ast},A_{x},A_{x}^{\ast},A_{y},A_{y}^{\ast})\label{Functional4}%
\end{equation}

Show that the Euler-Lagrange equations for this functional are%
\begin{align*}
\frac{\partial\mathcal{L}}{\partial A}-\partial_{t}\left(  \frac
{\partial\mathcal{L}}{\partial A_{t}}\right)  -\partial_{x}\left(
\frac{\partial\mathcal{L}}{\partial A_{x}}\right)  -\partial_{y}\left(
\frac{\partial\mathcal{L}}{\partial A_{y}}\right)   & =0\\
\frac{\partial\mathcal{L}}{\partial A^{\ast}}-\partial_{t}\left(
\frac{\partial\mathcal{L}}{\partial A_{t}^{\ast}}\right)  -\partial_{x}\left(
\frac{\partial\mathcal{L}}{\partial A_{x}^{\ast}}\right)  -\partial_{y}\left(
\frac{\partial\mathcal{L}}{\partial A_{y}^{\ast}}\right)   & =0
\end{align*}

Find the Euler-Lagrange equation corresponding to the functional%
\begin{equation}
I(A,A^{\ast})=\int_{D}dxdy\int_{t_{0}}^{t_{0}}dt\left(  A_{t}A_{t}^{\ast
}-c^{2}A_{x}A_{x}^{\ast}-c^{2}A_{y}A_{y}^{\ast}-mAA^{\ast}\right)
\label{Functional5}%
\end{equation}

This equation is called the complex Klein-Gordon equation and describe, among
other things, charged spin-less elementary particles. In this context $m$ is
the mass of the elementary particle.

\item A functional of the type (\ref{Functional4}) is said to be invariant
under variations of the form%
\begin{align}
A(t,x,y)  & \mapsto A(t,x,y)+\varepsilon\eta(t,x,y)\label{Variation3}\\
A^{\ast}(t,x,y)  & \mapsto A^{\ast}(t,x,y)+\varepsilon\eta^{\ast}(t,x,y)\nonumber
\end{align}

if there exists functions $F_{j}=F_{j}(t,x,y)$ for $j=1,2,3$ such that%
\begin{align*}
\frac{\partial\mathcal{L}}{\partial A}\eta+\frac{\partial\mathcal{L}}{\partial
A^{\ast}}\eta^{\ast}+\frac{\partial\mathcal{L}}{\partial A_{t}}\eta_{t}%
+\frac{\partial\mathcal{L}}{\partial A_{t}^{\ast}}\eta_{t}^{\ast}%
+\frac{\partial\mathcal{L}}{\partial A_{x}}\eta_{x}  & \\
+\frac{\partial\mathcal{L}}{\partial A_{x}^{\ast}}\eta_{x}^{\ast}%
+\frac{\partial\mathcal{L}}{\partial A_{y}}\eta_{y}+\frac{\partial\mathcal{L}%
}{\partial A_{y}^{\ast}}\eta_{y}^{\ast}  & =\partial_{t}F_{1}+\partial
_{x}F_{2}+\partial_{y}F_{3}%
\end{align*}

If $F_{1}=F_{2}=F_{3}=0$, we say that the Lagrangian is invariant. Show that
the components of the Noether current corresponding to the variation
(\ref{Variation3}) are%
\begin{align*}
j_{1}  & =F_{1}-\eta\frac{\partial\mathcal{L}}{\partial A_{t}}-\eta^{\ast
}\frac{\partial\mathcal{L}}{\partial A_{t}^{\ast}}\\
j_{2}  & =F_{2}-\eta\frac{\partial\mathcal{L}}{\partial A_{x}}-\eta^{\ast
}\frac{\partial\mathcal{L}}{\partial A_{x}^{\ast}}\\
j_{3}  & =F_{3}-\eta\frac{\partial\mathcal{L}}{\partial A_{y}}-\eta^{\ast
}\frac{\partial\mathcal{L}}{\partial A_{y}^{\ast}}%
\end{align*}

The conservation law for the Noether current is%
\[
\partial_{t}j_{1}+\partial_{x}j_{2}+\partial_{y}j_{3}=0
\]

\item Consider the following two functionals%
\begin{align}
I(\psi,\psi^{\ast})  & =\int_{D}dxdy\int_{t_{0}}^{t_{0}}dt(\frac{\hbar}%
{2}i(\psi^{\ast}\psi_{t}-\psi\psi_{t}^{\ast})\label{Functional6}\\
& -\frac{\hbar^{2}}{2m}(\psi_{x}\psi_{x}^{\ast}+\psi_{y}\psi_{y}^{\ast
})-V(x,y)\psi\psi^{\ast})\\
I(\psi,\psi^{\ast})  & =\int_{D}dxdy\int_{t_{0}}^{t_{0}}dt(\psi_{t}\psi
_{t}^{\ast}-c^{2}\psi_{x}\psi_{x}^{\ast}\label{Functional7}\\
& -c^{2}\psi_{y}\psi_{y}^{\ast}-m\psi\psi^{\ast})
\end{align}
The first functional is discussed on page \pageref{SchrodingerEquation} in the lecture notes on
variational calculus. There we proved that the Euler-Lagrange equation for
this functional is the quantum mechanical Schr\"{o}dinger equation. The second
functional we discussed in problem 10. There we proved that the corresponding
Euler-Lagrange equation is the complex Klein-Gordon equation. Show that the
Lagrangian for the functionals (\ref{Functional6}) and (\ref{Functional7}) are
invariant under an infinitesimal rotation of the complex phase of $\psi$.%
\begin{align}
\psi(t,x,y)  & \mapsto\psi(t,x,y)+i\varepsilon\psi(t,x,y)\label{Variation4}\\
\psi^{\ast}(t,x,y)  & \mapsto\psi^{\ast}(t,x,y)-i\varepsilon\psi^{\ast}(t,x,y)\nonumber
\end{align}

\begin{description}
\item[a)] Find the Noether current and its conservation law corresponding to
the infinitesimal phase variation (\ref{Variation4}) for the Schr\"{o}dinger
functional (\ref{Functional6}).

The Schr\"{o}dinger equation was derived by Erwin Schr\"{o}dinger in
1925. Initially it was not at all clear what the physical interpretation of
the wave function $\psi$ should be. Schr\"{o}dinger himself favored initially
an interpretation in terms of charge density, but could not make it work and
this interpretation was abandoned. The interpretation that lives on to this
day was given by Max Born in 1926. In this interpretation $\psi\psi^{\ast
}(t,x,y)$ is the probability for finding the electron at a point $(x,y)$ at
time $t$. Schr\"{o}dinger never accepted this interpretation, nether did Einstein.

In what way does the conservation law corresponding to the invariance of
the Lagrangian under the infinitesimal phase variation (\ref{Variation4})
support the interpretation introduced by Max Born?

\item[b)] Find the Noether current and its conservation law corresponding to
the infinitesimal phase variation (\ref{Variation4}) for the complex
Klein-Gordon functional (\ref{Functional7}).

The Klein-Gordon equation was in fact first derived by Erwin
Schr\"{o}dinger in 1925 as a quantum mechanical equation for the electron. In
his mind the Klein-Gordon equation was much more likely to be the right
equation for the electron than what we today call the Schr\"{o}dinger
equation. It is for example invariant under Lorentz transformations and thus
respect the fundamental rules of Einsteins special theory of relativity.
Today's Schr\"{o}dinger equation is an approximation to the Klein-Gordon
equation that is valid only for electrons moving slowly compared to the speed
of light, and it is not Lorentz invariant and thus does not respect Einsteins
special theory of relativity. However Schr\"{o}dinger was forced to abandon
his fully relativistic Klein-Gordon equation for the electron because it
turned out to be inconsistent with known atomic spectral data.

 Argue, using the the conservation law corresponding to the invariance of
the Lagrangian under the infinitesimal phase variation (\ref{Variation4}),
that it is not possible to generalize Max Born's interpretation to the
Klein-Gordon equation by constructing a probability density from $\psi$. This
is another reason why it was abandoned at the time.

 It turned out that abandoning the Klein-Gordon equation was premature.
It has since then been reintroduced as a quantum equation, but for charged, 
spin-less particles. In this context the interpretation is nowhere near the
original one introduced by Max Born. Today the wave function in the
Klein-Gordon equation is interpreted as a field of operators that create
particles from the void and return them to the void by annihilation.
\end{description}

One might be surprised at how much insight there is to be gained by using the
simple fact that the global phase of the quantum mechanical wave function is
arbitrary. Thus nothing change if we make the substitution%
\[
\psi\mapsto\psi e^{i\alpha}
\]

The conserved Noether currents corresponding to invariance with respect to
phase, that you have found in problem 12, is however not close to the end of
the story. When this phase invariance is coupled to the Gauge Principle, which
is closely related to Noether's theorem, the existence of the electromagnetic
field and the form of it's interaction with electrons is determined. And there
is more: The quantum mechanical wave equations that describe weakly
interacting particles like neutrinos, and strongly interacting ones like
quarks, also has an invariance with respect to rotation of the global phase.
However, for these equations the arbitrary phase is multidimensional and the
arbitrary rotation of phase involve 2 x 2 matrices for the weakly interacting
case, and 3 x 3 matrices for the strongly interacting case. When the Gauge
principle is applied to these two cases, the existence of the weak interaction
field and the strong interaction field, and their interaction with their
respective particles, like neutrinos and quarks, are determined, just like in
the case of electromagnetics. The field equations for weak interactions and
strong interactions are generalizations of the Maxwell equations involving
more than one vector potential. Also these generalized Maxwell-like equations
are nonlinear, not linear like the electromagnetic Maxwell equation.
\end{enumerate}

\setcounter{equation}{0}

\section{Dimensional analysis}

\subsection{Units and dimensions}

The basic aim of science is to establish functional relationships between physical quantities. 

Physical quantities are used to classify physical objects and events in terms of numbers. Physical quantities, however, are not all the same. \\
\ttx{Base} physical quantities, also called \textit{primary} physical quantities, are defined entirely in terms of physical operations. For such quantities, equality and addition are defined in physical terms. \\
Length is a familiar physical quantity that is primary.
Two sticks are of equal length if they cover each other perfectly when one is put on top of the other. Physical operations corresponding to addition of two lengths, $A$ , $B$ are defined in the familiar way. 

\begin{figure}[htbp]
\centering
\includegraphics{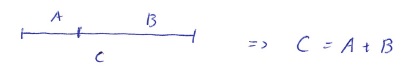}
\caption{Operation defining addition of lengths $A$ and $B$ giving the result $C$}
\label{fig1}
\end{figure}

\noindent Observe that there are no numbers involved in the equality and addition operations. They are defined entirely in physical terms. In order for a physical quantity to be primary the operation of equality and addition must satisfy the following familiar laws 
\begin{align}
A = B \; \text{and} \; B=C  \; \; \Ra \; \; A=C,\nonumber \\ 
A +B = B+A ,\nonumber\\ 
A+(B+C) = (A+B)+C, \nonumber \\ 
\text{If} \; A+B = C \; \; \text{then} \; \not\exists \; D \; \text{such that} \nonumber\\
A + B + D = C. \lbl{1.12}
\end{align}
If \rf{1} holds equality and addition can be used to define the following operations 
\begin{align}
A>B \; \; \;\;\;\;\Lra  &&\exists C\;\;\text{such that}\;\; B + C = A,\nonumber\\
A = C - B \; \; \Lra \; \; A + B = C, \nonumber\\
A = n \; B \; \; \Lra \; \; \underbrace{B+B+...+B}_{\text{n times}} = A,\nonumber\\ 
A = \inv{n} \; B \; \; \Lra \; \; B = n \; A.\lbl{2.12}
\end{align}
Other familiar base quantities are mass, time, area, volume, velocity and force. A base quantity that is perhaps less familiar is \ttx{cardinality} which is a measure of the number of discrete entities in a set of things. 
We use base quantities to assign numbers to objects and events in the familiar way. \\
We first choose a \ttx{unit}. This is a physical object or event displaying a particular instance of the primary quantity in question. \\
We now use this unit as a reference for assigning numerical values to physical objects and events by using \rf{1.12} and \rf{2.12}.

\begin{figure}[htbp]
\centering
\includegraphics{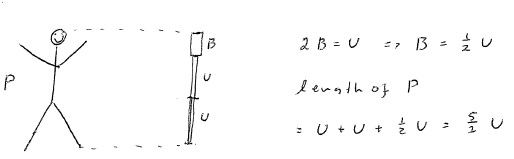}
\caption{A choice of a unit, u, makes it possible to classify objects and events using numbers}
\label{fig2}
\end{figure}

\noindent The numerical value for the length of $P$ is thus 2.5. This process should be familiar all the way back to grade school.

The numerical value assigned obviously depends on which units we use. My height is 
\begin{align*}
l &= 1.7 \; m,  \\ 
l&= 170 \; cm, \nonumber \\
l &= 1700000 \; \mu m. \nonumber 
\end{align*}
There are varying reasons for choosing a particular unit in a given situation. Often the unit is chosen so that the objects and events of interest are assigned numbers of magnitude close to one, if possible. \\
Thus for an astronomer studying nearby stars the natural length unit could be 
\begin{align*}
\text{lightyears} = 9.96.. \cdot 10 ^{15} \; m. 
\end{align*}
Using this unit, the distances to the 26 closest stars range from 4.2 - 11.7. \\ For a scientist studying atoms, a more reasonable unit would be micron or nanometer 
\begin{align*}
\text{micron} &= 1.0 \cdot 10^{-6} \; m,\nonumber\\ 
\text{nanometer} &= 1.0 \cdot 10^{-9}.  
\end{align*}
Using nanometers as our unit we find that most atoms have a diameter of around 0.1. 

Another motivation for choosing units is to simplify algebraic manipulations by getting rid of constant appearing in mathematical formulas. For example if we choose units for length, mass, time, charge and temperature to be 
\begin{align*}
\text{length} &= 1.61 \cdot \tp{-35} \; m,\nonumber\\ 
\text{mass} &=  1.18 \cdot \tp{-8}  \; kg, \nonumber \\
\text{time} &= 5.39\cdot \tp{-44} \; s, \nonumber \\ 
\text{charge} &= 1.88 \cdot \tp{-18} \; C, \nonumber \\ 
\text{temperature} &= 1.42 \cdot \tp{32} \; K.   
\end{align*}
Then the constants 
\begin{align*}
&\text{gravitational constant}  && G,\nonumber\\ 
&\text{Planck constant} && \hbar, \nonumber \\ 
&\text{speed of light} && c, \nonumber \\ 
&\text{Coulomb constant} && \inv{4 \pi \eps_0}, \nonumber \\ 
&\text{Boltzman constant} && K,    
\end{align*}
all get the numerical value 1. This choice leads to enormous simplifications in the algebraic manipulations that are required for predicting events in, for example,  hight energy physics. 

In addition to base, or primary quantities, we have \ttx{derived} quantities. Of these are two types; derived quantities of the first and the  second kind. \\ 
Let us start by discussing derived quantities of the first kind. These quantities appear from inserting numerical values corresponding to base quantities into mathematical formulas. 
\begin{align}
&\text{Base quantities} \begin{cases}
l &= 2 \; m, \\
t& = 60 s .
\end{cases} \nonumber\\ 
&\text{Derived quantities (first kind)} \begin{cases}
A = l \; t = 2\cdot60 = 120, \\ 
B = \frac{1}{2} \; l \; t^2 = \frac{1}{2}\cdot2\cdot 60^2 = 3600.
\end{cases} \lbl{8.12} 
\end{align}
However, not all mathematical formulas will produce a derived physical quantity in this way. \\
Observe that the base quantity, length, has the following property.
\\
If
\begin{align*}
l_1 &= 1 \; m,\nonumber \\ 
l_2 &= 2 \; m, 
\end{align*}
then 
\begin{align*}
\frac{l_1}{l_2} = \frac{1}{2}  
\end{align*}
If we change units to $cm$ we have 
\begin{align*}
l_1 &= 100 \; cm,\nonumber\\ 
l_2 &= 200 \; cm.   
\end{align*}
Both $l_1$ and $l_2$ have changed their numerical values when we introduced the new unit, but we still have 
\begin{align*}
\frac{l_1}{l_2} = \frac{100}{200} = \frac{1}{2}.  
\end{align*}
All base quantities have this property; ratios do not depend on the choice of unit.
Observe that both $A$ and $B$ from \rf{8.12} have the same property.

Let $C$ be defined by 
\begin{align}
C = e^{lt},  \lbl{9.12}
\end{align}
then, using the values of the base quantities $l$ and $t$ from \rf{8.12}, we get 
\begin{align*}
C = e^{2\cdot60} \approx 1.3 \cdot \tp{52}.  
\end{align*}
Thus formula \rf{9.12} certainly assign a numerical value to an event that is assigned numerical values 2 and 60 with respect to the base physical quantities $l$ length and time. However, it does not satisfy the ratio property that $A$ and $B$ and all other base quantities satisfy. 
\\
Let 
\begin{align*}
l_1 &= 1 \; m,\nonumber\\ 
l_2 &= 2 \; m, \nonumber \\ 
t_1 &= 1 \; s, \nonumber \\ 
t_2 &= 2 \; s,  
\end{align*}
Then
\begin{align}
\frac{C_1}{C_2} &= \frac{e}{e^4} \approx 0.049, \lbl{15.1}  
\end{align}
and changing units to $cm$ and minutes 
\begin{align*}
l_1 &= 100 \; cm,\nonumber\\ 
l_2 &= 200 \; cm, \nonumber \\ 
t_1 &= \frac{1}{60} \; min, \nonumber \\ 
t_2 &= 2 \cdot\inv{60} \; min,  
\end{align*}
 we have
 \begin{align*}
 \frac{C_1}{C_2} &= \frac{e^{\frac{100}{60}}}{e^{2 \frac{200}{60}}} \approx 0.0067. 
\end{align*}
P.W. Bridgeman was the first to elevate the invariance of ratios under change of units to a defining property for \ttx{any} physical quantity. \\
Using this property he proved that a mathematical formula
\begin{align*}
y = f(a_1,...,a_n), 
\end{align*}
where the $a_j$ are numerical values corresponding to base quantities, define a physical (derived) quantity \ttx{only} if $f$ is in the form of a monomial. 
\begin{align}
y = c\; a_1^{\alpha_1} \; a_2^{\alpha_2} ... a_n^{\alpha_n} && \alpha_j \in \mathbf{R} , \; \; c \in \mathbf{R} \lbl{18.12}
\end{align}
Thus \ttx{only} formulas of the type \rf{18.12} will define a derived \ttx{physical} quantity. 

In order to keep track of how the numerical values of a derived quantity change when we change units for the base quantities, we introduce the \ttx{dimension} for a derived quantity. 

First, all base physical quantity are assigned a letter chosen by convention. We have for example 
\begin{align*}
\text{length} &\ra L,\nonumber\\ 
\text{time} &\ra T, \nonumber \\ 
\text{mass} &\ra M, \nonumber \\ 
\text{force} &\ra F,   
\end{align*}
Next,  a derived physical quantity is assigned a monomial of letters based on the mathematical formula defining the quantity. We use the notation $[A]$ to denote the dimension of a physical quantity $A$. In general, if $A$ is a physical quantity defined by a monomial
\begin{align*}
A &= c \; a_1^{\alpha_1} \; a_2^{\alpha_2} ...a_n^{\alpha_n},\nonumber
\end{align*}
then it's dimension is given by
\begin{align*}
[A] &= [a_1]^{\alpha_1} \; [a_2]^{\alpha_2} ... [a_n]^{\alpha_n}.    
\end{align*}
We have for example 
\begin{align*}
A &= 2 \; l^2 \; \; \Ra \; [A] = L^2,\nonumber\\ 
V &= l_1 \; l_2 \; l_3 \;\Ra [V] = L^3, \nonumber \\
B &= 3 \; \frac{l}{t} \; \; \Ra [B] = LT^{-1}, \nonumber \\ 
C &= m \; l^2 \; \sqrt{t} \; \; \Ra [C] = M \; L^2 \; T^{\frac{1}{2}}.   
\end{align*}
The dimensions are used to keep track of how the numerical values of derived physical quantities change when we change units for the base quantities. You know how to do this.

\noindent Let an event be characterized by 
\begin{align*}
l&= 1 \; cm,\nonumber\\ 
m &= 1 \; kg, \nonumber \\ 
t &= 1 \; s,   
\end{align*}
and let $A$ be a derived quantity with dimensions 
\begin{align}
[A] = M^{\frac{1}{2}} \; L^2 \; T^3. \lbl{23.12} 
\end{align}
Let us change units for the base quantities to $mm, \; g$ and $hours$. \\ Then
\begin{align*}
l&= 10  \; mm,\nonumber\\ 
m&= 1000\; g, \nonumber \\ 
t &= \inv{3600} \; hours. \nonumber 
\end{align*}
According to \rf{23.12} the numerical value of $A$ will change by a factor 
\begin{align*}
(1000)^{\frac{1}{2}} \; (10)^2 \; (\inv{3600})^3 \approx 6.78 \cdot \tp{-8}.  
\end{align*}

\begin{figure}[htbp]\label{syst1}
\centering
\includegraphics{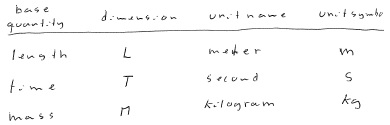}
\caption{System of units for mechanics.\label{fig1_3}}
\end{figure}

\noindent The number of base quantities and the choice of their units depends on what kind of objects and/or events are of interest. For this reason there are many such \ttx{systems of units} in use.\\
If one is mainly interested in mechanical systems one can use the system displayed in figure \ref{fig1_3}.

\noindent Using this system of units, the dimensions of some important derived quantities are
\begin{align}
\text{velocity} &= \dt{x}, && L\; T^{-1},\nonumber\\ 
\text{acceleration} &= \frac{d^2x}{dt^2}, && L \; T^{-2}, \nonumber \\
\text{area} &= \int dx \; dy, && L^2, \nonumber \\ 
\text{Force} &= m \; a, && M \; L \; T^{-2}, \nonumber \\
\text{volume} &= \int dx \; dy \; dz, && L^3 . \lbl{27.12} 
\end{align}
If two systems of units have the same base quantities, but different units, we say that they are of the \ttx{same type}. Thus shifting to units $cm, \; hours$ and $grams$ in the mechanical system of units displayed in figure \ref{fig1_3}, gives us a new system of units for mechanics that is of the same type as \ref{fig1_3}. \\
There is an aspect of the list of physical quantities \rf{27.12} that is somewhat confusing. We know that force is a primary physical quantity, it has the operations of equality and addition defined entirely in terms physical operations. 
But in \rf{27.12} it appears to be a derived quantity! It is defined as mass times acceleration. What is going on, is force a primary quantity or is it a derived quantity? 

The fact of the matter is that the identity 
\begin{align*}
F = m \; a, 
\end{align*}
that appears to tell us that force is a derived quantity of the first kind, is in fact a \ttx{physical law} first discovered by Newton. 
This law has a large, but not universal, domain of validity. It will for example not hold if the speed of objects approaches the speed of light. 

There are systems of units for mechanics that use length, time, mass \ttx{and} force as base quantities. One such system is the British Engineering System 
\begin{figure}[htbp]\label{syst1}
\centering
\includegraphics{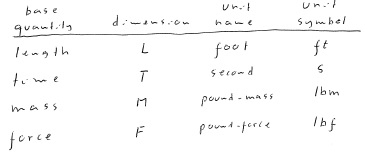}
\caption{The British Engineering System of mechanical units}
\label{fig4}
\end{figure}

\noindent In this system of units Newtons law takes the form 
\begin{align}
F = c \; m \; a, \lbl{30.12}
\end{align}
where $c$ is a \ttx{dimensional constant} with dimensions 
\begin{align*}
[c] = F \; M^{-1} \; L^{-1} \; T^2.
\end{align*}
In the British Engineering system $c$ has the numerical value 
\begin{align*}
c \approx 0.031.  
\end{align*}
By dispensing with force as a base quantity, like in the British Engineering System, and rather defining it as a derived quantity  with dimensions
\begin{align*}
[F] = M \; L \; T^{-2},  
\end{align*}
 as in the mechanical system of units from figure \ref{fig1_3}, we are essentially choosing units for force in such a way as to ensure that the constant of proportionality  in \rf{30.12} is equal to 1. 

\noindent This is certainly a good idea if we expect Newtons law to play a part in our investigations. However if it does not play a part, restricting the choice of units for force, like in the system from figure $\ref{fig1_3}$, is unnecessary, and in fact is detrimental to the utility of dimensional analysis. 

\noindent For example, if we are only concerned with situations where the forces are in balance 
\begin{align*}
\mathlarger{\sum}_i f_i = 0,  
\end{align*}
which is the case in the important subfield of mechanics called {\it statics}, then there is no point in using the system from figure \ref{fig1_3}, no dimensional constant appears if we use a system with force, mass, length and time as base quantities.
This is a good thing; the presence of dimensional constants is also detrimental to dimensional analysis as we will see. 

When primary quantities appear as derived quantities in a system of units through the existence of a physical law, we call them derived quantities of the \ttx{second kind}. \\ 
We can always remove derived quantities of the second kind from a system of units by shifting these quantities to the set of base quantities. We thereby extended the set of base quantities and thus define a new system of units. 
In this new system of units the physical laws, defining the original derived quantities of the second kind, will now include dimensional constants. However, if some of these laws do not play a role in our investigation, like Newtons law for statics, then the corresponding dimensional constant will not appear in our investigations and the power of dimensional analysis is enhanced.

\noindent Recall that volume is in fact a primary physical quantity. Thus,with respect to the mechanical system of units displayed in figure \ref{fig1_3}, volume is a derived quantity of the second kind. \\
The physical law behind this derived quantity is the fact that for a rectangular box with sides of length $l_1, \; l_2, \;l_3$ we have 
\begin{align}
{\cal V}= c \; l_1 \; l_2 \; l_3, \lbl{35.12} 
\end{align}
where $c$ is a dimensional constant whose dimensions are 
\begin{align*}
[c] = \frac{V}{L^3}, 
\end{align*}
and whose numerical value depends on choice of units for volume and length. Here $V=[{\cal V}]$ is by convention the dimension symbol for volume. 

Thus, if we do not use a system of units where volume is a derived quantity, a dimensional constant appears. However \ttx{if} the law \rf{35.12} does not play a part in our investigation, the dimensional constant will not appear and our dimensional analysis will be more powerful, as we will see. \\
Whenever some physical law is deemed relevant for the investigation of some situation, dimensional constants appears.
Two well known dimensional constant are the speed of light, $c$, and the Planck constant, $\hbar$. They appear in the physical laws 
\begin{align*}
E& = c^2  m,\nonumber\\
E&= \hbar \; \omega,   
\end{align*}
which determine the energy equivalence of any given amount of mass and the quantum of energy for electromagnetic radiation of frequency $\omega$.  These laws will play a role if we study systems involving  speeds close to the speed  of light and/or very weak electromagnetic fields. This is the domain covered by the most accurate theory constructed by man, quantum electrodynamics.

  The more such physical laws, that are deemed relevant for an investigation, the more dimensional constants appears and the weaker the dimensional analysis will tend to be. This will become clear when we shortly describe the main tool in dimensional analysis; the PI-theorem.

\begin{figure}[htbp]
\centering
\includegraphics{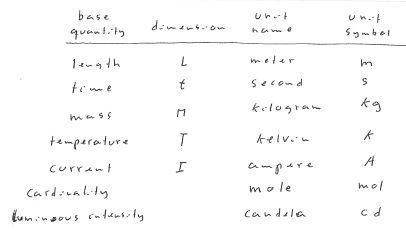}
\caption{The SI system of units.\label{fig1_5}}
\end{figure}

A well known system of units that is used all over the world is the SI - system. Its base quantities and units are displayed in figure \ref{fig1_5}.

Note that there  does not appear to be standardized dimension symbols for cardinality and luminous intensity. 
\\
\subsection{The PI-theorem}
After this review of systems of units, dimensions, base quantities etc, it is time to introduce the main tool of dimensional analysis. This is the famous PI - theorem. 

As stated in the very beginning of these notes, the basic aim of science is to establish functional relationships between physical quantities. Let us assume that a given system can be described by a set of physical quantities 
\begin{align}
Q_1,...,Q_n. \lbl{39.12}
\end{align}
In order to apply dimensional analysis to this situation we \ttx{must} ensure that 
\begin{enumerate}
\item The set of quantities are \ttx{independent}. This means that the values of each quantity, $Q_i$, can be adjusted arbitrarily without affecting the values of the other quantities. For example $Q_1, \; Q_2, \; Q_3$ are not independent if we have, by definition, 
\begin{align}
Q_3 = Q_1 \; Q_2 
\end{align}\label{item1}
\item The set of quantities must be \ttx{complete}. This means that there are no other quantities, beyond $Q_1, ..., Q_n$, that has any significant influence on the situation under study.\label{item2}
\end{enumerate}

\noindent Insuring that properties \ref{item1} and \ref{item2} holds clearly requires some insight into the physics of the situation we want to describe. 

Using $Q_1,...,Q_n$ we can say that the basic aim of science is to find a function $f$ such that
\begin{align}
f(Q_1, ... , Q_n) = 0 \lbl{41.12}
\end{align}
We want $f$ to be a \ttx{physical} relationship between the quantities $Q_1,...,Q_n$ and therefore require that the \ttx{form} of the function $f$ is independent of which units we happened to choose for $Q_1,...,Q_n$. This is a basic assumption of objectivity that applies everywhere in science and which is assumed to hold  in dimensional analysis. \\
For example, let $l$ be a length and $t$ a time and assume units are $m$ and $s$. Lets say that we observe a physical relationship  between these two physical quantities given by
\begin{align}
2\; l^2 + t = 0. \lbl{42.12} 
\end{align}
Define $Q_1 = 2 \; l^2 \; , \; \; Q_2 = t$ and $f(x,y) = x+y$. Then \rf{42.12} can be written as 
\begin{align*}
f(Q_1, Q_2) = 0.  
\end{align*}
However if we rather choose to use units $cm, \; hour$ the relationship we will observe is given by
\begin{align*}
2 \cdot \tp{4} \; l^2 + \inv{3600} \; t = 0.  
\end{align*}
If we define $g(x,y) = \tp{4} \; x + \inv{3600} \; y $ we have 
\begin{align*}
g(Q_1, Q_2) = 0.  
\end{align*}
Since the proposed functional relationship depends on the choice of units we will not accept \rf{42.12} as describing a physical relationship. 

It can however be made into a physical relationship by introducing two dimensional constants $c_1, \; c_2$ where 
\begin{align*}
[c_1] = L^{-2}, && [c_2] = T^{-1}. 
\end{align*}
The values of $c_1, \; c_2$ when units are $m$ and $s$ are 
\begin{align*}
c_1 = 1, && c_2 = 1. 
\end{align*}
We can now define $h(x,y,z,w) = z \; x + w \; y$. We then have a physical functional relationship 
\begin{align*}
h(Q_1, Q_2, c_1, c_2) &= 0,\nonumber\\ 
&\Updownarrow\nonumber\\
2 \; c_1 \; l^2 + c_2 \; t& = 0.   
\end{align*}
In this way, any proposed functional relationship that is not physical can be made physical by introducing enough dimensional constants. 

\noindent So, let us now assume that we have a complete, independent set of physical quantities \rf{39.12}, that satisfies the physical functional relation \rf{41.12}. 

From the set \rf{39.12} we pick a dimensionally \ttx{independent} and \ttx{complete} subset 
\begin{align}
Q_1, ... , Q_r. \lbl{49.12}
\end{align}
That \rf{49.12} is dimensionally complete, means that the dimension any $Q_i$ from \rf{39.12} can be written as a monomial in dimensions of the quantities from \rf{49.12}, and that \rf{49.12} is dimensionally independent means that the dimension of no quantity $Q_s$ in \rf{49.12} can be written as a monomial in the dimensions of the  remaining quantities in \rf{49.12}. \\ 
One might wonder how large $r$ can be?

  Let us assume that the physical quantities $Q_1,...,Q_n$ are expressed in a system of units consisting of $k$ base quantities with dimension symbols given by 
\begin{align*}
d_1, \; d_2, ..., \; d_k. 
\end{align*}
A quantity, $Q$, is said to be \ttx{dimensionless} if its dimensions are 
\begin{align}
[Q] =\Pi_{s=1}^k d_s^0= d_1^0 \; d_2^0 ...d_k^0. \lbl{51.12}
\end{align}
A dimensionless quantity is also called a \ttx{pure number}; its value is the same for all choices of units for the $k$ base quantities. \\
Let the dimensions of $Q_1,...,Q_n$ with respect to the chosen base quantities be 
\begin{align}
[Q_p] =\Pi_{s=1}^kd_s^{\alpha_{sp}}  && \alpha_{sp} \in \mathbf{R}. \lbl{53.12}
\end{align}
If $Q_1,...,Q_r$ are dimensionally independent there can exist no numbers 
\begin{align*}
x_1,...,x_r, 
\end{align*}
where at least one $x_{p_0} \ne 0$ such that 
\begin{align}
C = \Pi^r_{p=1} \; Q_p^{x_p}, \lbl{55.12} 
\end{align}
is dimensionless. This is true because if \rf{55.12} \ttx{did} hold, we would have 
\begin{align*}
Q_{p_0} = C^1 \; \Pi_{p \ne p_0} \; Q_p^{-\frac{x_p}{x_{p_0}}}, && C^1 = C^{-\inv{x_{p_0}}},  
\end{align*}
which would imply that $Q_1,...,Q_r$ are \ttx{not} dimensionally independent. Thus, the largest $r$ such that the quantity $C$  can not be made dimensionless for any choice of constants $x_1,...,x_r$ defines the larges possible value of $r$ such that $Q_1,...,Q_r$ is dimensionally independent.

Inserting the dimensions \rf{53.12} into \rf{55.12}we get 
\begin{align}
[C]=\Pi_{p=1}^r[Q_p]^{x_p}&=\Pi_{p=1}^r\Pi_{s=1}^kd_s^{\alpha_{sp}x_p}= \Pi_{s=1}^k\Pi_{p=1}^rd_s^{\alpha_{sp}x_p} \nonumber\\
&= \Pi_{s=1}^kd_s^{\sum_{p=1}^r\alpha_{sp}x_p}\lbl{50.12}
\end{align}
Thus, according to the definition of a dimensionless quantity \rf{51.12}, we conclude that the largest $r$ such that  $Q_1,...,Q_r$ is dimensionally independent is equal to the largest $r$ such that the homogeneous linear system
\begin{align}
\sum^r_{p=1} &\alpha_{sp} \; x_p = 0 && s = 1,...,k\;\;, \lbl{58.12}
\end{align}
only has the trivial solution.
From the theory of linear systems we know that \rf{58.12} will have non-zero solutions if $r>k$. Thus we conclude that we have the bound 
\begin{align}
r \leq k. \lbl{59.12}
\end{align}
We now return to our main argument. The fact that $Q_1,...,Q_r$ is dimensionally complete means that 
\begin{align}
[Q_{r+i}] = [Q_1^{x_{1(r+i)}} ...Q_r^{x_{r(r+i)}}] && i = 1,...,n-r, \lbl{60.12}
\end{align}
where $x_{1(r+i)}, ... , x_{r(r+i)}$ are some real numbers. \\
Using \rf{60.12} we define dimensionless quantities 
\begin{align}
\Pi_i = \frac{Q_{r+i}}{Q_1^{x_{1(r+i)}}...Q_r^{x_{r(r+i)}}} && i=1,...,n-r. \lbl{61.12} 
\end{align}
Using \rf{41.12}, \rf{60.12} and \rf{61.12} we have 
\begin{align}
f(Q_1,Q_2,...,Q_n)=0, \nonumber \\
\Updownarrow\nonumber\\
f(Q_1,...,Q_r,Q_1^{x_{1(r+1)}}, ..., Q_r^{x_{r(r+1)}}  \Pi_1, ... ,Q_1^{x_{1n}}...Q_r^{x_{rn}} \Pi_{n-r} ) = 0, \nonumber \\ 
\Updownarrow\nonumber\\
g(\Pi_1,...,\Pi_{n-r}, Q_1, ..., Q_r) = 0, \lbl{62.12} 
\end{align}
where $g$ has been defined in terms of $f$ in the obvious way.

  But $f$ and therefore $g$, should not depend on the choice of units of the base quantities. This is only possible if in fact $g$ in \rf{62.12} does \ttx{not} depend on $Q_1,...,Q_r$. 

Thus the conclusion is that \ttx{any} physical relationship, involving $Q_1,...,Q_n$, \ttx{must} be of the form
\begin{align*}
g(\Pi_1, ..., \Pi_{n-r}) = 0, 
\end{align*}
where the $\Pi_j$ are all the independent dimensionless quantities that can be constructed using $Q_1, ... ,Q_n$. \\
This is the PI-theorem. \\
Thus in order to write down all possible physical functional relationships involving $Q_1, ... ,Q_n$ we only need to find the dimensionless quantities $\Pi_j , j=1,...,n-r$. The number of such quantities is according to \rf{59.12} at least $n-k$. For a small number of $\Pi_j$'s, they can usually be constructed easily by manipulating the quantities $Q_1,...,Q_n$. However, for a large number of $\Pi_j$'s there is a systematic procedure that often is useful:\\
We want to find all monomials in $Q_1,..., Q_n$ that are dimensionless. \\ 
Arguing like in \rf{50.12} we have 
\begin{align}
[Q_1]^{x_1} ... [Q_n]^{x_n} &= d_1^0 ... d_k^0, \nonumber \\
\Updownarrow\nonumber\\d_1^{\mathlarger{\sum}^n_{j=1} \alpha_{1j}x_j} ... d_k^{\mathlarger{\sum}^n_{j=1} \alpha_{kj}x_j} &= d_1^0 ...d_k^0, \nonumber \\
\Updownarrow\nonumber\\
\mathlarger{\sum}^n_{j=1} \alpha_{sj}x_j &= 0 && s=1,...,k\;\;. \lbl{64.12}
\end{align}
Thus we only need to find the null-space for the $k\times n$ matrix $\alpha = (\alpha_{sj})$. This can be done using the methods of linear algebra. \\
If $\vb{x}_q = (x_{1q} , ..., x_{nq}),\;\; q =1 ,...,p$ is a basis for the null-space then the corresponding dimensionless quantities are 
\begin{align*}
\Pi_q = Q_1^{x_{1q}} ... Q_n^{x_{nq}} && q=1,...,p\;\;,  
\end{align*}
and where 
\begin{align*}
p \geq n-k.
\end{align*}
Let us now apply dimensional analysis to some simple examples. 

\subsection{Dimensional analysis 1.  No mathematical model is known}
We will first apply dimensional analysis in a situation where no mathematical model in known. This is the most elementary application of dimensional analysis, nothing is required beyond a list of the relevant physical quantities and a choice of a set of base physical quantities. When these choices have been made, the application of dimensional analysis is entirely mechanical. However, the choice of the relevant physical quantities and a set of base quantities is in general anything but elementary and require real insight into the physical systems or processes under investigation.
\begin{example}\label{D1}
Let us consider an object of mass $m$ that is hanging from a fixed point $P$, and is free to swing in one plane under the influence of gravity. \\
The string connecting the mass to the point $P$ is totally stiff and of length $l$. 
\begin{figure}[htbp]
\centering
\includegraphics{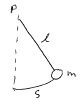}
\caption{A simple pendulum}
\label{fig6}
\end{figure}
We know that in the simple situation alluded to here, the mass moves along an arc of a circle of radius $l$ with center at $P$. We want to figure out how the time $t$, that  it takes the object to move a distance $s$ along the circle, depends on other physical quantities of interest. \\
We know that the force of gravity acting on the mass involves the well known dimensional constant $g$. 

Since this is a mechanical system, and we know that Newtons law must be a part of any modeling of the system, we choose a system of units with base quantities length, time and mass. The dimension symbols are L, T and M. Using our physical insight we conclude that 
\begin{align*}
t,\; l,\; s, \; m, \; g, 
\end{align*}
is a complete and independent system of physical quantities. \\
We need, according to the PI - theorem, dimensionless quantities $\Pi_q, \; q=1,...,p$ where 
\begin{align*}
p \geq 5-3=2.  
\end{align*}
Let us use the formal approach. We then need the matrix $\alpha$ which we list as a table in figure \ref{fig7.12}.

\begin{figure}[htbp]
\centering
\includegraphics{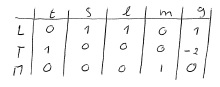}
\caption{Dimension matrix for the simple pendulum.\label{fig7.12}}
\end{figure}

\noindent Thus,using the formal approach outlined in \rf{64.12}, the quantity 
\begin{align*}
\Pi = t^x \; s^y \; l^z \; m^u \; g^v,  
\end{align*}
is dimensionless if and  only if 
\begin{align}
y + z + v &= 0,\nonumber \\
x - 2 \; v &=0, \nonumber \\
u &= 0.\lbl{72.12}
\end{align}
The general solution of \rf{72.12} is 
\begin{align*}
\mqty(x \\ y \\ z \\ u \\ v) = \mqty(2c_1 \\ - c_1 - c_2 \\ c_2 \\ 0 \\ c_1) = c_1 \; \mqty(2 \\ -1 \\ 0 \\ 0 \\ 1 ) + c_2 \; \mqty(0 \\ -1 \\ 1 \\ 0 \\ 0), 
\end{align*}
and therefore a basis for the null-space is
\begin{align*}
 \mqty(2 \\ -1 \\ 0 \\ 0 \\ 1 ) \;,\;\mqty(0 \\ -1 \\ 1 \\ 0 \\ 0). 
 \end{align*}
 Thus we get the following two dimensionless quantities 
\begin{align*}
\Pi_1 &= t^2 \; s^{-1} \; g = \frac{t^2 g}{s},\nonumber\\
\Pi_2 &= s^{-1} \; l = \frac{l}{s}.   
\end{align*}
Since we want to find $t$ as a function of the other variables we write the PI - theorem as 
\begin{align*}
\Pi_1 = f(\Pi_2), 
\end{align*}
where $f$ is an arbitrary function. Note that none of two dimensionless quantities depends on the mass of the pendulum. Thus our guess that the mass did matter in this problem turned out to be unwarranted. Using the expressions for $\Pi_1$ and $\Pi_2$ we have
\begin{align}
\Pi_1 &= f(\Pi_2), \nonumber \\
\Updownarrow\nonumber\\
\frac{t^2g}{s} &= f(\frac{l}{s}), \nonumber \\ 
\Updownarrow\nonumber\\
t^2 &= \frac{s}{g} \; f(\frac{l}{s}), \nonumber \\ 
\Updownarrow\nonumber\\
 t &= \sqrt{ \frac{s}{g} \; f(\frac{l}{s})}. \lbl{76.12}
\end{align}
Recall that by definition 
\begin{align*}
\theta = \frac{s}{l},  
\end{align*}
is the angle between the string and the vertical measured in radians. If we introduce $\theta$ in \rf{76.12} we get 
\begin{align}
t =\sqrt{  \frac{l}{g} \; \theta \; f(\inv{\theta})} = \sqrt{\frac{l}{g}} \; h(\theta). \lbl{77.12} 
\end{align}
Where $h(\theta)$ is an arbitrary function. 
\end{example}
\noindent Thus dimensional analysis tells us exactly how $t$ depends on $g,l$,  and it tells us that $t$ does not depend on the mass $m$. The only thing unknown at this point is how $t$ depends on $\theta$. Dimensional analysis cannot decide this. The function $h(\theta)$ can be found experimentally or through modeling. \\
Since the shape of $h$ does not depend on $l$ and $g$ we can do experiments on a small laboratory sized system and gets results that apply to system on which measurements are impractical. This is of course only true as long as the scaling does not introduce new physical variables beyond $t,s,l,m$ and $g$. This will happen if we scale far enough from our laboratory sized system. No person in his right mind would believe \rf{77.12} if we scale our system to atomic dimensions for example or to  kilometer sized dimensions. 

\begin{example}\label{D2}
Let us next consider small oscillations of a drop of liquid. We assume there is no gravitational field. The oscillations refer to shape change of the drop from spherical to ellipsoidal and back.
\begin{figure}[htbp]
\centering
\includegraphics{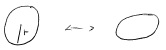}
\caption{An oscillating drop of liquid}
\label{fig8}
\end{figure}
We want to find how the oscillation period depends on the other quantities relevant to this problem. \\ 
Our physical insight informs us that the relevant quantities  are surface tensions  $s$, radius of the drop $r$, and density of the liquid $\rho$. \\
Recall that surface tension is a reaction force resisting the deformation of a liquid surface 
\begin{figure}[htbp]
\centering
\includegraphics{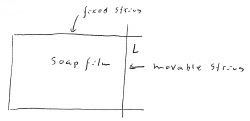}
\caption{Definition of surface tension}
\label{fig9}
\end{figure}

\noindent If we try to move the movable string to the right we have to overcome an elastic force of resistance from the film. The surface tension,$s$,  is a force density defined so that
\begin{align*}
s\vdot L.  
\end{align*} 
is the force acting on a string of length $L$. Thus surface tension is force per unit length. \\
Our dimension matrix, $\alpha$, is then
\begin{figure}[htbp]
\centering
\includegraphics{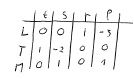}
\caption{Dimension matrix for the oscillating drop of liquid}
\label{fig10}
\end{figure}

\noindent For this case $n=4$ and $k=3$ so the number of dimensionless quantities are $p \geq n - k = 1$. \\
A quantity 
\begin{align*}
\Pi = t^x \; s^y \; r^z \; \rho^u,  
\end{align*}
is dimensionless only if
\begin{align}
z - 3u &= 0 \; \; \Ra \; \; z = 3 \; u,\nonumber \\ 
x - 2 \; y &= 0 \; \; \Ra \; \; x = 2 \; y = - 2 \; u, \nonumber \\
y+u &= 0 \; \; \Ra \; \; y = - u. \lbl{83.12} 
\end{align}
The general solution of \rf{83.12} is thus 
\begin{align*}
\mqty(x \\ y \\ z \\ u) = u \; \mqty(-2 \\ -1 \\ 3 \\1),  
\end{align*}
and a basis for the null-space is 
\begin{align*}
\mqty(-2 \\ -1 \\ 3 \\1).
\end{align*}
The single dimensionless quantity can therefore be chosen to be 
\begin{align*}
\Pi = t^{-2} \; s^{-1} \; r^3 \; \rho,  
\end{align*}
and the PI-theorem informs us that the most general physical relationship involving $s, \; r, \; \rho$ and $t$ is 
\begin{align*}
f(\Pi) &= 0, \nonumber \\ 
\Updownarrow\nonumber\\
\Pi &= c && c \in \mathbf{R} \; \; \; \text{arbitary}, \nonumber \\ 
\Updownarrow\nonumber\\
\frac{r^3 \rho}{s t^2 } &= c, \nonumber \\ 
\Updownarrow\nonumber\\
t &= c' \; \sqrt{\frac{\rho r^3}{s}} && c' \in \mathbf{R}.  
\end{align*}
Thus for this problem, dimensional analysis informed us how $t$ depends on \ttx{all} quantities of relevance to the problem! A single experiment can now determine the number $c'$. 
\end{example}

\begin{example}\label{D3}
We would like to find how the speed, $v$, of ocean surface waves depends on other relevant quantities. \\
In order for the problem not to become too complex, we focus on surface waves on deep water where the bottom topography plays no role. \\
Since surface waves are motions of sea water under the influence of gravity we expect the dimensional constant, $g$, to play a role. Also it is evident that the density of the liquid, $\rho$, must be part of the mix. \\
This appears to be a mechanical problem and since Newtons law evidently must play a role we use a system of units having base quantities, length, time and mass. \\
The dimension matrix $\alpha$ is now
\begin{figure}[htbp]
\centering
\includegraphics{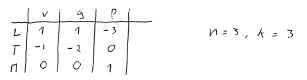}
\caption{Dimension matrix for surface waves, I.}
\label{fig11}
\end{figure}

\noindent There are $p \geq n-k = 0$ dimensionless quantities. 
A quantity
\begin{align*}
\Pi = v^x \; g^y \; \rho^z,  
\end{align*}
is dimensionless only if 
\begin{align*}
x+y -3z &= 0,\nonumber\\ 
-x -2y &= 0, \nonumber \\ 
z &= 0, \nonumber \\ 
\Downarrow\nonumber\\
x=y &= 0.  
\end{align*}
We have no dimensionless quantities and thus no possible physical relations involving $v,g$ and $\rho$. Dimensional analysis fails! \\
Actually dimensional analysis does \textit{not} fail, but we have failed in our application of dimensional analysis. We did not have enough insight into the physics of the situation and therefore did not  include all relevant quantities. In fact, anyone watching ocean waves for a while, will know that the speed of these waves also depends on their wave length, $\lambda$. \\
The dimension matrix is now
\begin{figure}[htbp]
\centering
\includegraphics{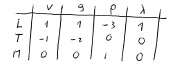}
\caption{Dimension matrix for surface waves, II.}
\label{fig12}
\end{figure}

\noindent A quantity
\begin{align*}
\Pi = v^x \; g^y \;\rho^z \; \lambda^w, 
\end{align*}
is dimensionless if 
\begin{align*}
x + y - 3\; z + w &= 0,\nonumber \\ 
-x - 2\; y &= 0, \nonumber \\ 
z &= 0. 
\end{align*}
The general solution is now 
\begin{align*}
\mqty(x \\ y \\ z \\ w) = c \; \mqty(1 \\ - \frac{1}{2} \\ 0 \\ -\frac{1}{2}).    
\end{align*}
Thus there is one dimensionless quantity 
\begin{align*}
\Pi = v \; g^{- \frac{1}{2}} \; \lambda^{-\frac{1}{2}} = \frac{v}{\sqrt{g \lambda}}, 
\end{align*}
and the only physical law is 
\begin{align*}
f(\Pi) &= 0, \nonumber \\ 
\Updownarrow\nonumber\\
\Pi &= c - \text{constant}, \nonumber \\
\Updownarrow\nonumber\\
v &= c \; \sqrt{\lambda g},  
\end{align*}
The dimensional analysis succeeds in determining wave speed in terms of all quantities of relevance to the problem! The constant, $c$, can be determined by a small-scale laboratory experiment. 
\end{example}

\subsection{Dimensional analysis 2.  A mathematical model is known}
Dimensional analysis, as we have applied it so far, has been physics on the cheap. No mathematical models were written down for the systems of interest. We merely listed the physical quantities relevant to the situation. Modeling only played a role to the extent that we made a decision with regard to what basic physical laws should be involved in such a modeling. This insight guided us in deciding which dimensional constant to include and which system of unit would be appropriate to use in the dimensional analysis. \\
Applying dimensional analysis at this level can be treacherous and requires considerable insight into the physics of the situation in order to succeed. This should be clear from the examples given, especially the last example. \\ 
However, a more common situation is that you, or somebody else, has derived an approximate model for the situation of interest and the challenge is to solve the equations defining the model. In such cases more is known about the system and dimensional analysis is easier to apply. 
\begin{example}\label{D4}
Let us return to the system from example \ref{D1}. The modeling of this system is done in any introductory class in mechanics.
\begin{figure}[htbp]
\centering
\includegraphics{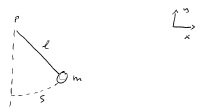}
\caption{Modelling the simple pendulum}
\label{fig13}
\end{figure}

\noindent Using a Cartesian coordinate system with origin at $p$ and axes oriented as indicated we have the model
\begin{align}
x &= l \; \sin \frac{s}{l},\nonumber \\ 
y&=- l \; \cos \frac{s}{l}, \nonumber \\
& s''(t) + g \; \sin (\frac{s(t)}{l}) = 0. \lbl{97.12} 
\end{align}
This model tells us immediately that the relevant physical quantities are 
\begin{align}
t, \; s, \; l, \; g. \lbl{98.12} 
\end{align}
In the original dimensional analysis we included the mass in the list of variables. The dimensional analysis informed us that any physical functional relationship should not include the mass. \\ 
Here we observe that the model does not include the mass and we therefore don't need to include it in our list \rf{98.12}. \\\
The modeling leading up to \rf{97.12} used Newtons law in the form 
\begin{align*}
F = m \; a,  
\end{align*}
with no dimensional constant. We thus have considered force to be a derived quantity of the second kind. We therefore use a system of units with base quantities length and time. \\
The dimension analysis is very simple but let us do it anyway. 
\begin{figure}[htbp]
\centering
\includegraphics{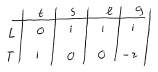}
\caption{Dimension matric for the simple pendulum model}
\label{fig14}
\end{figure}

\noindent A quantity 
\begin{align*}
\Pi = t^x \; s^y \; l^z \; g^w, 
\end{align*}
is dimensionless only if 
\begin{align*}
y + z + w &= 0,\nonumber\\ 
x - 2 \; w &= 0.   
\end{align*}
The general solution is 
\begin{align*}
\mqty(x \\ y \\ z \\ w) = c_1 \; \mqty(2 \\ 0 \\-1 \\ 1) + c_2 \; \mqty(0 \\ 1 \\ -1 \\ 0),     
\end{align*}
and we therefore get two dimensionless quantities 
\begin{align*}
\Pi_1 &= t^2 \; l^{-1} \; g = \frac{gt^2}{l},\nonumber \\ 
\Pi_2 &= s \; l^{-1} = \frac{s}{l}. 
\end{align*}
We want to express $s$ as a function of the other variables and therefore write the general law allowed by the PI - theorem in the form 
\begin{align*}
\Pi_2 &= f(\Pi_1),\nonumber\\
\Updownarrow\nonumber\\
\frac{s}{l} &= f(\frac{g t^2}{l}), \nonumber \\ 
\Updownarrow\nonumber\\
s &= l \; f(\frac{g}{l} \; t^2).  
\end{align*}
This is how far dimensional analysis can take us. However, we also know the actual model for this system, and can use this model to derive an equation for the unknown function $f$. \\
We have 
\begin{align*}
s'&= 2 \; g \; t \; f',\nonumber\\
\Downarrow\nonumber\\
s'' &= 2 \; g \; f' + 4 \; \frac{g^2 t^2}{l} \; f''.   
\end{align*}
Inserting this into \rf{97.12} gives us 
\begin{align}
2 \; g \; f' + 4 \; \frac{g^2 t^2}{l} \; f'' + g \; \sin f = 0. \lbl{108.12} 
\end{align}
Define 
\begin{align*}
\xi\equiv\Pi_1 = \frac{g t^2}{l},  
\end{align*}
then $f = f(\xi)$ and \rf{108.12} can be written in the form 
\begin{align}
4 \; \xi \; f'' + 2 \; f' + \sin f = 0. \lbl{110.12} 
\end{align}
This equation is now in the realm of pure mathematics; no dimensional quantities are involved. \\
The equation \rf{110.12} might look unfamiliar to people that has taken a course in mechanics. For the simple pendulum on would rather expect something like 
\begin{align*}
u'' + \sin u = 0.  
\end{align*}
The form of the equation \rf{110.12} depends on the choice of basis for the null-space of the dimensional matrix. If we rather used the basis 
\begin{align*}
\vb{x}_1 = (1, \; 0, \; - \frac{1}{2}, \frac{1}{2}),\nonumber\\
\vb{x}_2 = (0 , \; 1, \; -1 , \; 0 ),   
\end{align*}
we would get the general physical law in the form 
\begin{align*}
s = l \; h (\sqrt{\frac{g}{l} t}), 
\end{align*}
and inserting this into the model equation \rf{97.12} gives us
\begin{align}
h'' + \sin h = 0. \lbl{114.12} 
\end{align}
Of course, the two equations \rf{110.12} and \rf{114.12} are equivalent. They are connected through the change of variables 
\begin{align*}
f(\xi) = h(\sqrt{\xi}).  
\end{align*}
If the mathematical model for a given physical system is known, we can always do a dimensional analysis and then derive an equation for the unknown function that appears, like we did in example \ref{D4}. The resulting equation will contain no dimensional quantities but will frequently contain dimensionless parameters. 
\end{example}

\begin{example}\label{D5}
Let us include damping, for example caused by viscous effects, in the pendulum model from the previous example.\\
The model equation is now 
\begin{align}
s'' + \frac{\gamma}{m} \; s' + g \; \sin \frac{s}{l} = 0, \lbl{116.12} 
\end{align}
where $\gamma$ is the damping constant. It is a dimensional constant 
\begin{align*}
[\gamma] = \frac{M}{T}. 
\end{align*}
Let us redo the dimensional analysis. Our quantities are now 
\begin{align*}
s, \; t, \; l, \; g, \; m, \; \gamma,  
\end{align*}
and we use a system of units with base quantities length, time, mass. \newline The dimension matrix is displayed in figure \ref{fig15}
\begin{figure}[htbp]
\centering
\includegraphics{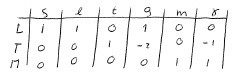}
\caption{Dimension matric for the damped pendulum model.}
\label{fig15}
\end{figure}

\noindent A quantity 
\begin{align*}
\Pi = s^x \; l^y \; t^z \; g^u \; m^v \; \gamma^w, 
\end{align*}
is dimensionless only if 
\begin{align}
x+y+u &= 0,\nonumber\\ 
z - 2 \; u - w &= 0, \nonumber \\ 
v+w &=0.  \lbl{121.12} 
\end{align}
The general solution to this system can be written 
\begin{align*}
\mqty(x \\ y \\ z \\ u \\ v\\ w) = c_1 \; \mqty(1 \\ -1 \\ 0 \\0 \\0\\0) + c_2 \; \mqty(0 \\ -\frac{1}{2} \\1\\ \frac{1}{2} \\ 0 \\ 0) + c_3 \; \mqty(0 \\ \inv{2} \\ 0 \\ - \inv{2} \\ -1 \\ 1),     
\end{align*}
and we therefore get the following three dimensionless quantities 
\begin{align*}
\Pi_1 &= s \; l ^{-1} = \frac{s}{l},\\ 
\Pi_2 &= l^{-\inv{2}} \; t \; g^{\inv{2}} = \sqrt{\frac{g}{l}} \; t, \nonumber \\ 
\Pi_3 &= l^{\inv{2}} \; g^{- \inv{2}} \; m^{-1} \; \gamma = \sqrt{\frac{l}{g}} \; \frac{\gamma}{m},
\end{align*}
and thus the most general physical law can be written 
\begin{align}
\Pi_1 &= f(\Pi_2, \Pi_3),\nonumber\\ 
\Updownarrow\nonumber\\
 s&= l \; f(\sqrt{\frac{g}{l}} \; t, \sqrt{\frac{l}{g}} \; \frac{\gamma}{m}). \lbl{125.12}
\end{align}
We insert \rf{125.12} into \rf{116.12} in order to get an equation for $f$ 
\begin{align*}
\dt{s} &= \sqrt{g \;l} f',\nonumber\\
\Downarrow\nonumber\\
\dst{s} &= g \; f'',  
\end{align*}
and \rf{116.12} implies that 
\begin{align}
g \; f'' + \frac{\gamma}{m} \; \sqrt{g \; l } \; f' + g \; \sin f &= 0, \nonumber \\
\Updownarrow\nonumber\\
\; \; \; f'' + \eps \; f' + \sin f &= 0, \lbl{127.12}
\end{align}
where $f=f(\xi, \eps)$ and the dimensionless variables are 
\begin{align*}
\eps& = \Pi_3, \;\nonumber\\
 \xi &= \Pi_2.  
\end{align*}
\end{example}
\noindent
In order to get \rf{127.12}, we chose a particular basis for the null space of the dimension matrix. 
For the previous example different choices of basis gave different ODE's that were all equivalent through a simple change of variables. 

For the current case we observe that the general solution of \rf{121.12} can also be written in the form
\begin{align}
\mqty(x\\y\\z\\u\\v\\w) = c_1 \; \mqty(1\\-1\\0\\0\\0\\0) + c_2 \; \mqty(0 \\ -\frac{1}{2} \\ 1 \\ \inv{2} \\ 0 \\ 0) + c_3 \; \mqty(0 \\ 0 \\ 1 \\ 0 \\-1 \\ 1), \lbl{129.12}     
\end{align}
and we get the dimensionless quantities 
\begin{align*}
\Pi_1 &= \frac{s}{l},\nonumber \\
\Pi_2 &= \sqrt{\frac{g}{l}} \; t, \nonumber \\ 
\Pi_3 &= \frac{\gamma}{m} \; t.  
\end{align*}
The most general physical law now takes  the form 
\begin{align*}
\Pi_1 &= f(\Pi_2, \Pi_3),\nonumber\\ 
\Updownarrow\nonumber\\
s &= l \; f(\sqrt{\frac{g}{l}} \; t , \frac{\gamma}{m} \; t). \nonumber 
\end{align*}
We now have 
\begin{align*}
\dt{s} &= \sqrt{g \; l} \; \prt{\xi} f + \frac{\gamma \; l}{m} \; \prt{\eta}f,\nonumber\\
\Downarrow\nonumber\\
\dst{s} &= g \; \prt{\xi \xi } f + 2 \; \sqrt{g \; l} \; \frac{\gamma}{m} \; \prt{\xi \eta} f + \frac{\gamma^2 \; l}{m^2} \; \prt{\eta \eta} f, 
\end{align*}
where we have defined
\begin{align*}
\xi &=\Pi_2= \sqrt{\frac{g}{l}} \; t, \nonumber\\
 \eta &= \Pi_3=\frac{\gamma}{m} \; t. \nonumber
 \end{align*}
 \noindent Inserting this into our equation \rf{116.12}, we get 
\begin{align*}
g \; \prt{\xi \xi } f + 2 \; \sqrt{g \; l} \; \frac{\gamma}{m} \; \prt{\xi \eta} f + \frac{\gamma^2 l}{m^2} \prt{\eta \eta} f &\nonumber \\
+ \frac{\gamma}{m} \; \sqrt{g \; l} \; \prt{\xi} f + \frac{\gamma^2 \;l}{m^2} \; \prt{\eta} f + g \; \sin f& = 0, \nonumber \\
& \Updownarrow\nonumber\\
 \prt{\xi \xi} f + 2 \; \sqrt{\frac{l}{g}} \; \frac{\gamma}{m} \; \prt{\xi \eta} f + \frac{\gamma^2 l}{g \; m^2} \; \prt{\eta \eta} f\\
+ \sqrt{\frac{l}{g}} \; \frac{\gamma}{m} \; \prt{\xi} f + \frac{\gamma^2 \; l}{g \; m^2} \; \prt{\eta} f + \sin f& = 0.   
\end{align*}
Defining the dimensionless constant $\eps$ to be 
\begin{align*}
\eps = \sqrt{\frac{l}{g}} \; \frac{\gamma}{m},  
\end{align*}
we get the equation 
\begin{align}
\prt{\xi \xi}f &+ 2 \; \eps \; \prt{\xi \eta} f + \eps^2 \; \prt{\eta \eta} f\nonumber\\
&+ \eps \; \prt{\xi} f + \eps^2 \; \prt{\eta} f + \sin f = 0. \lbl{135.12} 
\end{align}
This is a partial differential equation! Observe that
\begin{align*}
\begin{rcases}
\; \; \xi = \Pi_2\\
\; \; \eta = \Pi_3
\end{rcases} \Ra \frac{\eta}{\xi} = \sqrt{\frac{l}{g}} \; \frac{\gamma}{m} = \eps \Ra \eta = \eps \; \xi. 
\end{align*}
Using the basis \rf{129.12} to the null space we are thus lead to the following solution procedure.
\begin{enumerate}
\item  Find a solution $f(\xi, \eta)$ to  \rf{135.12}.
\item  Define $h(\xi) = f(\xi, \eps \; \xi)$.
\item  A solution to the damped pendulum  problem is $s(t) = l \; h(\sqrt{\frac{g}{l}}t)$.
\end{enumerate}
Since, in general,  it is much harder to solve a PDE than an ODE, the current choice of basis does not appear to be a smart choice. 
However, using partial differential equations to solve linear and non-linear oscillator problems is a key step in an analytic approximation method called \ttx{the method of multiple scales}. This belongs to perturbation methods and will be discussed in section five of these lecture notes.

\subsection{Scaling}
For the case when the mathematical model is known, there is a method which is more restrictive than dimensional analysis, but which is easier to apply. This method, like dimensional analysis, removes all dimensional quantities from the model. \\ 
The approach consists of introducing unknown dimensional constants for all dependent and independent variables, and then to choose values for the dimensional constants in a way that simplifies the equation and ensures that dependent and independent variables vary over a range of order 1, if possible. This approach is called \ttx{scaling}. In this context the dimensional constants are usually called scales.

Let us apply it to the undamped pendulum model from example \ref{D1}. 

\begin{example}\label{D6}
Our model equation is 
\begin{align}
s'' + g \; \sin \frac{s}{l} = 0. \lbl{137.12}
\end{align}
Introduce dimensional constants $c_s, \; c_t$ by 
\begin{align*}
t=c_t \; \tau, && s = c_s \; u,
\end{align*}
where 
\begin{align*}
[c_t] = [t] = T, &&  [c_s]=[s] = L,  
\end{align*}
so that $\tau$ and $u$ are dimensionless. Using the chain rule we have 
\begin{align}
\dt{s} &= c_s \; \dt{} u = \frac{c_s}{c_t} \; \frac{du}{d\tau},\nonumber \\
\Downarrow\nonumber\\
\dst{s} &= \frac{c_s}{c_t^2} \; \frac{d^2u}{d\tau^2}. \lbl{140.12} 
\end{align}
Inserting \rf{140.12} into \rf{137.12} we get 
\begin{align*}
\frac{c_s}{c_t^2} \; u'' + g \; \sin (\frac{c_s}{l} \; u) &= 0,\nonumber\\
&\Updownarrow\nonumber \\
 u'' + \frac{g \; c_t^2}{c_s} \; \sin (\frac{c_s}{l} \; u) &= 0, \nonumber \\ 
&\Updownarrow\nonumber \\
u'' + \eps_1 \; \sin (\eps_2 \; u) &= 0,  
\end{align*}
where $\eps_1, \; \eps_2$ are two dimensionless constants 
\begin{align*}
\eps_1 &= \frac{g \; c_t^2}{c_s},\nonumber\\ 
\eps_2 &= \frac{c_s}{l}.  
\end{align*}
The numerical values for $\eps_1, \; \eps_2$ are determined by the values we choose for $c_s$ and $c_t$.
The values we choose are motivated mainly by what kind of solution we expect will be of importance in the situation of interest. \\
For example, if we expect the pendulum to swing along an arc whose maximum length will be of order $l$,  it makes sense to choose 
\begin{align*}
c_s = l \; \; \; \Ra \; \; \; \eps_2 = 1.  
\end{align*}
Then for swings of this type $u_{max}$ is of order 1. If there is no compelling reason for choosing a particular value for $c_t$, we could choose it so that the dimensionless constant $\eps_1$ also is equal to 1 
\begin{align}
\frac{g \; c_t^2}{c_s} = 1 \; \; \Lra \; \; c_t = \sqrt{\frac{l}{g}}. \lbl{144.12}
\end{align}
We have now \ttx{chosen scales} for length and time and our equation is 
\begin{align*}
u'' + \sin (u) = 0.  
\end{align*}
This is the same dimensionless equation that we got from dimensional analysis.
\end{example}

\begin{example}\label{D7}
Let us assume that our model includes an initial condition
\begin{align*}
s'' + g \; \sin \frac{s}{l} &= 0,  \\
s(0) &= s_0. 
\end{align*}
Assuming that the solution of interest move along an arc of maximum length $s_0$, it makes sense to use $s_0$ as a new scale for length. For time we choose the same scale as in \rf{144.12}. 
\begin{align}
c_s &= s_0,\nonumber\\
c_t &= \sqrt{\frac{l}{g}}. \lbl{148.12} 
\end{align}
Our model is now
\begin{align*}
y'' + \inv{\eps} \; \sin \eps y = 0,\nonumber\\
y(0) = 1,  
\end{align*}
where $y$ is dimensionless and $\eps$ is a dimensionless number. They are defined by the expressions 
\begin{align*}
s &= s_0 \; y,\nonumber\\
\eps &= \frac{s_0}{l}.  
\end{align*}
We now consider the special case 
\begin{align*}
\eps << 1.  
\end{align*}
This condition makes sense because $\eps$ is a pure number whose value does not depend on the choice of units. The solutions of interest are, with respect to the units \rf{148.12}, of maximum size one. Since $\eps$ is small and $y$ of order one, we can use Taylor's theorem and have 
\begin{align*}
\inv{\eps} \; \sin \eps y = y - \inv{6} \; \eps^2 \; y^2 + ... \;\;.
\end{align*}
Thus our model can, for solutions of the assumed type, be approximated by 
\begin{align}
y'' + y = \inv{6} \; \eps^2 \; y^3, && y(0) = 1.\lbl{153.12} 
\end{align}
\end{example}
\noindent
Observe, that the validity of \rf{153.12} depends on the assumption that the maximum value of $s(t)$ is set by $s(0) = s_0$. Here, this is true, but applying approximations to mathematical models based on scaling assumptions can bring us into trouble. We must always be at the lookout for \ttx{breakdowns}. For example if we on the basis of \rf{153.12}, calculate that the solution eventually reach a size $y>>1$, we must reject this solution. There is nothing wrong with the solution as a solution to model \rf{153.12}, but it is not an approximation to a solution of the original model. \\
Our solution, whose size eventually reach $y>>1$, is an \ttx{unphysical} solution even if it is fine as a \ttx{mathematical} solution for \rf{153.12}. If we are not aware of the possibility of breakdowns, like the one I just described, we are not doing \ttx{applied} mathematics, we are rather doing \ttx{pure} mathematics. 

\begin{example}\label{D8}
Let us assume that we have a system that is modeled by the following initial value problem 
\begin{align}
s'(t) &= \alpha \; s^2(t) - \beta \; s^4(t),\nonumber\\
s(0) &= s_0, \lbl{154.12}  
\end{align}
where $s$ and $t$ are dimensional quantities whose dimensions are length and time 
\begin{align*}
[s] = L, && [t] =T,  
\end{align*}
and where $\alpha, \; \beta , s_0$ are positive dimensional constants
\begin{align*}
[s_0] &= L,\nonumber\\
[\alpha] &= T^{-1} \; L^{-1}, \nonumber \\
[\beta] &= T^{-1} \; L^{-3}. 
\end{align*}
We now introduce dimensional constants $c_s, \; c_t$ such that 
\begin{align}
s &= c_s \; y, && [c_s] = L,\nonumber \\
t &= c_t \; \xi, && [c_t] = T. \lbl{157.12}
\end{align}
The quantities $y$ and $\xi$ are then evidently dimensionless. Inserting \rf{157.12} into \rf{154.12} using the chain rule we get 
\begin{align*}
y'(\xi) &= \alpha \; c_s \; c_t \; y^2(\xi) - \beta \; c_s^3 \; c_t \; y^4(\xi),\nonumber \\
y(0) &= \frac{s_0}{c_s}.  
\end{align*}
Choose scales for time and space such that 
\begin{align*}
\alpha \; c_s \; c_t &= 1 \; \; \Lra \; \; c_t = \inv{\alpha \; s_0},\nonumber \\ 
\frac{s_0}{c_s} &= 1 \; \; \Lra \; \; c_s = s_0.  
\end{align*}
The we have 
\begin{align}
y'(\xi) &= y^2(\xi) - \eps \; y^4(\xi),\nonumber\\
y(0) &= 1, \lbl{160.12} 
\end{align}
where the dimensionless constant $\eps$ is 
\begin{align*} 
\eps = \frac{\beta \; s_0}{\alpha}.
\end{align*}
Let us assume that for the given values of the dimensional constants, $\alpha,\beta$ and for some choise of $s_0$, we have
\begin{align*}
\eps << 1. \nonumber
\end{align*}
Then, for a solution of \rf{160.12} whose size is of $O(1)$, we get the approximate model 
\begin{align}
y'(\xi) &= y^2(\xi),\nonumber\\
y(0) &=1. \lbl{162.12} 
\end{align}
This is a separable first order ODE which is easy to solve. The solution is 
\begin{align}
y(\xi) = \inv{1 - \xi}. \lbl{163.12}
\end{align}
\begin{figure}[htbp]
\centering
\includegraphics{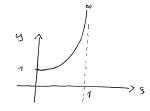}
\caption{Breakdown of solution in finite time}
\label{fig16}
\end{figure}

\noindent  The solution \rf{163.12} approaches infinity at the finite value $\xi = 1$. Thus the solution to the initial value problem \rf{162.12} does not exist beyond $\xi = 1$. However, long before we actually reach the value $\xi = 1$, the solution invalidates the assumption $y(\xi) = O(1)$ which is the justification for approximating \rf{160.12} using \rf{162.12}. Thus long before we have \ttx{mathematical} breakdown of the solution \rf{163.12} we have \ttx{physical} breakdown in the sense that the solution \rf{163.12} is no longer physical relevant. The actual physical solution exists for all $\xi$ and approaches 
\begin{align}
y= \inv{\sqrt{\eps}}, \lbl{164.12}
\end{align}
as $\xi \rightarrow \infty$. \\
When we use scaling to simplify the equations describing a system, breakdown is very common, the solutions we find tend to lose physical relevance in certain domains for the independent variables. In the previous example,  this was easy to see, in more realistic situations it can be everything from simple to  \ttx{very} hard. 
\end{example}

\begin{example}\label{D9}
The equations describing an incompressible non-ideal fluid are 
\begin{align}
\prt{t} \vb{v} + \vb{v} \vdot \grad_{\vb{x}}\vb{v} &= - \inv{\rho_0} \; \grad_{\vb{x}} p + \frac{\eta}{\rho_0} \; \laplacian_{\vb{x}} \vb{v},\nonumber \\
\grad_{\vb{x}} \vdot \vb{v} &= 0. \lbl{165.13} 
\end{align}
Let us introduce dimensional constants for all independent and dependent variables 
\begin{align}
\vb{v} &= c_v \; \vb{u},\nonumber\\ 
\vb{x} &= c_x \; \vb{\xi}, \nonumber \\ 
t &= c_t \; \tau, \nonumber \\
p &= c_p \; q,  \lbl{164}
\end{align}
where now $\vb{u}, \; \vb{\xi}, \; \tau$ and $q$ are dimensionless. Inserting \rf{164} into \rf{165.13} using the chain rule we find after division the system 
\begin{align*}
\prt{\tau} \vb{u} + \eps_1 \; \vb{u} \vdot \grad_{\vb{\xi}} \vb{u} &= - \eps_2 \; \grad_{\vb{\xi}} q + \eps_3 \; \laplacian_{\vb{\xi}} \vb{u},\nonumber\\
\grad_{\vb{\xi}} \vdot \vb{u} &= 0,  
\end{align*}
where the dimensionless constants $\eps_i$ are 
\begin{align*}
\eps_1 &= \frac{c_v \; c_t}{c_x}, \nonumber\\ 
\eps_2 &= \frac{c_p \; c_t}{c_v \; c_x \; \rho_0}, \nonumber \\ 
\eps_3 &= \frac{\eta\;c_t}{\rho_0 \; c_x^2}. \nonumber 
\end{align*}
We keep the velocity and space scales arbitrary but fix the time scale and pressure scale so that 
\begin{align*}
\eps_1 &= 1 \; \; \Lra \; \; \frac{c_v \; c_t}{c_x} =1 \; \; \Lra \; \; c_t = \frac{c_x}{c_v}, \nonumber \\ 
\eps_2 &= 1 \;\; \Lra \; \; \frac{c_p \; c_t}{\rho_0 \; c_v \; c_x} = 1 \; \; \Lra \; \; c_p = \rho_0 \; c_v^2. \nonumber 
\end{align*}
This gives us the system 
\begin{align}
\prt{\tau} \vb{u} + \vb{u} \vdot \grad_{\vb{\xi}} \vb{u} &= -\grad_{\vb{\xi}} q + \inv{\Re} \; \laplacian{}_{\vb{\xi}} \vb{u},\nonumber \\
\grad{_{\vb{\xi}}} \vdot \vb{u} &=0, \lbl{164.12}
\end{align}
where the dimensionless constant $\Re$ is 
\begin{align}
\Re = \frac{c_x \; c_v \; \rho_0}{\eta}. \lbl{165.12} 
\end{align}
\end{example}
\noindent This kind of scaling for the fluid equations was first introduced by G.G. Stokes in 1851. He noted the existence of the dimensionless number \rf{165.12} but the number was named the \textit{Reynolds number}, after Osborn Reynolds who popularized its use in 1883. In a given situation, the value of $\Re$ is set by making an assumption about the dominating fluid velocity and the scale over which it varies in space. When this assumption has been made one expects $\vb{u}$ and $p$ to be of order 1 and vary over time and space scales of order 1. \\
If these expectations holds true in a situation where
\begin{align}
\Re >> 1, \lbl{167.12}
\end{align}
we can approximate \rf{164.12} by 
\begin{align}
\prt{\tau} \vb{u} + \vb{u} \vdot \grad_{\vb{\xi}} \vb{u} &= - \grad_{\vb{\xi}} q,\nonumber \\ 
	\grad_{\vb{\xi}} \vdot \vb{u} &= 0. \lbl{168.12}
\end{align}
This are the Euler equations. Thus the Euler equations are a scaling approximation to the Navier-Stokes equations \rf{164.12}. After having seen example \ref{D8}, we should be wary of a scaling approximations like \rf{168.12}, it may break down in some space and/or time domains. And it does; solutions exist that satisfy natural initial/boundary data but which form singularities in finite time. \\
This does not mean that \rf{168.12} is useless. Many important fluid flows are well described by \rf{168.12}, potential flows is one of them, as we have seen. As I have mentioned in the fluid dynamics section of these lecture notes, the physical validity of the original Navier-Stokes equations under all conditions is still unresolved. However, this is mostly of academic interest, since \rf{164.12} does describe an enormous range of fluid phenomena in a stable way, without breakdown.

What is a fact, however,  is that solutions to \rf{164.12} tend to generate motions on a very large range of space and time scales. This makes them extremely hard to solve in realistic situations. This phenomena is called \ttx{turbulence}, and is an observed physical phenomena, not some mathematical artifact of the equations. 

\subsection{Exercises}

In all these problems I want you to use the formal method for finding the
dimensionless quantities.

\begin{enumerate}
\item
\begin{description}
\item[a)] Let a physical system consist of satellite of negligible mass orbiting
a mass $m$ at a distance $r$. We want to find how the orbital time $t$ of the
satellite depends on the mass and distance. We expect Newton's law of
universal gravitation to play a role here so the gravitational constant $G $
need to be included in our list of physical quantities. Use the formal method
of dimensional analysis to show that
\[
t=\alpha\frac{r^{\frac{3}{2}}}{G^{\frac{1}{2}}m^{\frac{1}{2}}}
\]

\item where $\alpha$ is a dimensionless number.

\item[b)] Now consider the problem of two masses $m_{1}$ and $m_{2}$ orbiting
each other at a distance $r$. Redo the dimensional analysis from a) for this case.
\end{description}

\item In this problem we will use dimensional analysis to derive a formula for
the electromagnetic mass $m$ of a small sphere of radius $r$ and charge $e$.
In the project report I want you to explain what electromagnetic mass actually
is. For this problem use the electrostatic system of units. Figuring out what
this system is, is part of the problem here.

\begin{description}
\item[a)] Show that dimensional analysis fails to find a physical functional
relationship involving only the quantities $m,r$ and $e$.

\item[b)] Since this is some sort of electromagnetic phenomenon it seems likely
that the dimensional constant $c$, which is the speed of light, plays a role.
Redo the dimensional analysis using the quantities $m,r,e$ and $c$. This time
the dimensional analysis succeeds and we get the formula%
\[
m=\alpha\frac{e^{2}}{rc^{2}}
\]

\item where $\alpha$ is a dimensionless constant. Note that I will be reusing
the symbol $\alpha$. The same symbol $\alpha$ in different formulas does not
mean that they denote the same dimensionless constant.
\end{description}

\item Let our system be a box of volume $V$ containing a liquid of density
$\rho$ hanging from a spring with spring constant $k$. The box is acted upon
by the force of gravity and we are required to find a formula for the
oscillation period $t$. \ Since the force of gravity is part of the problem we
expect that the dimensional constant $g$ ,which is the acceleration of
gravity, will be part of the problem.

\begin{description}
\item[a)] Use dimensional analysis based on the mechanical system of units
having base quantities length,time and mass to derive the formula%
\[
t=\sqrt{\frac{V\rho}{k}}f(\frac{k}{V^{\frac{2}{3}}g\rho})
\]

\item where $f(y)$ is an arbitrary function.

\item[b)] Observe that in this problem we are not actually using the physical
law expressing the volume of a box of sides $l_{1},l_{2}$ and $l_{3} $ through
the formula%
\[
V\propto l_{1}l_{2}l_{3}
\]

\item Because of this there is no reason to treat volume as a derived unit of
the second kind but rather treat it as one of the base quantities. Now redo
the dimensional analysis using a system of units with base quantities length,
time,mass and volume. Show that we now get the formula%
\[
t=\alpha\sqrt{\frac{V\rho}{k}}
\]

\item observe that $V\rho$ is the mass $m$ of the liquid in the box.

\item[c)] We consider the same problem as in a) and b) but now we will use a
system of units with base quantities length, time, mass, velocity and force.
Since both the law for moving bodies connecting velocity to length and time
and Newton's law are part of the problem, using velocity and force as base
quantities means that two dimensional constants $k_{1}$ and $k_{2}$ appears.
They are defined through the physical laws%
\begin{align*}
v  & =k_{1}\frac{l}{t}\\
f  & =k_{2}ma
\end{align*}

\item Thus physical quantities for this problem are $t,m,k,k_{1}$ and $k_{2}$.
We use the insight gained from b) to exclude the dimensional constant $g$ from
our list. Show that we now find the formula%
\[
t=\alpha\sqrt{\frac{mk_{1}k_{2}}{k}}
\]

\item which is in fact the same formula as in b) if we use units for velocity
and force such that $k_{1}=k_{2}=1$.
\end{description}

\item What is the formula determining how fast we can walk under ideal conditions?

Let $v$ be the maximum walking speed. We conjecture that $v$ should depend on
length of the walkers legs $l$ the acceleration of gravity $g$ and the mass
of the person. Show that we get the formula%
\[
v=\alpha\sqrt{gl}
\]

The dimensionless constant $\alpha$ can be found by walking tests in the
laboratory. One finds that $\alpha\approx1$. Use the formula to calculate how
fast you can possibly walk. How fast could you walk on the moon or on mars?
It can be argued that age should be relevant in this situation. Include the age as an additional variable and redo the dimensional analysis. Choose a basis for the null space in such a way that you get a law expressing the walking speed in terms of the other physical quantities. Try to approximate the unknown dimensionless function by fitting it to data collected on walking speed

\item Consider a sphere of radius $R$ exposed to an airflow flow of velocity
$v$. The air has density $\rho$ and viscosity $\mu$. The airflow will through
its action on the sphere create a drag force $F$. Use Dimensional analysis to
show that the drag force is given by the formula%
\[
F=\rho R^{2}v^{2}f(\frac{\mu}{Rv\rho})
\]

where $f$ is an arbitrary function that must be determined by mathematical
modelling or experiments.

\item If we put dominoes in a row and topple the first one, it will hit the
next and topple it and so on. We have all seen this on tv. The toppling
action will move through the row of dominoes at a certain speed $v$. It is
evident that this speed should depend on the spacing of the dominoes $d$ their
height $h$ and thickness $t$. Since they topple under the action of gravity we
expect the acceleration of gravity $g$ to be part of the mix. Use dimensional
analysis to show that we have the formula%
\[
v=\sqrt{gh}f(\frac{d}{h},\frac{t}{h})
\]

In the limiting case when $t<<h$ we get
\[
v=\sqrt{gh}f(\frac{d}{h})
\]

Experiments with several kinds of thin dominoes show that the function $f$ is
approximately constant and equal to the dimensionless value $1.5$. Thus we
have
\[
v\approx1.5\sqrt{gh}
\]

A typical domino have a height of $5$ $cm$. This give a toppling speed of
approximately $1$ $m/s$.

\item The first atomic bomb test was performed in New Mexico in 1945. The
energy released by the explosion was top secret at the time, but a series of
timelaps pictures of the exploding cloud was shortly thereafter published in a
popular magazine. This pictures showed that $0.006$ seconds after the
explosion the exploding cloud was approximately spherical of radius $80$
meter. A british physicist G. I. Taylor used this information to calculate
that the explosive yield was approximately $25$ kilo-tons of TNT.

Question: How did he do it?

Answer: He used dimensional analysis.

An explosion occur when a large amount of energy is quickly releases in a
small space. Under ideal conditions the expanding explosive cloud is
spherical. The physical quantities relevant to the situation then appears to
be energy $E$, radius of explosive cloud $R$ at time $t$ and the density of
the surrounding air $\rho$. Perform a dimensional analysis using a system of
units with primary quantities length, time and mass. Show that we get the
formula%
\[
E=\alpha\frac{R^{5}\rho}{t^{2}}
\]

Taylor then used a small scale laboratory test explosion to determine that the
dimensionless number $\alpha$ is approximately equal to one. Now calculate the
energy and show that you get indeed an energy release of approximately $25$ kilo-tons.

\item Consider a situation where two perfectly conducting parallel metal
plates are placed a distance $d$ apart in a perfect vacuum. The plates carry
no electric charge but your experimentalist friend report that he is
nevertheless measuring a negative pressure $p$ between the plates if they are
placed sufficiently close together. Negative pressure means that the plates
are pulled together. His experiment involved plates whose distance was less
than a micron. He is mystified since there is really noting at all between the
plates except a perfect vacuum, so what could be pulling them together? He
wants you to come up with a formula relating the pressure $p$ to the distance
$d$. You have a elementary course in quantum theory and know that a vacuum is
more than we used to think in the pre-quantum days. You know for example that
a quantum harmonic oscillator in its ground state is not stationary, its
position is subject to never ending fluctuations around the mean position which
is zero. You know that these fluctuations goes under the name zero-point
fluctuations and you also know that this funny behavior is generic for
quantum system. This is the content of the Heisenberg uncertainty principle.
Thus you suspect that the emptiness between the plates could be filled by some
quantum system in its ground state and that the negative pressure is caused by
its zero-point fluctuations. Since this something must fill the void between
the plates it is reasonable that it is some kind of field. The most common
such field is the electromagnetic one so you conjecture that the void is
filled by a electromagnetic field performing zero-point fluctuations. Use
dimensional analysis to derive the formula%
\[
p=\alpha\frac{\hslash c}{d^{4}}
\]

Hendrik Casimir used quantum field theory to calculate the constant $a$ in
1948 and found that
\[
\alpha=-\frac{\pi^{2}}{240}
\]

This formula has been verified up to an accuracy equal to a few percent.

\item The Fulling--Davies--Unruh effect shows that an observer experiencing an
acceleration $a$ with respect to empty space will measure his surroundings to
be at a temperature $T>0$. The effect was first described by Stephen Fulling in
1973, Paul Davies in 1975 and W. G. Unruh in 1976. Inspired by the previous
problem, we know that what appears to be empty space, from a classical point of
view, is not actually so from a quantum point of view. We know that the vacuum
is filled by a quantum electromagnetic field in its lowest energy state. Since
this state involve the continual creation and destruction of virtual photons
we perceive the possibility that some of these virtual photons could be
observed to be real by an accelerated observer. We therefore expect that the
speed of light $c$ and the Planck constant $\hslash$ to be part of the mix.
Also since temperature is one of the physical quantities, thermodynamics must
be operating behind the scene here. We therefore expect that the Boltzmann
constant $k_{B}$ will play a role. Note that acceleration is necessary for
this effect to occur. This is because the quantum vacuum is Lorentz invariant;
all observers moving at uniform speed with respect to each other will perceive
the quantum vacuum in exactly the same state. Use dimensional analysis to show
that we have the formula%
\[
T=\alpha\frac{\hslash a}{ck_{B}}
\]

Fulling, Davies and Unruh, pushing the boundaries of known physics, were able
to calculate the dimensionless number $\alpha$. They found%
\[
\alpha=\frac{1}{2\pi}
\]

This effect is strongly related to the Hawkings radiation from a black hole
that was predicted by Stephen Hawking around the same time. Both phenomena
depends on quantum effects and gravitationally effects to be of the same
order. No coherent physical theory involving both quantum effects and
gravitationally effects exists. Finding such a theory is the holy grail of
theoretical physics.

\item A black hole is a space-time singularity. It is the place where all our
physical theories go to die. A black whole is a solution of the Einstein
equations of general relativity. It was found by Karl Schwartschild in 1915.
The black hole is surrounded by a spherical region of space which is called the
event horizon. Anyone crossing this horizon is doomed, he can not escape and
will be crushed by the singularity. Let the event horizon have radius $R$. We
expect $R$ to depend on the mass of the hole $m$ and also the speed of light
$c$. Since this is a gravitational phenomenon we expect that the gravitational
constant $G$ is part of the mix. Use dimensional analysis to show that we have
the formula%
\[
R=\alpha\frac{Gm}{c^{2}}
\]

By solving the equations of general relativity Schwartschild found that
$\alpha=2$.

\item In 1972 Jacob Bekenstein suggested that black holes harbor an amount of  entropy that is proportional to the area of the hole. Stephen Hawking initially opposed the idea. This was because any thermodynamic system with nonzero entropy would be at a nonzero temperature, and all hot objects radiate energy. Hawing pointed out that this simple thermodynamic relationship implied that black holes could not have a nonzero entropy. This is because  a black hole does {\it not}  radiate energy. This was a very reasonable conclusion, after all, any matter or radiation and that cross the horizon of a black hole are sucked into the whole, never to be seen again! However, Bekenstein's argument was difficult to get around and in order to bolster his argument,  Hawking went ahead and calculated  the temperature of a black hole. This was not by any means a simple calculation to do. This is because Hawking realized that the only source of radiation emanating from the hole had to come from particle pairs that are spontaneously created from the vacuum at the horizon of the hole. The existence of such spontaneous creation of particles was a well established consequence of the theory of quantum fields, and the theory of quantum fields is gospel when it comes to describing quantum properties of fields like, for example, the electromagnetic field. A collection of such theories, called the Standard Model was created in the 1970s and has since then correctly predicted all experimental observations of particle interactions. It has never made a wrong prediction during all this time. 

However applying the theory of quantum fields in the presence of gravity meant going beyond known theory. Nobody knew how to combine these two theories into one, in a consistent way, neither does anybody today. However, Hawking, being a genius and all, plowed ahead and got an answer after a very lengthy and difficult calculation. His answer was given by the following beautiful formula
\begin{align}
T=\frac{\hbar c^3}{8\pi G M k_B}.\label{HawkingFormula}
\end{align}
\begin{description}
\item{a)}
From the discussion it is clear that since we are looking for the temperature of a black whole, it's temperature $T$ and mass $M$ must be taken into account.  Since gravity is involved in this problem, the speed of light $c$, and the gravitational constant $G$, must be part of the mix. Since this is a problem in thermodynamics  the Boltzmann constant $k_B$ must be important, as must the Planck constant $\hbar$, since Hawking conjectured that  the source of the radiation came from particle creation, which is described by quantum theory.

Given this, use dimensional analysis to show that the most general physical law connecting all these quantities must be of the form
\begin{align*}
T=\frac{c^2 M}{k_B}f(\frac{G M^2}{c \hbar}),
\end{align*}
where $f(\xi)$ is an arbitrary function of one variable. Find a simple  such function that give us Hawking's formula (\ref{HawkingFormula}) for the radiation from a black hole. 
\item{b)} At the center of our galaxy there is a super massive black hole called Sagittarius A$^*$. Calculate the temperature of this hole.
\item{c)} Since the black hole has a nonzero temperature it will radiate and thereby loose  energy. Through the Einstein relation $E=M c^2$ we can conclude that the whole will loose mass over time. But as the mass decrease, the formula for the temperature of the hole show that the temperature of the hole increase. Thus it radiates more, loose mass faster and increase the temperature even more. This is a runaway effect that will lead to the evaporation of the hole in finite time. If $t_{ev}$ is the time it takes a black hole to evaporate, and assuming that the evaporation time only depends on the speed of light $c$, the Planck constant $\hbar$, the gravitational constant $G$ and the mass of the hole $M$, show that the most general physical law connecting these quantities is
\begin{equation}
t_{ev}=\frac{\hbar^{\frac{1}{2}} G^{\frac{1}{2}}}
{c^{\frac{5}{2}}} f\left(\frac{GM^2}{c\hbar }\right)\label{GeneralEvaporationLaw}
\end{equation}
Using various kinds of arguments and approximations, it has been argued that the formula for the evaporation time is
\begin{equation}
t_{ev}=\frac{5120\pi G^2 M^3}{\hbar c^4}\label{SpecialEvaporationLaw}
\end{equation}
Show that the formula (\ref{SpecialEvaporationLaw}) is consistent with the one we derived using dimensional analysis (\ref{GeneralEvaporationLaw}), by finding a simple function $f$ that give us (\ref{SpecialEvaporationLaw}) starting with (\ref{GeneralEvaporationLaw}). 
\item{d)} Use the evaporation formula (\ref{SpecialEvaporationLaw}) to calculate the time it takes for Sagittarius A$^*$ to evaporate.
\item{e)} The Russian scientists Yakov Borisovich Zel'dovich and Igor Dmitriyevich Novikov proposed in 1966 the existence of {\it primordial } black holes. These are holes that does not come from the collapse of a star but which were created a short time after the Big Bang,  in the radiation dominated era of the universe. These holes can range in mass from $10^{-8} g$ and up. Find a lower limit on the mass of primordial black holes that can  exists at the present time.
 \end{description}
  Because of Hawking contribution to this problem, the thermal radiation from a black hole is called {\it Hawking radiation}.  Stephen Hawking died in 2018 and was buried in Westminster Abbey, an honor given to very few. He share this burial place with giants like Newton and Darwin. Carved into the stone covering his burial cite, you will find the formula for  Hawking radiation. His greatest work. 
\end{enumerate}

\setcounter{equation}{0}
\section{The method of multiple scales}\label{MMS}

Perturbation methods are aimed at finding approximate analytic solutions to
problems whose exact analytic solutions can not be found. The setting where
perturbation methods are applicable is where there is a family of equations,
 $\mathcal{P}(\varepsilon)$, depending on a parameter $\varepsilon<<1$, and
where $\mathcal{P}(0)$ has a known solution. Perturbation methods are designed
to construct solutions to $\mathcal{P}(\varepsilon)$ by adding small
corrections to known solutions of $\mathcal{P}(0)$. The singular aim of
perturbation methods is to calculate corrections to solutions of
$\mathcal{P}(0)$. Perturbation methods do not seek to prove that a solution of
$\mathcal{P}(0)$, with corrections added, is close to a solution of
$\mathcal{P}(\varepsilon)$ for $\varepsilon$ in some finite range with respect
to some measure of error. It's sole aim is to compute corrections and to make
sure that the first correction is small with respect to the chosen solution
of $\mathcal{P}(0)$, that the second correction is small with respect to the
first correction and so on, all in the limit when $\varepsilon$ approaches
zero. This formal nature and limited aim of is why we prefer to call it
\textit{perturbation methods} rather than \textit{perturbation theory}. A
mathematical theory is a description of proven mathematical relations among
well defined objects of human thought. Perturbation methods does not amount to
a mathematical theory in this sense. It is more like a very large bag of
tricks, whose elements have a somewhat vague domain of applicability, and where
the logical relations between the tricks are not altogether clear, to put it nicely.

After all this negative press you might ask why we should bother with this
subject at all, and why we should not rather stay with real, honest to God,
mathematics. The reason is simply this: If you want analytic solutions to
complex problems, it is the only game in town. In fact, for quantum theory,
which is arguably our best description of reality so far, perturbation methods
is almost always the first tool we reach for. For the quantum theory of
fields, like quantum electrodynamics, perturbation methods are essentially the
only tools available. These theories are typically only known in terms of
perturbation expansions. You could say that we don't actually fully  understand the
mathematical  structure of these very fundamental theories. But at the
same time, quantum theory of fields give some of the most accurate,
experimentally verified, predictions in all of science.

So clearly, even if perturbation methods are somewhat lacking in mathematical
justification, they work pretty well. And in the end that is the only thing
that really counts.

These lecture notes are not meant to be a general introduction to the wide
spectrum of perturbation methods that are used  all across science. Many textbooks
exists whose aim is to give such a broad overview, an overview that includes the most
commonly used perturbation methods\cite{Hinch},\cite{Nayfeh},\cite{Holmes}%
,\cite{Murdoc}. Our aim is more limited; we focus on one such method, which 
is widely used in many areas of applied science. This is the \textit{method of
multiple scales}. The method of multiple scales is described in all
respectable books on perturbation methods and there are also more specialized
books on singular perturbation methods where the method of multiple scales has
a prominent place\cite{Kevorkian},\cite{Johnson}. There are, however, quite
different views on how the method is to be applied, and what its limitations
are. Therefore, descriptions of the method appears quite different in the
various sources, depending on the views of the authors. In these lecture notes
we describe the method in a way that is different from most textbooks, but
which is very effective and makes it possible to take the perturbation
expansions to higher order in the small perturbation parameter that would otherwise be possible. The source that
is closest to our approach is \cite{NLO}.

We do not assume that the reader has had any previous exposure to perturbation
methods. These lecture notes therefore starts off by introducing the
basic ideas of asymptotic expansions and illustrate them  using
algebraic equations. The lecture notes then proceeds by introducing regular perturbation expansions for 
single ODEs, study the breakdown of these expansions, and show how to avoid the breakdown using the
method of multiple scales.  The method of multiple scales is then generalized to systems of ODEs, boundary layer problems for ODEs  and to PDEs.
In an appendix to these lecture notes, we further  illustrate the method of multiple scales by applying it to the Maxwell‘s equations; showing how the Nonlinear Schr\o dinger equation appears
as an approximation to the Maxwell equations in a situation where dispersion and nonlinearity balances. 
 Several exercises involving multiple scales for ODEs and PDEs are included in the lecture notes.

\subsection{Regular and singular problems.}

In this section we will introduce perturbation methods in the context of
algebraic equations. One of the main goals of this section is to introduce the
all-important distinction between regular and singular perturbation problems,
but we also use the opportunity to introduce the notion of a
\textit{perturbation hierarchy} and describe some of its general properties.

\subsubsection{A regularly perturbed quadratic equation}

Consider the polynomial equation
\begin{equation}
x^{2}-x+\varepsilon=0.\label{Eq2.1.12}%
\end{equation}
This is our \textit{perturbed problem }$\mathcal{P}(\varepsilon)$. The
\textit{unperturbed problem }$\mathcal{P}(0)$, is
\begin{equation*}
x^{2}-x=0.
\end{equation*}
This unperturbed problem is very easy to solve%
\begin{align*}
x^{2}-x  &  =0,\nonumber\\
&  \Updownarrow\nonumber\\
x_{0}  &  =0,\nonumber\\
x_{1}  &  =1.
\end{align*}
Let us focus on $x_{1}$ and let us assume that the \textit{perturbed problem}
has a solution in the form of a \textit{perturbation expansion}%
\begin{equation}
x(\varepsilon)=a_{0}+\varepsilon a_{1}+\varepsilon^{2}a_{2}+...\;\;.\label{Eq2.4.12}%
\end{equation}
where $a_{0}=1$. Our goal is to find the unknown numbers $a_{1},a_{2},..$\;.
These numbers should have a size of order $1$. This will ensure that
$\varepsilon a_{1}$ is a small correction to $a_{0}$, that $\varepsilon
^{2}a_{2}$ is a small correction to $\varepsilon a_{1}$and so on, all in the
limit of small $\varepsilon$. As we have stressed before, maintaining the
ordering of the perturbation expansion is the one and only unbreakable rule
when we do perturbation calculations. \ The perturbation method now proceeds
by inserting the expansion (\ref{Eq2.4.12}) into equation (\ref{Eq2.1.12}) and
collecting terms containing the same order of $\varepsilon$.%
\begin{gather*}
(a_{0}+\varepsilon a_{1}+\varepsilon^{2}a_{2}+...)^{2}-(a_{0}+\varepsilon
a_{1}+\varepsilon^{2}a_{2}+...)+\varepsilon=0,\nonumber\\
\Downarrow\nonumber\\
a_{0}^{2}+2\varepsilon a_{0}a_{1}+\varepsilon^{2}(a_{1}^{2}+2a_{0}a_{2}%
)-a_{0}-\varepsilon a_{1}-\varepsilon^{2}a_{2}-..+\varepsilon=0,\nonumber\\
\Downarrow\nonumber\\
a_{0}^{2}-a_0+\varepsilon(2a_{0}a_{1}-a_{1}+1)+\varepsilon^{2}(2a_{0}a_{2}%
+a_{1}^{2}-a_{2})+...=0.
\end{gather*}
Since $a_{1},a_{2},..$ are all assumed to be of order $1$ this equation will
hold in the limit when $\varepsilon$ approach zero only if
\begin{align*}
a_{0}^{2}-a_0  &  =0,\nonumber\\
2a_{0}a_{1}-a_{1}+1  &  =0,\nonumber\\
2a_{0}a_{2}+a_{1}^{2}-a_{2}  &  =0.
\end{align*}
We started with one nonlinear equation for $x$, and have ended up with three
coupled nonlinear equations for $a_{0}$, $a_{1}$ and $a_{2}$. Why should we
consider this to be progress? It seems like we have rather substituted one
complicated problem with one that is even more complicated!

The reason why this is progress, is that the coupled system of nonlinear
equations has a very special structure. We can rewrite it in the form%
\begin{align}
a_{0}(a_{0}-1)  &  =0,\nonumber\\
(2a_{0}-1)a_{1}  &  =-1,\nonumber\\
(2a_{0}-1)a_{2}  &  =-a_{1}^{2}.\label{Eq2.7.12}
\end{align}
The first equation is nonlinear but simpler than the perturbed equation
(\ref{Eq2.1.12}), the second equation is linear in the variable $a_{1}$ and 
the third equation is linear in the variable $a_{2}$ when $a_{1}$ has been
found. Moreover, the linear equations are all determined by the same linear
operator\ $\mathcal{L}(\cdot)=(2a_{0}-1)(\cdot)$. This reduction to a simpler
nonlinear equation and a sequence of linear problems determined by the same
linear operator is what makes (\ref{Eq2.7.12}) essentially simpler than the
original equation (\ref{Eq2.1.12}), which does not have this special structure.
The system (\ref{Eq2.7.12}) is called a \textit{perturbation hierarchy} for
(\ref{Eq2.1.12}). The special structure of the perturbation hierarchy is key to
any successful application of perturbation methods, whether it is for
algebraic equations, ordinary differential equations or partial differential equations.

The perturbation hierarchy (\ref{Eq2.7.12}) is easy to solve and we find
\begin{align*}
a_{0}  &  =1,\nonumber\\
a_{1}  &  =-1,\nonumber\\
a_{2}  &  =-1,
\end{align*}
and thus our perturbation expansion to second order in $\varepsilon$ is
\begin{equation*}
x(\varepsilon)=1-\varepsilon-\varepsilon^{2}+...
\end{equation*}
For this simple case we can solve the unperturbed problem directly using the
solution formula for a quadratic equation. Here are some numbers

\bigskip%

\begin{tabular}
[c]{lll}%
$\varepsilon$ & Exact solution & Perturbation solution\\
$0.001$ & $0.998999$ & $0.998999$\\
$0.01$ & $0.989898$ & $09989900$\\
$0.1$ & $0.887298$ & $0.890000$%
\end{tabular}

\bigskip

\noindent We see that our perturbation expansion is quite accurate even for
$\varepsilon$ as large as $0.1$.

Let us see if we can do better by finding an even more accurate approximation through
extension of the perturbation expansion to higher order in $\varepsilon$. In fact
let us take the perturbation expansion to infinite order in $\varepsilon$.
\begin{equation}
x(\varepsilon)=a_{0}+\epsilon a_{1}+\epsilon^{2}a_{2}+...=a_{0}+
{\displaystyle\sum_{n=1}^{\infty}}
\varepsilon^{n}a_{n}\label{Eq2.10.12}
\end{equation}
Inserting (\ref{Eq2.10.12}) into (\ref{Eq2.1.12}) and expanding we get
\begin{gather*}
(a_{0}+
{\displaystyle\sum_{n=1}^{\infty}}
\varepsilon^{n}a_{n})(a_{0}+
{\displaystyle\sum_{m=1}^{\infty}}
\varepsilon^{m}a_{m})-a_{0}-
{\displaystyle\sum_{n=1}^{\infty}}
\varepsilon^{n}a_{n}+\varepsilon=0,\nonumber\\
\Downarrow\nonumber\\
a_{0}^{2}-a_{0}+
{\displaystyle\sum_{p=1}^{\infty}}
\varepsilon^{p}(2a_{0}-1)a_{p}+
{\displaystyle\sum_{p=2}^{\infty}}
\varepsilon^{p}\left(
{\displaystyle\sum_{m=1}^{p-1}}
a_{m}a_{p-m}\right)  +\varepsilon=0,\nonumber\\
\Downarrow\nonumber\\
a_{0}^{2}-a_{0}+\varepsilon\left(  (2a_{0}-1)a_{1}+1\right)  +
{\displaystyle\sum_{p=2}^{\infty}}
\varepsilon^{p}\left(  (2a_{0}-1)a_{p}+
{\displaystyle\sum_{m=1}^{p-1}}
a_{m}a_{p-m}\right)  =0.
\end{gather*}
Therefore the complete perturbation hierarchy is
\begin{align*}
a_{0}(a_{0}-1)  &  =0,\nonumber\\
(2a_{0}-1)a_{1}  &  =-1,\nonumber\\
(2a_{0}-1)a_{p}  &  =-{\displaystyle\sum_{m=1}^{p-1}}a_{m}a_{p-m},\text{ \ \ \ }p\geqq2.
\end{align*}
The right-hand side of the equation for $a_{p}$ only depends on $a_{j}$ for
$j<p$. Thus the perturbation hierarchy is an infinite system of linear
equations that is coupled in such a special way that we can solve them one by
one. The perturbation hierarchy truncated at order $4$ is
\begin{align*}
(2a_{0}-1)a_{1}  &  =-1,\nonumber\\
(2a_{0}-1)a_{2}  &  =-a_{1}^{2},\nonumber\\
(2a_{0}-1)a_{3}  &  =-2a_{1}a_{2},\nonumber\\
(2a_{0}-1)a_{4}  &  =-2a_{1}a_{3}-a_{2}^{2}.
\end{align*}
Using $a_{0}=1$, the solution to the hierarchy is trivially found to be
\begin{align*}
a_{1}  &  =-1,\nonumber\\
a_{2}  &  =-1,\nonumber\\
a_{3}  &  =-2,\nonumber\\
a_{4}  &  =-5.
\end{align*}
For $\varepsilon=0.1$ the perturbation expansion gives
\begin{equation*}
x(0.1)=0.8875...\;\;,
\end{equation*}
whereas the exact solution is
\begin{equation*}
x(0.1)=0.8872...\;\;.
\end{equation*}
we are clearly getting closer. However, we did not get all that much in return
for our added effort.

Of course, we did not actually have to use perturbation methods in order to find
solutions to equation (\ref{Eq2.1.12}), since it is exactly solvable using the
formula for the quadratic equation. The example, however, illustrate many
general features of perturbation calculations that will appear again and again
in different guises.

\subsubsection{A regularly perturbed quintic equation}

Let us consider the equation
\begin{equation}
x^{5}-2x+\varepsilon=0.\label{Eq2.17.12}
\end{equation}
This is our perturbed problem, $\mathcal{P}(\varepsilon)$. For this case
perturbation methods are necessary, since there is no solution formula for
general polynomial equations of order higher than four. The unperturbed
problem, $\mathcal{P}(0)$, is
\begin{equation}
x^{5}-2x=0.\label{Eq2.18.12}
\end{equation}
It is easy to see that the unperturbed equation has a real solution
\begin{equation*}
x=\sqrt[4]{2}\equiv a_{0}.
\end{equation*}
We will now construct a perturbation expansion for a solution to
(\ref{Eq2.17.12}), starting with the solution $x=a_{0}$ of the unperturbed
equation (\ref{Eq2.18.12}). We therefore introduce the expansion
\begin{equation}
x(\varepsilon)=a_{0}+\varepsilon a_{1}+\varepsilon^{2}a_{2}+....\;\;.\label{Eq2.20.12}
\end{equation}
Inserting (\ref{Eq2.20.12}) into equation (\ref{Eq2.17.12}) and expanding we get
\begin{gather*}
(a_{0}+\varepsilon a_{1}+\varepsilon^{2}a_{2}+..)^{5}\nonumber\\
-2(a_{0}+\varepsilon a_{1}+\varepsilon^{2}a_{2}+..)+\varepsilon=0,\nonumber\\
\Downarrow\nonumber\\
a_{0}^{5}+5a_{0}^{4}(\varepsilon a_{1}+\varepsilon^{2}a_{2}+...)+10a_{0}
^{3}(\varepsilon a_{1}+...)^{2}+..\nonumber\\
-2a_{0}-2\varepsilon a_{1}-2\varepsilon^{2}a_{2}-...+\varepsilon=0,\nonumber\\
\Downarrow\nonumber\\
a_{0}^{5}-2a_{0}+\varepsilon(1+5a_{0}^{4}a_{1}-2a_{1})+\varepsilon^{2}
(5a_{0}^{4}a_{2}+10a_{0}^{3}a_{1}^{2}-2a_{2})+...=0.
\end{gather*}
Thus the perturbation hierarchy to order two in $\varepsilon$ is
\begin{align*}
a_{0}^{5}-2a_{0} &  =0,\nonumber\\
(5a_{0}^{4}-2)a_{1} &  =-1, \nonumber\\
(5a_{0}^{4}-2)a_{2} &  =-10a_{0}^{3}a_{1}^{2}.
\end{align*}
Observe that the first equation in the hierarchy for $a_{0}$ is nonlinear, 
whereas the equations for $a_{p}$ are linear in $a_{p}$ for $p>0$. All the
linear equations are defined in terms of the same linear operator
$\mathcal{L}(\cdot)=(5a_{0}^{4}-2)(\cdot)$. This is the same structure that we
saw in the previous example. If the unperturbed problem is linear, the first
equation in the hierarchy will also in general be linear.

  The perturbation hierarchy is easy to solve, and we find
\begin{align*}
a_{1} &  =-\frac{1}{5a_{0}^{4}-2}=\frac{-1}{8},\nonumber\\
a_{2} &  =-\frac{10a_{0}^{3}a_{1}^{2}}{5a_{0}^{4}-2}=-\frac{5\sqrt[4]{8}}{256}.
\end{align*}
The perturbation expansion to second order is then
\begin{equation*}
x(\varepsilon)=\sqrt[4]{2}-\frac{1}{8}\varepsilon-\frac{5\sqrt[4]{8}}
{256}\varepsilon^{2}+...\;\;.
\end{equation*}
Here are some numbers

\bigskip%

\begin{tabular}
[c]{lll}%
$\varepsilon$ & Exact solution & Perturbation solution\\
$0.001$ & $1.18908$ & $1.18908$\\
$0.01$ & $1.19795$ & $1.19795$\\
$0.1$ & $1.17636$ & $1.17638$%
\end{tabular}

\noindent Perturbation expansions for the other solutions to equation (\ref{Eq2.17.12}) can
be found by starting with the other four solutions of the equation
(\ref{Eq2.18.12}). In this way we get perturbation expansions for all the
solutions of (\ref{Eq2.17.12}), and the effort is not much larger than for the
quadratic equation.

If we can find perturbation expansions for all the solutions of a problem
$\mathcal{P}(\varepsilon)$, by starting with solutions of the unperturbed
problem $\mathcal{P}(0)$, we say that $\mathcal{P}(\varepsilon)$ is a
\textit{regular} perturbation of $\mathcal{P}(0)$. \ If the perturbation is
not regular it is said to be \textit{singular}. This distinction applies to
all kinds of perturbation problems whether we are looking at algebraic
equations, ordinary differential equations or partial differential equations.
Clearly, for polynomial equations a necessary condition for being a regular
perturbation problem is that $\mathcal{P}(\varepsilon)$ and $\mathcal{P}(0)$
have the same algebraic order. This is not always the case as the next example shows.

\subsubsection{A singularly perturbed quadratic equation.}

Let us consider the following equation
\begin{equation}
\varepsilon x^{2}+x-1=0.\label{Eq2.28.12}
\end{equation}
This is our perturbed problem $\mathcal{P}(\varepsilon)$. The unperturbed
problem $\mathcal{P}(0)$, is
\begin{equation}
x-1=0.\label{Eq2.29.12}
\end{equation}
There is only one solution to the unperturbed problem
\begin{equation}
x=1\equiv a_{0}.\label{Eq2.30.12}
\end{equation}
Let us find a perturbation expansion for a solution to (\ref{Eq2.28.12}) starting
with the solution (\ref{Eq2.30.12}) of the unperturbed problem
\begin{equation}
x(\varepsilon)=a_{0}+\varepsilon a_{1}+\varepsilon^{2}a_{2}+...\;\;.\label{Eq2.31.12}
\end{equation}
Inserting (\ref{Eq2.31.12}) into equation (\ref{Eq2.28.12}) and expanding we get
\begin{gather*}
\varepsilon(a_{0}+\varepsilon a_{1}+\varepsilon^{2}a_{2}+...)^{2}
+a_{0}+\varepsilon a_{1}+\varepsilon^{2}a_{2}+...-1=0,\nonumber\\
\Downarrow\nonumber\\
\varepsilon(a_{0}^{2}+2\varepsilon a_{0}a_{1}+...)+a_{0}+\varepsilon
a_{1}+\varepsilon^{2}a_{2}+...-1=0,\nonumber\\
\Downarrow\nonumber\\
a_{0}-1+\varepsilon(a_{1}+a_{0}^{2})+\varepsilon^{2}(a_{2}+2a_{0}
a_{1})+...=0.
\end{gather*}
The perturbation hierarchy, up to second order in $\varepsilon$ is thus
\begin{align*}
a_{0}  &  =1,\nonumber\\
a_{1}  &  =-a_{0}^{2},\nonumber\\
a_{2}  &  =-2a_{0}a_{1}.
\end{align*}
The solution of the perturbation hierarchy is
\begin{align*}
a_{0}  &  =1,\nonumber\\
a_{1}  &  =-1,\nonumber\\
a_{2}  &  =2,
\end{align*}
and the perturbation expansion for the solution to (\ref{Eq2.28.12}) starting
from the solution $x=1$ to the unperturbed problem (\ref{Eq2.29.12}) is
\begin{equation*}
x(\varepsilon)=1-\varepsilon+2\varepsilon^{2}+...\;\;.
\end{equation*}
In order to find a perturbation expansion for the other solution to the
quadratic equation (\ref{Eq2.28.12}), the unperturbed problem (\ref{Eq2.29.12}) is
of no help.

However, looking at equation (\ref{Eq2.28.12}) we learn something important: In
order for a solution different from $x=1$ to appear in the limit when
$\varepsilon$ approaches zero, the first term in (\ref{Eq2.28.12}) can not
approach zero. This is only possible if $x$ approaches infinity as
$\varepsilon$ approaches zero.

Inspired by this, let us introduce a change of variables
\begin{equation}
x=\varepsilon^{-p}y,\label{Eq2.36.12}
\end{equation}
where $p>0$. If $y$ is of order one, as $\varepsilon$ approaches zero, then $x $
will approach infinity in this limit, and will thus be the solution we lost in (\ref{Eq2.29.12}
). Inserting (\ref{Eq2.36.12}) into (\ref{Eq2.28.12}) gives us
\begin{align*}
\varepsilon(\varepsilon^{-p}y)^{2}+\varepsilon^{-p}y-1 &  =0,\nonumber\\
&  \Downarrow\nonumber\\
\varepsilon^{1-2p}y^{2}+\varepsilon^{-p}y-1 &  =0,\nonumber\\
&  \Downarrow\nonumber\\
y^{2}+\varepsilon^{p-1}y-\varepsilon^{2p-1} &  =0.
\end{align*}
The idea is now to pick a value for $p$, thereby defining a perturbed problem
$\mathcal{P}(\varepsilon)$, such that $\mathcal{P}(0)$ has a solution of order
one. For $p>1$ we get in the limit when $\varepsilon$ approaches zero the
problem%
\begin{equation*}
y^{2}=0,
\end{equation*}
which does not have any solution of order one. One might be inspired to choose
$p=\frac{1}{2}$. We then get the equation
\begin{equation*}
\sqrt{\varepsilon}y^{2}+y-\sqrt{\varepsilon}=0,
\end{equation*}
which in the limit when $\varepsilon$ approaches zero turns into
\begin{equation*}
y=0.
\end{equation*}
This equation clearly has no solution of order one. Another possibility is to
choose $p=1$. Then we get the equation
\begin{equation}
y^{2}+y-\varepsilon=0.\label{Eq2.41.12}
\end{equation}
In the limit when $\varepsilon$ approaches zero this equation turns into
\begin{equation}
y^{2}+y=0.\label{Eq2.42.12}
\end{equation}
This equation has a solution $y=-1$ which \textit{is} of order one. We
therefore proceed with this choise for $p$, and introduce a perturbation
expansion for the solution to (\ref{Eq2.41.12}) that starts at the solution
$y\equiv a_{0}=-1$ to the unperturbed equation (\ref{Eq2.42.12}).
\begin{equation}
y(\varepsilon)=a_{0}+\varepsilon a_{1}+\varepsilon^{2}a_{2}+...\;\;.\label{Eq2.43.12}
\end{equation}
Inserting the perturbation expansion (\ref{Eq2.43.12}) into equation
(\ref{Eq2.41.12}) and expanding we get
\begin{gather*}
(a_{0}+\varepsilon a_{1}+\varepsilon^{2}a_{2}+...)^{2}+a_{0}+\varepsilon
a_{1}+\varepsilon^{2}a_{2}+...-\varepsilon=0,\nonumber\\
\Downarrow\nonumber\\
a_{0}^{2}+a_{0}+\varepsilon((2a_{0}+1)a_{1}-1)+\varepsilon^{2}((2a
_{0}+1)a_{2}+a_{1}^{2})+...=0.
\end{gather*}
The perturbation hierarchy to second order in $\varepsilon$ is then
\begin{align*}
a_{0}^{2}+a_{0} &  =0,\nonumber\\
(2a_{0}+1)a_{1} &  =1,\nonumber\\
(2a_{0}+1)a_{2} &  =-a_{1}^{2}.
\end{align*}
We observe in passing, that the perturbation hierarchy has the special
structure we have seen earlier. The solution to the perturbation hierarchy is
\begin{align*}
a_{1} &  =-1,\nonumber\\
a_{2} &  =1,
\end{align*}
and the perturbation expansion to second order in $\varepsilon$ is
\begin{equation*}
y(\varepsilon)=-1-\varepsilon+\varepsilon^{2}+...\;\;.
\end{equation*}
Going back to the original coordinate $x$ we finally get
\begin{equation*}
x(\varepsilon)=-\varepsilon^{-1}-1+\varepsilon+...\;\;.
\end{equation*}
Even for $\varepsilon$ as large as $0.1$ the perturbation expansion and the
exact solution, $x_{E}(\varepsilon)$, are close
\begin{align*}
x(\varepsilon) &  =-\varepsilon^{-1}-1+\varepsilon+...\approx
-10.900..\;\;,\nonumber\\
x_{E}(\varepsilon) &  =\frac{-1-\sqrt{1+4\varepsilon}}{2\varepsilon}
\approx-10.916...\;\;.
\end{align*}

The perturbation problem we have discussed in this example is evidently a
singular problem. For singular problems, a coordinate transformation, like the
one defined by (\ref{Eq2.36.12}), must at some point be used to transform the
singular perturbation problem into a regular one.

At this point I need to be honest with you; there is really no general rule
for how to find the right transformations. Skill, experience, insight and
sometimes even dumb luck is needed to succeed. This is one of the reasons why
I prefer to call our subject perturbation methods and not perturbation theory.
Certain classes of commonly occurring singular perturbation problems have
however been studied extensively and rules for finding the correct
transformations have been designed. In general, what one observe, is that some
kind of \textit{scaling transformation}, like in (\ref{Eq2.36.12}), is almost
always part of the mix.

\subsection{Asymptotic sequences and series.}

When using perturbation methods, our main task is to investigate the behavior
of unknown functions $f(\varepsilon)$, in the limit when $\varepsilon$
approaches zero. This is what we did in examples one, two and three.

The way we approach this problem is to compare the unknown function $f\left(
\varepsilon\right)  $ to one or several known functions when $\varepsilon$
approaches zero. In example one and two we compared our unknown function to the known functions
$\{1,\varepsilon,\varepsilon^{2},...\}$ whereas in example three we used the
functions $\{\varepsilon^{-1},1,\varepsilon,...\}$. In order to facilitate
such comparisons, we introduce the "large-O" and "little-o" notation.

\subsubsection{Asymptotic ordering of functions}

\begin{definition}
Let $f(\varepsilon)$ be a function of $\varepsilon$. Then

\begin{description}
\item[i)] $f(\varepsilon)=O(g(\varepsilon))$ $\ ,\,\ \ \,\varepsilon
\rightarrow0$ $\ \Leftrightarrow\lim_{\varepsilon\rightarrow0}\left\vert
\frac{f(\varepsilon)}{g(\varepsilon)}\right\vert \neq0,$

\item[ii)] $f(\varepsilon)=o(g(\varepsilon))\ ,\,\ \ \,\varepsilon
\rightarrow0\ \ \ \ \Leftrightarrow\lim_{\varepsilon\rightarrow0}\left\vert
\frac{f(\varepsilon)}{g(\varepsilon)}\right\vert =0.$
\end{description}
\end{definition}

\noindent Thus,  $\ f(\varepsilon)=O(g(\varepsilon))\,$\ means that $f(\varepsilon)$ and
$g(\varepsilon)$ are of roughly the same size when $\varepsilon$ approaches
zero, and $f(\varepsilon)=o(g(\varepsilon))\,$\ means that $f(\varepsilon)$ is
much smaller than $g(\varepsilon)$ when $\varepsilon$ approaches zero.

We have for example that

\begin{enumerate}
\item $\sin(\varepsilon)=O(\varepsilon)$ $\ ,\,\ \ \,\varepsilon\rightarrow0,
$ \ \ \ because
\[
\lim_{\varepsilon\rightarrow0}\left\vert \frac{\sin(\varepsilon)}{\varepsilon
}\right\vert =1\neq0,
\]

\item $\sin(\varepsilon^{2})=o(\varepsilon)\ ,\,\ \ \,\varepsilon
\rightarrow0,$ \ \ because
\[
\lim_{\varepsilon\rightarrow0}\left\vert \frac{\sin(\varepsilon^{2}
)}{\varepsilon}\right\vert =\lim_{\varepsilon\rightarrow0}\left\vert
\frac{2\varepsilon\cos(\varepsilon^{2})}{1}\right\vert =0,
\]

\item $1-\cos(\varepsilon)=o(\varepsilon),\,\ \ \,\varepsilon\rightarrow0,$
\ \ because
\[
\lim_{\varepsilon\rightarrow0}\left\vert \frac{1-\cos(\varepsilon
)}{\varepsilon}\right\vert =\lim_{\varepsilon\rightarrow0}\left\vert
\frac{\sin(\varepsilon)}{1}\right\vert =0,
\]

\item $\ln(\varepsilon)=o(\varepsilon^{-1}),\,\ \ \,\varepsilon\rightarrow0,$ \ \ because

\[
\lim_{\varepsilon\rightarrow0}\left\vert \frac{\ln(\varepsilon)}
{\varepsilon^{-1}}\right\vert =\lim_{\varepsilon\rightarrow0}\left\vert
\frac{\varepsilon^{-1}}{\varepsilon^{-2}}\right\vert =\lim_{\varepsilon
\rightarrow0}\varepsilon=0.
\]
\end{enumerate}
When we apply perturbation methods, we usually use a whole sequence of
comparison functions. In examples one and two we used the sequence
\[
\{\delta_{n}(\varepsilon)=\varepsilon^{n}\}_{n=1}^{\infty},
\]
and in example three we used the sequence
\[
\{\delta_{n}(\varepsilon)=\varepsilon^{n}\}_{n=-1}^{\infty}.
\]
What is characteristic about these sequences is that
\begin{equation}
\delta_{n+1}(\varepsilon)=o(\delta_{n}(\varepsilon)),\,\ \ \,\varepsilon
\rightarrow0,\label{Eq3.7.12}
\end{equation}
for all $n$ in the range defining the sequences. Sequences of functions that
satisfy conditions (\ref{Eq3.7.12}) are called \textit{asymptotic sequences}.

Here are some asymptotic sequences

\begin{enumerate}
\item $\delta_{n}(\varepsilon)=\sin(\varepsilon)^{n}$,

\item $\delta_{n}(\varepsilon)=\ln(\varepsilon)^{-n}$,

\item $\delta_{n}(\varepsilon)=(\sqrt{\varepsilon})^{n}$.
\end{enumerate}

\noindent Using the notion of asymptotic sequences, we can define asymptotic expansion
analogous to the way infinite series are defined in elementary calculus

\begin{definition}
Let $\{\delta_{n}(\varepsilon)\}$ be an asymptotic sequence. Then a formal
series
\begin{equation*}
{\displaystyle\sum_{n=1}^{\infty}}a_{n}\delta_{n}(\varepsilon),
\end{equation*}

\end{definition}

is an asymptotic expansion for a function $f(\varepsilon)$ as $\varepsilon$
approaches zero if
\begin{equation*}
f(\varepsilon)-{\displaystyle\sum_{n=1}^{N}}a_{n}\delta_{n}(\varepsilon)=o(\delta_{N}(\varepsilon)),\text{ \ \ \ }\varepsilon\rightarrow0.
\end{equation*}
Observe that
\begin{align*}
f(\varepsilon)-a_{1}\delta_{1}(\varepsilon)  &  =o(\delta_{1}(\varepsilon
)),\text{ \ \ \ }\varepsilon\rightarrow0,\nonumber\\
&  \Downarrow\nonumber\\
\lim_{\varepsilon\rightarrow0}\left\vert \frac{f(\varepsilon)-a_{1}\delta
_{1}(\varepsilon)}{\delta_{1}(\varepsilon)}\right\vert  &  =0,\nonumber\\
&  \Downarrow\nonumber\\
\lim_{\varepsilon\rightarrow0}\left\vert a_{1}-\frac{f(\varepsilon)}
{\delta_{1}(\varepsilon)}\right\vert  &  =0,\nonumber\\
&  \Downarrow\nonumber\\
a_{1}  &  =\lim_{\varepsilon\rightarrow0}\frac{f(\varepsilon)}{\delta
_{1}(\varepsilon)}.
\end{align*}
In an entirely similar way we find that for all $m\geqq1$ that%
\begin{equation}
a_{m}=\lim_{\varepsilon\rightarrow0}\left\vert \frac{f(\varepsilon)-{\displaystyle\sum_{n=1}^{m-1}}a_{n}\delta_{n}(\varepsilon)}{\delta_{m}(\varepsilon)}\right\vert.\label{Eq3.13.12}
\end{equation}
This shows that for a fixed asymptotic sequence, the coefficients of the
asymptotic expansion for a function $f(\varepsilon)$ are determined by taking
limits. Observe that formula (\ref{Eq3.13.12}) does not require
differentiability for $f(\varepsilon)$ at $\varepsilon=0$. This is very
different from Taylor expansions which requires that $f(\varepsilon)$ is
infinitely differentiable at $\varepsilon=0$.

\noindent  This is a hint that asymptotic expansions are much more general than the usual
convergent expansions, for example power series, that we are familiar with from
elementary calculus. In fact, asymptotic expansions may well diverge, but this
does not make them less useful! The following example was first discussed by
Leonard Euler in 1754.

\subsubsection{Euler's example}

Let a function $f(\varepsilon)$ be defined by the formula
\begin{equation}
f(\varepsilon)=\int_{0}^{\infty}dt\frac{e^{-t}}{1+\varepsilon t}.\label{Eq3.14.12}
\end{equation}
The integral defining $f(\varepsilon)$ converge very fast, and because of this
$f(\varepsilon)$ is a very smooth function, in fact it is infinitely smooth
and moreover analytic in the complex plane where the negative real axis has been removed.

Using the properties of telescoping series we observe that for all $m\geqq0$
\begin{equation}
\frac{1}{1+\varepsilon t}={\displaystyle\sum_{n=0}^{m}}(-\varepsilon t)^{n}+\frac{(-\varepsilon t)^{m+1}}{1+\varepsilon t}.\label{Eq3.15.12}
\end{equation}
Inserting (\ref{Eq3.15.12}) into (\ref{Eq3.14.12}) we find that
\begin{equation*}
f(\varepsilon)=S_{m}(\varepsilon)+R_{m}(\varepsilon),
\end{equation*}
where
\begin{align}
S_{m}(\varepsilon) &  =
{\displaystyle\sum_{n=0}^{m}}(-1)^{n}n!\;\varepsilon^{n},\nonumber
\end{align}
For the quantity $R_{m}(\varepsilon)$ we have the estimate
\begin{equation*}
\left\vert R_{m}(\varepsilon)\right\vert \leqq\varepsilon^{m+1}\int_{0}^{\infty}dt\frac{t^{m+1}e^{-t}}{1+\varepsilon t}\leqq\varepsilon^{m+1}
\int_{0}^{\infty}dt\smallskip t^{m+1}e^{-t}=(m+1)!\varepsilon^{m+1},
\end{equation*}
from which it follows that
\begin{equation*}
\lim_{\varepsilon\rightarrow0}\left\vert \frac{R_{m}(\varepsilon)}{\varepsilon^{m}}\right\vert \leqq\lim_{\varepsilon\rightarrow0}(m+1)!\varepsilon=0.
\end{equation*}
Thus we have proved that an asymptotic expansion for $f(\varepsilon)$ is
\begin{equation}
f(\varepsilon)={\displaystyle\sum_{n=0}^{\infty}}(-1)^{n}n!\varepsilon^{n}.\label{Eq3.20.12}
\end{equation}
It is on the other hand trivial to verify that the formal power series
\begin{equation*}
{\displaystyle\sum_{n=0}^{\infty}}(-1)^{n}n!\varepsilon^{n},
\end{equation*}
diverge for all $\varepsilon\neq0$!

\begin{figure}[htbp]
\centering
\includegraphics[
height=2.501in,
width=3.8147in
]{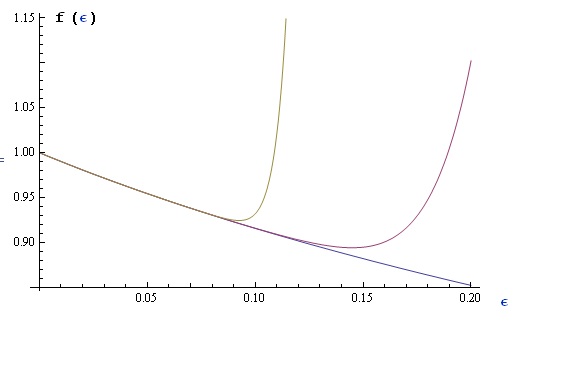}
\caption{Comparing the exact(blue) expression for $f(\varepsilon)$ with the
asymptotic expansion (\ref{Eq3.20.12}) containing ten(red) and twenty(yellow)
terms.\label{fig4_1}}
\end{figure}

\noindent In figure \ref{fig4_1}, we compare the function $f(\varepsilon)$ with what we get from the
asymptotic expansion for a range of $\varepsilon$ and several truncation
levels for the expansion. From this example we make the following two
observations that are quite generic with regards to the convergence or
divergence of asymptotic expansions.%

\noindent Firstly, the asymptotic expansion (\ref{Eq3.20.12}) is an accurate representation
of $f(\varepsilon)$ in the limit when $\varepsilon$ approaches zero even if
the expansions is divergent. Secondly, adding more terms to the expansion for
a fixed value of $\varepsilon$ makes the expansion less accurate.

In reality we are most of the time, because of algebraic complexity, only
able to calculate a few terms of an asymptotic expansion. Thus convergence
properties of the expansion are most of the time unknown. As this example
shows, convergence properties are also not relevant for what we are trying to
achieve when we solve problems using perturbation methods.

\subsection{Regular perturbation expansions for ODEs.}

It is now finally time to start solving differential equations using
asymptotic expansions. Let us start with a simple boundary value problem for a
first order ordinary differential equation.

\subsubsection{A weakly nonlinear boundary value problem.}\label{NonlinearBoundaryValueProblem}

Consider the following boundary value problem
\begin{align}
y^{\prime}(x)+y(x)+\varepsilon y^{2}(x)  &  =x,\text{ \ }0<x<1,\nonumber\\
y(1)  &  =1,\label{Eq3.22.12}
\end{align}
where $\varepsilon$ as usual is a small number. Since the differential
equation is nonlinear and non-separable, this is a nontrivial problem. The
unperturbed problem is
\begin{align*}
y^{\prime}(x)+y(x)  &  =x, \text{ \ }0<x<1,\nonumber\\
y(1)  &  =1.
\end{align*}
The unperturbed problem is easy to solve since the equation is a first order
linear equation. The general solution to the equation is
\begin{equation*}
y(x)=x-1+Ae^{-x}.
\end{equation*}
The arbitrary constant $A$ is determined from the boundary condition
\begin{align*}
y(1)  &  =1,\nonumber\\
&  \Downarrow\nonumber\\
1-1+Ae^{-1}  &  =1,\nonumber\\
&  \Downarrow\nonumber\\
A  &  =e.
\end{align*}
Thus the unique solution to the unperturbed problem is
\begin{equation}
y_{0}(x)=x-1+e^{1-x}.\label{Eq3.26.12}
\end{equation}
We now want to find an asymptotic expansion for the solution to the perturbed
problem (\ref{Eq3.22.12}), starting from the solution $y_{0}(x)$. We thus
postulate an expansion of the form
\begin{equation}
y(\varepsilon;x)=y_{0}(x)+\varepsilon y_{1}(x)+\varepsilon^{2}y_{2}
(x)+...\;\;.\label{Eq3.27.12}
\end{equation}
Inserting (\ref{Eq3.27.12}) into (\ref{Eq3.22.12}) and expanding we get
\begin{gather}
(y_{0}+\varepsilon y_{1}+\varepsilon^{2}y_{2}+...)^{\prime}+y_{0}+\varepsilon
y_{1}+\varepsilon^{2}y_{2}+...\nonumber\\
+\varepsilon(y_{0}+\varepsilon y_{1}+\varepsilon^{2}y_{2}+...)^{2}=x,\nonumber\\
\Downarrow\nonumber\\
y_{0}^{\prime}+\varepsilon y_{1}^{\prime}+\varepsilon^{2}y_{2}^{\prime
}+...+y_{0}+\varepsilon y_{1}+\varepsilon^{2}y_{2}+...\nonumber\\
\varepsilon(y_{0}^{2}+2\varepsilon y_{0}y_{1}+..)=x,\nonumber\\
\Downarrow\nonumber\\
y_{0}^{\prime}+y_{0}+\varepsilon(y_{1}^{\prime}+y_{1}+y_{0}^{2})+\varepsilon
^{2}(y_{2}^{\prime}+y_{2}+2y_{0}y_{1})+...=x.\label{Eq3.28.12}
\end{gather}
We must also expand the boundary condition
\begin{equation}
y_{0}(1)+\varepsilon y_{1}(1)+\varepsilon^{2}y_{2}(1)+...=1.\label{Eq3.29.12}
\end{equation}
From (\ref{Eq3.28.12}) and (\ref{Eq3.29.12}) we get the following perturbation
hierarchy
\begin{align*}
y_{0}^{\prime}(x)+y_{0}(x)  &  =x,\nonumber\\
y_{0}(1)  &  =1,\nonumber\\
& \nonumber\\
y_{1}^{\prime}(x)+y_{1}(x)  &  =-y_{0}^{2}(x),\nonumber\\
y_{1}(1)  &  =0,\nonumber\\
& \nonumber\\
y_{2}^{\prime}(x)+y_{2}(x)  &  =-2y_{0}(x)y_{1}(x),\nonumber\\
y_{2}(1)  &  =0.
\end{align*}
We observe that the perturbation hierarchy has the special structure that we
have noted earlier. All equations in the hierarchy are determined by the
linear operator $\mathcal{L}=\frac{d}{dx}+1$. The first boundary value problem
in the hierarchy has already been solved. The second equation in the hierarchy
is
\begin{equation}
y_{1}^{\prime}(x)+y_{1}(x)=-y_{0}^{2}(x).\label{Eq3.31.12}
\end{equation}
Finding a special solution to this equation is simple
\begin{gather*}
y_1^{p\;\prime}(x)+y_1^{p}(x)=-y_{0}^{2}(x),\nonumber\\
\Downarrow\nonumber\\
(y_1^{p}(x)e^{x})^{\prime}=-y_{0}^{2}(x)e^{x},\nonumber\\
\Downarrow\nonumber\\
y_1^{p}(x)=-e^{-x}\int_{0}^{x}dx^{\prime}e^{x^{\prime}}y_{0}^{2}(x^{\prime
}).
\end{gather*}
Adding a general solution to the homogeneous equation, we get the general solution to
equation (\ref{Eq3.31.12}) in the form
\begin{equation}
y_{1}(x)=A_{1}e^{-x}-e^{-x}\int_{0}^{x}dx^{\prime}e^{x^{\prime}}y_{0}^{2}(x^{\prime
}).\label{Eq3.32.1.12}
\end{equation}
Inserting the expression for $y_0(x)$ from (\ref{Eq3.26.12}) into (\ref{Eq3.32.1.12}), expanding and doing the required integrals, we find that, after applying the boundary condition $y_1(1)=0$, we have 
\begin{equation*}
y_{1}(x)=  -x^{2}+4x-5+(2x-x^{2})e^{1-x}+e^{2-2x}.
\end{equation*}
 We can thus conclude that the perturbation expansion to first order in
$\varepsilon$ is
\begin{equation}
y(\varepsilon;x)=x-1+e^{1-x}+\varepsilon\left(  -x^{2}+4x-5+(2x-x^{2}
)e^{1-x}+e^{2-2x}\right)  +...\;\;.\label{Eq3.34.12}
\end{equation}
The general solution to the third equation in the perturbation hierarchy is in
a similar way found to be
\begin{equation}
y_{2}(x)=A_{2}e^{-x}-2e^{-x}\int_{0}^{x}dx^{\prime}e^{x^{\prime}}
y_{0}(x^{\prime})y_{1}(x^{\prime}).\label{Eq3.35.12}
\end{equation}
 The integral in (\ref{Eq3.35.12})
will have eighteen terms that needs to be integrated. We thus see that even for this very simple
example the algebraic complexity grows quickly. \ 

Recall that we are only ensured that the correction $\varepsilon y_{1}(t)$ is
small with respect to the unperturbed solution $y_{0}(t)$ in the limit when
$\varepsilon$ approaches zero. The perturbation method does not say anything
about the accuracy for any finite value of $\varepsilon$. The hope is of
course that the perturbation expansion also gives a good approximation for
some range of $\varepsilon>0$.

\begin{figure}[htbp]
\centering
\includegraphics[
height=2.501in,
width=3.7144in
]{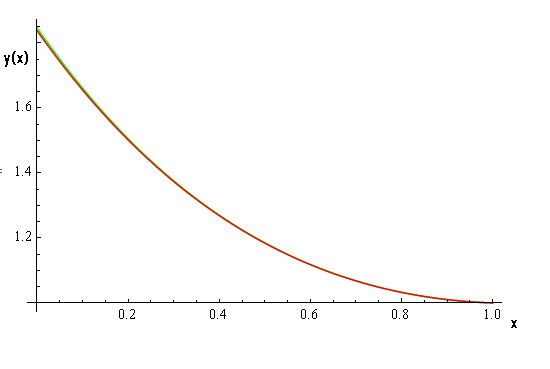}
\caption{Comparing direct perturbation expansion(red) to the exact solution(green),of the boundary value problem.\label{fig4_2}}
\end{figure}

Our original equation (\ref{Eq3.22.12}) is a Riccati equation and an exact solution to the boundary value problem can be found in terms of Airy functions.
In figure \ref{fig4_2} we compare our perturbation expansion(\ref{Eq3.34.12}) to the exact solution in the domain $0<x<1$. We observe that even for
$\varepsilon$ as large as $0.05$ our perturbation expansion give a very
accurate representation of the solution over the whole domain.

\noindent In general, we will not have an exact solution that can be used to
investigate the accuracy of the perturbation expansion for finite values of
$\varepsilon$. For example, if our original equation contained $y^3$ instead of 
$y^2$, an exact solution can not be found. This is the normal situation when we apply perturbation
methods. The only way to get at the accuracy of the perturbation expansion is
to compare it to an approximate solution found by some other, independent,
approximation scheme. Often this involve numerical methods, but it could also
be another perturbation method.

As the next example show, things does not always work out as well as in the current example.

\subsubsection{A weakly damped linear oscillator.}\label{WeaklyDampedLinearOscillator}

Consider the following initial value problem
\begin{align}
y^{\prime\prime}(t)+\varepsilon y^{\prime}(t)+y(t)  &  =0,\text{
\ }t>0,\nonumber\\
y(0)  &  =1,\nonumber\\
y^{\prime}(0)  &  =0.\label{Eq3.36.12}
\end{align}
This is our perturbed problem $\mathcal{P}(\varepsilon)$. The unperturbed
problem, $\mathcal{P}(0)$, is
\begin{align*}
y^{\prime\prime}(t)+y(t)  &  =0,\nonumber\\
y(0)  &  =1,\nonumber\\
y^{\prime}(0)  &  =0.
\end{align*}
The general solution to the unperturbed equation is evidently
\begin{equation*}
y_{0}(t)=A_{0}e^{it}+A_{0}^{\ast}e^{-it},
\end{equation*}
and the initial condition is satisfied if
\begin{align*}
A_{0}+A_{0}^{\ast}  &  =1,\nonumber\\
iA_{0}-iA_{0}^{\ast}  &  =0,
\end{align*}
which has the unique solution $A_{0}=\frac{1}{2}$. Thus the unique solution to
the unperturbed problem is
\begin{equation}
y_{0}(t)=\frac{1}{2}e^{it}+(\ast),\label{Eq3.40.12}
\end{equation}
where $z+(\ast)$ means $z+z^{\ast}$. This is a very common notation.

We want to find a perturbation expansion for the solution to the perturbed
problem, starting with the solution $y_{0}$ of the unperturbed problem. The
simplest approach is to use an expansion of the form
\begin{equation}
y(\varepsilon;t)=y_{0}(t)+\varepsilon y_{1}(t)+\varepsilon^{2}y_{2}(t)...\;\;.\label{Eq3.41.12}
\end{equation}
We now, as usual, insert (\ref{Eq3.41.12}) into the perturbed equation
(\ref{Eq3.36.12}) and expand
\begin{gather}
(y_{0}+\varepsilon y_{1}+\varepsilon^{2}y_{2}+...)^{\prime\prime}\nonumber\\
+\varepsilon(y_{0}+\varepsilon y_{1}+\varepsilon^{2}y_{2}+...)^{\prime}
+y_{0}+\varepsilon y_{1}+\varepsilon^{2}y_{2}+...=0,\nonumber\\
\Downarrow\nonumber\\
y_{0}^{\prime\prime}+y_{0}+\varepsilon(y_{1}^{\prime\prime}+y_{1}
+y_{0}^{\prime})+\varepsilon^{2}(y_{2}^{\prime\prime}+y_{2}+y_{1}^{\prime
})+...=0.\label{Eq3.42.12}
\end{gather}
We must in a similar way expand the initial conditions
\begin{align}
y_{0}(0)+\varepsilon y_{1}(0)+\varepsilon^{2}y_{2}(0)+...  &  =1,\nonumber
\\
y_{0}^{\prime}(0)+\varepsilon y_{1}^{\prime}(0)+\varepsilon^{2}y_{2}^{\prime
}(t)+...  &  =0.\label{Eq3.43.12}
\end{align}
From equations (\ref{Eq3.42.12}) and (\ref{Eq3.43.12}) we get the following
perturbation hierarchy%
\begin{align*}
y_{0}^{\prime\prime}+y_{0}  &  =0,\text{ \ }t>0,\nonumber\\
y_{0}(0)  &  =1,\nonumber\\
y_{0}^{\prime}(0)  &  =0,\nonumber\\
& \nonumber\\
y_{1}^{\prime\prime}+y_{1}  &  =-y_{0}^{\prime},\text{ \ }t>0,\nonumber\\
y_{1}(0)  &  =0,\nonumber\\
y_{1}^{\prime}(0)  &  =0,\nonumber\\
& \nonumber\\
y_{2}^{\prime\prime}+y_{2}  &  =-y_{1}^{\prime},\text{ \ }t>0,\nonumber\\
y_{2}(0)  &  =0,\nonumber\\
y_{2}^{\prime}(0)  &  =0.
\end{align*}
We note that the perturbation hierarchy has the special form discussed
earlier. Here the linear operator determining the hierarchy is $L=\frac{d^{2}
}{dt^{2}}+1$.

The first initial value problem in the hierarchy has already been solved. The
solution is (\ref{Eq3.40.12}). Inserting $y_{0}(t)$ into the second equation in
the hierarchy we get
\begin{equation}
y_{1}^{\prime\prime}+y_{1}=-\frac{i}{2}e^{it}+(\ast).\label{Eq3.45.12}
\end{equation}
Looking for particular solutions of the form
\[
y_{1}^{p}(t)=Ce^{it}+(\ast),
\]
will not work, here because the right-hand side of (\ref{Eq3.45.12}) is a solution
to the homogeneous equation. In fact (\ref{Eq3.45.12}) is a harmonic oscillator
driven on resonance. For such cases we must rather look for a special
solution of the form
\begin{equation}
y_{1}^{p}(t)=Cte^{it}+(\ast).\label{Eq3.46.12}
\end{equation}
By inserting (\ref{Eq3.46.12}) into (\ref{Eq3.45.12}) we find $C=-\frac{1}{4}$. The
general solution to equation (\ref{Eq3.45.12}) is then
\begin{equation*}
y_{1}(t)=A_{1}e^{it}-\frac{1}{4}te^{it}+(\ast).
\end{equation*}
Applying the initial condition for $y_{1}(t)$ we easily find that
$A_{1}=-\frac{i}{4}$. Thus the perturbation expansion to first order in
$\varepsilon$ is
\begin{equation*}
y(\varepsilon;t)=\frac{1}{2}e^{it}+\varepsilon\frac{1}{4}(i-t)e^{it}+(\ast).
\end{equation*}
Let $y_{E}(t)$ be a high precision numerical solution to the perturbed problem
(\ref{Eq3.36.12}). For $\varepsilon=0.01$ we get for increasing time

\bigskip

$
\begin{array}
[c]{ccc}
t & y_{E} & y\\
4 & -0.6444 & -0.6367\\
40 & -0.5426 & -0.5372\\
400 & -0.0722 & 0.5295
\end{array}
$

\bigskip

\begin{figure}[ptb]
\centering
\includegraphics[
height=3.4765in,
width=4.8801in
]{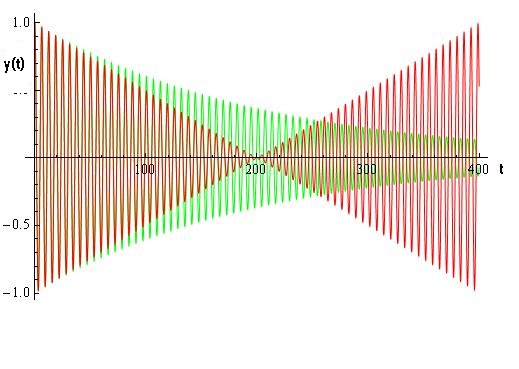}
\caption{Comparing the direct perturbation expansion(red) and a high precision
numerical solution(green).\label{fig4_3}}
\end{figure}

\noindent The solution starts out by being quite accurate, but as $t$ increases, the
perturbation expansion eventually looses any relation to the exact solution.
The true extent of the disaster is seen in figure \ref{fig4_3}.

So what is going on, why is the perturbation expansion such a bad
approximation in this example?

Observe that $y_{1}$ contain a term that is proportional to $t$. Thus as $t$
grows the size of $y_{1}$ also grows and when
\begin{equation*}
t\sim\frac{1}{\varepsilon}
\end{equation*}
the second term in the perturbation expansion become as large as the first
term. The ordering of the expansion breaks down and the first correction,
$\varepsilon y_{1}$, is of the same size as the solution to the unperturbed
problem, $y_{0}$.

The reason why the growing term, $y_{1}$, is a problem here, but was not a
problem in the previous example, is that here the domain for the independent
variable is unbounded.

\noindent Let us at this point introduce some standard terminology. The last two
examples involved perturbation expansions where the coefficients depended on a
parameter. In general such expansions takes the form
\begin{equation*}
f(\varepsilon;\mathbf{x})\sim{\displaystyle\sum_{n=1}^{\infty}}
a_{n}(\mathbf{x})\delta_{n}(\varepsilon),\text{ \ \ }\varepsilon\rightarrow 0
\end{equation*}
where the parameter, $\mathbf{x}$, ranges over some domain $V\subset
\mathbb{R}^{m}$ for some $m$. For the boundary value problem (\ref{Eq3.22.12}) , $V$ is the interval $[0,1]$ whereas for
the initial value problem (\ref{Eq3.36.12}), $V$ is the unbounded interval $(0,\infty)$.

With the introduction of a parameter dependence of the coefficients, a
breakdown of order in the expansion for some region(s) in $V$ becomes a
possibility. We saw how this came about for the case of the damped harmonic oscillator model (\ref{Eq3.36.12}).

And let me be clear about this; breakdown of order in parameter dependent
perturbation expansions is not some weird, rarely occurring, event. On the
contrary it is very common.

Thus methods has to be invented to handle this phenomenon, which is called
\textit{non-uniformity} of asymptotic expansions. The multiple scale method is
design to do exactly this.

\subsection{The method of multiple scales for ODE.}

In the previous section we saw that trying to represent the solution to the
problem
\begin{align}
y^{\prime\prime}(t)+\varepsilon y^{\prime}(t)+y(t)  &  =0,\text{
\ }t>0,\nonumber\\
y(0)  &  =1,\nonumber\\
y^{\prime}(0)  &  =0,\label{Eq4.1.12}
\end{align}
using a regular perturbation expansion
\begin{equation*}
y(\varepsilon;t)=y_{0}(t)+\varepsilon y_{1}(t)+\varepsilon^{2}y_{2}(t)...\;\;,
\end{equation*}
leads to a nonuniform expansion where ordering of the terms broke down for
$t\sim\frac{1}{\varepsilon}$. In order to understand how to fix this, let us
have a look at the exact solution to (\ref{Eq4.1.12}). The exact solution can be
found using characteristic polynomials. We get
\begin{equation}
y(t)=Ce^{-\frac{1}{2}\varepsilon t}e^{i\sqrt{1-\frac{1}{4}\varepsilon^{2}}t}+(\ast),\label{Eq4.3.12}
\end{equation}
where
\begin{equation*}
C=\frac{-\lambda^{\ast}}{\lambda-\lambda^{\ast}},\text{ \ \ }\lambda=-\frac
{1}{2}\varepsilon+i\sqrt{1-\frac{1}{4}\varepsilon^{2}}.
\end{equation*}
If we expand the square root in the exponent with respect to $\varepsilon$, we
get
\begin{equation}
y(t)\approx Ce^{it}e^{-\frac{1}{2}\varepsilon t}e^{-\frac{i}{8}\varepsilon
^{2}t}+(\ast).\label{Eq4.5.12}
\end{equation}
Observe that if $f(\xi)$ is a function whose derivative is of order one, then
the function
\begin{equation*}
g_{n}(t)=f(\varepsilon^{n}t),
\end{equation*}
satisfy
\begin{equation*}
\bigtriangleup g_{n}(t)=g_{n}(t+T)-g_{n}(t)\approx\varepsilon^{n}f^{\prime
}(\varepsilon^{n}t)T=O(1)\Longleftrightarrow T_n\sim\varepsilon^{-n}.
\end{equation*}
We express this by saying that the function $g_{n}(t)$ \textit{vary on the
time scale} $T_{n}=\varepsilon^{-n}$. If we now look at equation
(\ref{Eq4.5.12}), we see that the approximate solution (\ref{Eq4.5.12}) vary on three
separate time scales $T_{0}=\varepsilon^{0},T_{1}=\varepsilon^{-1}$ and
$T_{2}=\varepsilon^{-2}$. If we include more terms in the Taylor expansion
for the square root in (\ref{Eq4.3.12}) the resulting solution will depend on
even more time scales.

Inspired by this example we postulate the existence of a function
\begin{equation*}
h=h(t_{0},t_{1},t_{2},...),
\end{equation*}
such that%
\begin{equation}
y(t)=h(t_{0},t_{1},t_{2},...)|_{t_{j}=\varepsilon^{j}t},\label{Eq4.9.12}
\end{equation}
is a solution to problem (\ref{Eq4.1.12}). Using the chain rule we evidently
have
\[
\frac{dy}{dt}(t)=\left\{  (\partial_{t_{0}}+\varepsilon\partial_{t_{1}
}+\varepsilon^{2}\partial_{t_{2}}+...)h\right\}  |_{t_{j}=\varepsilon^{j}t},
\]
which we formally write as
\begin{equation}
\frac{d}{dt}=\partial_{t_{0}}+\varepsilon\partial_{t_{1}}+\varepsilon
^{2}\partial_{t_{2}}+...\;\;.\label{Eq4.10.12}
\end{equation}
The function $h$ is represented using a perturbation expansion of the form
\begin{equation}
h=h_{0}+\varepsilon h_{1}+\varepsilon^{2}h_{2}+...\;\;.\label{Eq4.11.12}
\end{equation}
The multiple scale method now proceed by substituting (\ref{Eq4.10.12}) and
(\ref{Eq4.11.12}) into the differential equation
\begin{equation*}
y^{\prime\prime}(t)+\varepsilon y^{\prime}(t)+y(t)=0,
\end{equation*}
and expanding everything in sight.
\begin{gather*}
(\partial_{t_{0}}+\varepsilon\partial_{t_{1}}+\varepsilon^{2}\partial_{t_{2}
}+...)(\partial_{t_{0}}+\varepsilon\partial_{t_{1}}+\varepsilon^{2}
\partial_{t_{2}}+...)\nonumber\\
(h_{0}+\varepsilon h_{1}+\varepsilon^{2}h_{2}+...)+\varepsilon(\partial
_{t_{0}}+\varepsilon\partial_{t_{1}}+\varepsilon^{2}\partial_{t_{2}
}+...)\nonumber\\
(h_{0}+\varepsilon h_{1}+\varepsilon^{2}h_{2}+...)+h_{0}+\varepsilon
h_{1}+\varepsilon^{2}h_{2}+...=0,\nonumber\\
\Downarrow\nonumber\\
(\partial_{t_{0}t_{0}}+\varepsilon(\partial_{t_{0}t_{1}}+\partial_{t_{1}t_{0}
})+\varepsilon^{2}(\partial_{t_{0}t_{2}}+\partial_{t_{1}t_{1}}+\partial
_{t_{2}t_{0}})+...)\nonumber\\
(h_{0}+\varepsilon h_{1}+\varepsilon^{2}h_{2}+...)+\varepsilon(\partial
_{t_{0}}+\varepsilon\partial_{t_{1}}+\varepsilon^{2}\partial_{t_{2}
}+...)\nonumber\\
(h_{0}+\varepsilon h_{1}+\varepsilon^{2}h_{2}+...)+h_{0}+\varepsilon
h_{1}+\varepsilon^{2}h_{2}+...=0,\nonumber\\
\Downarrow\nonumber\\
\partial_{t_{0}t_{0}}h_{0}+h_{0}+\varepsilon(\partial_{t_{0}t_{0}}h_{1}
+h_{1}+\partial_{t_{0}t_{1}}h_{0}+\partial_{t_{1}t_{0}}h_{0}+\partial_{t_{0}
}h_{0})\nonumber\\
+\varepsilon^{2}(\partial_{t_{0}t_{0}}h_{2}+h_{2}+\partial_{t_{0}t_{1}}
h_{1}+\partial_{t_{1}t_{0}}h_{1}+\partial_{t_{0}t_{2}}h_{0}+\partial
_{t_{1}t_{1}}h_{0}\nonumber\\
+\partial_{t_{2}t_{0}}h_{0}+\partial_{t_{1}}h_{0}+\partial_{t_{0}}
h_{1})+...=0,
\end{gather*}
which gives us the following perturbation hierarchy to second order in
$\varepsilon$
\begin{align*}
\partial_{t_{0}t_{0}}h_{0}+h_{0}  &  =0,\nonumber\\
& \nonumber\\
\partial_{t_{0}t_{0}}h_{1}+h_{1}  &  =-\partial_{t_{0}t_{1}}h_{0}
-\partial_{t_{1}t_{0}}h_{0}-\partial_{t_{0}}h_{0},\nonumber\\
& \nonumber\\
\partial_{t_{0}t_{0}}h_{2}+h_{2}  &  =-\partial_{t_{0}t_{1}}h_{1}
-\partial_{t_{1}t_{0}}h_{1}-\partial_{t_{0}t_{2}}h_{0}\nonumber\\
&  -\partial_{t_{1}t_{1}}h_{0}-\partial_{t_{2}t_{0}}h_{0}-\partial_{t_{1}
}h_{0}-\partial_{t_{0}}h_{1}.
\end{align*}
We observe, in passing, that the perturbation hierarchy has the special form we
have seen several times before. Here the common differential operator is
$L=\partial_{t_{0}t_{0}}+1$.

At this point a remark is in order. It is fair to say that there is not a full
agreement among the practitioners of the method of multiple scales about how
to perform these calculations. The question really hinges on whether to take
the multiple variable function $h(t_{0},t_{1},..)$ seriously or not. If you
do, you will be lead to a certain way of doing these calculation. This is the
point of view used in most textbooks on this subject. We will not follow this
path here. We will not take $h$ seriously as a multiple variable function and
never forget that what we actually want is not $h$, but rather $y$, which is
defined in terms of $h$ through equation (\ref{Eq4.9.12}). This point of view
will lead us to do multiple scale calculations in a different way from what
you see in most textbooks. This way is very efficient and will make it
possible to go to order $\varepsilon^{2}$ and beyond without being overwhelmed
by the amount of algebra that needs to be done.

  What I mean when I say that we
will not take $h$ seriously as a multiple variable function will become clear
as we proceed. One immediate consequence of this choise is already evident from the way I write the
perturbation hierarchy. Observe that I keep
\begin{equation*}
\partial_{t_{i}t_{j}}h_{k}\;\;\text{and \ }\partial_{t_{j}t_{i}}h_{k}\;\;,
\end{equation*}
as separate terms, I don't use the equality of cross derivatives to simplify
my expressions. This is the first rule we must follow when we do multiple
scale calculations in the way I am teaching you in these lecture notes. If we
took $h$ seriously as a multiple variable function we would put cross
derivatives equal. The second rule we must follow is to disregard the
initial values for the time being. We will fit the initial values at the very
end of our calculations rather than do it at each order in $\varepsilon$ like
we just did in section \ref{NonlinearBoundaryValueProblem} and \ref{WeaklyDampedLinearOscillator}.

Let us now proceed to solve the equations in the perturbation hierarchy. At
order $\varepsilon^{0}$ we have the equation
\begin{equation}
\partial_{t_{0}t_{0}}h_{0}+h_{0}=0.\label{Eq4.16.12}
\end{equation}
When we are applying multiple scales to ordinary differential equations we
always use the general solution to the order $\varepsilon^{0}$ equation. \ For
partial differential equations this will not be so, as we will see later. The
general solution to (\ref{Eq4.16.12}) is evidently
\begin{equation*}
h_{0}(t_{0},t_{1},..)=A_{0}(t_{1},t_{2},..)e^{it_{0}}+(\ast).
\end{equation*}
Observe that the equation only determines how $h_{0}$ depends on the fastest
time scale $t_{0}$, the dependence on the other time scales $t_{1},t_{2},..$,
is arbitrary at this point and this is reflected in the fact that the  integration
"constant" $A_{0}$ is actually a function depending on $t_{1},t_{2},..$\;.

We have now solved the order $\varepsilon^{0}$ equation. Inserting the
expression for $h_{0}$ into the order $\varepsilon$ equation, we get after
some simple algebra
\begin{equation}
\partial_{t_{0}t_{0}}h_{1}+h_{1}=-2i(\partial_{t_{1}}A_{0}+\frac{1}{2}
A_{0})e^{it_{0}}+(\ast).\label{Eq4.18.12}
\end{equation}
We now need a particular solution to this equation. Observe that since $A_{0}
$ only depends on the slow time scales $t_{1},t_{2},..$,  equation
(\ref{Eq4.18.12}) is in fact a harmonic oscillator driven on ressonance. It is
simple to verify that it has a particular solution of the form
\begin{equation}
h_{1}(t_{0},t_{1},..)=-t_{0}(\partial_{t_{1}}A_{0}+\frac{1}{2}A_{0})e^{it_{0}}.\label{Eq4.19.12}
\end{equation}
But this term is growing and will lead to breakdown of ordering for the
perturbation expansion (\ref{Eq4.11.12}) when $t_{0}\sim\varepsilon^{-1}$. This
breakdown was exactly what we tried to avoid using the multiple scales approach!

But everything is not lost, we now have freedom to remove the growing term by
postulating that
\begin{equation*}
\partial_{t_{1}}A_{0}=-\frac{1}{2}A_{0}.
\end{equation*}
With this choise, the order $\varepsilon$ equation simplifies into
\begin{equation*}
\partial_{t_{0}t_{0}}h_{1}+h_{1}=0.
\end{equation*}
Terms in equations leading to linear growth like in (\ref{Eq4.19.12}), are
traditionally called \textit{secular terms}. The name are derived from the
Latin word soeculum that means century and is used here because this kind of
non-uniformity was first observed on century time scales in planetary orbit calculations.

At this point we introduce the third rule for doing multiple scale
calculations in the particular way that I advocate in these lecture notes. The
rule is to disregard the general solution of the homogeneous equation for all
equations in the perturbation hierarchy except the first. \ We therefore
choose $h_{1}=0$ and proceed to the order $\varepsilon^{2}$ equation using
this choice. The equation for $h_{2}$ then simplifies into
\begin{equation*}
\partial_{t_{0}t_{0}}h_{2}+h_{2}=-2i(\partial_{t_{2}}A_{0}-\frac{i}{2}
\partial_{t_{1}t_{1}}A_{0}-\frac{i}{2}\partial_{t_{1}}A_{0})e^{it_{0}}
+(\ast).
\end{equation*}
We have a new secular term and in order to remove it we must postulate that
\begin{equation*}
\partial_{t_{2}}A_{0}=\frac{i}{2}\partial_{t_{1}t_{1}}A_{0}+\frac{i}
{2}\partial_{t_{1}}A_{0}.
\end{equation*}
Using this choice, our order $\varepsilon^{2}$ equation simplifies into
\begin{equation*}
\partial_{t_{0}t_{0}}h_{2}+h_{2}=0.
\end{equation*}
For this equation we use, according to the rules of the game, the special
solution $h_{2}=0$.

What we have found so far is then
\begin{equation}
h(t_{0},t_{1},t_{2},..)=A_{0}(t_{1},t_{2},..)e^{it_{0}}+(\ast)+O(\varepsilon
^{3}),\label{Eq4.25.12}
\end{equation}
where
\begin{align}
\partial_{t_{1}}A_{0}  &  =-\frac{1}{2}A_{0},\label{Eq4.26.12}\\
\partial_{t_{2}}A_{0}  &  =\frac{i}{2}\partial_{t_{1}t_{1}}A_{0}+\frac{i}
{2}\partial_{t_{1}}A_{0}.\label{Eq4.27.12}
\end{align}
At this point you might ask if we actually have done something useful. Instead
of one ODE we have ended up with two coupled partial differential equations,
and clearly, if we want to go to higher order we will get even more partial
differential equations.

Observe that if we use (\ref{Eq4.26.12}) we can simplify equation (\ref{Eq4.27.12})
by removing the derivatives on the right hand side. Doing this we get the
system
\begin{align}
\partial_{t_{1}}A_{0}  &  =-\frac{1}{2}A_{0},\label{Eq4.28.12}\\
\partial_{t_{2}}A_{0}  &  =-\frac{i}{8}A_{0}.\label{Eq4.29.12}
\end{align}
The first thing that should come to mind when we see a system like
(\ref{Eq4.28.12}) and (\ref{Eq4.29.12}), is that the count is wrong. There is one
unknown function, $A_{0}$, and two equations. The system is
\textit{overdetermined} and will get more so, if we extend our calculations to
higher order in $\varepsilon$. Under normal circumstances, overdetermined
systems of equations have no solutions, which for our setting means that under
normal circumstances the function $h(t_{0},t_{1},t_{2},..)$ does not exist!
This is what I meant when I said that we will not take the functions $h$
seriously as a multiple variable function. For systems of first order partial
differential equations like (\ref{Eq4.28.12}), (\ref{Eq4.29.12}) there is a simple
test we can use to decide if a solution actually does  exist. This is the cross
derivative test you know from elementary calculus. Taking $\partial_{t_{2}}$
of equation (\ref{Eq4.28.12}) and $\partial_{t_{1}}$ of equation (\ref{Eq4.29.12})
we get
\begin{align*}
\partial_{t_{2}t_{1}}A_{0}  &  =\partial_{t_{2}}\partial_{t_{1}}A_{0}
=-\frac{1}{2}\partial_{t_{2}}A_{0}=\frac{i}{16}A_{0},\nonumber\\
\partial_{t_{1}t_{2}}A_{0}  &  =\partial_{t_{1}}\partial_{t_{2}}A_{0}
=-\frac{i}{8}\partial_{t_{1}}A_{0}=\frac{i}{16}A_{0}.
\end{align*}
According to the cross derivative test the overdetermined system
(\ref{Eq4.28.12}), (\ref{Eq4.29.12}) is solvable. Thus in this case the function $h$
exists, at least as a two variable function. To make sure that it exists as a
function of three variables we must derive and solve the perturbation
hierarchy to order $\varepsilon^{3}$, and then perform the cross derivative
test. For the current example we will never get into trouble, the many
variable function $h$ will exist as a function of however many variables we
want. But I want you to reflect on how special this must be. We will at order
$\varepsilon^{n}$ have a system of $n$ partial differential equations for only
one unknown function ! In general we will not be so lucky as in the current example,and
the function $h(t_0,t_1,...)$ will not exist. This fact   is the
reason why we can not take $h$ seriously as a many variable function.

So, should we be disturbed by the nonexistence of the solution to
the perturbation hierarchy in the general
case? Actually no, and the reason is that we do not care about $h(t_{0}
,t_{1},..)$. What we care about is $y(t)$. \ 

Inspired by this let us define an \textit{amplitude}, $A(t)$, by
\begin{equation}
A(t)=A_{0}(t_{1},t_{2},..)|_{t_{j}=\varepsilon^{j}t}.\label{Eq4.31.12}
\end{equation}
Using this and equations (\ref{Eq4.9.12}) and  (\ref{Eq4.25.12}), our perturbation
expansion for $y(t)$ is
\begin{equation*}
y(t)=A(t)e^{it}+(\ast)+O(\varepsilon^{3}).
\end{equation*}
For the amplitude $A(t)$ we have, using equations (\ref{Eq4.10.12}),(\ref{Eq4.28.12}
),(\ref{Eq4.29.12}) and (\ref{Eq4.31.12})
\begin{gather*}
\frac{dA}{dt}(t)=\{(\partial_{t_{0}}+\varepsilon\partial_{t_{1}}
+\varepsilon^{2}\partial_{t_{2}}+...)A_{0}(t_{1},t_{2},...)\}|_{t_{j}
=\varepsilon^{j}t},\nonumber\\
\Downarrow\nonumber\\
\frac{dA}{dt}(t)=\{-\varepsilon\frac{1}{2}A_{0}(t_{1},t_{2},...)-\varepsilon
^{2}\frac{i}{8}A_{0}(t_{1},t_{2},...)\}|_{t_{j}=\varepsilon^{j}t},\nonumber\\
\Downarrow\nonumber\\
\frac{dA}{dt}=-\varepsilon\frac{1}{2}A-\varepsilon^{2}\frac{i}{8}
A.
\end{gather*}

\noindent This equation is our first example of an \textit{amplitude equation}. The
amplitude equation determines, through equation (\ref{Eq4.31.12}), the perturbation
expansion for our solution to the original equation (\ref{Eq4.1.12}). The
amplitude equation is of course easy to solve and we get
\begin{equation*}
y(t)=Ce^{-\frac{1}{2}\varepsilon t}e^{it}e^{-\frac{i}{8}\varepsilon^{2}
t}+(\ast)+O(\varepsilon^{3}).
\end{equation*}
The constant $C$ can be fitted to the initial conditions. What we get is equal
to the exact solution up to second order in $\varepsilon$, as we see by
comparing with (\ref{Eq4.5.12}).

Let us next apply the multiple scale method to some weakly nonlinear ordinary
differential equations. For these cases no exact solution is known, so the
multiple scale method will actually be useful!

\subsubsection{A cubic oscillator}

Consider the initial value problem
\begin{align}
\frac{d^{2}y}{dt^{2}}+y  &  =\varepsilon y^{3},\nonumber\\
y(0)  &  =1,\nonumber\\
\frac{dy}{dt}(0)  &  =0.\label{Eq4.35.12}
\end{align}
If we try do solve this problem using a regular perturbation expansion, we will
get secular terms that will lead to breakdown of ordering on a time scale
$t\sim\varepsilon^{-1}$. Let us therefore apply the multiple scale approach.
We introduce a function $h$ through
\begin{equation*}
y(t)=h(t_{0},t_{1},t_{2},...)|_{t_{j}=\varepsilon^{j}t},
\end{equation*}
and expansions
\begin{align*}
\frac{d}{dt}  &  =\partial_{t_{0}}+\varepsilon\partial_{t_{1}}+\varepsilon
^{2}\partial_{t_{2}}+...\;\;,\nonumber\\
h  &  =h_{0}+\varepsilon h_{1}+\varepsilon^{2}h_{2}+...\;\;.
\end{align*}
Inserting these expansions into (\ref{Eq4.35.12}), we get
\begin{gather*}
(\partial_{t_{0}}+\varepsilon\partial_{t_{1}}+\varepsilon^{2}\partial_{t_{2}
}+...)(\partial_{t_{0}}+\varepsilon\partial_{t_{1}}+\varepsilon^{2}
\partial_{t_{2}}+...)\nonumber\\
(h_{0}+\varepsilon h_{1}+\varepsilon^{2}h_{2}+...)+h_{0}+\varepsilon
h_{1}+\varepsilon^{2}h_{2}+...\nonumber\\
=\varepsilon(h_{0}+\varepsilon h_{1}+\varepsilon^{2}h_{2}+...)^{3},\nonumber\\
\Downarrow\nonumber\\
(\partial_{t_{0}t_{0}}+\varepsilon(\partial_{t_{0}t_{1}}+\partial_{t_{1}t_{0}
})+\varepsilon^{2}(\partial_{t_{0}t_{2}}+\partial_{t_{1}t_{1}}+\partial
_{t_{2}t_{0}})+...)\nonumber\\
(h_{0}+\varepsilon h_{1}+\varepsilon^{2}h_{2}+...)+h_{0}+\varepsilon
h_{1}+\varepsilon^{2}h_{2}+...\nonumber\\
=\varepsilon h_{0}^{3}+3\varepsilon^{2}h_{0}^{2}h_{1}+...\;\;,\nonumber\\
\Downarrow\nonumber\\
\partial_{t_{0}t_{0}}h_{0}+h_{0}+\varepsilon(\partial_{t_{0}t_{0}}h_{1}
+h_{1}+\partial_{t_{0}t_{1}}h_{0}+\partial_{t_{1}t_{0}}h_{0})\nonumber\\
+\varepsilon^{2}(\partial_{t_{0}t_{0}}h_{2}+h_{2}+\partial_{t_{0}t_{1}}
h_{1}+\partial_{t_{1}t_{0}}h_{1}+\partial_{t_{0}t_{2}}h_{0}+\partial
_{t_{1}t_{1}}h_{0}\nonumber\\
+\partial_{t_{2}t_{0}}h_{0})+...=\varepsilon h_{0}^{3}+3\varepsilon^{2}
h_{0}^{2}h_{1}+...\;\;,
\end{gather*}
which gives us the following perturbation hierarchy to second order in
$\varepsilon$
\begin{align*}
\partial_{t_{0}t_{0}}h_{0}+h_{0}  &  =0,\nonumber\\
& \nonumber\\
\partial_{t_{0}t_{0}}h_{1}+h_{1}  &  =h_{0}^{3}-\partial_{t_{0}t_{1}}
h_{0}-\partial_{t_{1}t_{0}}h_{0},\nonumber\\
& \nonumber\\
\partial_{t_{0}t_{0}}h_{2}+h_{2}  &  =3h_{0}^{2}h_{1}-\partial_{t_{0}t_{1}
}h_{1}-\partial_{t_{1}t_{0}}h_{1}-\partial_{t_{0}t_{2}}h_{0}\nonumber\\
&  -\partial_{t_{1}t_{1}}h_{0}-\partial_{t_{2}t_{0}}h_{0}.
\end{align*}
The general solution to the first equation in the perturbation hierarchy is
\begin{equation*}
h_{0}=A_{0}(t_{1},t_{2},...)e^{it_{0}}+(\ast).
\end{equation*}
Inserting this into the right hand side of the second equation in the
hierarchy and expanding, we get
\begin{equation*}
\partial_{t_{0}t_{0}}h_{1}+h_{1}=(3|A_{0}|^{2}A_{0}-2i\partial_{t_{1}}
A_{0})e^{it}+A_{0}^{3}e^{3it}+(\ast).
\end{equation*}
In order to remove secular terms we must postulate that
\begin{equation}
\partial_{t_{1}}A_{0}=-\frac{3i}{2}|A_{0}|^{2}A_{0}.\label{Eq4.43.12}
\end{equation}
This choice simplify the equation for $h_{1}$ into
\begin{equation*}
\partial_{t_{0}t_{0}}h_{1}+h_{1}=A_{0}^{3}e^{3it_{0}}+(\ast).
\end{equation*}
According to the rules of the game we now need a particular solution to this
equation. It is easy to verify that
\begin{equation*}
h_{1}=-\frac{1}{8}A_{0}^{3}e^{3it_{0}}+(\ast),
\end{equation*}
is such a particular solution.

We now insert $h_{0}$ and $h_{1}$ into the right hand side of the third
equation in the perturbation hierarchy and find
\begin{equation*}
\partial_{t_{0}t_{0}}h_{2}+h_{2}=(-\frac{3}{8}|A_{0}|^{4}A_{0}-2i\partial
_{t_{2}}A_{0}-\partial_{t_{1}t_{1}}A_{0})e^{it_{0}}+(\ast)+NST,
\end{equation*}
where $NST$ is an acronym for "nonsecular terms". Since we are not here
planning to go beyond second order in $\varepsilon$, we will at this order
only need the secular terms and group the rest into $NST$. In order to remove
the secular terms we must postulate that
\begin{equation}
\partial_{t_{2}}A_{0}=\frac{3i}{16}|A_{0}|^{4}A_{0}+\frac{i}{2}\partial
_{t_{1}t_{1}}A_{0}.\label{Eq4.47.12}
\end{equation}
As before, it make sense to simplify (\ref{Eq4.47.12}) using equation
(\ref{Eq4.43.12}). This leads to the following overdetermined system of equations
for $A_{0}$
\begin{align}
\partial_{t_{1}}A_{0}  &  =-\frac{3i}{2}|A_{0}|^{2}A_{0},\nonumber\\
\partial_{t_{2}}A_{0}  &  =-\frac{15i}{16}|A_{0}|^{4}A_{0}\label{Eq4.48.12}
\end{align}
Let us check solvability of this system using the cross derivative test
\begin{align*}
\partial_{t_{2}t_{1}}A_{0}  &  =-\frac{3i}{2}\partial_{t_{2}}(A_{0}^{2}
A_{0}^{\ast})\\
&  =-\frac{3i}{2}\left(  2A_{0}\partial_{t_{2}}A_{0}A_{0}^{\ast}+A_{0}
^{2}\partial_{t_{2}}A_{0}^{\ast}\right) \\
&  =-\frac{3i}{2}\left(  2A_{0}\left(  -\frac{15i}{16}|A_{0}|^{4}A_{0}\right)
A_{0}^{\ast}+A_{0}^{2}\left(  \frac{15i}{16}|A_{0}|^{4}A_{0}^{\ast}\right)
\right) \\
&  =-\frac{45}{32}|A_{0}|^{6}A_{0}.
\end{align*}

\begin{align*}
\partial_{t_{1}t_{2}}A_{0}  &  =-\frac{15i}{16}\partial_{t_{1}}\left(
A_{0}^{3}A_{0}^{\ast2}\right) \\
&  =-\frac{15i}{16}\left(  3A_{0}^{2}\partial_{t_{1}}A_{0}A_{0}^{\ast2}
+2A_{0}^{3}A_{0}^{\ast}\partial_{t_{1}}A_{0}^{\ast}\right) \\
&  =-\frac{15i}{16}\left(  3A_{0}^{2}\left(  -\frac{3i}{2}|A_{0}|^{2}
A_{0}\right)  A_{0}^{\ast2}+2A_{0}^{3}A_{0}^{\ast}\left(  \frac{3i}{2}
|A_{0}|^{2}A_{0}^{\ast}\right)  \right) \\
&  =-\frac{45}{32}|A_{0}|^{6}A_{0}.
\end{align*}
The system is compatible, and thus the function $h_{0}$ exists as a function of
two variables. Of course, whether or not $h_{0}$ exists is only of academic
interest for us, since our only aim is to find the solution of the original
equation $y(t)$.

  Defining an amplitude, $A(t)$ by
\begin{equation*}
A(t)=A_{0}(t_{1},t_{2},...)|_{t_{j}=\varepsilon^{j}t},
\end{equation*}
we find that the solution is
\begin{equation}
y(t)=A(t)e^{it}-\varepsilon\frac{1}{8}A^{3}e^{3it}+(\ast)+O(\varepsilon
^{2}),\label{Eq4.50.12}
\end{equation}
where the amplitude satisfy the equation
\begin{gather}
\frac{dA}{dt}(t)=\{(\partial_{t_{0}}+\varepsilon\partial_{t_{1}}
+\varepsilon^{2}\partial_{t_{2}}+...)A_{0}(t_{1},t_{2},...)\}|_{t_{j}
=\varepsilon^{j}t},\nonumber\\
\Downarrow\nonumber\\
\frac{dA}{dt}(t)=\{-\varepsilon\frac{3i}{2}|A_{0}|^{2}A_{0}(t_{1}
,t_{2},...)-\varepsilon^{2}\frac{15i}{16}|A_{0}|^{4}A_{0}(t_{1},t_{2}
,...)\}|_{t_{j}=\varepsilon^{j}t},\nonumber\\
\Downarrow\nonumber\\
\frac{dA}{dt}=-\varepsilon\frac{3i}{2}|A|^{2}A-\varepsilon^{2}\frac{15i}
{16}|A|^{4}A.\label{Eq4.51.12}
\end{gather}
Observe that this equation has a unique solution for a given set of initial
conditions regardless of whether the overdetermined system (\ref{Eq4.48.12}) has
a solution or not. Thus doing the cross derivative test was only motivated by
intellectual curiosity, we did not have to do it.

In summary, (\ref{Eq4.50.12}) and (\ref{Eq4.51.12}), determines a perturbation
expansion for $y(t)$ that is uniform for $t\lesssim\varepsilon^{-3}$.

At this point it is reasonable to ask in which sense we have made progress. We
started with one second order nonlinear ODE for a real function $y(t)$ and
have ended up with one first order nonlinear ODE for a complex function
$A(t)$.

  This question actually has two different answers.
The first one is that it is possible to get an analytical solution for
(\ref{Eq4.51.12}) whereas this is not possible for the original equation
(\ref{Eq4.35.12}). This possibility might however easily get lost as we proceed
to higher order in $\varepsilon$, since this will add more terms to the
amplitude equation. But even if we can not solve the amplitude equation
exactly, it is a fact that amplitude equations with the \textit{same
}mathematical structure will arise when we apply the multiple scale method to
many \textit{different} equations. Thus any insight into an amplitude equation
derived by some mathematical analysis has relevance for many different
situations. This is clearly very useful.

There is however a second, more robust, answer to the question of whether we
have made progress or not. From a numerical point of view, there is an
important difference between (\ref{Eq4.35.12}) and (\ref{Eq4.51.12}). If we solve
(\ref{Eq4.35.12}) numerically, the time step is constrained by the oscillation
period of the linearized system
\begin{equation}
\frac{d^{2}y}{dt^{2}}+y=0.\label{Eq4.52.12}
\end{equation}

\noindent which is of order $T\sim1$, whereas if we solve (\ref{Eq4.51.12}), numerically the
timestep is constrained by the period $T\sim\varepsilon^{-1}$. Therefore, if we
want to propagate out to a time $t\sim\varepsilon^{-2}$, we must take on the
order of $\varepsilon^{-2}$ time steps if we use (\ref{Eq4.52.12}) whereas we
only need on the order of $\varepsilon^{-1}$ time steps using (\ref{Eq4.51.12}).
If $\varepsilon$ is very small the difference in the number of time steps can
be highly significant. From this point of view, the multiple scale method is a
\textit{reformulation} that is the key element in a fast \textit{numerical} method
for solving weakly nonlinear ordinary and partial differential equation.

Let us next turn to the problem of fitting the initial conditions. Using
equation (\ref{Eq4.50.12}) and the initial conditions from (\ref{Eq4.35.12}) we get,
truncating at order $\varepsilon$, the following equations
\begin{align*}
A(0)-\varepsilon\frac{1}{8}A^{3}(0)+(\ast)  &  =1,\nonumber\\
iA(0)-\varepsilon(\frac{3i}{2}|A(0)|^{2}A(0)+\frac{3i}{8}A^{3}(0))+(\ast)  &
=0.
\end{align*}
The solution for $\varepsilon=0$ is
\begin{equation*}
A(0)=\frac{1}{2}.
\end{equation*}
For $\varepsilon>0$ we solve the equation by Newton iteration starting with
the solution for $\varepsilon=0$. This will give us the initial condition for
the amplitude equation correct to this order in $\varepsilon$

\begin{figure}[htbp]
\centering
\includegraphics[
height=2.501in,
width=3.9288in
]{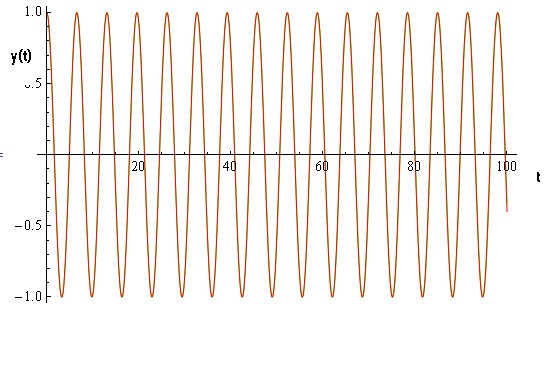}
\caption{Comparing the multiple scale solution, while keeping only the first
term in the amplitude equation(red), to a numerical solution(green) for
$t\lesssim\varepsilon^{-2}$.\label{fig4_4}}
\end{figure}

\begin{figure}[htbp]
\centering
\includegraphics[
height=2.501in,
width=4.0491in
]{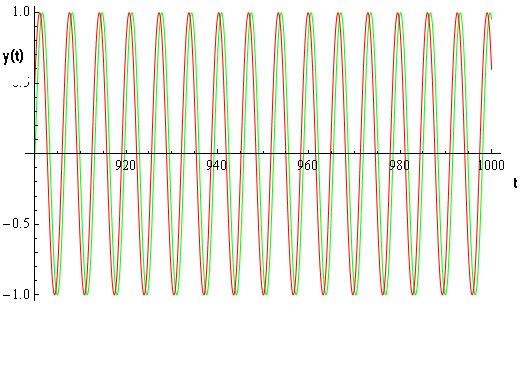}
\caption{Comparing the multiple scale solution, while keeping only the first
term in the amplitude equation(red), to a numerical solution(green) for
$t\lesssim\varepsilon^{-3}$.\label{fig4_5}}
\end{figure}

\begin{figure}[ptb]
\centering
\includegraphics[
height=2.501in,
width=3.6409in
]{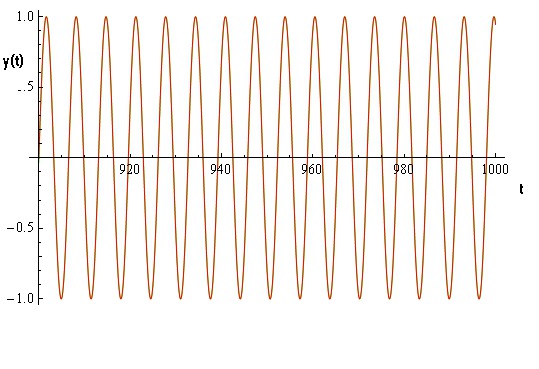}
\caption{Comparing the multiple scale solution, while keeping both terms in the
amplitude equation(red), to a numerical solution(green) for
$t\lesssim\varepsilon^{-3}$.\label{fig4_6}}
\end{figure}

In figure \ref{fig4_4} we compare the multiple scale solution, keeping only
the first term in the amplitude equation, to a high precision numerical
solution for $\varepsilon=0.1$ for
$t\lesssim\varepsilon^{-2}$. We see that the
perturbation solution is very accurate for this range of $t$. In figure
\ref{fig4_5} we do the same comparison as in figure (\ref{fig4_4}) but now for
$t\lesssim\varepsilon^{-3}$. As expected the multiple scale solution and the
numerical solution starts to deviate for this range of $t$. In figure
\ref{fig4_6} we make the same comparison as in figure (\ref{fig4_5}), but now
include both terms in the amplitude equation. We see that high accuracy is
restored for the multiple scale solution for $t\lesssim\varepsilon^{-3}$.

\subsubsection{A second order ODE with a quadratic nonlinearity}\label{QuadraticNonlinearity}

Let us consider the weakly nonlinear equation
\begin{equation}
\frac{d^{2}y}{dt^{2}}+\frac{dy}{dt}+\varepsilon y^{2}=0,\text{ \ \ }
t>0.\label{Eq4.55.12}
\end{equation}
We want to apply the multiple scale method, and introduce a function
$h(t_{0},t_{1},t_{2},..)$ such that
\begin{equation}
y(t)=h(t_{0},t_{1},t_{2},..)|_{t_{j}=\varepsilon^{j}t},\label{Eq4.56.12}
\end{equation}
is a solution to equation (\ref{Eq4.55.12}). As usual we have the formal
expansions
\begin{align}
\frac{d}{dt}  &  =\partial_{t_{0}}+\varepsilon\partial_{t_{1}}+\varepsilon
^{2}\partial_{t_{2}}+...\;\;,\label{Eq4.57.12}\\
h  &  =h_{0}+\varepsilon h_{1}+\varepsilon^{2}h_{2}+...\;\;.\label{Eq4.58.12}
\end{align}
Inserting (\ref{Eq4.56.12}),(\ref{Eq4.57.12}) and (\ref{Eq4.58.12}) into equation
(\ref{Eq4.55.12}) and expanding, we get
\begin{gather*}
(\partial_{t_{0}}+\varepsilon\partial_{t_{1}}+\varepsilon^{2}\partial_{t_{2}
}+...)(\partial_{t_{0}}+\varepsilon\partial_{t_{1}}+\varepsilon^{2}
\partial_{t_{2}}+...)\nonumber\\
(h_{0}+\varepsilon h_{1}+\varepsilon^{2}h_{2}+...)+(\partial_{t_{0}
}+\varepsilon\partial_{t_{1}}+\varepsilon^{2}\partial_{t_{2}}+...)\nonumber\\
(h_{0}+\varepsilon h_{1}+\varepsilon^{2}h_{2}+...)\nonumber\\
=-\varepsilon(h_{0}+\varepsilon h_{1}+\varepsilon^{2}h_{2}+...)^{2},\nonumber\\
\Downarrow\nonumber\\
(\partial_{t_{0}t_{0}}+\varepsilon(\partial_{t_{0}t_{1}}+\partial_{t_{1}t_{0}
})+\varepsilon^{2}(\partial_{t_{0}t_{2}}+\partial_{t_{1}t_{1}}+\partial
_{t_{2}t_{0}})+...)\nonumber\\
(h_{0}+\varepsilon h_{1}+\varepsilon^{2}h_{2}+...)+(\partial_{t_{0}
}+\varepsilon\partial_{t_{1}}+\varepsilon^{2}\partial_{t_{2}}+...)\nonumber\\
(h_{0}+\varepsilon h_{1}+\varepsilon^{2}h_{2}+...)=-\varepsilon h_{0}
^{2}-\varepsilon^{2}2h_{0}h_{1}+...\;\;,\nonumber\\
\Downarrow\nonumber\\
\partial_{t_{0}t_{0}}h_{0}+\partial_{t_{0}}h_{0}+\varepsilon(\partial
_{t_{0}t_{0}}h_{1}+\partial_{t_{0}}h_{1}+\partial_{t_{0}t_{1}}h_{0}
+\partial_{t_{1}t_{0}}h_{0}+\partial_{t_{1}}h_{0})\nonumber\\
+\varepsilon^{2}(\partial_{t_{0}t_{0}}h_{2}+\partial_{t_{0}}h_{2}
+\partial_{t_{0}t_{1}}h_{1}+\partial_{t_{1}t_{0}}h_{1}+\partial_{t_{0}t_{2}
}h_{0}+\partial_{t_{1}t_{1}}h_{0}\nonumber\\
+\partial_{t_{2}t_{0}}h_{0}+\partial_{t_{1}}h_{1}+\partial_{t_{2}}
h_{0})+...=-\varepsilon h_{0}^{2}-\varepsilon^{2}2h_{0}h_{1}+...\;\;,
\end{gather*}
which gives us the perturbation hierarchy
\begin{align*}
\partial_{t_{0}t_{0}}h_{0}+\partial_{t_{0}}h_{0}  &  =0,\nonumber\\
& \nonumber\\
\partial_{t_{0}t_{0}}h_{1}+\partial_{t_{0}}h_{1}  &  =-h_{0}^{2}
-\partial_{t_{0}t_{1}}h_{0}-\partial_{t_{1}t_{0}}h_{0}-\partial_{t_{1}}
h_{0},\nonumber\\
& \nonumber\\
\partial_{t_{0}t_{0}}h_{2}+\partial_{t_{0}}h_{2}  &  =-2h_{0}h_{1}
-\partial_{t_{0}t_{1}}h_{1}-\partial_{t_{1}t_{0}}h_{1}-\partial_{t_{0}t_{2}
}h_{0}\nonumber\\
&  -\partial_{t_{1}t_{1}}h_{0}-\partial_{t_{2}t_{0}}h_{0}-\partial_{t_{1}
}h_{1}-\partial_{t_{2}}h_{0}.
\end{align*}
The general solution to the first equation in the perturbation hierarchy is
\begin{equation}
h_{0}(t_{0},t_{1},t_{2},...)=A_{0}(t_{1},t_{2},..)+B_{0}(t_{1},t_{2}
,...)e^{-t_{0}},\label{Eq4.61.12}
\end{equation}
where $A_{0}$ and $B_{0}$ are real functions of their arguments. Inserting
$h_{0}$ into the second equation in the hierarchy we get
\begin{equation}
\partial_{t_{0}t_{0}}h_{1}+\partial_{t_{0}}h_{1}=-\partial_{t_{1}}A_{0}
-A_{0}^{2}+(\partial_{t_{1}}B_{0}-2A_{0}B_{0})e^{-t_{0}}-B_{0}^{2}e^{-2t_{0}
}.\label{Eq4.62.12}
\end{equation}
In order to remove secular terms we must postulate that
\begin{align}
\partial_{t_{1}}A_{0}  &  =-A_{0}^{2},\nonumber\\
\partial_{t_{1}}B_{0}  &  =2A_{0}B_{0}.\label{Eq4.63.12}
\end{align}
Equation (\ref{Eq4.62.12}) then  simplifies into
\begin{equation*}
\partial_{t_{0}t_{0}}h_{1}+\partial_{t_{0}}h_{1}=-B_{0}^{2}e^{-2t_{0}
},
\end{equation*}
which has a special solution
\begin{equation}
h_{1}(t_{0},t_{1},..)=-\frac{1}{2}B_{0}^{2}e^{-2t_{0}}.\label{Eq4.65.12}
\end{equation}
Inserting (\ref{Eq4.61.12}) and (\ref{Eq4.65.12}) into the third equation in the
perturbation hierarchy, we get
\begin{equation*}
\partial_{t_{0}t_{0}}h_{2}+\partial_{t_{0}}h_{2}=-\partial_{t_{2}}
A_{0}-\partial_{t_{1}t_{1}}A_{0}+(\partial_{t_{2}}B_{0}-\partial_{t_{1}t_{1}
}B_{0})e^{-t_{0}}+NST.
\end{equation*}
In order to remove secular terms we must postulate that%
\begin{align}
\partial_{t_{2}}A_{0}  &  =-\partial_{t_{1}t_{1}}A_{0},\nonumber\\
\partial_{t_{2}}B_{0}  &  =\partial_{t_{1}t_{1}}B_{0}.\label{Eq4.67.12}
\end{align}
We can as usual use (\ref{Eq4.63.12}) to simplify (\ref{Eq4.67.12}). We are thus
lead to the following overdetermined system for $A_{0}$ and $B_{0}$.
\begin{align}
\partial_{t_{1}}A_{0}  &  =-A_{0}^{2},\nonumber\\
\partial_{t_{1}}B_{0}  &  =2A_{0}B_{0},\nonumber\\
\partial_{t_{2}}A_{0}  &  =-2A_{0}^{3},\nonumber\\
\partial_{t_{2}}B_{0}  &  =2A_{0}^{2}B_{0}.\label{Eq4.68.12}
\end{align}
In order to satisfy our academic curiosity, let us do the cross derivative
test for solvability of (\ref{Eq4.68.12}).
\begin{align}
\partial_{t_{1}t_{2}}A_{0}  &  =-2\partial_{t_{1}}A_{0}^{3}=-6A_{0}
^{2}\partial_{t_{1}}A_{0}=6A_{0}^{4},\nonumber\\
\partial_{t_{2}t_{1}}A_{0}  &  =-\partial_{t_{2}}A_{0}^{2}=-2A_{0}
\partial_{t_{2}}A_{0}=4A_{0}^{4},\nonumber\\
& \nonumber\\
\partial_{t_{1}t_{2}}B_{0}  &  =2\partial_{t_{1}}(A_{0}^{2}B_{0}
)=4A_{0}\partial_{t_{1}}A_{0}B_{0}+2A_{0}^{2}\partial_{t_{1}}B_{0}
=0,\nonumber\\
\partial_{t_{2}t_{1}}B_{0}  &  =2\partial_{t_{2}}(A_{0}B_{0})=2\partial
_{t_{2}}A_{0}B_{0}+2A_{0}\partial_{t_{2}}B_{0}=0.\nonumber
\end{align}
We see that the test fails, so the system (\ref{Eq4.68.12}) has no solutions.
However the multiple scale method does \textit{not} fail since we are not
actually interested in the functions $A_{0}$ and $B_{0}$ that defines $h_{0}$,
but is rather interested in the function $y(t)$. Define two amplitudes $A(t)$
and $B(t)$ by%
\begin{align}
A(t)  &  =A_{0}(t_{1},t_{2},...)|_{t_{j}=\varepsilon^{j}t},\nonumber\\
B(t)  &  =B_{0}(t_{1},t_{2},...)|_{t_{j}=\varepsilon^{j}t},\label{Eq4.70.12}
\end{align}
then the solution to (\ref{Eq4.55.12}) is
\begin{equation}
y(t)=A(t)+B(t)e^{-t}-\varepsilon\frac{1}{2}B^{2}(t)e^{-2t}+O(\varepsilon
^{2}),\label{Eq4.71.12}
\end{equation}
where the amplitudes $A(t)$ and $B(t)$ satisfy the equations
\begin{align}
\frac{dA}{dt}(t)  &  =\{(\partial_{t_{0}}+\varepsilon\partial_{t_{1}
}+\varepsilon^{2}\partial_{t_{2}}+...)A_{0}(t_{1},t_{2},...)\}|_{t_{j}
=\varepsilon^{j}t},\nonumber\\
&  \Downarrow\nonumber\\
\frac{dA}{dt}(t)  &  =\{-\varepsilon A^{2}(t_{1},t_{2},...)-2\varepsilon
^{2}A(t_{1},t_{2},...)^{3}\}|_{t_{j}=\varepsilon^{j}t},\nonumber\\
&  \Downarrow\nonumber\\
\frac{dA}{dt}  &  =-\varepsilon A^{2}-2\varepsilon^{2}A^{3}.\label{Eq4.72.12}
\end{align}
and
\begin{align}
\frac{dB}{dt}(t)  &  =\{(\partial_{t_{0}}+\varepsilon\partial_{t_{1}
}+\varepsilon^{2}\partial_{t_{2}}+...)B_{0}(t_{1},t_{2},...)\}|_{t_{j}
=\varepsilon^{j}t},\nonumber\\
&  \Downarrow\nonumber\\
\frac{dB}{dt}(t)  &  =\{(2\varepsilon A_{0}(t_{1},t_{2},...)B_{0}(t_{1}
,t_{2},...)+2\varepsilon^{2}A_{0}^{2}(t_{1},t_{2},...)B_{0}(t_{1}
,t_{2},......)\}|_{t_{j}=\varepsilon^{j}t},\nonumber\\
&  \Downarrow\nonumber\\
\frac{dB}{dt}  &  =2\varepsilon AB+2\varepsilon^{2}A^{2}B.\label{Eq4.73.12}
\end{align}
Given the initial conditions for $A$ and $B$, equations (\ref{Eq4.72.12}) and
(\ref{Eq4.73.12}) clearly has a unique solution, and our multiple scale method
will ensure that the perturbation expansion (\ref{Eq4.71.12}) will stay uniform
for $t\lesssim\varepsilon^{-3}$. As for the previous example, the initial
conditions $A(0)$ and $B(0)$ are calculated from the initial conditions for
(\ref{Eq4.55.12}) by a Newton iteration. Thus we see again that the existence or
not of $h(t_{0},..)$ is irrelevant for constructing a uniform perturbation expansion.

The system (\ref{Eq4.72.12}) and (\ref{Eq4.73.12}) can be solved analytically in
terms of implicit functions. However, as we have discussed before,
\ analytical solvability is nice, but not robust. If we take the expansion to
order $\varepsilon^{3}$, more terms are added to the amplitude equations and
the property of analytic solvability can easily be lost. What \textit{is}
robust is that the presence of $\varepsilon$ in the amplitude equations makes
(\ref{Eq4.72.12}) and (\ref{Eq4.73.12}) together with (\ref{Eq4.71.12}) into a fast
numerical scheme for solving the ordinary differential equation (\ref{Eq4.55.12}%
). This property does \textit{not} go away if we take the perturbation
expansion to higher order in $\varepsilon$.

\subsubsection{Two coupled cubic oscillators}\label{CCO}

So far, we have only been applying the method of multiple scales to scalar
ODEs. This is not a limitation on the method, it may equally well be applied
to systems of ordinary differential equations. The mechanics of the method for
systems of equations is very similar to what we have seen for scalar
equations. The only major difference is how we decide which terms are secular
and must be removed. For systems, this problem is solved by using the Fredholm
Alternative theorem, this is in fact one of the major areas of application for
this theorem in applied mathematics.

Let us consider the following system of two coupled second order ODEs.
\begin{align}
x^{\prime\prime}+2x-y  & =\varepsilon xy^{2},\nonumber\\
y^{\prime\prime}+3y-2x  & =\varepsilon yx^{2},\label{Eq4.74.12}
\end{align}
where $\varepsilon\ll1$. We will solve the system using the method of multiple
scales and introduce therefore two functions $h=h(t_{0},t_{1},...)$ and
$k=k(t_{0},t_{1},...)$ such that
\begin{align}
x(t)  & =h(t_{0},t_{1},...)|_{t_{j}=\varepsilon^{j}t},\nonumber\\
y(t)  & =k(t_{0},t_{1},...)|_{t_{j}=\varepsilon^{j}t},\label{Eq4.75.12}
\end{align}
is a solution to equation (\ref{Eq4.74.12}). As usual we have
\begin{equation}
\frac{d}{dt}=\partial_{t_{0}}+\varepsilon\partial_{t_{1}}+...\;\;,\label{Eq4.76.12}
\end{equation}
and for $h$ and $k$ we introduce the expansions
\begin{align}
h  & =h_{0}+\varepsilon h_{1}+...\,\,,\nonumber\\
k  & =k_{0}+\varepsilon k_{1}+...\;\;.\label{Eq4.77.12}
\end{align}
Inserting (\ref{Eq4.75.12}),(\ref{Eq4.76.12}) and (\ref{Eq4.77.12}) into equation
(\ref{Eq4.74.12}), and expanding everything in sight to first order in
$\varepsilon$ we get, after some tedious algebra, the following perturbation
hierarchy
\begin{align}
\partial_{t_{0}t_{0}}h_{0}+2h_{0}-k_{0}  & =0,\nonumber\\
\partial_{t_{0}t_{0}}k_{0}+3k_{0}-2h_{0}  & =0,\label{Eq4.78.12}\\
& \nonumber\\
\partial_{t_{0}t_{0}}h_{1}+2h_{1}-k_{1}  & =-\partial_{t_{0}t_{1}}
h_{0}-\partial_{t_{1}t_{0}}h_{0}+h_{0}k_{0}^{2},\nonumber\\
\partial_{t_{0}t_{0}}k_{1}+3k_{1}-2h_{1}  & =-\partial_{t_{0}t_{1}}
k_{0}-\partial_{t_{1}t_{0}}k_{0}+k_{0}h_{0}^{2}.\label{Eq4.79.12}
\end{align}
Let us start by finding the general solution to the order $\varepsilon^{0} $\newline
equations (\ref{Eq4.78.12}). They can be written as the following linear system
\begin{equation}
\partial_{t_{0}t_{0}}\left(
\begin{array}
[c]{c}
h_{0}\\
k_{0}
\end{array}
\right)  =\left(
\begin{array}
[c]{cc}
-2 & 1\\
2 & -3
\end{array}
\right)  \left(
\begin{array}
[c]{c}
h_{0}\\
k_{0}
\end{array}
\right),\label{Eq4.80.12}
\end{equation}
Let us look for a solution of the form
\begin{equation}
\left(
\begin{array}
[c]{c}
h_{0}\\
k_{0}
\end{array}
\right)  =\mathbf{\alpha}e^{i\omega t_{0}},\label{Eq4.81.12}
\end{equation}
where $\mathbf{\alpha}$ is a unknown vector and $\omega$ an unknown real
number. Inserting (\ref{Eq4.81.12}) into the system (\ref{Eq4.80.12}) and canceling
a common factor we get the the following linear algebraic equation
\begin{equation}
\left(
\begin{array}
[c]{cc}
-2+\omega^{2} & 1\\
2 & -3+\omega^{2}
\end{array}
\right)  \mathbf{\alpha}=0.\label{Eq4.82.12}
\end{equation}
For there to be a nontrivial solution, the determinant of the matrix has to be
zero. This condition leads to the following polynomial equation for $\omega$
\begin{equation*}
\omega^{4}-5\omega^{2}+4=0,
\end{equation*}
which has four real solutions
\[
\omega_{1}=1,\omega_{2}=-1,\omega_{3}=2,\omega_{4}=-2.
\]
A basis for the solution space of (\ref{Eq4.82.12}) corresponding to
$\omega=\omega_{1},\omega_{2}$ is
\begin{equation*}
\mathbf{\alpha}=\left(
\begin{array}
[c]{c}
1\\
1
\end{array}
\right), 
\end{equation*}
and a basis corresponding to $\omega=\omega_{3},\omega_{4}$ is
\begin{equation*}
\mathbf{\beta}=\left(
\begin{array}
[c]{c}
1\\
-2
\end{array}
\right). 
\end{equation*}
It is then clear that a basis for the solution space for the linear system
(\ref{Eq4.80.12}) is
\[
\mathbf{\alpha e}^{\pm it_{0}},\mathbf{\beta}e^{\pm2it_{0}}.
\]
Therefore a general complex solution to (\ref{Eq4.80.12}) is
\begin{equation}
\left(
\begin{array}
[c]{c}
h_{0}\\
k_{0}
\end{array}
\right)  =A_{1}\mathbf{\alpha}e^{it_{0}}+A_{2}\mathbf{\alpha}e^{-it_{0}}%
+B_{1}\mathbf{\beta}e^{2it_{0}}+B_{2}\mathbf{\beta}e^{-2it_{0}}.\label{Eq4.86.12}
\end{equation}
However, we are looking for real solutions to the original system
(\ref{Eq4.74.12}), and in order to ensure reality for (\ref{Eq4.86.12}) we must
choose%
\begin{align*}
A_{1}  & =A_{0}^{\ast},\text{ \ \ }A_{2}=A_{0},\text{ \ \ \ }A_{0}=A_{0}
(t_{1},t_{2},...),\\
B_{1}  & =B_{0}^{\ast},\text{ \ \ }B_{2}=B_{0},\text{ \ \ \ }B_{0}=B_{0}
(t_{1},t_{2},...).
\end{align*}
Thus, a general real solution to (\ref{Eq4.80.12}) is
\begin{equation*}
\left(
\begin{array}
[c]{c}
h_{0}\\
k_{0}
\end{array}
\right)  =A_{0}\mathbf{\alpha}e^{-it_{0}}+B_{0}\mathbf{\beta}e^{-2it_{0}
}+(\ast).
\end{equation*}
In component form, the general real solution is
\begin{align}
h_{0}  & =A_{0}e^{-it_{0}}+B_{0}e^{-2it_{0}}+(\ast),\nonumber\\
k_{0}  & =A_{0}e^{it_{0}}-2B_{0}e^{-2it_{0}}+(\ast).\label{Eq4.88.12}
\end{align}
We now insert the expressions (\ref{Eq4.88.12}) into the order $\varepsilon$
equations (\ref{Eq4.79.12}). After a large amount of tedious algebra, Mathematica
can be useful here, we find that the order $\varepsilon$ equations can be
written in the form
\begin{gather}
\partial_{t_{0}t_{0}}\left(
\begin{array}
[c]{c}
h_{1}\\
k_{1}
\end{array}
\right)  +\left(
\begin{array}
[c]{cc}
2 & -1\\
-2 & 3
\end{array}
\right)  \left(
\begin{array}
[c]{c}
h_{1}\\
k_{1}
\end{array}
\right)  =\left(
\begin{array}
[c]{c}
4B_{0}^{3}\\
-2B_{0}^{3}
\end{array}
\right)  e^{-6it_{0}}+\left(
\begin{array}
[c]{c}
0\\
-3A_{0}B_{0}^{2}
\end{array}
\right)  e^{-5it_{0}}\nonumber\\
-\left(
\begin{array}
[c]{c}
3A_{0}^{\ast}A_{0}^{\ast}B_{0}\\
0
\end{array}
\right)  e^{-4it_{0}}+\left(
\begin{array}
[c]{c}
A_{0}^{2}\\
A_{0}^{3}-3A_{0}^{\ast}B_{0}^{2}
\end{array}
\right)  e^{-3it_{0}}\nonumber\\
+\left(
\begin{array}
[c]{c}
4i\partial_{t_{1}}B_{0}+12|B_{0}|^{2}B_{0}-6|A_{0}|^{2}B_{0}\\
-8i\partial_{t_{1}}B_{0}-6|B_{0}|^{2}B_{0}
\end{array}
\right)  e^{-2it_{0}}\nonumber\\
+\left(
\begin{array}
[c]{c}
2i\partial_{t_{1}}A_{0}+3|A_{0}|^{2}A_{0}\\
2i\partial_{t_{1}}A_{0}+3|A_{0}|^{2}A_{0}-6|B_{0}|^{2}A_{0}
\end{array}
\right)  e^{-it_{0}}+\left(
\begin{array}
[c]{c}
-\frac{3}{2}A_{0}^{\ast}A_{0}^{\ast}B_{0}\\
0
\end{array}
\right)  +(\ast).\label{Eq4.89.12}
\end{gather}
We are not going to go beyond order $\varepsilon$ so we don't actually need to
solve this equation. What we need to do, however, is to ensure that the
solution is bounded in $t_{0}$. We only need a special solution to
(\ref{Eq4.89.12}), and because it is a linear equation, such a special solution
can be constructed as a sum of solutions where each solution in the sum
corresponds to a single term from the right hand side of (\ref{Eq4.89.12}). What we
mean by this is that
\begin{equation}
\left(
\begin{array}
[c]{c}
h_{1}\\
k_{1}
\end{array}
\right)  =\sum_{n=1}^{7}\left(
\begin{array}
[c]{c}
u_{n}\\
v_{n}
\end{array}
\right)  +(\ast),\label{Eq4.90.12}
\end{equation}
where for example
\begin{align}
\partial_{t_{0}t_{0}}\left(
\begin{array}
[c]{c}
u_{1}\\
v_{1}
\end{array}
\right)  +\left(
\begin{array}
[c]{cc}
2 & -1\\
-2 & 3
\end{array}
\right)  \left(
\begin{array}
[c]{c}
u_{1}\\
v_{1}
\end{array}
\right)   & =\left(
\begin{array}
[c]{c}
4B_{0}^{3}\\
-2B_{0}^{3}
\end{array}
\right)  e^{-6it_{0}},\nonumber\\
\partial_{t_{0}t_{0}}\left(
\begin{array}
[c]{c}
u_{2}\\
v_{2}
\end{array}
\right)  +\left(
\begin{array}
[c]{cc}
2 & -1\\
-2 & 3
\end{array}
\right)  \left(
\begin{array}
[c]{c}
u_{2}\\
v_{2}
\end{array}
\right)   & =\left(
\begin{array}
[c]{c}
0\\
-3A_{0}B_{0}^{2}
\end{array}
\right)  e^{-5it_{0}},\label{Eq4.91.12}
\end{align}
and so on. For the first equation we look for a solution of the form
\begin{equation}
\left(
\begin{array}
[c]{c}
u_{1}(t_{0})\\
v_{1}(t_{0})
\end{array}
\right)  =\mathbf{\xi}e^{-6it_{0}},\label{Eq4.92.12}
\end{equation}
where $\mathbf{\xi}$ is a constant vector. Observe that any solution of the
form (\ref{Eq4.92.12}), is bounded in $t_{0}$. If we insert (\ref{Eq4.92.12}) into
the first equation in (\ref{Eq4.91.12}) and cancel the common exponential factor
we find that the unknown vector $\xi$ has to be a solution of the following
linear algebraic system
\[
\left(
\begin{array}
[c]{cc}
-34 & -1\\
-2 & -33
\end{array}
\right)  \mathbf{\xi}=\left(
\begin{array}
[c]{c}
4B_{0}^{3}\\
-2B_{0}^{3}
\end{array}
\right).
\]
The matrix of this system is clearly nonsingular and the solution is
\[
\mathbf{\xi}=\frac{1}{560}\left(
\begin{array}
[c]{c}
-67B_{0}^{3}\\
38B_{0}^{3}
\end{array}
\right),
\]
which gives us the following bounded solution
\[
\left(
\begin{array}
[c]{c}
u_{1}(t_{0})\\
v_{1}(t_{0})
\end{array}
\right)  =\frac{1}{560}\left(
\begin{array}
[c]{c}
-67B_{0}^{3}\\
38B_{0}^{3}
\end{array}
\right)  e^{-6it_{0}}.
\]
A similar approach works for all but the fifth and the sixth term on the
right hand side of equation (\ref{Eq4.89.12}). For these two terms we run into
trouble. For the fifth term we must solve the equation
\begin{eqnarray}
\partial_{t_{0}t_{0}}\left(
\begin{array}
[c]{c}
u_{5}\\
v_{5}
\end{array}
\right)  +\left(
\begin{array}
[c]{cc}
2 & -1\\
-2 & 3
\end{array}
\right)  \left(
\begin{array}
[c]{c}
u_{5}\\
v_{5}
\end{array}
\right)=\nonumber\\
\left(
\begin{array}
[c]{c}
4i\partial_{t_{1}}B_{0}+12|B_{0}|^{2}B_{0}-6|A_{0}|^{2}B_{0}\\
-8i\partial_{t_{1}}B_{0}-6|B_{0}|^{2}B_{0}
\end{array}
\right)  e^{-2it_{0}}.\label{Eq4.93.12}
\end{eqnarray}
A bounded trial solution of the form
\begin{equation*}
\left(
\begin{array}
[c]{c}
u_{5}(t_{0})\\
v_{5}(t_{0})
\end{array}
\right)  =\mathbf{\xi}e^{-2it_{0}},
\end{equation*}
leads to the following algebraic equation for $\mathbf{\xi}$
\begin{equation}
\left(
\begin{array}
[c]{cc}
-2 & -1\\
-2 & -1
\end{array}
\right)  \mathbf{\xi}=\left(
\begin{array}
[c]{c}
4i\partial_{t_{1}}B_{0}+12|B_{0}|^{2}B_{0}-6|A_{0}|^{2}B_{0}\\
-8i\partial_{t_{1}}B_{0}-6|B_{0}|^{2}B_{0}
\end{array}
\right). \label{Eq4.95.12}
\end{equation}
The matrix for this linear system is singular, and the system will in general
have no solution. It will only have a solution, which will lead to a bounded
solution for (\ref{Eq4.93.12}), if the right hand side of (\ref{Eq4.95.12}) satisfy a
certain constraint. This constraint we get from the Fredholm Alternative
Theorem. Recall that this theorem say that a linear system
\[
M\mathbf{x}=b\mathbf{_{0}},
\]
has a solution only if
\[
(\mathbf{f},\mathbf{b_{0}})=0,
\]
for all vectors $\mathbf{f}$ such that
\[
M^{\ast}\mathbf{f}=0,
\]
where $M^{\ast}$ is the adjoint of the matrix $M$. For a real matrix, like the
one we have, $M^{\ast}$ is just the transpose of $M$. For the matrix of the
system (\ref{Eq4.95.12}) we get
\[
\left(
\begin{array}
[c]{cc}
-2 & -2\\
-1 & -1
\end{array}
\right)  \mathbf{f}=0.
\]
A basis for the solution space of this homogeneous system can be taken to be
\[
\mathbf{f}=(1,-1).
\]
Thus in order to ensure solvability of the system (\ref{Eq4.95.12}) we must have
\begin{gather*}
(1,-1)\cdot\left(
\begin{array}
[c]{c}
4i\partial_{t_{1}}B_{0}+12|B_{0}|^{2}B_{0}-6|A_{0}|^{2}B_{0}\\
-8i\partial_{t_{1}}B_{0}-6|B_{0}|^{2}B_{0}
\end{array}
\right)  =0,\nonumber\\
\Updownarrow\nonumber\\
4i\partial_{t_{1}}B_{0}+12|B_{0}|^{2}B_{0}-6|A_{0}|^{2}B_{0}+8i\partial
_{t_{1}}B_{0}+6|B_{0}|^{2}B_{0}=0,\nonumber\\
\Updownarrow\nonumber\\
\partial_{t_{1}}B_{0}=\frac{i}{2}(3|B_{0}|^{2}-|A_{0}|^{2})B_{0}.
\end{gather*}
If this condition on the amplitudes is imposed on the original system
(\ref{Eq4.93.12}), it has a bounded solution. The sixth term in the sum
(\ref{Eq4.90.12}) must be treated in the same way. The equation that we must
solve is
\begin{eqnarray*}
\partial_{t_{0}t_{0}}\left(
\begin{array}
[c]{c}
u_{6}\\
v_{6}
\end{array}
\right)  +\left(
\begin{array}
[c]{cc}
2 & -1\\
-2 & 3
\end{array}
\right)  \left(
\begin{array}
[c]{c}
u_{6}\\
v_{6}
\end{array}
\right)  =\nonumber\\
\left(
\begin{array}
[c]{c}
2i\partial_{t_{1}}A_{0}+3|A_{0}|^{2}A_{0}\\
2i\partial_{t_{1}}A_{0}+3|A_{0}|^{2}A_{0}-6|B_{0}|^{2}A_{0}
\end{array}
\right)  e^{-it_{0}}.
\end{eqnarray*}
Using a bounded trial solution of the form
\begin{equation*}
\left(
\begin{array}
[c]{c}
u_{6}(t_{0})\\
v_{6}(t_{0})
\end{array}
\right)  =\mathbf{\xi}e^{-it_{0}},
\end{equation*}
leads to the following singular linear system
\begin{equation}
\left(
\begin{array}
[c]{cc}
1 & -1\\
-2 & 2
\end{array}
\right)  \mathbf{\xi}=\left(
\begin{array}
[c]{c}
2i\partial_{t_{1}}A_{0}+3|A_{0}|^{2}A_{0}\\
2i\partial_{t_{1}}A_{0}+3|A_{0}|^{2}A_{0}-6|B_{0}|^{2}A_{0}
\end{array}
\right). \label{Eq4.101.12}
\end{equation}
For this case we find that the null space of the transpose of the matrix in
(\ref{Eq4.101.12}) is spanned by the vector
\[
\mathbf{f}=(2,1),
\]
and the Fredholm Alternative gives us the solvability condition
\begin{gather*}
(2.1)\cdot\left(
\begin{array}
[c]{c}
2i\partial_{t_{1}}A_{0}+3|A_{0}|^{2}A_{0}\\
2i\partial_{t_{1}}A_{0}+3|A_{0}|^{2}A_{0}-6|B_{0}|^{2}A_{0}
\end{array}
\right)  =0,\nonumber\\
\Updownarrow\nonumber\\
4i\partial_{t_{1}}A_{0}+6|A_{0}|^{2}A_{0}+2i\partial_{t_{1}}A_{0}
+3|A_{0}|^{2}A_{0}-6|B_{0}|^{2}A_{0}=0,\nonumber\\
\Updownarrow\nonumber\\
\partial_{t_{1}}A_{0}=\frac{i}{2}(3|A_{0}|^{2}-2|B_{0}|^{2})A_{0}.
\end{gather*}
The solutions $h_{1}$ and $k_{1}$ are thus bounded if we impose the following
two conditions on the amplitudes $A_{0}$ and $B_{0}$
\begin{align*}
\partial_{t_{1}}A_{0}  & =\frac{i}{2}(3|A_{0}|^{2}-2|B_{0}|^{2})A_{0},\nonumber\\
\partial_{t_{1}}B_{0}  & =\frac{i}{2}(3|B_{0}|^{2}-|A_{0}|^{2})B_{0}.
\end{align*}
Returning to the original variables $x(t)$ and $y(t)$ in the usual way, we
have thus found the following approximate solution to our system
(\ref{Eq4.74.12})
\begin{align}
x(t)  & =A(t)e^{-it}+B(t)e^{-2it}+\mathcal{O}(\varepsilon),\nonumber\\
y(t)  & =A(t)e^{-it}-2B(t)e^{-2it}+\mathcal{O}(\varepsilon),\label{Eq4.104.12}
\end{align}
where the amplitudes $A(t)$ and $B(t)$ are defined by
\begin{align*}
A(t)  & =A_{0}(t_{1},t_{2},...)|_{t_{j}=\varepsilon^{j}t},\\
B(t)  & =B_{0}(t_{1},t_{2},...)|_{t_{j}=\varepsilon^{j}t},
\end{align*}
and satisfy the following equations
\begin{align}
\partial_{t}A  & =\varepsilon\frac{i}{2}(3|A|^{2}-2|B|^{2})A,\label{Eq4.105.12}\\
\partial_{t}B  & =\varepsilon\frac{i}{2}(3|B|^{2}-|A|^{2})B.\nonumber
\end{align}
The expansions (\ref{Eq4.104.12}) are uniform for $t\lesssim\varepsilon^{-2}$.

The amplitude equations (\ref{Eq4.105.12}) looks complicated, but they are
special in the sense that they can be solved exactly. We have noted before
that the amplitude equations that appears when we use the method of multiple
scale tends to be special. We will see more of this later when we apply the
method to partial differential equations.

Observe that
\begin{align*}
\partial_{t}|A|^{2}  & =\partial_{t}(AA^{\ast})=A^{\ast}\partial
_{t}A+A\partial_{t}A^{\ast}\\
& =A^{\ast}(\frac{i}{2}(3|A|^{2}-2|B|^{2})A)+A(-\frac{i}{2}(3|A|^{2}
-2|B|^{2})A^{\ast})\\
& =\frac{i}{2}(3|A|^{4}-2|B|^{2}|A|^{2}-3|a|^{4}+2|B|^{2}|A|^{2})=0.
\end{align*}
Thus $|A(t)|=|A(0)|$ for all $t$. In a similar way we find that
$|B(t)|=|B(0)|$. Therefore the amplitude equations can be written as
\begin{align*}
\partial_{t}A  & =\frac{i}{2}(3|A(0)|^{2}-2|B(0)|^{2})A,\\
\partial_{t}B  & =\frac{i}{2}(3|B(0)|^{2}-|A(0)|^{2})B.
\end{align*}
and this system is trivial to solve. We find
\begin{align}
A(t)  & =A(0)e^{\frac{i}{2}(3|A(0)|^{2}-2|B(0)|^{2})t},\nonumber\\
B(t)  & =B(0)e^{\frac{i}{2}(3|B(0)|^{2}-|A(0)|^{2})t}.\label{Eq4.106.12}
\end{align}
The formulas (\ref{Eq4.106.12}) together with the expansions (\ref{Eq4.104.12})
gives us an approximate analytic solution to the original system (\ref{Eq4.74.12}). The analytic solution is
\begin{align*}
x(t)  & =A(0)e^{-i\Omega_1 t}+B(0)e^{-i\Omega_2t}+\mathcal{O}(\varepsilon),\nonumber\\
y(t)  & =A(0)e^{-i\Omega_1t}-2B(0)e^{-i\Omega_2t}+\mathcal{O}(\varepsilon),
\end{align*}
where
\begin{align*}
\Omega_1 &=1+|B(0)|^2-\frac{3}{2}|A(0)|^2,\nonumber\\
\Omega_2 &=2+\frac{1}{2}|A(0)|^2-\frac{3}{2}|B(0)|^2,
\end{align*}
and where $A(0)$ and $B(0)$ are two arbitrary complex constants.
Thus, the dynamics of two linear oscillators that are subject to a weak nonlinear coupling  is composed of harmonic motions with respect to two different frequencies $\Omega_{1,2}$ just like for the two linear uncoupled oscillators. What is special about the nonlinearly coupled oscillators is that the frequencies of the two harmonic motions {\it depends on the initial data}. This is not the case for the uncoupled case where the two frequencies are $\Omega_1=1$ and $\Omega_2=2$, no matter what the initial data is. This effect of a weak nonlinear coupling between nonlinear oscillators  is universal.

\subsection{Boundary layer problems for ODEs}

Boundary layer problems first appeared in the theory of fluids. However, boundary layer problems are in no
way limited to fluid theory, but occurs in all areas of science and technology.

In these lecture notes, we will not worry about the physical context for these problems, but
will focus on how to apply the multiple scale method to solve a given problem
of this type. As usual we learn by doing examples.

\subsubsection{A singularly perturbed linear ODE}

Let us consider the following linear boundary value problem
\begin{align}
\varepsilon y^{\prime\prime}(x)+y^{\prime}(x)-y(x)  & =0,\text{ \ \ }
0<x<1,\nonumber\\
y(0)  & =1,\nonumber\\
y(1)  & =0.\label{Eq6.1.12}
\end{align}
We will assume that $\varepsilon\ll1$, and try to solve this problem using a
perturbation methods. The unperturbed problem is clearly
\begin{align*}
y^{\prime}(x)-y(x)  & =0,\text{ \ \ }0<x<1,\nonumber\\
y(0)  & =1,\nonumber\\
y(1)  & =0.\
\end{align*}
The general solution to the differential equation is
\[
y(x)=Ae^{x},
\]
and fitting the boundary condition at $x=0$ we find that
\[
y(x)=e^{x},
\]
but for this solution we have
\[
y(1)=e\neq0,
\]
so the unperturbed problem has no solution. Our perturbation approach fail at
the very first step; there is no unperturbed solution that we can start
calculating corrections to! What is going on?

What is going on is that equation (\ref{Eq6.1.12}) is a singular perturbation
problem. For $\varepsilon\neq0$, we have a second order ODE, whose general
solution has two free constants that can be fitted to the two boundary
conditions, whereas for $\varepsilon=0$ we have a first order ODE whose
general solution has only one free constant. This single constant can in
general not be fitted to two boundary conditions.

We have seen such singular perturbation problems before when we applied
perturbation methods to polynomial equations. For the polynomial case, the
unperturbed problem was of lower algebraic order than the perturbed problem.
Here the unperturbed problem is of lower differential order than the perturbed problem.

For the polynomial case we solved the singular perturbation problem by
transforming it into a regular perturbation problem using a change of
variables. We do the same here.

Let
\begin{equation*}
x=\varepsilon^{p}\xi,\text{ \ \ \ }y(x)=u(\frac{x}{\epsilon^{p}}),
\end{equation*}
then the function $u(\xi)$ satisfy the equation
\begin{equation*}
u^{\prime\prime}(\xi)+\varepsilon^{p-1}u^{\prime}(\xi)-\varepsilon^{2p-1}
u(\xi)=0.
\end{equation*}
This equation constitute a regular perturbation problem if we, for example,
choose $p=1$. We thus have the following regularly perturbed boundary value
problem
\begin{align}
u^{\prime\prime}(\xi)+u^{\prime}(\xi)-\varepsilon u(\xi)  & =0,\text{
\ \ }0<\xi<\frac{1}{\varepsilon},\nonumber\\
u(0)  & =1,\nonumber\\
u(\frac{1}{\varepsilon})  & =0.\label{Eq6.4.12}
\end{align}
Let us try to solve this problem using a perturbation expansion
\begin{equation}
u(\xi)=u_{0}(\xi)+\varepsilon u_{1}(\xi)+...\;\;.\label{Eq6.5.12}
\end{equation}
We will solve the problem by first finding $u_{0}$ and $u_{1}$ and then
fitting the boundary conditions. If we insert the perturbation expansion
(\ref{Eq6.5.12}) into the first equation from (\ref{Eq6.4.12}),  we find the following
perturbation hierarchy to first order in $\varepsilon$
\begin{align}
u_{0}^{\prime\prime}+u_{0}^{\prime}  & =0,\nonumber\\
u_{1}^{\prime\prime}+u_{1}^{\prime}  & =u_{0}.\label{Eq6.52.12}
\end{align}
The general solution to the first equation in the perturbation hierarchy (\ref{Eq6.52.12}), is
clearly
\begin{equation}
u_{0}(\xi)=A_{0}+B_{0}e^{-\xi}.\label{Eq6.6.12}
\end{equation}
If we insert the solution (\ref{Eq6.6.12}) into the second equation in the
perturbation hierarchy (\ref{Eq6.52.12}), we get
\begin{equation}
u_{1}^{\prime\prime}+u_{1}^{\prime}=A_{0}+B_{0}e^{-\xi}.\label{Eq6.7.12}
\end{equation}
Note, that we only need a particular solution to this equation, since the first
term in the perturbation expansion (\ref{Eq6.5.12}) already have two free
constants, and we only need two constants to fit the two boundary data.
Integrating equation (\ref{Eq6.7.12}) once we get
\[
u_{1}^{\prime}+u_{1}=A_{0}\xi-B_{0}e^{-\xi},
\]
and using an integrating factor we get the following particular solution
\[
u_{1}(\xi)=A_{0}(\xi-1)-B_{0}\xi e^{-\xi}.
\]
Thus our perturbation solution to first order in $\varepsilon$ is
\begin{equation*}
u(\xi)=A_{0}+B_{0}e^{-\xi}+\varepsilon\left(  A_{0}(\xi-1)-B_{0}\xi e^{-\xi
}\right)  +...\;\;.
\end{equation*}
The two constants are fitted to the boundary conditions using the following
two equations
\begin{align*}
u(0)  & =1\text{ \ \ }\Longleftrightarrow\text{ \ }A_{0}+B_{0}-\varepsilon
A_{0}=1,\text{\ }\\
u(\frac{1}{\varepsilon})  & =0\text{ \ \ }\Longleftrightarrow\text{ }
A_{0}+B_{0}e^{-\frac{1}{\varepsilon}}+\varepsilon\left(  A_{0}(\frac
{1}{\varepsilon}-1)-B_{0}\frac{1}{\varepsilon}e^{-\frac{1}{\varepsilon}
}\right)  =0.
\end{align*}
However at this point disaster strikes. When we evaluate the solution at the
right boundary $\xi=\frac{1}{\varepsilon}$, using the perturbation expansion,
the ordering of terms is violated. The first and the second term in the
expansion are of the same order. This can not be allowed. Our perturbation
method fails. The reason why the direct perturbation expansion (\ref{Eq6.5.12})
fails is similar to the reason why the direct perturbation expansion failed
for the weakly damped oscillator. In both cases the expansions failed because
they became nonuniform when we evaluated the respective functions at values of
the independent variable that was of order $\varepsilon^{-1}$.

We will resolve the problem with the direct expansion (\ref{Eq6.5.12}) by using
the method of multiple scales to derive a perturbation expansion for the
differential equation from (\ref{Eq6.4.12}), that is uniform for $\xi\lesssim
\varepsilon^{-2}$, and then use this expansion to satisfy the boundary
conditions at $x=0$ and $x=\varepsilon^{-1}$.

\begin{figure}[h]
\centering
\includegraphics[
height=2.2987in,
width=3.4817in
]
{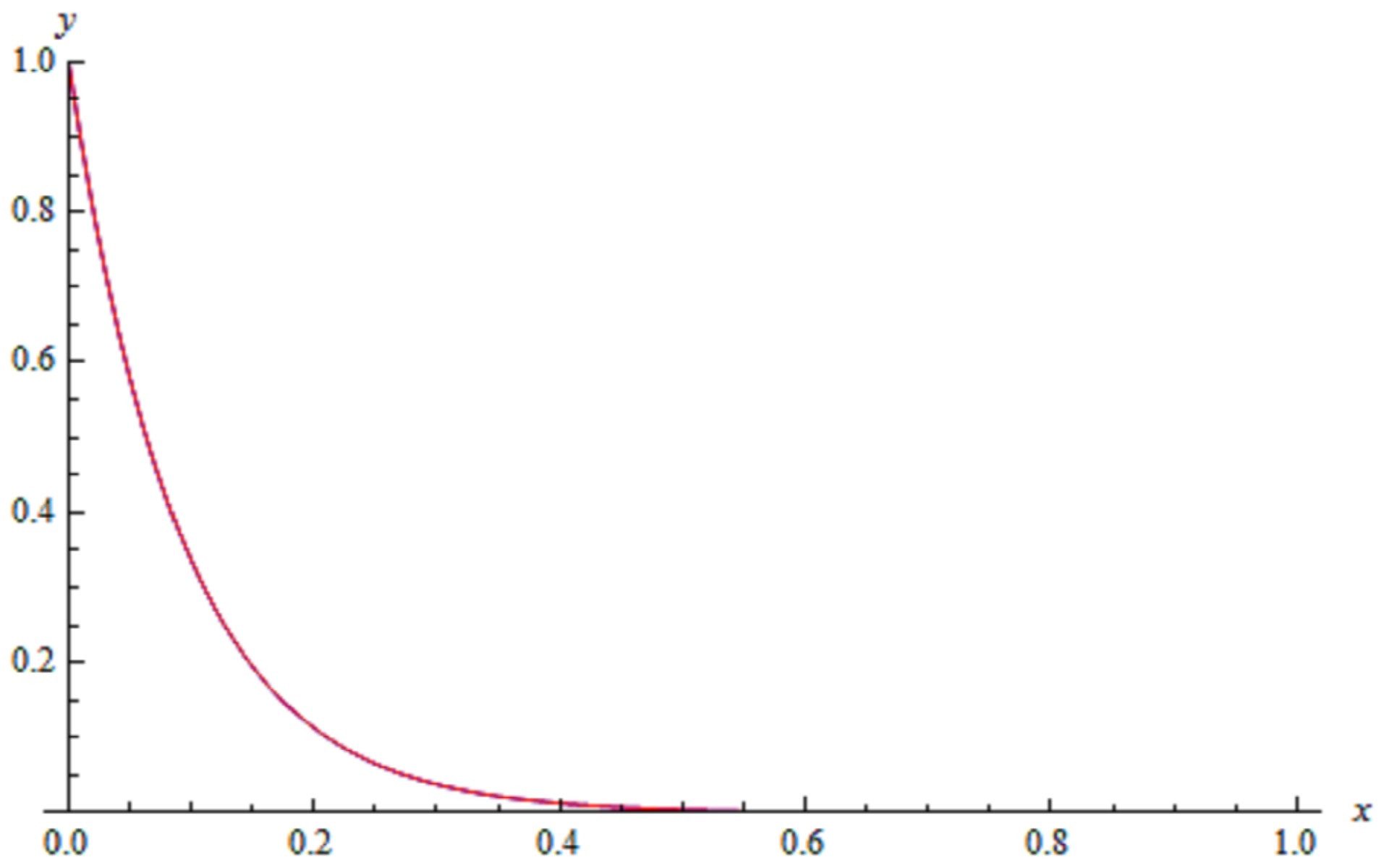}
\caption{Comparing the exact and approximate solution to the singularly
perturbed linear boundary value problem (\ref{Eq6.1.12}). \label{fig4_7}}
\end{figure}

We thus introduce a function $h=h(\xi_{0},\xi_{1},...)$, where $h$ is a function that will be designed 
to ensure that the function $u$,  defined by
\begin{equation}
u(\xi)=h(\xi_{0},\xi_{1},...)|_{\xi_{j}=\varepsilon^{j}\xi},\label{Eq6.8.12}
\end{equation}
is a solution to the equation (\ref{Eq6.4.12}). For the differential operator we
have in the usual way an expansion
\begin{equation}
\frac{d}{d\xi}=\partial_{\xi_{0}}+\varepsilon\partial_{\xi_{1}}+\varepsilon
^{2}\partial_{\xi_{2}}+...\;\;,\label{Eq6.9.12}
\end{equation}
and for the function $h$ we introduce the expansion
\begin{equation}
h=h_{0}+\varepsilon h_{1}+\varepsilon^{2}h_{2}+...\;\;.\label{Eq6.10.12}
\end{equation}
Inserting (\ref{Eq6.8.12}),(\ref{Eq6.9.12}) and (\ref{Eq6.10.12}) into the differential equation from 
(\ref{Eq6.4.12}), and expanding everything in sight to second order in
$\varepsilon$, we get after a small amount of algebra the following
perturbation hierarchy
\begin{align}
\partial_{\xi_{0}\xi_{0}}h_{0}+\partial_{\xi_{0}}h_{0}  & =0,\nonumber\\
& \nonumber\\
\partial_{\xi_{0}\xi_{0}}h_{1}+\partial_{\xi_{0}}h_{1}  & =h_{0}-\partial
_{\xi_{0}\xi_{1}}h_{0}-\partial_{\xi_{1}\xi_{0}}h_{0}-\partial_{\xi_{1}}
h_{0},\nonumber\\
& \nonumber\\
\partial_{\xi_{0}\xi_{0}}h_{2}+\partial_{\xi_{0}}h_{2}  & =h_{1}-\partial
_{\xi_{0}\xi_{1}}h_{1}-\partial_{\xi_{1}\xi_{0}}h_{1}\nonumber\\
& -\partial_{\xi_{0}\xi_{2}}h_{0}-\partial_{\xi_{1}\xi_{1}}h_{0}-\partial
_{\xi_{2}\xi_{0}}h_{0}\nonumber\\
& -\partial_{\xi_{1}}h_{1}-\partial_{\xi_{2}}h_{0}.\label{Eq6.11.12}
\end{align}
The general solution to the first equation in the perturbation hierarchy
(\ref{Eq6.11.12}) is
\begin{equation}
h_{0}(\xi_{0},\xi_{1},\xi_{2},...)=A_{0}(\xi_{1},\xi_{2},...)+B_{0}(\xi
_{1},\xi_{2},...)e^{-\xi_{0}}.\label{Eq6.12.12}
\end{equation}
We now insert this solution into the right hand side of the second equation in
the perturbation hierarchy. Thus the order $\varepsilon$ equation is of the
form
\[
\partial_{\xi_{0}\xi_{0}}h_{1}+\partial_{\xi_{0}}h_{1}=A_{0}-\partial_{\xi
_{1}}A_{0}+(\partial_{\xi_{1}}B_{0}+B_{0})e^{-\xi_{0}}.
\]
Both terms on the right hand side of the equation are secular and in order to
avoid non-uniformity in our expansion we must enforce the conditions
\begin{align}
\partial_{\xi_{1}}A_{0}  & =A_{0},\nonumber\\
\partial_{\xi_{1}}B_{0}  & =-B_{0}.\label{Eq6.13.12}
\end{align}
With these conditions in place, the equation for $h_{1}$ simplify into
\[
\partial_{\xi_{0}\xi_{0}}h_{1}+\partial_{\xi_{0}}h_{1}=0.
\]
and for this equation we choose the special solution
\begin{equation}
h_{1}=0.\label{Eq6.14.12}
\end{equation}
Inserting (\ref{Eq6.12.12}) and (\ref{Eq6.14.12}) into the third equation in the
perturbation hierarchy (\ref{Eq6.11.12}) we get
\[
\partial_{\xi_{0}\xi_{0}}h_{2}+\partial_{\xi_{0}}h_{2}=-\partial_{\xi_{2}
}A_{0}-\partial_{\xi_{1}\xi_{1}}A_{0}+(\partial_{\xi_{2}}B_{0}-\partial
_{\xi_{1}\xi_{1}}B_{0})e^{-\xi_{0}}.
\]
In order to avoid secular terms we enforce the conditions
\begin{align}
\partial_{\xi_{2}}A_{0}  & =-\partial_{\xi_{1}\xi_{1}}A_{0},\nonumber\\
\partial_{\xi_{2}}B_{0}  & =\partial_{\xi_{1}\xi_{1}}B_{0},\label{Eq6.15.12}
\end{align}
and with this choice the equation for $h_{2}$ simplify into
\[
\partial_{\xi_{0}\xi_{0}}h_{2}+\partial_{\xi_{0}}h_{2}=0,
\]
and for this equation we choose the special solution
\[
h_{2}=0.
\]
Using (\ref{Eq6.13.12}), equations (\ref{Eq6.15.12}) can be simplified into
\begin{align*}
\partial_{\xi_{2}}A_{0}  & =-A_{0},\nonumber\\
\partial_{\xi_{2}}B_{0}  & =B_{0}.
\end{align*}
Returning to the original variable $u(\xi)$ in the usual way, we have an
approximate solution to the differential equation from (\ref{Eq6.4.12}) of the form
\begin{equation}
u(\xi)=A(\xi)+B(\xi)e^{-\xi}+\mathcal{O}(\varepsilon^{3}),\label{Eq6.16.12}
\end{equation}
where the amplitudes $A$ and $B$ are defined by
\begin{align*}
A(\xi)  & =A_{0}(\xi_{1},\xi_{2},...)|_{\xi_{j}=\varepsilon^{j}\xi},\\
B(\xi)  & =B_{0}(\xi_{1},\xi_{2},...)|_{\xi_{j}=\varepsilon^{j}\xi},
\end{align*}
and satisfy the equations
\begin{align}
\frac{dA}{d\xi}  & =\varepsilon A-\varepsilon^{2}A,\nonumber\\
\frac{dB}{d\xi}  & =-\varepsilon B+\varepsilon^{2}B.\label{Eq6.17.12}
\end{align}
The amplitude equations (\ref{Eq6.17.12}) are easy to solve, the general solution
is
\begin{align}
A(\xi)  & =Ce^{(\varepsilon-\varepsilon^{2})\xi},\nonumber\\
B(\xi)  & =De^{(-\varepsilon+\varepsilon^{2})\xi},\label{Eq6.18.12}
\end{align}
where $C$ and $D$ are arbitrary real constants. If we insert the solution
(\ref{Eq6.18.12}) into (\ref{Eq6.16.12}) we get
\begin{equation}
u(\xi)=Ce^{(\varepsilon-\varepsilon^{2})\xi}+De^{(-\varepsilon+\varepsilon
^{2}-1)\xi}+\mathcal{O}(\varepsilon^{3}).\label{Eq6.19.12}
\end{equation}
We now determine the constants $C$ and $D$ such that (\ref{Eq6.19.12}) satisfy
the boundary conditions to order $\varepsilon^{2}$.
\begin{align*}
u(0)  & =1\text{ \ \ }\Longleftrightarrow\text{ \ }C+D=1,\\
u(\frac{1}{\varepsilon})  & =0\text{ \ \ }\Longleftrightarrow\text{
\ }Ce^{(1-\varepsilon)}+De^{(-1+\varepsilon-\frac{1}{\varepsilon})}=0.
\end{align*}
The linear system for $C$ and $D$ is easy to solve and we get
\begin{align*}
C  & =(1-e^{2-2\varepsilon+\frac{1}{\varepsilon}})^{-1},\\
D  & =(1-e^{-2+2\varepsilon-\frac{1}{\varepsilon}})^{-1},
\end{align*}
and the approximate solution to the original boundary value problem
(\ref{Eq6.1.12}) is
\begin{equation}
y(x)=(1-e^{2-2\varepsilon+\frac{1}{\varepsilon}})^{-1}e^{(1-\varepsilon
)x}+(1-e^{-2+2\varepsilon-\frac{1}{\varepsilon}})^{-1}e^{(-1+\varepsilon
-\frac{1}{\varepsilon})x}+\mathcal{O}(\varepsilon^{3}).\label{Eq6.20.12}
\end{equation}
In figure \ref{fig4_7} we compare a high precision numerical solution of (\ref{Eq6.1.12}) with the
approximate solution (\ref{Eq6.20.12}) for $\varepsilon=0.1$. The two solutions
are clearly very close over the whole domain.

\begin{figure}[htpb]
\centering
\includegraphics[
height=2.2987in,
width=3.4817in]
{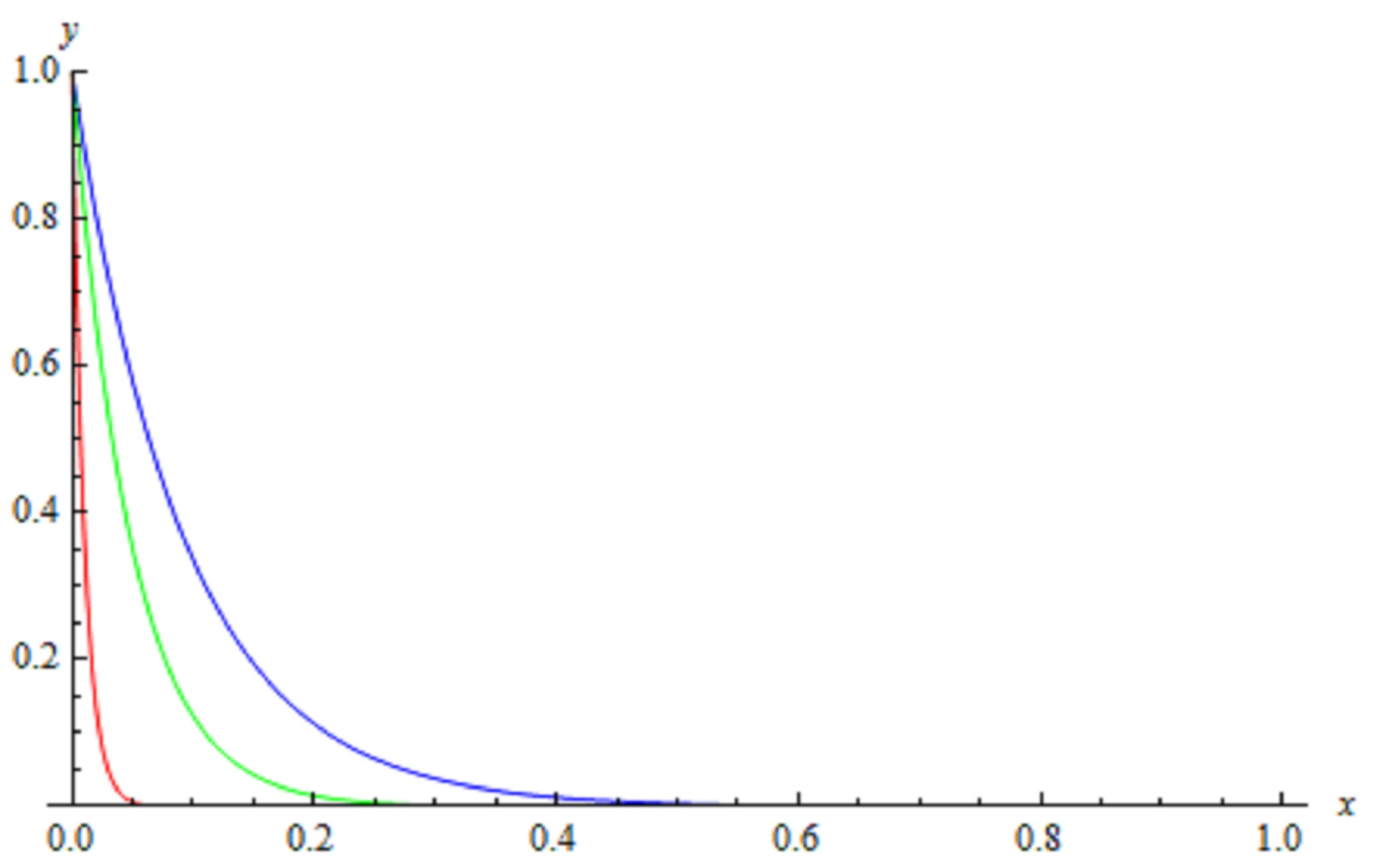}
\caption{A high precision numerical solution of the linear boundary value problem (\ref{Eq6.1.12}), for
$\varepsilon=0.1$ (Blue), $\varepsilon=0.05$ (Green) and $\varepsilon=0.01$
(Red).\label{fig4_8}}
\end{figure}

In figure \ref{fig4_8} we show a high precision numerical solution to the
 boundary value problem (\ref{Eq6.1.12}) for $\varepsilon=0.1$ (Blue),
$\varepsilon=0.05$ (Green) and $\varepsilon=0.01$ (Red).

We observe that the solution is characterized by a very fast variation close
to $x=0$. The domain close to $x=0$, where $y(x)$ experience a fast variation
is called a \textit{boundary layer}. It's extent is of the order of
$\varepsilon$.

In the context of fluids, the boundary layer is the part of the fluid where
the viscosity plays a role. Away from the boundary layer, the dynamics of the
fluid is to a good approximation described by the Euler equation.

\subsubsection{A singularly perturbed nonlinear ODE}

\begin{figure}[h]
\centering
\includegraphics[
height=1.6924in,
width=5.1387in
]
{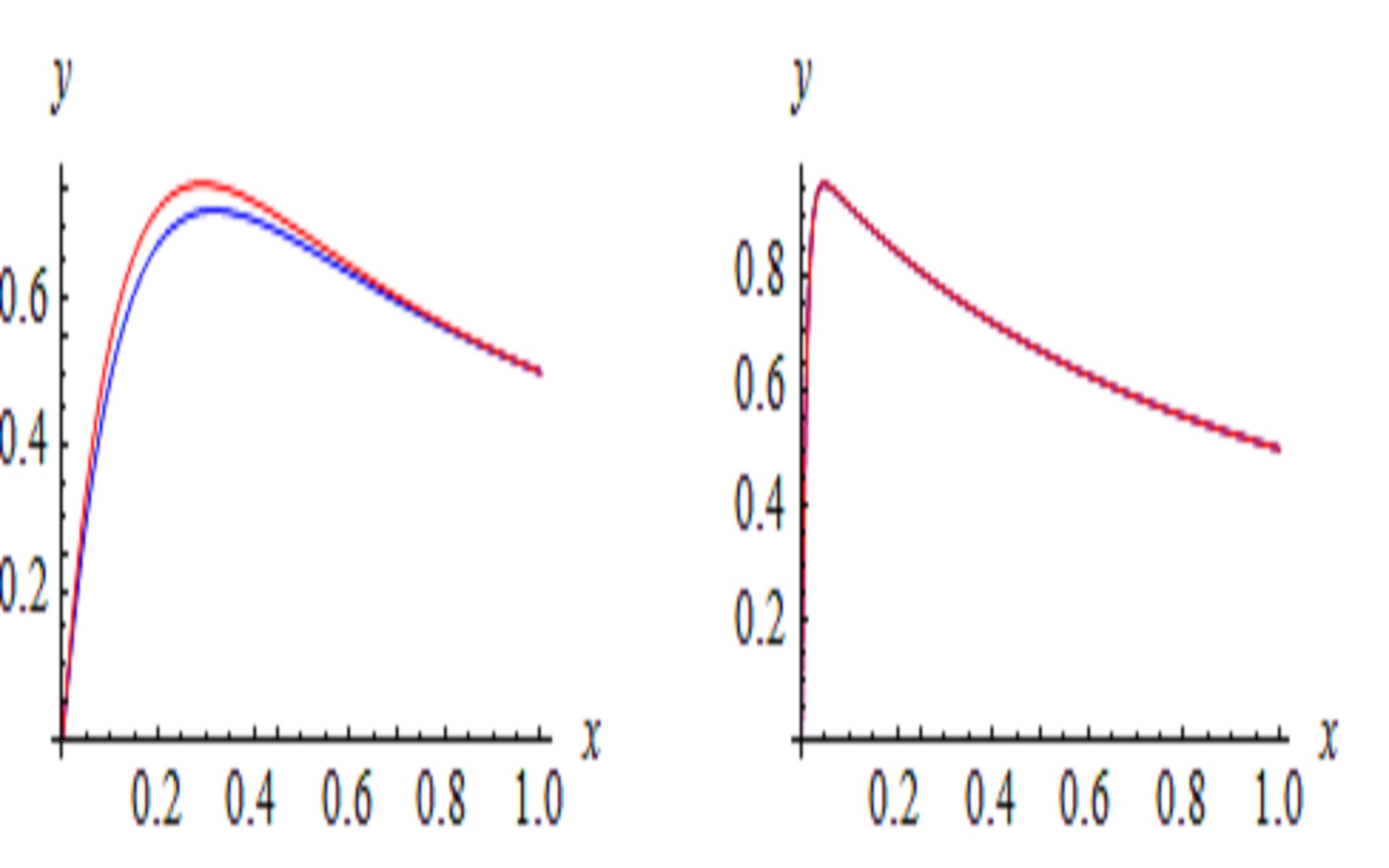}
\caption{Comparing a high precision numerical solution of the nonlinear  boundary value
problem (\ref{Eq6.21.12}), in (Red),  with the approximate multiple scale solution, in (Blue), for
$\varepsilon=0.1\,$(Left panel) and $\varepsilon=0.01$ (Right panel).\label{fig4_9}}
\label{fig9}
\end{figure}

Let us consider the following nonlinear boundary value problem
\begin{align}
\varepsilon y^{\prime\prime}+y^{\prime}+y^{2}  & =0,\text{ \ }
0<x<1,\nonumber\\
y(0)  & =0,\nonumber\\
y(1)  & =\frac{1}{2}.\label{Eq6.21.12}
\end{align}
We recognize that the differential equation in (\ref{Eq6.21.12}) is singularly
perturbed. The problem is transformed into a regularly perturbed problem using
the transformation
\begin{align*}
x  & =\varepsilon\xi,\nonumber\\
y(x)  & =u(\frac{x}{\varepsilon}).
\end{align*}
For the function $u(\xi)$ we get the following regularly perturbed boundary
value problem
\begin{align}
u^{\prime\prime}+u^{\prime}+\varepsilon u^{2}  & =0,\text{ \ \ }0<\xi<\frac
{1}{\varepsilon},\nonumber\\
u(0)  & =0,\nonumber\\
u(\frac{1}{\varepsilon})  & =\frac{1}{2}.\label{Eq6.23.12}
\end{align}
We have previously, in section \ref{QuadraticNonlinearity}, constructed an approximate
solution to the equation in (\ref{Eq6.23.12}), which is uniform for $\xi\lesssim\varepsilon^{-2}$.
\begin{equation}
u(\xi)=A(\xi)+B(\xi)e^{-\xi}-\varepsilon\frac{1}{2}B^{2}(\xi)e^{-2\xi
}+O(\varepsilon^{2}),\label{Eq6.24.12}
\end{equation}
where the amplitudes $A(\xi)$ and $B(\xi)$ satisfy the equations
\begin{align}
\frac{dA}{d\xi}  & =-\varepsilon A^{2}-2\varepsilon^{2}A^{3},\nonumber\\
\frac{dB}{d\xi}  & =2\varepsilon AB+2\varepsilon^{2}A^{2}B.\label{Eq6.25.12}
\end{align}
From the boundary conditions on $u(\xi)$, we get
\begin{align}
u(0)  & =0,\text{ \ \ }\Longleftrightarrow\text{ \ \ }A(0)+B(0)-\varepsilon
\frac{1}{2}B^{2}(0)=0,\text{ }\nonumber\\
u(\frac{1}{\varepsilon})  & =\frac{1}{2},\text{ \ \ }\Longleftrightarrow
\text{\ \ \ }A(\frac{1}{\varepsilon})+B(\frac{1}{\varepsilon})e^{-\frac
{1}{\varepsilon}}-\varepsilon\frac{1}{2}B^{2}(\frac{1}{\varepsilon}
)e^{-\frac{2}{\varepsilon}}=\frac{1}{2}.\label{Eq6.26.12}
\end{align}
The equations (\ref{Eq6.24.12}),(\ref{Eq6.25.12}) and (\ref{Eq6.26.12}) can now be used
to design an efficient numerical algorithm for finding the solution to the
boundary value problem. We do this by defining a function $F(B_{0})$ by
\[
F(B_{0})=A(\frac{1}{\varepsilon})+B(\frac{1}{\varepsilon})e^{-\frac
{1}{\varepsilon}}-\varepsilon\frac{1}{2}B^{2}(\frac{1}{\varepsilon}
)e^{-\frac{2}{\varepsilon}}-\frac{1}{2},
\]
where the functions $A(\xi)$ and $B(\xi)$ are calculated by solving the system
(\ref{Eq6.25.12}) with initial conditions
\begin{align}
A(0)  & =-B_{0}+\varepsilon\frac{1}{2}B_{0}^{2},\nonumber\\
B(0)  & =B_{0}.\label{Eq6.27.12}
\end{align}
Using Newton iteration, we find a value of $B_{0}$ such that
\[
F(B_{0})=0.
\]
Inserting this value of $B_{0}$ into the formulas for the initial conditions
(\ref{Eq6.27.12}), calculating the amplitudes $A(\xi)$, $B(\xi)$ from
(\ref{Eq6.25.12}) and inserting $A(\xi)$ and $B(\xi)$ into the formula (\ref{Eq6.24.12}), gives
us a solution to the initial value problem (\ref{Eq6.21.12}). In figure \ref{fig4_9} we compare
a high precision numerical solution of (\ref{Eq6.21.12}) with our approximate
multiple scale solution for $\varepsilon=0.1$ and $\varepsilon=0.01$.

\noindent Apart from being able to use the amplitude equations to construct an
efficient, purely numerical, algorithm for solving the boundary value problem,
it is also possible do quite a lot of analytic work on the amplitude equations (\ref{Eq6.25.12}). It is
fairly easy to find an explicit formula for $B$ as a function of $A$, it
involves nothing more fancy than using partial fractions. It is also possible
to find an implicit solution for the function $A$, also using partial fractions.

\subsection{The multiple scale method for weakly nonlinear PDEs.}

It is now finally time to start applying the multiple scale method to partial
differential equations. The partial differential equations that are of
interest in the science of linear and nonlinear wave motion are almost always
hyperbolic, dispersive and weakly nonlinear. We will therefore focus all our
attention on such equations. The multiple scale method is however in no way
restricted to equations of this type.

\subsubsection{A quadratically perturbed Klein-Gordon equation}

Let us consider the equation
\begin{equation}
\partial_{tt}u-\partial_{xx}u+u=\varepsilon u^{2}.\label{Eq5.1.12}
\end{equation}
Inspired by our work on ordinary differential equations, we introduce a
function $h(x_{0},t_{0},x_{1},t_{1},...)$ such that
\begin{equation}
u(x,t)=h(x_{0},t_{0},x_{1},t_{1},...)|_{t_{j}=\varepsilon^{j}t,x_{j}
=\varepsilon^{j}x},\label{Eq5.2.12}
\end{equation}
is a solution of (\ref{Eq5.1.12}). The derivatives turns into
\begin{align}
\partial_{t}  &  =\partial_{t_{0}}+\varepsilon\partial_{t_{1}}+\varepsilon
^{2}\partial_{t_{2}}+...\,\,,\nonumber\\
\partial_{x}  &  =\partial_{x_{0}}+\varepsilon\partial_{x_{1}}+\varepsilon
^{2}\partial_{x_{2}}+...\;\;,\label{Eq5.3.12}
\end{align}
and for $h$ we use the expansion
\begin{equation}
h=h_{0}+\varepsilon h_{1}+\varepsilon^{2}h_{2}+...\;\;.\label{Eq5.4.12}
\end{equation}
Inserting (\ref{Eq5.2.12}),(\ref{Eq5.3.12}) and (\ref{Eq5.4.12}) and expanding
everything in sight, we get
\begin{gather*}
(\partial_{t_{0}}+\varepsilon\partial_{t_{1}}+\varepsilon^{2}\partial_{t_{2}
}+...)(\partial_{t_{0}}+\varepsilon\partial_{t_{1}}+\varepsilon^{2}
\partial_{t_{2}}+...)\nonumber\\
(h_{0}+\varepsilon h_{1}+\varepsilon^{2}h_{2}+...)-\nonumber\\
(\partial_{x_{0}}+\varepsilon\partial_{x_{1}}+\varepsilon^{2}\partial_{x_{2}
}+...)(\partial_{x_{0}}+\varepsilon\partial_{x_{1}}+\varepsilon^{2}
\partial_{x_{2}}+...)\nonumber\\
(h_{0}+\varepsilon h_{1}+\varepsilon^{2}h_{2}+...)+(h_{0}+\varepsilon
h_{1}+\varepsilon^{2}h_{2}+...)\nonumber\\
=\varepsilon(h_{0}+\varepsilon h_{1}+\varepsilon^{2}h_{2}+...)^{2},\nonumber\\
\Downarrow\nonumber\\
(\partial_{t_{0}t_{0}}+\varepsilon(\partial_{t_{0}t_{1}}+\partial_{t_{1}t_{0}
})+\varepsilon^{2}(\partial_{t_{0}t_{2}}+\partial_{t_{1}t_{1}}+\partial
_{t_{2}t_{0}})+...)\nonumber\\
(h_{0}+\varepsilon h_{1}+\varepsilon^{2}h_{2}+...)-\nonumber\\
(\partial_{x_{0}x_{0}}+\varepsilon(\partial_{x_{0}x_{1}}+\partial_{x_{1}x_{0}
})+\varepsilon^{2}(\partial_{x_{0}x_{2}}+\partial_{x_{1}x_{1}}+\partial
_{x_{2}x_{0}})+...)\nonumber\\
(h_{0}+\varepsilon h_{1}+\varepsilon^{2}h_{2}+...)+(h_{0}+\varepsilon
h_{1}+\varepsilon^{2}h_{2}+...)\nonumber\\
=\varepsilon(h_{0}^{2}+2\varepsilon h_{0}h_{1}+...),\nonumber\\
\Downarrow\nonumber
\end{gather*}
\begin{gather*}
\partial_{t_{0}t_{0}}h_{0}+\varepsilon(\partial_{t_{0}t_{0}}h_{1}
+\partial_{t_{0}t_{1}}h_{0}+\partial_{t_{1}t_{0}}h_{0})+\nonumber\\
\varepsilon^{2}(\partial_{t_{0}t_{0}}h_{2}+\partial_{t_{0}t_{1}}h_{1}
+\partial_{t_{1}t_{0}}h_{1}+\partial_{t_{0}t_{2}}h_{0}+\partial_{t_{1}t_{1}
}h_{0}+\partial_{t_{2}t_{0}}h_{0})-...\nonumber\\
\partial_{x_{0}x_{0}}h_{0}-\varepsilon(\partial_{x_{0}x_{0}}h_{1}
+\partial_{x_{0}x_{1}}h_{0}+\partial_{x_{1}x_{0}}h_{0})-\nonumber\\
\varepsilon^{2}(\partial_{x_{0}x_{0}}h_{2}+\partial_{x_{0}x_{1}}h_{1}
+\partial_{x_{1}x_{0}}h_{1}+\partial_{x_{0}x_{2}}h_{0}+\partial_{x_{1}x_{1}
}h_{0}+\partial_{x_{2}x_{0}}h_{0})\nonumber\\
+h_{0}+\varepsilon h_{1}+\varepsilon^{2}h_{2}+...\nonumber\\
=\varepsilon h_{0}^{2}+2\varepsilon^{2}h_{0}h_{1}+...\;\;,
\end{gather*}
which gives us the perturbation hierarchy
\begin{gather}
\partial_{t_{0}t_{0}}h_{0}-\partial_{x_{0}x_{0}}h_{0}+h_{0}=0,\label{Eq5.6.12}\\
\nonumber\\
\partial_{t_{0}t_{0}}h_{1}-\partial_{x_{0}x_{0}}h_{1}+h_{1}=h_{0}^{2}
-\partial_{t_{0}t_{1}}h_{0}-\partial_{t_{1}t_{0}}h_{0}\nonumber\\
+\partial_{x_{0}x_{1}}h_{0}+\partial_{x_{1}x_{0}}h_{0},\label{Eq5.7.12}\\
\nonumber\\
\partial_{t_{0}t_{0}}h_{2}-\partial_{x_{0}x_{0}}h_{2}+h_{2}=2h_{0}
h_{1}-\partial_{t_{0}t_{1}}h_{1}-\partial_{t_{1}t_{0}}h_{1}\nonumber\\
-\partial_{t_{0}t_{2}}h_{0}-\partial_{t_{1}t_{1}}h_{0}-\partial_{t_{2}t_{0}
}h_{0}+\partial_{x_{0}x_{1}}h_{1}+\partial_{x_{1}x_{0}}h_{1}\nonumber\\
+\partial_{x_{0}x_{2}}h_{0}+\partial_{x_{1}x_{1}}h_{0}+\partial_{x_{2}x_{0}
}h_{0}.\label{Eq5.8.12}
\end{gather}
For ordinary differential equations, we used the general solution to the order
$\varepsilon^{0}$ equation. For partial differential equations we can not do
this. We will rather use a finite sum of linear modes. The simplest
possibility is a single linear mode which we use here
\begin{equation}
h_{0}(t_{0},x_{0},t_{1},x_{1},...)=A_{0}(t_{1},x_{1},...)e^{i(kx_{0}-\omega
t_{0})}+(\ast).\label{Eq5.9.12}
\end{equation}
Since we are not using the general solution, we will in not be able to satisfy
arbitrary initial conditions. However, in the theory of waves this is
perfectly alright, since most of the time the relevant initial conditions are
in fact finite sums of wave packets or even a single wave packet. Such initial
conditions can be included in the multiple scale approach that we discuss in
this section. For (\ref{Eq5.9.12}) to actually be a solution to (\ref{Eq5.6.12}) we
must have
\begin{equation*}
\omega=\omega(k)=\sqrt{1+k^{2}},
\end{equation*}
which we of course recognize as the dispersion relation for the linearized
version of (\ref{Eq5.1.12}). With the choice of signs used here, (\ref{Eq5.9.12})
will represent a right-moving disturbance.

Inserting (\ref{Eq5.9.12}) into (\ref{Eq5.7.12}) we get
\begin{gather*}
\partial_{t_{0}t_{0}}h_{1}-\partial_{x_{0}x_{0}}h_{1}+h_{1}=2|A_{0}
|^{2}\nonumber\\
+A_{0}^{2}e^{2i(kx_{0}-\omega t_{0})}+A_{0}^{\ast2}e^{-2i(kx_{0}-\omega
t_{0})}\nonumber\\
+(2i\omega\partial_{t_{1}}A_{0}+2ik\partial_{x_{1}}A_{0})e^{i(kx_{0}-\omega
t_{0})}\nonumber\\
-(2i\omega\partial_{t_{1}}A_{0}^{\ast}+2ik\partial_{x_{1}}A_{0}^{\ast
})e^{-i(kx_{0}-\omega t_{0})}.
\end{gather*}
In order to remove secular terms, we must postulate that
\begin{gather}
2i\omega\partial_{t_{1}}A_{0}+2ik\partial_{x_{1}}A_{0}=0,\nonumber\\
\Updownarrow\nonumber\\
\partial_{t_{1}}A_{0}=-\frac{k}{\omega}\partial_{x_{1}}A_{0}.\label{Eq5.12.12}
\end{gather}
Here we assume that the terms
\[
e^{2i(kx_{0}-\omega t_{0})},e^{-2i(kx_{0}-\omega t_{0})}\;\;,
\]
are \textit{not} solutions to the homogenous equation
\[
\partial_{t_{0}t_{0}}h_{1}-\partial_{x_{0}x_{0}}h_{1}+h_{1}=0.
\]
For this to be true we must have
\begin{equation}
\omega(2k)\neq2\omega(k),\label{Eq5.12.1.12}
\end{equation}
and this is in fact true for all $k$. This is however not generally true for
dispersive wave equations. Whether it is true or not will depend on the exact
form of the dispersion relation for the system of interest. In the theory of
interacting waves, equality in (\ref{Eq5.12.1.12}), is called \textit{phase
matching}, and is of utmost importance.

The equation for $h_{1}$ now simplify into
\begin{equation*}
\partial_{t_{0}t_{0}}h_{1}-\partial_{x_{0}x_{0}}h_{1}+h_{1}=2|A_{0}|^{2}
+A_{0}^{2}e^{2i(kx_{0}-\omega t_{0})}+A_{0}^{\ast2}e^{-2i(kx_{0}-\omega
t_{0})}.
\end{equation*}
According to the rules of the game we need a special solution to this
equation. It is easy to verify that
\begin{equation}
h_{1}=2|A_{0}|^{2}-\frac{1}{3}A_{0}^{2}e^{2i(kx_{0}-\omega t_{0})}-\frac{1}
{3}A_{0}^{\ast2}e^{-2i(kx_{0}-\omega t_{0})},\label{Eq5.14.12}
\end{equation}
is such a special solution. Inserting (\ref{Eq5.9.12}) and (\ref{Eq5.14.12}) into
(\ref{Eq5.8.12}), we get
\begin{align*}
\partial_{t_{0}t_{0}}h_{2}-\partial_{x_{0}x_{0}}h_{2}+h_{2}  &  =(2i\omega
\partial_{t_{2}}A_{0}+2ik\partial_{x_{2}}A_{0}-\partial_{t_{1}t_{1}}
A_{0}\\
&  +\partial_{x_{1}x_{1}}A_{0}+\frac{10}{3}|A_{0}|^{2}A_{0})e^{i(kx_{0}-\omega
t_{0})}+NST+(\ast).\nonumber
\end{align*}
In order to remove secular terms we must postulate that
\begin{equation}
2i\omega\partial_{t_{2}}A_{0}+2ik\partial_{x_{2}}A_{0}-\partial_{t_{1}t_{1}
}A_{0}+\partial_{x_{1}x_{1}}A_{0}+\frac{10}{3}|A_{0}|^{2}A_{0}=0.\label{Eq5.16.12}
\end{equation}
Equations (\ref{Eq5.12.12}) and (\ref{Eq5.16.12}) constitute, as usual, an overdetermined
system. In general it is not an easy matter to verify that an overdetermined
system of partial differential equations is solvable and the methods that do
exist to address such questions are mathematically quite sophisticated. For the
particular case discussed here it is however easy to verify that the system is
in fact solvable. But, as we have stressed several times in these lecture
notes, we are not really concerned with the solvability of the system
(\ref{Eq5.12.12}), (\ref{Eq5.16.12}) for the many variable function $A_{0} $. We
are rather interested in the function $u(x,t)$ which is a solution to
(\ref{Eq5.1.12}). With that in mind, we define an amplitude
\begin{equation}
A(x,t)=A_{0}(t_{1},x_{1},...)|_{t_{j}=\varepsilon^{j}t,x_{j}=\varepsilon^{j}
x}.\label{Eq5.17.12}
\end{equation}
The solution to (\ref{Eq5.1.12}) is then
\begin{align}
u(x,t)  &  =A(x,t)e^{i(kx-\omega t)}+\varepsilon(2|A|^{2}(x,t)-\frac{1}
{3}A^{2}(x,t)e^{2i(kx_{0}-\omega t_{0})}\nonumber\\
&  -\frac{1}{3}A^{\ast2}(x,t)e^{-2i(kx_{0}-\omega t_{0})})+O(\varepsilon^{2}),\label{Eq5.18.12}
\end{align}
where $A(x,t)$ satisfy a certain amplitude equation that we will now derive.

Multiplying equation (\ref{Eq5.12.12}) by $\varepsilon$ , equation (\ref{Eq5.16.12})
by $\varepsilon^{2}$ and adding the two expressions, we get
\begin{gather}
\varepsilon(2i\omega\partial_{t_{1}}A_{0}+2ik\partial_{x_{1}}A_{0})\nonumber\\
+\varepsilon^{2}(2i\omega\partial_{t_{2}}A_{0}+2ik\partial_{x_{2}}
A_{0}-\partial_{t_{1}t_{1}}A_{0}+\partial_{x_{1}x_{1}}A_{0}+\frac{10}{3}
|A_{0}|^{2}A_{0})=0,\nonumber\\
\Downarrow\nonumber\\
2i\omega(\partial_{t_{0}}+\varepsilon\partial_{t_{1}}+\varepsilon^{2}
\partial_{t_{2}})A_{0}+2ik(\partial_{x_{0}}+\varepsilon\partial_{x_{1}
}+\varepsilon^{2}\partial_{x_{2}})A_{0}\nonumber\\
-(\partial_{t_{0}}+\varepsilon\partial_{t_{1}}+\varepsilon^{2}\partial_{t_{2}
})^{2}A_{0}+(\partial_{x_{0}}+\varepsilon\partial_{x_{1}}+\varepsilon
^{2}\partial_{x_{2}})^{2}A_{0}+\varepsilon^{2}\frac{10}{3}|A_{0}|^{2}
A_{0}=0,\label{Eq5.19.12}
\end{gather}
where we have used the fact that $A_{0}$ does not depend on $t_{0}$ and
$x_{0}$ and where the equation (\ref{Eq5.19.12}) is correct to second order in
$\varepsilon$. \ If we now evaluate (\ref{Eq5.19.12}) at $x_{j}=\varepsilon
^{j}x,t_{j}=\varepsilon^{j}t$, using (\ref{Eq5.3.12}) and (\ref{Eq5.17.12}), we
get the amplitude equation

\begin{gather}
2i\omega\partial_{t}A+2ik\partial_{x}A-\partial_{tt}A+\partial_{xx}
A+\varepsilon^{2}\frac{10}{3}|A|^{2}A=0,\nonumber\\
\Updownarrow\nonumber\\
\partial_{t}A=-\frac{k}{\omega}\partial_{x}A-\frac{i}{2\omega}\partial
_{tt}A+\frac{i}{2\omega}\partial_{xx}A+\varepsilon^{2}\frac{5i}{3\omega
}|A|^{2}A.\label{Eq5.21.12}
\end{gather}
This equation appears to have a problem since it contains a second
derivative with respect to time. The initial conditions for (\ref{Eq5.1.12}) is
only sufficient to determine $A(x,0)$. However, in order to be consistent with
the multiple scale procedure leading up to (\ref{Eq5.21.12}) we can only consider
solutions such that
\begin{gather*}
\partial_{t}A\sim-\frac{k}{\omega}\partial_{x}A\sim\varepsilon,\nonumber\\
\Downarrow\nonumber\\
\partial_{tt}A\sim\left(  \frac{k}{\omega}\right)  ^{2}\partial_{xx}
A\sim\varepsilon^{2}.
\end{gather*}
Thus we can, to second order in $\varepsilon$, rewrite the amplitude equation
as
\begin{equation}
\partial_{t}A=-\frac{k}{\omega}\partial_{x}A+\frac{i}{2\omega^{3}}
\partial_{xx}A+\varepsilon^{2}\frac{5i}{3\omega}|A|^{2}A.\label{Eq5.24.12}
\end{equation}
This is now first order in time and has a unique solution for a given initial
condition $A(x,0)$.

The multiple scale procedure demands that the amplitude $A(x,t)$ vary slowly
on scales $L=\frac{2\pi}{k},T=\frac{2\pi}{\omega}$. This means that
(\ref{Eq5.18.12}) and (\ref{Eq5.24.12}) can be thought of as a fast numerical scheme
for \textit{wave packets} solutions to (\ref{Eq5.1.12}). If these are the kind of
solutions that we are interested in, and in the theory of waves this is often
the case, it is much more efficient to use (\ref{Eq5.18.12}) and (\ref{Eq5.24.12})
rather than having to resolve the scales $L$ and $T$ by integrating the
original equation (\ref{Eq5.1.12}).

The very same equation (\ref{Eq5.24.12}) appear as leading order amplitude
equation starting from a large set of nonlinear partial differential equations
describing a wide array of physical phenomena in fluid dynamics, climate
science, laser physics etc. The equation appeared for the first time more than
70 years ago, but it was not realized at the time that the Nonlinear
Schr\"{o}dinger equation (NLS), as it is called, is very special indeed.

V. Zakharov discovered in 1974 that NLS is in a certain sense completely
solvable. He discovered a nonlinear integral transform that decompose NLS into
an infinite system of uncoupled ODE's, that in many important cases are easy
to solve. This transform is called the \textit{Scattering Transform}.

Using this transform one can find explicit formulas for solutions of NLS that
acts like particles, they are localized disturbances in a wave field that does
not disperse and they collide elastically just like particles do. The NLS equation
has a host of interesting and beautiful properties. It has for example
infinitely many quantities that are conserved under the time evolution and is
the continuum analog of a \textit{completely integrable} system of ODE's.

Many books and $\infty$- many papers have been written about this equation. In
the process of doing this, many other equations having similar wonderful
properties has been discovered. They \textit{all} appear through the use of
the method of multiple scales. However, all these wonderful properties,
however nice they are, are not robust. If we want to propagate our waves for
$t\lesssim\varepsilon^{-4}$, the multiple scale procedure must be extended to
order $\varepsilon^{3}$, and additional terms will appear in the amplitude
equation. These additional terms will destroy many of the wonderful
mathematical properties of the Nonlinear Schr\"{o}dinger equation but it will
\textit{not} destroy the fact that it is the key element in a fast numerical
scheme for wave packet solutions to (\ref{Eq5.1.12}).

\subsubsection{A fourth order PDE with a cubic nonlinearity}

Let us consider the equation
\begin{equation}
u_{tt}+u_{xx}+u_{xxxx}+u=\varepsilon u^{3}.\label{Eq5.25.12}
\end{equation}
Introducing the usual tools for the multiple scale method, we have
\begin{align*}
u(x,t)  &  =h(x_{0},t_{0},x_{1},t_{1},...)|_{t_{j}=\varepsilon^{j}
t,x_{j}=\varepsilon^{j}x},\nonumber\\
\partial_{t}  &  =\partial_{t_{0}}+\varepsilon\partial_{t_{1}}+...\;\;,\nonumber\\
\partial_{x}  &  =\partial_{x_{0}}+\varepsilon\partial_{x_{1}}+...\;\;,\nonumber\\
h  &  =h_{0}+\varepsilon h_{1}+...\;\;.
\end{align*}
Inserting these expressions into (\ref{Eq5.25.12}) and expanding we get
\begin{gather*}
(\partial_{t_{0}}+\varepsilon\partial_{t_{1}}+...)(\partial_{t_{0}
}+\varepsilon\partial_{t_{1}}+...)(h_{0}+\varepsilon h_{1}+...)+\\
(\partial_{x_{0}}+\varepsilon\partial_{x_{1}}+...)(\partial_{x_{0}
}+\varepsilon\partial_{x_{1}}+...)(h_{0}+\varepsilon h_{1}+...)+\\
(\partial_{x_{0}}+\varepsilon\partial_{x_{1}}+...)(\partial_{x_{0}
}+\varepsilon\partial_{x_{1}}+...)\\
(\partial_{x_{0}}+\varepsilon\partial_{x_{1}}+...)(\partial_{x_{0}
}+\varepsilon\partial_{x_{1}}+...)(h_{0}+\varepsilon h_{1}+...)\\
+h_{0}+\varepsilon h_{1}+...=\varepsilon(h_{0}+...)^{3},\\
\Downarrow\\
(\partial_{t_{0}t_{0}}+\varepsilon(\partial_{t_{0}t_{1}}+\partial_{t_{1}t_{0}
})+...)(h_{0}+\varepsilon h_{1}+...)+\\
(\partial_{x_{0}x_{0}}+\varepsilon(\partial_{x_{0}x_{1}}+\partial_{x_{1}x_{0}
})+...)(h_{0}+\varepsilon h_{1}+...)+\\
(\partial_{x_{0}x_{0}}+\varepsilon(\partial_{x_{0}x_{1}}+\partial_{x_{1}x_{0}
})+...)(\partial_{x_{0}x_{0}}+\varepsilon(\partial_{x_{0}x_{1}}+\partial
_{x_{1}x_{0}})+...)\\
(h_{0}+\varepsilon h_{1}+...)+h_{0}+\varepsilon h_{1}+...=\varepsilon h_{0}^{3}+...\;\;,\\
\Downarrow\\
\partial_{t_{0}t_{0}}h_{0}+\varepsilon(\partial_{t_{0}t_{0}}h_{1}
+\partial_{t_{0}t_{1}}h_{0}+\partial_{t_{1}t_{0}}h_{0})+\\
\partial_{x_{0}x_{0}}h_{0}+\varepsilon(\partial_{x_{0}x_{0}}h_{1}
+\partial_{x_{0}x_{1}}h_{0}+\partial_{x_{1}x_{0}}h_{0})+\\
\partial_{x_{0}x_{0}x_{0}x_{0}}h_{0}+\varepsilon(\partial_{x_{0}x_{0}
x_{0}x_{0}}h_{1}+\partial_{x_{0}x_{0}x_{0}x_{1}}h_{0}+\partial_{x_{0}
x_{0}x_{1}x_{0}}h_{0}\\
+\partial_{x_{0}x_{1}x_{0}x_{0}}h_{0}+\partial_{x_{1}x_{0}x_{0}x_{0}}
h_{0})+...\\
+h_{0}+\varepsilon h_{1}+...=\varepsilon h_{0}^{3}+...\;\;,
\end{gather*}
which gives us the perturbation hierarchy
\begin{gather}
\partial_{t_{0}t_{0}}h_{0}+\partial_{x_{0}x_{0}}h_{0}+\partial_{x_{0}
x_{0}x_{0}x_{0}}h_{0}+h{0}=0,\nonumber\\
\nonumber\\
\partial_{t_{0}t_{0}}h_{1}+\partial_{x_{0}x_{0}}h_{1}+\partial_{x_{0}
x_{0}x_{0}x_{0}}h_{1}+h_{1}=h_{0}^{3}\label{Eq5.28.12}\\
-\partial_{t_{0}t_{1}}h_{0}-\partial_{t_{1}t_{0}}h_{0}-\partial_{x_{0}x_{1}
}h_{0}-\partial_{x_{1}x_{0}}h_{0}\nonumber\\
-\partial_{x_{0}x_{0}x_{0}x_{1}}h_{0}-\partial_{x_{0}x_{0}x_{1}x_{0}}
h_{0}-\partial_{x_{0}x_{1}x_{0}x_{0}}h_{0}+\partial_{x_{1}x_{0}x_{0}x_{0}
}h_{0}.\nonumber
\end{gather}
For the order $\varepsilon^{0}$ equation, we choose a wave packet solution
\begin{equation}
h_{0}(x_{0},t_{0},x_{1},t_{1},...)=A_{0}(x_{1},t_{1},...)e^{i(kx_{0}-\omega
t_{0})}+(\ast),\label{Eq5.29.12}
\end{equation}
where the dispersion relation is
\begin{equation}
\omega=\sqrt{k^{4}-k^{2}+1}.\label{Eq5.30.12}
\end{equation}
Inserting (\ref{Eq5.29.12}) into (\ref{Eq5.28.12}), we get after a few algebraic
manipulations
\begin{gather*}
\partial_{t_{0}t_{0}}h_{1}+\partial_{x_{0}x_{0}}h_{1}+\partial_{x_{0}
x_{0}x_{0}x_{0}}h_{1}=\nonumber\\
(2i\omega\partial_{t_{1}}A_{0}-2ik\partial_{x_{1}}A_{0}+4ik^{3}\partial
_{x_{1}}A_{0}+3|A_{0}|^{2}A_{0})e^{i(kx_{0}-\omega t_{0})}\nonumber\\
+A_{0}^{3}e^{3i(kx_{0}-\omega t_{0})}+(\ast).
\end{gather*}
In order to remove secular terms we must postulate that
\begin{equation}
2i\omega\partial_{t_{1}}A_{0}-2ik\partial_{x_{1}}A_{0}+4ik^{3}\partial_{x_{1}
}A_{0}+3|A_{0}|^{2}A_{0}=0.\label{Eq5.32.12}
\end{equation}
But using the dispersion relation (\ref{Eq5.30.12}), we have
\[
-2ik+4ik^{3}=2i\omega\omega^{\prime},
\]
so that (\ref{Eq5.32.12}) simplifies into
\begin{equation*}
2i\omega(\partial_{t_{1}}A_{0}+\omega^{\prime}\partial_{x_{1}}A_{0}
)+3|A_{0}|^{2}A_{0}=0.
\end{equation*}
Introducing an amplitude
\[
A(x,t)=A_{0}(x_{1},t_{1},...)|_{x_{j}=e^{j}x,t_{j}=\varepsilon^{j}t},
\]
we get, following the approach from the previous example, the amplitude
equation
\begin{equation}
2i\omega(\partial_{t}A+\omega^{\prime}\partial_{x}A)=-3|A|^{2}
A.\label{Eq5.34.12}
\end{equation}
This equation together with the expansion
\begin{equation}
u(x,t)=A(t)e^{i(kx-\omega t)}+(\ast)+O(\varepsilon),\label{Eq5.35.12}
\end{equation}
constitute a fast numerical scheme for wave packet solutions to (\ref{Eq5.25.12})
for $t\lesssim\varepsilon^{-2}$. Of course, this particular amplitude equation
can be solved analytically using the method of characteristics, but as stressed earlier, this property is not
robust and can easily be lost if we take the expansion to higher order in
$\varepsilon$.

There is however one point in our derivation that we need to look more closely
into. We assumed that the term
\begin{equation}
A_{0}^{3}e^{3i(kx_{0}-\omega t_{0})},\label{Eq5.36.12}
\end{equation}
was \textit{not} a secular term. The term \textit{is} secular if
\begin{equation*}
\omega(3k)=3\omega(k).
\end{equation*}
Using the dispersion relation (\ref{Eq5.30.12}) we have
\begin{align}
\omega(3k) &  =3\omega(k),\nonumber\\
&  \Updownarrow\nonumber\\
\sqrt{81k^{4}-9k^{2}+1} &  =3\sqrt{k^{4}-k^{2}+1},\nonumber\\
&  \Updownarrow\nonumber\\
81k^{4}-9k^{2}+1 &  =9k^{4}-9k^{2}+9,\nonumber\\
&  \Updownarrow\nonumber\\
k &  =\pm\frac{1}{\sqrt{3}}.\label{Eq5.38.12}
\end{align}
Thus the term (\ref{Eq5.36.12}) \textit{can} be secular if the wave number of the
wave packet is given by (\ref{Eq5.38.12}). This is another example of the
fenomenon that we in the theory of interacting waves call phase matching. As
long as we stay away from the two particular values of the wave numbers given
in (\ref{Eq5.38.12}), our expansion (\ref{Eq5.35.12}) with (\ref{Eq5.34.12}) is uniform
for $t\lesssim\varepsilon^{-2}$. However if the wave number takes on one of
the two values in (\ref{Eq5.38.12}), non-uniformities will make the ordering of
the expansion break down for $t\sim\varepsilon^{-1}$. However this does not
mean that the multiple scale method breaks down. We only need to include a
second amplitude at order $\varepsilon^{0}$ that we can use to remove the
additional secular terms at order $\varepsilon^{1}$. We thus, instead of
(\ref{Eq5.29.12}), use the solution
\begin{align*}
h_{0}(x_{0},t_{0},x_{1},t_{1},...) &  =A_{0}(x_{1},t_{1},...)e^{i(kx_{0}
-\omega t_{0})}\nonumber\\
&  +B_{0}(x_{1},t_{1},...)e^{3i(kx_{0}-\omega t_{0})}+(\ast),
\end{align*}
where $k$ now is given by (\ref{Eq5.38.12}). Inserting this expression for
$h_{0}$ into the order $\varepsilon$ equation (\ref{Eq5.28.12}) we get, after a
fair amount of algebra, the equation
\begin{gather*}
\partial_{t_{0}t_{0}}h_{1}+\partial_{x_{0}x_{0}}h_{1}+\partial_{x_{0}
x_{0}x_{0}x_{0}}h_{1}=\\
(2i\omega\partial_{t_{1}}A_{0}-2ik\partial_{x_{1}}A_{0}+4ik^{3}\partial
_{x_{1}}A_{0}\nonumber\\
+3|A_{0}|^{2}A_{0}+6|B_{0}|^{2}A_{0}+3A_{0}^{\ast2}B_{0})e^{i(kx_{0}-\omega
t_{0})}\nonumber\\
+(6i\omega\partial_{t_{1}}B_{0}-6ik\partial_{x_{1}}B_{0}+108ik^{3}
\partial_{x_{1}}B_{0}\nonumber\\
+3|B_{0}|^{2}B_{0}+6|A_{0}|^{2}B_{0}+A_{0}^{3})e^{3i(kx_{0}-\omega t_{0}
)}\nonumber\\
+NST+(\ast).\nonumber
\end{gather*}
In order to remove secular terms we must postulate the two equations
\begin{gather}
2i\omega\partial_{t_{1}}A_{0}-2ik\partial_{x_{1}}A_{0}+4ik^{3}\partial_{x_{1}
}A_{0}\nonumber\\
+3|A_{0}|^{2}A_{0}+6|B_{0}|^{2}A_{0}+3A_{0}^{\ast2}B_{0}=0,\nonumber\\
\nonumber\\
6i\omega\partial_{t_{1}}B_{0}-6ik\partial_{x_{1}}B_{0}+108ik^{3}
\partial_{x_{1}}B_{0}\nonumber\\
+3|B_{0}|^{2}B_{0}+6|A_{0}|^{2}B_{0}+A_{0}^{3}=0.\label{Eq5.41.12}
\end{gather}
Using the dispersion relation we have
\begin{align*}
& \\
-6ik+108ik^{3}  & =2i\omega(3k)\omega^{\prime}(3k).
\end{align*}
Inserting this into the system (\ref{Eq5.41.12}), simplifies it into
\begin{align*}
2i\omega(k)(\partial_{t_{1}}A_{0}+\omega^{\prime}(k)\partial_{x_{1}}A_{0}) &
=-3|A_{0}|^{2}A_{0}-6|B_{0}|^{2}A_{0}-3A_{0}^{\ast2}B_{0},\nonumber\\
2i\omega(3k)(\partial_{t_{1}}B_{0}+\omega^{\prime}(3k)\partial_{x_{1}}B_{0})
&  =-3|B_{0}|^{2}B_{0}-6|A_{0}|^{2}B_{0}-A_{0}^{3}.
\end{align*}
Introducing amplitudes
\begin{align*}
A(x,t) &  =A_{0}(x_{1},t_{1},...)|_{x_{j}=e^{j}x,t_{j}=\varepsilon^{j}
t},\nonumber\\
B(x,t) &  =B_{0}(x_{1},t_{1},...)|_{x_{j}=e^{j}x,t_{j}=\varepsilon^{j}
t},
\end{align*}
the asymptotic expansion and corresponding amplitude equations for this case
are found to be
\begin{align*}
u(x,t) &  =A(x,t)e^{i(kx-\omega t)}\nonumber\\
&  +B(x,t)e^{3i(kx-\omega t)}+(\ast)+O(\varepsilon),\nonumber\\
2i\omega(k)(\partial_{t}A+\omega^{\prime}(k)\partial_{x}A) &  =-3|A|^{2}
A-6|B|^{2}A-3A^{\ast2}B,\nonumber\\
2i\omega(3k)(\partial_{t}B+\omega^{\prime}(3k)\partial_{x}B) &  =-3|B|^{2}
B+6|A|^{2}B+A^{3}.
\end{align*}
The same approach must be used to treat the case when we do not have exact
phase matching but we still have
\[
\omega(3k)\approx3\omega(k)
\]
It should be apparent by now that the method of multiple scales is a method that can be applied in diverse situations where a naive approach using direct perturbation expansions leads to nonuniform expansions. All examples in these lecture notes has been included because they represent interesting generic features of nonlinear PDEs and ODEs and also because the amount of algebra required to construct and solve the perturbation hierarchy is manageable. In more real-life cases the amount of algebra can be challenging unless organized in an appropriate way. In order to illustrate these remarks we have in Appendix A included a derivation of the amplitude equation for linearly polarized  light pulses propagating in a dispersive medium. For these derivations the underlying system of equations  are  the full 3D Maxwell equations.

\subsubsection{Exercises}

\paragraph{Ordinary differential equations}

For the following initial value problems for ODEs, find asymptotic expansions that are
uniform for $t\lesssim\varepsilon^{-3}$. You thus need to take the expansions to
second order in $\varepsilon$. Compare your asymptotic solution to a high
precision numerical solution of the exact problem. Do the comparison for
several values of $\varepsilon$ and show that the asymptotic expansion and the
numerical solution of the exact problem deviates when $t\gtrsim\varepsilon^{-3}$.

\begin{enumerate}
\item
\begin{align*}
\frac{d^{2}y}{dt^{2}}+y  &  =\varepsilon y^{2},\\
y(0)  &  =1\\
\frac{dy}{dt}(0)  &  =0
\end{align*}

\item
\begin{align*}
\frac{d^{2}y}{dt^{2}}+y  &  =\varepsilon(1-y^{2})\frac{dy}{dt}\\
y(0)  &  =1,\\
\frac{dy}{dt}(0)  &  =0.
\end{align*}

\item
\begin{align*}
\frac{d^{2}y}{dt^{2}}+y  &  =\varepsilon(y^{3}-2\frac{dy}{dt}),\\
y(0)  &  =1,\\
\frac{dy}{dt}(0)  &  =0.
\end{align*}

\item Let the initial value problem
\begin{align}
\frac{d^{2}y}{dt^{2}}+\frac{dy}{dt}+\varepsilon y^{2}  &  =0,\text{
\ \ }t>0,\nonumber\\
y(0)  &  =1,\nonumber\\
y^{\prime}(0)  &  =1,\label{Eq4.731}
\end{align}
be given. Design a numerical solution to this problem based on the amplitude
equations (\ref{Eq4.72.12}),(\ref{Eq4.73.12}) and the expansion (\ref{Eq4.71.12}). Compare this
numerical solution to a high precision numerical solution of (\ref{Eq4.731})
for $t\lesssim\varepsilon^{-3}$. Use several different values of $\varepsilon$ and
show that the multiple scale solution and the high precision solution starts to deviate
when $t\gtrsim\varepsilon^{-3} $.
\end{enumerate}

\paragraph{Partial differential equations}

In the following problems for PDEs, use the methods from this section to find asymptotic
expansions that are uniform for $t\lesssim\varepsilon^{-3}$. Thus all
expansions must be taken to second order in $\varepsilon$.

\begin{enumerate}
\item
\[
u_{tt}-u_{xx}+u=\varepsilon^{2}u^{3},
\]

\item
\[
u_{tt}-u_{xx}+u=\varepsilon(u^{2}+u_{x}^{2}),
\]

\item
\[
u_{tt}-u_{xx}+u=\varepsilon(uu_{xx}-u^{2}),
\]

\item
\[
u_{t}+u_{xxx}=\varepsilon u^{2}u_{x},
\]

\item
\[
u_{tt}-u_{xx}+u=\varepsilon(u_{x}^{2}-uu_{xx}).
\]

\end{enumerate}

\setcounter{equation}{0}

\section{Green's functions} 

Green's functions were first introduced by the British mathematician George Green around 1830.  They can today be found everywhere in pure and applied mathematics and physics. They appear in many different guises and tend to have different names in different domains of science. \\
To a mathematician, Green's functions are the inverse of differential operators and he will tend to call them \ttx{fundamental solutions}.\\ To a solid state physicist, Green's functions are correlation coefficients for material parameters located at different space-time points. As such, Green's functions play a starring part in solid state physics to the extent that one can say that solid state physics \ttx{is} the theory of Green's functions. \\
To an elementary particle physicist, Green's functions describe the propagation of particles and antiparticles from one space-time location to another. They are associated with internal lines in \ttx{Feynman diagrams} which is the main computational engine in elementary particle physics. In this area of science Green's functions are called \ttx{propagators}. \\
Green's functions are the subject of  many textbooks. Most textbook authors, eager to quickly start discussing important nontrivial applications,  jumps right into the fray, discussing the main ideas of the theory in a fairly complicated setting. We will eventually get there,  but will approach the subject from a simpler setting where the main properties of Green's functions can be explained in a simple way. In this simple setting it will appear as if we solve simple problems in a complicated way. And we do, but the point is not to solve these simple problems but rather to introduce all the main constructions involving  Green's functions in the simplest setting possible. In a more complex and realistic setting, beloved by textbook authors, there are really no new ideas. Everything is just more complicated. 

\subsection{Green's functions for the operator $\vb{L = -\frac{d}{dx}}$}\label{LL1}

A Green's function, $k(x; \xi)$, for the operator $L = - \frac{d}{dx}$ is a solution to the equation 
\begin{align}
L \; k(x; \xi) = \delta(x - \xi). \lbl{1.12}
\end{align}
Recall that $\delta(x-\xi)$ is not a function, but a distribution. So a Green's function is not really a function either, but should be understood to be a distribution too. But how do we differentiate distributions? And what are distributions anyway? We will address these questions later in this section of the lecture notes, but for now we will proceed in a heuristic manner \'{a} la Dirac, or in other words we play with formulas. \\
Let $ I_{\epsilon} = (- \epsilon + \xi , \epsilon + \xi)$ be a small interval, centred on $x = \xi$. Integrate \rf{1} over $I_{\epsilon}$ 
\begin{align}
\int_{I_{\xi}} dx \; L \; k(x; \xi) &= \int_{I_\epsilon} dx \; \delta(x-\xi), \nonumber\\
&\Updownarrow\nonumber\\
 - \int^{\xi + \epsilon}_{\xi - \epsilon} dx \; k'(x; \xi) &= \int^{\xi + \epsilon}_{\xi - \epsilon} dx \; \delta(x - \xi) = 1, \nonumber\\ 
&\Updownarrow\nonumber\\
 k(\xi + \epsilon; \xi) - k(\xi - \epsilon;\xi)& = -1. \lbl{3.12} 
\end{align}
The last equation holds for all $ \epsilon > 0$. Taking the limit when $\epsilon$ approaches zero, we get 
\begin{align}
k_{+}(\xi; \xi) - k_{-}(\xi; \xi) = -1, \lbl{4.12} 
\end{align}
where by the definition 
\begin{align*}
k_{\pm} (\xi;\xi) = \lim\limits_{\epsilon \rightarrow 0} k(\xi \pm \epsilon; \xi)
\end{align*}
The right hand side identity in \rf{3.12} holds because $\delta(x- \xi)$ is concentrated infinitesimally close to $x = \xi$, so that the domain outside $I_{\epsilon}$ gives no contribution 
\begin{align*}
1 = \int^{\infty}_{-\infty} dx \; \delta(x - \xi) = \int_{I_{\epsilon}} dx \; \delta(x-\xi).
\end{align*}
I am of course just playing with formulas here... \\
Using \rf{4.12}, we can now say that $k(x; \xi)$ is a Green's function for $L = - \frac{d}{dx}$ if $ k(x; \xi)$ satisfies 
\begin{align}
- k'(x;\xi)& = 0, && x \ne \xi, \lbl{7.12} \\
k_{+}(\xi; \xi)  - k_{-}(\xi; \xi) &= -1. \lbl{8.12} 
\end{align}
Note that prime here means the derivative of the function $k(x; \xi)$ with respect to it's first argument, $x$.

\noindent This is a problem we can actually solve! From \rf{7.12} we get 
\begin{align*}
k(x; \xi) =\begin{cases}
 \; \; a(\xi) \; \; \; \; \; \; x > \xi \\ \; \; b(\xi) \; \; \; \; \; \; x < \xi\\ \end{cases}, 
\end{align*}
and \rf{8.12} imposes the condition 
\begin{align}
k_{+}(\xi;\xi) - k_{-}(\xi; \xi) &= -1, \nonumber \\
&\Updownarrow \nonumber \\
 a(\xi) - b(\xi) &= -1, \nonumber \\ 
&\Updownarrow \nonumber \\ b(\xi) &= 1 + a(\xi), \nonumber 
\end{align}
and thus $k(x; \xi)$ is a Green's function for $ L = - \frac{d}{dx} $ if only if it is of the form 
\begin{align}
k(x; \xi) =\begin{cases} \; \; a(\xi) \; \; \; \; \; \;\;\;\; x>\xi \\ \; \; 1 + a(\xi) \; \; \; x<\xi \end{cases}, \lbl{10.12}
\end{align}
where $a(\xi)$ is an arbitrary function. For example, if $a(\xi) = 0$, we get 
\begin{align*}
k(x; \xi) =\begin{cases} \; \; 0 \; \; \; \; \; \; x>\xi \\ \; \; 1  \; \; \;\;\;\; x<\xi \end{cases}.
\end{align*}

\begin{figure}[htbp]
\centering
\includegraphics{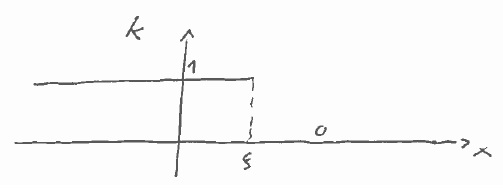}
\caption{A Green's functions for the operator $L=-\frac{d}{dx}$.}\label{L1}
\label{fig1}
\end{figure}

\noindent This function is displayed in figure (\ref{fig1}), and  is clearly a regular function and not a distribution, whatever that is! However, we only know that $ k(x; \xi)$ satisfies \rf{7.12} and \rf{8.12}, which we got from the original condition \rf{1.12} by playing with some formulas. And we certainly can't just substitute $k(x; \xi)$ from \rf{10.12} into \rf{1.12} in order to verify that it actually solves \rf{1.12}. The function $k(x; \xi)$ is not differentiable, it is not even continuous! \\
However,we will show later in these notes that there is a distribution corresponding to $k(x; \xi)$ and as a distribution it \ttx{can} be differentiated and it \ttx{will} in fact be a solution to \rf{1.12}. \\
For now we proceed heuristically and assume that the functions $k(x; \xi)$ from \rf{10.12} are solutions to \rf{1.12}, so that \rf{10.12} describes all Green's functions for the operator $L = - \frac{d}{dx}$. 

So, having  found the Green's functions for $L$, the next question is; what are they good for? 
In order to answer this question, we must introduce a certain integral identity associated with $L$. 

Let $[x_0 , x_1]$ be some interval on the real line, and let $\phi$ and $\psi$ be functions defined on the interval. Using integration by parts, we have 
$$ \intx dx \; L \; \phi \; \psi = - \intx dx \; \phi' \; \psi = - \{ \phi \; \psi \mathlarger{|}^{x_1}_{x_0} - \intx dx \; \phi \; \psi' \}. $$ 
Thus, if we define $L^{+} = \frac{d}{dx}$, we get the integral identity 
\begin{align}
\intx dx \; \{  L \; \phi \; \psi - \phi \; L^{+} \; \psi \} = - \phi \; \psi \mathlarger{|}^{x_1}_{x_0}. \lbl{11.12} 
\end{align}
This shows that $L^{+}$ is the \ttx{formal adjoint} of $L$. \\ 
We will now put \rf{11.12} to use, and with that goal in mind, let us consider the differential equation 
\begin{align}
L^{+} \; f(x) &= h(x) && x_0 < x < x_1,\nonumber \\ 
&\Updownarrow \nonumber \\ f'(x) &= h(x) && x_0 < x < x_1.  \lbl{12.12} 
\end{align}
Let now $\psi = f(x)$ be any solution of \rf{12.12} and let $ \phi = k(x; \xi)$, in the integral identity \rf{11.12}. Of course, we derived this identity using integration by parts, which makes assumptions about the smoothness of $\phi$ and $\psi$. We will disregard this fact and assume that the identity holds for any $\phi$ and $\psi$ we like. \\
From \rf{11.12} we get 
\begin{align*}
\intx dx \; \{  - k'(x; \xi) \; f(x) - k(x; \xi) \; h(x)  \} \nonumber \\
=-k(x; \xi) \; f(x) \mathlarger{|}^{x_1}_{x_0},
\end{align*}
which upon using \rf{12.12} gives us 
\begin{align}
f(\xi) = \intx dx \; k(x; \xi) \; h(x) - k(x_1;\xi) \; f(x_1) + k(x_0; \xi) \; f(x_0). \lbl{14.12} 
\end{align}
Formula \rf{14.12} introduces the first great theme in the theory of Green's functions: \\ \textit{Green's functions can be used to derive integral identities for solutions to differential equations that relate values of the solutions inside a domain to their values on the boundary of the domain.}\label{FirstGreatTheme} \\
Note that \rf{14.12} does not give us solutions to \rf{12.12}, it merely shows that values of solutions inside $[x_0, x_1]$ are related to values on the boundary, which in this case consists of two points $\{ x_0, x_1 \}$, in a particular way. \\
However \rf{14.12} can be used as a starting point for finding solutions to \rf{12.12} in two quite distinct ways. \\
Let us look for a solution that satisfy the boundary condition.
\begin{align*}
f(x_0) = f_0.
\end{align*}
For this kind of solution the integral identity \rf{14.12} gives us 
\begin{align}
f(\xi) = \intx dx \; k(x; \xi) \; h(x) - k(x_1; \xi) \; f(x_1) + k(x_0; \xi) \; f_0. \lbl{17.12} 
\end{align}
This identity \ttx{still} does not give us a solution because $f(x_1)$ on the right hand side is unknown. \\ 
However \ttx{if} we can find a Green's function that satisfy the condition 
\begin{align}
k (x_1; \xi) = 0, \lbl{18.12} 
\end{align}
\ttx{then} we do get a solution from \rf{17.12}. The solution is 
\begin{align}
f(\xi) = \intx dx \; k(x; \xi) \; h(x) + k (x_0; \xi) \; f_0. \lbl{19.12} 
\end{align}
From \rf{10.12} we see that \rf{18.12} holds if 
\begin{align*}
a(\xi) = 0, 
\end{align*}
and thus the required Green's function is 
\begin{align*}
k(x; \xi) =\begin{cases}
 \; \; 0 \; \; \; \; \; x>\xi \\ \; \; 1 \; \; \; \; \; x< \xi
\end{cases},
\end{align*}
and the solution \rf{19.12} is 
\begin{align*}
f(\xi) = \intxi dx \; h(x) + f_0.
\end{align*}
This is of course exactly what we would get if we applied the fundamental theorem of calculus to 
\begin{align*}
f'(x) &= h(x) && x_0 < x < x_1, \nonumber \\
f(x_0) &= f_0.
\end{align*}
because 
\begin{align}
\intxi dx \; f'(x) &= \intxi dx \; h(x), \nonumber \\
&\Updownarrow \nonumber \\ \;\;\;f(\xi) - f(x_0) &= \intxi dx \; h(x), \nonumber \\ 
&\Updownarrow \nonumber \\ \;\;\;f(\xi) &= \intxi dx \; h(x) + f_0. \nonumber 
\end{align} 
As another application of \rf{14.12} let us try to find a solution to \rf{12.12} that satisfy the condition 
\begin{align}
\frac{1}{2} \; (f(x_0) + f(x_1)) = \bar{f}, \lbl{24.12} 
\end{align}
where $\bar{f}$ is given. Introduce $a$ and $b$ through 
\begin{align*}
a &= \frac{1}{2} \; (f(x_0) + f(x_1)), \nonumber \\ 
b &= \frac{1}{2} \; (f(x_0) - f(x_1)), \nonumber \\ 
&\Updownarrow \nonumber \\
 f(x_0) &= a+b, \nonumber \\
f(x_1) &= a - b. 
\end{align*}
Inserting this into \rf{14.12} and rearranging terms we get 
\begin{align}
f(\xi) &= \intx dx \; k(x; \xi) \; h(x) + (k(x_0; \xi) - k(x_1; \xi)) \; a \nonumber \\ 
&+ (k(x_1; \xi)+ k(x_0; \xi))b. \lbl{26.12} 
\end{align}
The condition \rf{24.12} implies that $a = \bar{f}$, is given. Since $b$ is not given we need a Green's function that satisfy the condition 
\begin{align*}
k(x_1; \xi) + k(x_0;\xi) &= 0, \nonumber \\ 
&\Updownarrow \nonumber \\ \;\;\; a(\xi) + 1 + a(\xi) &= 0, \nonumber \\ 
&\Updownarrow \nonumber \\ \;\;\; a(\xi) &= - \frac{1}{2}. 
\end{align*}
The correct Green's function is  thus
\begin{align}
k(x; \xi) =\begin{cases} 
- \frac{1}{2} \; \; \;\;\;\;\;\;\; \; \; \; x> \xi \\ \; \; \frac{1}{2} \; \;\;\;\;\;\;\;\;\; \; \; \; x<\xi  
\end{cases}, \lbl{28.12}
\end{align}
and the solution of \rf{12.12} that satisfy \rf{24.12} is from \rf{26.12} 
\begin{align*}
f(\xi) = \frac{1}{2} \; \intxi dx \; h(x) - \frac{1}{2} \; \int^{x_1}_{\xi} dx \; h(x) + \bar{f},
\end{align*}
because for the Green's function \rf{28.12} we have 
\begin{align*}
k(x_0; \xi) - k(x_1; \xi) = \frac{1}{2} - (- \frac{1}{2}) = 1.	
\end{align*}
Thus we see that by making $k(x; \xi)$ satisfy the appropriate boundary condition we can find a solution of \rf{12.12} satisfying any chosen boundary conditions by using our integral identity \rf{14.12}. This is one of the ways we can use \rf{14.12} to find solutions to the differential equation \rf{12.12}. \\
For our chosen operator $L = - \frac{d}{dx}$, Green's functions that satisfy various boundary conditions are trivial to find. However, for more realistic and complex cases, finding the required Green's functions can be very hard. Mostly this must be done analytically,  because the presence of the Dirac delta in the differential equation implies that numerical methods are of limited use here. \\
This leads us to the second way we can use the integral identity \rf{14.12} to find solutions to the differential equation \rf{12.12}. Let us make the choice $a(\xi) = 0$ so that the Green's function is 
\begin{align*}
k(x; \xi) =\begin{cases}
\; \; 0 \; \; \;\;\;\;\;\;\; \; \; x>\xi \\ \; \; 1 \; \; \; \;\;\;\;\;\;\; \; x< \xi.
\end{cases}
\end{align*}
For this choice 
\begin{align*}
k (x_1 ; \xi) &= 0, \nonumber \\ 
k(x_0; \xi) & = 1,
\end{align*}
so that \rf{14.12} reads 
\begin{align}
f(\xi) = \intxi dx \; h(x) + f(x_0). \lbl{33.12} 
\end{align}
We now let $\xi$ approach $x_0$ from above and $x_1$ from below 
\begin{align}
x \rightarrow x_0 \; \; \Rightarrow \; \; f(x_0) &= f(x_0) \; \; - \text{trivially true},\nonumber \\ 
x \rightarrow x_1 \; \; \Rightarrow \; \; f(x_1) &= \intx dx \; h(x) + f(x_0). \lbl{35.12} 
\end{align}
Let us say we are looking for a solution that satisfy 
\begin{align*}
f(x_1) = f_1, 
\end{align*}
then \rf{35.12} is an \ttx{equation} whose solution determines the \ttx{unknown} boundary value $f(x_0)$ 
\begin{align}
f(x_0) + \intx dx \; h(x) &= f_1,\nonumber \\
\Updownarrow\nonumber\\
 f(x_0) &= f_1 - \intx dx h(x). \lbl{38.12}
\end{align}
Equation \rf{38.12} is called a \ttx{boundary integral equation}. The 'integral' part of the name will be clear when we move to a more realistic situation where the domain is 2D or 3D and the boundary 1D or 2D. The analogue to \rf{38.12} will in these cases actually be an integral equation for functions defined on the boundary of the domain. \\
We now insert the solution of the boundary integral equation \rf{38.12} into the integral identity \rf{33.12}. 
\begin{align*}
f(\xi) = \intxi dx \; h(x) + f_1 - \intx dx \; h(x),
\end{align*}
and get the solution $f(\xi)$ of \rf{12.12} that satisfies the boundary condition 
\begin{align*}
f(x_1) = f_1. 
\end{align*}
Observe that in this case, we did not need to pick a particular $k(x; \xi)$ by posing a boundary condition. In fact, we can choose almost any Green's function we want. The price we pay is that at some point we must solve a boundary integral equation. This is however easier to do for a complex boundary than trying to construct a Green's function satisfying some particular boundary condition. Here, of course, both approaches are trivial to deploy because of the simplicity of the operator $L = -\frac{d}{dx}$,  and the domain $[x_0, x_1]$. \\
The two different ways one can find solutions to the differential equation \rf{12.12} from the integral identity \rf{14.12} form the second and third great themes in the theory of Green's functions:\\
\textit{Green's functions satisfying particular boundary conditions can be used to find integral representations of solutions to initial/boundary value problems for differential equations}\label{SecondGreatTheme} \\
and 

\noindent \textit{Green's functions can be used to derive boundary integral equations whose solutions will give integral representations of solutions to initial/boundary value problems for differential equations.}\label{ThirdGreatTheme} \\ 
The choice of Green's function for a given differential equation is sometimes determined by the physical context of the equation. \\
Let us consider the problem 
\begin{align*}
\frac{dx}{dt} = v(t) && t_0 < t < t_1,
\end{align*}
where $t$ is time, $v(t)$ is the velocity of a particle and $x(t)$ it's position. The integral identity \rf{14.12} applies and we have 
\begin{align}
x(t) = \intt dt' \; k(t'; t) \; v(t') - k(t_1;  t) \; x(t_1) + k(t_0 ; t) \; x(t_0), \lbl{42.12} 
\end{align}
where $k(t'; t)$ is a Green's function for $L = -\frac{d}{dt'}$.

  Let us first choose 
\begin{align*}
k(t';t) =\begin{cases}
\; \; 0 \; \; \; \;\;\;\;\;\;\; \; t' > t \\ \; \; 1 \; \;\;\;\;\;\;\; \; \; \; t' < t
\end{cases}, 
\end{align*}
then \rf{42.12} turns into
\begin{align}
x(t) = \int_{t_0}^t dt' \; v(t') + x(t_0). \lbl{44.12}
\end{align}
This tells us that the current position of the particle depends on the current and past values of the velocity. \\
This makes physical sense; the past influences the future. Formula \rf{44.12} is an embodiment of \ttx{causality}. Green's functions, for time dependent ODE's and PDE's, that leads to formulas respecting causality are in general called \ttx{retarded} Green's functions. \\
Let us next choose 
\begin{align*}
k(t';t) =\begin{cases}
  -1 \; \; \; \; \; \; \;\;\qquad \qquad\; \;\;\;\;\;\;\; \; t' > t \\ 
\;\; 0 \; \; \; \; \; \; \;\qquad  \qquad\;\; \; \;\;\;\;\;\;\;\; t' < t 
\end{cases}.
\end{align*}
Then \rf{42.12} becomes 
\begin{align}
x(t) = - \int_{t}^{t_1} dt' \; v(t') + x(t_1), \lbl{46.12}
\end{align}
and this formula tells us that the current position of the particle depends on the current and future values of the velocity. This does not make physical sense; the future influences the past! The formula \rf{46.12} is an embodiment of \ttx{non-causality}. Green's functions that leads to non-causal formulas are called \ttx{advanced} Green's functions. Advanced Green’s functions play an important role in some areas of applied mathematics and theoretical physics, in particular they play a crucial role in the Standard Model of elementary particles.

Believe it or not, but the main ideas in the theory of Green's functions have now been introduced. What remains to do, is to solidify these ideas by looking at several, progressively more complicated, cases. During this work I will also introduce the main ideas from the theory of distributions.

\subsection{Green's functions for the operator $\vb{L=-\frac{d^2}{dx^2}}$}\label{L2}

A Green's function, $k(x; \xi)$, for the operator $L= - \frac{d^2}{dx^2}$ is a solution to the equation
\begin{align*}
L \; k(x; \xi) = \delta(x-\xi). 
\end{align*}
In order to derive an equation we can actually solve, we proceed heuristically like on page 2. 
\begin{align}
\int_{I_{\epsilon}} dx \; L \; k(x; \xi) &= \int_{I_{\epsilon}} dx \; \delta(x - \xi), \nonumber\\ 
&\Updownarrow\nonumber\\
- \int_{\xi - \epsilon}^{\xi + \epsilon} dx \; k''(x; \xi) &= \int_{\xi + \epsilon}^{\xi - \epsilon} dx \; \delta(x - \xi) = 1,\nonumber \\ 
&\Updownarrow\nonumber\\
k'(\xi + \epsilon; \epsilon) - k'(\xi - \epsilon; \xi) &= -1,\nonumber \\
&\Updownarrow\nonumber\\
k'_{+}(\xi; \xi) - k'_{-}(\xi; \xi) &= -1, \lbl{51}
\end{align}
where 
\begin{align*}
k_{\pm}'(\xi; \xi) = \lim\limits_{\epsilon \rightarrow 0} \; k'(\xi \pm \epsilon; \xi)
\end{align*}
Equation \rf{51} tells us that $k'(x; \xi) $ has a jump discontinuity at $x=\xi $. We postulate that $k (x;\xi )$ is continuous at $x= \xi$. Thus $k(x; \xi)$ is a Green's function of $L = -\dsx{} $ if 
\begin{align}
-k''(x; \xi) &= 0 && x\ne \xi, \nonumber\\
k(\xi, \xi) - k(\xi, \xi) &= 0, \lbl{54} \\
k'_+(\xi, \xi) - k_-'(\xi, \xi) &= -1. \lbl{55}
\end{align}
These equations we can now solve. From \rf{53} we get 
\begin{align}
k(x; \xi) =\begin{cases}
a(\xi) \; x + b(\xi) \; \; \; \; \; \qquad \qquad \; \; \; x>\xi \\
c(\xi) \; x + d(\xi) \; \; \; \;\qquad \qquad \; \; \; \; x< \xi 
\end{cases}. \lbl{56}
\end{align}
Equations \rf{54}, \rf{55} applied to the functions $k(x; \xi)$ in \rf{56} give 
\begin{align}
a(\xi) \; \xi + b(\xi) - c(\xi) \; \xi - d(\xi) &= 0, \nonumber \\
a(\xi) - c(\xi) &= -1. \lbl{58}
\end{align}
Equations \rf{58} are easy to solve, and we get 
\begin{align*}
c(\xi) = 1 + a(\xi), \nonumber \\ 
d(\xi) = b(\xi) - \xi,
\end{align*}
and thus Green's functions to $L = -\dst{}$, are of the form 
\begin{align}
k(x; \xi) =\begin{cases}
a(\xi) \; x + b(\xi) \; \; \;\;\; \; \qquad \qquad \; \; \; \; \; \; \; x>\xi \\
a(\xi) \; x + b(\xi) + x - \xi \; \; \; \; \; \; \; \qquad \; x<\xi
\end{cases},\lbl{61}
\end{align}
where $a(\xi)$ and $b(\xi)$ are arbitrary. For example if $ a(\xi) = b(\xi) = 0$ we get 
\begin{align*}
k(x; \xi) =\begin{cases}
0 \; \; \; \; \; \; \; \; \; \qquad \qquad \qquad \; \; x>\xi \\
x - \xi \; \; \; \; \; \;\;\;\; \qquad \qquad \; \; \; x<\xi
\end{cases}.
\end{align*}
This Green's function is illustrated in figure \ref{fig2}.

\begin{figure}[htbp]
\centering
\includegraphics{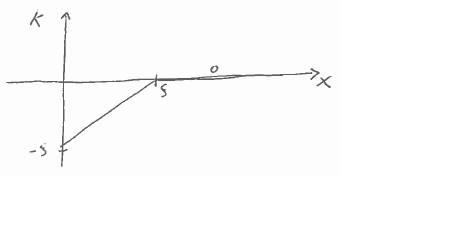}
\caption{A Green's functions for the operator $L=-\frac{d^2}{dx^2}$}
\label{fig2}
\end{figure}

\noindent In order to make use of the Green's function, we need it's associated integral identity.

  Let $[x_0, x_1]$ be some interval on the real line and let $\phi$ and $\psi$ be smooth functions. Using integration by parts we have 
\begin{align*} 
\intx dx \; L \; \phi \; \psi &= - \intx dx \; \phi'' \; \psi = -\phi' \; \psi \mathlarger{|}^{x_1}_{x_0} + \intx dx \; \phi' \; \psi'  \nonumber \\
&= - \phi' \; \psi \; \mathlarger{|}^{x_1}_{x_0} + \phi \; \psi'\mathlarger{|}^{x_1}_{x_0} - \intx dx \; \phi \; \psi'',
 \end{align*}
and thus we get the integral identity 
\begin{align}
\intx dx \; \{ L \; \phi \; \psi - \phi \; L \; \psi  \} = ( \phi \; \psi' - \phi' \; \psi)\mathlarger{|}^{x_1}_{x_0}. \lbl{63.12}
\end{align}
This identity reflects the fact that $L$ is formally self-adjoint. \\
We will now put the Green's functions \rf{61.12} to use, and introduce for this purpose the differential equation 
\begin{align}
L \; f(x) &= h(x) && x_0 < x < x_1, \nonumber \\ 
\Updownarrow \nonumber \\- f''(x) &= h(x) && x_0 < x < x_1. \lbl{64.12}
\end{align}
Inserting $\psi = f(x)$ and $ \phi = k(x; \xi)$ into the integral identity \rf{63.12}, gives us
\begin{align}
&\intx dx \{-k''(x; \xi) \; f(x) - k(x; \xi) \; h(x)\} \nonumber \\
&= ( k(x; \xi) \; f'(x) - k'(x; \xi) \; f(x) )\mathlarger{|}^{x_1}_{x_0}, \nonumber
\end{align}
which, upon using the fact that $k(x; \xi)$ is a Green's function for $L = -\dsx{} $, gives us the integral identity 
\begin{align}
f(\xi) &=  \intx dx \; k(x; \xi) \; h(x) \nonumber \\ 
&+ (k(x; \xi) \; f'(x) - k' (x; \xi) \; f(x)) \mathlarger{|}^{x_1}_{x_0}. \lbl{65.12} 
\end{align}
Like before, \rf{65.12} does not give us a solution to \rf{64.12}, but for any solution of \rf{64.12}, it gives us a relation between values of the solution inside the domain and values of the solution and its derivative on the boundary of the domain. The identity \rf{65.12} is evidently more complex than \rf{14.12} but the ideas behind them are are exactly the same. The only difference is the nature of $L$ in the two cases. \\
Equation \rf{64.12} is a second order differential equation and in order to get a unique solution, two pieces of boundary data must be specified. \\
Let us assume that Dirichlet data is given
\begin{align}
f(x_0) = f_0, && f(x_1) = f_1. \lbl{66.12}
\end{align}
If the unknown data $f'(x_0)$ and $f'(x_1)$ could be made to vanish, \rf{65.12} would give us a solution to the boundary value problem \rf{64.12}, \rf{66.12}. We can achieve this by assuming that the Green's function satisfy the boundary conditions 
\begin{align*}
&k(x_1;\xi) = 0, \nonumber\\
&k(x_0; \xi) = 0.  
\end{align*}
Using the general description of $k(x; \xi)$ from \rf{61.12} we get the equations 
\begin{align*}
& a(\xi) \; x_1 + b(\xi) = 0,\nonumber \\
& a(\xi) \; x_0 + b(\xi) + x_0 - \xi = 0,
\end{align*}
and solving these equations,  we find 
\begin{align*}
& a(\xi) = \frac{x_0 - \xi}{x_1 - x_0},\nonumber \\
& b(\xi) = -\frac{x_1 \; (x_0-\xi)}{x_1 - x_0}. 
\end{align*}
Inserting these into the formula for $k(x; \xi)$ from \rf{61} we get 
\begin{align}
k(x;\xi) =\begin{cases}
\frac{(x_0-\xi)(x-x_1)}{(x_1 - x_0)} \qquad \qquad \qquad x>\xi \\
\frac{(x_0 - x)(\xi - x_1)}{(x_1 - x_0)} \qquad \qquad \qquad x<\xi 
\end{cases}, \lbl{71.12}
\end{align}
and from \rf{71.12} we get 
\begin{align}
k'(x;\xi) =\begin{cases} 
-\frac{(\xi - x_0)}{(x_1 - x_0)} \qquad \qquad \qquad x>\xi \\
-\frac{(\xi - x_1)}{(x_1 - x_0)} \qquad \qquad \qquad x<\xi 
\end{cases}. \lbl{72.12}
\end{align}
Inserting \rf{71.12} and \rf{72.12} into \rf{65.12} will give us the solution of \rf{64.12} that satisfies the boundary conditions \rf{66.12}. Observe that 
\begin{align}
&\intx dx \; k(x;\xi) \; h(x) \nonumber \\
&= \intxi  dx \frac{(x_0 - x)(\xi - x_1)}{(x_1 - x_0)} \; h(x) + \intix dx \; \frac{(x_1 - x)(\xi - x_0)}{(x_1 - x_0)} \; h(x) \nonumber \\
&= \frac{(\xi - x_1)}{(x_1 - x_0)} \intxi dx \; (x_0 - x) \; h(x) + \frac{(\xi - x_0)}{(x_1 - x_0)} \; \intix dx \; (x_1 - x) \; h(x), \nonumber 
\end{align}
and from the expression \rf{72.12} we have
\begin{align}
k'(x_1; \xi) = - \frac{(\xi - x_0)}{(x_1 - x_0)},\nonumber \\ 
k'(x_0; \xi)  = - \frac{(\xi - x_1)}{(x_1 - x_0)}. \lbl{74.12}
\end{align}
From \rf{65.12} we thus get the solution
\begin{align}
f(\xi) &= \frac{(\xi - x_1)}{(x_1 - x_0)} \; \intxi dx \; (x_0 - x) \; h(x) + \frac{(\xi - x_0)}{(x_1 - x_0)} \; \intix dx \; (x_1 - x) \; h(x) \nonumber \\ 
&+ \frac{(\xi - x_0)}{(x_1 - x_0)} \; f_1 - \frac{(\xi - x_1)}{(x_1 - x_0)} \; f_0. \lbl{75.12} 
\end{align}
Since the calculations leading up to \rf{75.12} involved some play with formulas, it would be useful to verify directly that \rf{75.12} defines a function that solves the boundary value problem \rf{64.12},\rf{66.12}.

  The fact that it satisfies the boundary conditions is evident. Let us next verify that it also satisfies the differential equation \rf{64.12}
\begin{align}
f'(\xi) &= \frac{1}{L} \; \intxi dx (x_0 - x) \; h(x) + \frac{(\xi - x_1)}{L} \; (x_0 - \xi) \; h(\xi) \nonumber \\ 
&+ \frac{1}{L} \; \intix dx \; (x_1 - x) \; h(x) -\frac{(\xi - x_0)}{L} \; (x_1 - \xi) \; h(\xi) \nonumber \\
&= \frac{1}{L} \; \{ \; \intxi dx (x_0 - x) \; h(x) + \intix dx \; (x_1 - x) \; h(x) \; \}, \nonumber \\
\Downarrow \; \;\nonumber\\
 f''(\xi) &= \frac{1}{L} \; \{\; (x_0 - \xi) \; h(\xi) - (x_1 - \xi) \; h(\xi) \; \} \nonumber \\ 
&= \frac{(x_0 - x_1)}{L} \; h(\xi) = - h(\xi), \nonumber \\  
\Downarrow\nonumber\\
 - f''(\xi) &= h(\xi). \nonumber 
\end{align}
Note that we have here defined $L = x_1 - x_0$. 
As an example, for the special case 
\begin{align*}
h(x) = \bar{h} = \text{const}, 
\end{align*}
we get 
\begin{align}
\intxi dx \; (x_0 - x) \; h(x) &= -\frac{1}{2} \bar{h} \; (x_0 - x)^2\mathlarger{|}^{\xi}_{x_0} = -\frac{1}{2} \; \bar{h} \; (\xi - x_0)^2, \nonumber \\
\intix dx \; (x_1 - x) \; h(x) &= -\frac{1}{2} \bar{h} \; (x_1 - x)^2\mathlarger{|}_{\xi}^{x_1} = \frac{1}{2} \; \bar{h} \; (\xi - x_1)^2, \nonumber
\end{align}
and thus we have the explicit solution
\begin{align}
 f(\xi) &= - \frac{1}{2} \; \bar{h} \; \frac{(\xi - x_0)^2(\xi - x_1)}{(x_1 - x_0)} + \frac{1}{2} \; \bar{h} \; \frac{(\xi - x_1)^2(\xi - x_0)}{(x_1 - x_0)} \nonumber\\
 &+ \frac{(\xi - x_0)}{(x_1 - x_0)} \; f_1 - \frac{(\xi - x_1)}{(x_1 - x_0)} \; f_0. \nonumber
 \end{align}
 
\noindent Again, by direct differentiation, one can verify that this is an explicit solution of the boundary value problem \rf{64.12}, \rf{66.12}.

  Let us next assume that we have Cauchy data given at $x=x_0$ 
\begin{align}
f(x_0) = f_0, && f'(x_0) = g_0. \lbl{77.12}
\end{align}
Expression \rf{65.12} will give us a solution to the boundary value problem \rf{64.12},\rf{77.12} if the unknown boundary data $f(x_1), \; f'(x_1)$  is made to vanish. This we achieve by posing the following conditions on the Green's function
\begin{align*}
k(x_1; \xi) = 0,\nonumber \\ 
k'(x_1; \xi) = 0. 
\end{align*}
Using the general description of $k(x; \xi)$ from \rf{61} we now get 
\begin{align*}
a(\xi) \; x_1 + b(\xi) &= 0,\nonumber \\
a(\xi) &= 0,
\end{align*}
which leads us to the Green's function 
\begin{align}
k(x; \xi) =\begin{cases}
\; \; 0 \qquad \qquad \qquad \qquad x>\xi \\
x - \xi  \qquad \qquad \qquad\;\;\;  x<\xi 
\end{cases}.\lbl{80.12}
\end{align}
The expression \rf{65.12}, with Green's function given by \rf{80.12}, gives us the solution 
\begin{align}
f(\xi) &= \int_{x_0}^{x_1} dx \; k(x; \xi) \; h(x) - k(x_0;\xi) \; g_0\nonumber \\
&+ k'(x_0; \xi) \; f_0,  \lbl{81.12} 
\end{align}
and from \rf{80.12} we have
\begin{align*}
k'(x; \xi) &= \begin{cases}
\; \; 0 \qquad \qquad \qquad x>\xi \\
\; \; 1 \qquad \qquad \qquad x<\xi \\
\end{cases}.
\end{align*}
Therefore, for this case we have 
\begin{align*}
k(x_0;\xi) &= x_0 - \xi,\nonumber \\
k'(x_0; \xi) &= 1, 
\end{align*}
and the identity \rf{81.12} gives us the solution
\begin{align}
f(\xi) = \intxi dx (x- \xi) \; h(x) + (\xi - x_0) \; g_0 + f_0. \lbl{85.12} 
\end{align}
If $x(t)$ is position as a function of time, and $F(t)$ is the force acting on a particle,  we have according to Newton that 
\begin{align*}
m \; x''(t) = F(t). 
\end{align*}
This is equation \rf{64.12}, with $h = -\frac{F}{m} $. For this case the choice of Green's function \rf{80.12} can be written as 
\begin{align}
k(t';t) = \begin{cases}
\; \; 0 \qquad \qquad \qquad \qquad t' > t \\
\; \; t' - t \qquad \qquad \qquad \;t' < t 
\end{cases}, \lbl{87.12}
\end{align}
and the solution \rf{85.12} is 
\begin{align}
x(t) = \frac{1}{m} \; \int^{t}_{t_0} dt'(t - t') \; F(t') + (t-t_0) \; x'(t_0) + x(t_0). \lbl{88.12}
\end{align}
We observe that \rf{88.12} respects causality; the current position depends on the current and past values of the force. Thus \rf{87.12} is a retarded Green's function. \\
\indent It is evident from the description of all possible Green's function, \rf{61.12}, of $L = - \dsx{}$, that picking one that satisfies any given choice of boundary conditions is simple.
However in order to illustrate what one can do to find Green's functions that satisfy boundary conditions in more realistic and complex cases, we will construct the Green's function \rf{74.12} using the finite Fourier transform. This method applies in more complex cases also, where deriving a formula describing all possible function \rf{61.12} is impossible. \\ 
In order to simplify our exposition we choose $x_0 = 0, \; x_1 = l$. \\
The problem we seek to solve is the following one 
\begin{align}
- k''(x; \xi) &= \delta(x-\xi), \nonumber\\
k(0; \xi) &= k(l;\xi) = 0.  \lbl{89.12}
\end{align}
This is the unique Green's function that satisfy Dirichlet conditions at the boundary.
We want to use the finite Fourier transform and focus therefore on the eigenvalue problem 
\begin{align*}
-M''(x) = \lambda^2 \; M(x), \nonumber\\
M(0) = M(l) = 0.
\end{align*}
We solved this problem several times last semester. The solution is
\begin{align*}
M_k(x) &= \sqrt{\frac{2}{l}} \; \sin(\lambda_k x) && k=1,2...\;,\nonumber \\
\lambda_k &= \frac{\pi k}{l}.
\end{align*}
We then look for a solution to \rf{89.12} of the form 
\begin{align}
k(x,\xi)= \sqrt{\frac{2}{l}} \; \mathlarger{\sum}_{k=1}^{\infty} \; N_k(\xi) \; \sin(\lambda_k x), \lbl{92.12} 
\end{align}
where
\begin{align*}
N_k(\xi) = \sqrt{\frac{2}{l}} \; \int_{0}^l dx \; \sin(\lambda_k x) k(x,\xi).
\end{align*}
From \rf{89.12} we get 
\begin{align*}
- \sqrt{\frac{2}{l}} \; \int^l_0 dx \; \sin(\lambda_k x) \; k''(x; \xi) =  \sqrt{\frac{2}{l}} \; \int^l_0 dx \; \sin(\lambda_k x) \; \delta(x-\xi),
\end{align*}
and using integration by parts and the boundary conditions we get 
\begin{align*}
\lambda_k^2 \; N_k(\xi) &= \sqrt{\frac{2}{l}} \; \sin(\lambda_k \xi), \nonumber \\
&\Updownarrow \nonumber \\  N_k(\xi) &= \sqrt{\frac{2}{l}} \; \frac{\sin(\lambda_k \xi)}{\lambda_k^2}.
\end{align*}
Inserting this into \rf{92.12} we get the following formula for \rf{92.12}
\begin{align}
k(x,\xi) = \frac{2}{l} \; \mathlarger{\sum}_{k=1}^{\infty} \; \frac{\sin(\lambda_k \xi)\sin(\lambda_k x)}{\lambda_k^2}. \lbl{96.12} 
\end{align} 
If we introduce the orthogonal eigenfunctions $M_k(x)$, \rf{96.12} can be more compactly written 
\begin{align}
k(x,\xi) = \mathlarger{\sum}_{k=1}^{\infty} \frac{M_k(\xi) \; M_k(x)}{\lambda_k^2}. \lbl{97.12}
\end{align}
The structure of formula \rf{97.12} is very general, this is the kind of formula that we \ttx{always} get for the Green's function when we apply the finite Fourier transform. The only things that changes from case to case, whether they are 1D, 2D or 3D, are the nature of the eigenfunctions $M_k$ and the corresponding eigenvalues. \\
This approach also yields an interesting and useful formula for the Dirac delta function. Differentiating \rf{97.12} twice term by term and using \rf{89.12} we get 
\begin{align}
\mathlarger{\sum}^{\infty}_{k=1} M_k(\xi) \; M_k(x) = \delta(x-\xi). \lbl{98.12} 
\end{align}
Again \rf{98.12} is fully general and holds for any complete set of eigenfunctions whether they are 1D, 2D or 3D. 

In deriving \rf{96.12} we have certainly been playing with formulas, and given that \rf{97.12} and \rf{71.12} with $x_0 = 0, x_1 = l$ look formally very different it would be instructive to verify that $k(x,\xi)$ from \rf{97.12} is in fact equal to the expression from \rf{71.12}
\begin{align}
k(x,\xi) =\begin{cases}
\; \frac{\xi(l - x)}{l} \qquad \qquad \qquad \qquad \qquad x>\xi \\
\; \frac{x(l - \xi)}{l} \qquad \qquad \qquad \qquad \qquad x<\xi
\end{cases}.\lbl{99.12}
\end{align}
In order to do this we must calculate the Fourier coefficients of $k(x,\xi)$ in \rf{99.12} with respect to the orthogonal system 
\begin{align*}
M_k(x) &= \sqrt{\frac{2}{l}} \sin(\frac{\pi\; k}{l}x) \qquad \qquad \qquad \qquad \qquad k=1,2,3...\;, 
\end{align*}
We have
\begin{align}
N_k(\xi) &= \int^l_0 dx \; M_k(x) \; k(x,\xi) \nonumber \\
&= \sqrt{\frac{2}{l}} \;\int^l_0 dx  \; \sin(\frac{\pi \;k}{l}x) k(x,\xi) \nonumber \\
&= \sqrt{\frac{2}{l}} \; \frac{(l - \xi)}{l} \; \int_0^{\xi} dx \; \sin(\frac{\pi \;k}{l}x) \; x \nonumber \\ 
&+ \sqrt{\frac{2}{l}} \; \frac{\xi}{l} \; \int^l_{\xi} dx \; \sin(\frac{\pi \; k}{l}x) \; (l-x), \nonumber 
\end{align}
and 
\begin{align}
\int_0^{\xi} dx \; \sin(\frac{\pi \; k}{l}x) \; x &= - \frac{l}{\pi k} \; \cos(\frac{\pi \; k}{l} x) \; x \mathlarger{|}^{\xi}_0 \nonumber \\ 
&+ \frac{l}{\pi\; k} \; \int_{0}^{\xi} dx \; \cos(\frac{\pi \; k}{l} x) \nonumber \\
&= - \frac{\xi \; l}{\pi \; k} \; \cos(\frac{\pi \;  k}{l} \xi) + (\frac{l}{\pi \; k})^2 \; \sin(\frac{\pi \; k}{l} \; \xi), \nonumber \\ 
\int_{\xi}^{l} dx \; \sin(\frac{\pi \; k}{l}x) \; (l-x) &= - \frac{l}{\pi \; k} \; \cos(\frac{\pi \; k}{l} x) \; (l-x) \mathlarger{|}^l_{\xi} \nonumber \\ 
&- \frac{l}{\pi \; k} \; \int_{\xi}^l dx \; \cos(\frac{\pi \; k}{l} x) \nonumber \\ 
&= \frac{l(l-\xi)}{\pi \; k} \; \cos(\frac{\pi \; k}{l} \xi ) + (\frac{l}{\pi \; k})^2 \; \sin(\frac{\pi \; k}{l} \xi).  \nonumber 
\end{align}
Thus
\begin{align}
N_k(\xi) &= \sqrt{\frac{2}{l}} \; \frac{(l - \xi)}{l} \; \{ - \frac{l \; \xi}{\pi \; k} \; \cos(\frac{\pi \; k}{l} \xi ) + (\frac{l}{\pi \; k})^2 \; \sin(\frac{\pi \; k}{l} \xi) \} \nonumber \\
&+ \sqrt{\frac{2}{l}} \; \frac{\xi}{l} \; \{ (\frac{l(l-\xi)}{\pi \; k}) \; \cos(\frac{\pi \; k}{l} \xi )  + (\frac{l}{\pi \; k})^2 \; \sin(\frac{\pi \; k}{l} \xi) \} \nonumber \\
&= \sqrt{\frac{2}{l}} \; (\frac{l}{\pi \; k})^2 \; \sin(\frac{\pi \; k}{l} \xi) = \frac{M_k (\xi)}{\lambda_k^2}. \nonumber 
\end{align}
Let us next turn to the second way in which the integral identity \rf{65.12} can be used to find solutions to the differential equation \rf{64.12}. This is the way leading to boundary integral equations. 
\\
Let us consider the boundary value problem 
\begin{align}
-f''(x) &= h(x) && x_0 < x <x_1, \nonumber \\
f(x_0) &= f_0, \nonumber \\ 
f(x_1) &= f_1. \lbl{101.12}
\end{align}
Inserting the known boundary values into the general integral identity \rf{65.12} gives us 
\begin{align}
f(\xi) &= \intx dx \; k(x;\xi) \; h(x)\nonumber \\
&+ k(x_1; \xi) \; f'(x_1) - k(x_0; \xi) \; f'(x_0) \nonumber \\
&-k'(x_1; \xi) \; f_1 + k'(x_0;\xi)\; f_0,  \lbl{102.12} 
\end{align}
which holds for all $x_0 < \xi < x_1$. Let now $\xi$ approach $x_0$ from above and $x_1$ from below. We get 
\begin{align}
f_0 &= \intx dx \; k(x;x_0) \; h(x)\nonumber\\ 
&+ k(x_1; x_0) \; f'(x_1) - k(x_0; x_0) \; f'(x_0) \nonumber \\ 
&- k'(x_1; x_0) \; f_1 + k'(x_0; x_0) \; f_0, \nonumber \\
f_1 &= \intx dx \; k(x; x_1) \; h(x) \nonumber\\ 
&+ k(x_1; x_1) \; f'(x_1) - k(x_0; x_1) \; f'(x_0) \nonumber \\ 
&- k'(x_1; x_1) \; f_1 + k'(x_0; x_1) \; f_0,  \lbl{103.12}  
\end{align}
Where we must be careful to use limits when we evaluate $k'(x; x)$ in $x_0$ and $x_1$ since $k'(x;\xi)$ is discontinuous at $x=\xi$. The correct way to evaluate them are 
\begin{align*}
 k'(x_0; x_0)=\lim\limits_{\epsilon\rightarrow0}k'(x_0,x_0+\epsilon),\\
 k'(x_1; x_1)=\lim\limits_{\epsilon\rightarrow0}k'(x_1,x_1-\epsilon).
 \end{align*}
Observe that \rf{103.12} is a system of two equations for the two unknown boundary data $f'(x_0)$ and $f'(x_1)$. The system \rf{103.12} is the boundary integral equation for this situation. We can write the system as 
\begin{align}
\mqty(-k(x_0;x_0) && k(x_1;x_0) \\ -k(x_0; x_1) && k(x_1; x_1)) \; \mqty(f'(x_0) \\ f'(x_1)) = \mqty(b_0 \\ b_1), \lbl{105.12} 
\end{align}
where 
\begin{align*}
b_0 &= f_0 - \intx dx \; k(x; x_0) \; h(x) + k'(x_1;x_0) \; f_1, \nonumber \\
&- k'(x_0; x_0) \; f_0 \nonumber \\
b_1 &= f_1 - \intx dx \; k(x;x_1) \; h(x) + k'(x_1;x_1) \; f_1 \nonumber\\
&- k'(x_0;x_1) \; f_0.
\end{align*}
The only requirement on the Green's function is that the determinant of the matrix in \rf{105.12} is non-zero. This is not a very strict requirement on $k(x,\xi)$, most Green's functions will satisfy it. If you by some means have gotten hold of a Green's function, using that Green's function in the boundary integral equation will almost certainly be ok. 
\\
For the simple operator we are discussing here, all Green's functions are known, and the subset of Green's functions leading to a singularity of the boundary integral equations can be described precisely. From \rf{61.12} we get 
\begin{align*}
&k(x_1;x_0) = a(x_0) \; x_1 + b(x_0), \nonumber \\
&k(x_0;x_0) = a(x_0) \; x_0 + b(x_0), \nonumber \\
&k(x_1;x_1) = a(x_1) \; x_1 + b(x_1), \nonumber \\
&k(x_0;x_1) = a(x_1) \; x_0 + b(x_1)+x_0-x_1,
\end{align*}
and thus the condition for a singularity is 
\begin{align}
-k(x_0;x_0) \; k(x_1;x_1) + k(x_1; x_0) \; k(x_0; x_1)& = 0, \nonumber \\
&\Updownarrow \nonumber \\ -(a(x_0) \; x_0 + b(x_0)) \; (a(x_1)\; x_1 + b(x_1)) \nonumber \\ 
+ (a(x_0) \; x_1 + b(x_0)) \; (a(x_1)\; x_0 + b(x_1) + x_0 - x_1)&= 0, \nonumber \\ 
&\Updownarrow \nonumber \\ (a(x_0) \; b(x_1) - a(x_1) \; b(x_0) - a(x_0) \; x_1 - b(x_0))(x_1 - x_0) &= 0. \lbl{108.12} 
\end{align}
Choosing to use the Green's function 
\begin{align*}
k(x; \xi) =\begin{cases}
\; \; 0 \qquad \qquad \qquad \qquad \qquad x>\xi \\
x - \xi \qquad \qquad \qquad \qquad \; \; \; x<\xi 
\end{cases},
\end{align*}
corresponding to $a(\xi) = b(\xi) = 0$, clearly will not work, because then \rf{108.12} is satisfied. However, 
\begin{align}
k(x;\xi) =\begin{cases}
\; -x \qquad \qquad \qquad \qquad \qquad x>\xi \\
\; - \xi \qquad \qquad \qquad \qquad \qquad x<\xi 
\end{cases}, \lbl{110.12}
\end{align}
corresponding to $a(\xi) = -1, \; b(\xi) = 0 $, will work nicely because then \rf{108.12} is not satisfied.

  For this particular Green's function we have 
\begin{align*}
k(x_1;x_0) = -x_1,\nonumber\\
k(x_0;x_0) = -x_0,\nonumber \\
k(x_1;x_1) = -x_1, \nonumber \\
k(x_0;x_1) = -x_1, 
\end{align*}
and
\begin{align*}
k'(x; \xi) =\begin{cases}
\; -1 \qquad \qquad \qquad \qquad \qquad x>\xi \\
\; \; 0 \qquad \qquad \qquad \qquad \qquad x<\xi
\end{cases}.
\end{align*}
Thus
\begin{align*}
k'(x_1;x_0) &= -1,\nonumber\\
k'(x_0;x_0) & =\lim\limits_{\epsilon\rightarrow0}k'(x_0,x_0+\epsilon)= 0, \nonumber \\
k'(x_1;x_1)& = \lim\limits_{\epsilon\rightarrow0}k'(x_1,x_1-\epsilon)= -1, \nonumber \\
k'(x_0;x_1) &=  0.
\end{align*}
So our linear system is 
\begin{align*}
\mqty(x_0 && -x_1 \\ x_1 && -x_1 ) \; \mqty(f'(x_0) \\ f'(x_1)) = \mqty(b_0 \\ b_1),
\end{align*}
where now 
\begin{align*}
&b_0 = f_0 + \intx dx \; x \; h(x) - f_1, \nonumber \\
&b_1 = f_1 + x_1 \; \intx dx \; h(x) - f_1 = x_1 \; \intx dx \; h(x).
\end{align*}
The solution of the boundary integral equation for this case is then 
\begin{align*}
\mqty(f'(x_0) \\ f'(x_1)) =\frac{1}{x_1 \; (x_1 - x_0)} \; \mqty(-x_1 && x_1 \\ -x_1 && x_0) \; \mqty({f_0 - f_1 + \intx dx \; x \; h(x) \\ x_1 \; \intx dx \; h(x)}),
\end{align*}
or equivalently
\begin{align*}
f'(x_0) = \frac{1}{x_1 - x_0} \; \{ f_1 - f_0 + \intx dx \; (x_1 - x) \; h(x) \},\nonumber \\
f'(x_1) = \frac{1}{x_1 - x_0} \; \{ f_1 - f_0 + \intx dx \; (x_0 - x) \; h(x) \}. 
\end{align*}
We now insert these expressions for $f'(x_0)$ and $f'(x_1)$ together with $k(x; \xi)$ from \rf{110.12} into \rf{102.12}. This will give us the solution to the boundary value problem. \\
We have 
\begin{align*}
f(\xi) &= -\intxi dx \; \xi \; h(x) - \intix dx \; x \; h(x)\nonumber\\ 
&- x_1 \; f'(x_1) + \xi \; f'(x_0) + f_1 \nonumber \\ 
&= - \xi \; \intxi dx \; h(x) - \intix dx \; x\; h(x) \nonumber \\ 
&- x_1 \; \frac{1}{x_1 - x_0} \; \{f_1 - f_0 + \intx dx \; (x_0 - x) \; h(x)\} \nonumber \\ 
&+ \xi \; \frac{1}{x_1 - x_0} \; \{ f_1 - f_0 + \intx dx \; (x_1 - x) \; h(x) \} + f_1 \nonumber \\
&= \{1 - \frac{x_1}{x_1 - x_0} + \frac{\xi}{x_1 - x_0} \}\; f_1 \nonumber \\
&+ \{ \frac{x_1}{x_1 - x_0} - \frac{\xi}{x_1-x_0} \} \; f_0 - \xi \; \intxi dx \; h(x) \nonumber \\ 
&- \intix dx \; x \; h(x) - \frac{x_1}{x_1 - x_0} \; \intx dx \; (x_0 - x) \; h(x) \nonumber \\
&+ \frac{\xi}{x_1 - x_0} \; \intx dx \; (x_1 - x) \; h(x),   
\end{align*}
and thus we have 
\begin{align}
f(\xi) &= \frac{\xi - x_0}{x_1 - x_0} \; f_1 - \frac{\xi - x_1}{x_1 - x_0} \; f_0 - \xi \; \intxi dx \; h(x) - \intix dx \; x \; h(x)\nonumber\\
&- \frac{x_1}{x_1 - x_0} \; \intx dx \; (x_0 -x) \; h(x) + \frac{\xi}{x_1 - x_0} \; \intx dx \; (x_1 - x) \; h(x). \lbl{119.12}
\end{align}
This solution certainly looks very different from the solution \rf{75.12} that we found previously. However the solution to the boundary value problem \rf{101.12} is unique so \rf{75.12} and \rf{119.12} must really be the same. By rearranging the integrals in \rf{119.12} this can be proved (do it!). \\
A general (ordinary) differential operator of order 2 is of the form 
\begin{align*}
L = a(x) \; \dsx{} + b(x) \; \dx{} + c(x).
\end{align*}
A Green's function for $L$ is a function $k(x; \xi)$ such that 
\begin{align*}
L \; k(x; \xi) = \delta(x-\xi).
\end{align*}
Arguing heuristically like before we find that $k(x;\xi)$ is a Green's function for $L$ only if 
\begin{align*}
a(x) \; k''(x;\xi) + b(x) \; k'(x;\xi) + c(x) \; k(x, \xi) = 0  && \text{for}\; x \ne \xi, \nonumber \\ 
k^{+}(\xi; \xi) - k^{-}(\xi;\xi) = 0, \nonumber\\ 
k'^{+}(\xi;\xi) - k'^{-}(\xi;\xi) = \frac{1}{a(\xi)}. 
\end{align*}
The generalization to ordinary differential operators of order $n$ is straight forward and is left to the reader.
\\
We will now leave the theory of Green's functions for a while and spend some time discussing \ttx{distribution}.

\subsection{The theory of distributions}

Another, and perhaps even a better name, for the things we study in the theory of distributions are \ttx{generalized functions}. This name signifies that what we do is to extend and enlarge the set of functions.
\\
So what is the defining property of a function? It is simply this; a function is a rule that to each \ttx{number} in a set of numbers associate another number, possibly belonging to some different set of numbers. 
\begin{figure}[htbp]
\centering
\includegraphics{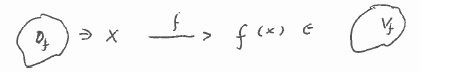}
\caption{}
\label{fig3}
\end{figure}
$f$ is the name of the function and $f(x)$ the \ttx{value} of $f$ for a given number $x$. $D_f$ is the \ttx{domain} of $f$ and $V_f$ the \ttx{range} or \ttx{codomain} of $f$. 
\\
{\it Generalized functions are functions whose domain consist of functions}. 

The first task is then to say what this domain of functions is. 
For any smooth function, $\phi\in C^{\infty}(\mathbf{R})$, on $\mathbf{R}$ we define the \ttx{support} of $\phi$ by 
\begin{align*}
\text{supp}(\phi) = \overline{\{ x \in \mathbf{R} \;\mathlarger{|} \; \phi(x) \ne 0 \}},
\end{align*}
\noindent where for any subset $A \subset \mathbf{R}, \; \bar{A}$ is the \ttx{closure} of $A$. Recall that taking the closure of $A$ consists of adding all \ttx{boundary points} of $A$. \\
Recall also that a closed and bounded subset of $\mathbf{R}$ is called \ttx{compact}. We are now ready to define the domain for our generalized functions.
\begin{align*}
D_0 = \{ \; \phi \in C^{\infty}(\mathbf{R}) \; \mathlarger{|} \; \text{supp}(\phi) \; \text{is compact} \}.
\end{align*}
It is not hard to imagine continuous functions of compact support. Here is one 
\begin{align*}
\phi(x) =\begin{cases}
\; \; 0 \qquad \qquad \qquad \qquad \qquad x>1 \\
\; 1-x \qquad \qquad \qquad \qquad \frac{1}{2} \leq x \leq 1 \\
\; \; x \qquad \qquad \qquad \qquad \; \; \; \; \; 0 \leq x \leq \frac{1}{2} \\
\; \; 0 \qquad \qquad \qquad \qquad \qquad x<0 
\end{cases}. 
\end{align*}

\noindent However $\phi(x)$ is not smooth. It is not differentiable at the points  $x=0, \; \frac{1}{2}, \; 1$.  In order to create an element of $D_0$, the function has, for example, to become zero at $x=0$ and $x=1$  in an infinitely smooth way. Can this even be done? Could it be that $D_0$ is in fact empty? 

Fortunately $D_0$ is not empty! Recall from Calculus, the following function 
\begin{align}
\phi(x) =\begin{cases}
\exp(- \frac{1}{x^2}) \qquad \qquad \; \;  x \geq 0 \\
\; \; 0 \qquad \qquad \qquad \qquad  x \leq 0 
\end{cases}. \lbl{128.12}
\end{align} 
Using elementary calculus, one can prove that $\phi(x)$ is infinitely differentiable at $x=0$ (do it!). By joining together functions like \rf{128.12}, we can create smooth functions of compact support. Here is an example of a smooth function whose compact support is $[-1,1]$
\begin{align*}
\phi(x) =\begin{cases}
\exp(- \frac{1}{1-x^2}) \qquad \qquad \; \;  |x|<1 \\
\; \; 0 \qquad \qquad \qquad \qquad \;\;\;\; |x|\geq1
\end{cases}.
\end{align*} 
 Thus $D_0$ is not empty, as a matter of fact,  it is very large. In a way that can be made precise it is in fact infinitely larger than any subset of real numbers. 

We now have a domain for our generalized functions. A function on $D_0$ is something that evaluates to a number for any $\phi \in D_0$. 
\begin{align*}
f(\phi) \in \mathbf{R} && \phi \in D_0.
\end{align*}
There are many such functions. Here are some examples. 
\begin{example}\label{DiracDelta}
Let $x_0 \in \mathbf{R}$. Define a function on $D_0$ by 
\begin{align*}
\delta(x-x_0) (\phi) = \phi(x_0) \in \mathbf{R} && \phi \in D_0.
\end{align*}
We will see that this function is nothing but the Dirac delta function (hence the notation).
\end{example}
\begin{example}\label{Tf}
Let $f: \mathbf{R} \rightarrow \mathbf{R}$ be a continuous function. Define a function $T_f$ on $D$ by 
\begin{align}
T_f(\phi) = \int^{\infty}_{-\infty} dx \; f(x) \; \phi(x) \; \; \in \mathbf{R}. \lbl{131.12}
\end{align}
Observe that $T_f$ is well defined because $\phi \in D_0$ has compact support so that the integral converges. In fact, $f$ does not have to be continuous in order for $T_f$ to be well defined. The very extensive class of \ttx{locally integrable} functions define functions on $D_0$ through \rf{131.12}.
\end{example}
\begin{example}\label{Heaviside}
Define a function on $D_0$ through 
\begin{align*}
H(\phi) = \int^{\infty}_0 dx \; \phi(x).
\end{align*}
$H$ is a very important function in the theory of distributions.
\end{example}
\noindent Observe that $D_0$ is a vector space over $\mathbf{R}$. Vector space operations are defined in the usual way 
\begin{align*}
(\phi_1 + \phi_2) (x) &= \phi_1(x) + \phi_2(x),\nonumber \\
(a \; \phi)(x) &= a(\phi(x)).
\end{align*}
\noindent In the theory of distributions we only consider functions on $D_0$ that are \ttx{linear} with respect to the vector space structure on $D_0$
\begin{align*}
f(\phi_1 + \phi_2) &= f(\phi_1) + f(\phi_2),\nonumber\\
f(a\; \phi) &= a \; f(\phi).
\end{align*}
\noindent The functions on $D_0$ defined in examples \ref{DiracDelta}, \ref{Tf} and \ref{Heaviside} are all linear. For example, for the one in \ref{DiracDelta} we have 
\begin{align*}
\delta(x-x_0) (\phi_1 + \phi_2) &= (\phi_1 + \phi_2) (x_0) = \phi_1 (x_0) + \phi_2(x_0)\nonumber\\
&= \delta(x-x_0) (\phi_1) + \delta(x-x_0) (\phi_2), \nonumber \\
\delta(x-x_0)(a \; \phi) &= (a \; \phi) (x_0) = a \; (\phi(x_0)) = a \; \delta(x-x_0) (\phi).
\end{align*}
Note that whereas the set of linear functions on $\mathbf{R}$ is very small, they must be of the form $f(x) = a \; x$ where $a$ is some real number, the set of linear functions on $D_0$ is very large. In fact without further restriction it is so large and varied that no general theory can be created for \ttx{all} linear functions on $D_0$. In order to get a useful theory we must restrict to a subclass of all linear functions on $D_0$. Like in calculus we do this by requiring that the functions on $D_0$ are continuous. Here,  $f$ defined on $D_0$ is continuous at $\phi_0 \in D$, if for \ttx{all} sequences $\{ \phi_n \}$ in $D_0$ with 
\begin{align}
\phi_n \rightarrow \phi_0, \lbl{136.12} 
\end{align}
we have 
\begin{align*}
f(\phi_n) \rightarrow f(\phi_0).
\end{align*}
Formally this definition of continuity is the same as the regular one from calculus. Of course, we have not actually defined continuity yet, since we have not given a meaning to the limit \rf{136.12}. For now however, let us assume that the limit \rf{136.12} has been given a precise meaning. 

Given this, let $f$ be a function on $D_0$ that is continuous at $\phi = 0$. Let $\phi_0 \in D_0$ be any element in $D_0$ and let $\{ \phi_n \}$ be a sequence in $D_0$ that converges to $\phi_0$ 
\begin{align*}
\phi_n \rightarrow \phi_0.
\end{align*}
Let $\{ \psi_n \}$ be the sequence 
\begin{align*}
\psi_n = \phi_n - \phi_0, 
\end{align*}
then 
\begin{align*}
\psi_n \rightarrow \phi_0 - \phi_0 = 0 && \text{when} \; n \rightarrow \infty,
\end{align*}
and since $f$ is continuous at $\phi_0 = 0$, we have 
\begin{align*}
f(\psi_n) \rightarrow f(0) = 0,
\end{align*}
using the linearity of $f$. From the linearity of $f$ we also get 
\begin{align*}
f(\phi_n) = f(\phi_0 + \phi_n - \phi_0) = f(\phi_0) + f(\psi_n) \rightarrow f(\phi_0), 
\end{align*}
and thus $f$ is continuous at $\phi_0 \in D_0$. The conclusion is that for a linear function on $D_0$, we only need to check continuity at $\phi_0 = 0$. We thus only need to specify precisely what it means for a sequence in $D_0$ to converge to zero. Such sequences are called \ttx{zero sequences}.
\begin{definition}
Let $\{ \phi_n \}$ be an infinite sequence in $D_0$. Then $\{ \phi_n\}$ is a zero sequence if and only if \\ \\
i)  There exists a bounded interval $I \subset \mathbf{R}$ such that 
\begin{align*}
\text{supp}(\phi_n) \subset I && \forall n,
\end{align*}
\\
ii) 
\begin{align*}
\lim\limits_{n\rightarrow\infty} \; \max\limits_{x \in \mathbf{R}} \; \abs{\frac{d^k \phi^n}{dx^k}} = 0 && \forall k \geq 0,
\end{align*}
or in other words the sequences $\{ \frac{d^k \phi^n}{dx^k} \}$ converge \ttx{uniformly} to zero on $\mathbf{R}$ for all $k \geq 0$.
\end{definition}
\noindent With these formulations out of the way we have
\begin{definition}\label{DistributionDef}
A generalized function, or distribution, is a continuous linear function on $D_0$. 
\end{definition}
\noindent Note that functions $f:D_0\rightarrow \mathbf{R}$ are often called \ttx{functionals} to distinguish them from regular calculus functions.

  We do not put a great emphasis on mathematical stringency in these lecture notes, and will usually assume that reasonably constructed linear functions on $D_0$ are in fact continuous, and thus define generalized functions. However, in order to solidify the definitions let us show that some of the previously defined linear functions on $D_0$ are in fact generalized functions according to definition 2. 
\begin{example}
In example \ref{example1} we defined the linear function $\delta(x-x_0)$ by 
\begin{align*}
\delta(x-x_0) (\phi) = \phi(x_0).
\end{align*} 
Let $\{ \phi_n \}$ be a zero sequence. Then according to the definition (\ref{DistributionDef})  we have 
\begin{align*}
\lim\limits_{n \rightarrow \infty} \; \max\limits_{x \in \mathbf{R}} \; \abs{\phi_n(x)} = 0,
\end{align*}
and thus
\begin{align*}
\abs{\delta(x-x_0) (\phi_n)} = \abs{\phi_n (x_0)} \leq \max\limits_{x \in \mathbf{R}} \; \abs{\phi_n(x)} \; \; \rightarrow 0 \qquad n \rightarrow \infty.
\end{align*}
Therefore, $\delta(x-x_0)$ is a generalized function.
\end{example}
\begin{example}\label{HeavisideProof}
In example \ref{example3} we defined the function $H$ on $D_0$ by 
\begin{align*}
H(\phi) = \int^{\infty}_0 dx \; \phi(x)
\end{align*}
$H$ is linear because 
\begin{align*}
H(c_1\phi_1 +c_2 \phi_2) &= \int_{0}^{\infty} dx \; (c_1\phi_1 + c_2\phi_2) (x)\\
& = \int^{\infty}_{0} dx (c_1\phi_1(x) + c_2\phi_2(x))  \\
&= c_1\int^{\infty}_{0} dx \; \phi_1(x) + c_2\int_0^{\infty} dx \; \phi_2(x)r\\
& =c_1 H(\phi_1) +c_2 H(\phi_2). 
\end{align*}
Let $\{ \phi_n \}$ be a zero sequence. Then there exists a finite interval $I$ such that 
\begin{align*}
\text{supp} (\phi_n) \subset I && \forall n,
\end{align*}
and
\begin{align*}
\max\limits_{x \in I} \; \abs{\phi_n (x)} = \max\limits_{x \in \mathbf{R}} \; \abs{\phi_n(x)} \rightarrow 0 \qquad \qquad \qquad n\rightarrow \infty.
\end{align*}
Therefore 
\begin{align*}
\abs{H (\phi_n)} &= \abs{\int^{\infty}_0 dx \; \phi_n(x)} \leq \int^{\infty}_0 dx \; \abs{\phi_n(x)}\nonumber\\ 
&\leq \int_I dx \; \abs{\phi_n(x)} \leq \abs{I} \; \max\limits_{x \in I} \; \abs{\phi_n(x)} \nonumber \\
&\rightarrow 0 && n \rightarrow \infty \nonumber \\
&\text{(Here $\abs{I}$ is the length of $I$)} 
\end{align*}
\end{example}
\noindent The proof that all functions $T_f$ from example \ref{Tf} are generalized functions is very similar to example \ref{HeavisideProof}.

The generalized functions of the form $T_f$ shows that any locally integrable function on $\mathbf{R}$ defines a corresponding generalized function. We can thus consider any regular calculus function to also be a generalized function. We evaluate the corresponding generalized function by integrating, like in example \ref{Tf}. \\
Since not all generalized functions are of the form $T_f$ for some $f$, $\delta(x-x_0)$ being the primary example, the set of generalized functions is a true extension of the notion of function as we know it from calculus. A generalized functions that is of the form $T_f$ for some locally integrable function $f$, is called \textit{regular}. All other generalized functions are called \textit{singular}.
\begin{figure}[htbp]
\centering
\includegraphics{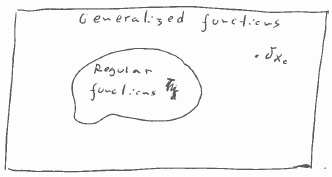}
\caption{The domain for generalized functions}
\label{fig4}
\end{figure}

\noindent The relation between locally integrable functions $f(x)$, and the corresponding generalized functions $T_f$ is not one to one. Clearly if $f(x), \; g(x)$ are equal except at a finite number of points we have 
\begin{align*}
T_f = T_g, 
\end{align*}
so they correspond to the \ttx{same} generalized function. Nevertheless, the association is so close that we should think of regular functions \ttx{as} generalized functions.\\
In fact we will encourage this identification by abandoning the notation $T_f$, just writing 
\begin{align*}
f(x)(\phi) = \int^{\infty}_{-\infty} dx \; f(x) \; \phi(x).
\end{align*}
In this way $\abs{x}$ is the generalized function defined by 
\begin{align*}
\abs{x} (\phi) &= \int^{\infty}_{-\infty} dx \; \abs{x} \; \phi(x) = - \int^{0}_{-\infty} dx \; x \; \phi(x)\nonumber\\
&+ \int^{\infty}_{0} dx \; x \; \phi(x),
\end{align*}
and the generalized function $H$ from example \ref{Heaviside} clearly corresponds to the locally integrable function 
\begin{align*}
H(x) =\begin{cases}
\; \; 1 \qquad \qquad \qquad \qquad \qquad  x \geq 0 \\ 
\; \; 0 \qquad \qquad \qquad \qquad \qquad  x < 0 
\end{cases},
\end{align*}
where the action of $H(x)$ on $D_0$ is defined by the formula
\begin{equation*}
H(x)(\phi)=\int_{-\infty}^{\infty}\;dx\;H(x)\phi(x)
\end{equation*}
The generalized function $H(x)$ is called the \textit{Heaviside} function. Note that I now think of the locally integrable functions $H(x)$ and $\abs{x}$ as generalized function. Here I introduce an abuse of notation that is common in calculus and which we bring to a new level for generalized functions.

  In calculus we often write $f(x)$ even if we really mean the function $f$, not the function value $f(x)$. This abuse of notation is very useful in calculus. Using this notation we say that the function is $f$ \ttx{and} that generic points in the domain will be denoted by $x$. We can then in a compact way introduce other functions using the notation $f(x)$ as a starting point 
\begin{align*}
g(x) = f(2x),\nonumber \\
h(x) = f(x-x_0).
\end{align*}
In a similar way we will use this abuse of notation to introduce new generalized functions from old ones.

   Let $\{f_n(x) \}$ be a sequence of generalized functions and let $f(x)$ be a generalized function. Then \label{WeakConvergence}
\begin{align*}
f_n(x) \rightarrow f(x) && n\rightarrow\infty,\nonumber \\ 
\Updownarrow \nonumber \\ f_n(x) (\phi) \rightarrow f(x) (\phi) && \forall \phi \in D , \; \; n\rightarrow\infty.
\end{align*}
Thus $f_n(x) \rightarrow f(x)$ means pointwise convergence as functions on $D_0$. This simple and natural definition is called \ttx{weak convergence}. There are many other notions of convergence for sequences of generalized functions, but we will not discuss them in these lecture notes.
\begin{example}
Let 
\begin{align*}
f_n(x) =\begin{cases}
\frac{1}{2} \; n \qquad \qquad \qquad \qquad \qquad \abs{x} < \frac{1}{n} \\
\; \; 0 \qquad \qquad \qquad \qquad \qquad \;  \; \abs{x} > \frac{1}{n} 
\end{cases},
\end{align*} 
then for any $\phi \in D_0$ we have 
\begin{align*}
f_n(x) (\phi) &= \frac{1}{2} \; n \; \int^{\frac{1}{n}}_{-\frac{1}{n}} dx \; \phi(x) = \frac{n \; \phi(\hat{x})}{2} \; \int^{\frac{1}{n}}_{-\frac{1}{n}} dx\nonumber \\ 
&= \phi(\hat{x}),
\end{align*}
where we have used the mean value theorem and $ - \frac{1}{n} < \hat{x} < \frac{1}{n} $. Observe that as $n\rightarrow \infty, \; \hat{x} \rightarrow 0$. The continuity of $\phi$ then gives us 
\begin{align*}
f_n(x)(\phi) = \phi(\hat{x}) \rightarrow \phi(0) && n\rightarrow \infty.
\end{align*}
But $\phi(0) = \delta(x) (\phi) \; \; \forall \phi \in D$. We have therefore proved that
\begin{align*}
f_n(x) \rightarrow \delta(x) \; \; \; \text{as} \; n \rightarrow \infty,
\end{align*}
and thus $\{ f_n(x) \}$ is a sequence of regular generalized functions that converge to the singular generalized function $\delta(x)$. Sequences that converge to $\delta(x)$ are important in the theory of distributions and are called \ttx{delta sequences}. \\
If $f_n(x)$ are like rational numbers, then $\delta(x)$ is like an irrational number and a delta sequence is an approximation of an irrational number in terms of rational numbers. This is \ttx{not} a superficial analogy, on the contrary, the analogy runs very deep.\\
Let $f(x)$ be regular. Then by definition $f(ax)$ is the generalized function 
\begin{align}
f(ax)(\phi) &= \int^{\infty}_{-\infty}dx \; f(ax) \; \phi(x) = \frac{1}{\abs{a}} \; \int^{\infty}_{-\infty} dy \; f(y) \; \phi_a(y)\nonumber \\ 
&= \frac{1}{\abs{a}} \; f(x) (\phi_a), \lbl{171.12} 
\end{align}
where $\phi_a(x) = \phi(\frac{x}{a}) \in D_0$ when $a \ne 0$. 
For any generalized function, not necessarily regular, we use \rf{171.12} to \textit{define} $f(ax)(\phi)$, 
\begin{align}
f(ax)(\phi) = \frac{1}{\abs{a}} \; f(x) (\phi_a). \lbl{172.12} 
\end{align}
For the particular case when $f(x) = \delta(x)$ we get from \rf{172.12} 
\begin{align*}
\delta(ax) (\phi) &= \frac{1}{\abs{a}} \; \delta(x) (\phi_a) = \frac{1}{\abs{a}} \; \phi_a(0)\nonumber \\
&= \frac{1}{\abs{a}} \; \phi(0)=\frac{1}{\abs{a}} \; \delta(x) (\phi),
\end{align*}
and thus we have the identity
\begin{align*}
\delta(ax) = \frac{1}{\abs{a}} \; \delta(x).
\end{align*}
Using $a=-1$ we get the interesting identity 
\begin{align*}
\delta(-x) = \delta(x),
\end{align*}
$\delta(x)$ is by definition an \ttx{even} generalized function. Observe how efficient our abuse of notation is! \\
Let $f(x)$ be a regular generalized function. Then by definition $f(x-a)$ is the generalized function
\begin{align*}
f(x-a) (\phi) &= \int^{\infty}_{-\infty} dx \; f(x-a) \; \phi(x)\nonumber \\
&= \int^{\infty}_{-\infty} dy \; f(y) \; \phi^a(y) = f(x) (\phi^a), 
\end{align*}
where $\phi^a(x) = \phi(x+a)$. For any generalized function $f(x)$ we \textit{define} $f(x-a)$ by the identity 
\begin{align*}
f(x-a)(\phi) = f(x)(\phi^a).
\end{align*}
For the particular case $f(x) = \delta(x)$ we get 
\begin{align}
\delta(x-a) (\phi) = \delta(x) (\phi^a) = \phi^a(0) = \phi(a).\lbl{178.12} 
\end{align}
For the particular case of $\delta(x)$ there is some further abuse of notation that is common. We write 
\begin{align}
\delta(x)(\phi) = \intinf dx \; \delta(x) \; \phi(x). \lbl{179.12} 
\end{align}
The right-hand side of \rf{179.12} is purely formal. Using this notation we have for \rf{178.12} 
\begin{align*}
\delta(x-a) (\phi) = \intinf dx \; \delta(x-a) \; \phi(x) = \phi(a).
\end{align*}
Let now $a(x)$ be a smooth function on $\mathbf{R}$ and let $f(x)$ be a regular generalized function. Then 
\begin{align}
(a(x)f(x)) (\phi) &= \intinf dx\; a(x) f(x) \; \phi(x)\nonumber \\
&= \intinf dx \; f(x)\;  a(x)  \phi(x)\nonumber\\
&= \intinf dx \; f(x)\; (a\phi)(x)\nonumber\\
&= f(x) (a \phi),  \lbl{181.12}
\end{align}
Where, by definition of products of functions, we have
\begin{align*}
(a\phi)(x) = a(x) \; \phi(x). 
\end{align*}
Observe that since $a(x)$ is smooth $a\phi \in D_0$ and \rf{181.12} make sense. We now use \rf{181.12} for any generalized function and \textit{define}
\begin{align}
(a(x) f(x)) (\phi) = f(x) (a \phi). \lbl{183.12} 
\end{align}
Let us again consider the special case $f(x) = \delta(x)$. We get 
\begin{align*}
(a(x)\delta(x)) (\phi) &= \delta(x) (a\phi) =(a \phi)(0)\nonumber \\
&= a(0) \; \phi(0) = a(0) \; \delta(x) (\phi),
\end{align*}
so we have the identity 
\begin{align*}
a(x)\; \delta(x) = a(0) \; \delta(x).
\end{align*}
For the special case $a(x) = x$ we get 
\begin{align*}
x \; \delta(x) = 0 \; \delta(x)=0.
\end{align*}
Let $f(x)$ be a regular generalized function with $f(x)$ differentiable. Thus $f'(x)$ is also a regular generalized function, and we have 
\begin{align}
f'(x)(\phi) &= \intinf dx \; f'(x) \; \phi(x) = f(x) \; \phi(x)\mathlarger{|}^{\infty}_{-\infty}\nonumber \\ 
&- \intinf dx \; f(x) \; \phi'(x) = - \intinf dx \; f(x) \; \phi'(x) \nonumber \\
&= - f(x) (\phi'),  \lbl{187.12}
\end{align}
where we have used the fact that all $\phi \in D_0$ have compact support so that $\phi(\pm \infty) = 0$. We now use \rf{187.12} to \textit{define} the derivative of any generalized function, $f(x)$, by
\begin{align}
f'(x) (\phi) = - f(x) (\phi'). \lbl{188.12}
\end{align}
Since $\phi \in D_0$ are smooth functions, $\phi' \in D_0$ and \rf{188.12} makes sense. Note that \rf{188.12} tells us that \ttx{all} generalized functions can be differentiated. In fact they are infinitely differentiable because \rf{188.12} can obviously be generalized into 
\begin{align*}
f^{(k)}(x) (\phi) = (-1)^k \; f(x) (\phi^{(k)}),
\end{align*}
by using repeated integration by parts in \rf{187.12}.
For the particular case $f(x) = \delta(x)$ \rf{188.12} gives us 
\begin{align*}
\delta'(x) (\phi) = - \delta(x) (\phi') = - \phi'(0).
\end{align*}
Let us consider the locally integrable function $H(x)$. It is certainly not smooth in the usual calculus sense, it is not even continuous. However it is (infinitely) differentiable as a generalized function!
\begin{align*}
H'(x)(\phi) &= -H(x) (\phi') = - \int^{\infty}_0 dx \; \phi'(x)\nonumber\\
&= - \phi(x)\mathlarger{|}^{\infty}_0 = \phi(0) = \delta(x)(\phi),
\end{align*}
and we therefore get the identity
\begin{align*}
H'(x) = \delta(x).
\end{align*}
Let us combine the derivative with multiplication by smooth functions from \rf{183.12}. This gives 
\begin{align*}
x \; \delta'(x) (\phi) &= \delta'(x)(x \phi) = - \delta(x) ((x \phi)')\nonumber\\
&= - \delta(x) (\phi + x \phi') = - \phi(0) - 0 \; \phi'(0) \nonumber \\
&= - \delta(x) (\phi), 
\end{align*}
and we get the identity 
\begin{align*}
x \; \delta'(x) + \delta(x) = 0.
\end{align*}
\end{example}
\noindent Thus the delta functions is a solution to the ODE
\begin{align}
xf'(x)+f(x)=0.\nonumber
\end{align}
Generalized functions, as it turns out, are very well suited for describing singular solutions to both ODEs and PDEs.

\noindent As a final example let us consider a function $f(x)$ that is smooth except for a point $x=a$ where it has a jump discontinuity. 
\begin{figure}[htbp]
\centering
\includegraphics{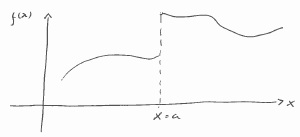}
\caption{A smooth function with a jump discontinuity at $x=a$}
\label{fig5}
\end{figure}

\noindent Let $[f(x)]_{x=a} = \lim\limits_{x \rightarrow a^{+}} f(x) - \lim\limits_{x \rightarrow a^{-}} f(x)$ be the jump in function value of $f$ as we pass through $x=a$.
Like all generalized functions,  $f(x)$ is differentiable 
\begin{align*}
f'(x)(\phi) &= -\intinf dx \; f(x) \; \phi'(x)\nonumber \\
&= - \int_{-\infty}^{a} dx \; f(x) \; \phi'(x) - \int_{a}^{\infty} dx \; f(x) \; \phi'(x) \nonumber \\
&= - f(x) \; \phi(x)\mathlarger{|}^{a}_{-\infty} + \int_{-\infty}^{a} dx \; f'(x) \; \phi(x) \nonumber \\
&- f(x) \; \phi(x)\mathlarger{|}^{\infty}_a +\int^{\infty}_a dx \; f'(x) \; \phi(x) \nonumber \\
&= \phi(a)(\lim\limits_{x \rightarrow a^{+}} f(x) - \lim\limits_{x \rightarrow a^{-}} f(x)) + \intinf dx \; f_c'(x) \; \phi(x) \nonumber \\
&=[f(x)]_{x=a} \; \delta(x -a)(\phi) + f'(x)_{c} (\phi), 
\end{align*}
and we get the identity 
\begin{align*}
f'(x)=[f(x)]_{x=a} \; \delta(x -a) + f'(x)_{c},
\end{align*}
where $f'(x)_{c}$ is a regular generalized function represented by the function $f'(x)$ for all $x \ne a$ and any value we want at $x =a$. Since $f'(x)_{c}$ acts through an integral the actual value of $f_c'(x)$ at $x=a$ does not matter. 

\subsubsection{Fourier transform of generalized functions}

Recall that the Fourier transform and inverse transform are defined as 
\begin{align}
F(\lambda) = \frac{1}{\sqrt{2 \pi}} \; \intinf dx \; \exp(i \lambda x) \; f(x),\nonumber \\
f(x) = \frac{1}{\sqrt{2 \pi}} \; \intinf d \lambda \; \exp(- i \lambda x) \; F(\lambda),  \lbl{200.12}
\end{align}
and also recall that we have the all-important \ttx{convolution theorem}: 
\begin{align}
\intinf d\lambda \; F(\lambda) \; G(\lambda) e^{-i\lambda x}= \intinf dt \; f(t) \; g(x-t). \lbl{201.12} 
\end{align}
Apply \rf{201.12} with $x=0$ and change the dummy variable to $x$ in the right hand side of the expression
\begin{align}
\intinf d\lambda \; F(\lambda) \; G(\lambda) = \intinf dx \; f(x) \; g(-x). \lbl{202.12} 
\end{align}
Let $G(\lambda) = \phi(\lambda)$ be an element of $D_0$. Observe that 
\begin{align}
\Phi(x) \equiv g(-x) = \frac{1}{\sqrt{2\pi}} \intinf d\lambda \; \exp(i \; \lambda x) \; \phi(\lambda). \lbl{203.12}
\end{align}
Thus $\Phi(x)$ is the Fourier transform of the test function $\phi(\lambda)$. Let now $f(x)$ be a regular generalized function with Fourier transform $F(\lambda)$. Then \rf{202.12}, \rf{203.12} shows that 
\begin{align}
F(\lambda) (\phi)=f(x) (\Phi), \lbl{204.12} 
\end{align}
where $\Phi(x)$ is the Fourier transform of $\phi$. It is now tempting to use \rf{204.12} to define the Fourier transform of any generalized function. However \rf{204.12} has a problem. If we take the Fourier transform of a $\phi \in D_0$ it is not always the case that its Fourier transform, $\Phi$, is in $D_0$ and if it's not, the right hand side of \rf{204.12} is not defined. What can go wrong, is that $\Phi$ might not have compact support even if $\phi$ does. \\ 
In order to make Fourier transforms possible, we must enlarge the domain of our generalized functions to include certain functions that do not have compact support. We now rather require that the domain for our distributions should consist of smooth functions that decay rapidly at $\pm \infty$. In fact we require that $\phi(x)$ and all its derivatives decrease faster than any inverse power of $x$ as $\abs{x} \rightarrow \infty$. We call these \ttx{test functions of rapid decay}, or say that they are in the \textit{Schwartz class}. The collection of all Schwartz class functions is denoted by $D_s$. The archetypical example of such a function is the Gaussian $e^{-x^2}$. Clearly any function of compact support belongs to this class and it can be proved that if $\phi$ belongs to this class then its Fourier transform does too. The set of generalized functions based on this new and larger domain, $D_s$, is smaller than the one based on $D_0$. Basically, a regular generalized function corresponding to a function $f(x)$ will be of our new and more restrictive type of generalized function only if $f(x)$ grows slower than exponential at $\pm \infty$. We call this new class \ttx{generalized functions of slow growth}.\label{ft} \\
The function 
\begin{align}
f(x) = 1, \lbl{205.12}
\end{align} 
certainly does not have a Fourier transform in the ordinary sense, because the Fourier transform integrals \rf{200.12} diverges. However, the regular generalized function defined by \rf{205.12} through the association $f \leftrightarrow T_f$ has a Fourier transform. 
\begin{align*}
F(\lambda)(\phi) &= f(x) (\Phi) = \intinf dx \; \Phi(x)\nonumber \\
&=\isp \; \intinf dx \; \intinf d\lambda \; \exp(i \lambda x) \; \phi(\lambda) \nonumber \\ 
&= \sqrt{2 \; \pi} \{ \frac{1}{2\pi} \; \intinf dx \; \intinf d\lambda \exp(-ix(t-\lambda)) \; \phi(\lambda) \} \mathlarger{|}_{t=0}, 
\end{align*}
and from \rf{200.12}, replacing dummy variables, we have 
\begin{align*}
\phi(t) &= \isp \intinf dx \; \exp(-itx) \; \Phi(x)\nonumber \\
&= \frac{1}{2 \; \pi} \; \intinf dx \; \exp(-itx) \; \int_{-\infty}^{\infty} d\lambda \; \exp(i \lambda x) \; \phi(\lambda) \nonumber \\ 
&= \frac{1}{2 \; \pi} \; \intinf dx \; \intinf d\lambda \; \exp(-ix(t-\lambda)) \; \phi(\lambda).
\end{align*}
Thus 
\begin{align*}
F(\lambda)(\phi) = \sqrt{2 \; \pi} \; \phi(0) = \sqrt{2 \; \pi} \; \delta(\lambda) (\phi),
\end{align*}
or 
\begin{align*}
F(\lambda) = \sqrt{2 \; \pi} \; \delta(\lambda).
\end{align*}
As another example let us find the Fourier transform, $D(\lambda)$, of $\delta$ 
\begin{align*}
D(\lambda) (\phi) &= \delta(x) (\Phi) = \Phi(0)\nonumber \\
&= \isp \; \intinf d\lambda \; \phi(\lambda) = \intinf d\lambda \; (\isp) \; \phi(\lambda) \nonumber \\ 
&= (\isp) (\phi), 
\end{align*}
and thus 
\begin{align*}
D(\lambda) = \isp
\end{align*}

\subsubsection{Sequences, series and derivatives}

We have previously, on page \pageref{WeakConvergence},  defined the convergence of a sequence of generalized functions by 
\begin{align*}
&f_n(x) \rightarrow f(x)  &&n\rightarrow\infty,\nonumber \\ 
&\Updownarrow \nonumber \\& f_n(x)(\phi) \rightarrow f(x)(\phi) && n \rightarrow \infty \; \; \; \; \; \; \forall \phi \in D, 
\end{align*}
where $D$ consists of smooth functions of compact support \textit{or} functions decaying fast to zero at $\pm \infty$ as discussed on the previous page. 

Let now a sequence of generalized functions $ \{ f_n(x) \} $ approach a generalized function $f(x)$. Both $f_n(x)$ and $f(x)$ are differentiable since \ttx{all} generalized functions are differentiable. 

But then we have for all $\phi \in D$ 
\begin{align*}
f_n' (x)(\phi) = -f_n(x)(\phi') \rightarrow - f(x)(\phi') = f'(x)(\phi),
\end{align*}
so 
\begin{align}
f_n'(x) \rightarrow f'(x). \lbl{214.12}
\end{align}
This is a very strong statement, the corresponding statement for ordinary functions and derivatives is \ttx{not} true in general. Even if each element of a sequence of functions, $\{ f_n(x) \}$ is smooth the limiting function $f(x) = \lim\limits_{n \rightarrow \infty} f_n(x)$ does not even have to be continuous. The standard example here is the sequence of continuous functions defined by
\begin{align*}
f_n(x) = x^n && 0 \leq x \leq 1,
\end{align*}
whose limit is the discontinuous function 
\begin{align*}
f(x) =\begin{cases}
\; \; 0 \qquad \qquad \qquad \qquad \qquad \qquad \qquad 0 \leq x \leq 1 \\ 
\; \; 1 \qquad \qquad \qquad \qquad \qquad \qquad \qquad x = 1 
\end{cases}. 
\end{align*}
Generalized functions are thus very well behaved with respect to limits.

  The convergence of infinite series of generalized functions is defined in the obvious way. An infinite series 
\begin{align*}
\mathlarger{\sum}^{\infty}_{n=1} f_n(x),
\end{align*}
converges to $f(x)$ iff the sequence of partial sums
\begin{align*}
S_N = \sum^N_{n=1} f_n(x),
\end{align*}
converge to $f(x)$. Thus 
\begin{align*}
f(x) &= \mathlarger{\sum}^{\infty}_{n=1} f_n(x), \nonumber \\
&\Updownarrow \nonumber \\ f(x) &= \lim\limits_{N \rightarrow \infty} S_N(x)\;\;\;\;\;\;\text{where} && S_N(x) = \sum^N_{n=1} f_n(x). 
\end{align*}
The statement \rf{214.12} implies that infinite series of generalized functions can be differentiated term-wise 
\begin{align*}
f(x) = \mathlarger{\sum}^{\infty}_{n=1} f_n(x),\nonumber\\
\Rightarrow f'(x) = \mathlarger{\sum}^{\infty}_{n=1} f_n'(x) 
\end{align*}
This is also a very strong statement that does not hold for ordinary derivatives. 
Let now 
\begin{align}
f(x) = \mathlarger{\sum}^{\infty}_{n=1} f_n(x), \lbl{221.12} 
\end{align}
be an infinite series of functions that converge uniformly in any bounded set. Let $D_0$ be the set of test functions of compact support and let $S_N$ be the sequence of generalized functions defined by 
\begin{align*}
S_N(x) =\sum ^N_{n=1} f_n(x). 
\end{align*}
By assumption, the sequence $S_N(x)$ converges, \textit{as functions on $\mathbf{R}$},  uniformly on any bounded region to some function $f(x)$. But then for all $\phi \in D_0$ we have 
\begin{align}
\lim\limits_{N\rightarrow\infty} S_N(x)(\phi) &= \lim\limits_{N\rightarrow\infty} \intinf dx \; S_N(x) \; \phi(x)\nonumber\\ 
&= \intinf dx \; \lim\limits_{N \rightarrow \infty} S_N(x) \; \phi(x) = \intinf dx \; f(x) \; \phi(x) \nonumber \\
&= f(x) (\phi),  \lbl{221.1.12}  
\end{align}
and thus there exists a generalized function $f(x)$, that is the sum of the infinite series. We have therefore proved that the identity
\begin{align*}
f(x) = \mathlarger{\sum}^{\infty}_{n=1} f_n(x), 
\end{align*}
holds in the sense of generalized functions. Observe that we used uniform convergence when we interchanged limits and integrals in \rf{221.1.12}. Thus we can reinterpret \rf{221.12} in the sense of generalized functions only if we have uniform convergence. 

This reinterpretation of \rf{221.12} in terms of generalized functions is very useful when it can be done. Application of this idea leads to an interpretation of all sorts of very singular series in terms of generalized functions. Consider the series 
\begin{align}
f(x) =\mathlarger{\sum}^{\infty}_{n=1} \; \frac{a_n}{n^2} \; \sin(nx)\label{223.12}
\end{align}
Where $\abs{a_n} < M < \infty \; \; \forall n$. The series (\ref{223.12}) converges uniformly according to the Weierstrass $M$-test. According to our derivative on the previous page we can interpret (\ref{223.12}) in the sense of generalized functions. In this sense the series can be differentiated two times term-wise and we get 
\begin{align}
 f''(x) = -\mathlarger{\sum}^{\infty}_{n=1} a_n \; \sin(nx) \label{224.12}
\end{align}
The series (\ref{224.12}) does not in general converge in the ordinary sense, but it does converge in the sense of generalized functions and in fact represents the second derivative of $f(x)$. As a concrete example of this construction, consider the function defined by 
\begin{align*}
f(x) = \mathlarger{\sum}^{\infty}_{n=1} 1/n^2\sin(n\pi x),
\end{align*}
which is displayed in figure \ref{fig10}. As a generalized function, $f(x)$ is smooth and its second derivative is represented by the series
\begin{align*}
f(x) =- \mathlarger{\sum}^{\infty}_{n=1} \sin(n\pi x), 
\end{align*}
which certainly does not converge in the conventional sense, but which \textit{does} converge in the sense of distributions.

\begin{figure}[htbp]
\centering
\includegraphics{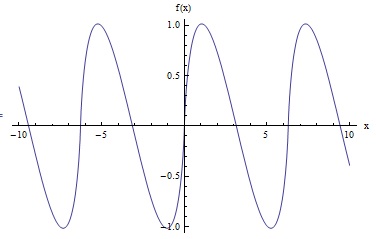}
\caption{A graph of the periodic function $f(x)$}
\label{fig10}
\end{figure}

\noindent Generalized functions are well suited for representing  singular solutions to differential equations. 
\begin{example}
Let us consider the Heaviside function $H(x)$. We have shown that it is differentiable and its derivative is $\delta(x)$. Therefore, for any $\xi$ we have 
\begin{equation*}
 H'(x - \xi) = \delta(x - \xi),
 \end{equation*}
 Thus, the Heaviside generalized function gives us a solution to the equation
 \begin{equation*}
 LG(x;\xi)=\delta(x-\xi),
 \end{equation*}
 where $L$ is the differential operator $L=\frac{d}{dx}$. Thus, the generalized function, $G(x;\xi)=H(x-\xi)$,  is a Green's function to the operator $L$. In this way we can check that proposed functions are Green's functions to an operator by direct substitution. 
\end{example}
\begin{example}
Let a generalized function $U(x,t)$ be given by
\begin{align*}
u(x,t) = H(t - \frac{x}{c}), 
\end{align*}
taking generalized derivatives we have 
\begin{align*}
u_t &= \delta(t - \frac{x}{c}),\nonumber \\ 
u_{tt} &= \delta'(t - \frac{x}{c}), \nonumber \\ 
u_x &= - \inv{c} \; \delta(t - \frac{x}{c}), \nonumber \\
u_{xx} &= \inv{c^2} \; \delta'(t - \frac{x}{c}), \nonumber \\ 
\Downarrow\nonumber\\
 u_{tt} - c^2 \; u_{xx} &= \delta'(t - \frac{x}{c}) - c^2 \; (\inv{c^2} \; \delta'(t - \frac{x}{c})) \nonumber \\
&= \delta'(t - \frac{x}{c}) - \delta'(t - \frac{x}{c}) = 0. 
\end{align*}
Thus $u(x,t)$ is a solution of the wave equation 
\begin{figure}[htbp]
\centering
\includegraphics{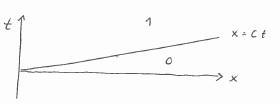}
\caption{A singular solution to the wave equation}
\label{fig6}
\end{figure}
\end{example}

\subsubsection{Properties of the Dirac delta function}

The Dirac delta function is arguably the most important generalized function. It satisfies many interesting identities and have many generalized cousins. We will now discuss a few of these. 

There is one common operation involving ordinary functions that does not extend to generalized functions: Given two functions $f(x), \; g(x)$ we can form the \ttx{product} of $f$ and $g$
\begin{align*}
h(x) = f(x)g(x) \;
\end{align*}
It is, however, not in general possible to form products of generalized functions. Let the following locally integrable function
\begin{align*}
f(x) = \frac{1}{\sqrt{x}},
\end{align*}
be given. Since it is locally integrable, it will  define a generalized function through the association $f \leftrightarrow T_f$.
\begin{align}
(\frac{1}{\sqrt{x}}) (\phi) = \intinf dx \; \inv{\sqrt{x}} \; \phi(x). \lbl{232.12}
\end{align}
The integral \rf{232.12} is well defined as an ordinary improper integral. However the product 
\begin{align*}
h(x) =f(x) f(x) = \inv{x},
\end{align*}
is not a locally integrable function and does not define a generalized function through the association $f \leftrightarrow T_f$ since 
\begin{align*}
h(x) (\phi) = \intinf dx \; \inv{x} \; \phi(x),
\end{align*}
is not meaningful as an ordinary improper integral. There is another way to make $\inv{x}$ into a generalized function using a more general kind of integral called a Cauchy principle value, but the point here is that if you are compelled to write down products of generalized functions, you should be very careful. Your constructions will not automatically make sense and it is up to you to ensure that what you write down is sensible.

Products of Dirac delta functions tend to occur in applications and they sometimes can be given a sense in the theory of generalized functions. For example 
\begin{align*}
\delta(x) \; \delta(x-a) && x \ne a, 
\end{align*}
make sense. In fact we observe that
\begin{align*}
\delta(x) \; \delta(x-a) = 0.
\end{align*}
We can also define the product of two Dirac delta with different variables
\begin{align}
\delta(x) \; \delta(y). \lbl{237.12}
\end{align}
This generalized function now acts on smooth functions of compact support in $\mathbf{R}^2$. This is how \rf{237.12} is given meaning 
\begin{align*}
\delta(x) \; \delta(y) (\phi) = \phi(0,0).
\end{align*}
Thus, whereas $\delta(x)$ and $\delta(y)$ have  domains consisting of functions of one variable, $\delta(x) \; \delta(y)$ has a domain consisting of functions of two variable. By definition,  \rf{237.12} is called the \ttx{outer} product of $\delta(x)$ and $\delta(y)$. \\
This can obviously be generalized to more dimensions. 
\begin{align*}
\delta(x) \; \delta(y) \; \delta(z) \qquad \qquad \qquad \text{etc.}\;\;.
\end{align*}
Another important operation involving ordinary functions is \ttx{composition}. This is possible, to a limited extent, also for generalized functions. We will only consider composition of a generalized function $f(x)$  and an smooth function $g(x)$. We will assume that this smooth function is invertible with inverse $g^{-1}(x)$. This implies that $g'(x)\neq 0$ for all $x$. Let first us assume that $g'(x) > 0$ for all $x$. Then we have 
\begin{align}
f(g(x))(\phi) &= \intinf dx f((g(x)) \; \phi(x)\nonumber \\ 
&= \intinf dy \; f(y)\; \frac{\phi(g^{-1}(y))}{g'(g^{-1}(y))} \nonumber \\ 
&=f(x)(\phi_g),  \lbl{240.12}
\end{align}
where we have used a change of coordinates $y = g(x), \; dy = g'(x)\;dx$. The function $\phi_g(x)$ is defined by 
\begin{equation}
\phi_g(x)=\frac{\phi(g^{-1}(y))}{g'(g^{-1}(y))}\label{fig}.
\end{equation}
If we repeat the above calculations for the case when $g'(x)<0$, we get a formula like (\ref{fig}), excepts for  a minus sign. Both cases can be subsumed into one formula using the absolute value sign
\begin{equation*}
\phi_g(x)=\frac{\phi(g^{-1}(y))}{|g'(g^{-1}(y))|}.
\end{equation*}
From the assumptions we have made about the function $g(x)$ we can conclude that $\phi_g(x)$ is in $D_0$ so that for any generalized function it makes sense to define composition using the formula
\begin{equation*}
f(g(x))(\phi)=f(x)(\phi_g).
\end{equation*}
We now apply this formula to the case when $f(x)=\delta(x)$. For this case we get
\begin{equation*}
\delta(g(x))(\phi)=\delta(x)(\phi_g)=\phi_g(0)=\frac{\phi(g^{-1}(0))}{|g'(g^{-1}(0))|}.
\end{equation*}
Let us now assume that the function $g(x_0)=0$ where $x_0=g^{-1}(0)$. This gives us
\begin{equation*}
\delta(g(x))(\phi)=\frac{\phi(x_0)}{|g'(x_0)|}=\frac{1}{|g'(x_0)|}\delta(x-x_0)(\phi),
\end{equation*}
which leads to the extremely useful identity 
\begin{align*}
\delta(g(x)) = \frac{1}{|g'(x_0)|} \; \delta(x-x_0).
\end{align*}
If $g(x)$ has several isolated zeroes $g(x_n) = 0$ we can repeat \rf{240.12} locally around each zero and get the general formula 
\begin{align*}
\delta(g(x)) = \sum_n \inv{\abs{g'(x_n)}} \; \delta(x-x_n).
\end{align*}
More properties of the Dirac delta can be found in handbooks of mathematical formulas or on the web.

\subsection{Green's functions for the Laplace operator}

A Green's function for the Laplace operator,  $L = - \laplacian{}$,  in $\mathbb{R}^n$, is a function $k(\vb{x};\vb{\xi})$ that solves the equation
\begin{align}
-\laplacian{} k(\vb{x}; \vb{\xi}) = \delta(\vb{x} - \vb{\xi}) && \vb{x} , \; \vb{\xi} \in \mathbf{R}^n. \lbl{250.12} 
\end{align}
Here we are using the n-dimensional Dirac delta generalized function. It can be written as an outer product of $n$ one-dimensional Dirac delta generalized functions 
\begin{align*}
\delta(\vb{x}) = \delta(x_1)\cdots\delta(x_n). 
\end{align*}

In dimension one, $L$ is simply the ordinary differential operator  $L= - \dsx{}$. We have previously, in section \ref{L2}, found all Green's functions for this case. Here, we will concentrate on the case of dimensions higher than one. For these cases we can not find closed form formulas for all possible Green's functions. 

We will start our investigation by finding the appropriate integral identity for the Laplace operator. \\
First observe that we have the formula 
\begin{align}
\div{(\phi \; \grad{\psi})} &= \grad{\phi} \cdot \grad{\psi} + \phi \; \laplacian{\psi},\nonumber \\
\Updownarrow\nonumber\\
 - \phi \; \laplacian{\psi} &= \grad{\phi} \cdot \grad{\psi} - \div{(\phi \; \grad{\psi})}. \lbl{245.12}
\end{align}
Interchanging $\phi$ and $\psi$ in \rf{245.12} we get 
\begin{align*}
\grad{\phi} \cdot \grad{\psi} =\grad{\psi} \cdot \grad{\phi} &=   \div{(\psi \; \grad{\phi})} - \psi \; \laplacian{\phi}, 
\end{align*}
Using these identities, we have for any domain $V \subset \mathbf{R}^n$ with boundary $S$
\begin{align}
\int_V dV \; \phi \; L \; \psi &= \int_V dV \; (-\phi \; \laplacian{\psi})\nonumber \\
&= \int_V dV \; \{  \grad{\phi} \cdot \grad{\psi} - \div{(\phi \; \grad{\psi})} \} \nonumber \\
&= - \int_S dA \; \phi \; \grad{\psi} \cdot \vb{n} + \int_V dV \; \grad{\phi} \vdot \grad{\psi} \nonumber\\
&= - \int_S dA \; \phi \; \grad{\psi} \cdot \vb{n}\nonumber \\
&+ \int_V dV \;\{ \div{(\psi \; \grad{\phi})} - \psi \; \laplacian{\phi}\}\nonumber\\
&= - \int_S dA \; \phi \; \grad{\psi} \cdot \vb{n}+ \int_S dA \; \psi \; \grad{\phi} \cdot \vb{n}\nonumber\\
&+\int_VdV\; \psi \; L{\phi}.  \lbl{246.12}
\end{align}
Thus we get the following fundamental integral identity for the Laplace operator
\begin{align}
\int_V dV \; \{ \phi \; L \; \psi -\psi\; L \; \phi \} = \int_S dA \; \{ \psi \; \prt{\vb{n}} \phi - \phi \; \prt{\vb{n}} \psi \} \lbl{249.12} 
\end{align}
Observe that, just as for the one-dimensional case in section \ref{L2}, the integral identity is derived by (generalized) integration by parts. Integral identities associated with differential operators are \ttx{always} derived using $n$ dimensional generalizations of integration by parts.

Let $\phi(\vb{x})$ be any solution to the equation 
\begin{align}
- \laplacian{\phi(\vb{x})} = F(\vb{x}). \lbl{252.12}
\end{align}
Recall that this is called Poisson's equation. Inserting such a $\phi$ and $\psi(\vb{x}) = k(\vb{x}; \vb{\xi})$ into the integral identity \rf{249.12} gives us 
\begin{align*}
\int_V dV_{\vb{x}} \; \{ \phi(\vb{x}) \; \delta(\vb{x} - \vb{\xi}) - k(\vb{x}; \vb{\xi})\; F(\vb{x}) \}\nonumber \\ 
= \int_S dA_{\vb{x}} \; \{ k(\vb{x}; \vb{\xi}) \; \prt{\vb{n}} \phi(\vb{x}) - \phi(\vb{x}) \; \prt{\vb{n}} k(\vb{x}; \vb{\xi}) \}. 
\end{align*}
Using the fundamental property of the delta function we get
\begin{align}
 \phi(\xi) &= \int_V dV_{\vb{x}} \; k(\vb{x}; \vb{\xi}) \; F(\vb{x})\nonumber\\ 
&+ \int_S dA_{\vb{x}} \; \{ k(\vb{x}; \vb{\xi})\; \prt{\vb{n}} \; \phi(\vb{x}) - \phi(\vb{x}) \; \prt{\vb{n}} k(\vb{x}; \vb{\xi}) \}. \lbl{254.12} 
\end{align}
As previously, the identity \rf{254.12} does not give us a solution to the equation \rf{252.12}, it is merely an integral identity relating values of solutions to the equation \rf{252.12} inside $V$ and on the boundary of $V$. This is the first great theme in the theory of Green's functions introduced on page \pageref{FirstGreatTheme}. As before, \rf{254.12} can be used to find solutions to \rf{252.12} in two distinct ways. We will discuss both approaches in the same order as we did for the simpler operators $L = - \dx{}$ and $L = - \dsx{}$ in section \ref{LL1} and \ref{L2} . \\
Starting with the first approach, let  us look for a solution to \rf{252.12} that satisfies Dirichlet conditions at the boundary of the domain. 
\begin{align}
\phi(\vb{x}) = f(\vb{x}), && \vb{x} \in S = \partial V. \lbl{255.12} 
\end{align}
For this case the integral identity \rf{254.12} is 
\begin{align}
\phi(\vb{\xi}) &= \int_V dV_{\vb{x}} \; k(\vb{x}; \vb{\xi}) \; F(\vb{x})\nonumber \\ 
&+ \int_S dA_{\vb{x}} \; \{ \prt{\vb{x}} \phi(\vb{x}) \; k(\vb{x}; \vb{\xi}) - f(\vb{x}) \; \prt{\vb{n}} k(\vb{x}, \vb{\xi}) \}. \lbl{256.12} 
\end{align}
We next choose a Green's function that satisfies the boundary condition 
\begin{align}
k(\vb{x}; \vb{\xi}) = 0 && \vb{x} \in S. \lbl{257.12} 
\end{align}
Then the unknown boundary data vanish from \rf{256.12} and we get the unique solution to the boundary value problem \rf{252.12}, \rf{255.12} in the form 
\begin{align}
\phi(\vb{\xi}) &= \int_V dV_{\vb{x}} \; k(\vb{x}; \vb{\xi}) \; F(\vb{x})\nonumber \\
&- \int_S dA_{\vb{x}} \; f(\vb{x}) \; \prt{\vb{n}} k(\vb{x}; \vb{\xi}). \lbl{258.12} 
\end{align}
This solution was found using the second great theme in the theory of Green's functions introduced on page \pageref{SecondGreatTheme}. The work remaining is to actually construct the Green's function satisfying the boundary condition \rf{257.12}. 

There are several ways of doing this, depending on the shape and dimension of the domain $V$. Let us first use Fourier series. \\
Let us start by considering the eigenvalue problem.
\begin{align}
- \laplacian{M_k (\vb{x})} &= \lambda_k \; M_k(\vb{x}) && \vb{x} \in V,\nonumber \\
M_k(\vb{x}) &=0 && \vb{x} \in S = \partial V.  \lbl{259.12}
\end{align}
The operator $L =-\laplacian{}$, subject to the given boundary condition,  is self adjoint and positive. For the standard type of boundary conditions used here, the spectral theory of the Laplace operator is well known. Using a source like \cite{textbook}, we conclude that the eigenvalues of \rf{259.12} are real and non-negative 
\begin{align*}
0 \leq \lambda_1 \leq \lambda_2 \leq ...\;\;. 
\end{align*}
The eigenvalues can be enumerated in such a way that to each $\lambda_k$ there is one independent eigenfunction $M_k(\vb{x})$,  and $M_k(\vb{x})$ is orthogonal to  $\; M_j(\vb{x})$ for $k \neq j$. The eigenfunctions are also assumed to be normalized. 

For the Laplace operator we have from \rf{246.12} 
\begin{align*}
\int_V dV \; \phi(- \laplacian{\psi}) = - \int_S dA \; \phi \; \grad{\psi} \vdot \vb{n} + \int_V dV \; \grad{\phi} \vdot \grad{\psi}. 
\end{align*}
Let $\phi = \psi = M_k$. Then we get 
\begin{align}
\int_V dV \; \lambda_k \; M_k^2(\vb{x})& = - \int_S dA \; M_k \; \prt{\vb{n}} M_k + \int_V dV \; \grad{M_k^2},\nonumber \\ 
&\Downarrow\nonumber\\
\int_V dV \; \lambda_k \; M_k^2(\vb{x})& = \int_V dV \; \grad{M_k^2}.  \lbl{262.1.12} 
\end{align}
Assume that the smallest eigenvalue is actually zero, $\lambda_1 = 0$. 
Then, from \rf{262.1.12}, we get 
\begin{align*}
\int_V dV \; \grad{M_1^2}(\vb{x}) &= 0,\nonumber\\ 
&\Updownarrow \nonumber \\ \grad{M_1}(\vb{x}) &= 0, \nonumber \\ 
&\Updownarrow \nonumber \\ M_1 (\vb{x}) &= c.  
\end{align*}
But $M_1(\vb{x}) = 0$ on the boundary. Thus $c = 0 \Rightarrow M_1(\vb{x}) = 0 \; \; \forall \vb{x} \in V$. Therefore we conclude that the smallest eigenvalue is \ttx{not} zero and we have 
\begin{align*}
0 < \lambda_1 \leq \lambda_2 \leq \lambda_3 ...\;\;. 
\end{align*}
We now write the Green's function $k(\vb{x}; \vb{\xi})$ as a Fourier series 
\begin{align*}
k(\vb{x}; \vb{\xi}) = \mathlarger{\sum}_{k=1}^{\infty}\; N_k(\vb{\xi}) \; M_k (\vb{x}).  
\end{align*}
Then $k(\vb{x}; \vb{\xi})$ satisfy the boundary condition \rf{257.12}. We now multiply the  equation for Green's functions \rf{250.12} with the eigenfunction $M_k$ and integrate over the domain $V$. This gives us
\begin{align*}
-\int_{V_{\vb{x}}}dV M_k(\vb{x})\grad_{\vb{x}} ^2\; k(\vb{x}; \vb{\xi}) &= \int_{V}dV_{\vb{x}}M_k(\vb{x})\delta(\vb{x} - \vb{\xi})=M_k(\vb{\xi}),\nonumber \\
&\Updownarrow \nonumber \\ \lambda_k \; N_k(\vb{\xi}) &  = M_k (\vb{\xi}),  
\end{align*}
and thus we must have 
\begin{align*}
N_k(\vb{\xi}) = \frac{M_k(\vb{\xi})}{\lambda_k}, 
\end{align*}
which gives us the following formula for the Green's function
\begin{align}
k(\vb{x}; \vb{\xi}) = \mathlarger{\sum}_{k=1}^{\infty} \; \frac{M_k(\vb{x}) \; M_k (\vb{\xi}))}{\lambda_k}. \lbl{263.12}
\end{align}
We have seen this type of formula for the Green's function appearing in the one-dimensional case in equation \rf{97.12}. For that case we also had a closed form solution for the Green's function that did not involve an infinite sum. Here, the infinite sum representation is the only one we get for a general domain. \\ Observe that by substituting the formula \rf{263.12} for the Green's function into equation \rf{250.12}, and formally applying the Laplace operator term wise, we get the following useful representation for the Dirac delta function
\begin{equation}
\delta(\vb{x}-\vb{\xi})= \mathlarger{\sum}_{k=1}^{\infty} \; M_k(\vb{x}) \; M_k (\vb{\xi})). \lbl{263.1.12}
\end{equation}
If we use \rf{263.12} in \rf{258.12} we get the solution to the boundary value problem \rf{252.12}, \rf{255.12}. 

However, using \rf{258.12} to find numerical values of $\phi$ at selected points, is in general not a small matter. In fact, for general domains, finding approximations for $\lambda_k$ and $M_k(\vb{x})$ is in itself not a small matter.\\
Chapter 8 in the book \cite{textbook} describe some of the standard methods for doing this. \\
We will not pursue these types of methods here,  but will rather illustrate the general theory using symmetric domains where exact formulas for $\lambda_k$ and $M_k(\vb{x})$ can be found. We will exclusively focus on the two-dimensional case in these examples in order for the algebra not to get out of hand.

\begin{example}\label{UnitSquare}
Let $V$ be the unit square in $\mathbf{R}^2$.
\begin{figure}[htbp]
\centering
\includegraphics{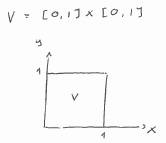}
\caption{The unit square domain for the Laplace equation}
\label{fig7}
\end{figure}

\noindent Our task is to solve the eigenvalue problem 
\begin{align*}
- M_{xx}(x,y) - M_{yy}(x,y) &= \lambda \; M(x,y), \nonumber\\
M(0,y) = M(1,y) &= 0, \nonumber \\
M(x,0) = M(x,1) &= 0. 
\end{align*}
We separate variables using
\begin{align*}
M(x,y) = X(x)\; Y(y). 
\end{align*}
Then the equations can be written 
\begin{align*}
\lambda = - \frac{X''}{X} - \frac{Y''}{Y}. 
\end{align*}
We thus get the following pair of uncoupled equations \label{SimpleBoundaryProblem}
\begin{align}
X'' = - \mu \; X, && X(0) = X(1) = 0,\nonumber \\
Y'' = - \gamma \; Y, && Y(0) = Y(1) = 0, \nonumber \\ 
\lambda = \mu + \gamma. \lbl{268.1.12} 
\end{align}
The boundary value problems \rf{268.1.12} are entirely standard and their solution are
\begin{align*}
X_n(x) = a_n \; \sin(n \pi x) && \mu_n = n^2 \; \pi^2,\nonumber \\
Y_m(y) = b_m \; \sin(m \pi y) && \gamma_m = m^2 \; \pi ^2, \nonumber \\
\lambda_{nm} = n^2 \; \pi^2 + m^2 \; \pi^2, 
\end{align*}
where $n,m=1,2...\;$. 

The appropriate inner product for this problem is
\begin{align}
(\phi , \psi) = \int_0^1 \; \int^1_0 dx \; dy \; \phi(x,y) \; \psi(x,y). \lbl{270.12} 
\end{align}
With respect to this inner product, defined for functions on the unit square with Dirichlet boundary, the Laplace operator is self adjoint. Thus, the set of eigenfunctions  for the Laplace operator $L= - \laplacian{}$, after normalization, forms an orthonormal set of functions given explicitly by
\begin{align*}
M_{nm} (x,y) = 2 \; \sin(n \pi x) \; \sin(m \pi y). 
\end{align*}
The general formula \rf{263.12} then give us the Green's function in the form 
\begin{align}
k(x,y; \xi , \eta) = 4 \; \mathlarger{\sum}^{\infty}_{n,m=1} \frac{\sin(n\pi x) \; \sin(m\pi y) \; \sin(n \pi \xi) \;  \sin(m \pi \eta)}{\pi^2 \; n^2 + \pi^2 \; m^2}. \lbl{272.12}
\end{align}

\end{example}
\begin{example}
Let $V$ be the unit disk.

\begin{figure}[htbp]
\centering
\includegraphics{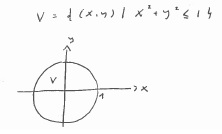}
\caption{The unit disk domain for the Laplace equation}
\label{fig8}
\end{figure}
\noindent Our task is to solve the following boundary problem for the Laplace equation in the unit disk
\begin{align*}
-\laplacian{M(x,y)} &= \lambda, \; M(x,y) && x^2 + y^2 < 1,\nonumber\\
M(x,y) &= 0, && x^2 + y^2 =1. 
\end{align*}
For this problem it is useful to introduce polar coordinates 
\begin{align*}
- \inv{r} \; \prt{r} (r \; \prt{r} M) - \inv{r^2} \; \prt{\theta\theta} M &= \lambda\; M,\nonumber\\
M(1, \theta) &= 0. 
\end{align*}
We separate variables according to
\begin{align*}
M(r, \theta) = R(r) \; \Theta(\theta), 
\end{align*}
and find that the equation can be written as 
\begin{align*}
- \frac{r(r \; R')'}{R} - \lambda \; r^2 = \frac{{\Theta''}}{\Theta}, 
\end{align*}
which leads to the system 
\begin{align}
\Theta'' = -\mu \; \Theta, \lbl{278.12} \\
r(r\; R')' + \lambda \; r^2 \;R = \mu \; R. \lbl{279.12} 
\end{align}
For \rf{278.12} we have the boundary condition 
\begin{align}
\Theta(\theta + 2 \pi ) = \Theta(\theta), && \forall \; \Theta \in [ 0,2\pi ], \lbl{280.12}
\end{align}
and for \rf{279.12} we have 
\begin{align*}
R(1) = 0,\nonumber \\
R(r) \;  \text{bounded as } \; r \rightarrow 0.  
\end{align*}
Let us first focus on the boundary value problem for $\Theta$. We find that there are no eigenvalues in the range $\mu< 0$. For $\mu>0$ the general solution of \rf{278.12} is 
\begin{align}
\theta(\theta) = A \; \cos(\sqrt{\mu} \theta) + B\; \sin(\sqrt{\mu}  \theta), \lbl{282.12} 
\end{align}
and the boundary condition \rf{280.12} implies that 
\begin{align}
\Theta(0) = \Theta (2 \pi),\nonumber\\
\Theta'(0) = \Theta'(2 \pi).   \lbl{283.12} 
\end{align}
\rf{282.12} and \rf{283.12} lead us to the following linear system for $A$ and $B$ 
\begin{align}
\mqty[1 - \cos(2 \pi \sqrt{\mu}) && - \sin(2 \pi \sqrt{\mu}) \\ 
\sin(2 \pi \sqrt{\mu}) && 1 - \cos(2 \pi \sqrt{\mu})] \; \mqty[A \\B] = 0. \lbl{284.12}  
\end{align}
Non-trivial solutions exist only is the determinant of the matrix in \rf{284.12} is zero. 
\begin{align*}
(1 - \cos(2 \pi \sqrt{\mu}))^2 + \sin^2(2 \pi \sqrt{\mu}) &= 0,\nonumber \\
&\Updownarrow \nonumber \\ \cos(2 \pi \sqrt{\mu}) &= 1, \nonumber \\
&\Updownarrow \nonumber \\ \mu_n &= n^2, && n= 1,2,..\;.
\end{align*}
For $\mu = \mu_n$ the linear system \rf{284.12} takes the form 
\begin{align*}
\mqty[0&& 0 \\ 0 && 0] \; \mqty[A \\ B] = 0.  
\end{align*}
The solution space for this system is  two-dimensional . For each $n$, a basis for the solution space is 
\begin{align*}
\{ \cos(n \theta),\sin(n \theta)\}.    
\end{align*}
It is easy to verify that $\mu=0$ is also an eigenvalue with a corresponding basis for the eigenspace given by $\Theta_0(\theta)=1$.  For each $n$, \rf{279.12} turns into the equation 
\begin{align*}
r^2 \; R'' + r \; R' + (\lambda \; r^2 - n^2) \; R = 0. 
\end{align*}
Let 
\begin{align*}
R(r) = \phi(\sqrt{\lambda} r). 
\end{align*}
Then $\phi = \phi(x)$ satisfies the equation 
\begin{align*}
x^2 \; \phi'' + x \; \phi' + (x^2 - n^2) \; \phi = 0. 
\end{align*}
This is Bessel's equation. The space of solutions is two-dimensional  and is spanned by the two Bessel functions $J_n(x),Y_n(x)$. However, only $J_n(x)$ is  bounded at $x=0$ and thus the space of solutions of Bessel's equation that are bounded at the origin is one-dimensional and is  spanned by the Bessel function $J_n(x)$. Thus we get 
\begin{align*}
R(r) \propto J_n (\sqrt{\lambda} r), 
\end{align*}
and the boundary condition at $r=1$ now leads to 
\begin{align*}
J_n(\sqrt{\lambda}) = 0.  
\end{align*}
The infinite set of zeroes of the Bessel function $J_n(x)$  has been tabulated. Let us denote them by 
\begin{align*}
0<\alpha_{n 0} < \alpha_{n 1} < ...\;\;,  
\end{align*}
which finally gives us the eigenvalues 
\begin{align*}
\lambda_{nm} = \alpha_{nm}^2, && n, m = 0, 1, 2, ... \;\;. 
\end{align*}
The corresponding eigenfunctions are 
\begin{align*}
M_{nm} (r, \theta) =\begin{cases}
\; \; a_{nm} \; \cos(n \theta) \; J_n(\alpha_{nm} r) \\
\; \; b_{nm} \; \sin(n \theta) \; J_n(\alpha_{nm} r)
\end{cases}, 
\end{align*}
where $a_{nm}, \; b_{nm}$ are normalization constants that we will now determine. The appropriate inner product for this problem is 
\begin{align*}
(\phi, \psi) = \int^{2\pi}_0 d\theta \; \int^1_0 dr \;r \; \phi(r,\theta) \; \psi(r, \theta).  
\end{align*}
It is then evident that 
\begin{align*}
\; \; a_{nm} \; \cos(n \theta) \; J_n(\alpha_{nm} r),\nonumber \\
\; \; b_{nm} \; \sin(n \theta) \; J_n(\alpha_{nm} r),  
\end{align*}
are orthogonal. From the theory of Bessel functions we have the identity 
\begin{align*}
\int^1_0 dr \; r \; J_n(\alpha_{nm} r_) \; J_n(\alpha_{nm} r) = \frac{1}{2} \; [J_{n+1}(\alpha_{nm})]^2 . 
\end{align*}
Using this identity, we have 
\begin{align*}
& (\cos(n \theta) \; J_n(\alpha_{nm} r), \;\cos(n \theta) \; J_n(\alpha_{nm} r))\nonumber\\
&= \int^{2\pi}_0 d\theta \; \cos^2(n \theta) \; \int^1_0 dr \; r \; J_n (\alpha_{nm} r)^2 \nonumber \\
&= \frac{1}{2} \; \pi ( J_{n+1}(\alpha_{nm}))^2,   
\end{align*}
and
\begin{align}
& (\sin(n \theta) \; J_n(\alpha_{nm} r), \;\sin(n \theta) \; J_n(\alpha_{nm} r))\nonumber \\
&= \frac{1}{2} \; \pi ( J_{n+1}(\alpha_{nm}))^2. \lbl{300.12}
\end{align}
The orthonormal set of eigenfunctions is then
\begin{align*}
M_{nm}(r,\theta) =\begin{cases}
\; \; c_{nm} \; \cos(n \theta) \; J_n(\alpha_{nm} r) \\
\; \; c_{nm} \; \sin(n \theta) \; J_n(\alpha_{nm} r)
\end{cases}, 
\end{align*}
where 
\begin{align*}
c_{nm} = \frac{\sqrt{2}}{\sqrt{\pi} J_{n+1} (\alpha_{nm})}, 
\end{align*}
and the corresponding  Green's function is 
\begin{align*}
k(r,\theta;r',\theta')\nonumber \\
&= \mathlarger{\sum}^{\infty}_{n,m=0} \{ \frac{c_{nm}^2 \; \cos(n \theta) \; \cos(n \theta') \; J_n (\alpha_{nm} r) \; J_n(\alpha_{nm} r')}{\alpha_{nm}^2} \nonumber \\ 
&+ \frac{c_{nm}^2 \; \sin(n \theta) \; \sin(n \theta') \; J_n(\alpha_{nm} r) \; J_n(\alpha_{nm} r')}{\alpha_{nm}^2} \}. 
\end{align*}
Using this in a numerical context is obviously not a simple matter.

  For domains that are generalized cylinders, series that are faster to evaluate can be found using the finite Fourier transform. 
\end{example}
\begin{example}
Let us redo example \ref{UnitSquare} using the Finite Fourier transform. Recall that the domain is the unit square 

\begin{align}
V = [ 0, 1] \times [0,1],\nonumber
\end{align}
and the equation for the Green's function in Cartesian coordinates is 
\begin{align}
 - \laplacian{}k(x,y;\xi,\eta) = \delta (x - \xi) \; \delta(y - \eta), \lbl{305.12} 
\end{align}
where $k=k(x,y;\xi, \eta)$.We want to construct the Green's functions that satisfy Dirichlet conditions at the boundary, $S$, of the square
\begin{eqnarray*}
&k(0,y;\xi,\eta)=0, \;k(1,y;\xi,\eta)=0,\\
&k(x,0;\xi,\eta)=0,\; k(x,1;\xi,\eta)=0.
\end{eqnarray*}
For this purpose, we introduce the ordinary differential operator
\begin{align*}
L = -\prt{xx}, 
\end{align*}
and consider the following eigenvalue problem for this operator
\begin{eqnarray*}
&L M = \lambda M, \\
&M(0) = M(1) = 0.
\end{eqnarray*}
This eigenvalue problem is entirely standard\cite{textbook}. The eigenvalues and normalized eigenfunctions are 
\begin{align*}
\lambda_k = k^2 \; \pi^2 && k =1,2, ...\;,\nonumber \\
M_k(x) = \sqrt{2} \; \sin(k \pi x). 
\end{align*}
The Green's function, expressed using the inverse Finite Fourier transform based on the orthonormal system $\{ M_k \}$, is 
\begin{align}
k(x, y; \xi , \eta) = \mathlarger{\sum}_{k=1}^{\infty} N_k(y; \xi, \eta) \; M_k(x), \lbl{310.12} 
\end{align}
where 
\begin{align*}
N_k(y; \xi, \eta) = \sqrt{2} \; \int^1_0 dx \; \sin(k \pi x) \; k(x , y; \xi , \eta). 
\end{align*}
Multiplying \rf{305.12} by $M_k(x)$ and integrating over the variable $x$, we get 
\begin{align*}
& - \sqrt{2} \; \int^1_0 dx \; \sin(k \pi x) \; \prt{yy}  k(x, y; \xi , \eta)\nonumber\\ 
& - \sqrt{2} \; \int^1_0 dx \; \sin(k \pi x) \; \prt{xx} k(x , y; \xi , \eta) \nonumber \\ 
&= \sqrt{2} \; \int^1_0 dx \; \sin(k \pi x) \; \delta(x - \xi) \; \delta(y - \eta). 
\end{align*}
Using integration by part and the boundary values gives us
\begin{align*}
\prt{yy} \; N_k - k^2 \; \pi^2 \; N_k = - \sqrt{2} \sin(k \pi \xi) \delta(y - \eta). 
\end{align*}
We can rewrite this equation into 
\begin{align}
G_k'' - k^2 \; \pi^2 \; G_k = - \delta(y - \eta), \lbl{314.12} 
\end{align}
where we have introduced 
\begin{align*}
G_k(y; \xi , \eta) = \frac{N_k(y; \xi,\eta)}{\sqrt{2} \; \sin(k \pi \xi)}. 
\end{align*}
Equation \rf{314.12} is subject to the boundary conditions 
\begin{align*}
G_k (0; \eta) = G_k(1; \eta) = 0,
\end{align*}
where we have suppressed the dependence of $G_k$ on the parameter $\xi$, as it plays no active role in the current calculations. Using the approach developed for one-dimensional Green's functions we are lead to the problem 
\begin{align}
G_k'' - k^2 \; \pi^2 \; G_k &= 0, && y \ne \eta, \lbl{317.12} \\ 
G_k(0; \eta) = G_k(1;\eta) &= 0, \lbl{319.12}\\
G_{k\text{+}} (\eta;\eta) - G_{k\text{-}}(\eta;\eta) &= 0,\nonumber \\ 
G'_{k\text{+}} (\eta; \eta) - G'_{k\text{-}} (\eta, \eta) &= -1. \lbl{318.12}
\end{align}
Equation \rf{317.12} implies that 
\begin{align*}
G_k(y;\eta) &= a(\eta) \; \cosh(k \pi (1-y)) + b(\eta) \; \sinh(k\pi(1-y)), && y > \eta,\\ 
G_k(y; \eta) &= c(\eta) \; \cosh(k \pi y) + d(\eta) \; \sinh(k \pi y), && y< \eta.   
\end{align*}
Boundary condition \rf{319.12} implies that 
\begin{align*}
a(\eta) = c(\eta) = 0,  
\end{align*}
and \rf{318.12} leads to the system 
\begin{align*}
b(\eta) \; \sinh(k \pi (1-\eta)) - d(\eta) \; \sinh(k \pi \eta) &= 0,\nonumber \\ 
-k \; \pi \; b(\eta) \; \cosh(k \pi (1-\eta)) - k \; \pi \; d(\eta) \; \cosh(k \pi \eta) &= -1, \nonumber \\
&\Updownarrow \nonumber \\ \mqty[\sinh(k \pi (1-\eta)) && - \sinh(k \pi \eta) \\ 
\cosh(k \pi (1-\eta)) && \cosh(k \pi \eta)] \; \mqty[b(\eta) \\ d(\eta)] = \mqty[0 \\ \inv{k \pi}].  
\end{align*}
The determinant of the matrix is 
\begin{align*}
D &= \sinh(k \pi (1-\eta)) \; \cosh(k \pi \eta) + \sinh(k \pi \eta) \; \cosh(k \pi (1 - \eta))\nonumber\\
&= \sinh(k \pi),   
\end{align*}
and we find 
\begin{align*}
b(\eta) &= \frac{\sinh(k \pi \eta)}{k \pi \sinh(k \pi)},\nonumber\\ 
d(\eta) &= \frac{\sinh(k \pi (1-\eta))}{k \pi \sinh(k \pi)}. 
\end{align*}
Thus the Green's function is 
\begin{align}
N_k(y;\xi, \eta) = \begin{cases}
\sqrt{2} \; \sin(k \pi \xi) \; \frac{\sinh(k \pi (1-y)) \; \sinh(k \pi \eta)}{k \; \pi \; \sinh(k \pi)},  \qquad \qquad y > \eta \\
\sqrt{2} \; \sin(k \pi \xi) \; \frac{\sinh(k \pi y) \; \sinh(k \pi (1-\eta))}{k \; \pi \; \sinh(k \pi)}, \qquad \qquad y < \eta 
\end{cases}.\lbl{326.12}
\end{align}
This formula applies $\forall \; k \geq 1 $. Using \rf{326.12} in \rf{310.12} gives us the Green's function for $L=-\laplacian{}$ on the unit square in the form of a single infinite sum, not a double infinite sum as in \rf{272.12}. It will be much faster to evaluate than formula the original formula \rf{272.12} which contains a doubly infinite sum.

  Formulas like \rf{263.12} can be derived for a wide class of operators and boundary conditions. This might be mathematically complex to do, but in principle it can be done. This is however only if $\lambda = 0$ is not an eigenvalue for the operator. If this happens, and it easily can, the procedure must be modified.

\end{example}
\begin{example}\label{example11}
Let $k(x,y; \xi, \eta)$ be the Green's function for the operator $L= - \laplacian{}$ on the unit square that satisfies Neumann conditions on the boundary 
\begin{align*}
- \laplacian{} k(x,y;\xi,\eta) &= \delta(x-\xi) \; \delta(y - \eta), && (x,y) \in V,\nonumber\\
\prt{\vb{n}} k(x,y;\xi,\eta) &= 0, && (x,y) \in \prt{}V.  
\end{align*}
We will try to solve this using eigenfunctions like in the previous examples.
The relevant eigenvalue problem is 
\begin{align*}
- \laplacian{M(x,y)} &= \lambda \; M(x,y), && 0< x<1 , \; \; 0 < y < 1,\nonumber \\
M_x(0,y) &= M_x(1,y) = 0, \nonumber \\ 
M_y(x,0) &= M_y(x,1) = 0.  
\end{align*}
We use separation of variables 
\begin{align*}
M(x,y) = X(x) \; Y(y),  
\end{align*}
which gives us the equation 
\begin{align*}
\lambda = - \frac{X''}{X} - \frac{Y''}{Y}.  
\end{align*}
The separated equations and boundary conditions are 
\begin{align*}
X'' = -\mu \; X, && X'(0) = X'(1) = 0,  \\
Y'' = -\alpha \; Y, && Y'(0) = Y'(1) = 0,  
\end{align*}
and $\lambda=\mu+\alpha$. We have solved boundary value problems like these several times before. The eigenvalues and corresponding eigenfunctions are 
\begin{align*}
\mu_n &= n^2  \pi^2, && X_n(x) = a_n \; \cos(n \pi x),\nonumber \\
\mu_0 &= 0, && X_0 (x) = a_0, \\
\alpha_m &= m^2  \pi^2, && Y_m(y) = b_m \; \cos(m \pi y),\nonumber\\ 
\alpha_0 &= 0, && Y_0(y) = b_0.   
\end{align*}
Normalizing using the inner product \rf{270.12}, we get the eigenvalues and eigenfunctions 
\begin{align*}
& M_{mn}(x,y) = \begin{cases}
\; \; 1, \qquad \qquad \qquad \qquad \; \; \;  \qquad \qquad \qquad n = m = 0 \\
\; \sqrt{2} \; \cos(n \pi x), \qquad   \qquad  \qquad \qquad \qquad m= 0, n =1,2,... \\
\; \sqrt{2} \; \cos(m \pi y), \qquad   \qquad  \qquad \qquad \qquad n= 0, m =1,2,... \\
\; 2 \; \cos(n \pi x) \; \cos(m \pi y),  \;   \qquad \qquad \qquad n, m = 1,2,... 
\end{cases},\nonumber\\
& \lambda_{mn} = n^2 \; \pi^2 + m^2 \; \pi^2. 
\end{align*}
Since the eigenvalue $\lambda_{00} = 0$,  we can \ttx{not} use formula \rf{263.12} to construct the Green's function. \\
Let us leave this specific example for now, and consider a general situation where the first eigenvalue $\lambda_0 = 0$. Let the corresponding eigenfunction be $M_0 (\vb{x})$, and define a function $\hat{k}(\vb{x};\vb{\xi})$ by 
\begin{align*}
\hat{k}(\vb{x}; \vb{\xi}) = \mathlarger{\sum}^{\infty}_{k=1} \frac{M_k(\vb{x}) \; M_k(\vb{\xi})}{\lambda_k}. 
\end{align*}
Observe that 
\begin{align*}
L \; \hat{k}(\vb{x};\vb{\xi}) &= \mathlarger{\sum}^{\infty}_{k=1} \frac{L \; M_k(\vb{x}) \; M_k(\vb{\xi})}{\lambda_k}\nonumber \\ 
&= \mathlarger{\sum}^{\infty}_{k=1}  M_k(\vb{x}) \; M_k(\vb{\xi}) \nonumber \\ 
&= \mathlarger{\sum}^{\infty}_{k=0}  M_k(\vb{x}) \; M_k(\vb{\xi}) - M_0(\vb{x}) \; M_0(\vb{\xi}) \nonumber \\ 
&= \delta(\vb{x} - \vb{\xi}) - M_0(\vb{x}) \; M_0(\vb{\xi}), 
\end{align*}
where we have used the formal representation of the Dirac delta function introduced in \rf{263.1.12}. 

Thus $\hat{k}(\vb{x}; \vb{\xi})$ is \textit{not} a Green's function for $L= -\laplacian{}$, but is rather a solution to the equation 
\begin{align*}
L \; \hat{k} = \delta(\vb{x} - \vb{\xi}) - M_0(\vb{x}) \; M_0(\vb{\xi}).  
\end{align*}
By definition, $\hat{k}(\vb{x}; \vb{\xi})$ is a \ttx{modified Green's function} for the operator $L= - \laplacian{}$. Returning to  example \ref{example11}, we see that the modified Green's function satisfying Neumann conditions on the boundary of the unit square is given by
\begin{align*}
& \hat{k}(\vb{x}; \vb{\xi})\nonumber\\ 
&= 2 \; \mathlarger{\sum}^{\infty}_{k=1} \frac{\cos(n \pi x) \; \cos(n \pi \xi)}{n^2 \; \pi^2} + 2 \; \mathlarger{\sum}^{\infty}_{k=1} \frac{\cos(m \pi y) \; \cos(m \pi \eta)}{m^2 \; \pi^2} \nonumber \\ 
&+ 4 \; \mathlarger{\sum}^{\infty}_{m, n = 1} \frac{\cos(n \pi x) \; \cos(n \pi \xi) \; \cos(m \pi y) \; \cos(m \pi \eta)}{n^2 \; \pi^2 + m^2 \; \pi^2}.   
\end{align*}

\noindent Using $\psi = \hat{k}$ in the general integral identity \rf{249.12} we get,  instead of the integral identity \rf{254.12}, the following  identity 
\begin{align*}
\phi(\vb{\xi}) &= \int_V dV_{\vb{x}} \; \hat{k}(\vb{x};\vb{\xi}) \; F(\vb{x})\nonumber\\
&+ M_0(\vb{\xi})\;\int_V dV_{\vb{x}}M_0(\vb{x})\phi(\vb{x}) \nonumber \\ 
&+ \int_S dA_{\vb{x}} \; \{ \prt{\vb{n}} \phi(\vb{x}) \; \hat{k}(\vb{x}; \vb{\xi}) - \phi(\vb{x}) \; \prt{\vb{n}} \hat{k}(\vb{x}; \vb{\xi}) \}.
\end{align*}
If we proceed like before, starting with this formula, we do not get explicit solutions by fixing boundary conditions but will rather end up with integral equations for the solutions. This is certainly less satisfying than before, but if $L$ has a zero eigenvalue this is the best we can do. 

In order to gain a better insight into modified Green's functions, let us return to the one-dimensional Laplace operator $L= - \dsx{}$. 
\end{example}
\begin{example}
Recall that the integral identity associated with $L= - \dsx{}$ is 
\begin{align}
\intx dx \{\phi \;  L \; \psi - L \; \phi \; \psi \} = ( \psi\;\phi' -\psi' \;\phi )\mathlarger{|}^{x_1}_{x_0}. \lbl{342.12} 
\end{align}
Let $f(x)$ be a solution to the boundary value problem 
\begin{align*}
-f''(x) &= h(x), && x_0 < x < x_1,\nonumber \\
f'(x_0) &= a_0, \nonumber \\ 
f'(x_1) &= a_1. 
\end{align*}
Let $k(x; \xi)$ be a Green's function for $L= - \dsx{}$, and use $\phi = f(x), \; \psi = k(x; \xi)$ in \rf{342.12}. This gives us 
\begin{align*}
f(\xi) &=\intx dx \; k(x; \xi) \; h(x)\nonumber\\
&+ k(x_1;\xi) \; a_1 - k(x_0;\xi) \; a_0 \nonumber \\
&-k'(x_1; \xi) \; f(x_1) + k'(x_0; \xi) \; f(x_0). 
\end{align*}
In order to find a solution we must use a Green's function solving the problem
\begin{align}
- k''(x;\xi) &= \delta(x - \xi),\nonumber \\ 
k'(x_0; \xi) &= k'(x_1;\xi) =0. \lbl{345.12} 
\end{align}
We now solve \rf{345.12} using Fourier series  and therefore introduce the eigenvalue problem 
\begin{align}
- M''(x) &= \lambda \; M(x), && x_0 < x <x_1,\nonumber\\
M'(x_0) &= M'(x_1) = 0.  \lbl{346.12}  
\end{align}
For $\lambda > 0$ the general solution to \rf{346.12} is 
\begin{align}
M(x) &= A \; \cos(\sqrt{\lambda} x) + B \; \sin(\sqrt{\lambda} x),\lbl{348.12} \\
\Downarrow\nonumber\\
 M'(x) &= - \sqrt{\lambda} \; A \; \sin(\sqrt{\lambda} x) + \sqrt{\lambda} \; B \; \cos(\sqrt{\lambda} x). \nonumber
 \end{align}
 Applying the boundary conditions we get
 \begin{align*}
M'(x_0) &= 0 \; \; \Rightarrow \; \; -A \; \sin(\sqrt{\lambda} x_0) + B \; \cos(\sqrt{\lambda} x_0) = 0, \nonumber \\
M'(x_1) &= 0 \; \; \Rightarrow \; \; -A \; \sin(\sqrt{\lambda} x_1) + B \; \cos(\sqrt{\lambda} x_1) = 0, \nonumber 
\end{align*}
and we thus have the homogeneous linear system 
\begin{align*}
\mqty[-\sin(\sqrt{\lambda} x_0) && \cos(\sqrt{\lambda} x_0) \\ -\sin(\sqrt{\lambda} x_1)  
&& \cos(\sqrt{\lambda} x_1)] \; \mqty[A \\ B] = 0.   
\end{align*}
In order to have non-trivial solutions, the determinant of the matrix must be zero 
\begin{align*}
- \cos(\sqrt{\lambda} x_1) \; \sin(\sqrt{\lambda} x_0) + \sin(\sqrt{\lambda} x_1) \; \cos(\sqrt{\lambda} x_0) &= 0,\nonumber \\ 
&\Updownarrow \nonumber \\ \sin(\sqrt{\lambda} (x_1 - x_0)) &= 0, \nonumber \\ 
&\Updownarrow \nonumber \\ \lambda_n &= (\frac{n \pi}{l})^2 && n=1,2,...\;\;,  
\end{align*}
where $l=x_1 - x_0$. For these values of $\lambda$ the system reduces to a single independent equation 
\begin{align*}
-A \; \sin(\sqrt{\lambda_n} x_0) + B \; \cos(\sqrt{\lambda_n} x_0) &= 0,\nonumber \\ 
&\Updownarrow \nonumber \\ B &= \frac{\sin(\sqrt{\lambda_n} x_0)}{\cos(\sqrt{\lambda_n} x_0)} \; A.  
\end{align*}
Inserting this into \rf{348.12} gives us, after using trigonometric addition formulas, 
\begin{align*}
M_n(x) = C_n \; \cos(\frac{n\pi}{l}(x-x_0)).  
\end{align*}
The appropriate inner product for this case is 
\begin{align*}
(\phi, \psi) = \intx dx \; \phi(x) \; \psi(x).  
\end{align*}
Using this we get an orthonormal system of eigenfunctions 
\begin{align}
M_n(x) &= \sqrt{\frac{2}{l}} \;\cos(\frac{n\pi}{l}(x-x_0)) && n=1,2,...\;\;,\nonumber \\ 
\lambda_n &= (\frac{n \pi}{l})^2. \lbl{354.12}
\end{align}
The system \rf{354.12} is not complete because $\lambda = 0$ is also an eigenvalue with a corresponding eigenfunction $M_0(x)$ given by 
\begin{align*}
M_0(x) = \frac{1}{\sqrt{l}}.  
\end{align*}
Because of the zero eigenvalue, we cannot write down a Green's function using the general formula \rf{263.12}. The modified Green's function for $L = - \dsx{}$ is however given by
\begin{align*}
\hat{k}({x}; {\xi}) = \mathlarger{\sum}^{\infty}_{n=1} \frac{2\;l\;\cos(\frac{n\pi}{l}(x-x_0))\;\cos(\frac{n\pi}{l}(\xi-x_0))}{n^2 \; \pi^2}.  
\end{align*}
The modified Green's function satisfies the equation 
\begin{align*}
-\hat{k}''(x; \xi) = \delta(x - \xi) - \frac{1}{l}.  
\end{align*}
If we now use $\phi = f(x), \; \psi = \hat{k}(x; \xi)$ in the integral identity \rf{342.12} we get 
\begin{align}
f(\xi) &=  \intx dx \; \hat{k}(x; \xi) \; h(x) + \frac{1}{l} \; \intx dx \; f(x)\nonumber \\ 
&+ \hat{k}(x_1; \xi) \; a_1 - \hat{k}(x_0;\xi) \; a_0, \label{integralEq} 
\end{align}
where we have used the fact that the modified Green's function $\hat{k}({x}; {\xi})$ satisfies the boundary conditions 
\begin{align*}
\hat{k}'(x_1;\xi) = \hat{k}(x_0; \xi) = 0.  
\end{align*}
Evidently (\ref{integralEq}) is an integral equation for $f(x)$ which can be written 
\begin{align}
f(\xi) - \inv{l} \; \intx d\xi \; f(\xi) = g(\xi),  \lbl{360.12}
\end{align}
where 
\begin{align*}
g(\xi) &=  \intx dx \; \hat{k} (x; \xi) \; h(x)\nonumber\\ 
&+ \hat{k}(x_1; \xi) \; a_1 - \hat{k} (x_0;\xi) \; a_0.   
\end{align*}
Let $f(\xi) = f_0$. Then we have 
\begin{align*}
f(\xi) - \inv{l} \; \intx d\xi \; f(\xi) = f_0 - \inv{l} \; f_0 \; \intx dx = f_0 - f_0 = 0.  
\end{align*}
Thus the integral equation \rf{360.12} is \ttx{singular}, $f(\xi) = f_0$ is in the kernel of the integral operator. This is in general true for the integral equations that appear when we are working with modified Green's functions and  is a complication which  means that equation \rf{360.12} will in general have no solution. \\ 
In order to find out exactly when \rf{360.12} has a solution, we will introduce a piece of linear algebra that is of great utility in applied mathematics. In section \ref{CCO} we used it to find the solvability conditions for the perturbation hierarchy corresponding to a system of  two weakly coupled cubic oscillators. Here it appears again, in a very different mathematical context. It is probably the most useful piece of linear algebra that you don't learn in a standard course in linear algebra, and it is well worth, for a second time in these lecture notes, to  sum up it's most salient features.

\subsubsection{The Fredholm alternative} 

Let $V$ be a vector space, which may be of infinite dimension, and let $A$ be a linear operator. In the infinite dimensional case $A$, might be an integral or differential operator. \\
We would like to know when the linear system 
\begin{align}
A \; x = b, \lbl{366.12} 
\end{align}
has a solution. Here, $b$ is some vector in $V$. Both differential and integral equations can be written in the form \rf{366.12} and is thus covered by the Fredholm alternative. We will assume that $V$ is an inner product space with an inner product denoted by $(x,y)$. Recall that $x, y$ might be functions for the case when $V$ is infinite dimensional. \\
The \ttx{adjoint} of $A$, denoted by $A^{\dagger}$, is the unique linear operator such that 
\begin{align*}
(A \; x, y) = (x, A^{\dagger} \; y) && \forall \; x, \;y \in V. 
\end{align*}
In the infinite dimensional case one should really worry about the domain of definition for $A$ and $A^{\dagger}$. They are in general not defined on the whole $V$ unless they are bounded. Integral operators are often bounded, whereas differential operators are always unbounded. \\
Pursuing these kinds of issues really belongs in a class in mathematical analysis, and I will not talk more about them here. Here, I will concentrate on the algebra, not the analysis. Let us assume that $A^\dagger$ is not invertible. Then there exists vectors $x^* \in V$ such that 
\begin{align*}
A^{\dagger} \; x^* = 0.  
\end{align*} 
For any such $x^*$ we have 
\begin{align*}
(x^*, b ) = (x^*, A\; x) = (A^{\dagger} \; x^*, x) = (0,x) =0.  
\end{align*}
Thus a necessary condition for the system
\begin{align}
A \; x = b, \lbl{367.12} 
\end{align}
to have a solution is that $ \forall \; x^* $ such that 
\begin{align*}
A^{\dagger} \; x^* = 0, 
\end{align*}
we must have 
\begin{align*}
(x^*, b) = 0.  
\end{align*}
This is the \ttx{Fredholm alternative}. To prove that it is also sufficient in the infinite dimensional case requires mathematical analysis. However, here,  we will not worry about this, and just assume that the Fredholm alternative is also sufficient for solvability of \rf{367.12}. 

After this piece of very useful linear algebra, we now return to our example. The vector space $V$ is some reasonable space of functions defined on $[ x_0, x_1 ]$ and the inner product is 
\begin{align*}
(\phi, \psi) = \intx dx \; \phi(x) \; \psi(x).  
\end{align*}
The operator $A$ is 
\begin{align*}
A(f) = f(\xi) - \inv{l} \; \intx dx \; f(\xi).  
\end{align*}
We need the adjoint of $A$ 
\begin{align*}
(A \; \phi, \psi) &= \intx d\xi \; A \; \phi(\xi) \; \psi(\xi)\nonumber \\ 
&= \intx d \xi \; (\phi(\xi) - \inv{l} \; \intx d \xi' \; \phi(\xi') ) \; \psi(\xi) \nonumber\\ 
&= \intx d\xi \; \phi(\xi) \; \psi(\xi) - \inv{l} \; \intx d\xi \; \psi(\xi) \; \intx d\xi' \; \phi(\xi') \nonumber\\
&= \intx d\xi \; \phi(\xi) \; \psi(\xi) - \inv{l}\; \intx d\xi \; \phi(\xi) \; \intx d\xi' \; \psi (\xi') \nonumber \\ 
&= \intx d\xi \; \phi(\xi) \;(\psi(\xi) - \inv{l} \; \intx d\xi' \; \psi(\xi')) \nonumber \\
&= (\phi, A\; \psi). 
\end{align*}
Thus $A$ is self-adjoint, $A^{\dagger} = A$. In order to apply the Fredholm alternative we must now find the kernel of $A^{\dagger} = A$. This amounts to finding all solutions to the equation 
\begin{align}
f(\xi) - \inv{l} \; \intx d\xi \; f(\xi) &= 0,\nonumber \\
\nonumber\\
&\Updownarrow \nonumber \\ 
\nonumber\\
f(\xi) &= \inv{l} \; \intx d\xi \; f(\xi) = \text{const}. \lbl{373.12} 
\end{align}
Thus any solution of \rf{373.12} must be constant, $f(\xi) = f_0$. But we have already proved that any such constant \ttx{is} a solution. Thus $f(\xi)$ is in the kernel of $A^{\dagger}$ iff 
\begin{align*}
f(\xi) = f_0 && \forall \; \xi \in [x_0, x_1]. 
\end{align*}
The Fredholm alternative now gives us the single {\it solvability} condition 
\begin{align*}
& (1,g) = 0,\nonumber\\
&\Updownarrow \nonumber \\ & \intx d\xi \; \{ \hat{k}(x_1; \xi) \; a_1 - \hat{k} (x_0; \xi) \; a_0 +\intx dx \; \hat{k} (x; \xi) \; h(x) \} = 0. 
\end{align*}

\noindent We have now seen some of the complications that can arise when we try to apply the second great theme in the theory of Green's function to solve boundary value problems for differential operators. The message is that finding Green's functions satisfying specific boundary conditions, is not by any means easy and straight forward, even for an operator as simple as $L = -\laplacian{}$. 

\noindent We will now try to deploy the third great theme in the theory of Green's functions for the operator $L = - \laplacian{}$. Within this theme our choice of Green's function is much less constrained, for all practical purposes any Green's function will do. This is very important, because then we can simplify the equation for the Green's function using symmetry. 

Recall that Green's functions for the Laplace operator $L= - \laplacian{}$ are solutions to the equation 
\begin{align}
-\laplacian{k(\vb{x}; \vb{\xi})} = \delta(\vb{x} - \vb{\xi}). \lbl{376.12}
\end{align}
Observe that if we can find a solution to the equation 
\begin{align}
- \laplacian{k(\vb{\eta})} = \delta(\vb{\eta}), \lbl{377.12}
\end{align}
then 
\begin{align}
k(\vb{x} - \vb{\xi}), \lbl{378.12} 
\end{align}
will be a solution to \rf{376.12}. We will therefore focus on equation \rf{377.12}. First observe that any solution of \rf{377.12} satisfies the equation 
\begin{align}
 - \laplacian{k(\vb{\eta})} = \vb{0}, && \vb{\eta} \ne \vb{0}. \lbl{379.12}
\end{align}
The statement \rf{379.12} and all following statements can be justified through the theory of generalized functions but here I prefer to proceed heuristically. Thus, the Dirac delta is assumed to be a function satisfying 
\begin{align*}
\delta(\vb{x}) = 0, && \vb{x} \ne \vb{0},\nonumber \\
\int_V dV \; \delta(\vb{x}) = 1, && \text{if} \; \; \vb{0} \in V. 
\end{align*}
In the following we will focus on the two-dimensional case. \\
Let $S_{\epsilon}$ be a circular disk of radius $\epsilon$ centered on $\vb{\eta} = \vb{0}$. Integrating \rf{377.12} over $S_{\epsilon}$ gives us
\begin{align*}
- \int_{S_{\epsilon}} dA \; \div{(\grad{k})}(\vb{\eta}) = \int_{S_{\epsilon}} dA \; \delta(\vb{\eta}) &= 1,\nonumber\\
&\Updownarrow\nonumber\\
 \int_{C_{\epsilon}} dl \; \grad{k(\vb{\eta})} \cdot \vb{n}& = -1, 
\end{align*}
where $C_{\epsilon}$ is a circle of radius $\epsilon$ and centered on $\vb{\eta} = 0 $. Taking the limit as $\epsilon$ approach zero we get the following constraint satisfied by all Green's function of $L = - \laplacian{}$ 
\begin{align}
\lim\limits_{\epsilon \rightarrow 0} \int_{C_{\epsilon}} dl \; \grad{k(\vb{\eta})} \cdot \vb{n} = -1. \lbl{384.12} 
\end{align}
We will now try to find a solution of \rf{377.12} that is rotationally invariant
\begin{align}
k = k(r). \lbl{385.12} 
\end{align}
For such a function \rf{384.12} simplifies into 
\begin{align}
\lim\limits_{\epsilon \rightarrow 0} \int^{2 \pi}_0 d\theta \; \epsilon \; \prt{r} k(\epsilon) &= -1,\nonumber\\ 
&\Updownarrow \nonumber \\ \lim\limits_{\epsilon \rightarrow 0} \epsilon \; \prt{r} k(\epsilon) &= - \inv{2 \pi}. \lbl{387.12} 
\end{align}
Writing the equation \rf{379.12} in polar coordinates and using \rf{385.12}, we get the equation
\begin{align}
\inv{r} \; \prt{r} (r \; \prt{r} k) &= 0, && r \ne 0,\nonumber \\ 
&\Updownarrow \nonumber \\ r \; \prt{r} k &= c, \nonumber \\
&\Updownarrow \nonumber \\ \prt{r} k &= \frac{c}{r}. \lbl{389.12}
\end{align}
One solution of \rf{389.12} is 
\begin{align}
k(r) = c \; \ln(r). \lbl{390.12} 
\end{align}
Applying the condition \rf{387.12} for \rf{390.12} we get 
\begin{align*}
\lim\limits_{\epsilon \rightarrow 0} \epsilon \; \frac{c}{\epsilon} &= - \inv{2 \pi},\nonumber \\
&\Updownarrow \nonumber \\ c &= -\inv{2 \pi}.  
\end{align*}
Thus a  rotationally invariant solution to \rf{377.12} is 
\begin{align*}
k(r) = - \inv{2 \pi} \; \ln(r), 
\end{align*}
which in Cartesian coordinates is 
\begin{align*}
k(\vb{\eta}) = -\inv{2 \pi} \; \ln(\norm{\vb{\eta}}).  
\end{align*}
Thus using \rf{378.12} we get the following Green's function for $L = - \laplacian{}$ in two dimensions 

\begin{align}
k(x; \vb{\eta}) = - \inv{2 \pi} \; \ln(\norm{\vb{x} - \vb{\xi}}). \lbl{399.1.12}
\end{align}
Let us now return to the challenge of solving the boundary value problem 
\begin{align}
- \laplacian{\phi(\vb{x})} &= F(\vb{x}), && \vb{x} \in V \subset \mathbf{R}^2,\nonumber\\
\phi(\vb{x}) &= f(\vb{x}), && \vb{x} \in S = \prt{} V. \lbl{395.12}  
\end{align}
For any Green's function for the Laplace operator, in particular for \rf{399.1.12}, we have the integral identity 
\begin{align}
\phi(\vb{\xi}) &= \int_V dA_{\vb{x}} \; k(\vb{x}; \vb{\xi}) \; F(\vb{x})\nonumber\\ 
&+ \int_S dl_{\vb{x}} \; \{ \prt{\vb{n}} \phi(\vb{x}) \; k(\vb{x}; \vb{\xi}) - f(\vb{x}) \; \prt{\vb{n}} k(\vb{x}; \vb{\xi})  \}, \lbl{396.12} 
\end{align}
where $\phi(\vb{\xi})$ is the unique solution of \rf{395.12}. Observe that \rf{396.12} does not give us an explicit solution to \rf{395.12} since the boundary data $\prt{\vb{n}} \phi(\vb{x})$ is only known when the unique solution to \rf{395.12} is known. \\
The idea is now to get a closed equation for the boundary data $\prt{\vb{n}} \phi(\vb{x})$ by evaluating \rf{396.12} for $\vb{\xi}$ on $S = \prt{} V$. There is however a complication,
if we substitute a $\xi \in S$ into the curve integral over $S$ on the right hand side of \rf{396.12} we will end up having to evaluate
\begin{align*}
k(\vb{\xi};\vb{\xi}) = - \inv{2 \pi} \; \ln(\norm{\vb{\xi} - \vb{\xi}}) = - \inv{2 \pi} \; \ln(0), 
\end{align*}
which does not make sense. The way to resolve this problem is to evaluate \rf{396.12} on the boundary through a limit process. \\
There are many ways of doing this, but they all give the same equation, so I just pick the simplest one.
\begin{figure}[htbp]
\centering
\includegraphics{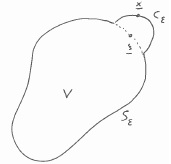}
\caption{Regularizing the boundary integral using a boundary deformation approach}
\label{fig9}
\end{figure}

\noindent We pick a point on the boundary $S$ and place $\vb{\xi}$ at this point. We then deform the boundary into a small semicircle centered on $\vb{\xi}$ and of radius $\epsilon$. This is $C_{\epsilon}$ in figure \ref{fig9}. The part of $S$ with the dotted piece removed we call $S_{\epsilon}$. We now evaluate \rf{396.12} at $\vb{\xi} \in S$ through the following limit process
\begin{align*}
f(\vb{\xi}) &= \lim\limits_{\epsilon \rightarrow 0} \int_{S_{\eps} \cup C_{\eps}} dl_{\vb{x}} \; \{ \prt{\vb{n}} \phi(\vb{x}) \; k(\vb{x};\vb{\xi})\nonumber \\  
&- f(\vb{x}) \; \prt{\vb{n}} k(\vb{x};\vb{\xi})  \} \nonumber \\
&+ \int_V dA_{\vb{x}} \; k(\vb{x}; \vb{\xi}) \; F(\vb{x}).
\end{align*}
Introduce the notion of a \textit{Cauchy Principal Value} integral using
\begin{align}
PV_{\vb{\xi}} \; \int_S dl_{\vb{x}} \; \{ \; \} = \lim\limits_{\eps \ra 0} \int_{S_{\eps}} dl_{\vb{x}} \; \{ \; \}, \lbl{399.12}
\end{align}
and also introduce 
\begin{align*}
R_{\eps} = \int_{C_{\eps}} dl_{\vb{x}} \; \{ \prt{\vb{n}} \phi(\vb{x}) \; k(\vb{x}; \vb{\xi}) - f(\vb{\xi}) \; \prt{\vb{n}} k(\vb{x}; \vb{\xi}) \}. 
\end{align*}
Using this notation we have 
\begin{align}
f(\xi) &= \int_V dA_{\vb{x}} \; k(x; \vb{\xi}) \; F(\vb{x})\nonumber \\
&+ PV_{\vb{\xi}} \; \int_S dl_{\vb{x}} \; \{ \prt{\vb{n}} \phi(\vb{x}) \; k(\vb{x}, \vb{\xi}) - f(\vb{x}) \; \prt{\vb{n}} k(\vb{x}, \vb{\xi})  \} \nonumber \\ 
&+ \lim\limits_{\eps \ra 0} R_{\eps}. \lbl{401.12}
\end{align}
We will shortly compute the last limit, but before that let me say a few words about Cauchy principal value integrals 
\\
\subsubsection{Cauchy principal value integrals} 

Let $f$ be a continuous function on a closed interval $[a,b]$. Then the usual Riemann integral of $f$ exist 
\begin{align*}
\int^b_a dx \; f(x) \in \mathbf{R}.
\end{align*}
If $f$ has a vertical asymptote at one or more points in $[a,b]$ the usual Riemann integral will not exist. \\
Let us for simplicity assume that there is a single vertical asymptote at $x_0$ with $a<x_0<b$. The integral from, $a$ to $b$ of $f$ is now defined by 
\begin{align}
\int^b_a dx \; f(x) = I_1 + I_2, \lbl{403.12}
\end{align}
where
\begin{align}
I_1 &= \lim\limits_{\eps \ra 0} \int_a^{x_0 - \eps} dx \; f(x),\nonumber\\ 
I_2 &= \lim\limits_{\delta \ra 0} \int_{x_0 + \delta}^{b} dx \; f(x). \nonumber 
\end{align}
The integral of $f$ from $a$ to $b$ exists as an \ttx{improper} integral if \ttx{both} $I_1$ and $I_2$ exist and then the value of the improper integral is given by \rf{403.12} 
\end{example}
\begin{example}
Let $a=-1, \; b=1$ and 
\begin{align*}
f(x) = \inv{\sqrt{\abs{x}}} && x \ne 0.  
\end{align*}
Then 
\begin{align*}
I_1 &= \lim\limits_{\eps \ra 0} \int_{-1}^{-\eps} dx \; \inv{\sqrt{-x}} =- \lim\limits_{\eps \ra 0} 2 \; \sqrt{-x}\mid^{-\eps}_{-1}\nonumber \\
&= \lim\limits_{\eps \ra 0} (2-2 \; \sqrt{\eps}) = 2, \nonumber \\ 
I_2 &= \lim\limits_{\delta \ra 0} \int_{\delta}^{1} dx \; \inv{\sqrt{x}} = \lim\limits_{\delta \ra 0} 2 \; \sqrt{x}\mid_{\delta}^1 \nonumber\\
&= \lim\limits_{\delta \ra 0} (2-2 \; \sqrt{\delta}) = 2. 
\end{align*}
Thus the improper integral $\int^1_{-1} dx \; \inv{\sqrt{\abs{x}}}$ exists and has the value 
\begin{align*}
\int^1_{-1} dx \; \inv{\sqrt{\abs{x}}} = 2 +2 = 4. 
\end{align*}
\end{example}
\begin{example}\label{Example47.1}
Let $a=-1, \; b=1$ and 	
\begin{align*}
f(x) = \inv{x}, && x \ne 0.  
\end{align*}
Then 
\begin{align*}
I_2 &= \lim\limits_{\delta \ra 0} \int_{\delta}^{1} dx \inv{x} = \lim\limits_{\delta \ra 0} \ln(x)\mathlarger{|}_{\delta}^{1}\nonumber \\ 
&= \lim\limits_{\delta \ra 0} (-\ln \delta ).   
\end{align*}
This limit does not exist so $\int^1_{-1} dx \; \inv{x} $ does not exist as an improper integral. 
Integrals like the one in this example can however be given a meaning as a Cauchy principal value integral. 

In general let $f:[a,b] \ra \mathbf{R}$ be a function with a vertical asymptote at $x_0, \; a< x_0 <b$. Then if the limit 
\begin{align}
I = \lim\limits_{\eps \ra 0} \{ \int^{x_0 - \eps}_a dx \; f(x) + \int^b_{x_0 + \eps} dx \; f(x) \}, \lbl{411.12} 
\end{align}
exist then $\int^b_a dx \; f(x)$ exists as a Cauchy principal value integral which we write as 
\begin{align*}
PV_{x_0} \; \int^b_a dx \; f(x) \equiv I.  
\end{align*}
\end{example}
\begin{example}
Let us return to the function from example \ref{Example47.1}. Using the limit \rf{411.12} we have
\begin{align*}
I &= \lim\limits_{\eps \ra 0} \{ \int^{-\eps}_{-1} dx \; \inv{x} + \int^1_{\eps} dx \; \inv{x}  \}\nonumber\\
&= \lim\limits_{\eps \ra 0} \{ - \int^1_{\eps} dy \; \inv{y} + \int^1_{\eps} dx \; \inv{x}   \} \nonumber\\ 
&= \lim\limits_{\eps \ra 0} \{0\} = 0,  
\end{align*} 
so the integral exists as a Cauchy principal value integral and has the value zero 
\begin{align*}
PV_0 \int_{-1}^{1} dx \; \inv{x} = 0.  
\end{align*}
Observe that the difference between the definition of an improper integral \rf{403.12} and a Cauchy principal value integral is that the limit is taken in a \ttx{symmetric} way for Cauchy principal value integrals. This allows for the possibility of canceling infinite terms that occurs with opposite signs. 
\end{example}
\begin{example}
Let $a=-1, \; b=1$ and 
\begin{align*}
f(x) = \inv{\abs{x}}.  
\end{align*}
We then have 
\begin{align*}
I &= \lim\limits_{\eps \ra 0} \{ \int^{-\eps}_{-1} dx \; \inv{\abs{x}} + \int^1_{\eps} dx \inv{\abs{x}} \}\nonumber \\
&= \lim\limits_{\eps \ra 0} \{ - \int^{-\eps}_{-1} dx \; \inv{x} + \int^1_{\eps} dx \; \inv{x} \} \nonumber \\ 
&= \lim\limits_{\eps \ra 0} \{\int^1_{\eps} dx \; \inv{x} + \int^1_{\eps} dx \; \inv{x} \} \nonumber \\
&=- 2 \; \lim\limits_{\eps \ra 0} \ln(\eps)= \infty.
\end{align*}
\end{example}
\noindent Thus the integral $\int^1_{-1} dx \; \inv{\abs{x}}$ does not exist as a Cauchy principal value integral either. In a sense it has the wrong kind of singular behavior. There is an even more general notion of singular integral called a Hadamard integral that can take care of some integrals that do not exist as Cauchy principal value integrals. We will not pursue this topic any further here. Observe the limit defined in \rf{399.12} is a symmetric limit, and when $S_{\eps}$ is parametrized, the limit will exactly define a Cauchy principal value integral. 

Let us now return to the evaluation of $R_{\eps}$ in the limit $\eps \ra 0$. Recall that 
\begin{align*}
k(\vb{x}; \vb{\xi}) = - \inv{2 \pi} \; \ln \norm{ \vb{x} - \vb{\xi}}.  
\end{align*}
The normal derivative is 
\begin{align*}
\prt{\vb{n}} k(\vb{x}; \vb{\xi}) &= \vb{n}(\vb{x}) \cdot (- \inv{2 \pi} \; \inv{\norm{ \vb{x}- \vb{\xi}}}) \; \grad{\norm{\vb{x} - \vb{\xi}}}\nonumber\\
&= - \vb{n}(\vb{x}) \vdot \frac{\vb{x} - \vb{\xi}}{2 \pi \norm{\vb{x} - \vb{\xi}}^2}.  
\end{align*}
But since $S_{\eps}$ is a semi-circle of radius $\eps$ centred on $\vb{\xi}$ we have 
\begin{align*}
\vb{n} (\vb{x}) = \frac{\vb{x} - \vb{\xi}}{\norm{\vb{x} - \vb{\xi}}}. 
\end{align*}
For $\vb{x}$ on the semi-circle of radius $\eps$ and center located at $\vb{\xi}$, we have
\begin{align*}
\norm{\vb{x} - \vb{\xi}} = \eps,  
\end{align*}
and therefore 
\begin{align*}
\prt{\vb{n}} k(\vb{x};\vb{\xi}) = - \inv{2 \pi \eps}.  
\end{align*}
The semi-circle $C_{\eps}$ is parametrized by 
\begin{align*}
\vb{\gamma}(t) &= \vb{\xi} + \eps ( \cos \theta, \sin \theta), && 0 \leq \theta \leq \pi 
\end{align*}
and we thus have
\begin{align*}
\vb{\gamma}'(t) &= \eps (-\sin \theta , \cos \theta ), \nonumber \\
\Downarrow\nonumber\\ 
dl_{\vb{x}} &= \eps \; d\theta.   
\end{align*}
Therefore 
\begin{align*}
\int_{C_{\eps}} dl_{\vb{x}} \; f(\vb{x}) \; \prt{\vb{n}} k(\vb{x}; \vb{\xi}) &= \int^{\pi}_0 d\theta \; \eps (-\inv{2 \pi \eps}) \; f(\vb{\xi} + \eps(\cos\theta,sin\theta))\nonumber\\
&\sim - \inv{2 \pi} \; f(\vb{\xi}) \int^{\pi}_0 d\theta = - \frac{1}{2} \; f(\vb{\xi}) \;\;\;\text{when}\;\;\;\eps \ra 0,  
\end{align*}
and also
\begin{align*}
\int_{C_{\eps}} dl_{\vb{x}} \; \prt{\vb{n}} \phi(\vb{x}) \; k(\vb{x}; \vb{\xi}) &= \int^{\pi}_0 d\theta \; \eps (-\inv{2 \pi } \; \ln \eps) \;\prt{\vb{n}}  \phi(\vb{\xi} + ...)\nonumber\\
&\sim - \frac{1}{2} \; \prt{\vb{n}} \phi(\vb{\xi}) \; \eps \; \ln \eps \ra 0\;\;\;\text{when}\;\;\; \eps \ra 0.  
\end{align*}
Our boundary integral equation is then from \rf{401.12}\label{BoundaryintegralEquation}
\begin{align}
PV_{\vb{\xi}} \; \int_S dl_{\vb{x}} \; k(\vb{x}; \vb{\xi}) \; \prt{\vb{n}} \phi(\vb{x}) = b(\vb{\xi}), && \vb{\xi} \in S, \lbl{426.12} 
\end{align}
where $b(\vb{\xi})$ is a known function given by 
\begin{align*}
b(\vb{\xi}) &= - \int_V dA_{\vb{x}} \; k(\vb{x}; \vb{\xi}) \; F(\vb{x}) + \frac{1}{2} \; f(\vb{\xi})\nonumber\\ 
&+ PV_{\vb{\xi}} \; \int_S dl_{\vb{x}} \; \prt{\vb{n}} k(\vb{x}, \vb{\xi}) \; f(\vb{x}). 
\end{align*}
After we have used \rf{426.12} to calculate the unknown boundary data $\prt{\vb{n}} \phi(\vb{x})$ we can use \rf{396.12} to calculate the solution to the boundary value problem \rf{395.12} at any chosen point. \\
In \rf{426.12} there is no restriction on $S$, like it being nice and symmetric. However solving \rf{426.12} must be done numerically, I know of no closed form solutions to equation \rf{426.12}. There exists,  however,  very efficient ways to solve \rf{426.12} numerically. 

Analytic methods for calculating Green's functions form a large body of mathematics. However, beyond eigenfunction expansions, the level of generality of these methods are low. They usually only apply to special operators and/or special domains. If you ever need these methods you must dive into the research literature. One of these special methods is \ttx{the method of images}. It only works for very special geometries and mostly only for the Laplace operator. Because of the importance of these few cases, it is however important to be somewhat familiar with this method.

\subsection{Computational projects}
\subsubsection{The Helmholtz equation}

\noindent Let the operator%
\[
\mathcal{L}=\frac{d^{2}}{dx^{2}}+n^{2}(x)
\]

\noindent be given. The funtion $n(x)$ is piecewise constant and given by
\[
n(x)=\left\{
\begin{tabular}
[c]{l}%
$n_{0},$ \ \ \ \ \ \ $a<x<L$\\
$n_{1},$ \ \ \ \ \ $-a<x<a$\\
$n_{0},$ \ \ \ \ \ \ $-L<x<-a$%
\end{tabular}
\right.
\]

\begin{figure}[htbp]
\centering
\includegraphics[
height=2.001in,
width=3.8147in
]{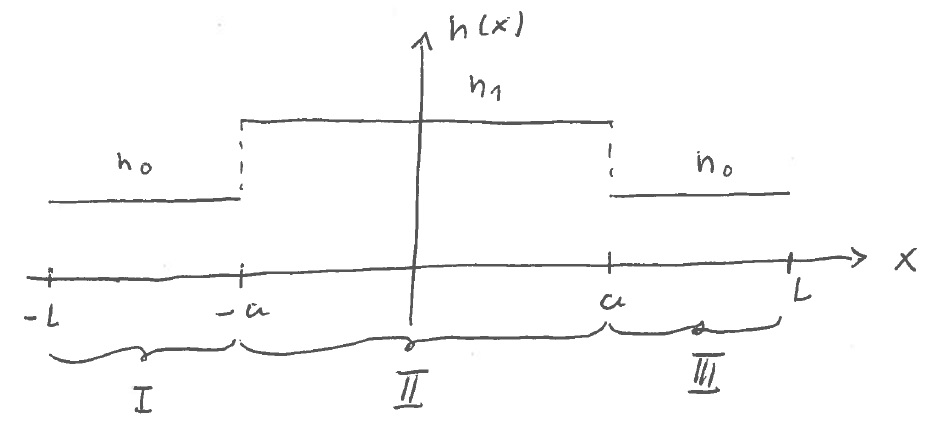}
\caption{Graph of the function n(x)}
\label{Fig1}
\end{figure}

\noindent The following boundary value problem is given%
\begin{align}
\mathcal{L}f(x)  & =h(x),\text{ \ \ \ }-L<x<L\label{BVP}\\
f(-L)  & =0\nonumber\\
f(L)  & =0\nonumber
\end{align}

\noindent where we also stipulate that $f(x)$ and $f^{\prime}(x)$ are continuous at
$x=-a$ and $x=a$.

\begin{description}
\item[a)] Solve the boundary value problem (\ref{BVP}) directly:

\begin{description}
\item[Step1] In each of regions $I,II$ and $III$ find the general solution of
the differential equation. You can find a solution to the inhomogeneous
equation using variation of parameters.

\item[Step2] After step 1 you have 6 free constants. Applying the boundary
conditions at $-L,-a,a$ and $L$ will give a system of 6 linear equations for
the 6 unknown constants.

\item[Step3] Solve the linear system. Mathematica can be very helpful here or
you can choose to do it numerically.

\item[Step4] Insert the formulas from step 3 into the functions from step 1.
These are your solution formulas for the boundary value problem.
\end{description}

\item Test your solution formulas using the artificial source test. I remind
you of the key steps in this type of test.

\begin{itemize}
\item Pick a function, $f(x)$, that is continuously differentiable in the
interval $[-L,L]$.

\item Insert this function into the left side of the differential equation and
find a formula for the function $h(x)$.

\item Use this function in the solution formulas for the boundary value problem.

\item Plot $f(x)$ and the function produced by your solution formulas in the
interval $[-L,L]$. They should overlap perfectly if your solutions formulas
are correct.
\end{itemize}

\item Pick now a particular function $h(x)$ and plot your solution of the
boundary value problem and the function $n(x)$ in the same plot. I will
suggest a Gaussian function with support in the domain denoted by $III$ in the
figure.%
\[
h(x)=e^{-\gamma(x-x_{0})^{2}}
\]

\item where $x_{0}$ is in the center of region $III$.

\item[b)] Solve the boundary value problem (\ref{BVP}) using Green's functions
to $\mathcal{L}$ satisfying particular boundary conditions:

\begin{description}
\item[Step1] Derive an integral identity involving $\mathcal{L}$ on the
interval $[-L,L]$. We have done this several times in the lectures.

\item[Step2] Apply the integral identity using a Green's function for
$\mathcal{L}$ and the solution of the boundary value problem.

\item[Step3] Pick boundary conditions for the Green's function so that the
integral identity from step 2 give an explicit solution to the boundary value
problem in terms of the Green's function and $h(x)$.

\item[Step4] Find the Green's function in the following way.

\begin{description}
\item[1.] Assume $\xi$ is in the sub domain $I$. Write down the general
solution to the homogeneous equation%
\[
K^{\prime\prime}(x;\xi)+n(x)^{2}K(x;\xi)=0
\]

\item for $x$ in each of the regions $[-L,\xi],[\xi,-a],II$ and $III$.

\item[2.] From point $1$ you have expressions containing $8$ arbitrary
functions of $\xi$. We assume that $K(x;\xi)$ is continuous and has continuous
first derivative at the points $a$ and $-a$. By definition of a Green's
function we have the two usual boundary conditions at $x=\xi$. Finally we have
the two boundary conditions at $-L$ and $L$ from step $3$. We have now $8$
boundary conditions.

\item[3.] Apply the boundary conditions from point $2$ to the solution
formulas from $1$. This give a linear system of $8$ equations for the $8$
unknown functions of $\xi$.

\item[4.] Solve the linear system from step $3$ symbolically or numerically.
Inserting the resulting formulas into the solutions from point $1$ will give
you a solution formula for the Green's function when $\xi$ is in region $I$.

\item[5.] Repeat points $1-4$ for $\xi$ in regions $II$ and $III$.
\end{description}
\end{description}

\item Pick the same $h(x)$ as in problem $a)$ and plot the solution found in
step 3 and the solution found in $a)$ together. They should overlap perfectly
if you have done everything right.

\item[c)] Solve the boundary value problem (\ref{BVP}) using boundary integral equations:

\begin{description}
\item[Step1] Find a Green's function for the operator $\mathcal{G}$ defined by%
\[
\mathcal{G=}\frac{d^{2}}{dx^{2}}+n^{2}
\]

\item where $n$ is some constant.

\item[Step2] Put $n=n_{0}$ and derive an integral identity for the operator
$\mathcal{G}$ on the interval $I$. Apply this integral identity using your
Green's function from step 1 and the solution to the boundary value problem
(\ref{BVP}). This give you a integral identity relating values of $f$ in the
interior of the interval $I$ and the unknown boundary values $f^{\prime
}(-L),f(-a)$ and $f^{\prime}(-a)$.

\item[Step3] Repeat step 2 for the regions $II$ and $III$ using $n=n_{1}$ and
$n=n_{0}$.

\item[Step4] In the integral identities from step 1 and step 2 take the limit
as $\xi$ approach lower and upper endpoints of the intervals $I,II$ amd $III$.
Taken together with the assumed continuity of $f$ and $f^{\prime} $ at the
points $-a$ and $a$ you will now have a linear system of $6$ equations for the
$6$ unknown values $f^{\prime}(-L),f(-a),f^{\prime}(-a),f(a),f^{\prime}(a)$
and $f^{\prime}(L)$.

\item[Step5] Solve the linear system from step 4. Inserting the solution into
the integral identities from step1 and step 2 give you explicit formulas for
the solution, $f(\xi),$ of the boundary value problem (\ref{BVP})~for $\xi$ in
the three intervals $I,II$ and $III$.

\item Using the same choice of $h(x)$ as in part $a)$ plot the solution from
step 5 and the solution found in $a)$ together in the same plot. They should
overlap perfectly if you have done everything right.
\end{description}

\item[d)] Solve the boundary value problem (\ref{BVP}) using finite difference methods:

\begin{description}
\item[Step1] Discretize the boundary value problem on a uniform grid%
\[
x_{j}=-L+dx\ast j,\text{ \ \ \ \ }j=0,1,...,N+1
\]

\item where $dx=\frac{2L}{N+1}$ . Use the center difference for the second derivative.

\item[Step2] Apply the boundary conditions at $-L=x_{0}$ and $L=x_{N+1}$to the
equations from step 1. You will have a linear system of $N$ equations for the
$N$ unknown functions values.%
\[
f(x_{1}),f(x_{1}),....,f(x_{N})
\]

\item[Step3] Choose some large $N$ and solve the linear system from step 2.
This you can do using Mathematica, Matlab,C......
\end{description}

\item Now plot the numerical solution found in step 3 and the solution from
$a)$ on the grid from step 1 in the same plot. They should overlap if you have
done everything right.
\end{description}

\subsubsection{The Laplace equation}

In this project we are going to solve the Laplace equation on the unit square
using several different methods. Our problem is
\begin{align}
-\nabla^{2}u(x,y)  & =\rho(x,y),\text{ \ \ \ \ }0<x<1,\text{ \ }%
0<y<1\label{Poissont}\\
u(x,y)  & =f(x,y)\text{ \ \ \ on the boundary of the unit square}\nonumber
\end{align}

\begin{description}
\item[a)] Solve the problem directly using a generalized Fourier series based
on the eigenvalue problem%
\begin{align*}
-\nabla^{2}M(x,y)  & =\lambda M(x,y),\text{ \ \ \ \ }0<x<1,\text{ \ }0<y<1\\
M(x,y)  & =0\text{ \ \ \ on the boundary of the unit square}%
\end{align*}

\item[b)] Introduce a 2D grid for the unit square, $(x_{i},y_{j})$ where
\begin{align*}
x_{i}  & =ih,\text{ \ \ }i=0,1,...,N+1\\
y_{j}  & =jk,\text{ \ \ }j=0,1,...,M+1
\end{align*}
where $M$ and $N$ are some fixed integer determining the number of points,
$MN$, in the grid and where the grid parameters $h$ and $k$ are
\begin{align*}
h  & =\frac{1}{N+1}\\
k  & =\frac{1}{M+1}%
\end{align*}
We are now going to test the solution from a) by using an artificial source.
Let $u(x,y)$ be a Gaussian function of the form%
\begin{equation}
u(x,y)=ae^{-b(x-x_{0})^{2}-c(y-y_{0})^{2}}\label{SourceSolution}%
\end{equation}
where $a,b$ and $c$ are some positive real numbers and where $(x_{0},y_{0})$
is a point inside the unit square. For a given choice of $a,b,c,x_{0}$ and
$y_{0}$ define a function $g$ on the boundary of the unit square by%
\[
g(x,y)=u(x,y)
\]
and a function $G$ inside the unit square by
\[
G(x,y)=-u_{xx}-u_{yy}
\]
Use the Fourier solution from a) with $\rho=G$ and $f=g$ to calculate the
solution to Poisson's equation at the grid points $(x_{i},y_{j})$. Plot the
numerically computed solution and the exact solution (\ref{SourceSolution}) in
the same figure.Show both a 3D plot and selected 2D slices that in a
convincing way show that the numerical solution and exact solution
(\ref{SourceSolution}) are the same. Show pictures for several choices of the
parameters $a,b,c,x_{0}$ and $y_{0}$.

\item[c)] Find a formula for the Green's function satisfying
\begin{align*}
-\nabla^{2}K(\mathbf{x};\mathbf{\xi})  & =\delta(\mathbf{x}-\mathbf{\xi
),}\text{ \ \ \ }\mathbf{x},\text{in the unit square}\\
K(\mathbf{x};\mathbf{\xi})  & =0,\text{ \ \ \ \ \ \ \ }\mathbf{x}\text{ on the
boundary of the unit square}%
\end{align*}
using a finite Fourier transform based on the eigenvalue problem
\begin{align*}
-M^{\prime\prime}(x)  & =\lambda M(x),\text{ \ }0<x<1\\
M(0)  & =M(1)=0
\end{align*}
I have done this derivation on page \pageref{SimpleBoundaryProblem} in these lecture notes,  but I
want you do redo it and include all details in your project report.

\item[d)] Use the Green's function from c) and the general integral identity
for the Laplace operator to find a formula for the solution of the boundary
value problem%
\begin{align}
-\nabla^{2}u(x,y)  & =0,\text{ \ \ \ \ }0<x<1,\text{ \ }0<y<1\label{laplace}\\
u(x,y)  & =f(x,y)\text{ \ \ \ on the boundary of the unit square}\nonumber
\end{align}

\item[e)] Pick a function $f$. You can for example let $f$ be a Gaussian of
the type (\ref{SourceSolution}) restricted to the boundary. \ Evaluate the
formula from d) on the grid from b) and show by using 3D and 2D plots that
your solution from d) and the solution calculated using the formula from a)
are the same.

\item[f)] Solve the equation (\ref{laplace}) from d) by using the boundary
integral equation from page \pageref{BoundaryintegralEquation} in these lecture notes. After you
have calculated the unknown boundary data using the boundary integral equation
you use the general integral identity to calculate the solution on the grid
from b). As boundary data you can use the function from b) or something of
your own choice. Compare the boundary integral solution with the solution
calculated in a) by making 3D and 2D plots.

Setting up and solving the boundary integral equations require you to go
through the following steps.

\begin{enumerate}
\item Parametrize the boundary of the square. The boundary consists of four
straight lines.%
\[
C=C_{1}\cup C_{2}\cup C_{3}\cup C_{4}
\]
A simple parametrization for the line $C_{1}$%
\[
C_{1}=\{(x,0);\text{ }0<x<1\}
\]

is for example
\[
\mathbf{\gamma}_{1}(t)=(t,0),\text{ \ \ \ \ \ }0<t<1
\]

and similar parameterizations, $\gamma_{k}$, for the other pieces of the
boundary. Your integral equation will now have the general structure%
\[%
{\displaystyle\sum_{l=1}^{4}}
\int_{0}^{1}dt^{\prime}A^{kl}(t,t^{\prime})v^{l}(t^{\prime})=B^{k}(t),\text{
\ \ }k=1,2..,4
\]

where by definition%
\begin{align*}
v^{l}(t^{\prime})  & =\partial_{n(\mathbf{\gamma}_{l}(t^{\prime}))}%
\varphi(\mathbf{\gamma}_{l}(t^{\prime}))\\
A^{kl}(t,t^{\prime})  & =K(\mathbf{\gamma}_{l}(t^{\prime});\mathbf{\gamma}%
_{k}(t))\\
B^{k}(t)  & =b(\mathbf{\gamma}_{k}(t))
\end{align*}

\item Discretize the boundary of the unit square using the parametrization
from step 1: Introduce intervals%
\[
I_{i}=(\alpha_{i},\alpha_{i-1})
\]

where
\[
\alpha_{i}=ih,\text{ \ \ \ }i=0,1,...,N+1
\]

and where the grid parameters $h$ is
\[
h=\frac{1}{N+1}
\]

Let $s_{i}$ be the midpoint of the interval $I_{i}$%
\[
s_{i}=\frac{\alpha_{i}+\alpha_{i-1}}{2}=(i-\frac{1}{2})h,\text{ \ \ \ }%
i=1,2,...,N
\]

Using this discretization your integral equation will have the general form%
\[%
{\displaystyle\sum_{l=1}^{4}}
{\displaystyle\sum_{j=1}^{N}}
\int_{I_{j}}dt^{\prime}A^{kl}(t,t^{\prime})v^{l}(t^{\prime})=B^{k}(t),\text{
\ \ }k=1,2..,4
\]

which we approximate by the linear algebraic system of equations
\[%
{\displaystyle\sum_{l=1}^{4}}
{\displaystyle\sum_{j=1}^{N}}
A_{ij}^{kl}v_{j}^{l}=B_{i}^{k},\text{ \ \ }k=1,2..,4
\]

where by definition%
\begin{align*}
A_{ij}^{kl}  & =\int_{I_{j}}dt^{\prime}A^{kl}(s_{i},t^{\prime})\\
v_{j}^{l}  & =v^{l}(s_{j})\\
B_{i}^{k}  & =B^{k}(s_{i})
\end{align*}

\item Observe that the integrand in the integrals defining the matrix elements
$A_{ij}^{kl}$ has a singularity in the domain of integration only when
$(k,i)=(l,j)$. Otherwise the integrands are smooth functions. Inspired by this
we approximate the coefficients $A_{ij}^{kl}$ for $(k,i)\neq(l,j)$ in the
following way%
\[
A_{ij}^{kl}=hA^{kl}(s_{i},s_{j})
\]

We are thus using the midpoint rule for evaluating the integrals. The matrix
elements $A_{ii}^{kk}$ are evaluated by using principle value integrals. The
coefficients $B_{i}^{k}$ are evaluated using the same approximations.

\item You now have a system of $4N$ equations for the $4N$ unknowns $v_{l}%
^{k}$. Solve this linear system using a linear system solver. (You don't have
to write your own)

\item Discretize the integral formula that express the solution inside the
unit square in terms of the boundary data in the same way as in points 2 and 3
and use it to calculate the solution of equation (\ref{laplace}) from d) on
the grid from b). Compare your boundary integral solution with the solution
from a) by using 3D and 2D plots.
\end{enumerate}

\item[g)] Solve the equation (\ref{laplace}) from d) by using finite
differences on the grid from b).  As boundary data you
can use the function from b) or something else of your choice. Compare your
finite difference solution with the boundary integral solution from f) by
using 3D and 2D plots.
\end{description}

\section{Acknowledgment}
The master level class in applied mathematics that these lecture notes has been written  for, has had a long history at the department of Mathematics and Statistics at the Arctic University of Norway. Its roots stretch back all the way to the 1970s. The architect for the first version of this class was  Einar Mj\o lhus, which at the time was associate professor of applied mathematics at the university. Several written accounts of the topics covered in the various iterations of the class appeared over the years, most of them now  lost to history, but all of them have influenced the author while preparing these lecture notes. The author has also been inspired many of the text books and lecture notes written by other authors on the topics that form the five sections of these lecture notes. Some of these written accounts has not been cited directly in the text but all sources that significantly influenced the author of the current lecture notes can be found in the list of references at the end og this document. 

  Over the years, the nature of the topics, and their number, fluctuated somewhat, with contributions from Kristian B. Dysthe and Tor Fl\aa, both at the time professors in applied mathematics at the university. Overall,  mathematical topics included in the class has been remarkably stable over all the years the class has been running. This said,however,  some new topics has been included in the current set of  lecture notes that were not present in the original accounts of the class, and many of the topics that were part of the original accounts, has been greatly expanded upon. 

\section{Appendix A}
\setcounter{equation}{0}
\subsection{ The multiple scale method for Maxwell's equations}

In optics the equations of interest are of course Maxwell's equations. For a
situation without free charges and currents they are given by
\begin{align}
\partial_{t}\mathbf{B+\nabla\times E}  & =0,\nonumber\\
\partial_{t}\mathbf{D-\nabla\times H}  & =0,\nonumber\\
\nabla\cdot\mathbf{D}  & =0,\nonumber\\
\nabla\cdot\mathbf{B}  & =0.\label{Eq8.1}
\end{align}
At optical frequencies materials of interest are almost always nonmagnetic. It is thus  appropriate to assume that
\begin{align}
\mathbf{H}  & =\frac{1}{\mu}\mathbf{B},\nonumber\\
\mathbf{D}  & =\varepsilon_{0}\mathbf{E}+\mathbf{P}.\label{Eq8.2}
\end{align}
The polarization is in general a sum of a term that is linear in $\mathbf{E}
$, and one that is nonlinear in $\mathbf{E}$. We have
\begin{equation}
\mathbf{P}=\mathbf{P}_{L}+\mathbf{P}_{NL},\label{Eq8.3}
\end{equation}
where the term linear in $\mathbf{E}$ has the general form
\begin{equation}
\mathbf{P}_{L}(\mathbf{x},t)=\varepsilon_{0}\int_{-\infty}^{t}dt^{\prime}
\chi(t-t^{\prime})\mathbf{E}(\mathbf{x},t^{\prime}).\label{Eq8.4}
\end{equation}
Thus the polarization at a time $t$ depends on the electric field at all times
previous to $t$. This memory effect is what we in optics call \textit{temporal
dispersion}. The presence of dispersion in Maxwell equations spells trouble
for the integration of the equations in time; we can not solve them as a
standard initial value problem. This is of course well known in optics and
various, more or less ingenious, methods have been designed for getting around
this problem. In optical pulse propagation, one gets around the problem by
solving Maxwell's equations as a boundary value problem rather
than as an initial value problem. A very general version of this approach is
the well known UPPE \cite{UPPE}\cite{per1} propagation scheme. In
these lecture notes we will, using the multiple scale method, derive
an approximation to Maxwell's equations that can be solved as an initial value problem.

In the explicit calculations that we do, we will assume that the nonlinear
polarization is generated by the Kerr effect. Thus we will assume that
\begin{equation}
\mathbf{P}_{NL}=\varepsilon_{0}\eta\mathbf{E}\cdot\mathbf{EE},\label{Eq8.5}
\end{equation}
where $\eta$ is the Kerr coefficient. This is a choice we make just to be
specific, the applicability of the multiple scale method to Maxwell's
equations in no way depend on this particular choice for the nonlinear response.

Before we proceed with the multiple scale method we will introduce a more
convenient representation of the dispersion. Observe that we have
\begin{align}
\mathbf{P}_{L}(\mathbf{x},t)  & =\varepsilon_{0}\int_{-\infty}^{t}dt^{\prime
}\chi(t-t^{\prime})\mathbf{E}(\mathbf{x},t^{\prime})\nonumber\\
& =\varepsilon_{0}\int_{-\infty}^{\infty}d\omega\widehat{\chi}(\omega
)\widehat{\mathbf{E}}(\mathbf{x},\omega)e^{-i\omega t}\nonumber\\
& =\varepsilon_{0}\int_{-\infty}^{\infty}d\omega\left(
{\displaystyle\sum_{n=0}^{\infty}}
\frac{\widehat{\chi}^{(n)}(0)}{n!}\omega^{n}\right)  \widehat{\mathbf{E}
}(\mathbf{x},\omega)e^{-i\omega t}\nonumber\\
& =\varepsilon_{0}
{\displaystyle\sum_{n=0}^{\infty}}
\frac{\widehat{\chi}^{(n)}(0)}{n!}\left(  \int_{-\infty}^{\infty}d\omega
\omega^{n}\widehat{\mathbf{E}}(\mathbf{x},\omega)e^{-i\omega t}\right)
\nonumber\\
& =\varepsilon_{0}
{\displaystyle\sum_{n=0}^{\infty}}
\frac{\widehat{\chi}^{(n)}(0)}{n!}\left(  \int_{-\infty}^{\infty}
d\omega(i\partial_{t})^{n}\widehat{\mathbf{E}}(\mathbf{x},\omega)e^{-i\omega
t}\right) \nonumber\\
& =\varepsilon_{0}
{\displaystyle\sum_{n=0}^{\infty}}
\frac{\widehat{\chi}^{(n)}(0)}{n!}(i\partial_{t})^{n}\left(  \int_{-\infty
}^{\infty}d\omega\widehat{\mathbf{E}}(\mathbf{x},\omega)e^{-i\omega t}\right)
\nonumber\\
& =\widehat{\chi}(i\partial_{t})\mathbf{E}(\mathbf{x},t),\nonumber
\end{align}
where $\widehat{\chi}(\omega)$ is the Fourier transform of $\chi(t)$. These
manipulations are of course purely formal; in order to make them into honest
mathematics we must dive into the theory of \textit{pseudo differential
operators}. In these lecture notes we will not do this as our focus is on
mathematical methods rather than mathematical theory.

Inserting (\ref{Eq8.2}),(\ref{Eq8.3}),(\ref{Eq8.4}) and (\ref{Eq8.5}) into
(\ref{Eq8.1}), we get Maxwell's equations in the form
\begin{align}
\partial_{t}\mathbf{B+\nabla\times E}  & =0,\nonumber\\
\partial_{t}\mathbf{E-c}^{2}\mathbf{\nabla\times B+\partial}_{t}\widehat{\chi
}(i\partial_{t})\mathbf{E}  & =-c^{2}\mu_{0}\partial_{t}\mathbf{P}
_{NL},\nonumber\\
\nabla\cdot\left(  \mathbf{E}+\widehat{\chi}(i\partial_{t})\mathbf{E}\right)
& =-\frac{1}{\varepsilon_{0}}\nabla\cdot\mathbf{P}_{NL},\nonumber\\
\nabla\cdot\mathbf{B}  & =0.\label{Eq8.7}
\end{align}

\subsubsection{TE scalar wave packets}

Let us first simplify the problem by only considering solutions of the form
\begin{align}
\mathbf{E}(x,y,z,t)  & =E(x,z,t)\mathbf{e}_{y},\nonumber\\
\mathbf{B}(x,y,z,t)  & =B_{1}(x,z,t)\mathbf{e}_{x}\mathbf{+}B_{2}
(x,z,t)\mathbf{e}_{z}.\label{Eq8.8}
\end{align}
For this simplified case, Maxwell's equations takes the form
\begin{align}
\partial_{t}B_{1}-\partial_{z}E  & =0,\nonumber\\
\partial_{t}B_{2}+\partial_{x}E  & =0,\nonumber\\
\partial_{t}E-c^{2}(\partial_{z}B_{1}-\partial_{x}B_{2})+\partial
_{t}\widehat{\chi}(i\partial_{t})E  & =-\partial_{t}P_{NL},\nonumber\\
\partial_{x}B_{1}+\partial_{z}B_{2}  & =0,\label{Eq8.9}
\end{align}
where
\begin{equation}
P_{NL}=\eta E^{3}.\label{Eq8.10}
\end{equation}
It is well known that this vector system is fully equivalent to the following
scalar equation
\begin{equation}
\partial_{tt}E-c^{2}\nabla^{2}E+\partial_{tt}\widehat{\chi}(i\partial
_{t})E=-\partial_{tt}P_{NL},\label{Eq8.11}
\end{equation}
where we have introduced the operator
\begin{equation}
\nabla^{2}=\partial_{xx}+\partial_{zz}.\label{Eq8.12}
\end{equation}
Equation (\ref{Eq8.11}) will be the staring point for our multiple scale
approach, but before that I will introduce the notion of a \textit{formal}
perturbation parameter. For some particular application of equation
(\ref{Eq8.11}) we will usually start by making the equation dimension-less by
picking some scales for space, time, and $E$ relevant for our particular
application. Here we don't want to tie our calculations to some particular
choice of scales and introduce therefore a formal perturbation parameter in
the equation multiplying the nonlinear polarization term. Thus we have
\begin{equation}
\partial_{tt}E-c^{2}\nabla^{2}E+\partial_{tt}\widehat{\chi}(i\partial
_{t})E=-\varepsilon^{2}\eta\partial_{tt}E^{3}.\label{Eq8.13}
\end{equation}
Starting with this equation, we will proceed with our perturbation calculations
assuming that $\varepsilon<<1$ and in the end we will remove $\varepsilon$ by
setting it equal to $1$. What is going on here is that $\varepsilon$ is a
"place holder" for the actual small parameter that will appear in front of the
nonlinear term in the equation when we make a particular choice of scales.
Using such formal perturbation parameters is very common.

You might ask why I use $\varepsilon^{2}$ instead of $\varepsilon$ as formal
perturbation parameter? I will not answer this question here but will say
something about it at the very end of the lecture notes. We proceed with the
multiple scale method by introducing the expansions
\begin{align}
\partial_{t}  & =\partial_{t_{0}}+\varepsilon\partial_{t_{1}}+\varepsilon
^{2}\partial_{t_{2}}+...\;\;,\nonumber\\
\nabla & =\nabla_{0}+\varepsilon\nabla_{1}+\varepsilon^{2}\nabla
_{2}+...\;\;,\nonumber\\
e  & =e_{0}+\varepsilon e_{1}+\varepsilon^{2}e_{2}+...\;\;,\nonumber\\
E(\mathbf{x},t)  & =e(\mathbf{x}_{0},t_{0},\mathbf{x}_{1},t_{1},...)|_{t_{j}
=\varepsilon^{j}t,\mathbf{x}_{j}=\varepsilon^{j}\mathbf{x}},\label{Eq8.14}
\end{align}
where
\begin{equation}
\nabla_{j}=(\partial_{x_{j}},\partial_{z_{j}}),\label{Eq8.15}
\end{equation}
is the gradient with respect to $\mathbf{x}_{j}=(x_{j},z_{j})$. We now insert
(\ref{Eq8.14}) into (\ref{Eq8.13}) and expand everything in sight
\begin{align*}
& (\partial_{t_{0}}+\varepsilon\partial_{t_{1}}+\varepsilon^{2}\partial
_{t_{2}}+...)(\partial_{t_{0}}+\varepsilon\partial_{t_{1}}+\varepsilon
^{2}\partial_{t_{2}}+...)\\
& (e_{0}+\varepsilon e_{1}+\varepsilon^{2}e_{2}+...)-\\
& c^{2}(\nabla_{0}+\varepsilon\nabla_{1}+\varepsilon^{2}\nabla_{2}
+...)\cdot(\nabla_{0}+\varepsilon\nabla_{1}+\varepsilon^{2}\nabla_{2}+...)\\
& (e_{0}+\varepsilon e_{1}+\varepsilon^{2}e_{2}+...)+\\
& (\partial_{t_{0}}+\varepsilon\partial_{t_{1}}+\varepsilon^{2}\partial
_{t_{2}}+...)(\partial_{t_{0}}+\varepsilon\partial_{t_{1}}+\varepsilon
^{2}\partial_{t_{2}}+...)\\
& \widehat{\chi}(i\partial_{t_{0}}+i\varepsilon\partial_{t_{1}}+i\varepsilon
^{2}\partial_{t_{2}}+...)(e_{0}+\varepsilon e_{1}+\varepsilon^{2}e_{2}+...)\\
& =-\varepsilon^{2}\eta(\partial_{t_{0}}+\varepsilon\partial_{t_{1}
}+\varepsilon^{2}\partial_{t_{2}}+...)(\partial_{t_{0}}+\varepsilon
\partial_{t_{1}}+\varepsilon^{2}\partial_{t_{2}}+...)\\
& (e_{0}+\varepsilon e_{1}+\varepsilon^{2}e_{2}+...)^{3}\\
& \Downarrow\\
& (\partial_{t_{0}t_{0}}+\varepsilon(\partial_{t_{0}t_{1}}+\partial
_{t_{1}t_{0}})+\varepsilon^{2}(\partial_{t_{0}t_{2}}+\partial_{t_{1}t_{1}
}+\partial_{t_{2}t_{0}})+...)\\
& (e_{0}+\varepsilon e_{1}+\varepsilon^{2}e_{2}+...)-\\
& c^{2}(\nabla_{0}^{2}+\varepsilon(\nabla_{1}\cdot\nabla_{0}+\nabla_{0}
\cdot\nabla_{1})+\varepsilon^{2}(\nabla_{2}\cdot\nabla_{0}+\nabla_{1}
\cdot\nabla_{1}+\nabla_{0}\cdot\nabla_{2})+...)\\
& (e_{0}+\varepsilon e_{1}+\varepsilon^{2}e_{2}+...)+\\
& (\partial_{t_{0}t_{0}}+\varepsilon(\partial_{t_{0}t_{1}}+\partial
_{t_{1}t_{0}})+\varepsilon^{2}(\partial_{t_{0}t_{2}}+\partial_{t_{1}t_{1}
}+\partial_{t_{2}t_{0}})+...)\\
& (\widehat{\chi}(i\partial_{t_{0}})+\varepsilon\widehat{\chi}^{\prime
}(i\partial_{t_{0}})i\partial_{t_{1}}+\varepsilon^{2}(\widehat{\chi}^{\prime
}(i\partial_{t_{0}})i\partial_{t_{2}}-\frac{1}{2}\widehat{\chi}^{\prime\prime
}(i\partial_{t_{0}})\partial_{t_{1}t_{1}})+...)\\
& (e_{0}+\varepsilon e_{1}+\varepsilon^{2}e_{2}+...)\\
& =-\varepsilon^{2}\partial_{t_{0}t_{0}}e_{0}^{3}+...\;\;,\\
\nonumber\\
& \Downarrow
\end{align*}

\begin{align*}
& \partial_{t_{0}t_{0}}e_{0}+\varepsilon(\partial_{t_{0}t_{0}}e_{1}
+\partial_{t_{0}t_{1}}e_{0}+\partial_{t_{1}t_{0}}e_{0})\\
& +\varepsilon^{2}(\partial_{t_{0}t_{0}}e_{2}+\partial_{t_{0}t_{1}}
e_{1}+\partial_{t_{1}t_{0}}e_{1}+\partial_{t_{0}t_{2}}e_{0}+\partial
_{t_{1}t_{1}}e_{0}+\partial_{t_{2}t_{0}}e_{0})+...\\
& -c^{2}\nabla_{0}^{2}e_{0}-\varepsilon c^{2}(\nabla_{0}^{2}e_{1}+\nabla
_{1}\cdot\nabla_{0}e_{0}+\nabla_{0}\cdot\nabla_{1}e_{0})\\
& -\varepsilon^{2}c^{2}(\nabla_{0}^{2}e_{2}+\nabla_{1}\cdot\nabla_{0}
e_{1}+\nabla_{0}\cdot\nabla_{1}e_{1}\\
& +\nabla_{2}\cdot\nabla_{0}e_{0}+\nabla_{1}\cdot\nabla_{1}e_{0}+\nabla
_{0}\cdot\nabla_{2}e_{0})+...\\
& +\partial_{t_{0}t_{0}}\widehat{\chi}(i\partial_{t_{0}})e_{0}+\varepsilon
(\partial_{t_{0}t_{0}}\widehat{\chi}(i\partial_{t_{0}})e_{1}+\partial
_{t_{0}t_{0}}\widehat{\chi}^{\prime}(i\partial_{t_{0}})i\partial_{t_{1}}
e_{0}\\
& +\partial_{t_{0}t_{1}}\widehat{\chi}(i\partial_{t_{0}})e_{0}+\partial
_{t_{1}t_{0}}\widehat{\chi}(i\partial_{t_{0}})e_{0})+\varepsilon^{2}
(\partial_{t_{0}t_{0}}\widehat{\chi}(i\partial_{t_{0}})e_{2}\\
& +\partial_{t_{0}t_{0}}\widehat{\chi}^{\prime}(i\partial_{t_{0}}
)i\partial_{t_{1}}e_{1}+\partial_{t_{0}t_{1}}\widehat{\chi}(i\partial_{t_{0}
})e_{1}+\partial_{t_{1}t_{0}}\widehat{\chi}(i\partial_{t_{0}})e_{1}\\
& +\partial_{t_{0}t_{0}}\widehat{\chi}^{\prime}(i\partial_{t_{0}}
)i\partial_{t_{2}}e_{0}-\frac{1}{2}\partial_{t_{0}t_{0}}\widehat{\chi}
^{\prime\prime}(i\partial_{t_{0}})\partial_{t_{1}t_{1}}e_{0}+\partial
_{t_{1}t_{0}}\widehat{\chi}^{\prime}(i\partial_{t_{0}})i\partial_{t_{1}}
e_{0}\\
& +\partial_{t_{0}t_{1}}\widehat{\chi}^{\prime}(i\partial_{t_{0}}
)i\partial_{t_{1}}e_{0}+\partial_{t_{2}t_{0}}\widehat{\chi}(i\partial_{t_{0}
})e_{0}+\partial_{t_{1}t_{1}}\widehat{\chi}(i\partial_{t_{0}})e_{0}\\
& +\partial_{t_{0}t_{2}}\widehat{\chi}(i\partial_{t_{0}})e_{0})+...\\
& =-\varepsilon^{2}\partial_{t_{0}t_{0}}e_{0}^{3}+...\;\;,
\end{align*}
which gives us the perturbation hierarchy
\begin{gather}
\partial_{t_{0}t_{0}}e_{0}-c^{2}\nabla_{0}^{2}e_{0}+\partial_{t_{0}t_{0}
}\widehat{\chi}(i\partial_{t_{0}})e_{0}=0,\label{Eq8.16}\\
\nonumber\\
\partial_{t_{0}t_{0}}e_{1}-c^{2}\nabla_{0}^{2}e_{1}+\partial_{t_{0}t_{0}
}\widehat{\chi}(i\partial_{t_{0}})e_{1}=\nonumber\\
-\partial_{t_{0}t_{1}}e_{0}-\partial_{t_{1}t_{0}}e_{0}-c^{2}\nabla_{1}
\cdot\nabla_{0}e_{0}-c^{2}\nabla_{0}\cdot\nabla_{1}e_{0}\nonumber\\
-\partial_{t_{0}t_{0}}\widehat{\chi}^{\prime}(i\partial_{t_{0}})i\partial
_{t_{1}}e_{0}-\partial_{t_{0}t_{1}}\widehat{\chi}(i\partial_{t_{0}}
)e_{0}-\partial_{t_{1}t_{0}}\widehat{\chi}(i\partial_{t_{0}})e_{0},\label{Eq8.17}\\
\nonumber\\
\partial_{t_{0}t_{0}}e_{2}-c^{2}\nabla_{0}^{2}e_{2}+\partial_{t_{0}t_{0}
}\widehat{\chi}(i\partial_{t_{0}})e_{2}=\nonumber\\
-\partial_{t_{0}t_{1}}e_{1}-\partial_{t_{1}t_{0}}e_{1}-\partial_{t_{0}t_{2}
}e_{0}-\partial_{t_{1}t_{1}}e_{0}-\partial_{t_{2}t_{0}}e_{0}\nonumber\\
-c^{2}\nabla_{1}\cdot\nabla_{0}e_{1}-c^{2}\nabla_{0}\cdot\nabla_{1}e_{1}
-c^{2}\nabla_{2}\cdot\nabla_{0}e_{0}-c^{2}\nabla_{1}\cdot\nabla_{1}
e_{0}\nonumber\\
-c^{2}\nabla_{0}\cdot\nabla_{2}e_{0}-\partial_{t_{0}t_{0}}\widehat{\chi
}^{\prime}(i\partial_{t_{0}})i\partial_{t_{1}}e_{1}-\partial_{t_{0}t_{1}
}\widehat{\chi}(i\partial_{t_{0}})e_{1}\nonumber\\
-\partial_{t_{1}t_{0}}\widehat{\chi}(i\partial_{t_{0}})e_{1}-\partial
_{t_{0}t_{0}}\widehat{\chi}^{\prime}(i\partial_{t_{0}})i\partial_{t_{2}}
e_{0}+\frac{1}{2}\partial_{t_{0}t_{0}}\widehat{\chi}^{\prime\prime}
(i\partial_{t_{0}})\partial_{t_{1}t_{1}}e_{0}\nonumber\\
-\partial_{t_{1}t_{0}}\widehat{\chi}^{\prime}(i\partial_{t_{0}})i\partial
_{t_{1}}e_{0}-\partial_{t_{0}t_{1}}\widehat{\chi}^{\prime}(i\partial_{t_{0}
})i\partial_{t_{1}}e_{0}-\partial_{t_{2}t_{0}}\widehat{\chi}(i\partial_{t_{0}
})e_{0}\nonumber\\
-\partial_{t_{1}t_{1}}\widehat{\chi}(i\partial_{t_{0}})e_{0}-\partial
_{t_{0}t_{2}}\widehat{\chi}(i\partial_{t_{0}})e_{0}-\partial_{t_{0}t_{0}}
e_{0}^{3}.\label{Eq8.18}
\end{gather}
For the order $\varepsilon^{0}$ equation we choose the wave packet solution
\begin{equation}
e_{0}(\mathbf{x}_{0},t_{0},\mathbf{x}_{1},t_{1},..)=A_{0}(\mathbf{x}_{1}
,t_{1},...)e^{i\theta_{0}}+(\ast),\label{Eq8.19}
\end{equation}
where
\begin{align}
\mathbf{x}_{j}  & =(x_{j},z_{j}),\nonumber\\
\theta_{0}  & =\mathbf{k}\cdot\mathbf{x}_{0}-\omega t_{0},\label{Eq8.20}
\end{align}
and where $\mathbf{k}$ is a plane vector with components $\mathbf{k}=(\xi
,\eta)$. In (\ref{Eq8.20}), $\omega$, is a function of $k=||\mathbf{k}||$ that
satisfy the dispersion relation
\begin{equation}
\omega^{2}n^{2}(\omega)=c^{2}k^{2},\label{Eq8.21}
\end{equation}
where the refractive index, $n(\omega)$, is defined by
\begin{equation}
n^{2}(\omega)=1+\widehat{\chi}(\omega).\label{Eq8.22}
\end{equation}
We now must now calculate the right-hand side of the order $\varepsilon$ equation.

Observe that
\begin{align}
\partial_{t_{1}t_{0}}e_{0}  & =-i\omega\partial_{t_{1}}A_{0}e^{i\theta_{0}
}+(\ast),\nonumber\\
\partial_{t_{0}t_{1}}e_{0}  & =-i\omega\partial_{t_{1}}A_{0}e^{i\theta_{0}
}+(\ast),\nonumber\\
\nabla_{1}\cdot\nabla_{0}e_{0}  & =ik\nabla_{1}A_{0}\cdot\mathbf{u}
e^{i\theta_{0}}+(\ast),\nonumber\\
\nabla_{0}\cdot\nabla_{1}e_{0}  & =ik\nabla_{1}A_{0}\cdot\mathbf{u}
e^{i\theta_{0}}+(\ast),\nonumber\\
\partial_{t_{0}t_{0}}\widehat{\chi}^{\prime}(i\partial_{t_{0}})i\partial
_{t_{1}}e_{0}  & =-i\omega\widehat{\chi}^{\prime}(\omega)\partial_{t_{1}}
A_{0}e^{i\theta_{0}}+(\ast),\nonumber\\
\partial_{t_{0}t_{1}}\widehat{\chi}(i\partial_{t_{0}})e_{0}  & =-i\omega
\widehat{\chi}(\omega)\partial_{t_{1}}A_{0}e^{i\theta_{0}}+(\ast),\nonumber\\
\partial_{t_{1}t_{0}}\widehat{\chi}(i\partial_{t_{0}})e_{0}  & =-i\omega
\widehat{\chi}(\omega)\partial_{t_{1}}A_{0}e^{i\theta_{0}}+(\ast),\label{Eq8.23}
\end{align}
where $\mathbf{u}$ is a unit vector in the direction of $\mathbf{k}$.
Inserting (\ref{Eq8.23}) into (\ref{Eq8.17}) we get
\begin{gather}
\partial_{t_{0}t_{0}}e_{1}-c^{2}\nabla_{0}^{2}e_{1}+\partial_{t_{0}t_{0}
}\widehat{\chi}(i\partial_{t_{0}})e_{1}=\nonumber\\
-\{-2i\omega\partial_{t_{1}}A_{0}-2ic^{2}k\mathbf{u}\cdot\nabla_{1}
A_{0}\nonumber\\
-i\omega^{2}\widehat{\chi}^{\prime}(\omega)\partial_{t_{1}}A_{0}
-2i\omega\widehat{\chi}(\omega)\partial_{t_{1}}A_{0}\}e^{i\theta_{0}}
+(\ast).\label{Eq8.24}
\end{gather}
In order to remove secular terms we must postulate that
\begin{gather}
-2i\omega\partial_{t_{1}}A_{0}-2ic^{2}k\mathbf{u}\cdot\nabla_{1}A_{0}
-i\omega^{2}\widehat{\chi}^{\prime}(\omega)\partial_{t_{1}}A_{0}
-2i\omega\widehat{\chi}(\omega)\partial_{t_{1}}A_{0}=0,\nonumber\\
\Updownarrow\nonumber\\
\omega(2n^{2}+\omega\widehat{\chi}^{\prime}(\omega))\partial_{t_{1}}
A_{0}-2ic^{2}k\mathbf{u}\cdot\nabla_{1}A_{0}=0.\label{Eq8.25}
\end{gather}
Observe that from the dispersion relation (\ref{Eq8.21}) we have
\begin{align*}
\omega^{2}n^{2}(\omega)  & =c^{2}k^{2},\\
& \Updownarrow\\
\omega^{2}(1+\widehat{\chi}(\omega))  & =c^{2}k^{2},\\
& \Downarrow\\
2\omega\omega^{\prime}n^{2}(\omega)+\omega^{2}\widehat{\chi}^{\prime}
(\omega)\omega^{\prime}  & =2c^{2}k,\\
& \Downarrow\\
\omega(2n^{2}+\omega\widehat{\chi}^{\prime}(\omega))\omega^{\prime}  &
=2c^{2}k.
\end{align*}
Thus (\ref{Eq8.25}) can be written in the form
\begin{equation}
\partial_{t_{1}}A_{0}+\mathbf{v}_{g}\cdot\nabla_{1}A_{0}=0,\label{Eq8.26}
\end{equation}
where $\mathbf{v}_{g}$ is the group velocity
\begin{equation}
\mathbf{v}_{g}=\omega^{\prime}(k)\mathbf{u}.\label{Eq8.27}
\end{equation}
The order $\varepsilon$ equation simplifies into
\begin{equation}
\partial_{t_{0}t_{0}}e_{1}-c^{2}\nabla_{0}^{2}e_{1}+\partial_{t_{0}t_{0}
}\widehat{\chi}(i\partial_{t_{0}})e_{1}=0.\label{Eq8.29}
\end{equation}
According to the rules of the game we choose the special solution
\begin{equation}
e_{1}=0,\label{Eq8.30}
\end{equation}
for (\ref{Eq8.29}). We now must compute the right-hand side of the order
$\varepsilon^{2}$ equation. Observe that

\begin{align}
\partial_{t_{2}t_{0}}e_{0}  & =-i\omega\partial_{t_{2}}A_{0}e^{i\theta_{0}
}+(\ast),\nonumber\\
\partial_{t_{1}t_{1}}e_{0}  & =\partial_{t_{1}t_{1}}A_{0}e^{i\theta_{0}}
+(\ast),\nonumber\\
\partial_{t_{0}t_{1}}e_{0}  & =-i\omega\partial_{t_{2}}A_{0}e^{i\theta_{0}
}+(\ast),\nonumber\\
\nabla_{2}\cdot\nabla_{0}e_{0}  & =ik\mathbf{u}\cdot\nabla_{2}A_{0}
e^{i\theta_{0}}+(\ast),\nonumber\\
\nabla_{1}\cdot\nabla_{1}e_{0}  & =\nabla_{1}^{2}A_{0}e^{i\theta_{0}}
+(\ast),\nonumber\\
\nabla_{0}\cdot\nabla_{2}e_{0}  & =ik\mathbf{u}\cdot\nabla_{2}A_{0}
e^{i\theta_{0}}+(\ast),\nonumber\\
\partial_{t_{0}t_{0}}\widehat{\chi}^{\prime}(i\partial_{t_{0}})i\partial
_{t_{2}}e_{0}  & =-i\omega^{2}\widehat{\chi}^{\prime}(\omega)\partial_{t_{2}
}A_{0}e^{i\theta_{0}}+(\ast),\nonumber\\
\frac{1}{2}\partial_{t_{0}t_{0}}\widehat{\chi}^{\prime\prime}(i\partial
_{t_{0}})\partial_{t_{1}t_{1}}e_{0}  & =-\frac{1}{2}\omega^{2}\widehat{\chi
}^{\prime\prime}(\omega)\partial_{t_{1}t_{1}}A_{0}e^{i\theta_{0}}
+(\ast),\nonumber\\
\partial_{t_{1}t_{0}}\widehat{\chi}^{\prime}(i\partial_{t_{0}})i\partial
_{t_{1}}e_{0}  & =\omega\widehat{\chi}^{\prime}(\omega)\partial_{t_{1}t_{1}
}A_{0}e^{i\theta_{0}}+(\ast),\nonumber\\
\partial_{t_{0}t_{1}}\widehat{\chi}^{\prime}(i\partial_{t_{0}})i\partial
_{t_{1}}e_{0}  & =\omega\widehat{\chi}^{\prime}(\omega)\partial_{t_{1}t_{1}
}A_{0}e^{i\theta_{0}}+(\ast),\nonumber\\
\partial_{t_{2}t_{0}}\widehat{\chi}(i\partial_{t_{0}})e_{0}  & =-i\omega
\widehat{\chi}(\omega)\partial_{t_{2}}A_{0}e^{i\theta_{0}}+(\ast),\nonumber\\
\partial_{t_{1}t_{1}}\widehat{\chi}(i\partial_{t_{0}})e_{0}  & =\widehat{\chi
}(\omega)\partial_{t_{1}t_{1}}A_{0}e^{i\theta_{0}}+(\ast),\nonumber\\
\partial_{t_{0}t_{2}}\widehat{\chi}(i\partial_{t_{0}})e_{0}  & =-i\omega
\widehat{\chi}(\omega)\partial_{t_{2}}A_{0}e^{i\theta_{0}}+(\ast),\nonumber\\
\partial_{t_{0}t_{0}}e_{0}^{3}  & =-3\omega^{2}\eta|A_{0}|^{2}A_{0}e^{i\theta
}+NST+(\ast).\label{Eq8.31}
\end{align}
Inserting (\ref{Eq8.30}) and (\ref{Eq8.31}) into the right-hand side of the
order $\varepsilon^{2}$ equation we get
\begin{gather}
\partial_{t_{0}t_{0}}e_{2}-c^{2}\nabla_{0}^{2}e_{2}+\partial_{t_{0}t_{0}
}\widehat{\chi}(i\partial_{t_{0}})e_{2}=\nonumber\\
-\{-2i\omega\partial_{t_{2}}A_{0}+\partial_{t_{1}t_{1}}A_{0}-2ic^{2}
k\mathbf{u}\cdot\nabla_{2}A_{0}-c^{2}\nabla_{1}^{2}A_{0}\nonumber\\
-i\omega^{2}\widehat{\chi}^{\prime}(\omega)\partial_{t_{2}}A_{0}+\frac{1}
{2}\omega^{2}\widehat{\chi}^{\prime\prime}(\omega)\partial_{t_{1}t_{1}}
A_{0}+2\omega\widehat{\chi}^{\prime}(\omega)\partial_{t_{1}t_{1}}
A_{0}\nonumber\\
-2i\omega\widehat{\chi}(\omega)\partial_{t_{2}}A_{0}+\widehat{\chi}
(\omega)\partial_{t_{1}t_{1}}A_{0}-3\omega^{2}\eta|A_{0}|^{2}\}e^{i\theta_{0}
}+NST+(\ast).\label{Eq8.32}
\end{gather}
In order to remove secular terms we must postulate that
\begin{gather}
-2i\omega\partial_{t_{2}}A_{0}+\partial_{t_{1}t_{1}}A_{0}-2ic^{2}
k\mathbf{u}\cdot\nabla_{2}A_{0}-c^{2}\nabla_{1}^{2}A_{0}-i\omega
^{2}\widehat{\chi}^{\prime}(\omega)\partial_{t_{2}}A_{0}\nonumber\\
+\frac{1}{2}\omega^{2}\widehat{\chi}^{\prime\prime}(\omega)\partial
_{t_{1}t_{1}}A_{0}+2\omega\widehat{\chi}^{\prime}(\omega)\partial_{t_{1}t_{1}
}A_{0}-2i\omega\widehat{\chi}(\omega)\partial_{t_{2}}A_{0}+\widehat{\chi
}(\omega)\partial_{t_{1}t_{1}}A_{0}\nonumber\\
-3\omega^{2}\eta|A_{0}|^{2}=0.\label{Eq8.33}
\end{gather}
Using the dispersion relation (\ref{Eq8.22}), equation (\ref{Eq8.33}) can be
simplified into
\begin{equation}
\partial_{t_{2}}A_{0}+\mathbf{v}_{g}\cdot\nabla_{2}A_{0}-i\beta\nabla_{1}
^{2}A_{0}+i\alpha\partial_{t_{1}t_{1}}A_{0}-i\gamma|A_{0}|^{2}A_{0}
=0,\label{Eq8.34}
\end{equation}
where
\begin{align}
\alpha & =\omega^{\prime}\frac{n^{2}+2\omega\widehat{\chi}^{\prime}
(\omega)+\frac{1}{2}\omega^{2}\widehat{\chi}^{\prime\prime}(\omega)}{2c^{2}
k},\nonumber\\
\beta & =\frac{\omega^{\prime}}{2k},\nonumber\\
\gamma & =\frac{3\eta\omega^{2}\omega^{\prime}}{2c^{2}k}.\nonumber
\end{align}
Defining an amplitude $A(\mathbf{x},t)$ by
\begin{equation}
A(\mathbf{x},t)=A_{0}(\mathbf{x}_{1},t_{1},...)|_{t_{j}=e^{j}t,\mathbf{x}
_{j}=\varepsilon^{j}\mathbf{x}},\label{Eq8.36}
\end{equation}
and proceeding in the usual way, using (\ref{Eq8.26}) and (\ref{Eq8.34}), we
get the following amplitude equation
\begin{equation}
\partial_{t}A+\mathbf{v}_{g}\cdot\nabla A-i\beta\nabla^{2}A+i\alpha
\partial_{tt}A-i\gamma|A|^{2}A=0,\label{Eq8.37}
\end{equation}
where we have put the formal perturbation parameter equal to $1$. From what we
have done it is evident that for
\begin{equation}
E(\mathbf{x},t)=A(\mathbf{x},t)e^{i(\mathbf{k}\cdot\mathbf{x}-\omega t)}
+(\ast),\label{Eq8.38}
\end{equation}
to be an approximate solution to (\ref{Eq8.13}) we must have
\begin{align}
\gamma|A|^{2}  & \sim\beta\nabla^{2}A\sim\alpha\partial_{tt}A\sim
O(\varepsilon^{2}),\nonumber\\
\partial_{t}A  & \sim\mathbf{v}_{g}\cdot\nabla A\sim O(\varepsilon),\label{Eq8.39}
\end{align}
where $\varepsilon$ is a number much smaller than $1$. Under these
circumstances (\ref{Eq8.37}),(\ref{Eq8.38}) is the key elements in a fast
numerical scheme for wave packet solutions to (\ref{Eq8.13}). Because of the
presence of the second derivative with respect to time, equation
(\ref{Eq8.37}) can not be solved as a standard initial value problem. However,
because of (\ref{Eq8.39}) we can remove the second derivative term by
iteration
\begin{align}
\partial_{t}A  & =-\mathbf{v}_{g}\cdot\nabla A\sim O(\varepsilon
),\nonumber\\
& \Downarrow\nonumber\\
\partial_{tt}A  & =(\mathbf{v}_{g}\cdot\nabla)^{2}A\sim O(\varepsilon
^{2}),\label{Eq8.40}
\end{align}
which leads to the equation
\begin{equation}
\partial_{t}A+\mathbf{v}_{g}\cdot\nabla A-i\beta\nabla^{2}A+i\alpha
(\mathbf{v}_{g}\cdot\nabla)^{2}A-i\gamma|A|^{2}A=0,\label{Eq8.41}
\end{equation}
which \textit{can} be solved as a standard initial value problem.

  In deriving
this equation we assumed that the terms proportional to
\[
e^{\pm3i(\mathbf{k}\cdot\mathbf{x}-\omega t)},
\]
where nonsecular. For this to be true we must have
\begin{equation}
\omega(3k)\neq3\omega(k),\label{Eq8.42}
\end{equation}
where $\omega(k)$ is a solution to (\ref{Eq8.21}). If an equality holds in
(\ref{Eq8.42}) we have \textit{phase matching} and the multiple scale
calculation has to be redone, starting at (\ref{Eq8.19}), using a sum of two
wave packets with the appropriate center wave numbers and frequencies instead
of the single wave packet we used in the calculation leading to (\ref{Eq8.37}).
It could also be the case that we are modeling a situation where several wave
packets are interacting in a Kerr medium. For such a case we would instead of
(\ref{Eq8.19}) use a finite sum of wave packets
\begin{equation}
e_{0}(\mathbf{x}_{0},t_{0},\mathbf{x}_{1},t_{1},..)=
{\displaystyle\sum_{j=0}^{N}}
A_{j}(\mathbf{x}_{1},t_{1},...)e^{i\theta_{j}}+(\ast).\label{Eq8.43}
\end{equation}
Calculations analogous to the ones leading up to equation (\ref{Eq8.37})
\ will now give a separate equation of the type (\ref{Eq8.37}) for each wave
packet, \textit{unless} we have phase matching. These phase matching
conditions appears from the nonlinear term in the order $\varepsilon^{2}$
equation and takes the familiar form
\begin{align}
\mathbf{k}_{j}  & =s_{1}\mathbf{k}_{j_{1}}+s_{2}\mathbf{k}_{j_{2}}
+s_{3}\mathbf{k}_{j_{3}},\nonumber\\
\omega(k_{j})  & =s_{1}\omega(k_{j_{1}})+s_{2}\omega(k_{j_{2}})+s_{3}
\omega(k_{j_{3}}),\label{Eq8.44}
\end{align}
where $s=\pm1$. The existence of phase matching leads to coupling of the
amplitude equations. If (\ref{Eq8.44}) holds, the amplitude equation for
$A_{j}$ will contain a coupling term proportional to
\begin{equation}
A_{j_{1}}^{s_{1}}A_{j_{2}}^{s_{2}}A_{j_{3}}^{s_{3}}\label{Eq8.46}%
\end{equation}
where by definition $A_{j}^{+1}=A_{j}$ and $A_{j}^{-1}=A_{j}^{\ast}$.

We have seeen that assuming a scaling of $\varepsilon$ for space and time
variables and $\varepsilon^{2}$ for the nonlinear term leads to an amplitude
equation where second derivatives of space and time appears at the same order
as the cubic nonlinearity. This amplitude equation can thus describe a
situation where diffraction, group velocity dispersion and nonlinearity are of
the same size. Other choices of scaling for space,time and nonlinearity will
lead to other amplitude equations where other physical effects are of the same
size. Thus, the choice of scaling is determined by what kind of physics we want to describe.

\subsubsection{Linearly polarized vector wave packets}

Up til now all applications of the multiple scale method  PDEs  has involved scalar
equations. The multiple scale method is not limited to scalar equations, but is
equally applicable to vector equations. However, for vector equations we need
to be more careful than for the scalar case when it comes to eliminating
secular terms. We will here use Maxwell's equations (\ref{Eq8.7}) to
illustrate how the method is applied to vector PDEs in general. 

Assuming,
as usual, a polarization response induced by the Kerr effect, our basic
equations are
\begin{align}
\partial_{t}\mathbf{B+\nabla\times E}  & =0,\nonumber\\
\partial_{t}\mathbf{E-c}^{2}\mathbf{\nabla\times B+\partial}_{t}\widehat{\chi
}(i\partial_{t})\mathbf{E}  & =-\varepsilon^{2}\eta\partial_{t}(\mathbf{E}
^{2}\mathbf{E)},\nonumber\\
\nabla\cdot\mathbf{B}  & =0,\nonumber\\
\nabla\cdot\mathbf{E}+\widehat{\chi}(i\partial_{t})\nabla\cdot\mathbf{E}  &
=-\varepsilon^{2}\eta\nabla\cdot(\mathbf{E}^{2}\mathbf{E)},\label{Eq8.47}
\end{align}
where we have introduced a formal perturbation parameter in front of the
nonlinear terms. We now introduce the usual machinery of the multiple scale method.

Let $\mathbf{e}(\mathbf{x}_{0},t_{0},\mathbf{x}_{1},t_{1},...)$ and
$\mathbf{b}(\mathbf{x}_{0},t_{0},\mathbf{x}_{1},t_{1},...)$ be functions such
that
\begin{align}
\mathbf{E}(x,t)  & =\mathbf{e}(\mathbf{x}_{0},t_{0},\mathbf{x}_{1}
,t_{1},...)|_{\mathbf{x}_{j}=\varepsilon^{j}\mathbf{x},t_{j}=\varepsilon^{j}
t},\nonumber\\
\mathbf{B}(x,t)  & =\mathbf{b}(\mathbf{x}_{0},t_{0},\mathbf{x}_{1}
,t_{1},...)|_{\mathbf{x}_{j}=\varepsilon^{j}\mathbf{x},t_{j}=\varepsilon^{j}
t},\label{Eq8.48}
\end{align}
and let
\begin{align}
\partial_{t}  & =\partial_{t_{0}}+\varepsilon\partial_{t_{1}}+\varepsilon
^{2}\partial_{t_{2}}+...\;\;,\nonumber\\
\nabla\times & =\nabla_{0}\times+\varepsilon\nabla_{1}\times+\varepsilon
^{2}\nabla_{2}\times+...\;\;,\nonumber\\
\nabla\cdot & =\nabla_{0}\cdot+\varepsilon\nabla_{1}\cdot+\varepsilon
^{2}\nabla_{2}\cdot+...\;\;,\nonumber\\
\mathbf{e}  & =\mathbf{e}_{0}+\varepsilon\mathbf{e}_{1}+\varepsilon
^{2}\mathbf{e}_{2}+...\;\;,\nonumber\\
\mathbf{b}  & =\mathbf{b}_{0}+\varepsilon\mathbf{b}_{1}+\varepsilon
^{2}\mathbf{b}_{2}+...\;\;.\label{Eq8.49}
\end{align}
We now insert (\ref{Eq8.49}) into (\ref{Eq8.47}) and expand everything in
sight to second order in $\varepsilon$. Putting each order of $\varepsilon$ to
zero separately gives us the perturbation hierarchy. At this point you should
be able to do this on your own so I will just write down the elements of the
perturbation hierarchy when they are needed.

The order $\varepsilon^{0}$ equations, which is the first element of the
perturbation hierarchy, is of course
\begin{align}
\partial_{t_{0}}\mathbf{b}_{0}\mathbf{+\nabla}_{0}\mathbf{\times e}_{0}  &
=0,\nonumber\\
\partial_{t_{0}}\mathbf{e}_{0}\mathbf{-c}^{2}\mathbf{\nabla}_{0}\mathbf{\times
b}_{0}\mathbf{+\partial}_{t_{0}}\widehat{\chi}(i\partial_{t_{0}}
)\mathbf{e}_{0}  & =0,\nonumber\\
\nabla_{0}\cdot\mathbf{b}_{0}  & =0,\nonumber\\
\nabla_{0}\cdot\mathbf{e}_{0}+\widehat{\chi}(i\partial_{t_{0}})\nabla_{0}
\cdot\mathbf{e}_{0}  & =0.\label{Eq8.50}
\end{align}
For the order $\varepsilon^{0}$ equations, we chose a linearly polarized wave
packet solution. It must be of the form
\begin{align}
\mathbf{e}_{0}(\mathbf{x}_{0},t_{0},\mathbf{x}_{1},t_{1},...)  & =\omega
A_{0}(\mathbf{x}_{1},t_{1},...)\mathbf{qe}^{i\theta_{0}}+(\ast),\nonumber
\\
\mathbf{b}_{0}(\mathbf{x}_{0},t_{0},\mathbf{x}_{1},t_{1},...)  &
=kA_{0}(\mathbf{x}_{1},t_{1},...)\mathbf{te}^{i\theta_{0}}+(\ast),\label{Eq8.51}
\end{align}
where
\begin{equation}
\theta_{0}=\mathbf{k}\cdot\mathbf{x}_{0}-\omega t_{0},\label{Eq8.52}
\end{equation}
and where
\[
\omega=\omega(k),
\]
is a solution to the dispersion relation
\begin{equation}
\omega^{2}n^{2}(\omega)=c^{2}k^{2}.\label{Eq8.53}
\end{equation}
The orthogonal unit vectors $\mathbf{q}$ and $\mathbf{t}$ span the space
transverse to $\mathbf{k}=k\mathbf{u}$, and the unit vectors $\{\mathbf{q}
,\mathbf{t},\mathbf{u}\}$ define a positively oriented frame for $\mathbb{R}^{3}$.

The order $\varepsilon$ equations are
\begin{gather}
\partial_{t_{0}}\mathbf{b}_{1}\mathbf{+\nabla}_{0}\mathbf{\times e}
_{1}=-\partial_{t_{1}}\mathbf{b}_{0}-\nabla_{1}\times\mathbf{e}_{0},\nonumber\\
\partial_{t_{0}}\mathbf{e}_{1}\mathbf{-c}^{2}\mathbf{\nabla}_{0}\mathbf{\times
b}_{1}\mathbf{+\partial}_{t_{0}}\widehat{\chi}(i\partial_{t_{0}}
)\mathbf{e}_{1}=\nonumber\\
-\partial_{t_{1}}\mathbf{e}_{0}+c^{2}\nabla_{1}\times\mathbf{b}_{0}
-\partial_{t_{1}}\widehat{\chi}(i\partial_{t_{0}})\mathbf{e}_{0}
-i\partial_{t_{0}}\widehat{\chi}^{\prime}(i\partial_{t_{0}})\partial_{t_{1}
}\mathbf{e}_{0},\nonumber\\
\nabla_{0}\cdot\mathbf{b}_{1}=-\nabla_{1}\cdot b_{0},\nonumber\\
\nabla_{0}\cdot\mathbf{e}_{1}+\widehat{\chi}(i\partial_{t_{0}})\nabla_{0}
\cdot\mathbf{e}_{1}=\nonumber\\
-\nabla_{1}\cdot\mathbf{e}_{0}-\widehat{\chi}(i\partial_{t_{0}})\nabla
_{1}\cdot\mathbf{e}_{0}-i\widehat{\chi}^{\prime}(i\partial_{t_{0}}
)\partial_{t_{1}}\nabla_{0}\cdot\mathbf{e}_{0}.\label{Eq8.54}
\end{gather}
Inserting (\ref{Eq8.51}) into (\ref{Eq8.54}) we get
\begin{gather}
\partial_{t_{0}}\mathbf{b}_{0}\mathbf{+\nabla}_{0}\mathbf{\times e}
_{0}=-\{k\partial_{t_{1}}A_{0}\mathbf{t}+\omega\nabla_{1}A_{0}\times
\mathbf{q\}}e^{i\theta_{0}}+(\ast),\nonumber\\
\partial_{t_{0}}\mathbf{e}_{0}\mathbf{-c}^{2}\mathbf{\nabla}_{0}\mathbf{\times
b}_{0}\mathbf{+\partial}_{t_{0}}\widehat{\chi}(i\partial_{t_{0}}
)\mathbf{e}_{0}=-\{(\omega n^{2}(\omega)+\omega^{2}\widehat{\chi}^{\prime
}(\omega))\partial_{t_{1}}A_{0}\mathbf{q}\nonumber\\
-c^{2}k\nabla_{1}A_{0}\times\mathbf{t}\}e^{i\theta_{0}}+(\ast),\nonumber\\
\nabla_{0}\cdot\mathbf{b}_{0}=-\{k\nabla_{1}A_{0}\cdot\mathbf{t\}}
e^{i\theta_{0}}+(\ast),\nonumber\\
\nabla_{0}\cdot\mathbf{e}_{0}+\widehat{\chi}(i\partial_{t_{0}})\nabla_{0}
\cdot\mathbf{e}_{0}=-\{\omega n^{2}(\omega)\nabla_{1}A_{0}\cdot\mathbf{q\}}
e^{i\theta_{0}}+(\ast).\label{Eq8.55}
\end{gather}  
If we can find a special solution to this system that is bounded, we will get a
perturbation expansion that is uniform for $t\lesssim\varepsilon^{-1}$. We
will look for solutions of the form
\begin{align}
\mathbf{e}_{1}  & =\mathbf{a}e^{i\theta_{0}}+(\ast),\nonumber\\
\mathbf{b}_{1}  & =\mathbf{b}e^{i\theta_{0}}+(\ast),\label{Eq8.56}
\end{align}
where $\mathbf{a}$ and $\mathbf{b}$ are constant vectors. Inserting
(\ref{Eq8.56}) into (\ref{Eq8.55}), we get the following linear algebraic
system of equations for the unknown vectors $\mathbf{a}$ and $\mathbf{b}$
\begin{gather}
-i\omega\mathbf{b}+ik\mathbf{u}\times\mathbf{a}=-\{k\partial_{t_{1}}
A_{0}\mathbf{t}+\omega\nabla_{1}A_{0}\times\mathbf{q\}},\label{Eq8.57}\\
-i\omega n^{2}(\omega)\mathbf{a}-ic^{2}k\mathbf{u}\times\mathbf{b}=-\{(\omega
n^{2}(\omega)+\omega^{2}\widehat{\chi}^{\prime}(\omega))\partial_{t_{1}}
A_{0}\mathbf{q}\nonumber\\
-c^{2}k\nabla_{1}A_{0}\times\mathbf{t}\},\label{Eq8.58}\\
ik\mathbf{u}\cdot\mathbf{b}=-k\nabla_{1}A_{0}\cdot\mathbf{t},\label{Eq8.59}\\
ikn^{2}(\omega)\mathbf{u}\cdot\mathbf{a}=-\omega n^{2}(\omega)\nabla_{1}
A_{0}\cdot\mathbf{q}\;.\label{Eq8.60}
\end{gather}
Introduce the longitudinal and transverse parts of $\mathbf{a}$ and
$\mathbf{b}$ through
\begin{align}
a_{\Vert}  & =(\mathbf{u}\cdot\mathbf{a)u},\text{ \ \ \ \ \ \ \ }a_{\bot
}=\mathbf{a}-a_{\Vert},\nonumber\\
b_{\Vert}  & =(\mathbf{u}\cdot\mathbf{b})\mathbf{u},\text{ \ \ \ \ \ \ \ }
b_{\bot}=\mathbf{b}-b_{\Vert}.\label{Eq8.61}
\end{align}
Then from (\ref{Eq8.59}) and (\ref{Eq8.60}) we get
\begin{align}
a_{\Vert}  & =(i\frac{\omega}{k}\nabla_{1}A_{0}\cdot\mathbf{q)u},
\label{Eq8.62}\\
b_{\Vert}  & =(i\nabla_{1}A_{0}\cdot\mathbf{t})\mathbf{u}.\label{Eq8.63}
\end{align}
However, the longitudinal part of (\ref{Eq8.57}) and (\ref{Eq8.58}) will also
determine $a_{\Vert}$ and $b_{\Vert}$. These values must be the same as the
ones just found in (\ref{Eq8.62}),(\ref{Eq8.63}). These are
\textit{solvability conditions}. Taking the longitudinal part of
(\ref{Eq8.57}) we get
\begin{gather}
-i\omega\mathbf{u}\cdot\mathbf{b}=-\omega\mathbf{u}\cdot(\nabla_{1}A_{0}
\times\mathbf{q}),\nonumber\\
\Updownarrow\nonumber\\
\mathbf{u}\cdot\mathbf{b}=i\nabla_{1}A_{0}\cdot\mathbf{t},\label{Eq8.64}
\end{gather}
which is consistent with (\ref{Eq8.63}). Thus this solvability condition is
automatically satisfied. Taking the longitudinal part of (\ref{Eq8.58}) we get
\begin{gather}
-i\omega n^{2}(\omega)\mathbf{u}\cdot\mathbf{a}\mathbf{=}c^{2}k\mathbf{u\cdot
(\nabla}_{1}A_{0}\times\mathbf{t}),\nonumber\\
\Updownarrow\nonumber\\
\mathbf{u}\cdot\mathbf{a}=i\frac{\omega}{k}\nabla_{1}A_{0}\cdot\mathbf{q},\label{Eq8.65}
\end{gather}
which is consistent with (\ref{Eq8.62}). Thus this solvability condition is
also automatically satisfied. 

The transversal part of (\ref{Eq8.57}) and
(\ref{Eq8.58}) are
\begin{align}
-i\omega b_{\bot}+ik\mathbf{u}\times a_{\bot}  & =-\{k\partial_{t_{1}}
A_{0}+\omega\nabla_{1}A_{0}\cdot\mathbf{u\}t},\label{Eq8.66}\\
-i\omega n^{2}(\omega)a_{\bot}-ic^{2}k\mathbf{u}\times b_{\bot}  &
=-\{\omega(n^{2}(\omega)+\omega\widehat{\chi}^{\prime}(\omega))\partial
_{t_{1}}A_{0}+c^{2}k\nabla_{1}A_{0}\cdot\mathbf{u}\}\mathbf{q},\nonumber
\end{align}
and this linear system is singular; the determinant is zero because of the
dispersion relation (\ref{Eq8.53}). It can therefore only be solved if the
right-hand side satisfy a certain solvability condition. The most effective
way to find this condition is to use the \textit{Fredholm Alternative}. It say
that a linear system
\[
A\mathbf{x}=\mathbf{c},
\]
has a solution if and only if
\[
\mathbf{f}\cdot\mathbf{c}=0,
\]
for all vectors $\mathbf{f}$, such that
\[
A^{\dag}\mathbf{f}=0,
\]
where $A^{\dag}$ is the adjoint of $A$.

The matrix for the system (\ref{Eq8.66}) is
\[
M=\left(
\begin{tabular}
[c]{ll}
$ik\mathbf{u}\times$ & $-i\omega$\\
$-i\omega n^{2}$ & $-ic^{2}k\mathbf{u}\times$
\end{tabular}
\right).
\]
The adjoint of this matrix is clearly

\begin{equation}
M^{\dag}=\left(
\begin{tabular}
[c]{ll}
$-ik\mathbf{u}\times$ & $-i\omega n^{2}$\\
$-i\omega$ & $ic^{2}k\mathbf{u}\times$
\end{tabular}
\right), \label{Eq8.67}
\end{equation}
and the null space of the adjoint is thus determined by
\begin{align}
-ik\mathbf{u}\times\mathbf{\alpha}-i\omega n^{2}\mathbf{\beta}  &
=0,\nonumber\\
-i\omega\mathbf{\alpha}+ic^{2}k\mathbf{u}\times\mathbf{\beta}  & =0.\label{Eq8.68}
\end{align}
A convenient basis for the null space is
\begin{equation}
\left\{  \left(
\begin{tabular}
[c]{l}
$-c^{2}k\mathbf{q}$\\
$\omega\mathbf{t}$
\end{tabular}
\right)  ,\left(
\begin{tabular}
[c]{l}
$c^{2}k\mathbf{t}$\\
$\omega\mathbf{q}$
\end{tabular}
\right)  \right\} \label{Eq8.69}
\end{equation}
The first basis vector gives a trivial solvability condition, whereas the
second one gives a nontrivial condition, which is
\begin{gather}
c^{2}k\{k\partial_{t_{1}}A_{0}+\omega\nabla_{1}A_{0}\cdot\mathbf{u\}+\omega
}\{\omega(n^{2}(\omega)+\omega\widehat{\chi}^{\prime}(\omega))\partial_{t_{1}
}A_{0}+c^{2}k\nabla_{1}A_{0}\cdot\mathbf{u}\}=0,\nonumber\\
\Updownarrow\nonumber\\
\omega^{2}(2n^{2}+\omega\widehat{\chi}^{\prime}(\omega))\partial_{t_{1}}
A_{0}+2c^{2}k\omega\mathbf{u}\cdot\nabla_{1}A_{0}=0.\label{Eq8.70}
\end{gather}
Observe that from the dispersion relation (\ref{Eq8.53}) we have
\begin{gather}
\omega^{2}n^{2}(\omega)=\omega^{2}(1+\widehat{\chi}(\omega))=c^{2}
k^{2},\nonumber\\
\Downarrow\nonumber\\
2\omega\omega^{\prime}n^{2}+\omega^{2}\widehat{\chi}^{\prime}(\omega
)\omega^{\prime}=2c^{2}k,\nonumber\\
\Downarrow\nonumber\\
\omega(2n^{2}+\omega\widehat{\chi}^{\prime}(\omega))\omega^{\prime}
=2c^{2}k.\label{Eq8.71}
\end{gather}
Using (\ref{Eq8.71}) in (\ref{Eq8.70}) the solvability condition can be
compactly written as
\begin{equation}
\partial_{t_{1}}A_{0}+\mathbf{v}_{g}\cdot\nabla_{1}A_{0}=0,\label{Eq8.72}
\end{equation}
where $\mathbf{v}_{g}$ is the \textit{group velocity}
\begin{equation}
\mathbf{v}_{g}=\frac{d\omega}{dk}\mathbf{u}.\label{Eq8.73}
\end{equation}
The system (\ref{Eq8.66}) is singular but consistent. We can therefore
disregard the second equation, and look for a special solution of the form
\begin{align}
a_{\bot}  & =a\mathbf{q},\nonumber\\
b_{\bot}  & =0.\label{Eq8.74}
\end{align}
Inserting (\ref{Eq8.74}) into the first equation in (\ref{Eq8.66}) we easily get
\begin{equation}
a_{\bot}=i\left\{  \partial_{t_{1}}A_{0}+\frac{\omega}{k}\mathbf{u}\cdot
\nabla_{1}A_{0}\right\}  \mathbf{q}.\label{Eq8.75}
\end{equation}
From (\ref{Eq8.62}),(\ref{Eq8.63}),(\ref{Eq8.74}) and (\ref{Eq8.75}), we get
the following bounded special solution to the order $\varepsilon$ equations
\begin{align}
\mathbf{e}_{1}  & =\{i(\partial_{t_{1}}A_{0}+\frac{\omega}{k}\mathbf{u}
\cdot\nabla_{1}A_{0})\mathbf{q}+i(\frac{\omega}{k}\mathbf{q}\cdot\nabla
_{1}A_{0})\mathbf{u\}e}^{i\theta_{0}}+(\ast),\nonumber\\
\mathbf{b}_{1}  & =\{i(\mathbf{t}\cdot\nabla_{1}A_{0})\mathbf{u\}e}
^{i\theta_{0}}+(\ast).\label{Eq8.76}
\end{align}
The order $\varepsilon^{2}$ equations are
\begin{gather}
\partial_{t_{0}}\mathbf{b}_{2}\mathbf{+\nabla}_{0}\mathbf{\times e}
_{2}=-\{\partial_{t_{1}}\mathbf{b}_{1}+\nabla_{1}\times\mathbf{e}_{1}
+\partial_{t_{2}}\mathbf{b}_{0}+\nabla_{2}\times\mathbf{e}_{0}\},\nonumber
\\
\nonumber\\
\partial_{t_{0}}\mathbf{e}_{2}\mathbf{-c}^{2}\mathbf{\nabla}_{0}\mathbf{\times
b}_{2}\mathbf{+\partial}_{t_{0}}\widehat{\chi}(i\partial_{t_{0}}
)\mathbf{e}_{2}=-\{\partial_{t_{1}}\mathbf{e}_{1}-c^{2}\nabla_{1}
\times\mathbf{b}_{1}+\partial_{t_{2}}\mathbf{e}_{0}\nonumber\\
-c_{2}\nabla_{2}\times\mathbf{b}_{0}+\partial_{t_{1}}\widehat{\chi}
(i\partial_{t_{0}})\mathbf{e}_{1}+i\partial_{t_{0}}\widehat{\chi}^{\prime
}(i\partial_{t_{0}})\partial_{t_{1}}\mathbf{e}_{1}\nonumber\\
\partial_{t_{2}}\widehat{\chi}(i\partial_{t_{0}})\mathbf{e}_{0}+i\partial
_{t_{1}}\widehat{\chi}^{\prime}(i\partial_{t_{0}})\partial_{t_{1}}
\mathbf{e}_{0}+i\partial_{t_{0}}\widehat{\chi}^{\prime}(i\partial_{t_{0}
})\partial_{t_{2}}\mathbf{e}_{0}\nonumber\\
-\frac{1}{2}\partial_{t_{0}}\widehat{\chi}^{\prime\prime}(i\partial_{t_{0}
})\partial_{t_{1}t_{1}}\mathbf{e}_{0}+\eta\partial_{t_{0}}\mathbf{e}_{0}
^{2}\mathbf{e}_{0}\},\nonumber\\
\nonumber\\
\nabla_{0}\cdot\mathbf{b}_{2}=-\{\nabla_{1}\cdot\mathbf{b}_{1}+\nabla_{2}
\cdot\mathbf{b}_{0}\},\nonumber\\
\nonumber\\
\nabla_{0}\cdot\mathbf{e}_{2}+\widehat{\chi}(i\partial_{t_{0}})\nabla_{0}
\cdot\mathbf{e}_{2}=-\{\nabla_{1}\cdot\mathbf{e}_{1}+\nabla_{2}\cdot
\mathbf{e}_{0}+\widehat{\chi}(i\partial_{t_{0}})\nabla_{1}\cdot\mathbf{e}
_{1}\nonumber\\
+i\widehat{\chi}^{\prime}(i\partial_{t_{0}})\partial_{t_{1}}\nabla_{0}
\cdot\mathbf{e}_{1}+\widehat{\chi}(i\partial_{t_{0}})\nabla_{2}\cdot
\mathbf{e}_{0}+i\widehat{\chi}^{\prime}(i\partial_{t_{0}})\partial_{t_{2}
}\nabla_{0}\cdot\mathbf{e}_{0}\nonumber\\
+i\widehat{\chi}^{\prime}(i\partial_{t_{0}})\partial_{t_{1}}\nabla_{1}
\cdot\mathbf{e}_{0}-\frac{1}{2}\widehat{\chi}^{\prime\prime}(i\partial_{t_{0}
})\partial_{t_{1}t_{1}}\nabla_{0}\cdot\mathbf{e}_{0}+\eta\nabla_{0}
\cdot(\mathbf{e}_{0}^{2}\mathbf{e}_{0})\}.\label{Eq8.77}
\end{gather}

We now insert (\ref{Eq8.51}) and (\ref{Eq8.76}) into (\ref{Eq8.77}). This
gives us
\begin{gather}
\partial_{t_{0}}\mathbf{b}_{2}\mathbf{+\nabla}_{0}\mathbf{\times e}
_{2}=-\{i(\partial_{t_{1}}\nabla_{1}A_{0}\cdot\mathbf{t)u}+i\mathbf{\nabla
}_{1}\partial_{t_{1}}A_{0}\times\mathbf{q},\nonumber\\
+i\frac{\omega}{k}\nabla_{1}(\nabla_{1}A_{0}\cdot\mathbf{u)\times q}
+i\frac{\omega}{k}\nabla_{1}(\nabla_{1}A_{0}\cdot\mathbf{q)\times
u}+k\mathbf{\partial}_{t_{2}}A_{0}\mathbf{t}\nonumber\\
\mathbf{+\omega\nabla}_{2}A_{0}\times\mathbf{q\}e}^{i\theta_{0}}
+(\ast),\nonumber\\
\nonumber\\
\partial_{t_{0}}\mathbf{e}_{2}\mathbf{-}c^{2}\mathbf{\nabla}_{0}\mathbf{\times
b}_{2}\mathbf{+\partial}_{t_{0}}\widehat{\chi}(i\partial_{t_{0}}
)\mathbf{e}_{2}=-\{iF(\omega)\partial_{t_{1}t_{1}}A_{0}\mathbf{q}\nonumber\\
+i\mathbf{G(\omega)(\partial}_{t_{1}}\nabla_{1}A_{0}\cdot\mathbf{u)q}
+iG(\omega)(\partial_{t_{1}}\nabla_{1}A_{0}\cdot\mathbf{q)u}-ic^{2}\nabla
_{1}(\nabla_{1}A_{0}\cdot\mathbf{t})\times\mathbf{u}\nonumber\\
-c^{2}k\nabla_{2}A_{0}\times\mathbf{t}+H\mathbf{(\omega)\partial}_{t_{2}}
A_{0}\mathbf{q}-3i\eta\omega^{4}|A_{0}|^{2}A_{0}\}e^{i\theta_{0}}
+(\ast)\nonumber,\\
\nonumber\\
\nabla_{0}\cdot\mathbf{b}_{2}=-\{i\nabla_{1}(\nabla_{1}A_{0}\cdot
\mathbf{t})\cdot\mathbf{u}+k\mathbf{\nabla}_{2}A_{0}\cdot\mathbf{t\}e}
^{i\theta_{0}}+(\ast),\nonumber\\
\nonumber\\
\nabla_{0}\cdot\mathbf{e}_{2}+\widehat{\chi}(i\partial_{t_{0}})\nabla_{0}
\cdot\mathbf{e}_{2}=-\{in^{2}\nabla_{1}\partial_{t_{1}}A_{0}\cdot
\mathbf{q}+in^{2}\frac{\omega}{k}\nabla_{1}(\nabla_{1}A_{0}\cdot
\mathbf{u})\cdot\mathbf{q}\nonumber\\
+in^{2}\frac{\omega}{k}\nabla_{1}(\nabla_{1}\cdot\mathbf{q})\cdot
\mathbf{u}+\omega n^{2}\nabla_{2}A_{0}\cdot\mathbf{q}\}e^{i\theta_{0}}
+(\ast),\label{Eq8.78}
\end{gather}
where we have defined
\begin{align}
F(\omega)  & =n^{2}+2\omega\widehat{\chi}^{\prime}(\omega)+\frac{1}{2}
\omega^{2}\widehat{\chi}^{\prime\prime}(\omega),\nonumber\\
G(\omega)  & =\frac{\omega}{k}(n^{2}+\omega\widehat{\chi}^{\prime}
(\omega)),\nonumber\\
H(\omega)  & =\omega(n^{2}+\omega\widehat{\chi}^{\prime}(\omega)).\label{Eq8.79}
\end{align}
Like for the order $\varepsilon$ equations, we will look for bounded solutions
of the form
\begin{align}
\mathbf{e}_{2}  & =\mathbf{a}e^{i\theta_{0}}+(\ast),\nonumber\\
\mathbf{b}_{2}  & =\mathbf{b}e^{i\theta_{0}}+(\ast).\label{Eq8.80}
\end{align}
Inserting (\ref{Eq8.80}) into (\ref{Eq8.78}) we get the following linear
system of equations for the constant vectors $\mathbf{a}$ and $\mathbf{b}$
\begin{gather}
-i\omega\mathbf{b}+ik\mathbf{u}\times\mathbf{a}=-\{i\mathbf{(}\partial_{t_{1}
}\nabla_{1}A_{0}\cdot\mathbf{t})\mathbf{u}+i\mathbf{\nabla}_{1}\partial
_{t_{1}}A_{0}\times\mathbf{q}\nonumber\\
+i\frac{\omega}{k}\nabla_{1}(\nabla_{1}A_{0}\cdot\mathbf{u)\times q}
+i\frac{\omega}{k}\nabla_{1}(\nabla_{1}A_{0}\cdot\mathbf{q)\times
u}+k\mathbf{\partial}_{t_{2}}A_{0}\mathbf{t}\nonumber\\
\mathbf{+\omega\nabla}_{2}A_{0}\times\mathbf{q\}},\label{Eq8.81}\\
\nonumber\\
-i\omega n^{2}(\omega)\mathbf{a}-ic^{2}k\mathbf{u}\times\mathbf{b}
=-\{iF(\omega)\partial_{t_{1}t_{1}}A_{0}\mathbf{q}\nonumber\\
+iG\mathbf{(\omega)(\partial}_{t_{1}}\nabla_{1}A_{0}\cdot\mathbf{u)q}
+iG\mathbf{(\omega)(\partial}_{t_{1}}\nabla_{1}A_{0}\cdot\mathbf{q)u}
\nonumber\\
-ic^{2}\nabla_{1}(\nabla_{1}A_{0}\cdot\mathbf{t})\times\mathbf{u}-c^{2}
k\nabla_{2}A_{0}\times\mathbf{t}\nonumber\\
+H\mathbf{(\omega)\partial}_{t_{2}}A_{0}\mathbf{q}-3i\eta\omega^{4}|A_{0}
|^{2}A_{0}\},\label{Eq8.82}\\
\nonumber\\
ik\mathbf{u}\cdot\mathbf{b}=-\{i\nabla_{1}(\nabla_{1}A_{0}\cdot\mathbf{t}
)\cdot\mathbf{u}+k\mathbf{\nabla}_{2}A_{0}\cdot\mathbf{t\}},\label{Eq8.83}\\
\nonumber\\
ikn^{2}\mathbf{u}\cdot\mathbf{a}=-\{in^{2}\nabla_{1}\partial_{t_{1}}A_{0}
\cdot\mathbf{q}+in^{2}\frac{\omega}{k}\nabla_{1}(\nabla_{1}A_{0}
\cdot\mathbf{u})\cdot\mathbf{q}\nonumber\\
+in^{2}\frac{\omega}{k}\nabla_{1}(\nabla_{1}\cdot\mathbf{q})\cdot
\mathbf{u}+\omega n^{2}\nabla_{2}A_{0}\cdot\mathbf{q}\}.\label{Eq8.84}
\end{gather}
We introduce the  longitudinal and transversal vector components for $\mathbf{a}$
and $\mathbf{b}$ like before, and find from (\ref{Eq8.83}) and (\ref{Eq8.84})
that%
\begin{align}
a_{\Vert}  & =(-\frac{1}{k}\nabla_{1}\partial_{t_{1}}A_{0}\cdot\mathbf{q}
-\frac{\omega}{k^{2}}\nabla_{1}(\nabla_{1}A_{0}\cdot\mathbf{u})\cdot
\mathbf{q}\nonumber\\
& -\frac{\omega}{k^{2}}\nabla_{1}(\nabla_{1}A_{0}\cdot\mathbf{q}
)\cdot\mathbf{u}+i\frac{\omega}{k}\nabla_{2}A_{0}\cdot\mathbf{q})\mathbf{u},\label{Eq8.85}\\
b_{\Vert}  & =(i\nabla_{2}A_{0}\cdot\mathbf{t}-\frac{1}{k}\nabla_{1}
(\nabla_{1}\cdot\mathbf{t})\cdot\mathbf{u})\mathbf{u}.\label{Eq8.86}
\end{align}
The longitudinal part of (\ref{Eq8.81}) is
\begin{equation}
\mathbf{u}\cdot\mathbf{b}=\frac{1}{\omega}\{\partial_{t_{1}}\nabla_{1}
A_{0}\cdot\mathbf{t-\nabla}_{1}\partial_{t_{1}}A_{0}\cdot\mathbf{t}
-\frac{\omega}{k}\nabla_{1}(\nabla_{1}A_{0}\cdot\mathbf{u})\cdot
\mathbf{t}+i\omega\mathbf{\nabla}_{2}A_{0}\cdot\mathbf{t}\},\label{Eq8.87}
\end{equation}
and in order for (\ref{Eq8.87}) to be consistent with (\ref{Eq8.86}), we find that
the following solvability condition must hold
\begin{equation}
\partial_{t_{1}}\nabla_{1}A_{0}\cdot\mathbf{t}=\nabla_{1}\partial_{t_{1}}
A_{0}\cdot\mathbf{t}.\label{Eq8.88}
\end{equation}
The longitudinal part of (\ref{Eq8.82}) is
\begin{equation}
\mathbf{u}\cdot\mathbf{a=}\frac{1}{\omega n^{2}}\{G(\omega)\partial_{t_{1}
}\nabla_{1}A_{0}\cdot\mathbf{q}+ic^{2}k\nabla_{2}A_{0}\cdot\mathbf{q}
\},\label{Eq8.89}
\end{equation}
and in order for (\ref{Eq8.89}) to be consistent with (\ref{Eq8.85}) we find,
after a little algebra, that the solvability condition
\begin{gather}
\frac{\omega}{k}n^{2}(\omega)\nabla_{1}\partial_{t_{1}}A_{0}\cdot
\mathbf{q}+G(\omega)\mathbf{\partial}_{t_{1}}\nabla_{1}A_{0}\cdot
\mathbf{q=}\nonumber\\
-c^{2}\nabla_{1}(\nabla_{1}A_{0}\cdot\mathbf{q)\cdot u}-c^{2}\nabla_{1}
(\nabla_{1}A_{0}\cdot\mathbf{u})\cdot\mathbf{q},\label{Eq8.90}
\end{gather}
must hold. 

The transverse parts of (\ref{Eq8.81}) and (\ref{Eq8.82}) are
\begin{gather}
-i\omega b_{\bot}+ik\mathbf{u}\times a_{\bot}=-\{i\nabla_{1}\partial_{t_{1}
}A_{0}\cdot\mathbf{u}+i\frac{\omega}{k}\nabla_{1}(\nabla_{1}A_{0}
\cdot\mathbf{u})\cdot\mathbf{u}\nonumber\\
-i\frac{\omega}{k}\nabla_{1}(\nabla_{1}A_{0}\cdot\mathbf{q})\cdot
\mathbf{q}+k\partial_{t_{2}}A_{0}+\omega\nabla_{2}A_{0}\cdot\mathbf{u\}t}
-\{i\frac{\omega}{k}\nabla_{1}(\nabla_{1}A_{0}\cdot\mathbf{q)\cdot
t\}q},\nonumber\\
\nonumber\\
-i\omega n^{2}a_{\bot}-ic^{2}k\mathbf{u}\times b_{\bot}=-\{iF(\omega
)\partial_{t_{1}t_{1}}A_{0}+iG(\omega)\partial_{t_{1}}\nabla_{1}A_{0}
\cdot\mathbf{u}\nonumber\\
-ic^{2}\nabla_{1}(\nabla_{1}A_{0}\cdot\mathbf{t})\cdot\mathbf{t}+c^{2}
k\nabla_{2}A_{0}\cdot\mathbf{u}+H(\omega)\partial_{t_{2}}A_{0}-3\eta
i\omega^{4}|A_{0}|^{2}A_{0}\}\mathbf{q}\nonumber\\
-\{ic^{2}\nabla_{1}(\nabla_{1}A_{0}\cdot\mathbf{t})\cdot\mathbf{q}
\}\mathbf{t}.\label{Eq8.91}
\end{gather}

\noindent  The matrix for this linear system is the same as for the order $\varepsilon$
case, (\ref{Eq8.66}), so that the two solvability conditions are determined,
through the Fredholm Alternative, by the vectors (\ref{Eq8.69}). The
solvability condition corresponding to the first of the vectors in
(\ref{Eq8.69}) is
\begin{gather}
(-c^{2}k)(-i\frac{\omega}{k}\nabla_{1}(\nabla_{1}A_{0}\cdot\mathbf{q)\cdot
t)}+\omega\mathbf{(-}ic^{2}\nabla_{1}(\nabla_{1}A_{0}\cdot\mathbf{t}
)\cdot\mathbf{q})=0,\nonumber\\
\Updownarrow\nonumber\\
\nabla_{1}(\nabla_{1}A_{0}\cdot\mathbf{q})\cdot\mathbf{t}=\nabla_{1}
(\nabla_{1}\cdot\mathbf{t})\cdot\mathbf{q},\label{Eq8.92}
\end{gather}
and the solvability condition corresponding to the second vector in
(\ref{Eq8.69}) is
\begin{gather}
c^{2}k(-\{i\nabla_{1}\partial_{t_{1}}A_{0}\cdot\mathbf{u}+i\frac{\omega}
{k}\nabla_{1}(\nabla_{1}A_{0}\cdot\mathbf{u})\cdot\mathbf{u}\nonumber\\
-i\frac{\omega}{k}\nabla_{1}(\nabla_{1}A_{0}\cdot\mathbf{q})\cdot
\mathbf{q}+k\partial_{t_{2}}A_{0}+\omega\nabla_{2}A_{0}\cdot\mathbf{u\})}
+\omega(-\{iF(\omega)\partial_{t_{1}t_{1}}A_{0}\nonumber\\
+iG(\omega)\partial_{t_{1}}\nabla_{1}A_{0}\cdot\mathbf{u}-ic^{2}\nabla
_{1}(\nabla_{1}A_{0}\cdot\mathbf{t})\cdot\mathbf{t}+c^{2}k\nabla_{2}A_{0}
\cdot\mathbf{u}\nonumber\\
+H(\omega)\partial_{t_{2}}A_{0}-3\eta i\omega^{4}|A_{0}|^{2}A_{0}
\}\mathbf{q})=0,\nonumber\\
\Updownarrow\nonumber\\
\partial_{t_{2}}A_{0}+\mathbf{v}_{g}\cdot\nabla_{2}A_{0}+i\delta_{1}\nabla
_{1}\partial_{t_{1}}A_{0}\cdot\mathbf{u}+i\delta_{2}\partial_{t_{1}}\nabla
_{1}A_{0}\cdot\mathbf{u}\nonumber\\
-i\beta(\nabla_{1}(\nabla_{1}A_{0}\cdot\mathbf{q)\cdot q}+\nabla_{1}
(\nabla_{1}A_{0}\cdot\mathbf{t})\cdot\mathbf{t-\nabla}_{1}(\nabla_{1}
A_{0}\cdot\mathbf{u})\cdot\mathbf{u})\nonumber\\
+i\alpha\partial_{t_{1}t_{1}}A_{0}-i\gamma|A_{0}|^{2}A_{0}=0,\label{Eq8.93}
\end{gather}
where we have defined
\begin{align}
\alpha & =\frac{\omega^{\prime}F(\omega)}{2c^{2}k},\nonumber\\
\beta & =\frac{\omega^{\prime}}{2k},\nonumber\\
\gamma & =\frac{3\eta\omega^{\prime}\omega^{4}}{2c^{2}k},\nonumber\\
\delta_{1}  & =\frac{\omega^{\prime}}{2\omega},\nonumber\\
\delta_{2}  & =\frac{\omega^{\prime}G(\omega)}{2c^{2}k}.\nonumber
\end{align}
We have now found all solvability conditions. These are (\ref{Eq8.88}
),(\ref{Eq8.90}),(\ref{Eq8.92}) and (\ref{Eq8.93}). 

We now, as usual, define an amplitude
$A(\mathbf{x},t)$ by
\[
A(\mathbf{x},t)=A_{0}(\mathbf{x}_{1},t_{1},...)|\mathbf{x}_{j}=\varepsilon
^{j}\mathbf{x},t_{j}=\varepsilon^{j}t,
\]
and derive the amplitude equations from the solvability conditions in the
usual way. This gives us the following system
\begin{gather}
\partial_{t}\nabla A\cdot\mathbf{t=\nabla}\partial_{t}A\cdot\mathbf{t},
\label{Eq8.95}\\
\nonumber\\
\frac{\omega}{k}n^{2}(\omega)\nabla\partial_{t}A\cdot\mathbf{q}
+G\mathbf{(\omega)\partial}_{t}\nabla A\cdot\mathbf{q=}\nonumber\\
-c^{2}\nabla(\nabla A\cdot\mathbf{q)\cdot u}-c^{2}\nabla(\nabla A\cdot
\mathbf{u})\cdot\mathbf{q},\label{Eq8.96}\\
\nonumber\\
\nabla(\nabla A_{0}\cdot\mathbf{q})\cdot\mathbf{t}=\nabla(\nabla
A\cdot\mathbf{t})\cdot\mathbf{q},\label{Eq8.97}\\
\nonumber\\
\partial_{t}A+\mathbf{v}_{g}\cdot\nabla A+i\delta_{1}\nabla\partial_{t}
A\cdot\mathbf{u}+i\delta_{2}\partial_{t}\nabla A\cdot\mathbf{u}\nonumber
\\
-i\beta(\nabla(\nabla A\cdot\mathbf{q)\cdot q}+\nabla(\nabla A\cdot
\mathbf{t})\cdot\mathbf{t-\nabla}(\nabla A\cdot\mathbf{u})\cdot\mathbf{u}
)\nonumber\\
+i\alpha\partial_{tt}A-i\gamma|A|^{2}A=0,\label{Eq8.98}
\end{gather}
where we as usual have set the formal perturbation parameter equal to $1$.
Equations (\ref{Eq8.95}) and (\ref{Eq8.97}) are automatically satisfied since
$A(\mathbf{x},t)$ is a smooth function of space and time. We know that only
amplitudes such that
\begin{equation}
\partial_{t}A\sim-\mathbf{v}_{g}\cdot\nabla A=\omega^{\prime}\nabla
A\cdot\mathbf{u},\label{Eq8.99}
\end{equation}
can be allowed as solutions. This is assumed by the multiple scale method. If
we insert (\ref{Eq8.99}) into (\ref{Eq8.96}), assume smoothness and use the
dispersion relation, we find that (\ref{Eq8.96}) is automatically satisfied.
The only remaining equation is then (\ref{Eq8.98}) and if we insert the
approximation (\ref{Eq8.99}) for the derivatives with respect to time in the
second and third term of (\ref{Eq8.98}) we get, using the dispersion relation,
that (\ref{Eq8.98}) simplify into
\begin{equation}
\partial_{t}A+\mathbf{v}_{g}\cdot\nabla A-i\beta\nabla^{2}A+i\alpha
\partial_{tt}A-i\gamma|A|^{2}A=0,\label{Eq8.100}
\end{equation}
where we have also used the fact that
\[
\mathbf{qq}+\mathbf{tt}+\mathbf{uu}=I.
\]
Equation (\ref{Eq8.100}) is the celebrated 3D
nonlinear Schr\o dinger equation, including group velocity dispersion, and is a key equation in the field of optical pulse propagation in dispersive media. As we have
seen before, an equation like this can be solved as an ordinary initial value
problem if we first use (\ref{Eq8.99}) to make the term containing a second
derivative\ with respect to time into one containing only a first derivative with respect to time.
\begin{equation}
\partial_{t}A+\mathbf{v}_{g}\cdot\nabla A-i\beta\nabla^{2}A+i\alpha(\mathbf{v}_{g}\cdot\nabla)^2
A-i\gamma|A|^{2}A=0.\label{Eq8.100}
\end{equation}
The amplitude $A$ determines the electric and magnetic fields through the
identities
\begin{align}
\mathbf{E}(\mathbf{x},t)  & \approx\{(\omega A+i(\frac{\omega}{k}
-\omega^{\prime})\mathbf{u}\cdot\nabla A\mathbf{)q}\nonumber\\
& +i(\frac{\omega}{k}\mathbf{q}\cdot\nabla A)\mathbf{u}\}e^{i(\mathbf{k}
\cdot\mathbf{x}-\omega t)}+(\ast)\nonumber\\
\mathbf{B}(x,t)  & \approx\{kA\mathbf{t}+i(\mathbf{t}\cdot\nabla
A)\mathbf{u}\}e^{i(\mathbf{k}\cdot\mathbf{x}-\omega t)}+(\ast).\label{Eq8.101}
\end{align}
The equations (\ref{Eq8.100}) and (\ref{Eq8.101}) are the key elements in a
fast numerical scheme for linearly polarized wave packet solutions to
Maxwell's equations. Wave packets of circular polarization or arbitrary
polarization can be treated in an entirely similar manner, as can sums of
different polarized wave packets and materials with nontrivial magnetic response.

The derivation of the nonlinear Scrodinger equation for linearly polarized
wave packets I have given in this section, is certainly not the simplest one possible.
However, the  aim in this section has been to illustrate
how to apply the multiple scale method to vector PDEs in general, not to do
it in the most effective way possible, for the particular case of linearly
polarized electromagnetic wave packets in non-magnetic materials. If the material
has a significant magnetic response, a derivation along the lines given is necessary.

All the essential elements we need in order to apply the method of multiple
scales to problems in optics and laser physics, and other areas of science
too, are at this point known. There are no new tricks to learn. Using the
approach described in these lecture notes, amplitude equations can be derived
for most situations of interest. Applying the method is mechanical, but for
realistic systems there can easily be a large amount of algebra involved. This
is unavoidable; solving nonlinear partial differential equations, even
approximately, is hard.

In these lecture notes we have focused on applications of the multiple scale
method for time-propagation problems. The method was originally developed for
these kind of problems and the mechanics of the method is most transparent for
such problems. However the method is by no means limited to time propagation problems.

Many pulse propagation schemes are most naturally formulated as a boundary
value problem where the propagation variable is a space variable. A very
general scheme of this type is the well known UPPE\cite{UPPE} propagation
scheme. More details on how the multiple scale method is applied for these
kind of schemes can be found in \cite{per1} and \cite{newell}.

\setcounter{equation}{0}
\section{Appendix B}

\subsection{The maximum entropy principle for classical systems}

Let $x_1, ... , x_n$ be random variables with an associated probability distribution $\rho(x_1,...,x_n)$. Let $f_1(x_1,...,x_n),...,f_p(x_1,...,x_n)$ be functions defined on the space of random variables $\Omega=\{(x_1,x_2,...,x_n)\}$ where the variables $x_n$ can run over a finite set, an infinite discrete set, for example a set indexed by a finite set of integers, or the variables can run over the real numbers. We will usually think about the real number case and will therefore write integrals instead of sums.
The functions $f_j$ are our \ttx{observables}. Their $\it{expectation}$ values are as usual defined by 
\begin{align}
\left<f_j\right> \; = \int_{\mathbf{R}^n} dV \; f_j(x_1,...,x_n) \; \rho(x_1,...,x_n). \lbl{75} 
\end{align}
The expectation value of a given observable of course depends on which probability distribution, $\rho$, we use. The challenge in statistics is to figure out which probability distribution one should use in any given situation. Let us say that we for some reason, (expert knowledge,guesswork, hearsay, ...) believe that a probability distribution $\rho_0$  accurately represents what we currently know about a given system. The probability distribution $\rho_0$ is called the \textit{prior distribution}, or just the \textit{prior}.

Let us next assume that we measure the mean values of the observables $f_1, ...,f_p$ and find the values $c_1,...,c_p$. If 
\begin{align}
\left<f_j\right>_0 \;  = \int_{\mathbf{R}^n} dV \; f_j \; (x_1,...,x_n) \; \rho_0 (x_1,...,x_n) = c_j, \lbl{76} 
\end{align}
we are satisfied with our choice of prior. It predicts exactly the mean values that are observed. \\
But we might not be so lucky. Perhaps 
\begin{align}
\left<f_j\right>_0 \; \ne c_j, \lbl{77} 
\end{align}
for at least one $j$. Our selected $\rho_0$ is then not the correct one, it predicts expectation values that are not observed. The challenge is to modify $\rho_0$ into a new distribution $\rho$ that is consistent with \ttx{all} the observed mean values. 

For this purpose we define a functional $S(\rho)$ by 
\begin{align}
S(\rho) = - \int_{\mathbf{R}^n} dV \; \rho \; \ln (\frac{\rho}{\rho_0}). \lbl{78} 
\end{align}
$S$ is by definition the \ttx{relative entropy} of the probability distribution $\rho$ with respect  to $\rho_0$. We will see later that our use of the word entropy here is consistent with its usage in thermodynamics. 

The \ttx{maximum entropy} principle states that one should choose the probability distribution that maximizes the functional 
\begin{align}
S(\rho) = - \int_{\mathbf{R}^n} dV \; \rho \; \ln (\frac{\rho}{\rho_0}), \lbl{79} 
\end{align}
subject to the constraints 
\begin{align}
\left<f_j\right> \; = \int_{\mathbf{R}^n} dV \; f_j \; \rho = c_j,\;\; j=1,2,...,p. \lbl{80}
\end{align}
\\
\subsubsection{The general thermodynamical formalism}
In this section we will solve the maximum principle stated in the previous section using the calculus of variations. The problem will initially be solved in the general setting described in the previous section,  but we will eventually specialize to the case of statistical mechanics. 

In order to proceed we must first recognize that in additional to the $p$ constraints \rf{80}, we have one more constraint that simply expresses the fact that $\rho$ is a probability distribution.
\begin{align}
\left<1\right> \; = \int_{\mathbf{R}^n} dV \; \rho (x_1,...,x_n) = 1, \lbl{206} 
\end{align}
and we thus have $p+1$ constraints and therefore introduce an extended functional 
\begin{align}
T(\rho) = S(\rho) - \lambda_0 \left<1\right> - \mathlarger{\sum}^p_{j=1} \lambda_j \; \left<f_j\right>, \lbl{207} 
\end{align}
Note that we could have written 
\begin{align}
T (\rho) = S(\rho) - \lambda_0 \; ( \left<1\right> - 1) - \mathlarger{\sum}^p_{j=1} \lambda_j \; ( \left<f_j\right> - c_j), \lbl{208} 
\end{align}
in order to make the values of the constraints explicit. However, all constant terms vanish when we take variational derivative, so we might as well drop the constant terms. Also note that our choice of minus sign in front of the Lagrange multiplier terms in \rf{207} and \rf{208} is a convention inspired by the application of this formalism to the case of statistical mechanics.

\noindent The integral density corresponding to the extended functional $T (\rho)$ is 
\begin{align}
\mathcal{L} = - \rho \; \ln (\frac{\rho}{\rho_0}) - \lambda_0 \; \rho - \mathlarger{\sum}^p_{j=1} \lambda_j \; f_j \; \rho. \lbl{209} 
\end{align}
Observe that $\mathcal{L}$ does not depend on any derivatives of $\rho$. The Euler-Lagrange equation for $T$ is therefore simply 
\begin{align}
\frac{\partial \mathcal{L}}{\partial \rho} &= 0, \lbl{210} \\ 
&\Updownarrow \nonumber \\
- \ln (\frac{\rho}{\rho_0}) - 1 - \lambda_0 - \mathlarger{\sum}^p_{j=1} \lambda_j f_j&=0 , \nonumber
\end{align}
whose solution is
\begin{align}
\rho &= \frac{\rho_0}{Z} \; \exp{-\mathlarger{\sum}_j \lambda_j \; f_j}, \nonumber 
\end{align}
where we have defined $Z=\exp{(1+ \lambda_0 )}$. In order for the constraint $\left<1\right> = 1$ to be satisfied, we must have 
\begin{align}
\left<1\right> &= 1, \nonumber \\
&\Updownarrow \nonumber \\  \intRn dV \; \frac{\rho_0}{Z} \; \exp{-\mathlarger{\sum}_j \lambda_j \; f_j} &= 1, \nonumber \\ 
&\Updownarrow \nonumber \\ Z = Z(\lambda_1,...,\lambda_p) &= \intRn dV \; \rho_0 \; \exp{-\mathlarger{\sum}_j \; \lambda_j \; f_j}, \lbl{211} 
\end{align}
and the stationary distribution is 
\begin{align}
\rho(x_1,...,x_n) = \frac{\rho_0(x_1,...,x_n)}{Z(\lambda_1,...,\lambda_n)} \; \exp{-\mathlarger{\sum}^p_{j=1} \lambda_j \; f_j (x_1, ...,x_n)}. \lbl{212} 
\end{align}
$\rho$ is called the \ttx{maximum entropy distribution} and $Z$ is the $partition function$. Note that we have not proved that the distribution \rf{212} in fact gives a maximum value  for $S$, but this can be done\cite{Jaynes4}. \\
The Lagrange multipliers $\lambda_1,...,\lambda_p$ are chosen so that all the constraints are satisfied 
\begin{align}
\left<f_j\right> \; = \intRn dV \; f_j (x_1,...,x_n) \; \rho(x_1,...,x_n) = c_j && j=1,...,p\;. \lbl{213} 
\end{align}
The system of equations \rf{213} consists of  $p$ equations for the $p$ unknown quantities $\lambda_j$. 

As it turns out, we almost never need to know the distribution $\rho$ from \rf{212}, it is enough to know the partition function. Observe that 
\begin{align}
\left<f_j\right> \; &= \intRn dV \; f_j \; \rho\nonumber\\ 
&= \inv{Z} \; \intRn dV \; f_j \; \rho_0 \; \exp(-\mathlarger{\sum}^p_{i=1} \lambda_{i} f_{i}) \nonumber \\ 
&= -\inv{Z} \; \intRn dV \; \partial_{\lambda_j}\{ \rho_0 \; \exp(-\mathlarger{\sum}^p_{i=1} \lambda_{i} f_{i})\} \nonumber \\ 
&= -\inv{Z} \; \partial_{\lambda_j}  \intRn dV \;  \rho_0 \; \exp(-\mathlarger{\sum}^p_{i=1} \lambda_{i} f_{i}) \nonumber \\
&= -\inv{Z} \; \prt{\lambda_j} Z =- \prt{\lambda_j} \ln Z \nonumber \\
&\Downarrow\nonumber\\
\left<f_j\right> \; &=- \prt{\lambda_j} \ln Z,\lbl{214} 
\end{align}
and thus we can find the mean of all the quantities $f_j$ by taking partial derivatives of the natural logarithm of the partition function with respect to the Lagrangian multipliers. Moreover, we also have 
\begin{align}
\prt{\lambda_j \lambda_k}\ln Z &= \prt{\lambda_j} (\inv{Z} \; \prt{\lambda_k} Z) \lbl{215} \\
&= - \inv{Z^2} \; \prt{\lambda_j} Z \; \prt{\lambda_k} Z + \inv{Z} \; \prt{\lambda_j \lambda_k} Z \nonumber \\
&= - \prt{\lambda_j} \ln Z \; \prt{\lambda_k} \ln Z + \inv{Z} \; \intRn dV \; f_j \; f_k \; \rho_0 \; \exp(-\mathlarger{\sum}_i \lambda_i \; f_i) \nonumber \\
&= - \prt{\lambda_j}\ln Z \; \prt{\lambda_k} \ln Z + \left<f_j \; f_k\right>. \nonumber 
\end{align}
Thus 
\begin{align}
\left<f_j \; f_k\right> \; = \prt{\lambda_j} \ln Z \; \prt{\lambda_k} \; \ln Z + \prt{\lambda_j \; \lambda_k} \ln Z \lbl{216} 
\end{align}
In a similar way \ttx{all} correlation coefficients $\left<f_1^{n_1} ... f_p^{n_p}\right>$ can be expressed through partial derivatives of the partition function.

Inserting the maximum entropy distribution \rf{212} into the entropy functional \rf{79} gives us the following expression for the maximal value of the entropy
\begin{equation}
S=\ln Z+\mathlarger{\sum}_j \; \lambda_j \;\left< f_j\right>.\lbl{MaxEntropy}
\end{equation}
From a mathematical point of view we now have two sets of variables $\{\left< f_1\right>,...,\left< f_p\right>\}$ and $\{\lambda_1,...,\lambda_p\}$. Geometrically we imagine that these two pairs of variables, together with $S$, defines a space $\Omega$  of odd dimension $2p+1$ with coordinates $\{S,\left< f_1\right>,...,\left< f_p\right>,\lambda_1,...,\lambda_p\}$. The $p$ identities \rf{213} defines a $p+1$ dimensional surface $\Lambda$ in $\Omega$.

\noindent Taking the differential of the identity \rf{MaxEntropy} we get
\begin{equation*}
dS=\mathlarger{\sum}_j\; \frac{\partial\ln Z}{\partial\lambda_j}d\lambda_j+\mathlarger{\sum}_j \;\{\left< f_j\right>\;d\lambda_j+\lambda_jd\left< f_j\right>\}.
\end{equation*}
Restricting this differential to the surface $\Lambda$, and thus using the identities \rf{214}, gives us the following expression for the differential $dS$ restricted to the surface $\Lambda$
\begin{equation}
dS=\mathlarger{\sum}_j \;\lambda_jd\left< f_j\right>\lbl{dS1}.
\end{equation}
The identity \rf{MaxEntropy} defines the entropy as a function depending on all $2p$ variables in $\Omega$. We therefore have
\begin{equation}
dS=\mathlarger{\sum}_j\;\frac{\partial S}{\partial\lambda_j}d\lambda_j+\mathlarger{\sum}_j\;\frac{\partial S}{\partial\left< f_j\right>}d\left< f_j\right>\lbl{dS2}.
\end{equation}
Comparing \rf{dS1} and \rf{dS2} we conclude that on the surface $\Lambda$ we must have the identities 
\begin{align}
\frac{\partial S}{\partial\lambda_j}&=0,\nonumber\\
\frac{\partial S}{\partial{\left< f_j\right>}}&=\lambda_j.\lbl{SRelation1}
\end{align}
\noindent
Thus, on the surface $\Lambda$, the entropy depends only on the variables $\{\left< f_1\right>,...,\left< f_p\right>\}$ and the derivative with respect to these variables determines the values of the Lagrange multipliers in terms of the data $\{c_1,...,c_p\}$ of the problem.

It is frequently the case that in addition to the variables $\{x_1,...,x_n\}$, the observables depends on parameters. For notational simplicity, let us assume that there is only one parameter denoted by the symbol $\alpha$. Thus we have observables $\{f_1(x_1,...x_n;\alpha),...,f_p(x_1,...x_n;\alpha)\}$. The presence of the parameter does not change the argument leading up to the maximum entropy distribution \rf{212} and thus we have the formulas
\begin{align}
\rho(x_1,...,x_n;\alpha)& = \frac{\rho_0(x_1,...,x_n)}{Z(\lambda_1,...,\lambda_p;\alpha)} \; \exp{-\mathlarger{\sum}^p_{j=1} \lambda_j \; f_j (x_1, ...,x_n;\alpha)},\nonumber\\
 Z(\lambda_1,...,\lambda_p;\alpha) &= \intRn dV \; \rho_0 \; \exp{-\mathlarger{\sum}_j \; \lambda_j \; f_j(x_1,...,x_n;\alpha)}.\lbl{parameterZ}
\end{align}
 \noindent Differentiation of the partition function \rf{parameterZ} with respect to the parameter $\alpha$ gives us the expression
 \begin{align}
 \frac{\partial Z}{\partial \alpha}&=-\intRn dV \; \rho_0 \;{\mathlarger{\sum}_j \; \lambda_j \;\frac{\partial  f_j}{\partial \alpha}} \exp{-\mathlarger{\sum}_j \; \lambda_j \; f_j},\nonumber\\
 &\Downarrow\nonumber\\
  \frac{\partial \ln Z}{\partial \alpha}&=-\mathlarger{\sum}_j \; \lambda_j \;\left<\frac{\partial  f_j}{\partial \alpha}\right>.\lbl{DParameterZ}
 \end{align}
 If we repeat the calculation leading from \rf{MaxEntropy} to \rf{dS1} for the case when the observables depends on a parameter $\alpha$, we now get instead of \rf{dS1} the following more general expression for the differential of the entropy
 \begin{equation}
dS=\mathlarger{\sum}_j \;\lambda_jd\left< f_j\right>-\mathlarger{\sum}_j \;\lambda_j \;\left<\frac{\partial  f_j}{\partial \alpha}\right>d\alpha\lbl{dS3},
\end{equation}
where we have used the identity \rf{DParameterZ}.
 
 Note that this differential identity can be written in the form
  \begin{align}
dS=\mathlarger{\sum}_j \;\lambda_jdQ_j\lbl{dS3},
\end{align}
where we have introduced the quantities $dQ_j$ representing {\it generalized heat} associated with the observables
\begin{align}
dQ_j&=d\left< f_j\right>-\left<\frac{\partial  f_j}{\partial \alpha}\right>d\alpha\nonumber\\
&=d\left< f_j\right>-\left<\frac{\partial  f_j}{\partial \alpha}d\alpha\right>\nonumber\\
&=d\left< f_j\right>-\left<df_j\right>.\lbl{GeneralizedHeat}
\end{align}
Formula \rf{GeneralizedHeat} tells us what heat actually represents. Physical systems on the human scale, these are evidently the ones of most immediate interest to us, consists of an immense number of elementary subsystems. The detailed configurationl variables for all these elementary systems defines the {\it  microscopic} degrees of freedom of the human scale system. The state of these microscopic degrees of freedom are unknown to us and our ability manipulate then directly is entirely lacking. The few degrees of freedom of the system whose state we {\it  can} know and which we have the means to manipulate defines the {\it macroscopic} degrees of freedom for the system. In our description of thermodynamics these are the observables $f_j$. A change in the  mean value of a macroscopic degrees of freedom,$d\left<f_j\right>$, comes from two sources. The first source is a change in the observable representing the said macroscopic degree of freedom, this is the kind of change that we have the ability to induce by direct manipulation. This quantity is represented by $\left<df_j\right>$ in formula \rf{GeneralizedHeat}. When this quantity is subtracted from $d\left<f_j\right>$ , what remains is the second source of change of the mean. This second source  is a change in the underlying probability distribution which represents a change in our information about the microscopic state of the system. When our ignorance about the microscopic state of a system increase the system grows ``hotter'', corresponding to an increase in $dQ_j$.

As is usual in thermodynamics, the formalism is misleading in the sense that $dQ_j$ merely denote an infinitesimal amount of generalized heat and is {\it not} the differential of some function $Q_j$. No such function exists. The proper mathematical way to think about the identity \rf{dS3} is that $dS$ and $dQ_j$ are differential forms where $dS$ is an exact differential forms, meaning it is the differential of a function, and $dQ_j$ are inexact differential forms and thus not the differential of a function. However, the mathematical formalism of differential forms must be  introduced in the very large context of differential geometry and we will not digress into this area of mathematics.

The above explanation of the nature of heat,  referred to the original application of the thermodynamical formalism, where the systems has an immense number of microscopic degrees of freedom,  which are in principle knowable and controllable, but  as practical matter, not. We however know that the thermodynamical formalism can be applied to any situation where systems has more degrees of freedom than the ones we chose to observe. This might be because the underlying degrees of freedom are unknown but it could also be the case that they are known but that we for various reasons choose to ignore them. In both cases the argument above stands and the existence of the unknown or ignored degrees of freedom manifest as heat in the theory.

 \noindent We will now derive a generalized version of  identity \rf{DParameterZ} that plays a crucial role when the thermodynamical formalism is applied to the special case for which the underlying space is the state space of a physical system. The system could be a classical mechanical system consisting of a finite number of particles, a system of classical  fields or even the Fock state space for a quantum mechanical many particle system.
 
 In all these cases, one consider systems that are confined to a bounded spatial domain $D$ which is defined by its bounding surface $\Gamma$. Thus all observables for the system will typically depend on the bounding surface $\Gamma$, $f_j=f_j(x_1,...x_n;\Gamma)$. We will now consider a small deformation, $\delta\Gamma$, of the bounding surface $\Gamma$. Thus $\Gamma\longrightarrow\Gamma+\delta\Gamma$. This deformation leads to variations
 \begin{align}
 \delta_\Gamma f_j(x_1,...x_n;\Gamma)=f_j(x_1,...x_n;\Gamma+\delta\Gamma)-f_j(x_1,...x_n;\Gamma),\nonumber\\
 \delta_\Gamma Z(\lambda_1,...,\lambda_p; \Gamma)=Z(\lambda_1,...,\lambda_p; \Gamma+\delta\Gamma)-Z(\lambda_1,...,\lambda_p; \Gamma).\lbl{variations}
 \end{align}
 Arguing exactly like we did for the simple case of a single parameter we now find the important identity
 \begin{align}
 \delta_\Gamma \ln Z=-\mathlarger{\sum}_j \; \lambda_j \;\left< \delta_\Gamma f_j\right>.\lbl{FundamentalVariationalIdentity}
  \end{align}
  \noindent This identity will, for the special cases mentioned above, lead to the definition of the thermodynamic pressure and related quantities. 
  Corresponding to the {\it differential} identity for the entropy \rf{dS3} we now get the following more general {\it variational} identity
\begin{align}
\delta S=\mathlarger{\sum}_j \;\lambda_jd\left< f_j\right>-\mathlarger{\sum}_j \;\lambda_j \;\left< \delta_\Gamma f_j\right>\lbl{dS5}.
\end{align}

\subsubsection{The thermodynamic formalism in statistical physics} 
 
   Let us now consider the special case when our underlying space is the classical state space for a mechanical system with $n$ degrees of freedom. This could for example consist of $n$ mass points. We will assume that the system is confined to a bounded domain $D$ in $\mathbf{R}^3$ defined by a bounding surface $\Gamma$. The state space is thus a subset of the euclidean space $\mathbf{R}^{6n}$ with coordinates $(\vb{q},\vb{p})=(\vb{q}_1,...,\vb{q}_n,\vb{p}_1,...,\vb{p}_n)$,  since we need 3 position coordinates and 3 momentum coordinates for each particle in order to uniquely specify the state of the system. \\
Let $\mathcal{H} = \mathcal{H}(\vb{q},\vb{p})$ be the Hamiltonian for the system of mass  points. Recall that the value of the Hamiltonian on any given state, $(\vb{q},\vb{p})$, is the energy, $E$,  of that state. \\
When $n$ is large it is very hard, and also mostly useless, to try to track the exact state $(\vb{q}(t),\vb{p}(t))$ of a system of mass points. \\
For such a large system it is more useful to consider a probability distribution $\rho(\vb{q},\vb{p})$ on the state-space. This is the point of view introduced by Gibbs. 
We will first consider the simplest, and by far the most common situation, where the Hamiltonian, $\mathcal {H}=\mathcal {H}(\vb{q},\vb{p})$ is the only observable. The maximum entropy distribution for this case is
\begin{align}
\rho (\vb{q},\vb{p}) = \frac{\rho_0(\vb{q},\vb{p})}{Z} \; \exp(-\frac{\mathcal {H}(\vb{q},\vb{p})}{kT}), \lbl{217.1} 
\end{align}
where the partition function is given by
\begin{align}
Z=Z(T)=\int_{\mathcal{R}^{6n}}d\vb{q}d\vb{p}\;\rho_0(\vb{q},\vb{p})\exp(-\frac{\mathcal {H}(\vb{q},\vb{p})}{kT}),\lbl{statmechZ1}
\end{align}
and where we have redefined the single Lagrange multiplier using
\begin{align}
\lambda = \inv{k \; T}. \lbl{218} 
\end{align}
In this formula, $k$ is the Boltzmann constant and $T$ is a new parameter which by definition is the thermodynamic temperature.
The parameter $T$ is determined by 
\begin{align}
 E&= \left<\mathcal{H}\right>, \nonumber \\
&\Updownarrow \nonumber \\ 
E&=k \; T^2 \; \prt{T} \ln Z ,\lbl{220}
\end{align}
where we have used the chain rule 
\begin{align}
\prt{\lambda} =- k \; T^2 \; \prt{T}, \lbl{221} 
\end{align}
in the general formula \rf{214}.

  Formula \rf{220} is in statistical mechanics and thermodynamics called the {\it equation of state}, and all thermodynamic statements that can be made about the system flows from this formula. The formula for the equation of state may look innocent, in order to find it you merely need to take the derivative of the partition function, and partition function also looks innocent, after all it is just a function of one variable, the kind of function we study in first year calculus. However, in order to actually find an expression for this single variable function one needs to do the integral in formula \rf{statmechZ1}, and this is a multiple integral involving something like $10^{27}$ integration variables in typical situations! Clearly, an exact formula for the partition function can rarely be found. Approximate expressions where the large number of particles are used to ones advantage can more frequently be found, but pushing through calculations like these are as a rule extremely technical. More than one Nobel price has been handed out for developing feasible schemes for calculating the partition function. Given the level of complexity involved in calculating the partition function from the defining formula \rf{statmechZ1}, and the fact that the partition function simply is a function of one or a few variables,  it should come as no surprise that the most common approach to finding the equation of state is to fit parametrized functions to experimental data.

The maximum entropy distribution \rf{217.1} is recognized to essentially be the \ttx{Gibb's Canonical ensemble} from statistical physics.

The Gibb's ensemble is the foundation of statistical physics. All results in statistical physics flows from formula \rf{217.1}. Statistical physics is also the foundation of thermodynamics so all conclusions from that subject also flow from the Gibb's ensemble \rf{217.1}. In the thermodynamics context,  \rf{220} is, as we have already remarked, nothing but the \ttx{equation of state}. 

An interesting insight here is that the temperature of a thermodynamic system is in fact a Lagrange multiplier!! This is a profound insight that to this day has not been fully understood or explored.

From this example, it appears  useful to think of any application of the maximal entropy principle as an extension of the methods of statistical mechanics to systems that has absolutely nothing to do with the motion of mass points. 

This wide general applicability of the methods of statistical physics has  lead to deep questions and insights into the nature and significance of the assumption of equilibrium that appears to underline the application of the Gibb's ensemble in statistical physics. 

There is also the intriguing fact that the very same functional \rf{79} used in the maximum entropy principle,  is also the foundation of information theory which was discovered by Shannon in 1948. This connection between information theory and statistical mechanics (and thermodynamics) has lead to deep insights into the role of information in our fundamental physical theories. \\
As already discussed in the introduction, the general nature and wide applicability of the maximum entropy principle has been described well by E.T. Jaynes in many papers and the (unfinished) monumental book "Probability theory: The Logic of Science". 

\noindent As if all this is not impressive enough for one single principle, it is also a very intriguing fact that when one looks deep into the heart of fundamental physics, in the form of quantum field theory, one again finds an appropriately generalized form of  the partition function \rf{211}. The whole computational engine in the theory of quantum fields revolve around a generalized Gibb's ensemble! 

What on earth is going on...

\subsubsection{The problem of prior}

Note that formula \rf{217.1} does not uniquely define Gibbs ensemble because of the presence of the prior distribution $\rho_0$. The actual Gibbs ensemble corresponds to the choice $\rho_0=1$. When using the information theoretical approach to statistical mechanics and thermodynamics, like we do here, one should be very wary when it comes to the choice of the prior distribution. It is simply the most contentious  issue in the whole theory. One should ask pointed questions of justification for any proposed choice. What kind of information about the system is it based on, and is it the correct embodiment of said information? 

In fact, if one study expositions of statistical mechanics and thermodynamics, which are based on the traditional objective dynamical approach to the subject,  one finds that the choice of what from the information theoretical point of view is the prior distribution, is much discussed. The reasons for choosing $\rho_0=1$ that have appeared through these discussions are, in our humble opinion, not very convincing.

The problem of determining the prior distribution has been at the center of probability theory
and statistics from the very start. The general rules of
probability theory tells us how to compute probabilities for
derived events from probabilities of primary events. The problem
of prior is concerned with the problem of assigning probabilities
to primary events. The assignment is supposed to reflect an
observers state of knowledge about the primary events. The
assignment should be the same for different observers with the
same state of knowledge but can be different for observers with
different states of knowledge \cite{Jaynes3}. In this sense
probability assignments are subjective
 \cite{jeffreys},\cite{cox1},\cite{cox2}. The problem of the prior
is how to turn states of knowledge into probability assignments.
The first solution to this problem was used by the very founders
of probability theory (Bernoulli and Laplace). If the observers
only knowledge of the primary events are their number, then a
uniform probability assignment should be used. This idea was later
named the principle of indifference by J. M. Keynes. Generalizing
this idea to countably infinite or even continuous spaces of
primary events has turned out to be very problematic. Laplace
himself used such a generalization is his work on probability
theory. His probability distribution was uniform and not
normalizable since it was defined on the whole real line. Using a
uniform distribution for representing indifference about a random
variable on a finite interval on the real line would seem to be
more reasonable, at least it is normalizable. However even in this
case serious problems arise as the well known Bertrand's paradox
shows. Problems and paradoxes arising from the various
generalizations of the principle of indifference to continuous
random variables played no small part in the creation
and for a long time complete dominance of the frequency interpretation\cite%
{Fisher} of probability theory.

The principle of maximum entropy appears first in the writings of
W. Gibbs  \cite{Gibbs} on thermodynamics and statistical physics
and later in the fundamental work on information theory by Shannon
\cite{Shannon}. However it was E. T. Jaynes \cite{Jaynes1} who
realized the real importance and general nature of the principle
of maximum entropy. In his hands it turned into a general method
for turning prior knowledge in the form of mean values for observables defined on finite state spaces, into prior probability assignments.

 Let us consider this simplest case in more detail. Let $\Omega =\{x_{1},x_{2},....,x_{n}\}$ be a finite
space of primary events. The algebra of possible events is the set of all subsets of $
\Omega$. A probability assignment on the set of primary events is a finite set of
numbers $p=\{p_{i}\}$ such that $0\leq p_{i}\leq 1$ and $\sum_{i=1}^{n}p_{i}=1$%
. Let $f_{1},...,f_{k}$ be real valued functions on $\Omega $. The
principle of maximum entropy states that if the means of the
functions $f_{1},..,f_{k}$ are known, $\left<f_{i}\right>=c_{i}$, one should,
among all probability assignments that satisfy the constraints,
pick the one that maximizes the entropy $S=-\sum_{i=1}^{n}p_{i}\ln
p_{i}$. The solution to this constrained maximization problem is, as we have seen, the maximum entropy distribution
\begin{align}
p=\frac{1}{Z(\lambda _{1},...,\lambda
_{k})}\exp\left(-\sum_{j=1}^{k}\lambda _{j}f_{j}\right),
\end{align}
where $Z$ is the partition function and is given by
\begin{align}
Z(\lambda _{1},...,\lambda _{k})=\sum_{i=1}^{n}\exp\left(-\sum_{j=1}^{k}\lambda _{j}f_{j}(i)\right).
\end{align}
Observe that for the particular situation where there are no constraints,
the principle gives $Z=n$ and the maximum entropy distribution is uniform
\begin{align}
p_{i}=\frac{1}{n}.
\end{align}
Thus, for an observer that only know that there are $n$ possible primary events, the maximum entropy distribution is exactly the one suggested by the principle of indifferent! The conclusion appeared to be that not only could the maximum entropy distribution tell us how to choose the best distribution in the presence of observed means of a finite number of observables, it could also tell us which distribution to choose when our ignorance is so profound that the only thing we know about a situation is the number of possible primary events. This distribution is of course exactly what we have called the prior distribution.  For
a time it looked as if the problem of prior was essentially
solved.  However continuous valued random variables again turned
out to be the Achilles heel. For finite spaces of events the
principle will give a unique probability assignment, but when
generalizing it to continuous random variables by taking a continuum limit of the finite discrete expression for the entropy, an unknown
probability density appears. The density appears because the continuum limit is not unique.
Different limiting expressions are found depending on how one approach the continuum through a countable set of discrete spaces.
The unknown probability distributions that appears essentially depends on how the discrete points bunch up in the limit. The meaning of this probability density became
clear when it was realized that it is the maximum entropy
distribution corresponding to no constraints. Thus it was
understood that in order to apply the principle of maximum entropy
one must start with a prior distribution. The principle of maximum
entropy could not determine the prior, it could only tell us how
to modify an already existing prior in order to satisfy
constraints in the form of mean values. It seemed as if one were
back to square one.

There {\it does} however exist  a systematic way to turn prior information on means of observables into prior
distributions, and it {\it does} involve the maximum entropy principle,  but not in the direct way just described.
In fact, after a certain reformulation  it will become evident that the problem of
 selecting a prior is not merely a side issue that has to be resolved in order to proceed with the real work of applying
  the maximum entropy principle, the problem of prior is the {\it only} issue as far as the maximum entropy principle is concerned.

In order to describe this reformulation of the principle of maximum entropy, we will return to the special case of statistical mechanics.
In the previous section we discussed the problem of specifying the prior in the context of statistical mechanics and expressed our doubt as to the justifications for making the standard choice $\rho_0(\vb{q},\vb{p})=1$. Even if we are doubtful about the justification for this particular choice,  it is clear that when we apply the maximum entropy principle  in statistical mechanics there is a physical context that certainly makes some choices of the prior more reasonable than othersr. By picking the Hamiltonian function as our observable we must also acknowledge that the system evolve according to the corresponding Hamiltonian equations. It is always the case that the the Hamiltonian function, $\mathcal{H}(\vb{q},\vb{p})$, which represents the energy, is a constant of the motion. Depending on the symmetries of the interaction,  Hamiltonian systems of equations may also have other conserved quantities. The generic situation is however that the energy is the only conserved quantity. We will assume that this is the case and let the corresponding Hamiltonian function $\mathcal{H}$ be our only observable. The maximum entropy distribution is now given by expression \rf{217} where $\rho_0$ is the prior distribution. It  is in the current context reasonable to impose the condition that the prior is a stationary solution to the corresponding Liouville  equation. But this means that the prior distribution, $\rho_0=\rho_0(\vb{q},\vb{p})$,  is a conserved quantity for  the Hamiltonian system and since the Hamiltonian is the only independent conserved quantity for our generic Hamiltonian system we must have
\begin{align}
\rho_0(\vb{q},\vb{p})=f_0(\mathcal{H}(\vb{q},\vb{p})),\nonumber
\end{align}
where $f_0$ is an arbitrary function defined on the positive real line.
Using this fact we have from \rf{statmechZ1}
\begin{align}
Z(T)&=\int_{\mathcal{R}^{6n}}d\vb{q}\;d\vb{p}\rho_0(\vb{q},\vb{p})\exp(-\frac{\mathcal {H}(\vb{q},\vb{p}}{kT})\nonumber\\
&=\int_{\mathcal{R}^{6n}}d\vb{q}\;d\vb{p}f_0(\mathcal{H}(\vb{q},\vb{p}))\exp(-\frac{\mathcal {H}(\vb{q},\vb{p}}{kT})\nonumber\\
&=\int_0^{\infty} dE\;\exp(-\frac{E}{kT})\;f_0(E)\;\int_{\mathcal{H}=E}d\vb{q}\;d\vb{p}\nonumber\\
&=\int_0^{\infty} dE\;\exp(-\frac{E}{kT})\;\rho_0(E),\lbl{statmechZ2}
\end{align}
where we have defined 
\begin{align}
\rho_0(E)=f_0(E)\;\int_{\mathcal{H}=E}d\vb{q}\;d\vb{p}.\lbl{macroPrior}
\end{align}
The constraints on the {\it microscopic} prior distribution $\rho_0(\vb{q},\vb{p})$ has reduced our original maximum entropy principle on the extremely high dimensional space $\mathbb{R}^{6n}$ with the Hamiltonian as our observable,  to a maximum entropy problem on the real line where the coordinate on the line, $E$, is the observable and the {\it macroscopic} prior is given by \rf{macroPrior}. The maximum entropy distribution for this case is
\begin{align}
\rho(E)=\frac{\rho_0(E)}{Z(T)}\exp(-\frac{E}{kT}),\lbl{MaxEntOnLine}
\end{align}
where the partition function is given by
\begin{align}
Z(T)=\int_0^{\infty}dE\;\rho_0(E)\exp(-\frac{E}{kT}).\lbl{StatMechReducedPartitionFunction}
\end{align}
This simple situation where we apply the entropy principle to a low dimensional state space $\mathbb{R}^p$ and where the observables are the 
coordinate functions, $x_1,...,x_p$ on the space is not special at all, in fact this is the most common situation when we apply the maximum entropy principle and  other applications can almost always be reduced to this situation using an approach similar to the reduction from $\mathbb{R}^{6n}$ to $\mathbb{R}$ described for the case of statistical mechanics.

In most applications of probability theory in statistics there is
no underlying high dimensional space of primary events, $\Omega $,  like in statistical mechanics and other areas of physics,  and the random
variables are not some functions,  like the Hamiltonian, defined on this space. 

Thus in the typical case one can assume that $\Omega =\mathbb{R}^{p}$, where $p$ is a fairly small number,  and that the random variables are just the coordinate function on $\mathbb{R}^{p}$. 
The prior probability distribution is then a function, $\rho_0=\rho_0(x_1,...,x_p)$, on $\mathbb{R}^p$, and  the partition function is given by the formula
\begin{align}
Z(\lambda _{1},..,\lambda _{p})=\int_{\mathbb{R}^p}dx_1dx_2...dx_p\;\rho_0(x_1,...,x_p)\exp(-\sum_{j=1}^{p}\lambda _{j}x_{j}).\lbl{FromRho0ToZ}
\end{align}
The partition function is thus nothing but the multi dimensional Laplace
transform of the prior distribution. This relation can be inverted,
using analytical continuation and the multidimensional Fourier transform on the imaginary $\lambda_j$ axes, and thereby expressing the prior in terms of the partition function
\begin{align}
\rho _{0}(x_{1},..x_{p})=\frac{1}{(2\pi )^{p}}\int_{\mathbb{R}^{p}}d\lambda _{1}..d\lambda _{p}Z(i\lambda _{1},..,i\lambda _{p})\exp(i\sum_{j=1}^{p}\lambda
_{j}x_{j}).\lbl{FromZToRho0}
\end{align}The whole content of the maximum entropy principle is contained in the integral transforms \rf{FromRho0ToZ} and \rf{FromZToRho0} connecting the partition function and the prior distribution. This is the promised reformulation of the maximum entropy principle, and we understand now that the prior distribution is not merely a bit player in this drama, it is the {\it only} player. In the next section we will show how this reformulation of the maximum entropy principle gives us a method for solving the problem of prior, which we now see is the only remaining fundamental problem in the statistical modeling of natural or artificial systems.

\subsubsection{Solving the problem of prior using stochastic relations}

In probability theory and statistics, random variables are often grouped into
statistical quantities. These are certain algebraic
combinations of means of functions of the random variables. A
large set of such statistical quantities are in use, some simple
examples are
\begin{align}
&\left<x\right>\;\;\;\;\;\;\;\;\;\;\;\;\;\;\;\;\;\;\;\;\;\;\;\;\;\;\;\;\;\;\;\;\;\;\;\;\;\;\;\;\;\;\;\;\text{ The mean of x.}\\
&\left<x^{2}\right>-\left<x\right>^{2}\;\;\;\;\;\;\;\;\;\;\;\;\;\;\;\;\;\;\;\;\;\;\;\;\;\;\;\;\;\;\; \text{ The variance of x.}\\
&\left<x^{3}\right>-3\left<x\right>\left<x^{2}\right>+2\left<x\right>^{3}\;\;\;\;\;\;\;\;\;\;\;\text{The third cumulant.} \\
&\left<xy\right>-\left<x\right>\left<y\right>\;\;\;\;\;\;\;\;\;\;\;\;\;\;\;\;\;\;\;\;\;\;\;\;\;\;\;\;\;\text{The cross variance of x and y.}
\end{align}
All such quantities can systematically be expressed as functions of the form
$F(q_{1},..,q_{k})$ where the variables $q_{j}$ are means of monomials in
the random variables. We will define {\it stochastic relations} to be systems of
equations for the quantities $q_{j}$.%
\begin{align}
F_{i}(q_{1},..,q_{k}) =0\text{ \ \ \ \ }i=1,...,s.\lbl{StochasticRelation}
\end{align}
Such relations are common in probability and statistics. Examples
are zero mean, fixed variance, uncorrelated variables and
identities expressing higher order cumulants in terms of lower
ones. Identities such as the last ones in the previous list  are the fundamental tools used to construct theories of turbulence
in fluid, gases and elsewhere. They are also, in their quantum incarnations,the key tools used to find viable simplified models in solid state physics and material science.

  In the previous section we have seen that the maximum
entropy principle defines a Laplace transform that map the prior
distribution to a partition function. As a direct consequence of this
transformation, we can express means of monomials in the random
variables in
terms of partial derivatives of the partition function. For example we have
\begin{align}
\left<f_i\right>&=-\frac{1}{Z}\partial_{\lambda_i}Z,\nonumber\\
\mathop{var}(f_i)&=\left<f_i^2\right>-\left<f_i\right>^2=-\frac{1}{Z^2}\partial_{\lambda}Z^2+\frac{1}{Z}\partial_{\lambda\lambda}Z.\nonumber
\end{align}

\noindent This means that the maximum entropy principle turns stochastic
relations into systems of partial differential equations for the
partition function and therefore imposes constraints on
the prior distribution.

    The problem is now how to describe the space of
solutions of these systems of partial differential equations. In
general, not all solutions to the equations can correspond to prior
probability distributions. From the definition of the partition function it is
for example clear that $Z(0)=1$ must hold for any acceptable solution.
 Finding necessary and sufficient conditions for functions to be the Laplace
transform of a probability distribution, and thus be acceptable solutions of
the systems of differential equations corresponding to stochastic relations,
 is not a simple matter, but
some results are known \cite {Rothaus}. We will not discuss this
problem but rather try to explicitly construct the solution space
or to say something useful about the structure of the solution space using
methods from the formal theory of differential equations .
Typically, the solution space is not a linear space and even when
it is, the dimension could easily be infinite. However, depending
on the number and types of stochastic relations the solution space
can end up being parametrized by a finite set of parameters or
even be a single point. In this last situation the stochastic
relations determine the prior uniquely. Note that in ordinary
(parametric) statistics finite parameter families of probability
distributions (Gaussian, Poisson, Bernoulli, t-distribution, etc)
are assumed to apply in given situations. From the point of view
discussed in these notes, this means that in ordinary statistics,
stochastic relations constrain the solution space enough for it to
be parameterized in terms of a finite number of parameters.
Nonparametric statistics correspond to the situation when the
solution space is so weakly constrained that it can not be
parameterized in terms of a finite number of parameters. Methods
from the theory of partial differential equations can in some
cases parameterize such weakly constrained solution spaces, not in
terms of real numbers, but in terms of arbitrary functions.
However for such weakly constrained solution spaces there is
another powerful tool available. This is the formal theory of
partial differential equations. The main object of study in this
theory is the infinitely prolonged hierarchy of a given systems of
differential equations. Thus one studies the infinite set of all
differential consequences of a given system of equations. Each
such differential consequence can be converted back into a
stochastic relation by using the relation between mean of
monomials and partial derivatives in reverse. One therefore gets
the corresponding hierarchy of stochastic relations that are
consequences of the original relations induced by the maximum
entropy principle and implemented through the Laplace transform.

In the remaining part of this subsection we will discuss several examples that
illustrate the method that has been outlined.

\paragraph{Stochastic relations for one random variable}
Essentially all families of distribution in use in parametric
statistics can be derived from simple stochastic relations
involving the mean, variance and skewness. In this section we show
some examples that support this statement.

\subparagraph{Delta distribution}
Let us consider the stochastic relation corresponding to a fixed
mean. It is
\begin{align}
\left<x\right>-q=0.
\end{align}
The Laplace transform convert this into the ordinary differential equation
\begin{align}
Z_{\lambda }=-qZ.
\end{align}
For this simple stochastic relation our system of partial differential
equations is a single linear ordinary differential equation. The solution
space is linear and parameterized by a single parameter
\begin{align}
Z(\lambda )=ae^{-q\lambda }.
\end{align}

\noindent The condition $Z(0)=1$ fixes the parameter $a$ to be one and we
have a unique solution. It is a simple matter to apply the inverse
transform \rf{FromZToRho0} to find the corresponding prior distribution
\begin{align}
\rho_0(x)& =\frac{1}{2\pi}\int_{-\infty}^{\infty}d\lambda Z(i\lambda)\exp(i\lambda x)\nonumber\\
&=\frac{1}{2\pi}\int_{-\infty}^{\infty}d\lambda\;\exp(-iq\lambda)\exp(i\lambda x)\nonumber\\
&=\frac{1}{2\pi}\int_{-\infty}^{\infty}d\lambda\;\exp(i(x-q)\lambda)\nonumber\\
&=\delta (x-q).
\end{align}

\subparagraph{Normal distribution}

The stochastic relation corresponding to constant variance  is
\begin{align}
\mathop{var}(x)=q,
\end{align}
and the corresponding differential equation is
\begin{align}
ZZ_{\lambda \lambda }-Z_{\lambda }^{2}-qZ^{2}=0.
\end{align}
This is a second order nonlinear ordinary differential equation.
The general solution of the nonlinear equation that satisfies the
requirement $Z(0)=1$ is
\begin{align}
Z(\lambda )=e^{-a\lambda +\frac{1}{2}q\lambda ^{2}}\text{, \ \ \ \
\ \ \ }a\in\mathbb{R}.
\end{align}
Using this partition function we can predict the mean of the random variable $x$ to be
\begin{align}
\left<x\right)&=-\frac{1}{Z(\lambda)}\frac{\partial Z(\lambda)}{\partial \lambda}\nonumber\\
&=a
\end{align}
and the corresponding prior distribution is found, using \rf{FromZToRho0}, to be

\begin{align}
\rho_0 (x)&=\frac{1}{2\pi }\int_{\mathbb{R}}d\lambda\;Z(i\lambda )\exp(i\lambda\;x)\nonumber\\
&=\frac{1}{2\pi}\int_{-\infty}^{\infty}d\lambda\;\exp(-ia\lambda  -\frac{1}{2}q\lambda ^{2})\exp(i\lambda\;x)\nonumber\\
&=\frac{1}{\sqrt{2\pi q}}e^{-\frac{(x-a)^{2}}{2q}}.
\end{align}
which is the normal distribution.

\subparagraph{Poisson distribution}

Let us consider the stochastic relation
\begin{align}
\mathop{var}(x)=<x>.
\end{align}

The corresponding differential equation is
\begin{align}
ZZ_{\lambda \lambda }-Z_{\lambda }^{2}+ZZ_{\lambda }=0.
\end{align}
This equation and most equations derived from stochastic relations simplify
considerably if we introduce a new function $\varphi $ through $Z=e^{\varphi
}$. The equation for $\varphi $ is
\begin{align}
\varphi _{\lambda \lambda }=-\varphi _{\lambda }.
\end{align}
This equation is easy to solve and the corresponding family of partition
functions satisfying, as always, the constraint $Z(0)=1$ is
\begin{align}
Z(\lambda )=e^{a(e^{-\lambda }-1)}.
\end{align}
The corresponding prior distribution is found using \rf{FromZToRho0} to be supported on $\Omega =\{0,1,2,....\}$ and
is of the form
\begin{align}
\rho_0 (k)=\frac{e^{-a}a^{k}}{k!}.
\end{align}
This is the Poisson distribution.

\subparagraph{Gamma distribution}

Let us consider a stochastic relation
\begin{align}
\mathop{var}(x)=\frac{1}{k}{\left<x\right>}^{2}\text{ \ \ }k>0.
\end{align}
Expressed in terms of $\varphi $ the corresponding differential equation is
\begin{align}
\varphi _{\lambda \lambda }=\frac{1}{k}\varphi _{\lambda }^{2}.
\end{align}
The general solution of this equation gives the following family of
partition functions
\begin{align}
Z(\lambda )=(1-a\lambda )^{-k}\text{ \ \ }a>0.
\end{align}
The corresponding prior  distribution is supported on $\Omega =(0,\infty )$ and is given by 
\begin{align}
\rho_0 (x)=x^{k-1}\frac{e^{-\frac{x}{a}}}{a^{k}\Gamma (k)}.
\end{align}
This is the Gamma distribution

\subparagraph{Bernoulli and Binomial distribution}

Let the variance be the following quadratic function of the mean
\begin{align}
\mathop{var}(x)=\left<x\right>(1-\left<x\right>).
\end{align}
The corresponding differential equation for $\varphi $ is
\begin{align}
\varphi _{\lambda \lambda }=-\varphi _{\lambda }(1+\varphi _{\lambda }).
\end{align}
The solution of the equation gives the \ following family of partition
functions
\begin{align}
Z(\lambda )=p+qe^{-\lambda }\text{ \ \ \ \ \ }p+q=1.
\end{align}
The corresponding distribution is supported on $\Omega =\{0,1\}$
and is given by $\rho (0)=p$, $\rho (1)=q$. This is the Bernoulli
distribution. If we generalize the stochastic relation to
\begin{align}
\mathop{var}(x)=\left<x\right>(1-\frac{1}{n}\left<x\right>).
\end{align}
where $n$ is a natural number we get the differential equation
\begin{align}
\varphi _{\lambda \lambda }=-\varphi _{\lambda }(1+\frac{1}{n}\varphi
_{\lambda }).
\end{align}
The partition function is found to be
\begin{align}
Z(\lambda )=(p+qe^{-\lambda }\text{ )}^{n}\text{\ \ \ \ \ }p+q=1.
\end{align}
The corresponding prior distribution is now found to be supported on $\Omega =\{0,1,...n\}$ and is
on this domain given by
\begin{align}
\rho_0 (k)=\binom{n}{k}p^{k}q^{n-k}.
\end{align}
This is the Binomial distribution.

\bigskip

\paragraph{Stochastic relations for more than one random variable}

When the number of random variables become larger than one,
stochastic relations in general leads to systems of nonlinear
partial differential equations. Unless the number and type of
relations is right, it is impossible to describe the solution
space in terms of a finite number of parameters. This lead us into
the domain of nonparametric statistics. This is the domain where
the methods from the formal theory of differential equations comes
into play. It is not possible to give nontrivial applications of
the theory here and  we will therefore limit ourselves to
two simple examples.

\subparagraph{The Multinomial distribution}

Let $x_{1},...x_{n}$ be $n$ random variables and consider the
following
system of stochastic relations
\begin{align}
\mathop{var}(x_{i}) &=\left<x_{i}\right>(1-\frac{1}{n}\left<x_{i}\right>)\text{ \ }
i=1,..n\;, \nonumber\\
\mathop{cov}(x_{i},x_{j}) &=-\frac{1}{n}\left<x_{i}\right>\left<x_{j}\right>\text{ \ \ }
i,j=1,...n,\text{ \ \thinspace }i\neq j.
\end{align}
The corresponding system of partial differential equations is%
\begin{align}
\varphi _{\lambda _{i}\lambda _{i}} &=-\varphi _{\lambda _{i}}(1+\frac{1}{n}
\varphi _{\lambda _{i}}), \nonumber\\
\varphi _{\lambda _{i}\lambda _{j}} &=-\frac{1}{n}\varphi _{\lambda
_{i}}\varphi _{\lambda _{j}}.
\end{align}
The second part of the  system of equations has general solutions of the
form $\varphi =n\ln (\theta )$ where $\theta (\lambda _{1},..,\lambda
_{n})=\sum_{i=1}^{n}\theta _{i}(\lambda _{i})$. \ Inserted into the first
part of the system this form of $\varphi $ easily gives the partition
function corresponding to the multinomial distribution. This system of
relations thus constrained the space of solutions so much that it could be
describes in terms of a finite number of parameters.

\subparagraph{Stochastic relations for the mean}

For a single random variable, stochastic relations involving only
the mean gives distributions located on a finite number of points.
For more than one random variable such relations gives rise to
nonparametric statistics, or solution spaces parameterized by
functions. The theory of partial differential equations can be
used to give a full description of these solution spaces. As an
example of such a relation consider the case of two
random variables whose means are constrained to be on a circle of radius $r$.
\begin{align}
{\left<x\right>}^{2}+{\left<y\right>}^{2}=r^{2}.
\end{align}
The corresponding partial differential equation is i terms of $\varphi $
\begin{align}
\varphi _{\lambda }^{2}+\varphi _{\mu }^{2}=r^{2},
\end{align}
and is known from optics as the {\it Eiconal} equation.The following $Z$ is in the solution space
\begin{align}
Z=e^{r\sqrt{\lambda ^{2}+\mu ^{2}}}.
\end{align}
This partition functions predicts that the following stochastic relation
should hold%
\begin{align}
\mathop{var}(x)=\left( \frac{\left<y\right>}{\left<x\right>}\right) ^{2}\mathop{var}(y).
\end{align}
The partial differential equation has, however, infinitely many solutions.
The method of characteristics can be used to describe the complete solution
space. In order to derive stochastic relations that holds for all $Z$ in the
solution space, these are the ones that can be said to be consequences of
the of the circle constrain, we should consider differential prolongations
of the original differential equation. The first prolongation is the system
\begin{align}
\varphi _{\lambda }^{2}+\varphi _{\mu }^{2} &=r^{2}, \\
\varphi _{\lambda }\varphi _{\lambda \lambda }+\varphi _{\mu }\varphi _{\mu
\lambda } &=0, \\
\varphi _{\lambda }\varphi _{\lambda \mu }+\varphi _{\mu }\varphi _{\mu \mu
} &=0,
\end{align}
and this system implies that%
\begin{align}
\varphi _{\lambda \lambda }=\left( \frac{\varphi _{\mu }}{\varphi _{\lambda }%
}\right) ^{2}\varphi _{\mu \mu }.
\end{align}

Translated into stochastic relations this is exactly the one we derived for
the special solution $\varphi =r\sqrt{\lambda ^{2}+\mu ^{2}}$ and it thus
holds for all solutions. It is of considerable interest to find a finite set
of basic stochastic relations that through some construction procedure
implies all consequences of some given system of stochastic relations. This
is exactly the kind of question addressed in the formal theory of partial
differential equations and the tools developed there can now through the
maximum entropy principle be brought into the area of nonparametric
statistics.

\subsubsection{Thermodynamic pressure and its cousins}\lbl{pressure}
We will now investigate an important consequence of the fundamental variational identity \rf{FundamentalVariationalIdentity} for the thermodynamic case when the total energy is the only observable. For this special case the variational identity turns into
 \begin{align}
\left< \delta_\Gamma\mathcal{ H} \right>&= -kT\;\delta_\Gamma \ln Z.\lbl{FundamentalVariationalIdentity1}
  \end{align}
  The force acting at a boundary point $\vb{x}$ can, taking into account the fact that the state of the system is determined by the position and momenta of all the $n$ particles comprising the system, be given by a function
  \begin{align}
  \vb{F}&=\vb{F}(\vb{p},\vb{q},\vb{x}),\;\;\;\;\;\vb{x}\in\Gamma,\;(\vb{p},\vb{q})\in\mathbb{R}^{6n}.\lbl{BoundaryForce}
  \end{align}
  A small deformation of the boundary is determined by an infinitesimal deformation vector  field $d\vb{r}_{\Gamma}$ defined on the boundary $\Gamma$. The change in total energy induced by this deformation is given by
  \begin{align}
\delta_{\Gamma}\mathcal{H}(\vb{p},\vb{q};\Gamma) &=-\int_{\Gamma}\;d\vb{x}\;\vb{F}(\vb{p},\vb{q},\vb{x})\cdot d\vb{r}_{\Gamma}.
  \end{align}
  The fundamental variational identity \rf{FundamentalVariationalIdentity} now gives
 \begin{align}
\int_{\Gamma}\;d\vb{x}\;\left<\vb{F}(\vb{p},\vb{q},\vb{x}) \right>\cdot d\vb{r}_{\Gamma}&= kT\;\delta_\Gamma \ln Z.\lbl{FundamentalVariationalIdentity2}
  \end{align}  
  We will in the following only consider the common situation defined by
  \begin{align}
  \left<\vb{F}(\vb{p},\vb{q},\vb{x}) \right>&=p(\vb{x})\;\vb{n},\lbl{PressureAssumption}
  \end{align}
 where $\vb{n}$ is the unit normal for the surface $\Gamma$.  Note that by definition, $p$ is now the pressure for the system. For this case \rf{FundamentalVariationalIdentity2} turns into the identity
 \begin{align}
\int_{\Gamma}\;d\vb{x}\;p(\vb{x})\;\vb{n}\cdot d\vb{r}_{\Gamma}&= kT\;\delta_\Gamma\ln Z.\lbl{FundamentalVariationalIdentity2}
\end{align}
Let us first consider the case of a smooth surface, and for this kind of surface,  let us consider an infinitesimal variation of the surface that is a pure expansion or contraction. This means that $d\vb{r}_{\Gamma}=\vb{n}ds$. For this kind of variation the fundamental variational identity \rf{FundamentalVariationalIdentity2} takes the form
 \begin{align}
ds\;\int_{\Gamma}\;d\vb{x}\;p(\vb{x})&= kT\;\delta_\Gamma\ln Z.\lbl{FundamentalVariationalIdentity3}
\end{align}
Using the fact that the volume spanned by the deformation is $d_{\Gamma}V=A(\Gamma)ds$, where $A(\Gamma)$ is the area of the surface, we have 
\begin{align}
\left<p\right>_{\Gamma}= kT\;\frac{\delta_\Gamma\ln Z}{d_{\Gamma}V},\lbl{FundamentalVariationalIdentity3.5}
\end{align}
where $\left<p\right>_{\Gamma}$ is the average of the pressure over the surface of the cavity.

For some important cases the partition function depends on the surface only through the volume. For this situation we have
\begin{align}
\delta_\Gamma\ln Z=\frac{\partial\ln Z}{\partial\;V}d_{\Gamma}V,\nonumber
\end{align}
so that 
 \begin{align}
\left<p\right>_{\Gamma}= kT\;\frac{\partial\ln Z}{\partial\;V}.\lbl{FundamentalVariationalIdentity4}
\end{align}
This is the standard formula for the thermodynamic pressure that one finds in any textbook. It is very frequently true that, independently of the shape, the partition function for large cavities depends only on the volume of the cavity. This may however not be the case for smaller cavities and for such cases we must retreat to the more general identity \rf{FundamentalVariationalIdentity3.5}. It is easy to verify that the indentity \rf{FundamentalVariationalIdentity4} holds for any surface smooth or not. If the deformation is a pure expansion or contraction of a part of the surface defined by $\Gamma_0\subset\Gamma$, we also get the identity \rf{FundamentalVariationalIdentity4}, but now with $\Gamma\rightarrow\Gamma_0$.

The general variational identity for the entropy \rf{dS5} takes for the particular case discussed in this section the form
\begin{align}
\delta S&=\frac{1}{kT}d\left<\mathcal{H}\right>-\frac{1}{kT}\left< \delta_\Gamma\mathcal{ H} \right>\nonumber\\
&=\frac{1}{kT}d\left<\mathcal{H}\right>+\delta_\Gamma \ln Z\nonumber\\
&=\frac{1}{kT}d\left<\mathcal{H}\right>+\frac{1}{kT}\left<p\right>_{\Gamma}d_{\Gamma}V,\nonumber
\end{align}
which can be rewritten as
\begin{align}
kTdS=dE+\left<p\right>_{\Gamma}d_{\Gamma}V,\lbl{dS6}
\end{align}
where we have used \rf{FundamentalVariationalIdentity3.5} and where now $E=\left<\mathcal{H}\right>$ is the energy of the system. We recognize \rf{dS6} as one of the fundamental formulas from conventional thermodynamics.

In this section we have done the derivation of the formulas for the Thermodynamical pressure for the case of a classical system. However, the derivation of the pressure formula for the case of quantum systems leads the exact same formulas. If there are more observables in addition to the energy, for example total momentum and/or total angular momentum, the pressure formulas must be generalized. The derivation of the generalizations follow the pattern laid down in this section.

\subsection{The maximum entropy principle for quantum systems}
We have in the first section of these notes introduced the thermodynamical formalism in the context of classical physics and classical observables. It involved a state space that was finite or at least finite dimensional, and the challenge was to determine which probability distribution on the state space is the best to use, given the means of a finite number of observables of the system. The solution to this problem was to choose the probability distribution that maximized the entropy functional \rf{78} under the constraints determined by the given means. At the face of it, for a quantum system, the situation appears to be very different. For this case the state space is an infinite dimensional Hilbert space and the full information that an observer has is encoded in the \ttx{density operator} for the system. This is a self-adjoint positive operator $\hat\rho$, on the Hilbert space with trace equal to one.
\begin{align}
\mathrm{Tr}(\hat\rho)=1.\nonumber
\end{align}
The expectation value of any quantum observable, $\hat A$ is by definition
\begin{align}
\left<\hat A\right>=Tr(\hat\rho\hat A).\lbl{QuantumMean}
\end{align} 
The question one poses is which density operator should be used if we only know the expectation value of a finite number of quantum observables $\hat A_i$.
\begin{align}
\left<\hat A_i\right>&=a_i,\;\;\;i=1,2,...p.\lbl{Constraints}
\end{align}
Even though there are real differences between the classical and the quantum case, much is also the same. 

The analog of the Gibbs entropy measure \rf{79} is the Von Neumann entropy measure for density operators given by
\begin{align}
S(\hat\rho)=-Tr(\hat\rho\ln\hat\rho).\lbl{NeumanEntropy}
\end{align}
 The solution to the question posed on the previous page proposed by the Maximum entropy method is to choose the density operator that maximize the Von Neumann entropy while satisfying the constraints \rf{Constraints}. It is a simple exercise to show that the unique solution to this maximization problem is the following density operator
 \begin{align}
 \hat\rho=\frac{1}{Z(\lambda_1,...,\lambda_p)}\exp(-\sum_i\;\lambda_i\hat A_i).\lbl{MaxEntDensityOperator}
 \end{align}
 This operator is the \ttx{maximum entropy density operator}. The function $Z$ is the partition function and is given by
 \begin{align}
 Z(\lambda_1,...,\lambda_p)=Tr\left\{\exp(-\sum_i\;\lambda_i\hat A_i)\right\}.\lbl{PartitionFunction}
 \end{align}
 Arguing like in the classical case we find that much of the thermodynamic formalism is the same as before. Specifically we have
 \begin{align}
 S&=\ln Z+\sum_{i=1}^p\;\lambda_i\left<\hat A_i\right>,\nonumber\\
 \left<\hat A_i\right>&=-\frac{\partial \ln Z}{\partial \lambda_i},\nonumber\\
 \lambda_i&=\frac{\partial S}{\partial \left<\hat A_i\right>},\\
\sum_i\lambda_i \left<\frac{\partial \hat A_i}{\partial \alpha}\right>&=-\frac{\partial \ln Z}{\partial \alpha},\lbl{QThermodynamicFormalism}
\end{align}
where in the last identity we assume that all the observables depend on some parameter $\alpha$. Like for the classical case, it is common that the system of interest  is restricted to some spatial domain $D$. In this case we get a quantum analog to the classical formula \rf{FundamentalVariationalIdentity}
 \begin{align}
 \delta_\Gamma\ln Z=-\mathlarger{\sum}_j \; \lambda_j \;\left< \delta_\Gamma \hat A_j\right>.\lbl{QuantumFundamentalVariationalIdentity}
  \end{align}
  As we can see, much of the thermodynamic formalism is the same for the classical and the quantum case. However some things are different, or they at least appear to be different. In the classical case we can find correlations between different observables by computing mixed partial derivatives of the partition function as shown in \rf{215} and \rf{216}. In the quantum case this is problematic  unless the operators commute. For the case of two observables $\hat A$ and $\hat B$ we have for example
  \begin{align}
  \frac{\partial^2 Z}{\partial_{\mu}\partial_{\lambda}}= \frac{\partial^2 Z}{\partial_{\lambda}\partial_{\mu}}=2\left<\hat A\hat B+\hat B\hat A\right>-\left<\hat A\right>\left<\hat B\right>.\nonumber
  \end{align}
 In a sense this should not come as a surprise. The reason for this is that $\hat A\hat B$ is not in general Hermitian even if both $\hat A$ and $\hat B$ are. Thus $\left<\hat A\hat B\right>$ is not something that you could ever measure, so it does not make sense to try to predict it. However, $\hat A\hat B+\hat B\hat A$ $\text{\it is}$ a Hermitian operator and thus $\left<\hat A\hat B+\hat B\hat A\right>$ is something one would want to predict. And,  this is exactly what you would be able to predict using the thermodynamical formalism.
  
  There is another important way in which the classical and quantum cases are different; namely the question of how to include prior information about the system into the thermodynamical formalism. In the classical case this was taken care of by using the entropy measure \rf{79} that included the prior distribution $\rho_0$. We have seen how maximization the entropy in the context of statistical mechanics leads to the distribution 
\begin{align}
\rho (\vx_1,...,\vx_n, \vb{p}_1,..., \vb{p}_n) = \frac{\rho_0}{Z} \; \exp(-\frac{\mathcal{H}}{kT}), \lbl{217} 
\end{align}
which we recognized to be the Gibbs canonical ensemble. This is however not entirely correct, the canonical ensemble corresponds to the case when we have a uniform prior. We could have gotten this distribution directly from a maximization of the functional
\begin{align}
S(\rho) = - \int_{\mathbf{R}^n} dV \; \rho \; \ln (\rho). \lbl{GibbsEntropy} 
\end{align}
This is in fact the functional used by Gibbs in his foundational studies of thermodynamics. It is this functional that is called the Gibbs entropy measure. The functional we introduced in \rf{79} measure by definition the relative entropy of $\rho$ with respect to $\rho_0$. It is also denoted by other names in the research literature.

  The Von Neumann entropy introduced in \rf{NeumanEntropy} is the direct analog of the Gibbs entropy measure \rf{GibbsEntropy}. However, in contrast to the classical case, there does not exists a measure of information for quantum systems that is universally recognized to be the best measure to use in all cases where there is  prior information to take into account. Many different kind of measures has been studied by physicists and mathematicians over the years. Today these questions are intensely pursued in the topical field of quantum information theory. 

Note that the classical relative entropy measure can be written in the form
\begin{align}
S(\rho|\;\rho_0) = - \int_{\mathbf{R}^n} dV \; \rho \; (\ln (\rho-\ln\rho_0)). \lbl{RelativeEntropy} 
\end{align}
One approach to a quantum version of relative entropy is to directly generalize \rf{RelativeEntropy} to the quantum case. This gives us for any pair of density operators $\hat\rho$ and $\hat\rho_0$ the relative quantum entropy in the form
\begin{align}
S_Q(\hat \rho|\;\hat\rho_0) = - Tr( \hat \rho \; (\ln \hat \rho-\ln\hat\rho_0)). \lbl{QRealtiveEntropy} 
\end{align}
One could now guess that the corresponding maximal entropy distribution, when the only observable is the total energy, will take the form
\begin{align}
\hat\rho = \frac{\hat\rho_0}{Z} \; \exp(-\frac{\mathcal{\hat H}}{kT}), \lbl{QMaxEntDensityOperator} 
\end{align}
However, a formula like this can not possibly  be correct because the right hand side of \rf{QMaxEntDensityOperator} is not a Hermitian operator unless the prior $\hat\rho_0$ commutes with the total energy operator $\mathcal{\hat H}$.

  Here we will not pursue these matters, but we should keep in mind that there is a real issue here concerning the general validity of the thermodynamical formalism for quantum systems in cases when there is prior information present.

\nocite{Landau}
\nocite{Whitham}
\nocite{Bridgman}
\nocite{Sonin}
\nocite{Jones}

\pagebreak

\end{document}